\documentclass[a4paper, 10pt,openright, twoside,oldfontcommands]{memoir}
\usepackage[english]{babel}

\usepackage{fancyhdr}
\usepackage{setspace}

\usepackage{amsmath}
\usepackage{graphicx,epsfig,color}
\usepackage[color]{xypic}
\xyoption{rotate}
\input xy
\xyoption{all}
\xyoption{2cell}
\UseAllTwocells
\usepackage{verbatim}
\usepackage{amssymb}
\usepackage{eucal}
\usepackage{enumerate}
\usepackage[shortlabels]{enumitem}
\usepackage{multicol}
\usepackage{manfnt}
\usepackage{inputenc}
\usepackage{amsthm}
\usepackage{amsfonts}
\usepackage{mathrsfs}
\usepackage[all,2cell]{xy}
\usepackage{float}
\usepackage{fancybox}
\usepackage{a4wide}
\usepackage{hyperref}
\usepackage{wrapfig}
\usepackage{wasysym}

\usepackage{stmaryrd} 

\usepackage{fixltx2e}
\MakeRobust{\overrightarrow}
\MakeRobust{\overleftarrow}

\newcommand{\Mat}[1]{\mathsf{Mat}_{\scriptscriptstyle #1}}

\mathcode`:="003A  
\mathcode`;="003B  
\mathcode`?="003F  
\mathcode`|="026A  

\mathchardef\ls="213C    
\mathchardef\gr="213E    
\mathchardef\uparrow="0222  
\mathchardef\downarrow="0223  

\newcommand{\pippo}[1]{\marginpar{{\bf Trotski:} #1 }}    
\newcommand{\fabio}[1]{\marginpar{{\bf Lenin:} #1 }}    

\newtheoremstyle{italic}
  {3pt}
  {3pt}
  {\normalfont}
  {}
  {\scshape}
  {.}
  { }
  {}

\theoremstyle{plain}

\newtheorem{theorem}{Theorem}[chapter]
\newtheorem{proposition}[theorem]{Proposition}
\newtheorem{corollary}[theorem]{Corollary}
\newtheorem{lemma}[theorem]{Lemma}

\theoremstyle{definition}

\newtheorem{convention}[theorem]{Convention}
\newtheorem{definition}[theorem]{Definition}
\newtheorem{remark}[theorem]{Remark}

\newtheorem{example}[theorem]{Example}

\theoremstyle{italic}

\renewenvironment{proof}[1][Proof]{\begin{trivlist}
\item[\hskip \labelsep {\bfseries #1}]}{\hfill $\square$\end{trivlist}}

\newenvironment{todo}{\begin{itemize}\item \footnotesize \scshape}{\end{itemize}}
\newenvironment{Iff-RL}{\textbf{($\Rightarrow$)} }{\bigskip}
\newenvironment{Iff-LR}{\textbf{($\Leftarrow$)} }{}

\newcommand{\m}{\mathit}
\def \: {\colon}

\newcommand{\mf}{\mathbf}

\newcommand\id{\m{id}}

\newcommand\funF{\mathcal{F}}

\newcommand\G{\mathcal{G}}
\newcommand\D{\mathcal{D}}

\newcommand\FH{\mathcal{H}}
\newcommand\T{\mathcal{T}}
\newcommand\J{\mathcal{J}}

\def \catC {\mathbb{C}}
\def \Set  {\mathbf{Set}}
\def \Cat {\mathbf{Cat}}
\def \Rel  {mathbb{Rel}}

\def \N {\mathbb{N}}
\def \Z {\mathbb{Z}}
\def \Bool {\mathsf{Bool}}
\def \T {\mathbb{T}} 
\def \PS {\mathbb{S}} 
\def \PR {\mathbb{R}} 

\def \PROP {\mathbf{PROP}} 
\def \PRO {\mathbf{PRO}} 
\newcommand{\Mod}[1]{\mf{Mod}(#1)} 
\def \Perm {\mathsf{P}} 
\def \F {\mathsf{F}}
\def \Fop {\F^{\m{op}}}

\def \AB {\mathbb{AB}}

\newcommand{\Span}[1]{\mathsf{Span}(#1)}
\newcommand{\Cospan}[1]{\mathsf{Cospan}(#1)}
\newcommand{\CospanO}{\mathsf{Cospan}}

\def \wmon {\mf{Mon}^w}

\def \bcom {\mf{Com}^b}

\def \tns {\oplus}

\newcommand{\eql}[1]{\ \overset{\text{#1}}=\ }

\def \SV {\mathsf{SV}} 
\newcommand{\sem}[1]{\mathcal{Sem}_{#1}} 

\newcommand{\strsem}[1]{\langle\!\langle #1 \rangle\!\rangle} 
\newcommand{\strsemO}{\strsem{\cdot}} 
\newcommand{\strsemHA}[1]{\overrightarrow{\strsem{#1}}} 
\newcommand{\strsemHAO}{\strsemHA{\cdot} } 

\newcommand{\initVect}{!} 

\DeclareMathSymbol{\reversedExclMark}{\mathord}{operators}{"3C}
\newcommand{\finVect}{\reversedExclMark} 

\newcommand\Gmult{\lower5pt\hbox{$\includegraphics[width=20pt]{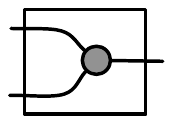}$}}
\newcommand\Gcomult{\lower5pt\hbox{$\includegraphics[width=20pt]{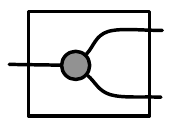}$}}
\newcommand\Gunit{\lower5pt\hbox{$\includegraphics[width=16pt]{graffles/Gunit.pdf}$}}
\newcommand\Gcounit{\lower5pt\hbox{$\includegraphics[width=16pt]{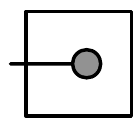}$}}
\newcommand\twoGcounit{\lower5pt\hbox{$\includegraphics[width=16pt]{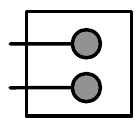}$}}

\newcommand\Bmult{\lower5pt\hbox{$\includegraphics[width=20pt]{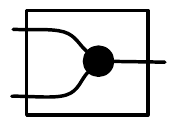}$}}
\newcommand\Bcomult{\lower5pt\hbox{$\includegraphics[width=20pt]{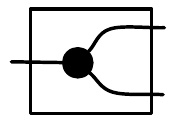}$}}
\newcommand\Bunit{\lower5pt\hbox{$\includegraphics[width=16pt]{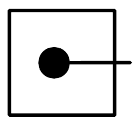}$}}
\newcommand\Bcounit{\lower5pt\hbox{$\includegraphics[width=16pt]{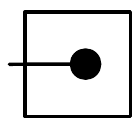}$}}

\newcommand\Wmult{\lower5pt\hbox{$\includegraphics[width=20pt]{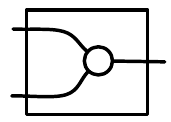}$}}
\newcommand\Wcomult{\lower5pt\hbox{$\includegraphics[width=20pt]{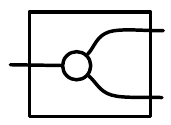}$}}
\newcommand\Wunit{\lower5pt\hbox{$\includegraphics[width=16pt]{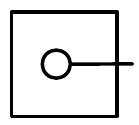}$}}
\newcommand\Wcounit{\lower5pt\hbox{$\includegraphics[width=16pt]{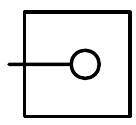}$}}

\newcommand\ASepUno{\lower3pt\hbox{$\includegraphics[width=28pt]{graffles/ASep1.pdf}$}}
\newcommand\ASepDue{\lower3pt\hbox{$\includegraphics[width=24pt]{graffles/ASep2.pdf}$}}
\newcommand\Idnet{\lower3pt\hbox{$\includegraphics[width=20pt]{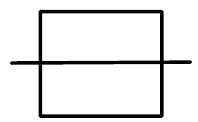}$}}
\newcommand\symNet{\lower3pt\hbox{$\includegraphics[width=20pt]{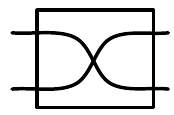}$}}

\newcommand{\pullbackcorner}[1][d]{\save*!/#1+1.0pc/#1:(1,-1)@^{|-}\restore}
\newcommand{\pushoutcorner}[1][d]{\save*!/#1-1.0pc/#1:(-1,1)@^{|-}\restore}
\newcommand\BSepOne{\lower3pt\hbox{\includegraphics[width=30pt]{graffles/BlackSep1.pdf}}}
\newcommand\BlackX{\lower5pt\hbox{\includegraphics[width=30pt]{graffles/BlackX.pdf}}}
\newcommand\WSep{\lower5pt\hbox{\includegraphics[width=35pt]{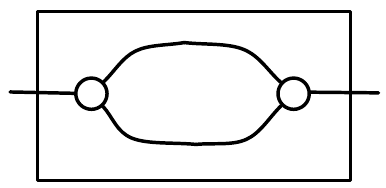}}}

\newcommand{\matrixOnebis}{1 }

\newcommand{\matrixZerobis}{0}

\newcommand{\matrixOneOne}{{\scriptsize\left(%
                \begin{array}{c}
                 \!\! 1 \!\!\\
                 \!\! 1 \!\!
                \end{array}\right)}}

\newcommand{\matrixOneZero}{{\scriptsize\left(%
                \begin{array}{c}
                 \!\! 1 \!\!\\
                 \!\! 0 \!\!
                \end{array}\right)}}

\newcommand{\matrixZeroOne}{{\scriptsize\left(%
                \begin{array}{c}
                 \!\! 0 \!\!\\
                 \!\! 1 \!\!
                \end{array}\right)}}

\newcommand{\matrixOneOneFlat}{\left(%
                \begin{array}{cc}
                 \!\! 1 \!&\! 1\!\!
                \end{array}\right)}

\newcommand{\matrixOne}{\left(%
                \begin{array}{c}
                 \!\!\! 1 \!\!\!
                \end{array}\right)}

\newcommand{\matrixZero}{\left(%
                \begin{array}{c}
                \!\!\!  0 \!\!\!
                \end{array}\right)}

\newcommand{\matrixK}{\left(%
                \begin{array}{c}
                \!\!\!  k \!\!\!
                \end{array}\right)}

\newcommand{\matrixNull}{\left(%
                \begin{array}{c}
                \!\!\!\!\!
                \end{array}\right)}

\newcommand{\tr}[1]{\xrightarrow{#1}}    
\newcommand{\tl}[1]{\xleftarrow{#1}}    

\newcommand{\op}{\mathit{op}}

\def \PID {\mathsf{R}}
\def \poly {\mathsf{k}[x]} 
\def \frpoly {\mathsf{k}(x)} 
\def \laur {\mathsf{k}((x))} 
\def \fps {\mathsf{k}[[x]]} 
\def \field {\mathsf{k}} 
\def \frPID {\mathsf{k}} 
\def \ratio {\mathsf{k}\langle x\rangle} 

\newcommand \LSum[1] {\sum_{i= #1}^{\infty}}

\newcommand \wstream[1] {\widehat{#1}}
\newcommand \wlrn[1] {\widetilde{#1}}

\newcommand \stream[1] {\hat{#1}}
\newcommand \lrn[1] {\tilde{#1}}

\def \embpolyfps {\stream{\cdot}}
\def \embpolyfrpoly {\delta}
\def \embfpsfls {\nu}
\def \embfrpolyfls {\lrn{\cdot}}

\newcommand \PROPRing[1] {\mathbb{#1}}
\def \PROPR {\PROPRing{R}} 

\newcommand \HA[1] {\mathbb{HA}_{\scriptscriptstyle #1}}

\def \ABR {\HA{\PID}}
\def \ABRop {\ABR^{\m{op}}} %
\def \ABpoly {\HA{\poly}}
\def \ABpolyop {\HA{\poly}^{\m{op}}}

\newcommand \IH[1]{\mathbb{IH}_{\scriptscriptstyle #1}}
\def \IBR {\IH{\PID}} 
\def \IBpoly {\IH{\poly}}

\def \Matpoly {\Mat{\poly}}
\def \Matpolyop {\Mat{\poly}^{\m{op}}}

\def \Matfrpoly {\Mat{\frpoly}}

\def \Matfps {\Mat{\fps}}
\def \Matfpsop {\Mat{\fps}^{\m{op}}}

\def \Matlaur {\Mat{\laur}}

\def \Matratio {\Mat{\ratio}}

\def \VectR {\Mat{\PID}}
\def \VectRop {\Mat{\PID}^{\m{op}}}

\newcommand \Subspace[1]{\mathsf{SV}_{\scriptscriptstyle #1}}
\def \SVR {\Subspace{\frPID}}
\def \SVpoly {\Subspace{\frpoly}}
\def \SVfps {\Subspace{\laur}}

\def \Mon {\mathbb{M}}
\def \Com {\mathbb{C}}
\def \wmon {\Mon}

\def \bcom {\Com}
\def \SupSpan {\mathsf{Sp}}
\def \SupCospan {\mathsf{Cp}}
\def \IBRw {\IBR^{\scriptscriptstyle \SupSpan}} 
\def \IBRb {\IBR^{\scriptscriptstyle \SupCospan}} 
\def \To {\Rightarrow}
\DeclareMathOperator{\RModule}{\mathsf{Mod}}
\DeclareMathOperator{\FRModule}{\mathsf{FMod}}
\DeclareMathOperator{\Kernel}{Ker}

\newcommand \RMod[1] {\RModule_{\scriptscriptstyle #1}} 
\newcommand \FRMod[1] {\FRModule_{\scriptscriptstyle #1}} 
\newcommand \Ker[1] {\Kernel(#1)} 
\newcommand \minusB {\scalebox{0.75}[1.0]{\( - \)}B}
\newcommand \KerAB {\Ker{A \mid \minusB}}
\def \VectSpFr {\FRMod{\frPID}}
\def \poi {\,\ensuremath{;}\,}
\newcommand \coc[1] {{#1}^{\star}}
\newcommand \refl[1] {{#1}^{R}}
\def \df {\ \ensuremath{:\!\!=}\ }
\def \dfop {\ \ensuremath{=\!\!:}\ }
\def \minus {\ensuremath{\!-\!\!}}
\newcommand \restr[2] {#1_{\upharpoonright #2}} 
\newcommand \tra {\ensuremath{\mathcal{T}}} 
\newcommand \pn {\ensuremath{\mathcal{N}}} 

\newcommand\scalar{\!\lower5pt\hbox{$\includegraphics[width=20pt]{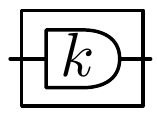}$}\!}
\newcommand\coscalar{\!\lower5pt\hbox{$\includegraphics[width=20pt]{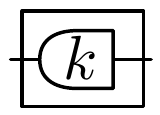}$}\!}
\newcommand\unitscalar{\lower3pt\hbox{$\includegraphics[width=24pt]{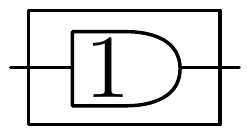}$}}
\newcommand\scalarminusone{\lower4pt\hbox{$\includegraphics[width=30pt]{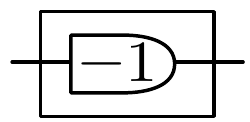}$}}
\newcommand\scalarminusoneop{\lower4pt\hbox{$\includegraphics[width=30pt]{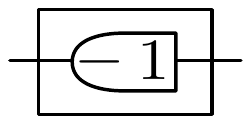}$}}
\newcommand\antipode{\lower4pt\hbox{$\includegraphics[width=20pt]{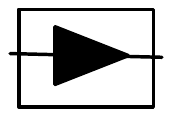}$}}
\newcommand\antipodeop{\lower4pt\hbox{$\includegraphics[width=20pt]{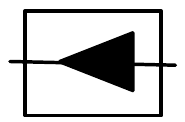}$}}
\newcommand\antipodesquare{\lower3pt\hbox{$\includegraphics[width=22pt]{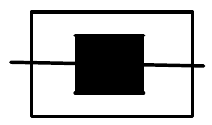}$}}
\newcommand\circuitAdots{\lower6pt\hbox{$\includegraphics[width=30pt]{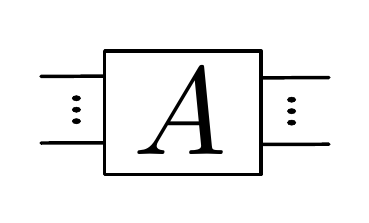}$}}

\newcommand\wcounitn{\lower5pt\hbox{$\includegraphics[width=25pt]{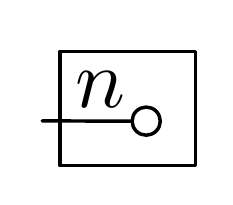}$}}
\newcommand\bcounitn{\lower5pt\hbox{$\includegraphics[width=25pt]{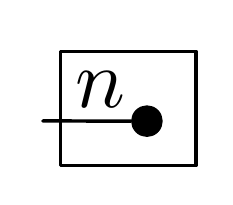}$}}
\newcommand\lccn{\lower5pt\hbox{$\includegraphics[width=25pt]{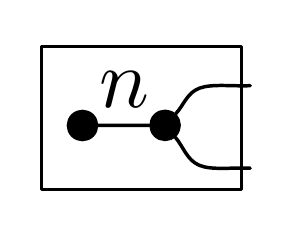}$}}
\newcommand\rccn{\lower5pt\hbox{$\includegraphics[width=25pt]{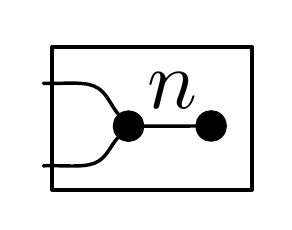}$}}
\newcommand\idncircuit{\lower5pt\hbox{$\includegraphics[width=25pt]{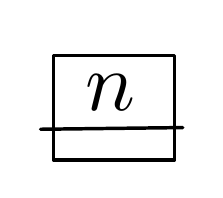}$}}
\newcommand\circuitrbcounits{\lower5pt\hbox{$\includegraphics[width=25pt]{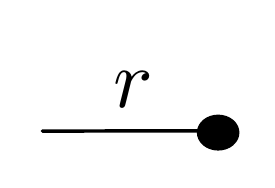}$}}
\newcommand\lccB{\lower5pt\hbox{$\includegraphics[width=25pt]{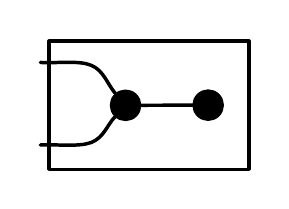}$}}
\newcommand\rccB{\lower5pt\hbox{$\includegraphics[width=25pt]{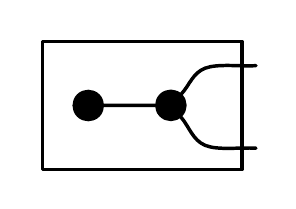}$}}
\newcommand\IdBcounitc{\lower5pt\hbox{$\includegraphics[width=20pt]{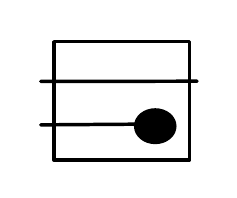}$}}
\newcommand\BcounitId{\lower5pt\hbox{$\includegraphics[width=20pt]{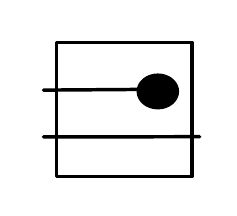}$}}
\newcommand\symNetTwoOne{\lower7pt\hbox{$\includegraphics[width=25pt]{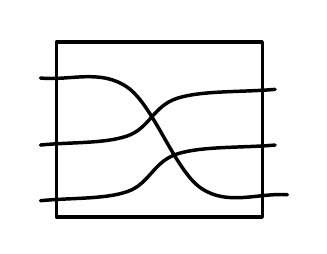}$}}
\newcommand\nscalar{\!\!\lower6pt\hbox{$\includegraphics[width=35pt]{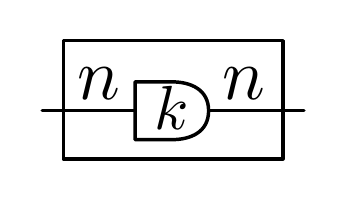}$}\!\!}
\newcommand\Wmultstar{\!\lower5pt\hbox{$\includegraphics[width=20pt]{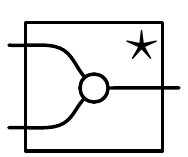}$}\!}
\newcommand\scalarstar{\!\!\lower7pt\hbox{$\includegraphics[width=35pt]{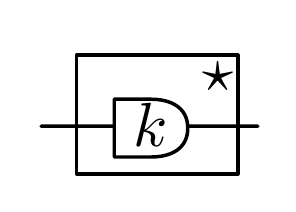}$}\!\!}
\newcommand\twoBcounit{\lower5pt\hbox{$\includegraphics[width=15pt]{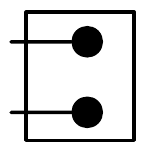}$}}

\newcommand\delay{\!\lower6pt\hbox{$\includegraphics[width=25pt]{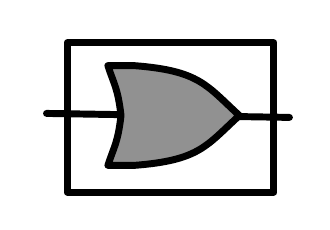}$}\!}
\newcommand\circuitX{\!\lower4pt\hbox{$\includegraphics[width=22pt]{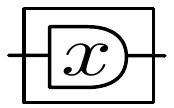}$}\!}
\newcommand\scalarp{\!\lower4pt\hbox{$\includegraphics[width=22pt]{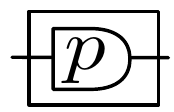}$}\!}
\newcommand\scalarpstar{\!\lower3pt\hbox{$\includegraphics[width=25pt]{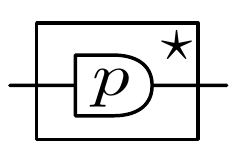}$}\!}
\newcommand\scalarpop{\!\lower4pt\hbox{$\includegraphics[width=22pt]{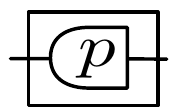}$}\!}
\newcommand\nscalarp{\!\lower3pt\hbox{$\includegraphics[width=30pt]{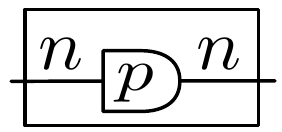}$}\!}
\newcommand\circuitfibr{\!\lower3pt\hbox{$\includegraphics[width=44pt]{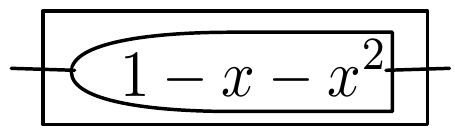}$}\!}
\newcommand\circuitXop{\!\lower4pt\hbox{$\includegraphics[width=24pt]{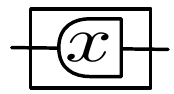}$}\!}
\newcommand\ncircuitX{\!\lower3pt\hbox{$\includegraphics[width=30pt]{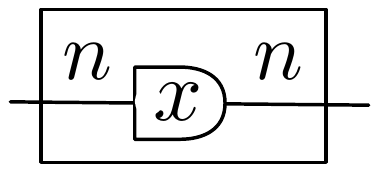}$}\!}
\newcommand\scalarpone{\!\lower3pt\hbox{$\includegraphics[width=22pt]{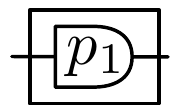}$}\!}
\newcommand\scalarptwoop{\!\lower3pt\hbox{$\includegraphics[width=22pt]{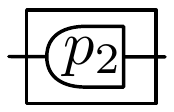}$}\!}

\newcommand{\xx}{\vectorfont{x}}
\newcommand{\yy}{\vectorfont{y}}
\newcommand{\zz}{\vectorfont{z}}
\newcommand{\uu}{\vectorfont{u}}
\newcommand{\vv}{\vectorfont{v}}
\newcommand{\ww}{\vectorfont{w}}
\newcommand{\ee}{\vectorfont{e}} 

\newcommand{\IdnetT}{\Idnet}
\newcommand{\symNetT}{\symNet}
\newcommand{\WunitT}{\Wunit}
\newcommand{\BcounitT}{\Bcounit}
 \newcommand{\BcomultT}{\Bcomult}
\newcommand{\WmultT}{\Wmult}
\newcommand{\scalarT}{\scalar}
\newcommand{\circuitXT}{\circuitX}
\newcommand{\BunitT}{\Bunit}
\newcommand{\WcounitT}{\Wcounit}
\newcommand{\WcomultT}{\Wcomult}
\newcommand{\BmultT}{\Bmult}
\newcommand{\scalaropT}{\coscalar}
\newcommand{\circuitXopT}{\circuitXop}
\newcommand{\ZeronetT}{\lower4pt\hbox{$\includegraphics[width=14pt]{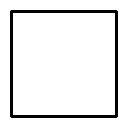}$}}

\newcommand{\circuitminusone}{\!\lower5pt\hbox{$\includegraphics[width=20pt,height=15pt]{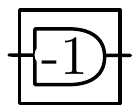}$}\!}
\newcommand{\circuitminusoneop}{\!\lower5pt\hbox{$\includegraphics[width=20pt,height=15pt]{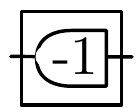}$}\!}

\newcommand{\cospanXTzero}{\circuitX^{\labelSep 0} ; \circuitXop^{\labelSep 0}}
\newcommand{\cospanXTkl}{\circuitX^{\labelSep k} ; \circuitXop^{\labelSep l}}

\newcommand{\spanXTzero}{\circuitXop^{\labelSep 0} ; \circuitX^{\labelSep 0}}
\newcommand{\spanXTk}{\circuitXop^{\labelSep k_1} ; \circuitX^{ \labelSep k_1} }
\newcommand{\spanXTktwo}{\circuitXop^{\labelSep k_2} ; \circuitX^{ \labelSep k_2} }
\newcommand{\spanXTkthree}{\circuitXop^{\labelSep k_3} ; \circuitX^{ \labelSep k_3} }

\newcommand{\typ}{\mathrel{:}}

\newcommand{\lbbd}{\mathopen{[\![}}
\newcommand{\rbbd}{\mathclose{]\!]}}
\newcommand{\dsem}[1]{\lbbd{#1}\rbbd}
\newcommand{\dsemO}{\dsem{\cdot}}
\newcommand{\dsemHA}[1]{\overrightarrow{\dsem{#1}}}
\newcommand{\dsemHAO}{\dsemHA{\cdot}}
\newcommand{\dsemHAop}[1]{\overleftarrow{\dsem{#1}}}
\newcommand{\dsemHAOop}{\dsemHAop{\cdot}}

\newcommand{\lbbo}{\mathopen{\langle}}
\newcommand{\rbbo}{\mathclose{\rangle}}
\newcommand{\osem}[1]{\lbbo #1 \rbbo}
\newcommand{\osemO}{\osem{\cdot}}

\newcommand{\pair}[2]{(#1,#2)}

\newcommand{\vectorfont}[1]{\mathbf{#1}}
\newcommand{\vlist}[1]{\overrightarrow{#1}}

\newcommand\idzcircuit{\lower5pt\hbox{$\includegraphics[width=20pt]{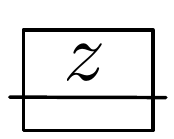}$}}
\newcommand{\circuitXspan}{\!\lower4pt\hbox{$\includegraphics[width=40pt]{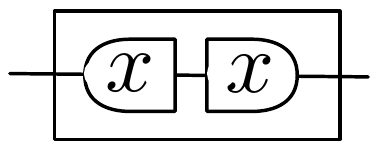}$}\!}
\newcommand{\circuitXcospan}{\!\lower5pt\hbox{$\includegraphics[width=40pt]{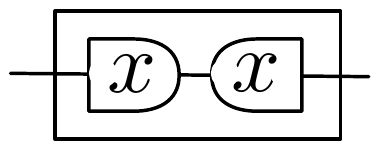}$}\!}
\newcommand\zeroscalar{\lower3pt\hbox{$\includegraphics[width=20pt]{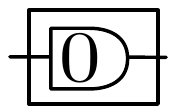}$}}
\newcommand\zeroscalarr{\lower3pt\hbox{$\includegraphics[width=25pt]{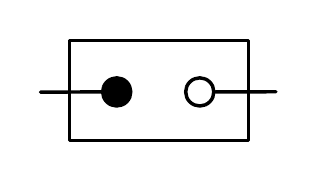}$}}
\newcommand\rationalcircuit{\lower5pt\hbox{$\includegraphics[width=35pt]{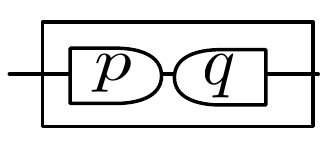}$}}
\newcommand\Wccl{\lower5pt\hbox{$\includegraphics[width=22pt]{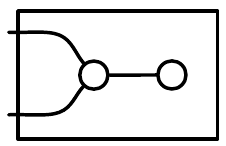}$}}
\newcommand{\circuitkkop}{\!\lower4pt\hbox{$\includegraphics[width=32pt]{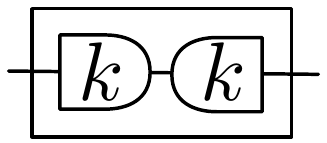}$}\!}

\usepackage{centernot}

\newcommand{\Impl}[1]{\ \overset{\text{#1}}\Rightarrow\ }

\newcommand{\mn}{{-}}

\newcommand{\dtransF}[2]{\underset{#2}{\scriptstyle\xrightarrow{#1}}}

\newcommand\circuitUnoMinusX{\lower4.5pt\hbox{$\includegraphics[width=32pt]{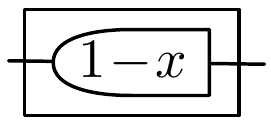}$}}
\newcommand\circuitUnoMinusXSquare{\lower5pt\hbox{$\includegraphics[width=35pt]{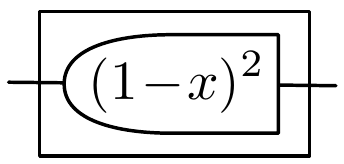}$}}

\usepackage{prooftree}

\newcommand{\ruleLabel}[1]{#1}
\newcommand{\labelSep}{\,}

\newcommand{\bnfEq}{\; ::= \;}
\newcommand{\bnfSep}{\;\; | \;\;}

\makeatletter
\def\moverlay{\mathpalette\mov@rlay}
\def\mov@rlay#1#2{\leavevmode\vtop{%
\baselineskip\z@skip \lineskiplimit-\maxdimen
\ialign{\hfil$#1##$\hfil\cr#2\crcr}}}
\makeatother

\newcommand{\typeJudgment}[3]{{ {#2} \,\typ\, {#3}}}
\newcommand{\sort}[2]{\ensuremath{(#1,\,#2)}}

 \newcommand{\derivationRule}[3]{{\prooftree{ #1}\justifies{ #2}\using\ruleLabel{#3}\endprooftree}}
\newcommand{\reductionRule}[2]{{\prooftree{\scriptstyle #1}\justifies{\scriptstyle #2}\endprooftree}}

\newcommand\twarr[2]{%
\mathrel{\mathop{\moverlay{\scriptstyle\xrightarrow{\,#1\,}\cr{\lower.2em\hbox{$\scriptstyle{}_{#2}$}}}}}}
\newcommand\twarrw[2]{%
\mathrel{\mathop{\moverlay{\scriptstyle\Longrightarrow\cr{\lower-.6em\hbox{$\scriptstyle{}_{#1}$}}
\cr{\lower.3em\hbox{$\scriptstyle{}_{#2}$}}}}}}

\newcommand{\dtrans}[2]{\hbox{$\;\twarr{#1}{#2}\;$}}
\newcommand{\dtransw}[2]{\raise1pt\hbox{$\;\twarrw{#1}{#2}\;$}}

\newcommand{\ftr}[1]{ft(#1)}
\newcommand{\itr}[1]{it(#1)}

\def \CD {\mathsf{Circ}}
\newcommand{\FCeq}{\mathsf{C}\protect\overrightarrow{\hspace{-.02cm}\scriptstyle\mathsf{irc}}}
\newcommand{\FCopeq}{\mathsf{C}\overleftarrow{\hspace{-.02cm}\scriptstyle\mathsf{irc}}}
\newcommand{\FC}{\mathsf{C}\protect\overrightarrow{\hspace{-.1cm}\scriptstyle\mathsf{irc}}}
\newcommand{\FCop}{\mathsf{C}\overleftarrow{\hspace{-.1cm}\scriptstyle\mathsf{irc}}}

\usepackage{bbm}
\newcommand \Rl[1]{{\mathsf{Rel}}_{\scriptscriptstyle #1}}
\def \Relpoly {\Rl{\poly}}
\def \Relfps {\Rl{\fps}}

\def \Sat {\mathcal{F}}
\def \Forget {\mathcal{U}}

\newcommand{\axiom}[1]{\scriptsize{\mathrm{#1}}}
\newcommand{\eqrule}[1]{\stackrel{\axiom{#1}}{=}}
\newcommand{\rev}[1]{#1^{\op}}
\newcommand{\eqIH}{\stackrel{\tiny \IH{} }{=}}
\newcommand{\eqHA}{\stackrel{\tiny \HA{} }{=}}
\newcommand{\eqHAop}{\stackrel{\tiny \HA{}^{\op} }{=}}

\def\SFGform{\mathsf{SF}}
\def\SFG{\mathbb{SF}} 
\newcommand\Tr[1]{\mathsf{Tr}^{#1}} 

\newcommand\Rwl[1]{\mathsf{L}_{#1}} 
\newcommand\Rwr[1]{\mathsf{R}_{#1}} 

\newcommand{\E}{\mathcal{E}}

\newcommand{\Inj}{\mathsf{In}}  
\newcommand{\Injop}{\Inj^{\op}}  
\newcommand{\Surj}{\mathsf{Su}} 
\newcommand{\Injplus}{\Inj + \Injop} 
\newcommand{\ER}{\mathsf{ER}}  
\newcommand \eqr[1] {\lfloor{#1}\rfloor} 
\newcommand \inv[1] {{#1}^{-1}}
\newcommand{\PER}{\mathsf{PER}}  
\newcommand \peqr[1] {\lfloor\!\lfloor{#1}\rfloor\!\rfloor} 

\usepackage{bbm}
  \def \UNIT {\mathbb{U}} 
   \def \MULT {\mathbb{M}\mathbbm{u}} 
   \def \COUNIT {\mathbb{C}\mathbbm{u}} 
    \def \PF {\mathsf{PF}} 
        \def \PFROB {\mathbb{PF}\mathbbm{r}} 

 \newcommand{\B}{\mathbb{B}}  
\def \FROB {\mathbb{F}\mathbbm{r}} 
\def \IFROB {\mathbb{IF}\mathbbm{r}} 
\def \IPFROB {\mathbb{IPF}\mathbbm{r}} 

\newcommand{\relcomp}{\!\ast\!}

\newcommand{\PMon}{\COUNIT \bicomp{\Perm} \MULT \bicomp{\Perm} \UNIT} 

\newcommand{\copier}{\lower4pt\hbox{$\includegraphics[height=15pt]{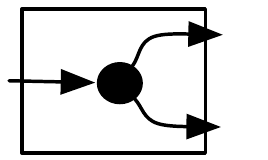}$}}
\newcommand{\adder}{\lower4pt\hbox{$\includegraphics[height=15pt]{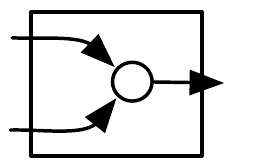}$}}
\newcommand{\amplifier}{\lower4pt\hbox{$\includegraphics[height=15pt]{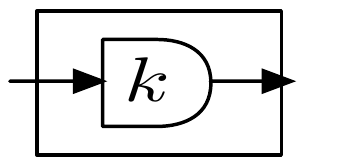}$}}
\newcommand{\register}{\lower4pt\hbox{$\includegraphics[height=15pt]{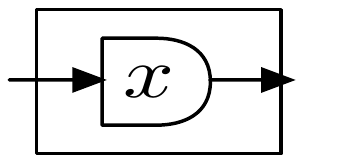}$}}

\newcommand\ord[1] {\overline{#1}}
\newcommand{\Prof}[1]{\mathsf{Prof}(#1)}
\newcommand{\MonB}[1]{\mathsf{Mnd}(#1)} 
\newcommand{\MonC}{\mathbf{Mon}} 

\newcommand{\Ar}[1]{\mathit{Ar}_{\scriptscriptstyle #1}}
\newcommand{\Ob}[1]{\mathit{Ob}_{\scriptscriptstyle #1}}
\newcommand{\ObO}{\Ob{}}
\newcommand{\ArO}{\Ar{}}
\newcommand{\dom}[1]{\mathit{dom}_{\scriptscriptstyle #1}}
\newcommand{\cod}[1]{\mathit{cod}_{\scriptscriptstyle #1}}
\newcommand{\domO}{\dom{}}
\newcommand{\codO}{\cod{}}

\def \catD {\mathbb{D}}
\def \catE {\mathbb{E}}
\def \bicat {\mathfrak{B}}
\newcommand\bicomp[1] {\hspace{-.2pt}{\otimes}_{\scriptscriptstyle #1}\hspace{-.2pt}}

\def \PJ {\mathbb{J}} 

\def\lambdapb {\lambda_{\scriptscriptstyle \mathit{pb}}}
\def\lambdapo {\lambda_{\scriptscriptstyle \mathit{po}}}

\newcommand\starredtext[2]{\begin{center}\parbox{10cm}{#1}\hspace*{2cm}(#2)\end{center}}

\newcommand{\WBboneText}{\lower5pt\hbox{$\includegraphics[height=.5cm]{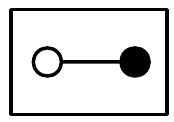}$}}

\newcommand\Ocounit{\lower5pt\hbox{$\includegraphics[width=16pt]{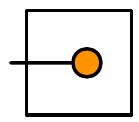}$}}
\newcommand\Ounit{\lower5pt\hbox{$\includegraphics[width=16pt]{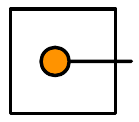}$}}
\newcommand\node{\lower3pt\hbox{$\includegraphics[width=22pt]{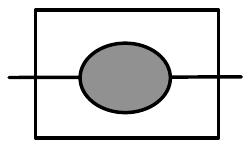}$}}

\newcommand \CospanToER {\Pi}
\newcommand \SpanToER {\Upsilon}
\newcommand \ERToPER {\Xi}
\newcommand \CospanToPCospan {\Lambda}

\newcommand \PermToT[1] {\text{A}_{\scriptscriptstyle #1}}

\newcommand \lefta[1] {\tau^{\scriptscriptstyle #1}}
\newcommand \righta[1] {\rho^{\scriptscriptstyle #1}}

\newcommand \LwA[1] {\mathbb{L}_{\scriptscriptstyle #1}}
\newcommand \Lw {\LwA{}}
\newcommand \var[1] {\mathit{ovar}(#1)}
\newcommand \size[1] {|#1|}
\newcommand\OpDiag{\lower4pt\hbox{$\includegraphics[width=24pt]{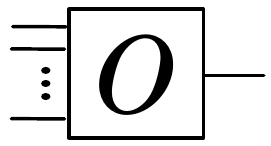}$}}
\newcommand\diagD{\lower4pt\hbox{$\includegraphics[width=24pt]{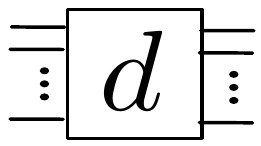}$}}
\newcommand \tpl[1] {{\langle#1\rangle}} 
\newcommand \Mquot {\theta} 

\newcommand \xmapsto[1] {\overset{#1}{\mapsto}}

\newcommand\idcircuit{\lower3pt\hbox{$\includegraphics[width=20pt,height=13pt]{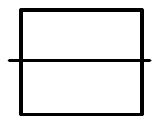}$}}
\newcommand\BcomultSingleAntipode{\lower9pt\hbox{$\includegraphics[width=21pt]{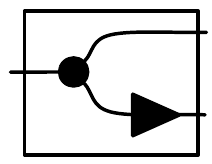}$}}
\newcommand \SVH[1] {\mathbb{SV}_{\scriptscriptstyle #1}}
\def \Q {\mathbb{Q}} 
\newcommand{\spaceFull}{\lower4pt\hbox{$\includegraphics[height=15pt]{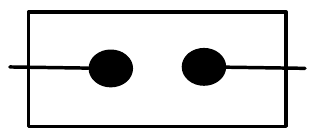}$}}
\newcommand{\spaceZero}{\lower4pt\hbox{$\includegraphics[height=15pt]{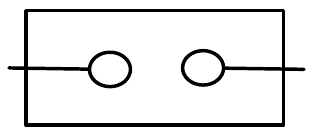}$}}
\newcommand{\spaceYaxis}{\lower4pt\hbox{$\includegraphics[height=15pt]{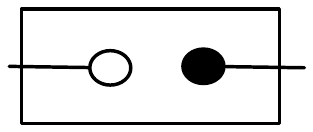}$}}
\newcommand{\spaceXaxis}{\lower4pt\hbox{$\includegraphics[height=15pt]{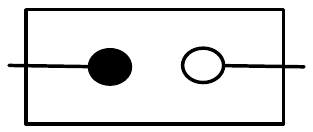}$}}
\newcommand{\spaceXaxisTiny}{\lower3pt\hbox{$\includegraphics[height=12pt]{graffles/spaceXaxis.pdf}$}}
\newcommand{\scalarktwo}{\lower5.5pt\hbox{$\includegraphics[height=18pt]{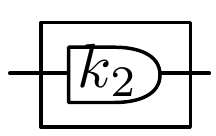}$}}
\newcommand{\spacekonektwo}{\lower5pt\hbox{$\includegraphics[height=18pt]{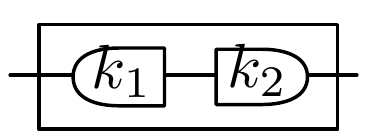}$}}

\newcommand \SpanMat {\Mat{\PID}^{\op} \bicomp{\PJ} \Mat{\PID}}
\newcommand \CospanMat {\Mat{\PID} \bicomp{\PJ} \Mat{\PID}^{\op}}

\newcommand\WFrobR{\includegraphics[width=18pt]{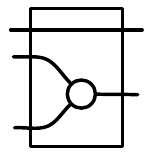}}
\newcommand\WFrobL{\includegraphics[width=18pt]{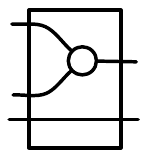}}
\newcommand\BFrobR{\includegraphics[width=18pt]{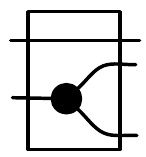}}
\newcommand\BFrobL{\includegraphics[width=18pt]{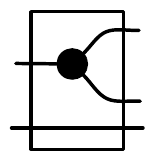}}

\newcommand \MatToSpan {\xi_1}
\newcommand \MatopToSpan {\xi_2}
\newcommand \MatToCospan {\pi_1}
\newcommand \MatopToCospan {\pi_2}
\newcommand \HAToIHw {\upsilon_1}
\newcommand \HAopToIHw {\upsilon_2}

\newcommand \trasp[1] {{#1}^{\scriptscriptstyle T}} 
\newcommand \contrid[1] {\coc{#1}} 

\newcommand \feq {=} 

\newcommand\kernelA{\lower3pt\hbox{$\includegraphics[width=40pt]{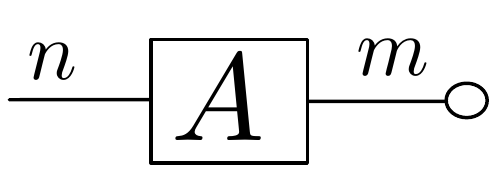}$}}
\newcommand\emptyspace{\lower1pt\hbox{$\includegraphics[width=20pt]{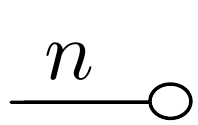}$}}
\newcommand\Astar{\lower4pt\hbox{$\includegraphics[height=.5cm]{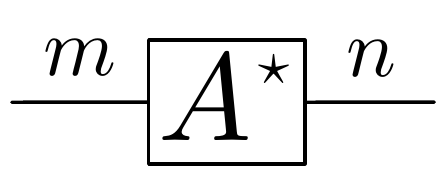}$}}
\newcommand\circuitA{\lower4pt\hbox{$\includegraphics[height=.5cm]{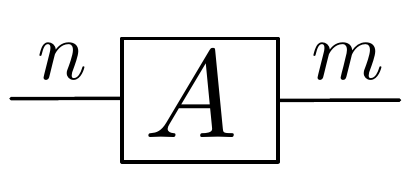}$}}

\newcommand\circuitn{\lower3pt\hbox{$\includegraphics[width=24pt]{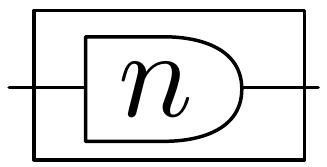}$}}

\newcommand\matrixTwoROneC[2]{{\tiny\begin{pmatrix}  #1 \\ #2 \end{pmatrix}}}
\newcommand{\BboneText}{\lower5pt\hbox{$\includegraphics[height=.5cm]{graffles/Bbone.pdf}$}}

\newcommand \lrnwide[1] {\widetilde{#1}}
\newcommand \FCtoHA {I}
\newcommand \CDtoIH {I'}
\newcommand \SpanMatpoly {\Matpolyop \bicomp{\PJ} \Matpoly}
\newcommand \CospanMatpoly {\Matpoly \bicomp{\PJ} \Matpolyop}
\newcommand \SpanMatfps {\Matfpsop \bicomp{\PJ} \Matfps}
\newcommand \CospanMatfps {\Matfps \bicomp{\PJ} \Matfpsop}

\newcommand \zerov {\vectorfont{0}}  

\makeatletter
\newcommand{\xmapsfrom}[2][]{%
  \ext@arrow3095\leftarrowfill@{#1}{#2}\mapsfromchar
}
\makeatother

\usepackage{datatool}
\usepackage[nopostdot,toc]{glossaries}
\newglossarystyle{mystyle}{%
{\vspace{-1.7em}
\begin{longtable}{@{}p{0.10\textwidth}@{}p{0.05\textwidth}@{}p{0.80\textwidth}@{}p{0.15\textwidth}@{}}}%
{\end{longtable}}%
\renewcommand*{\glsgroupheading}[1]{}%
\ifstrequal{\glsentryuseri{##1}}{-}{

}{
}

}

\makeglossaries

\glsdisablehyper

\newglossaryentry{emptyset}
{
name={\ensuremath{\emptyset}},
description={empty set \dotfill}
}

\newglossaryentry{UnitTns}
{
name={\ensuremath{I}},
description={unit object of the monoidal product in a monoidal category \dotfill}
}

\newglossaryentry{ordn}
{
name={\ensuremath{\ord{n}}},
description={the set $\{1,\dots,n\}$ \dotfill}
}

\newglossaryentry{DisjointUnion}
{
name={\ensuremath{\uplus}},
description={disjoint union of sets \dotfill}
}

\newglossaryentry{Inclusion}
{
name={\ensuremath{\subseteq}},
description={inclusion of sets \dotfill}
}

\newglossaryentry{setminus}
{
name={\ensuremath{\setminus}},
description={difference of sets \dotfill}
}

\newglossaryentry{coslice}
{
name={\ensuremath{x/\catC}},
description={coslice category of $\catC$ under $x \in \catC$ \dotfill},
symbol={\ensuremath{/}}
}

\newglossaryentry{poi}
{
name={\ensuremath{\poi}},
description={sequential composition of arrows in a category \dotfill}
}

\newglossaryentry{tns}
{
name={\ensuremath{\tns}},
description={monoidal product in a monoidal category \dotfill}
}

\newglossaryentry{sigmaxy}
{
name={\ensuremath{\sigma_{x,y}}},
description={symmetry associated with objects $x$ and $y$ in a symmetric monoidal category \dotfill}
}

\newglossaryentry{PROPMon}
{
name={\ensuremath{\Mon}},
description={PROP of the theory of commutative monoids \dotfill}
}

\newglossaryentry{PROPCom}
{
name={\ensuremath{\Com}},
description={PROP of the theory of commutative comonoids \dotfill}
}

\newglossaryentry{PROPB}
{
name={\ensuremath{\B}},
description={PROP of the theory of bialgebras \dotfill}
}

\newglossaryentry{PROPFrob}
{
name={\ensuremath{\FROB}},
description={PROP of the theory of special Frobenius algebras \dotfill}
}

\newglossaryentry{PROPFunction}
{
name={\ensuremath{\F}},
description={PROP of functions \dotfill}
}

\newglossaryentry{PROPPerm}
{
name={\ensuremath{\Perm}},
description={PROP of permutations \dotfill}
}

\newglossaryentry{PermToT}
{
name={\ensuremath{\PermToT{\T}}},
description={Unique PROP morphism from $\Perm$ to a PROP $\T$ \dotfill}
}

\newglossaryentry{bicompD}
{
name={\ensuremath{\bicomp{\D}}},
description={composition of 1-cells along the object $\D$ in a bicategory of bimodules \dotfill}
}

\newglossaryentry{Bicatxx}
{
name={\ensuremath{\bicat(x,x)}},
description={monoidal category of 1-cells $x\to x$ in a bicategory $\bicat$ \dotfill}
}

\newglossaryentry{MonBC}
{
name={\ensuremath{\MonB{\catC}}},
description={category of monoids in a monoidal category $\catC$ \dotfill}
}

\newglossaryentry{MonC}
{
name={\ensuremath{\MonC}},
description={category of monoids and monoid homomorphisms \dotfill}
}

\newglossaryentry{Set}
{
name={\ensuremath{\Set}},
description={category of sets and functions \dotfill}
}

\newglossaryentry{Cat}
{
name={\ensuremath{\Cat}},
description={category of small categories and functors \dotfill}
}

\newglossaryentry{SpanC}
{
name={\ensuremath{\Span{\catC}}},
description={bicategory of objects in $\catC$, spans and span morphisms \dotfill}
}

\newglossaryentry{Cospan}
{
name={\ensuremath{\Cospan{\catC}}},
description={bicategory of objects in $\catC$, cospans and cospan morphisms \dotfill},
symbol={\ensuremath{\CospanO}}
}

\newglossaryentry{ModBicat}
{
name={\ensuremath{\Mod{\bicat}}},
description={bicategory of monads in $\bicat$, bimodules and bimodule morphisms \dotfill}
}

\newglossaryentry{PROPMult}
{
name={\ensuremath{\MULT}},
description={PROP of the theory of multiplication \dotfill}
}

\newglossaryentry{PROPUnit}
{
name={\ensuremath{\UNIT}},
description={PROP of the theory of unit \dotfill}
}

\newglossaryentry{PROPSurj}
{
name={\ensuremath{\Surj}},
description={PROP of surjective functions \dotfill}
}

\newglossaryentry{PROPInj}
{
name={\ensuremath{\Inj}},
description={PROP of injective functions \dotfill}
}

\newglossaryentry{PROPCore}
{
name={\ensuremath{\PJ}},
description={PROP of isomorphisms in a given PROP $\T$ (also called the core of $\T$) \dotfill}
}

\newglossaryentry{LawSigmaE}
{
name={\ensuremath{\LwA{\Sigma,E}}},
description={Lawvere theory generated by a signature $\Sigma$ and a set $E$ of equations \dotfill}
}

\newglossaryentry{PROPCounit}
{
name={\ensuremath{\COUNIT}},
description={PROP of the theory of counit \dotfill}
}

\newglossaryentry{equivCat}
{
name={\ensuremath{\simeq}},
description={equivalence of categories \dotfill}
}

\newglossaryentry{iso}
{
name={\ensuremath{\cong}},
description={isomorphism of objects \dotfill}
}

\newglossaryentry{homset}
{
name={\ensuremath{\catC[x,y]}},
description={set of arrows from $x$ to $y$ in a small category $\catC$ \dotfill}
}

\newglossaryentry{matrixNull}
{
name={\ensuremath{(\,)}},
description={unique element of the $0$-dimensional vector space \dotfill}
}

\newglossaryentry{PROPER}
{
name={\ensuremath{\ER}},
description={PROP of equivalence relations \dotfill}
}

\newglossaryentry{PROPPER}
{
name={\ensuremath{\PER}},
description={PROP of partial equivalence relations \dotfill}
}

\newglossaryentry{PROPIFROB}
{
name={\ensuremath{\IFROB}},
description={PROP of the theory of irredundant special Frobenius algebras \dotfill}
}

\newglossaryentry{PROPIPFROB}
{
name={\ensuremath{\IPFROB}},
description={PROP of the theory of irredundant partial special Frobenius algebras \dotfill}
}

\newglossaryentry{PROPPFROB}
{
name={\ensuremath{\PFROB}},
description={PROP of the theory of partial special Frobenius algebras \dotfill}
}

\newglossaryentry{PROPPF}
{
name={\ensuremath{\PF}},
description={PROP of partial functions \dotfill}
}

\newglossaryentry{PROPCD}
{
name={\ensuremath{\CD}},
description={PROP of circuit diagrams \dotfill}
}

\newglossaryentry{PROPCDleftright}
{
name={\ensuremath{\FC}},
description={PROP of ``functional'' circuit diagrams \dotfill}
}

\newglossaryentry{poly}
{
name={\ensuremath{\poly}},
description={ring of polynomials \dotfill}
}

\newglossaryentry{laur}
{
name={\ensuremath{\laur}},
description={field of formal Laurent series \dotfill}
}

\newglossaryentry{frpoly}
{
name={\ensuremath{\frpoly}},
description={field of fractions of polynomials \dotfill}
}

\newglossaryentry{fps}
{
name={\ensuremath{\fps}},
description={ring of formal power series (streams) \dotfill}
}

\newglossaryentry{ratio}
{
name={\ensuremath{\ratio}},
description={ring of rational fractions of polynomials \dotfill}
}

\newglossaryentry{PROPSFGform}
{
name={\ensuremath{\SFGform}},
description={PROP of orthodox signal flow graphs \dotfill}
}

\newglossaryentry{dsemO}
{
name={\ensuremath{\dsemO}},
description={polynomial semantics $\CD \to \SV{\frpoly}$ \dotfill}
}

\newglossaryentry{dsemHAO}
{
name={\ensuremath{\dsemHAO}},
description={polynomial semantics $\FC \to \Matpoly$ \dotfill}
}

\newglossaryentry{dsemHAOop}
{
name={\ensuremath{\dsemHAOop}},
description={polynomial semantics $\FCop \to \Matpoly^{\op}$ \dotfill}
}

\newglossaryentry{PROPCDrightleft}
{
name={\ensuremath{\FCop}},
description={PROP of ``co-functional'' circuit diagrams (opposite of $\FC$) \dotfill}
}

\newglossaryentry{strsemHAO}
{
name={\ensuremath{\strsemHAO}},
description={denotational stream semantics $\FC \to \SVfps$ \dotfill}
}

\newglossaryentry{strsemO}
{
name={\ensuremath{\strsemO}},
description={denotational stream semantics $\CD \to \SVfps$ \dotfill}
}

\newglossaryentry{osemO}
{
name={\ensuremath{\osemO}},
description={observable behaviour $\CD \to \Relpoly \times \Relfps$ \dotfill}
}

\newglossaryentry{TrO}
{
name={\ensuremath{\Tr{}(\cdot)}},
description={trace in $\IBpoly$ defined by feedback loops guarded by a register \dotfill}
}

\newglossaryentry{ftrC}
{
name={\ensuremath{\ftr{c}}},
description={finite traces of the operational semantics of a circuit $c$ of $\CD$ \dotfill},
symbol={\ensuremath{\ftr{c}}}
}

\newglossaryentry{itrC}
{
name={\ensuremath{\itr{c}}},
description={infinite traces of the operational semantics of a circuit $c$ of $\CD$ \dotfill},
symbol={\ensuremath{\itr{c}}}
} 

\chapterstyle{ger}

\begin{document}

\begin{center}
\vspace{10mm}

{\large\bf TH\`{E}SE}

\vspace{5mm}

{\em en vue d'obtenir le grade de}

\vspace{2.5mm}

{\bf Docteur de l'Universit\'{e} de Lyon}\\[2mm]
{\bf d\'{e}livr\'{e} par l'\'{E}cole Normale Sup\'{e}rieure de Lyon}\\[10mm]
{\bf Discipline~: Informatique}\\[2mm]
{\bf Laboratoire de l'Informatique du Parall\'{e}lisme}\\[2mm]
{\bf \'{E}cole Doctorale en Informatique et Math\'{e}matiques de Lyon}\\[2mm]

\vspace{8mm}
{\em pr\'{e}sent\'{e}e et soutenue publiquement le 5 Octobre 2015}\\
{\em par Monsieur Fabio ZANASI}\\

\vspace{20mm}

{
    \sffamily\bfseries\LARGE
    \begin{center}
    \line(1,0){250}

    \vspace{2mm}
      \textbf{\Large Interacting Hopf Algebras} \\
      \textsc{\normalsize the theory of linear systems}
    \vspace{1mm}

          \line(1,0){250}
\end{center}%
  }%

\vspace{3cm}
\small

\end{center}

\begin{tabular}{lllll}
\textit{Directeurs de th\`{e}se~:} 
&M. &Filippo & BONCHI & \\
&M. &Daniel & HIRSCHKOFF & \\[2mm]
\textit{Apr\`{e}s l'avis de~:}  & M. &Samson &  ABRAMSKY & \\
                & M. &Pierre-Louis &  CURIEN &\\
                & M. &Peter &  SELINGER & \\[2mm]
\textit{Devant le jury compos\'{e}e de~:}      &M. &Samson &  ABRAMSKY &\textit{Rapporteur}\\
                & M. &Filippo &  BONCHI &\textit{Directeur}\\
                & M. &Pierre-Louis &  CURIEN &\textit{Rapporteur}\\
		&M. &Daniel  & HIRSCHKOFF &\textit{Directeur}\\
                &M. &Samuel  & MIMRAM &\textit{Examinateur}\\
                &M. & Prakash &  PANANGADEN &\textit{Examinateur}
\end{tabular}
\thispagestyle{empty}
\makeatother

\newpage\null\thispagestyle{empty}\newpage


\newpage
\thispagestyle{empty}
\section*{Acknoweldgements}

I am deeply grateful to Filippo Bonchi for the amazing amount of time, energy and passion that he invested in me. I think it is very rare to find such a dedicated supervisor and I am very lucky to have met him. I also thank him for making me work on beautiful topics and teach me to seek elegant solutions and stay away from convoluted ones.

I thank Daniel Hirschkoff for his guidance through French lifestyle, regulations and their mysteries. Life in and outside the university would have been much harder without his support. Daniel's self-control and positive attitude really helped me carrying on during bad periods.

Even if it was not officially my supervisor, Pawel Sobocinski played a key role for this thesis. He first disclosed to me the beauties of ``Australian'' category theory and influenced me with his radical views on concurrency and circuit theory. Also, he co-authored the articles that formed this thesis and he has always been extremely available for questions and discussion. 
I sincerely thank him for his time and his great teaching. 

I wish to thank Samson Abramsky, Pierre-Louis Curien and Peter Selinger for writing a report about my thesis and for the huge amount of feedback they sent me, which helped immensely in improving the manuscript. I also thank Samuel Mimram 
 and Prakash Panangaden for accepting of being part of the committee and bringing their insightful perspective on my work.


Thanks to my co-authors Facundo Carreiro, Alessandro Facchini, Stefan Milius, Alexandra Silva and Yde Venema: working with them was a pleasant and enriching experience. I want to also thank Alexandra, as well as Tom Hirschowitz, Matteo Mio and Damien Pous, for the support and the precious advices they have been giving me during my PhD. 

Working in the \emph{Plume} team was a very enjoyable experience. I wish to thank all the members that have been working at the lab during my stay, as well as the staff, for the nice atmosphere they have been creating and the interesting discussions.

I thank my parents for their constant support --- both moral and substantial, with provisions of balsamic vinegar and other goods that made me feel less homesick. My last and speechless thank is for Laura: this thesis is dedicated to her.
\newpage
\thispagestyle{empty}
\begin{abstract}
Scientists in diverse fields use diagrammatic formalisms to reason about various kinds of networks, or compound systems. Examples include electrical circuits, signal flow graphs, Penrose and Feynman diagrams, Bayesian networks, Petri nets, Kahn process networks, proof nets, UML specifications, amongst many others.
Graphical languages provide a convenient abstraction of some underlying mathematical formalism, which gives meaning to diagrams. For instance, signal flow graphs, foundational structures in control theory, are traditionally translated into systems of linear equations. This is typical: diagrammatic languages are used as an interface for more traditional mathematics, but rarely studied per se.

Recent trends in computer science analyse diagrams as first-class objects using formal methods from programming language semantics. In many such approaches, diagrams are generated as the arrows of a PROP --- a special kind of \emph{monoidal category} --- by a two-dimensional syntax and equations. The domain of interpretation of diagrams is also formalised as a PROP and the (compositional) semantics is expressed as a functor preserving the PROP~structure.

The first main contribution of this thesis is the characterisation of
$\SVR$,
the PROP of linear subspaces over a field $\field$. This is an
important domain of interpretation for diagrams appearing in diverse
research areas, like the signal flow graphs mentioned above. We
present by generators and equations the PROP $\IH{}$ of string
diagrams whose free model is $\SVR$. The name $\IH{}$
stands for \emph{interacting Hopf algebras}: indeed, the equations of
$\IH{}$ arise by distributive laws between Hopf algebras, which
we obtain using Lack's technique for composing PROPs.
The significance of the result is two-fold. On the one hand, it offers
a canonical \emph{string diagrammatic syntax} for linear algebra: linear
maps, kernels, subspaces and the standard linear algebraic
transformations are all faithfully represented in the graphical
language. On the other hand, the equations of $\IH{}$ describe
familiar algebraic structures --- Hopf algebras and Frobenius algebras
--- which are at the heart of graphical formalisms as seemingly
diverse as quantum circuits, signal flow graphs, simple electrical
circuits and Petri nets. Our characterisation enlightens the
provenance of these axioms and reveals their linear algebraic nature.

Our second main contribution is an application of $\IH{}$ to the
\emph{semantics of signal processing circuits}. We develop a formal
theory of \emph{signal flow graphs}, featuring a string diagrammatic syntax
for circuits, a structural operational semantics and a denotational
semantics. We prove soundness and completeness of the equations of
$\IH{}$ for denotational equivalence. Also, we study the full
abstraction question: it turns out that the purely operational picture
is too concrete --- two graphs that are denotationally equal may
exhibit different operational behaviour. We classify the ways in which
this can occur and show that any graph can be realised --- rewritten,
using the equations of $\IH{}$, into an executable form where
the operational behaviour and the denotation coincide. This
realisability theorem --- which is the culmination of our developments
--- suggests a reflection about the role of causality in the semantics
of signal flow graphs and, more generally, of computing devices.
\end{abstract} 
\tableofcontents

\chapter{Introduction}
\section{Background}
Scientists in diverse fields use diagrammatic formalisms to reason about various kinds of networks, or compound systems. Examples include electrical circuits, signal flow graphs, Penrose and Feynman diagrams, proof nets, Bayesian networks, Petri nets, Kahn process networks, UML specifications, amongst many others.

These diagrams are formalised to various extent and the mathematics that lies behind the intended \emph{meaning} of diagrams in several such families is, by now, well-understood. An illustrative example are signal flow graphs, foundational structures widely used in control theory and engineering since the 1950's, which are traditionally translated into systems of equations and then solved using standard techniques. This perspective is influenced by physics, where a system is typically modeled by a continuous state-space and the interactions that may occur in it are expressed as continuous state-space transformations, e.g. using differential equations.

Computer science has a rather different approach to modeling. Rather than on global behaviour, the focus is on local, rule-based interactions --- typically, occurring in a \emph{discrete} state-space. The formal semantics of programming languages rests on cornerstones such as \emph{compositionality}, \emph{types} and the use of methods from algebra and logic. In recent years, these principles have started to be fruitfully transferred from one-dimensional syntax to the analysis of diagrammatic languages. \emph{Monoidal categories} have been widely recognised~\cite{BaezRosetta,Abramsky2004,Baez2014,Pavlovic13} as the right mathematical setting in which diagrammatic notations can be studied in a compositional, resource sensitive fashion. Arrows of a monoidal category enjoy a graphical rendition as \emph{string diagrams}~\cite{Joyal1991,Selinger2009} and the two ways --- composition and monoidal product --- of combining arrows are represented pictorially, respectively, by horizontal and vertical juxtaposition of diagrams.

The main actors of our developments are \emph{PROPs} (\textbf{Pro}duct and \textbf{P}ermutation categories~\cite{MacLane1965}), which are symmetric monoidal categories with objects the natural numbers. PROPs can serve both as a \emph{syntax} and as a \emph{semantics} for graphical languages. Also, similarly to \emph{Lawvere theories}~\cite{LawvereOriginalPaper,hyland2007category}, they naturally support the expression of an algebraic structure describing equivalence of string diagrams.

We mention two illustrative examples of this approach. The first concerns \emph{concurrency theory}: in this area coexist traditional graphical formalisms, like Petri nets~\cite{Peterson:1977:PetriNets}, and the more recent process calculi, like CCS~\cite{Milner:1982:CCS}, CSP~\cite{Hoare:1978:CSP} and the $\pi$-calculus~\cite{Sangiorgi:2001:PiCalculus}. In the last two decades, some approaches~\cite{Bruni2013,Soboci'nski2010,Bruni2006} attempted to merge the benefits of the two worlds by modeling Petri nets in a \emph{compositional} way, as graphical process algebras formally described in the framework of PROPs. A proposal that naturally fits this picture is the \emph{Petri calculus}~\cite{Soboci'nski2010}. The syntax is given by a PROP $\mathbb{P}\mathbbm{etri}$ whose arrows $n \to m$ are bounded Petri nets with $n$ ports on the left and $m$ on the right, freely constructed starting from a small set of connectors. The meaning of these diagrams is given in terms of transition systems whose transitions have two labels, intuitively corresponding to left and right boundary of a Petri net: these systems also form a PROP $\mathsf{2LTS}$. The compositional semantics is given as a PROP functor taking a Petri net to its state graph.
\[ \mathbb{P}\mathbbm{etri} \to \mathsf{2LTS} \]

The equations between string diagrams which axiomatise this semantics are subject of ongoing work~\cite{Sobocinski2013a}. Interestingly, the identified algebraic theory is not far removed from those appearing in compositional approaches to quantum information, like the \emph{ZX-calculus}~\cite{Coecke2008,CoeckeDuncanZX2011}. This is our second motivating example of diagrammatic formalism, originated in the research programme of \emph{categorical quantum mechanics}~\cite{Abramsky2004,Abramsky2008:CQM}, whose aim is to develop high-level methods --- informed by the formal semantics of programming languages --- for quantum physics. The ZX-calculus is an algebra of interacting quantum observables, which can be presented as a PROP $\mathbb{ZX}$ whose string diagrams represent physical processes. The equations of $\mathbb{ZX}$ describe the interplay of familiar structures such as \emph{Frobenius algebras} and \emph{Hopf algebras}, which will also appear in our developments. The meaning of diagrams of $\mathbb{ZX}$ is given by linear maps between finite-dimensional Hilbert spaces, forming a PROP~$\mathsf{HS}$.
\[ \mathbb{ZX} \to \mathsf{HS} \]



\section{Content of the Thesis}
The first main contribution of this thesis is a characterisation of the PROP $\SVR$ whose arrows $n \to m$ are \emph{linear subspaces} of $\field^n \times \field^m$, for a field $\field$, and composition is relational. This is a particularly important domain of interpretation for many diagrammatic languages: the meaning of well-behaved classes of systems --- like the signal flow graphs and certain families of Petri nets and quantum processes --- can be typically expressed in terms of linear subspaces.
Our result is a presentation by generators and equations of the PROP $\IH{}$ of string diagrams whose free model is $\SVR$. That means, there is an interpretation of the diagrams of $\IH{}$ as subspaces of $\SVR$, which is also a (symmetric monoidal) isomorphism
\begin{equation*} \IH{} \tr{\cong} \SVR. \end{equation*}

The significance of the result is two-fold. On the one hand, we contend that $\IH{}$ is a \emph{canonical syntax for linear algebra}. Traditional linear algebra abounds in different encodings of the same entities: for instance, spaces are described as a collection of basis elements or as the solution set to a system of equations; matrices, and matrix-related concepts are used ubiquitously as stopgap, common notational conveniences. $\IH{}$ provides an \emph{uniform} description for linear maps, spaces, kernels, etc. based on a small set of simple string diagrams as primitives. Standard methods like Gaussian elimination can be faithfully mimicked in the graphical language, resulting in an alternative, often insightful perspective on the subject matter.

On a different viewpoint, we believe that the equational theory of $\IH{}$ is of independent interest, as it describes fundamental algebraic structures --- Hopf algebras and Frobenius algebras --- which are at the heart of graphical formalisms as seemingly diverse as categorical quantum mechanics, signal flow graphs, simple electrical circuits and Petri nets. Our characterisation enlightens the provenance of these axioms and reveals their linear algebraic nature.

The name $\IH{}$ stands for \textbf{i}nteracting \textbf{H}opf algebras. Indeed, we construct $\IH{}$ \emph{modularly}, starting from the PROP $\HA{}$ --- freely generated by the equations of \textbf{H}opf \textbf{a}lgebras --- and its opposite PROP $\HA{}^{\op}$. Using Lack's technique for composing PROPs~\cite{Lack2004a}, we define two distributive laws that describe different ways of letting $\HA{}$ and $\HA{}^{\op}$ interact. $\IH{}$ is the result of merging the equational theories generated by the two distributive laws. This modular account of $\IH{}$ is actually crucial in constructing the isomorphism $\IH{} \cong \SVR$ --- both with an inductive definition and a universal property --- and will be useful in a number of other ways in our developments. More abstractly, our analysis gives new insights on the interplay of Frobenius and Hopf algebras: for instance, while the authors of the ZX-calculus initially regarded the Frobenius structures as more fundamental, our modular construction reveals that the constituting blocks are Hopf algebras, and the Frobenius equations arise by their composition. In fact, $\IH{}$ axiomatise the \emph{phase-free fragment} of $\mathbb{ZX}$.

\medskip

Our second main contribution is the use of $\IH{}$ to develop a formal theory of signal processing in which \emph{circuits} are first-class citizens. We introduce the \emph{signal flow calculus} and analyse it using the standard methods of programming language theory. The calculus is based on a string diagrammatic syntax, whose terms are meant to represent signal processing circuits. A key feature which makes our language different from similar proposals is that there is no primitive for recursion: feedbacks are a derived notion. Moreover, the wires in our circuits are non-directed and thus there are no assumptions about causal direction of signal flow, allowing us to forego traditional restrictions such as connecting ``inputs'' to ``outputs''. This motivates our formulation of the denotational semantics in terms of linear \emph{relations} rather than functions. Circuit diagrams form a PROP $\CD$ and the (compositional) semantics of a circuit $c \: n \to m$ is given by a functor
\[ \strsemO{} \: \CD \to \SVfps \]
where we regard the subspace $\strsem{c} \subseteq \laur^n \times \laur^m$ as a relation between $\laur^n$ and $\laur^m$. Here $\laur$ is the field of \emph{Laurent series}, a generalised notion of stream typical in algebraic approaches~\cite{Barnabei19983} to signal processing. We are able to characterise ordinary signal flow graphs --- with information flowing from inputs on the left to outputs on the right --- as a certain subclass of $\CD$, whose semantics are precisely the \emph{rational} behaviours in $\SVfps$.

Our design choices make the syntax of $\CD$ abstract enough to enable the use of $\IH{}$ to reason about equivalence of circuits. We prove that the equations of $\IH{}$ are a \emph{sound and complete axiomatisation} for the denotational semantics. This result supports our claim that signal flow graphs are \emph{first-class citizens} of our theory: contrary to traditional approaches, there is a completely graphical way of reasoning about graph transformations and their properties, without the need of translating them first into systems of equations.

\medskip

A fully fledged theory of signal flow graphs demands an operational understanding of circuit diagrams in $\CD$ as executable state-machines. For this purpose, we equip the signal flow calculus with a \emph{structural operational semantics} and study the \emph{full-abstraction} question: how denotational and operational equivalence compare.  Interestingly, it turns out that, in our approach, it is the purely operational picture to be too concrete -- two circuits that are denotationally equal may exhibit different operational behaviour. The problem lies in the generosity of our syntax, which allows for the formation of circuits in which flow directionality cannot be coherently determined. This is not problematic for the denotational semantics, which simply describes a \emph{relation} between ports, but it is for the operational semantics, which is instead deputed to capture the execution of circuits. We classify the ways in which the operational semantics may be less abstract than the denotational semantics, and prove full-abstraction for all the circuits that are free of \emph{deadlocks} and of \emph{initialisation steps}. Interestingly, our argument relies on a syntactic characterisation of these properties, which reveals a connection with a duality that can be elegantly described using the modular character of $\IH{}$.

Because the semantics is not fully abstract for the whole signal flow calculus, one may wonder about the status of all those circuit diagrams --- featuring deadlocks or initialisation steps --- which do not have a clear operational status. Our answer is that they do \emph{not} contribute by any means to the expressivity of the calculus: we prove that, for any behaviour $\strsem{c}$ denoted by a circuit $c$, there exists a circuit $d$, for which the operational semantics is fully abstract, that properly \emph{realises} $\strsem{c}$, that is, $\strsem{d} = \strsem{c}$. In the spirit of the diagrammatic approach, we formulate this result as a procedure effectively transforming $c$ into $d$, using the equations of $\IH{}$ as the rewriting steps.

This \emph{realisability theorem} is the culmination of our work. It makes us able to crystallise what we believe is the main conceptual contribution of the signal flow calculus: a fully fledged operational theory of signal flow graphs as mathematical objects is possible without relying on primitives for flow directionality. Discarding the concept of causality is \emph{harmless}, because the realisability theorem guarantees that any diagram can be transformed into a proper circuit, for which the operational semantics describes the step-by-step execution of a state machine. Moreover, it is \emph{beneficial}, because it is only by forgetting flow that we disclose the beautiful algebraic landscape $\IH{}$ underlying signal flow graphs.

We believe that this lesson can be fruitfully applied to the categorical modeling of other dynamical systems, like electrical circuits and Kahn process networks. Hopefully, the modular techniques that we used to shape $\IH{}$ will contribute to a uniform methodology to axiomatise various kinds of behaviour, thus shedding light on the algebraic structure of a wider spectrum of computing devices, as well as connecting them with existing approaches in quantum and concurrency theory.


\section{Plan of the Thesis and Original Contributions}

We give an overview of the structure of the thesis and pointers to the main contributions. The reader may find at the beginning of each chapter a more detailed introduction and a synopsis.

\paragraph{Chapter \ref{sec:background}} introduces the basics of PROPs (\S~\ref{sec:props}) and PROP operations: sum (\S~\ref{sec:coproduct}), composition (\S~\ref{sec:composingprops}) and fibered sum (\S~\ref{sec:pushout}). Each operation is illustrated with several examples. The heart of the chapter is the technique of PROP composition, which we illustrate by recalling part of the formal theory of monads~\cite{Street_MonadsI} and the work of Lack~\cite{Lack2004a}. This background section also contains new material, whose aim is to demonstrate the pervasiveness of modular techniques and to develop useful tools for the next chapters. The following table gives pointers to the main original contributions.
\begin{center}
\begin{tabular}{m{11cm}|m{2.5cm}}
  \hline
  Generalised distributive laws of PROPs by pullback and pushout & Proposition~\ref{prop:distrLawPbPo} \\
  Distributive laws of PROPs yielding Lawvere theories & Theorem~\ref{Th:LawvereCompositePROP}  \\
Modular characterisation of the PROP of equivalence relations & Theorem~\ref{th:IFROB=ER} \\
  Modular characterisation of the PROP of partial equivalence relations & Theorem~\ref{th:PIFROB=PER} \\
  Modular characterisation of the PROP of partial functions & Example~\ref{ex:partialfunctions} \\
\hline
\end{tabular}
\end{center}

\paragraph{Chapter \ref{chapter:hopf}} uses the techniques introduced in Chapter~\ref{sec:background} to develop the theory of interacting Hopf algebras. Our starting point is the PROP $\HA{}$ of Hopf algebras: we give a novel proof of the fact that it characterises PROPs of matrices, based on PROP composition (\S~\ref{sec:theorymatr}). The technical core of the chapter is the study of distributive laws between $\HA{}$ and $\HA{}^{\op}$ (\S~\ref{sec:ibrw}). We prove that a first distributive law, defined by pullback of matrices, has a characterisation by generators and relations as the PROP $\IBRw$. By duality, it follows a presentation $\IBRb$ by generators and relations also for the PROP resulting from a second distributive law, defined by pushout of matrices. We then merge $\IBRb$ and $\IBRw$ into the theory $\IH{}$ and prove that it characterises PROPs of linear subspaces (\S~\ref{sec:cubetop}). The modular construction yields two factorisation properties for $\IH{}$, in terms of spans and of cospans of $\HA{}$-diagrams. An important aspect of our methodology is the rendition of standard linear algebraic transformations as equational reasoning in the graphical theory. We shall give several demonstrations of this approach in the proofs of the above statements, as well as in the conclusive part of the chapter, where we prove some facts about matrices and subpaces using string diagrams (\S~\ref{sec:graphicallinearlagebra}) and describe the theory $\IH{}$ for linear subspaces over the field $\Q$ of rationals (\S~\ref{sec:instances}). The table below give pointers to the main results.

\begin{center}
\begin{tabular}{m{11.5cm}|m{2.5cm}}
  \hline
  Modular characterisation of the PROP $\HA{}$ & Proposition~\ref{prop:ab=vect} \\
  Axiomatisation of the distributive law between $\HA{}$ and $\HA{}^{\op}$ by pullback & Theorem~\ref{th:Span=IBw} \\
  Axiomatisation of the distributive law between $\HA{}$ and $\HA{}^{\op}$ by pushout & Theorem~\ref{th:IBRb=Cospan} \\
Span and cospan factorisation properties of $\IH{}$ & Theorem~\ref{Th:factIBR} \\
  Isomorphism between $\IH{}$ and $\SVR$ & Theorem~\ref{th:IBR=SVR} \\
  Example: $\IH{}$ as an equational theory of rational subspaces & \S~\ref{sec:instances} \\
  \hline
\end{tabular}
\end{center}
This chapter is based on the following papers.
\begin{itemize}[itemsep=1pt,topsep=1pt,parsep=1pt,partopsep=1pt]
\item F.Bonchi, P.Soboci\'nski, F.Zanasi - \emph{Interacting Bialgebras are Frobenius} - FoSSaCS'14.
\item F.Bonchi, P.Soboci\'nski, F.Zanasi - \emph{Interacting Hopf Algebras} - \url{http://arxiv.org/abs/1403.7048}.
\end{itemize}

\paragraph{Chapter \ref{chapter:SFG}} introduces the signal flow calculus. We present its syntax, the structural operational semantics (\S~\ref{sec:SFcalculus}) and the denotational semantics (\S~\ref{sec:polysem}-\ref{sec:stream}). Circuits of the signal flow calculus can be interpreted as string diagrams of $\IH{}$: we use this observation to prove that the equations of $\IH{}$ are a sound and complete axiomatisation for denotational equivalence. Then we recover traditional signal flow graphs as a sub-class of our circuits and prove that they characterise the rational behaviours of the denotational semantics. This result is well-known in control theory~\cite{Lahti}, but our approach, based on a syntax and a complete set of axioms, allows to formulate it as a Kleene's theorem (\S~\ref{sec:SFG}). The second part of the chapter focuses on the comparison between the operational and the denotational picture. We investigate the two design flaws --- deadlocks and initialisation steps --- making the operational semantics less abstract and give syntactic characterisations for them: this lead us to prove full abstraction for deadlock and initialisation free circuits (\S~\ref{sec:fullabstract}). We then show that any circuit can be realised --- rewritten, using the equations of $\IH{}$, into an executable form where the operational behaviour and the denotation coincides (\S~\ref{sec:realisability}). We conclude our exposition with a formal explanation of the fact that direction of flow is a derivative notion of our theory (\S~\ref{sec:types}). The following table summarises the main contributions of the chapter.
\begin{center}
\begin{tabular}{m{9.8cm}|m{4.5cm}}
  \hline
  Soundness and completeness of $\IH{}$ for the denotational semantics & Theorem~\ref{cor:completeness} \\
  Kleene's theorem for rational stream subspaces &  Theorem~\ref{th:SFGcharactRationals} \\
  Compositionality of the operational semantics & Proposition~\ref{prop:OsemFunctor} \\
  Span form prevents deadlocks & Theorem~\ref{thm:spandeadlock}\\
  Cospan form prevents initialisation steps & Theorem~\ref{thm:cospaninit}\\
  Full abstraction & Corollary~\ref{cor:fullabstractInitDeadFree} \\
  Realisability theorem & Theorem~\ref{thm:realisability},~Corollary~\ref{cor:rea} \\
  \hline
\end{tabular}
\end{center}
This chapter is based on the following papers.
\begin{itemize}[itemsep=1pt,topsep=1pt,parsep=1pt,partopsep=1pt]
\item F.Bonchi, P.Soboci\'nski, F.Zanasi - \emph{A Categorical Semantics for Signal Flow Graphs} - \mbox{CONCUR'14}.
\item F.Bonchi, P.Soboci\'nski, F.Zanasi - \emph{Full Abstraction for Signal Flow Graphs} - PoPL'15.
\end{itemize}

\paragraph{Chapter~\ref{chapter:conclusion}} illustrates some research directions that we propose for future work.




\section{Related Work}

String diagrams originally came to the fore in the study of monoidal categories because they clear away swathes of cumbersome coherence bureaucracy, thereby dramatically simplifying algebraic arguments. Inspired by the seminal paper of Penrose~\cite{Penrose-tensornotation}, there is a tradition of works using string diagrams for characterising \emph{free} monoidal categories, beginning with Joyal and Street~\cite{Joyal1991} --- a comprehensive guide to the state of art is given by Selinger's survey~\cite{Selinger2009}. Our methodology heavily relies on Lack's approach to composing PROPs~\cite{Lack2004a}. Another source of inspiration was Cheng's works on composition of Lawvere theories~\cite{ChengDistrLawsLT} and iterated distributive laws~\cite{Cheng_IteratedLaws}.

The use of string diagrams as compositional syntax of interacting systems is increasingly widespread among computer scientists. We confine ourselves to mentioning some approaches which are particularly close and motivate our developments. In concurrency theory, we mention the algebra of $\mathsf{Span}(\mathsf{Graph})$~\cite{Katis1997a}, the calculus of stateless connectors~\cite{Bruni2006} and the algebra of Petri nets with boundaries~\cite{Soboci'nski2010,Bruni2013}. Frobenius algebras and Hopf algebras appear ubiquitously in these research lines, often interacting as part of the same theory~\cite{Bruni2006,Sobocinski2013a}: $\IH{}$ describes this interaction in a particularly well-behaved setting, in which all behaviours are linear homogeneous. This is of much relevance for the aforementioned approaches, although it leaves out some phenomena that are particularly interesting for concurrency theorists, such as mutual exclusion~\cite{Bruni2006}. 


The programme of categorical quantum mechanics~\cite{Abramsky2004,Abramsky2008:CQM} is another source of inspiration for our approach: in particular, we share the idea of giving an alternative foundation, informed by computer science, category theory and logic, to a subject which is traditionally studied with non-compositional methods. Our theory $\IH{}$ is particularly relevant for one of the most studied formalisms in categorical quantum mechanics, namely the ZX-calculus~\cite{Coecke2008,CoeckeDuncanZX2011}. The equations of $\IH{}$ are at the core of the ZX-calculus, which essentially only adds the properly quantum features such as phase operators. 

In this thesis we give presentations by generations and relations of various PROPs whose arrows are well-known mathematical objects, such as (partial) functions, equivalence relations, matrices and subspaces. This kind of characterisation has been studied for different purposes in diverse areas. We want to mention in particular the research thread on two-dimensional rewriting~\cite{Burroni1993,Lafont95-equationalReasoningTwoDimDiagrams,Lafont2003,Mimram14} where presentations for PROPs of matrices~\cite{Lafont2003}, functions~\cite{Burroni1993} and relations~\cite{Lafont95-equationalReasoningTwoDimDiagrams} are derived in a uniform way by the study of normal forms. Our work relies on a rather different methodology, being based on distributive laws instead of rewriting systems. Actually, there are points of contact between the two approaches, which could be fruitfully combined: we comment more extensively on this in the conclusions (Chapter~\ref{chapter:conclusion}).

Closely related to rewriting approaches is the formalism of interaction nets~\cite{Lafont90_interaction-nets}, a diagrammatic language which generalises proof nets~\cite{Girard87-linearLogic,DanosReigner-multiplicativesProofNets} and is adapted to the encoding various computational models such as Turing machines and cellular automata~\cite{Lafont97-interaction_combinators}. Apparently, $\IH{}$ cannot be reproduced using interaction nets: the form of interaction that it expresses is of a more general kind, featuring diagrams that communicate on \emph{multiple} ports.

The earliest reference for signal flow graphs that we are aware of is Shannon's 1942 technical report~\cite{Shannon1942}. They appear to have been independently rediscovered by Mason in the 1950s~\cite{mason1953feedback} and subsequently gained foundational status in electrical engineering, signal processing and control theory. Our vision of signal flow graphs is inspired by Willems' \emph{behavioural approach}~\cite{Willems2007,Willems-linearsystems}, which is the attempt to, in part, reexamine the central concepts of control theory without giving definitional status to derivable causal information such as direction of flow. Interestingly, signal flow graphs recently attracted coalgebraic modeling \cite{DBLP:journals/tcs/Rutten05,Rutten08_rationalstreamscoalgebraically,Prak2014}. This line of research analyses the coincidence of signal flow graphs, rational streams and a certain class of finite weighted automata using coinduction and the theory of coalgebras. The main difference with these works is that we give a formal syntax for circuits and a sound and complete axiomatisation for semantic equivalence. These features are also present in the work of Milius~\cite{Milius_streamaxiom}, but its syntax is one-dimensional and diagrams are just used for notational convenience. Also, the circuit language is of a rather different flavour; most notably, it features primitives for recursion, which are not necessary in our approach.

Another recent approach to signal flow graphs is Baez and Erbele's manuscript~\cite{BaezErbele-CategoriesInControl}, which appeared on arXiv shortly after our works \cite{BialgAreFrob14,interactinghopf} and the submission of \cite{Bonchi2014b}. In~\cite{BaezErbele-CategoriesInControl}, the authors independently give an equational presentation for PROPs of linear subspaces, which is equivalent to our theory $\IH{}$ --- this paper is inserted in Baez's programme of \emph{network theory}~\cite{Baez2014}, which aims at uniformly describing various kinds of networks used by engineers, ecologists and other scientists using methods from (higher) category theory. A major difference with~\cite{BaezErbele-CategoriesInControl} is in the use of distributive laws of PROPs, which is pervasive in our work and enables a number of analyses that are hampered by a monolithic approach, most notably the characterisation of the isomorphism $\IH{} \tr{\cong} \SVR$ as a universal arrow and the span/cospan factorisation for $\IH{}$. The modular account of $\IH{}$ also means a different choice of primitives: in our approach, feedback is a derivative notion, being constructible by combining the generators of the building blocks $\HA{}$ and $\HA{}^{\op}$ of $\IH{}$; instead, in~\cite{BaezErbele-CategoriesInControl} the ``cup'' and ``cap'' forming a feedback loop appear among the generators. Another significant difference with~\cite{BaezErbele-CategoriesInControl} is that we give a formal operational semantics, which allows us to study full abstraction and realisability, and make a statement about the role of causality in signal flow theory. 

\section{Prerequisites and Notation}
We assume familiarity with the basics of category theory (see e.g.~\cite{mclane,Borceux:1994a}), the definition of symmetric strict monoidal category~\cite{mclane,Selinger2009} (which we often abbreviate as SMC) and of bicategory~\cite{Borceux:1994a,BenabouBicategories}. We write $\catC^{\op}$ for the opposite of a category $\catC$ and $x / \catC$ for the coslice category of $\catC$ under $x \in \catC$. Composition of arrows $f \: x \to y$, $g\: y \to z$ is indicated with $f \gls{poi} g \: x \to z$. We write $\gls{homset}$ for the set of arrows from $x$ to $y$ in a small category $\catC$. It will be sometimes convenient to indicate an arrow $f \: x \to y$ of $\catC$ as $x\tr{f}y$ or $x\tr{f \in \catC}y$. When naming objects and arrows is unnecessary we simply write $\tr{\in \catC}$ or $\tr{}$ if $\catC$ is clear from the context. For $\catC$ symmetric monoidal, we use $\gls{tns}$ for the monoidal product, $\gls{UnitTns}$ for the unit object and $\gls{sigmaxy} \: x \tns y \to y \tns x$ for the symmetry associated with $x,y \in \catC$. 
For a natural number $n>0$, $\gls{ordn}$ is the set $\{1,\dots,n\}$ and $\ord{0} = \emptyset$. We reserve bold letters $\xx,\yy,\zz,\vv,\uu,\ww$ for vectors over a field $\field$. We write $\zerov$ for the zero vector (the length will typically be clear from the context) and $[\vv_1,\dots,\vv_n]$ for the space spanned by vectors $\vv_1,\dots,\vv_n$. Also, $\gls{matrixNull}$ is the unique element of the space with dimension zero.

\chapter{PROPs and their Composition}\label{sec:background}

\section{Overview}

This chapter introduces the basics of the theory of PROPs, focusing on operations to combine PROPs to form richer structures.

PROPs --- an abbreviation of \textbf{pro}duct and \textbf{p}ermutation category --- are symmetric monoidal categories with objects the natural numbers. They made their first appearance in~\cite{MacLane1965} as a means to describe one-sorted algebraic theories. There is a close analogy between PROPs and \emph{Lawvere theories}~\cite{LawvereOriginalPaper,hyland2007category}, with the former being strictly more general. Lawvere theories describe the algebraic structure borne on an object of a \emph{cartesian} category, whereas PROPs fulfill the same purpose in arbitrary \emph{symmetric monoidal} categories. We will further explore the relation between the two notions in \S~\ref{sec:lawvere}.

PROPs share the ability to describe non-cartesian contexts with \emph{operads}~\cite{leinster2004higher}, another family of categories adapted to the study of universal algebra. However, whereas operads are restricted to operations with coarity $1$, PROPs can describe operations with arbitrary arity and coarity. For instance, the level of generality of PROPs is required to express Frobenius algebras and Hopf algebras, which are central in our developments.

Just as Lawvere theories and operads, PROPs allow natural constructions that arise in universal algebra: in this chapter we focus on three of them. The first is the \emph{sum} of theories, which simply takes the disjoint union of the generators and of the equations. We also study the \emph{fibered sum}, in which some structure in common between the summed theories may be identified. The main focus of our developments will be on a third kind of construction: the \emph{composition} of theories by means of a distributive law. This operation, which for PROPs has been developed by Lack~\cite{Lack2004a}, is helpful to describe the modular nature of many algebraic structures. To explain the core intuition, a simple motivating example is the one of a ring, presented by equations:
\begin{multicols}{3} \noindent
\begin{eqnarray*}
 (a + b) + c &=& a + (b + c) \\
 a + b &=& b + a \\
 a + 0 &=& a \\
 a + (-a) &=& 0
  \end{eqnarray*}
    \begin{eqnarray*}
 (a \cdot b) \cdot c &=& a \cdot (b \cdot c) \\
 a \cdot 1 &=& a \\
 1 \cdot a &=& a
  \end{eqnarray*}
  \begin{eqnarray*}
 a \cdot (b + c) &=& (a \cdot b) + (a \cdot c) \\
 (b + c) \cdot a &=& (b \cdot a) + (c \cdot a).
  \end{eqnarray*}
 \end{multicols}
 The idea is to read these equations according to the following pattern: the first column defines an \emph{abelian group}, the second a \emph{monoid} and the third the \emph{distributivity} of the monoid over the group. One can make this formal by expressing the monoid and the abelian group as \emph{monads}; then, orienting left-to-right the equations in the third column defines a distributive law of monads in the sense of Beck~\cite{Beck_distributivelaws1969}. This law yields a new monad, presented by all the above equations: thus rings arise by the composition of monoids with abelian groups.

Note that, differently from sum and fibered sum, a distributive law yields \emph{new} equations expressing the interaction of the theories involved.
We will see in a number of examples that PROP composition, combined with sum and fibered sum, is a powerful heuristics to ease the analysis of complex algebraic structure, allowing to understand them \emph{modularly}, similarly to the case of rings.

This methodology will be applied to the PROPs of commutative monoids, of bialgebras and of special Frobenius algebras. All these examples are also included in~\cite{Lack2004a}. We will also show, as original contributions, the modular understanding of the PROP of partial functions (Example~\ref{ex:partialfunctions}), of equivalence relations (\S~\ref{sec:ER}) and of partial equivalence relations (\S~\ref{sec:PER}). Our analysis will produce a presentation by generators and equations for each of these PROPs. For our purposes, it will be also of importance to develop some ramifications of the composing PROP technique: in particular, we show how Lack's definition of composition can be extended to include distributive laws by pullback and pushout (\S~\ref{sec:distrLawPullback}); we recast in the setting of PROPs some basic operations on distributive laws such as composition, quotient and dual (\S~\ref{sec:iteratedDistrLaws}); finally, we study a family of distributive laws yielding Lawvere theories as the result of composition (\S~\ref{sec:lawvere}). These contributions are also original, when not stated otherwise. They are included to demonstrate the pervasiveness of the modular approach, as well as to give a series of useful techniques for the developments of the next chapter.

\paragraph{Synopsis} The chapter is organised as follows.
\begin{itemize}
\item \S~\ref{sec:props} introduces PROPs and their graphical language of string diagrams. We describe the generation of a PROP by a signature and equations.
\item \S~\ref{sec:coproduct} introduces the operation of PROP sum. 
\item \S~\ref{sec:composingprops} illustrates the operation of PROP composition. We first explain this form of composition in the simpler case of plain categories: categories can be thought as monads (\S~\ref{sec:composingCats}) and composed by distributive laws (\S~\ref{sec:distrlawCats}). We then describe this approach for the case of PROPs: \S~\ref{sec:PROPsMonads} shows how PROPs can be thought as monads and \S~\ref{sec:distrLawPROPs} introduces distributive laws of PROPs.

    In the second part we investigate some ramifications of this technique. In \S~\ref{sec:distrLawPullback} we show how to define distributive laws by pullback and pushouts. \S~\ref{sec:iteratedDistrLaws} explains some basic operations on distributive laws: composition, quotient and dual. Finally, in \S~\ref{sec:lawvere} we investigate a family of distributive law of PROPs yielding Lawvere theories as the result of composition.
\item \S~\ref{sec:pushout} discusses the operation of fibered sum of PROPs. We give a detailed example of how fibered sum, along with PROP sum and composition, can be used to give a presentation by generators and equations to the PROP of equivalence relations (\S~\ref{sec:ER}) and of partial equivalence relations (\S~\ref{sec:PER}).

    We remark that the material presented in \S~\ref{sec:ER}-\ref{sec:PER} is not needed in the sequel, thus it can be safely skipped on a first reading. Nonetheless, those sections offer warm-up examples of the ``cube'' construction that will be pivotal in Chapter~\ref{chapter:hopf}.
\end{itemize}

\label{sec:introProps}

\section{PROPs}
Our exposition is founded on categories called
PROPs (\textbf{pro}duct and \textbf{p}ermutation categories~\cite{MacLane1965}).

\begin{definition} A \emph{PROP} is a symmetric strict monoidal category with objects the natural numbers, where $\tns$ on objects is addition. Morphisms between PROPs are strict symmetric monoidal functors that are identity on objects: PROPs and their morphisms form the category $\PROP$.

We call a \emph{sub-PROP} a sub-category of a PROP $\T$ which is also a PROP.
\end{definition}


PROPs are adapted to the study of universal algebra in a symmetric monoidal setting. Within this perspective, a typical way of defining a PROP is as the free construction on a given set of generators and equations. We express these data in the form of a (one-sorted) \emph{symmetric monoidal theory} (SMT).
%
\begin{definition} A \emph{symmetric monoidal theory} (SMT) is a pair $(\Sigma, E)$ consisting of a \emph{signature} $\Sigma$ and a set of \emph{equations} $E$. The signature $\Sigma$ is a set of \emph{generators} $o \: n\to m$ with \emph{arity} $n$ and \emph{coarity} $m$. The set of \emph{$\Sigma$-terms} is obtained by composing generators in $\Sigma$, the unit $\id \: 1\to 1$ and the symmetry $\sigma_{1,1} \: 2\to 2$ with $;$ and $\tns$. This is a purely formal process: given $\Sigma$-terms $t \: k\to l$, $u \: l\to m$, $v \: m\to n$, one constructs new $\Sigma$-terms $t \mathrel{;} u \: k\to m$ and $t \tns v \: k+n \to l+n$. The set $E$ of \emph{equations} contains pairs $(t,t' \: n\to m)$ of $\Sigma$-terms with the same arity and coarity.
\end{definition}

\begin{figure}[t]
$$
\begin{array}{c}
(t_1 \poi t_3) \tns (t_2 \poi t_4) = (t_1 \tns t_2) \poi (t_3 \tns t_4)\end{array} $$ $$
\begin{array}{rcl}
(t_1 \poi t_2) \poi t_3= t_1 \poi (t_2 \poi t_3) & &
id_n \poi c = c = c\poi id_m\\
(t_1 \tns t_2) \tns t_3 = t_1 \tns (t_2\tns t_3) & &
id_0 \tns t = t  = t \tns id_0\\
\sigma_{1,1}\poi\sigma_{1,1}=id_2 & &
(t \tns id_z) \poi \sigma_{m,z} = \sigma_{n,z} \poi (id_z \tns t)
\end{array} $$
\caption{Axioms of symmetric strict monoidal categories for a PROP $\T$.}\label{fig:axSMC}
\end{figure}

  Now, given an SMT $(\Sigma,E)$, one (freely) obtains a PROP $\T$ by letting the arrows $n\to m$ be the set of $\Sigma$-terms $n\to m$ taken modulo the laws of symmetric strict monoidal categories --- Fig.~\ref{fig:axSMC} --- and the smallest congruence (with respect to $\poi$ and $\tns$) containing the equations $t=t'$ for any $(t,t')\in E$.

   There is a natural graphical representation of these terms as string diagrams, which we now sketch referring to~\cite{Selinger2009} for the details. A $\Sigma$-term $n \to m$ is pictured as a box with $n$ ports on the left and $m$ ports on the right, to which we shall refer with top-bottom enumerations $1,\dots,n$ and $1,\dots,m$. Composition via $\poi$ and $\tns$ are rendered graphically by horizontal and vertical juxtaposition of boxes, respectively.
    \begin{eqnarray}\label{eq:graphlanguage}
t \poi s \text{ is drawn }
\lower7pt\hbox{$\includegraphics[height=.6cm]{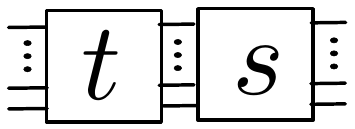}$}
\quad
t \tns s \text{ is drawn }
\lower13pt\hbox{$\includegraphics[height=1.1cm]{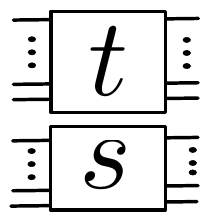}$}.
\end{eqnarray}
    In any SMT there are specific $\Sigma$-terms generating the underlying symmetric monoidal structure: these are $\id_1 \: 1 \to 1$, represented as $\Idnet$, the symmetry $\sigma_{1,1} \: 1+1 \to 1+1$, represented as $\symNet$, and the unit object for $\tns$, that is, $\id_0 \: 0 \to 0$, whose representation is an empty space $\ZeronetT$. Graphical representation for arbitrary identities $\id_n$ and symmetries $\sigma_{n,m}$ are generated according to the pasting rules in~\eqref{eq:graphlanguage}.

The axioms of symmetric strict monoidal categories (Fig.~\ref{fig:axSMC}) are naturally displayed in the graphical language. Compatibility of $\tns$ and $\poi$ is already implicit in the representation of $(t \poi s) \tns (t' \poi s')$ and $(t \tns s) \poi (t' \tns s')$ as the same string diagram:
\[\includegraphics[height=1.2cm]{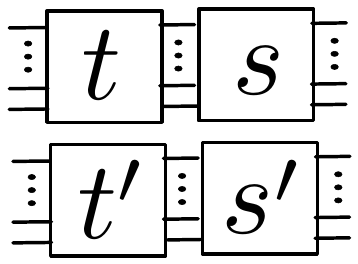}\]
Similarly, associativity of $\poi$, $\tns$ and compatibility of $\tns$ with the unit $\ZeronetT$ are also implicit in the graphical representation. We then have two sliding axioms yielding compatibility of $\poi$ with the identity and naturality of symmetry:

\noindent\begin{minipage}{.5\linewidth}
\begin{equation}\label{eq:sliding1}\tag{SM1}
\lower5pt\hbox{$\includegraphics[height=.6cm]{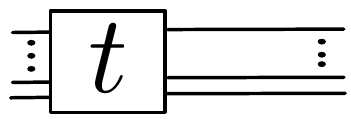}$}
=
\lower5pt\hbox{$\includegraphics[height=.6cm]{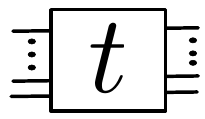}$}
=
\lower5pt\hbox{$\includegraphics[height=.6cm]{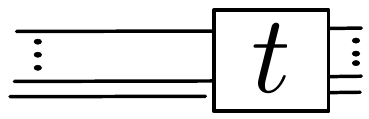}$}.
\vspace{.3cm}
\end{equation}
\end{minipage}
\begin{minipage}{.5\linewidth}
        \begin{equation}\label{eq:sliding2}\tag{SM2}
\lower9pt\hbox{$\includegraphics[height=1cm]{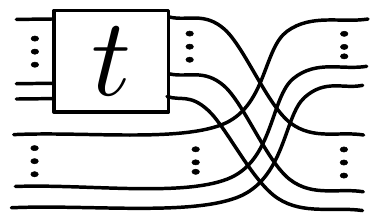}$}
=
\lower9pt\hbox{$\includegraphics[height=1cm]{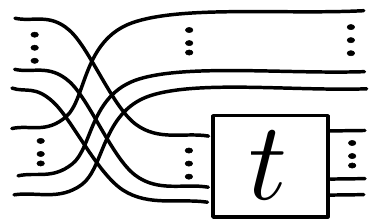}$}.
\end{equation}
\end{minipage}

\medskip
\noindent Finally, we have that $\sigma_{1,1}$ is self-inverse, that is,
        \begin{equation}\label{eq:SymIso}\tag{SM3}
\lower6pt\hbox{$\includegraphics[height=.6cm]{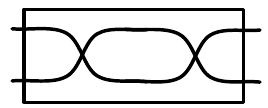}$}
=
\lower6pt\hbox{$\includegraphics[height=.6cm]{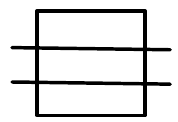}$}.
\end{equation}
As expected, graphical reasoning is \emph{sound and complete}, in the sense that an equality between arrows of a PROP follows from the axioms in Fig.~\ref{fig:axSMC} if and only if it can be derived in the graphical language by using \eqref{eq:sliding1}-\eqref{eq:SymIso} --- \emph{cf.}~\cite{Joyal1991,Selinger2009}.

\begin{convention} In equational reasoning, we will often \emph{orient} equations of SMTs: the notation $c_1 \Rightarrow c_2$ means the use of the equation $c_1 = c_2$ to rewrite a string diagram $c_1$ into $c_2$.
\end{convention}

  \begin{example}~ \label{ex:equationalprops}
  \begin{itemize}
  \item We write $(\Sigma_M,E_M)$ for the SMT of \emph{commutative monoids}. The signature $\Sigma_M$ contains two generators:  \emph{multiplication} --- which we depict as the string diagram $\Wmult \: 2 \to 1$ --- and \emph{unit}, represented as $\Wunit \: 0 \to 1$.
 Equations $E_M$ assert associativity \eqref{eq:wmonassoc}, commutativity~\eqref{eq:wmoncomm} and unitality~\eqref{eq:wmonunitlaw}.
 \begin{multicols}{3}\noindent
 \begin{equation}
\label{eq:wmonassoc}
\tag{A1}
\lower12pt\hbox{$\includegraphics[height=1cm]{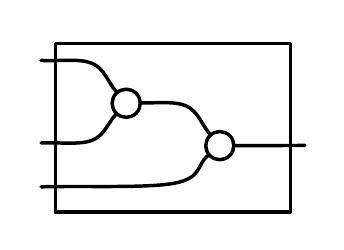}$}
\!\!\!
=
\!\!\!
\lower12pt\hbox{$\includegraphics[height=1cm]{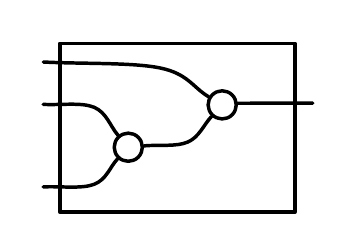}$}
\end{equation}
\begin{equation}
\label{eq:wmoncomm}
\tag{A2}
\lower5pt\hbox{$\includegraphics[height=.6cm]{graffles/Wmult.pdf}$}
\!
=
\!\!\!\!
\lower11pt\hbox{$\includegraphics[height=1cm]{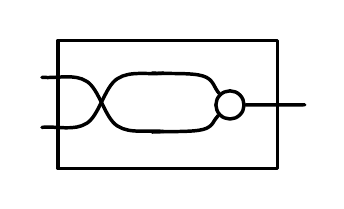}$}
\end{equation}
\begin{equation}
\label{eq:wmonunitlaw}
\tag{A3}
\lower11pt\hbox{$\includegraphics[height=1cm]{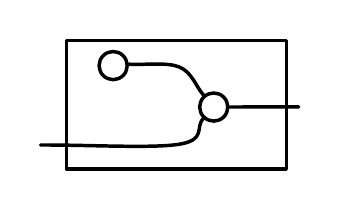}$}
\!\!\!
=\!
\lower5pt\hbox{$\includegraphics[height=.6cm]{graffles/idcircuit.pdf}$}
\end{equation}
\end{multicols}
We call $\gls{PROPMon}$ the PROP freely generated by the SMT $(\Sigma_M,E_M)$.
\item We also introduce the SMT $(\Sigma_C, E_C)$ of \emph{cocommutative comonoids}. The signature $\Sigma_C$ consists of a \emph{comultiplication} $\Bcomult \: 1 \to 2$ and a \emph{counit} $\Bcounit \: 1\to 0$. $E_C$ is the following set of equations.
\begin{multicols}{3}\noindent
\begin{equation}
\label{eq:bcomonassoc}
\tag{A4}
\lower11pt\hbox{$\includegraphics[height=1cm]{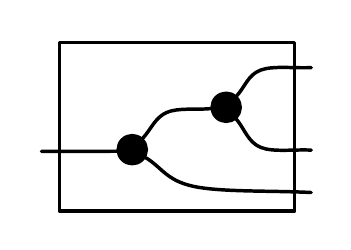}$}
\!\!\!
=
\!\!\!
\lower11pt\hbox{$\includegraphics[height=1cm]{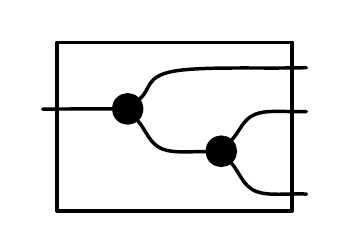}$}
\end{equation}
\begin{equation}
\label{eq:bcomoncomm}
\tag{A5}
\lower5pt\hbox{$\includegraphics[height=.6cm]{graffles/Bcomult.pdf}$}
\!
=
\!\!\!
\lower11pt\hbox{$\includegraphics[height=1cm]{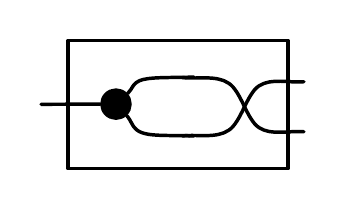}$}
\end{equation}
\begin{equation}
\label{eq:bcomonunitlaw}
\tag{A6}
\lower11pt\hbox{$\includegraphics[height=1cm]{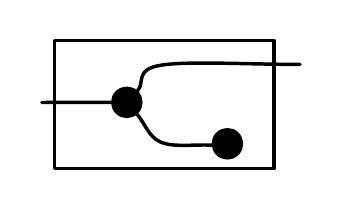}$}
\!\!\!
=
\!
\lower6pt\hbox{$\includegraphics[height=.6cm]{graffles/idcircuit.pdf}$}
\end{equation}
\end{multicols}
 We call $\gls{PROPCom}$ the PROP freely generated by $(\Sigma_C, E_C)$. Modulo the white vs. black colouring, the (string diagrams representening the) arrows of $\Com$ can be seen as those of $\Mon$ ``reflected about the $y$-axis''. This observation yields that $\Com \gls{iso} \Mon^{op}$. More generally, for $\T$ a freely generated PROP, $\T^{op}$ can be presented by generators and equations which are those of $\T$ reflected about the $y$-axis.
\item The PROP $\gls{PROPB}$ of (commutative/cocommutative) \emph{bialgebras} is generated by the theory $(\Sigma_M \uplus \Sigma_C, E_M \uplus E_C \uplus B)$, where $B$ is the following set of equations.
    \begin{multicols}{2}
\noindent
\begin{equation}
\label{eq:unitsl}
\lower2pt\hbox{$
\tag{A7}
\lower5pt\hbox{$\includegraphics[height=.6cm]{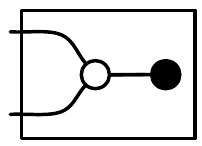}$}
=
\lower5pt\hbox{$\includegraphics[height=.6cm]{graffles/lunitsr.pdf}$}
$}
\end{equation}
\begin{equation}
\label{eq:unitsr}
\tag{A9}
\lower5pt\hbox{$\includegraphics[height=.6cm]{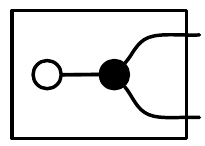}$}
=
\lower5pt\hbox{$\includegraphics[height=.6cm]{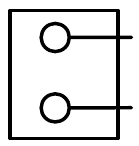}$}
\end{equation}
\begin{equation}
\label{eq:bialg}
\tag{A8}
\lower6pt\hbox{$\includegraphics[height=.6cm]{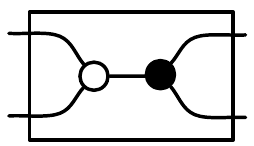}$}
=
\lower11pt\hbox{$\includegraphics[height=.9cm]{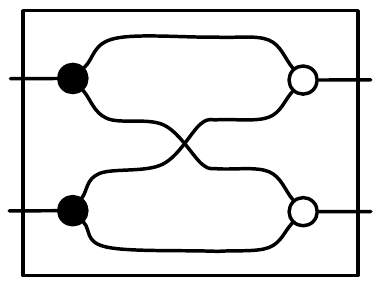}$}
\end{equation}
\begin{equation}
\label{eq:bwbone}
\tag{A10}
\lower4pt\hbox{$\includegraphics[height=.5cm]{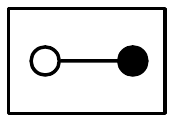}$}
=
\lower4pt\hbox{$\includegraphics[width=15pt]{graffles/idzerocircuit.pdf}$}
\end{equation}
\end{multicols}
One can read \eqref{eq:unitsl}-\eqref{eq:bwbone} as saying that the operations of the comonoid preserve the monoid structure.
\item The PROP $\gls{PROPFrob}$ of \emph{special Frobenius algebras}~\cite{Carboni1987} is generated by the theory $(\Sigma_M \uplus \Sigma_C, E_M \uplus E_C \uplus F)$, where $F$ is the following set of equations.
    \begin{multicols}{2}\noindent
    \begin{equation}\label{eq:BWFrob}\tag{F1}
\lower11pt\hbox{$\includegraphics[height=1cm]{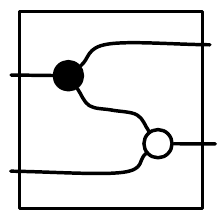}$}
=
\lower8pt\hbox{$\includegraphics[height=.7cm]{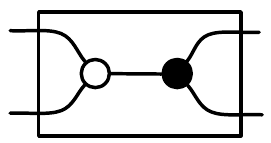}$}
=
\lower11pt\hbox{$\includegraphics[height=1cm]{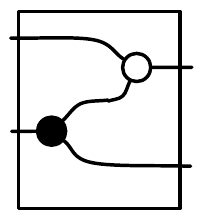}$}
\end{equation}
\begin{equation}\label{eq:BWSep}\tag{F2}
\lower8pt\hbox{$\includegraphics[height=.7cm]{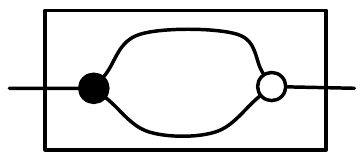}$}
=
\lower6pt\hbox{$\includegraphics[height=.6cm]{graffles/idcircuit.pdf}$}
\end{equation}
\end{multicols}
Intuitively, \eqref{eq:BWFrob}-\eqref{eq:BWSep} allow for any continuous deformation of diagrams. In other words, arrows $f$ of $\FROB$ can be defined by their topology only: the relevant information in any diagrammatic representation of $f$ is which ports on the left and on the right are linked.

    Bialgebras and special Frobenius algebras express two different ways of combining a monoid and a comonoid. We shall see later (\S~\ref{sec:composingprops}) how the equations describing such interaction can be seen as arising by an operation on the corresponding PROPs.
\end{itemize}
  \end{example}

\begin{remark}[Models of a PROP] \label{rmk:models}
The statement that $(\Sigma_M, E_M)$ is the SMT of commutative monoids --- and, similarly for those of comonoids, bialgebras and Frobenius algebras --- can be made rigorous through the notion of \emph{model} (sometimes also called algebra) of a PROP. For any symmetric strict monoidal category $\catC$, a model of a PROP $\T$ in $\catC$ is a symmetric strict monoidal functor $\funF \: \T \to \catC$. There is a category $\mathsf{Model}(\T,\catC)$ whose objects are the models of $\T$ in $\catC$. Now, turning to the example of commutative monoids, we can form a category $\mathsf{Monoid}(\catC)$ whose objects are the commutative monoids in $\catC$, i.e., objects $x \in \catC$ equipped with arrows $x \tns x \to x$ and $I \to x$ satisfying commutativity, associativity and unitality.
Given any model $\funF \: \Mon \to \catC$, one can easily prove that $\funF(1)$ is a commutative monoid in $\catC$: this yields a functor $\mathsf{Model}(\Mon, \catC) \to \mathsf{Monoid}(\catC)$. Saying that $(\Sigma_M, E_M)$ is the SMT of commutative monoids amounts to saying that this functor is an equivalence natural in $\catC$. 

We shall not go into more details about models as they are not necessary for our developments. We refer the reader to~\cite{Lack2004a} for more information.
\end{remark}

Example~\ref{ex:equationalprops} only shows PROPs freely generated from an algebraic specification. However, one can also define PROPs in a more direct manner, without relying on SMTs. We give two basic examples (using distinct typesetting to emphasize the different flavour):
\begin{itemize}
\item the PROP $\gls{PROPFunction}$ whose arrows $n \to m$ are functions from $\ord{n}$ to $\ord{m}$;
\item the PROP $\gls{PROPPerm}$ whose arrows $n \to m$ are bijections from $\ord{n}$ to $\ord{m}$. Note that arrows $n \to m$ exist only if $n = m$, in which case they are the permutations on $\ord{n}$.
\end{itemize}
This kind of definition is often useful to give a different, more concrete perspective on symmetric monoidal theories. For instance, the PROP $\F$ is \emph{presented} by the theory of commutative monoids, in the sense that there is an isomorphism between $\F$ and the PROP $\Mon$ freely generated by that theory. The correspondence is given by considering a string diagram $t \in \Mon[n,m]$ as the graph of a function of type $\{1,\dots, n\}\to \{1, \dots, m\}$. For instance, $\Wmult \tns \Wunit \: 2 \to 2$ describes the function $f \: \{1,2\} \to \{1,2\}$ mapping both elements to $1$. By duality, $\Com \cong \Fop$, that is, $\Fop$ is presented by the theory of commutative comonoids.

Similarly, $\Perm$ provides a concrete description of the theory $(\gls{emptyset}, \emptyset)$ with empty signature and no equations. To see this, note that arrows of the free PROP over $(\emptyset, \emptyset)$ are constructed by tiling together only $\Idnet$ and $\symNet$. Up-to the laws of SMCs, a term $n \to n$ of this kind uniquely represents a permutation of the elements of $\ord{n}$. For instance,
$$\includegraphics[height=1cm]{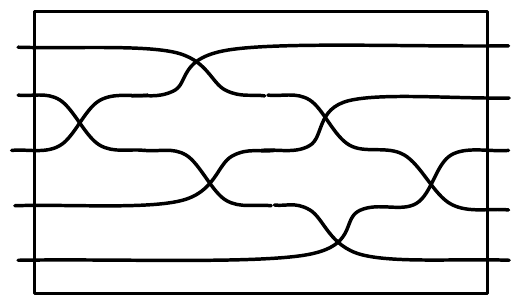}$$
describes the permutation on $\{1,2,3,4,5\}$ exchanging the first two elements with the last three.

One can also provide PROPs giving a concrete description of the theories of bialgebras and of special Frobenius algebras. Interestingly, these PROPs can be obtained \emph{modularly}, by composing together the ones already given for the theories of commutative monoids and of comonoids. This reflects our observation that bialgebras and Frobenius algebras are different ways in which a monoid and a comonoid interact. The next section will introduce the necessary tools to understand such compositions formally and give the desired characterisation for the two SMTs.

\paragraph{PROPs as Functors} For our developments it is useful to recall Lack's observation~\cite{Lack2004a} that PROPs are objects of a certain coslice category. To state this, we first need to recall the non-symmetric version of a PROP, called PRO (\textbf{pro}duct category).

\begin{definition}  A \emph{PRO} is a strict monoidal category with objects the natural numbers, where $\tns$ on objects is addition. Morphisms between PROs are strict monoidal functors that are identity on objects: PROs and their morphisms form the category $\PRO$.
\end{definition}

Roughly, a PROP $\T$ can be described as a PRO that contains a copy of $\Perm$, which forms its symmetry structure. This is made precise by observing that $\Perm$ is the \emph{initial object} in the category $\PROP$. The unique PROP morphism $\gls{PermToT} \: \Perm \to \T$ can be inductively defined starting from the assignment of the symmetry $\sigma_{1,1} \: 2 \to 2$ to the permutation $p_{1,1} \in \Perm[2,2]$ which interchanges the two elements of $\ord{2} = \{1,2\}$ --- all the other permutations in $\Perm$ are obtained from $p_{1,1}$ and the identities via $\poi$ and $\tns$.
Now, by regarding $\PermToT{\T} \: \Perm \to \T$ as a PRO morphism, one can define a functor from $\PROP$ to the coslice category $\Perm \glssymbol{coslice} \PRO$, which maps $\T$ to $\PermToT{\T} \: \Perm \to \T$. By initiality of $\Perm$, this functor is fully faithful and thus exhibits $\PROP$ as a full subcategory of $\Perm \glssymbol{coslice} \PRO$.

Conversely, it is worth spelling out why \emph{not} all the objects of $\Perm / \PRO$ are PROPs: starting from an arbitrary PRO morphism $\Theta \: \Perm \to \PS$, one could be tempted to \emph{define} the symmetry $\sigma_{n,m} \: n+m \to m+n$ in the PRO $\PS$ as the image under $\Theta$ of the permutation $p_{n,m}  \: n+m \to m+n$ which interchanges the first $n$ elements with the last $m$, as to make $\PS$ a PROP. However, in principle there is no reason why $\sigma_{n,m}$ should be natural in $n$ and $m$ as required.


\begin{example} We give a simple counterexample. Let $\PS$ be the PRO whose arrows are obtained by composing the identity $\id_1 \: 1 \to 1$ and a generator $\Delta \: 2 \to 1$  via $\tns$ and $\poi$, and then quotienting by the laws of strict monoidal categories. One can then form the coproduct $\Perm + \PS$ --- this is the PRO whose arrows are obtained by composing those of $\Perm$ and $\PS$ via $\tns$ and $\poi$, identifying the identities and quotienting by the laws of strict monoidal categories. Observe that $\Perm + \PS$ is an object of $\Perm / \PRO$: there is a PRO morphism $\iota_1 \: \Perm \to \Perm + \PS$ given by the coproduct injection. However, it is not a PROP, because the permutations do not yield a symmetry structure in $\Perm + \PS$. For instance,
\[
\xymatrix@=15pt{
2 + 1 \ar[d]_{\iota_1(p_{2,1})} \ar[rr]^{\Delta \tns \id_1} && 1+1 \ar[d]^{\iota_1(p_{1,1})} \\
1+2 \ar[rr]^{\id_1 \tns \Delta} && 1+1
}
\]
does not commute.
\end{example}

Following the above observations, we can fix the relationship between $\PROP$ and $\Perm \glssymbol{coslice} \PRO$.

\begin{proposition}\label{prop:coslicePROP} 
$\PROP$ is isomorphic to the full subcategory of $\Perm \glssymbol{coslice} \PRO$ whose objects are PRO morphisms $\Theta \: \Perm \to \PS$ such that the family of arrows $(\Theta(p_{n,m}) \: n +m \to m+n)_{n,m \in \PS}$ forms a symmetry in $\PS$. \end{proposition}

\begin{remark} In~\cite{Lack2004a} Lack \emph{defines} $\PROP$ as $\Perm / \PRO$. Instead, we chose to stick to the (more restrictive) definition of PROPs as symmetric monoidal categories: ours is a rather standard formulation, also given in MacLane's original paper~\cite{MacLane1965}, which we find more intuitive and simple for applications.
A shortcoming in not following Lack's approach is that we lose the 1-1 correspondence between PROPs and monads in a certain bicategory (\emph{cf.} Proposition~\ref{cor:PROPsAreMonads}). 
However, this does not affect the definition of PROP composition by distributive laws of monads --- see Remark~\ref{rmk:PROPsclosedunderdistrlaw}.
\end{remark}

\medskip


In the rest of the chapter we shall present three different ways of combining PROPs: sum (coproduct), sequential composition and fibered sum (pushout). These operations will allow us to understand PROPs \emph{modularly}, as the result of the interaction of simpler components. This perspective is crucial in investigating more sophisticated and interesting examples of SMTs, like the theories of matrices and linear subspaces that are the theme of the next chapter. \label{sec:props}


\section{PROP Sum}\label{sec:coproduct}
Given PROPs $\T$ and $\PS$, one can calculate their coproduct $\T + \PS$ in $\PRO$ by identifying their symmetry structures. First, following Proposition~\ref{prop:coslicePROP}, we associate with $\T$ and $\PS$ PRO morphisms $\PermToT{\T} \: \Perm \to \T$ and $\PermToT{\PS} \: \Perm \to \PS$. Then, let $\T + \PS$ be given by the following pushout in $\PRO$:
$$\xymatrix{ \Perm \ar[r]^{\PermToT{\PS}} \ar[d]_{\PermToT{\T}} & \PS \ar[d] \\ \T \ar[r] & \T + \PS}$$

\begin{proposition} $\T + \PS$ is the coproduct of $\T$ and $\PS$ in $\PROP$. \end{proposition}
 \begin{proof} We check that $\T + \PS$ is a PROP. Pushouts in $\PRO$ may be calculated as in $\Cat$: that means, arrows of $\T + \PS$ are given by (1) combining the arrows of $\T$ and $\PS$ via $\tns$ and $\poi$, and (2) identifying the permutations, i.e. the arrows $\tr{\in \T}$ and $\tr{\in \PS}$ in the image of the same arrow $\tr{\in \Perm}$. 
 PRO morphisms $\T \to \T + \PS \tl{} \PS$ simply interpret arrows of $\T$ and $\PS$ as arrows of $\T + \PS$.

We define the symmetry $\sigma_{n,m} \: n+m \to m+n$ in $\T + \PS$ to be the image under $\PermToT{\T}$ (equivalently, under $\PermToT{\PS}$) of the permutation in $\Perm$ which interchanges the first $n$ elements with the last $m$. This arrow is a symmetry (i.e., a natural isomorphism) in $\T$ by definition of $\PermToT{\T}$, and also in $\PS$ by definition of $\PermToT{\PS}$. Since arrows in $\T + \PS$ are just combinations of arrows of $\T$ and $\PS$, it follows that $\sigma_{n,m}$ is an isomorphism natural in $n$ and $m$ also in $\T + \PS$. Therefore, $\T + \PS$ is a symmetric monoidal category and thus a PROP.

Since $\T$, $\PS$ and $\T + \PS$ are PROPs and $\PROP$ is a full subcategory of $\Perm / \PRO$ (Proposition~\ref{prop:coslicePROP}), it follows that arrows ${\T \to \T + \PS \tl{} \PS}$ in the above diagram are PROP morphisms: we let them be the coproduct injections. With an analogous reasoning it is straitghtforward to check that the universal property of $\T + \PS$ as pushout in $\PRO$ yields the one as coproduct in $\PROP$.
\end{proof}

 When $\T$ and $\PS$ are freely generated PROPs, the above description provides a simple recipe for a presentation of $\T + \PS$.
\begin{proposition}\label{prop:SMCforSum} Suppose that $\T$ and $\PS$ are PROPs  freely generated by SMTs $(\Sigma_1, E_1)$ and $(\Sigma_2, E_2)$ respectively. Then $\T + \PS$ is freely generated by the sum of theories $(\Sigma_1 \gls{DisjointUnion} \Sigma_2, E_1 \uplus E_2)$.
 \end{proposition}
By Proposition~\ref{prop:SMCforSum}, arrows $n \to m$ of $\T + \PS$ are $\Sigma_1 \uplus \Sigma_2$-terms quotiented by $E_1 \uplus E_2$. We can always represent these arrows as sequences
\begin{equation}\label{eq:seqT1+T2}
n\tr{\in \T}\tr{\in \PS}\tr{\in \T} \dots \tr{\in \PS}\tr{\in \T}m\end{equation}
 of $\Sigma_1$- and $\Sigma_2$-terms modulo $E_1$ and $E_2$. To see this, recall that $\Sigma_1 \uplus \Sigma_2$-terms are constructed by composing the generators of $\Sigma_1 \uplus \Sigma_2$, $\id \: 1 \to 1$ and $\sigma_{1,1} \: 2 \to 2$ with $\poi$ and $\tns$. Then, functoriality of $\tns$ --- \emph{cf.}~Fig.\ref{fig:axSMC} --- allows to put any term $f \tns g$ consisting of a $\Sigma_1$-term $f$ and a $\Sigma_2$-term $g$ into the shape $(f \tns \id) \poi (\id \tns g)$ of a $\Sigma_1$-term followed by a $\Sigma_2$-term, and similarly for $g \tns f$. It follows that any $\Sigma_1 \uplus \Sigma_2$-terms is equal modulo the equations of Fig.~\ref{fig:axSMC} to a sequence as in~\eqref{eq:seqT1+T2}.

\begin{example}[Directed Acyclic Graphs] \label{ex:dags} In~\cite{Fiore2013} the sum of PROPs is used to characterise \emph{directed acyclic graphs} (dags). A dag is a graph with directed edges in which there are no cycles\footnote{Connectivity of dags considered in~\cite{Fiore2013} is relational, i.e. there is at most one edge between every two nodes.}. \emph{Interface}-dags (idags) are directed acyclic graphs extended with a left interface $\ord{n}$ and a right interface $\ord{m}$: edges can have elements of $\ord{n}$ as sources and elements of $\ord{m}$ as targets. Below are two examples with interfaces $\ord{2}/\ord{3}$ and $\ord{3}/\ord{1}$ respectively, taken from~\cite{Fiore2013}.
    \begin{eqnarray}\label{eq:idags}
\lower7pt\hbox{$\includegraphics[height=1.3cm]{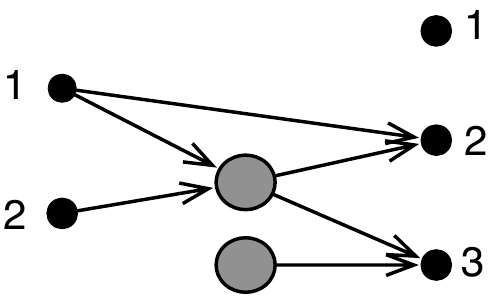}$}
& \qquad \qquad &
\lower7pt\hbox{$\includegraphics[height=1.3cm]{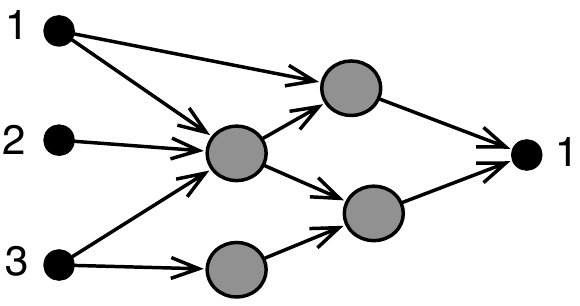}$}.
\end{eqnarray}
One can form a PROP $\mathsf{D}$ whose arrows $n \to m$ are idags with left interface $\ord{n}$ and right interface $\ord{m}$. We refer to~\cite{Fiore2013} for a precise definition of $\mathsf{D}$. Informally, the monoidal product is given by putting two idags side by side. Composition $n \tr{g_1} z \tr{g_2} m$ works by gluing the common interface $\ord{z}$ and redirecting edges of $g_1$ to nodes of $g_2$ accordingly. Here is the composite of the two idags in~\eqref{eq:idags}.
    \begin{eqnarray}\label{eq:idagsconcat}
\lower7pt\hbox{$\includegraphics[height=1.3cm]{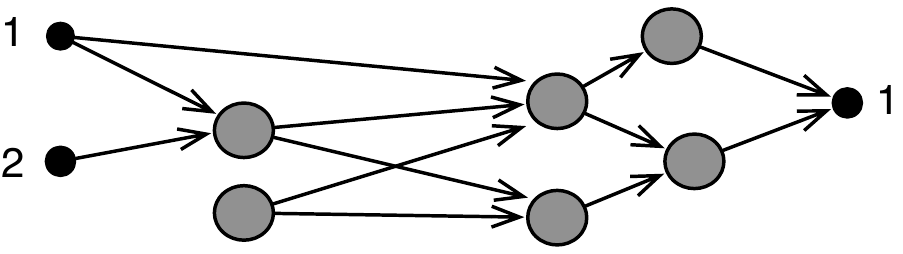}$}
\end{eqnarray}
The PROP $\mathsf{D}$ enjoys a presentation by generators and equations as the sum $\mathbb{SB} + \mathbb{N}$. Here $\mathbb{SB}$ is the PROP of \emph{special bialgebras}, obtained by quotienting $\B$ (Example~\ref{ex:equationalprops}) by~\eqref{eq:BWSep}. $\mathbb{N}$ is the PROP freely generated by the signature consisting of one ``node'' $\node \: 1 \to 1$ and no equations.

The underlying idea is that $\Bcomult$, $\Wmult$, $\Bcounit$ and $\Wunit$ give the branching structure of edges and $\node$ is used to represent nodes. For instance,~\eqref{eq:idags} becomes:
    \begin{eqnarray*}
\lower7pt\hbox{$\includegraphics[height=1.4cm]{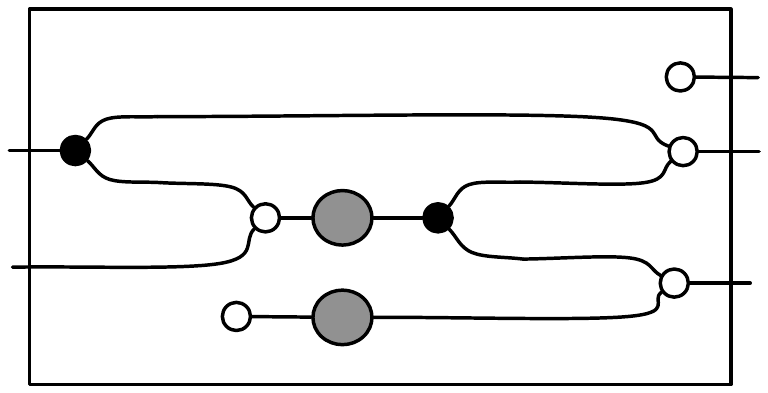}$}
& \qquad \qquad &
\lower7pt\hbox{$\includegraphics[height=1.4cm]{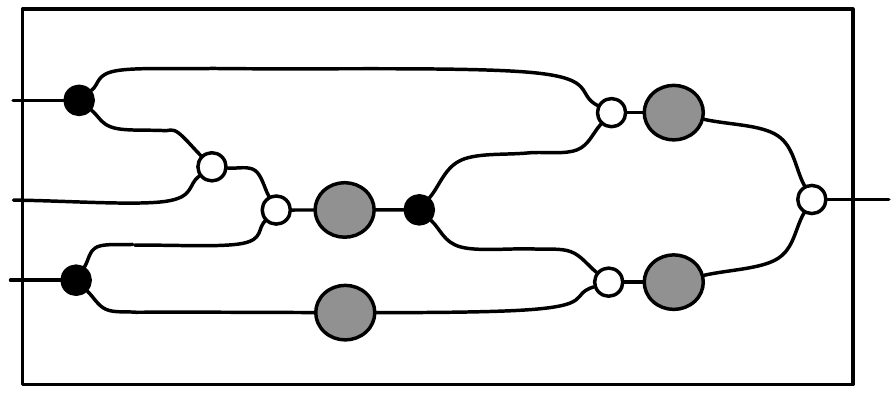}$}.
\end{eqnarray*}
The equations of $\mathbb{SB}$ allow to model composition of idags in the desired way. 

As a concluding note, we mention that by tweaking $\mathbb{SB} + \mathbb{N}$ one can characterise other familiar classes of structures. For instance, the quotient of $\mathbb{SB} + \mathbb{N}$ by $\Bcomult \poi (\node \tns \idcircuit) \poi \Wmult = \node$ characterises finite posets --- see~\cite{MimramThesis,Mimram15posets}.
\end{example}


\section{PROP Composition}\label{sec:composingprops}
The sum $\T+\PS$ is the least interesting way of combining PROPs, because there are no equations that express compatibility conditions between $\T$ and $\PS$ when ``interacting'' in $\T+\PS$. Such interactions are common in algebra: for instance, a ring is given by a monoid and an abelian group, subject to equations telling how the former structure distributes over the latter. Similarly, bialgebras and Frobenius algebras (Example~\ref{ex:equationalprops}) describe two different ways of combining a commutative monoid and a commutative comonoid. An example of a different flavour is the one of a function, which is always decomposable as a surjection followed by an injection.


In~\cite{Lack2004a} Lack shows how these phenomena can be uniformly described as the operation of composing PROPs. The conceptual switch is to understand PROPs as certain monads, which allows to define their composition as a distributive law. We will present this approach in steps, first presenting the simpler case of composition of plain categories (\S~\ref{sec:composingCats}-\ref{sec:distrlawCats}) and then adding the symmetric monoidal structure to the picture (\S~\ref{sec:PROPsMonads}-\ref{sec:distrLawPROPs}).

\subsection{Categories are Monads}\label{sec:composingCats}

As shown in the classical work \cite{Street_MonadsI} of Street, the theory of monads can be developed in an arbitrary bicategory $\bicat$~\footnote{Actually, Street worked in a 2-category, but the same theory can be developed in any bicategory with relatively minor modifications~\cite[\S 3.1]{Lack2004a}.}.

\begin{definition} A \emph{monad} on an object $x$ of $\bicat$ is a 1-cell $\funF \: x \to x$ with 2-cells $\eta^{\funF} \: \id_x \to \funF$ and $\mu^{\funF} \: \funF \poi \funF \to \funF$ (called the \emph{unit} and the \emph{multiplication} respectively) making the following diagrams commute.
\begin{multicols}{2}\noindent
\begin{eqnarray}
\vcenter{
    \xymatrix{
    \funF \ar[dr]_{\id} \ar[r]^{\funF \eta^{\funF}} & \funF \poi \funF \ar[d]_{\mu^{\funF}} & \ar[l]_-{\eta^{\funF} \funF} \ar[dl]^{\id} \funF\\
    & \funF &
    }
}
\label{diag:lawsmonads1}
\end{eqnarray}
\begin{eqnarray}
\label{diag:lawsmonads2}
\vcenter{
    \xymatrix{
    \funF \poi \funF \poi \funF \ar[d]_{\mu^{\funF} \funF} \ar[rr]^{\funF \mu^{\funF}} && \funF \poi \funF \ar[d]^{\mu^{\funF}} \\
    \funF \poi \funF \ar[rr]^{\mu^{\funF}} && \funF
    }
}
\end{eqnarray}
\end{multicols}

\label{def:monadquotient}
A \emph{morphism} between monads $x \tr{\funF} x$ and $x \tr{\G} x$ is a 2-cell $\Mquot \: \funF \to \G$ making the following diagrams commute\footnote{A notion of morphism can be defined also between monads on different objects, like in~\cite{Street_MonadsI}. We will not need that level of generality here.}.
\begin{multicols}{2}\noindent
\begin{eqnarray}
\vcenter{
    \xymatrix@C=18pt@R=20pt{
    \id_x \ar[d]_{\eta^{\funF}}\ar[dr]^{\eta^{\G}} & \\
    \funF \ar[r]^-{\Mquot} & \G
   }
}
\end{eqnarray}
\begin{eqnarray}
\vcenter{
    \xymatrix@C=20pt@R=18pt{
    \funF\poi\funF\ar[d]_{\mu^{\funF}}\ar[r]^-{\Mquot \Mquot} & \G \poi \G \ar[d]^{\mu^{\G}}\\
    \funF \ar[r]^{\Mquot} & \G\\
       }
   }
\end{eqnarray}
\end{multicols}
An epimorphic monad morphism is called a \emph{monad quotient}.
\end{definition}
For $\bicat = \Cat$, the above definition yields the standard notion of monad as an endofunctor with a pair of natural transformations. Something interesting happens for the case of the bicategory $\bicat = \Span{\gls{Set}}$, defined below.

\begin{definition} Let $\catC$ be a category with pullbacks. The bicategory $\gls{SpanC}$ of \emph{spans} on $\catC$ is given as follows:
\begin{itemize}[itemsep=1pt,topsep=1pt,parsep=1pt,partopsep=1pt]
\item objects are the objects of $\catC$
\item a 1-cell from $x$ to $y$ is a span $x \tl{f} z \tr{g} y$ in $\catC$.
\item a 2-cell from $x \tl{f} z \tr{g} y$ to $x \tl{f'} z' \tr{g'} y$ is a span morphism, that is, an arrow $h \: z \to z'$ in $\catC$ making the following diagram commute.
$$ \xymatrix@=8pt{ & \ar[dl]_{f} z \ar[dd]^{h} \ar[dr]^{g} & \\ x && y \\ & \ar[ul]^{f'} z' \ar[ur]_{g'} }$$
We shall call $h$ a \emph{span isomorphism} whenever it is invertible, i.e. there is $h^{-1}$ such that $h\poi h^{-1}=\id = h^{-1}\poi h$.
\item composition of 1-cells is by pullback; that is, the composite of $x \tl{f} z \tr{g} y$ and $y \tl{f'} z' \tr{g'} x'$ is $x \tl{f}\tl{p_1} \tr{p_2}\tr{g} y$ given by
 \begin{equation*}
 \xymatrix@R=10pt@C=20pt{
           &         & \ar[dl]_{p_1} . \pushoutcorner \ar[dr]^{p_2}  &  & \\
         & \ar[dl]_{f}z \ar[dr]^{g} &         & \ar[dl]_{f'}z' \ar[dr]^{g'}& \\
 x &        & y &        & x'
 }
 \end{equation*}
    We assume that there is a specific choice of pullback span for each pair of arrows. Thus composition of 1-cells is only weakly associative and unital, making $\Span{\bicat}$ a bicategory.
 \item Composition of 2-cells along objects (that is, horizontal composition) is given by universal property of pullback:
  $$ \text{given 2-cells $h$ and $h'$   }
  \vcenter{
      \xymatrix@R=8pt@C=14pt{
             & \ar[dl]_{} \ar[dd]^{h} \ar[dr]^{} &         & \ar[dl]_{} \ar[dd]^{h'} \ar[dr]^{}& \\
                 &                            &     &                              &  \\
              & \ar[ul]_{}\ar[ur]^{} &         & \ar[ul]_{}\ar[ur]^{}&
     }
 }
  \quad \text{their composite $i$ is } \quad
     \vcenter{
         \xymatrix@R=8pt@C=16pt{
               &         & \ar[dl]_{} . \pushoutcorner \ar[dr]^{} \ar@{-->}@/^12.5pt/[dddd]^{i} &  & \\
             & \ar[dl]_{}  \ar[dd]^{h} \ar[dr]^{} &         & \ar[dl]_{} \ar[dd]^{h'} \ar[dr]^{}& \\
                  &                            &          &                              &   \\
              & \ar[ul]_{} \ar[ur]^{} &         & \ar[ul]_{}\ar[ur]^{}& \\
               &         & \ar[ul]_{} . \pullbackcorner \ar[ur]^{} &  &
     }
 }.
 $$
     Composition of 2-cells along 1-cells (i.e., vertical composition) is simply composition in $\catC$:
     $$
          \vcenter{
          \xymatrix@R=15pt@C=25pt{
                 & \ar[dl]_{} \ar[d]^{h} \ar[dr]^{} &         \\
                  & \ar[l]_{} \ar[d]^{h'} \ar[r]^{} &   \\
                       & \ar[ul]_{}\ar[ur]^{}&
         }
     }
    $$
\end{itemize}
\end{definition}

The interest for the bicategory of spans stems from the following folklore observation.
\begin{proposition}\label{prop:CatAsMonads} Small categories are precisely the monads in $\Span{\Set}$. \end{proposition}
\begin{proof} Let $\funF$ be a monad on an object $\ObO$ of $\Span{\Set}$. A category $\catC$ can be recovered from $\funF$ as follows.
\begin{itemize}[noitemsep,topsep=0pt,parsep=0pt,partopsep=0pt]
\item the objects of $\catC$ are the elements of $\ObO$.
 \item the arrows of $\catC$ are given by $\funF$. Indeed, $\funF$ is a span $\ObO \tl{\domO} \ArO \tr{\codO} \ObO$ and arrows $f \: x \to y$ of $\catC$ are the elements $f \in \ArO$ such that $\domO(f) = x$ and $\codO(f) = y$.

 \item the composition of arrows $f \: x \to y$ and $g \: y \to z$ is handled by the multiplication of the monad. By definition, $\mu$ is the following span morphism
 \begin{equation*}
 \xymatrix@R=10pt@C=20pt{
           &         & \ar[dl]_{p_1} P \pushoutcorner \ar@/^2.5pc/[dddd]^{\mu}  \ar[dr]^{p_2}  &  & \\
         & \ar[dl]_{\domO}\ArO \ar[dr]^{\codO} &         & \ar[dl]_{\domO}\ArO \ar[dr]^{\codO}& \\
 \ObO &        & \ObO &        & \ObO \\
 &&&&\\
     && \ar[uull]^{\domO} \ArO \ar[uurr]_{\codO} &&
 }
 \end{equation*}
     where the innermost square is a pullback.
     Objects of $P$ are pairs of arrows $(f,g)$ which are composable, that is, $\codO(f) = \domO(g)$. $p_1$ and $p_2$ are the two projections. $\mu$ associates with $(f,g) \in P$ their composition $f \poi g \in \ArO$.
     Commutativity of the diagram guarantees that $f \poi g$ is an arrow of type $x \to z$ in $\catC$. Commutativity of \eqref{diag:lawsmonads2} yields associativity of composition.
  \item the identity arrow on $x \in \ObO$ is the image of $x$ under the span morphism $\eta$:
            $$ \xymatrix@=8pt{ &  \ar[dl]_{\id} \ObO \ar[dd]^{\eta} \ar[dr]^{\id} & \\ \ObO && \ObO \\ & \ar[ul]^{\domO} \ArO \ar[ur]_{\codO} }.$$
       Commutativity of \eqref{diag:lawsmonads1} amounts to the standard identity laws.
 \end{itemize}
The converse transformation from $\catC$ to a monad $\funF$ follows the same ideas. \end{proof}


\subsection{Distributive Laws of Categories}\label{sec:distrlawCats}

Now that we have an understanding of categories as monads, we can compose categories (with the same objects) via \emph{distributive laws}~\cite{Beck_distributivelaws1969}. First, we give the abstract definition.

\begin{definition} Let $(\funF, \eta^{\funF},\mu^{\funF})$, $(\G, \eta^{\G},\mu^{\G})$ be monads in a bicategory $\bicat$ on the same object. A \emph{distributive law} of $\funF$ over $\G$ is a 2-cell $\lambda \: \funF \poi \G \to \G \poi \funF$ in $\bicat$ making the following diagrams commute.
{%
\begin{multicols}{2}\noindent
\begin{eqnarray}
\vcenter{
    \xymatrix@C=18pt@R=20pt{
    \funF\ar[d]_{\funF\eta^\G}\ar[dr]^{\eta^{\G} \funF} & \\
    \funF\poi\G\ar[r]^-{\lambda} & \G\poi\funF \\
    \G\ar[u]^{\eta^{\funF} \G}\ar[ur]_{\G\eta^{\funF}}
    &
   }
}
\label{eq:distrlawequationsUnit}
\end{eqnarray}
\begin{eqnarray}
\vcenter{
    \xymatrix@C=20pt@R=18pt{
    \funF\poi\G\poi\G\ar[d]_{\funF\mu^{\G}}\ar[r]^-{\lambda \G} &
       \G\poi\funF\poi\G\ar[r]^-{\G\lambda} &
       \G\poi\G\poi\funF\ar[d]^{\mu^{\G} \funF} \\
    \funF\poi\G\ar[rr]^-{\lambda} & & \G\poi\funF\\
    \funF\poi\funF\poi\G\ar[u]^{\mu^{\funF} \G}\ar[r]_{\funF \lambda} &
       \funF\poi\G\poi\funF\ar[r]_{\lambda \funF} &
       \G\poi\funF\poi\funF\ar[u]_{\G\mu^{\funF} }
       }
   }
   \label{eq:distrlawequationsMult}
\end{eqnarray}
\end{multicols}
}%
\end{definition}

A distributive law $\lambda \: \funF \poi \G \to \G \poi \funF$ yields a monad $\G \poi \funF$ with the following unit and multiplication:
\begin{equation}\label{eq:defcomposedmonad}
\begin{aligned}
\eta^{\G\poi\funF} \ \: \ & \id \xrightarrow{\eta^{\funF}} \funF \xrightarrow{\eta^{\G}\funF} \G \poi \funF \\
\mu^{\G\poi\funF} \ \: \ & \G \poi \funF \poi \G \poi \funF \xrightarrow{\G \lambda \funF} \G \poi \G \poi \funF \poi \funF \xrightarrow{\mu^{\G}\funF\funF} \G \poi \funF \poi \funF\xrightarrow{\G\mu^{\funF}} \G \poi \funF
\end{aligned}
\end{equation}

Let us verify how the abstract definition works for the case of categories. Pick categories $\catC$ and $\catD$ with the same set $\ObO$ of objects, seen as monads $\ObO \tl{\dom{\catC}} \Ar{\catC}\tr{\cod{\catC}} \ObO$ and $\ObO \tl{\dom{\catD}} \Ar{\catD}\tr{\cod{\catD}} \ObO$ in $\Span{\Set}$. A distributive law $\lambda \: \catC \poi \catD \to \catD \poi \catC$ is a span morphism
 \begin{equation*}
 \xymatrix@R=10pt@C=25pt{
           &         & \ar[dl] . \pushoutcorner \ar@/^3pc/[dddd]^{\lambda}  \ar[dr]  &  & \\
         & \ar[dl]_{\dom{\catC}}\Ar{\catC} \ar[dr]^{\cod{\catC}} &         & \ar[dl]_{\dom{\catD}}\Ar{\catD} \ar[dr]^{\cod{\catD}}& \\
 \ObO &        & \ObO &        & \ObO \\
         & \ar[ul]^{\dom{\catD}}\Ar{\catD} \ar[ur]_{\cod{\catD}} &         & \ar[ul]^{\dom{\catC}}\Ar{\catC} \ar[ur]_{\cod{\catC}}& \\
           &         & \ar[ul]_{} . \pullbackcorner  \ar[ur]^{}  &  &
 }
 \end{equation*}
mapping composable pairs $x\tr{\in \catC}\tr{\in \catD}y$ to composable pairs $x\tr{\in \catD}\tr{\in \catC}y$. As described in~\eqref{eq:defcomposedmonad}, $\lambda$ allows to define a monad structure on $\catD \poi \catC$. That means, $\lambda$ yields a category $\catD \poi \catC$ whose arrows $x\to y$ are composable pairs $x\tr{\in \catD}\tr{\in \catC}y$ of arrows of $\catD$, $\catC$ and
$$\text{the composite of } x\tr{f\in \catD}\tr{g\in \catC}y \text{ and }y\tr{f'\in \catD}\tr{g'\in \catC}z \text{ is }x\tr{f\in \catD}\lambda(\tr{g\in \catC}\tr{f'\in \catD})\tr{g'\in \catC}z.$$

\begin{remark} Distributive laws in the above sense have an equivalent description in terms of factorisation systems: a category is expressible as a composite $\catC \poi \catD$ precisely when each arrow $x \to y$ has a \emph{unique} factorisation as $x \tr{\in \catC}\tr{\in \catD} y$~\cite{RosebrRWood_fact}. This observation reveals that distributive laws of categories hardly arise naturally, as they require factorisations to be unique on-the-nose. For this reason we shall postpone any example to when more relaxed notions of distributive laws --- for which factorisations have weaker uniqueness conditions --- are introduced in the next sections. Also, we shall not illustrate further the relation between distributive laws and factorisation systems, because it is not central for our exposition. The interested reader is referred to~\cite[\S 4]{ChengDistrLawsLT} for an overview of the topic.
\end{remark}

\subsection{PROPs are Monads}\label{sec:PROPsMonads}

In this and the next section we detail how PROPs can be seen as monads in a bicategory and composed together via distributive laws. At first glance, one could be tempted of working within the framework of \S~\ref{sec:composingCats}-\ref{sec:distrlawCats}: being categories, PROPs yield monads in $\Span{\Set}$. However, this approach does not take into account the symmetric monoidal structure carried by PROPs: composing PROPs via a distributive law in $\Span{\Set}$ would yield a category that is not necessarily a PROP.

The idea is then to refine the bicategory of interest. First, rather than considering spans in $\Set$, we take spans in the category $\gls{MonC}$ of monoids and monoid homomorphisms. Intuitively, this takes into account the monoidal structure and we obtain a variation of Proposition~\ref{prop:CatAsMonads}.

\begin{proposition}\label{prop:MonCatsasMonads}
Small strict monoidal categories are precisely monads in $\Span{\MonC}$.
\end{proposition}
\begin{proof} The construction is analogous to the one of Proposition~\ref{prop:CatAsMonads}. Given a monad with underlying span $\ObO \tl{\domO} \ArO \tr{\codO} \ObO$, its unit $\eta$ and multiplication $\mu$ define respectively the identities and composition in the corresponding category. The fact that $\ArO$ and $\ObO$ are monoids yields a monoidal product $\tns$ with unit object the unit of the monoid $\ObO$. Also, $\tns$ obeys the laws of Fig.~\ref{fig:axSMC} because $\eta$ and $\mu$ are arrows in $\MonC$.
\end{proof}

Following Proposition~\ref{prop:MonCatsasMonads}, monads in $\Span{\MonC}$ over the monoid $(\N, +, 0)$ are precisely PROs. Similarly to above, we could now try defining composition of PROPs as composition of the underlying PROs  $\T$ and $\PS$ via a distributive law $\lambda \: \PS \poi \T \to \T \poi \PS$ in $\Span{\MonC}$. As expected, this notion of composition is still ill-behaved as it does not take correctly into account the symmetry structure. The problem is that $\T \poi \PS$ contains \emph{two} copies of $\Perm$, one given by $\PermToT{\T} \: \Perm \to \T\to \T \poi \PS$ and the other by $\PermToT{\PS} \: \Perm \to \PS \to \T \poi \PS$, which do not necessarily agree.

The correct approach is to make explicit the symmetry structure of any PROP $\PR$ in the form of a \emph{left} and a \emph{right action} $\lefta{\PR} \: \Perm \poi \PR \to \PR$ and $\righta{\PR} \: \PR \poi \Perm \to \PR$, yielded by $\PermToT{\PR} \: \Perm \to \PR$. 
Then, we shall define the composite $\T \bicomp{\Perm} \PS$ of PROPs $\T$ and $\PS$ as a coequaliser in $\PRO$
\begin{equation}\label{eq:coequalizer}
\xymatrix{ \T \poi \Perm \poi \PS \ar@<0.66ex>[r]^{\righta{\T} \PS} \ar@<-0.66ex>[r]_{\PS_1 \lefta{\PS}} & \T \poi \PS \ar[r] & \T \bicomp{\Perm} \PS}
\end{equation}
which, intuitively, is responsible for identifying the two copies of $\Perm$ in $\T \poi \PS$. 

This account of PROPs is actually reminiscent of the familiar notion of \emph{bimodule}, which in algebra designates abelian groups with both a left and a right action over a ring; the construction~\eqref{eq:coequalizer} corresponds to the usual tensor product of bimodules.


 This suggests the idea to express PROPs as monads in $\Span{\MonC}$ with a bimodule structure and compose them using~\eqref{eq:coequalizer}. To make this formal, we first define the bicategory of bimodules in a given bicategory $\bicat$. We will then focus on bimodules in $\Span{\MonC}$ to capture PROPs. 


\begin{definition}\label{def:mod} Given a bicategory $\bicat$ with coequalisers, $\gls{ModBicat}$ is the bicategory of \emph{bimodules} in $\bicat$:
\begin{itemize}[itemsep=1pt,topsep=1pt,parsep=1pt,partopsep=1pt]
\item objects are the monads in $\bicat$
\item 1-cells are bimodules; that is, given monads $x \xrightarrow{\funF} x$ and $y \xrightarrow{\G} y$ in $\bicat$, a 1-cell in $\Mod{\bicat}$ from $\funF$ to $\G$ is a 1-cell $x \xrightarrow{\FH} y$ in $\bicat$ equipped with 2-cells $\tau \: \FH \poi \G\to \FH$ and $\rho \: \funF \poi \FH \to \FH$ in $\bicat$ called \emph{left} and \emph{right action} respectively. They satisfy compatibility conditions expressed by commutativity of the following diagrams in $\bicat$.
    \begin{eqnarray}\label{eq:actionscompatibilityconds}
    \vcenter{
    \xymatrix{
    \funF \poi \funF \poi \FH \ar[d]_{\funF \rho} \ar[r]^{\mu^{\funF} \FH} & \funF \poi\FH \ar[d]_{\funF \rho} & \ar[l]_{\eta^{\funF}} \FH \ar[dl]^{\id} \\
    \funF \poi \FH \ar[r]^{\rho} & \FH
    }
    }\quad
    \vcenter{
        \xymatrix{
    \FH \poi \G \poi \G \ar[d]_{\tau \G} \ar[r]^{\FH\mu^{\G}} & \FH \poi\G \ar[d]_{\tau \G} & \ar[l]_{\eta^{\G}} \FH \ar[dl]^{\id} \\
    \FH \poi \G\ar[r]^{\tau} & \FH
    }
    }
     \quad
    \vcenter{
    \xymatrix{
    \funF \poi \FH \poi \G \ar[r]^{\funF \tau} \ar[d]_{\rho \G} & \funF \poi \FH \ar[d]^{\rho} \\
    \FH \poi \G \ar[r]^{\tau} & \FH
    }
    }
    \end{eqnarray}
\item 2-cells are bimodule morphisms, that is, given 1-cells $\funF \xrightarrow{\FH} \G$ (with actions $\tau$ and $\rho$) and $\funF \xrightarrow{\FH'} \G$ (with actions $\tau'$ and $\rho'$), a 2-cell in $\Mod{\bicat}$ from $\FH$ to $\FH'$ is a 2-cell $\FH \xrightarrow{\psi} \FH'$ in $\bicat$ compatible with left and right actions:
      \begin{eqnarray}\label{eq:bimodulemorphism}
      \vcenter{
      \xymatrix{
      \funF \poi \FH \ar[r]^{\rho} \ar[d]_{\funF \psi} & \FH \ar[d]_{\psi} & \ar[l]_{\tau} \FH \poi \G \ar[d]^{\psi \G} \\
       \funF \poi \FH' \ar[r]^{\rho'} & \FH' & \ar[l]_{\tau'} \FH' \poi \G
      }
      }
      \end{eqnarray}
\item the identity 1-cell on an object $x \tr{\funF} x$ of $\Mod{\bicat}$ is $\funF$ itself, with left and right actions $\funF \poi \funF \to \funF$ given by the multiplication of $\funF$ as a monad in $\bicat$. The identity 2-cells are the same as in $\bicat$.
\item composition of 1-cells is by coequaliser. Suppose that $\funF \xrightarrow{\FH} \D$ (with actions $\rho$ and $\tau$) and $\D \xrightarrow{\FH'} \G$ (with actions $\rho'$ and $\tau'$) are 1-cells in $\Mod{\bicat}$ given by 1-cells
    $$x \xrightarrow{\funF} x \xrightarrow{\FH} z \xrightarrow{\D} z \xrightarrow{\FH'} y \xrightarrow{\G} y$$
    in $\bicat$. The composite $\FH \gls{bicompD} \FH' \: \funF \to \G$ is defined by coequaliser in $\bicat$:
    \begin{equation}\label{eq:compBimodules} \vcenter{
    \xymatrix{ \FH \poi \D \poi \FH' \ar@<0.66ex>[r]^{\rho \FH'} \ar@<-0.66ex>[r]_{\FH \tau'} & \FH \poi \FH' \ar[r] & \FH \bicomp{\D} \FH'}
    }
    \end{equation}
\item given that 2-cells in $\Mod{\bicat}$ are also 2-cells in $\bicat$, horizontal and vertical composition of 2-cells in $\Mod{\bicat}$ is defined as in $\bicat$.
\end{itemize}
\end{definition}
The same construction of Definition~\ref{def:mod} is used in~\cite{ChengDistrLawsLT} to give an account of Lawvere theories as monads in a bicategory. Interestingly, it also appears in topological field theory to describe orbifold completion --- see~\cite[Def. 4.1]{OrbifoldTFT}.

We now focus on our main application. Since $\Span{\MonC}$ has coequalisers~\cite{SobocinskiHeindel-VanKempenUniversal2011}, one can form the bicategory $\Mod{\Span{\MonC}}$ of bimodules in $\Span{\MonC}$. The next example details how Definition~\ref{def:mod} instantiates for this case. We shall later verify that PROPs are monads in $\Mod{\Span{\MonC}}$.

\begin{example}\label{ex:ModSpan} Objects in $\Mod{\Span{\MonC}}$ are monads in $\Span{\MonC}$, that is, by Proposition~\ref{prop:MonCatsasMonads}, monoidal categories. Fix any two of these objects, say categories $\catC$ and $\catD$ whose corresponding monads in $\Span{\MonC}$ have underlying spans
$$\vcenter{\xymatrix@R=8pt{ & \ar[dl]_{\dom{\catC}} \Ar{\catC} \ar[dr]^{\cod{\catC}} & \\ \Ob{\catC} && \Ob{\catC}}} \qquad \text{ and } \qquad \vcenter{\xymatrix@R=8pt{ & \ar[dl]_{\dom{\catD}} \Ar{\catD} \ar[dr]^{\cod{\catD}} & \\ \Ob{\catD} && \Ob{\catD}}}.$$
By definition, a 1-cell $\catC \tr{} \catD$ in $\Mod{\Span{\MonC}}$ is a span $\Ob{\catC} \tl{\dom{\FH}} \Ar{\FH} \tr{\cod{\FH}} \Ob{\catD}$ equipped with left and right actions 2-cells $\tau \: \catC \poi \FH \to \FH$ and $\rho \:  \FH \poi \catD \to \FH$ in $\Span{\MonC}$.
$$
 \xymatrix@R=10pt@C=17pt{
           &         & \ar[dl]_{} .\pushoutcorner \ar@/^2.5pc/[dddd]^{\tau}  \ar[dr]^{}  &  & \\
         & \ar[dl]_{\dom{\catC}}\Ar{\catC} \ar[dr]^{\cod{\catC}} &         & \ar[dl]_{\dom{\FH}}\Ar{\FH} \ar[dr]^{\cod{\FH}}& \\
  \Ob{\catC}  &        &  \Ob{\catC}  &        &  \Ob{\catD}  \\
 &&&&\\
     && \ar[uull]^{\dom{\FH}} \Ar{\FH} \ar[uurr]_{\cod{\FH}} &&
 } \qquad
  \xymatrix@R=10pt@C=17pt{
           &         & \ar[dl]_{} .\pushoutcorner \ar@/^2.5pc/[dddd]^{\rho}  \ar[dr]^{}  &  & \\
         & \ar[dl]_{\dom{\FH}}\Ar{\FH} \ar[dr]^{\cod{\FH}} &         & \ar[dl]_{\dom{\catD}}\Ar{\catD} \ar[dr]^{\cod{\catD}}& \\
  \Ob{\catC}  &        &  \Ob{\catD}  &        &  \Ob{\catD}  \\
 &&&&\\
     && \ar[uull]^{\dom{\FH}} \Ar{\FH} \ar[uurr]_{\cod{\FH}} &&
 }
 $$
One should think of $\Ar{\FH}$ as a set of arrows with source an object of $\catC$ and target one of $\catD$. Note that, a priori, $\Ob{\catC} \tl{} \Ar{\FH} \tr{} \Ob{\catD}$ does not define a category, as it is not supposed to carry a monad structure. Nonetheless, it will make notation easier to designate with $\tr{\in \FH}$ an element of $\Ar{\FH}$, as we do for arrows of a category.

Composites $\catC \poi \FH$ and $\FH \poi \catD$, source of $\tau$ and of $\rho$ respectively, are defined by pullback: as suggested in the proof of Proposition~\ref{prop:CatAsMonads}, the pullback object of $\catC \poi \FH$ should be regarded as the set of composable pairs $\tr{\in\catC}\tr{\in \FH}$ of arrows, and similarly for $\FH \poi \catD$. The action $\tau \: \catC \poi \FH \to \FH$ gives a way of pre-composing arrows of $\FH$ with arrows of $\catC$. Symmetrically, $\rho \: \FH \poi \catD \to \FH$ defines post-composition of arrows of $\FH$ with arrows of $\catD$. Compatibility conditions~\eqref{eq:actionscompatibilityconds} express that $\tau$ respects composition in $\catC$, $\rho$ respects composition in $\catD$ and that pre-/post-composition described by the two actions is associative.

To explain composition of 1-cells, fix an object $\catE = (\Ob{\catE} \tl{} \Ar{\catE} \tr{} \Ob{\catE})$ of $\Mod{\Span{\MonC}}$ and a 1-cell $\catD \tr{\FH'} \catE$, say with actions $\tau'$ and $\rho'$. Following the recipe~\eqref{eq:compBimodules}, the composite $\FH \bicomp{\catD} \FH'$ of $\FH$ and $\FH'$ is defined by the following coequaliser in $\Span{\MonC}$:
    \begin{equation*} \vcenter{
    \xymatrix{ \FH \poi \catD \poi \FH' \ar@<0.66ex>[r]^{\rho \FH'} \ar@<-0.66ex>[r]_{\FH \tau'} & \FH \poi \FH' \ar[r] & \FH \bicomp{\catD} \FH'}
    }
    \end{equation*}
    It is helpful to explain the definition of $\FH \bicomp{\catD} \FH'$ in terms of arrows. First, the carrier of $\FH \poi \catD \poi \FH'$ can be seen as the set of composable triples $\tr{\in \FH}\tr{\in \catD} \tr{\in \FH'}$. The action $\rho \FH'$ composes $\tr{\in \FH}$ and $\tr{\in \catD}$ to form an arrow of $\FH$, whereas the action $\FH \tau'$ composes $\tr{\in \catD}$ and $\tr{\in \FH'}$ to form an arrow of $\FH'$. Either ways we obtain a composable pair $\tr{\in \FH}\tr{\in \FH'}$. Equalizing these two actions amounts to quotient the set of pairs $\tr{\in \FH}\tr{\in \FH'}$ by the equivalence generated by the following relation:
\begin{equation}\label{eq:coeqt2t1} \tr{h}\tr{h'}\  {\equiv}_{\catD}\  \tr{g}\tr{g'} \quad \text{iff} \quad \text{there exist $\tr{d \in \catD}$ such that $\tr{h}\ =\rho(\tr{g}\tr{d})$ and $\tr{g'}\ =\tau'(\tr{d}\tr{h'})$.}\end{equation}
Therefore the 1-cell $\FH \bicomp{\catD} \FH'$ will be a span $\Ob{\catC} \tl{} \tr{} \Ob{\catE}$, whose carrier is the set of ${\equiv}_{\catD}$-equivalence classes of composable pairs $\tr{\in \FH}\tr{\in \FH'}$. We shall use the notation $[\tr{f}\tr{g}]_{{\equiv}_{\catD}}$ for the equivalence class with witness $\tr{f}\tr{g}$.
 \end{example}

\begin{remark}[Unit Laws] \label{rmk:unitlawMod} Let $\funF \xrightarrow{\FH} \G$ be a 1-cell in $\Mod{\Span{\MonC}}$. Since composition is weakly unital, there are isomorphisms
\begin{multicols}{2}\noindent
\begin{equation}\label{eq:idiso1}\FH \cong \funF \bicomp{\funF} \FH\end{equation}
\begin{equation}\label{eq:idiso2}\FH \cong \FH \bicomp{\G} \G\end{equation}
\end{multicols}
\noindent involving the identity 1-cells $\funF \xrightarrow{\funF} \funF$ and $\G \xrightarrow{\G} \G$. For later use it is useful to make explicit how these isomorphisms are defined. We focus on~\eqref{eq:idiso1}, the description of~\eqref{eq:idiso2} being analogous. Using the arrow view of 1-cells $\FH$ and $\funF \bicomp{\catD} \FH$ , the two directions of the iso~\eqref{eq:idiso1} are given by
\begin{eqnarray*}
\tr{f \in \FH} \quad \mapsto \quad [\tr{\id \in \funF}\tr{f\in \FH}]_{\equiv_{\funF}} & \qquad\qquad & [\tr{h \in \funF}\tr{g\in \FH}]_{\equiv_{\funF}}\quad \mapsto \quad\rho(\tr{h}\tr{g})
\end{eqnarray*}
 where $\rho$ is the right action of the bimodule $\FH$. The right-to-left direction is well-defined: the definition of $\equiv_{\funF}$ is given according to~\eqref{eq:coeqt2t1}, with left and right actions of $\funF$ being given by the multiplication $\mu^{\funF}$ of the monad $\funF$ in $\Span{\MonC}$. Then, compatibility of $\rho$ with $\mu^{\funF}$ guarantees that composable pairs which are equal modulo $\equiv_{\funF}$ are mapped into the same value by $\rho$.
We now check that the two mappings are invertible. First,
\[ \tr{f \in \FH} \quad \mapsto \quad [\tr{\id \in \funF}\tr{f\in \FH}]_{\equiv_{\funF}} \quad \mapsto \quad \rho(\tr{\id}\tr{f}) \ =\  \tr{f \in \FH} \]
because $\rho$ is compatible with the unit $\eta^{\funF}$ of the monad $\funF$, which is given by insertion of the identity arrow $\tr{\id}$ (see proof of Proposition~\ref{prop:CatAsMonads}). Conversely,
\[ [\tr{h \in \funF}\tr{g\in \FH}]_{\equiv_{\funF}} \quad\mapsto\quad \rho(\tr{h}\tr{g}) \quad\mapsto\quad [\tr{\id \in \funF}\tr{\rho(\tr{h}\tr{g})\in \FH}]_{\equiv_{\funF}} \ =\  [\tr{h \in \funF}\tr{g\in \FH}]_{\equiv_{\funF}}\]
because $\tr{h \in \funF}$ witnesses condition~\eqref{eq:coeqt2t1} for $\tr{h \in \funF}\tr{g\in \FH}$ and $\tr{\id \in \funF}\tr{\rho(\tr{h}\tr{g})}$, meaning that they are in the same $\equiv_{\funF}$- equivalence class.
\end{remark}

We now prove that any PROP yields a monad in $\Mod{\Span{\MonC}}$. We fist present an abstract approach (the same given in~\cite{Lack2004a}) and then also sketch a more direct argument. Let us write $\gls{Bicatxx}$ for the monoidal category of 1-cells $x\to x$ in $\bicat$ and $\gls{MonBC}$ for the category of monoids in a monoidal category $\catC$. It is useful to recall the following standard result (see e.g.~\cite{Lack2004a,ChengDistrLawsLT}).

\begin{proposition}\label{prop:purelyformalreasons} Fix a bicategory $\bicat$, $x \in \bicat$ and a monad $x \xrightarrow{\funF} x$. There is an equivalence
\begin{equation*}\MonB{\Mod{\bicat}(\funF,\funF)} \gls{equivCat} \funF / \MonB{\bicat(x,x)}.\end{equation*}
\end{proposition}
We instantiate Proposition~\ref{prop:purelyformalreasons} to the case in which $\bicat = \Span{\MonC}$, $x = \N$ and $\funF = \Perm$.
\begin{equation}\label{eq:propsaremonads}\MonB{\Mod{\Span{\MonC}}(\Perm,\Perm)} \simeq \Perm / \MonB{\Span{\MonC}(\N,\N)}.\end{equation}
We verify that PROPs live in the category on the right hand side.
Objects in $\MonB{\Span{\MonC}(\N,\N)}$ are monoids in ${\Span{\MonC}(\N,\N)}$, which are monad on $\N$ in $\Span{\MonC}$ and thus, by Proposition~\ref{prop:MonCatsasMonads} are precisely PROs. Morphisms in $\MonB{\Span{\MonC}(\N,\N)}$ are monoid homomorphisms in $\Span{\MonC}(\N,\N)$, thus are identity-on-objects monoidal functors between PROs, that is, PRO morphisms. Therefore, objects of the coslice $\Perm / \MonB{\Span{\MonC}(\N,\N)}$ are PRO morphisms with source $\Perm$. We can thereby conclude by Proposition~\ref{prop:coslicePROP} that PROPs are objects of the right hand side of \eqref{eq:propsaremonads}.

The left hand side of \eqref{eq:propsaremonads} tells us that PROPs are also monoids in $\Mod{\Span{\MonC}}(\Perm,\Perm)$, equivalently:

\begin{corollary}\label{cor:PROPsAreMonads} PROPs are monads on $\Perm$ in the bicategory $\Mod{\Span{\MonC}}$.\end{corollary}

\begin{remark} \label{rmk:MonadsNotPROPs} Differently from the case of small categories (Proposition~\ref{prop:CatAsMonads}), Corollary~\ref{cor:PROPsAreMonads} does not give a complete characterisation for PROPs: there are monads on $\Perm$ in $\Mod{\Span{\MonC}}$ which do \emph{not} correspond to any PROP. This is because, as we noticed at the end of \S~\ref{sec:props}, not all the objects of $\Perm / \PRO$ are PROPs.
\end{remark}

Although Corollary~\ref{cor:PROPsAreMonads} immediately follows by Proposition~\ref{prop:purelyformalreasons}, it is illuminating to sketch a direct argument for its statement.
\begin{proof}[Proof of Corollary~\ref{cor:PROPsAreMonads}]  Being a monoidal category with set of objects $\N$, the PRO $\Perm$ yields a monad in $\Span{\MonC}$ on $\N$, as shown in Proposition~\ref{prop:MonCatsasMonads}.
$$\xymatrix@R=8pt{ & \ar[dl]_{\dom{\Perm}} \Ar{ \Perm} \ar[dr]^{\cod{\Perm}} & \\ \N  && \N }$$
Starting now from a PROP $\T$, we shall define a monad on $\Perm$ in $\Mod{\Span{\MonC}}$. The underlying 1-cell is a span $\N \tl{} \Ar{\T} \tr{}\N$
  $$\xymatrix@R=8pt{
  & \ar[dl]_{\dom{\Perm}} \Ar{\Perm} \ar[dr]^{\cod{\Perm}} & & \ar[dl]_{\dom{\T}} \Ar{\T} \ar[dr]^{\cod{\T}} & & \ar[dl]_{\dom{\Perm}} \Ar{\Perm} \ar[dr]^{\cod{\Perm}} &\\
  \N  && \N && \N && \N
  }$$
 whose carrier $\Ar{\T}$ is the set of arrows of $\T$ and $\dom{\T}$, $\cod{\T}$ are the monoid homomorphisms for source and target. For convenience, we shall call $\T$ also the span $\N \tl{} \Ar{\T} \tr{}\N$. To be a 1-cell in $\Mod{\Span{\MonC}}$, $\T$ should carry a bimodule structure, meaning that it is equipped with 2-cells $\tau$ and $\rho$, respectively left and right action, making the following diagrams commute. 
 $$
 \xymatrix@R=10pt@C=17pt{
           &         & \ar[dl]_{} .\pushoutcorner \ar@/^2.5pc/[dddd]^{\tau}  \ar[dr]^{}  &  & \\
         & \ar[dl]_{\dom{\Perm}}\Ar{\Perm} \ar[dr]^{\cod{\Perm}} &         & \ar[dl]_{\dom{\T}}\Ar{\T} \ar[dr]^{\cod{\T}}& \\
  \N  &        &  \N  &        &  \N  \\
 &&&&\\
     && \ar[uull]^{\dom{\T}} \Ar{\T} \ar[uurr]_{\cod{\T}} &&
 } \qquad
  \xymatrix@R=10pt@C=17pt{
           &         & \ar[dl]_{} .\pushoutcorner \ar@/^2.5pc/[dddd]^{\rho}  \ar[dr]^{}  &  & \\
         & \ar[dl]_{\dom{\T}}\Ar{\T} \ar[dr]^{\cod{\T}} &         & \ar[dl]_{\dom{\Perm}}\Ar{\Perm} \ar[dr]^{\cod{\Perm}}& \\
  \N  &        &  \N  &        &  \N  \\
 &&&&\\
     && \ar[uull]^{\dom{\T}} \Ar{\T} \ar[uurr]_{\cod{\T}} &&
 }
 $$
  Since any PROP contains the permutations (as symmetries), we can regard $\Ar{\Perm}$ as a subset of $\Ar{\T}$ and define $\tau \: \Perm \poi \T \to \T$ and $\rho \: \T \poi \Perm \to \T$ respectively by pre- and post-composition of arrow in $\Ar{\T}$ with arrows in $\Ar{ \Perm}$. The compatibility conditions~\eqref{eq:actionscompatibilityconds} correspond to pre-/post-composition being compatible with composition in $\Perm$ and being associative.

 It remains to equip the bimodule $\T$ with a monad structure.
For the multiplication, note that, differently from the case of plain categories (Proposition~\ref{prop:CatAsMonads}), $\mu \: \T \bicomp{\Perm} \T \to \T$ acts on \emph{equivalence classes} of pairs of composable arrows $\tr{\in \T}\tr{\in \T}$, because of the way the composite $\T \bicomp{\Perm} \T$ is defined (Example~\ref{ex:ModSpan}). Nonetheless, defining $\mu$ by composition in $\T$ still yields a 2-cell: indeed, in the notation of \eqref{eq:coeqt2t1},
\begin{eqnarray*}
\tr{f \in \T}\tr{g \in \T} \ {\equiv}_{\Perm} \ \tr{f' \in \T}\tr{g' \in \T} & \Rightarrow & \exists \tr{p \in \Perm}\text{ s.t. } \tr{f \in \T}\ = \ \rho(\tr{f' \in \T}\tr{p \in \Perm}) \text{ and }\tr{g' \in \T} \ = \ \tau(\tr{p \in \Perm}\tr{g \in \T}) \\
& \Rightarrow & \exists \tr{p \in \Perm} \text{ s.t. }  \tr{f \in \T} \ = \ \tr{f' \poi p\  \in \T} \text{ and }\tr{g' \in \T}\ = \ \tr{p \poi g\  \in \T} \\
& \Rightarrow & \tr{f \poi g \  \in \T} \ = \ \tr{f' \poi p \poi g \  \in \T} \ = \ \tr{f' \poi g'\  \in \T} \\
& \Rightarrow & \mu(\tr{f\in \T}\tr{g \in \T}) = \mu(\tr{f'\in \T}\tr{g' \in \T}).
\end{eqnarray*}
With this definition of $\mu$, the monad law~\eqref{diag:lawsmonads2} holds by associativity of composition.

 The case of unit is more subtle. Indeed, whereas in $\Span{\MonC}$ identity 1-cells are given by identity spans, in $\Mod{\Span{\MonC}}$ the identity 1-cell on the object $\N \tl{} \Ar{\Perm} \tr{} \N$ is the object itself, now regarded as a bimodule with actions given by composition in $\Perm$. The unit $\eta \: \Perm \to \T$ will then be a span morphism
 $$
 \xymatrix@=6pt{
 & \ar[dl] \Ar{\Perm}\ar[dr] \ar[dd]^{\eta} & \\
 \N && \N \\
& \ar[ul] \Ar{\T} \ar[ur] &
} $$
defined by interpreting a permutation as an arrow of the PROP $\T$ (in fact, this is the mapping given by initiality of $\Perm$ in $\PROP$, see~\S\ref{sec:props}). This definition of $\eta$ satisfies the monad law~\eqref{diag:lawsmonads1}. 
\end{proof}

\subsection{Distributive Laws of PROPs}\label{sec:distrLawPROPs}

Now that we have an understanding of PROPs as monads we can compose them via distributive laws. Fix PROPs $\T_1$, $\T_2$, seen as monads in $\Mod{\Span{\MonC}}$, say with actions $\lefta{}_1, \righta{}_1$ and $\lefta{}_2$, $\righta{}_2$ respectively. Let $\lambda \: \T_2 \bicomp{\Perm} \T_1 \to \T_1 \bicomp{\Perm} \T_2$ be a distributive law between them. Recall that, by definition, the composite $\T_2 \bicomp{\Perm} \T_1$ is a 1-cell whose carrier as a span has elements composable pairs $\tr{f\in \T_2}\tr{g\in \T_1}$ of arrows subject to the following equivalence relation ${\equiv}_{\Perm}$, obtained by instantiating \eqref{eq:coeqt2t1} to 1-cells $\T_2$, $\Perm$ and $\T_1$\footnote{Note that, for a generic $\catD$ as in~\eqref{eq:coeqt2t1}, we quotient by the equivalence relation generated by the relation $\equiv_{\catD}$. For $\equiv_{\Perm}$, the two coincide by self-duality of $\Perm$.}:
\begin{equation}
\begin{aligned}\label{eq:quotPerm}
\tr{f}\tr{g} \ {\equiv}_{\Perm} \ \tr{f'}\tr{g'} \  \ & \text{ iff } & \text{there is $\tr{p \in \Perm}$ such that $\tr{f}\ =\righta{}_2(\tr{f'}\tr{p})$ and $\tr{g'}\ =\lefta{}_1(\tr{p}\tr{g})$} \\
 & \text{ iff } & \text{there is $\tr{p \in \Perm}$ making }
\vcenter{
            \xymatrix@R=13pt@C=23pt{
                & \ar[dr]^{g} & \\
               \ar[ur]^{f} \ar[r]_{f'} &\ar[u]_-{p} \ar[r]_{g'} &
            }
}
\text{ commute.}
\end{aligned}
\end{equation}
 A perhaps more illuminating way of phrasing condition~\eqref{eq:quotPerm} is by saying that, in presence of a triple $\tr{\in \T_2}\tr{\in \Perm}\tr{\in \T_1}$, the choices of letting $\tr{\in \Perm}$ be part of $\T_2$ or of $\T_1$ determine the same element of $\T_2 \bicomp{\Perm} \T_1$. An analogous description applies to $\T_1 \bicomp{\Perm} \T_2$. Therefore, we can present $\lambda$ as a mapping of arrows $\tr{\in \T_2}\tr{\in \T_1}$ to arrows $\tr{\in \T_1}\tr{\in \T_2}$: condition~\eqref{eq:quotPerm} expresses that $\lambda$ does not discriminate between $(\tr{\in \T_2}\tr{\in \Perm})\tr{\in \T_1}$, where the middle arrow $\tr{\in \Perm}$ is considered as part of $\T_2$, and $\tr{\in \T_2}(\tr{\in \Perm}\tr{\in \T_1})$, where it is considered as part of $\T_1$. 

We remark that $\lambda$ respects the PROP structure by definition. First, it preserves identity and composition by~\eqref{eq:distrlawequationsUnit}-\eqref{eq:distrlawequationsMult}. Compatibility with the monoidal product is guaranteed by $\lambda$ being a morphism between spans in $\MonC$. Finally, $\lambda$ behaves well with respect to the symmetry structure of $\T_1$ and $\T_2$: this is because, being a 2-cell in $\Mod{\Span{\MonC}}$, $\lambda$ is compatible with left and right action of the bimodules $\T_1$ and $\T_2$ (see~\eqref{eq:bimodulemorphism}). 

In conclusion, $\lambda$ yields a PROP $\T_1 \bicomp{\Perm} \T_2$ defined as follows.
\begin{itemize}[noitemsep,topsep=0pt,parsep=0pt,partopsep=0pt]
\item Arrows of $\T_1 \bicomp{\Perm} \T_2$ are composable pairs $\tr{\in \T_1}\tr{\in \T_2}$, identified when they are equal up-to permutation in the way described by \eqref{eq:quotPerm}.
 \item Following~\eqref{eq:defcomposedmonad}, the composite of $\tr{f\in \T_1}\tr{g\in \T_2}$ and $\tr{h\in \T_1}\tr{i\in \T_2}$ is $\tr{f \poi h' \in \T_1}\tr{g' \poi i \in \T_2}$, where $\tr{h'\in \T_1}\tr{g'\in \T_2}$ has been obtained by applying $\lambda$ to $\tr{g\in \T_2}\tr{h\in \T_1}$.
\end{itemize}

\begin{remark}\label{rmk:PROPsclosedunderdistrlaw} As we noticed in Remark~\ref{rmk:MonadsNotPROPs}, not all the monads in $\Mod{\Span{\MonC}}$ on $\Perm$ are PROPs. Therefore, to define PROP composition in a sensible way one should guarantee that a distributive law of PROPs yields a monad which is again a PROP. It is not hard to check that this is indeed the case: the key observation is that the permutations in $\T_1 \bicomp{\Perm}\T_2$ are exactly those in $\T_1$ and $\T_2$, equalized via~\eqref{eq:coequalizer}. Now, permutations satisfy the naturality requirements w.r.t. to the arrows in $\T_1$ and $\T_2$, because those two are PROPs. It follows that they satisfy the same property w.r.t. to the arrows $\tr{\in \T_1}\tr{\in\T_2}$ of $\T_1 \bicomp{\Perm}\T_2$.
\end{remark}

\begin{example}~ \label{ex:composingSemPROP}
    \begin{enumerate}[(a)]
    \item We describe the PROP $\F$ of functions as the composite of PROPs for surjections and injections. Let $\gls{PROPInj}$ be the PROP whose arrows $n \to m$ are injective functions from $\ord{n}$ to $\ord{m}$. The PROP $\gls{PROPSurj}$ of surjective functions is defined analogously. Epi-mono factorisation of functions gives a mapping of composable pairs $\tr{\in \Inj}\tr{\in \Surj}$ to composable pairs $\tr{\in \Surj}\tr{\in \Inj}$. This mapping extends to $\equiv_{\Perm}$-equivalence classes: if there is $p\in\Perm$ making
        \[
        \vcenter{
            \xymatrix@R=13pt@C=23pt{
                & \ar[dr]^{s} & \\
               \ar[ur]^{i} \ar[r]_{i'} &\ar[u]_-{p} \ar[r]_{s'} &
            }
            }
            \]
        commute, then $\tr{i}\tr{s}$ and $\tr{i'}\tr{s'}$ have the same epi-mono factorisation up-to permutation (that means, all their factorisations $\tr{\in \Surj}\tr{\in \Inj}$ are in the same $\equiv_{\Perm}$-class).

        This mapping yields a 2-cell $\lambda \: \Inj \bicomp{\Perm} \Surj \to \Surj \bicomp{\Perm} \Inj$ satisfying the equations of distributive laws~\cite{Lack2004a}.
        The resulting PROP $\Surj \bicomp{\Perm} \Inj$ is isomorphic to $\F$ because any function in $\F$ can be uniquely factorised (up-to permutation) as a surjection followed by an injection. From a different perspective, this result tells us that $\F$ can be \emph{decomposed} into simpler PROPs $\Surj$ and $\Inj$. \label{ex:composingSemPROP1}
    \item We describe the PROP of spans in $\F$ as the result of a distributive law defined by pullback in $\F$. First, note that a composable pair $\tr{f \in \F}\tr{g\in \Fop}$ is the same thing as a cospan $\tr{f} \tl{g}$ in $\F$. Dually, pairs $\tr{\in \Fop}\tr{\in \F}$ yield spans in $\F$. Define a 2-cell $\lambdapb \: \F \bicomp{\Perm} \Fop \to \Fop \bicomp{\Perm} \F$ as the mapping of a cospan $\tr{f} \tl{g}$ to its pullback span $\tl{g'} \tr{f'}$. This definition respects $\equiv_{\Perm}$-equivalence. Indeed:
        \begin{itemize}
        \item $\tr{f \in \F}\tr{g \in \Fop}\ \equiv_{\Perm}\ \tr{f' \in \F}\tr{g' \in \Fop}$ means that $\tr{f}\tl{g}$ and $\tr{f'}\tl{g'}$ are isomorphic cospans and thus they are pulled back by isomorphic spans.
        \item Isomorphisms in $\F$ coincide with permutations in $\Perm$.
        \end{itemize}
        Moreover, $\lambdapb$ satisfies the equations of distributive laws~\cite{Lack2004a}. This yields a PROP $\Fop \bicomp{\Perm} \F$ whose arrows $n \to m$ are equivalence classes of spans $n\tl{g'} \tr{f'}m$ in $\F$. Following~\eqref{eq:quotPerm}, two spans $n\tl{g_1} z\tr{f_1}m$ and $n\tl{g_2}z \tr{f_2}m$ are identified as arrows of $\Fop \bicomp{\Perm} \F$ whenever there is a permutation $p$ (i.e., an isomorphism in $\F$) making the following diagram commute
         $$
         \xymatrix@R=6pt@C=12pt{
         & \ar[dl]_{g_1} z\ar[dr]^{f_1} \ar[dd]^{p} & \\
         n && m \\
        & \ar[ul]^{g_2} z \ar[ur]_{f_2} &
        } $$
         that means, when $\tl{g_1} \tr{f_1}$ and $\tl{g_2} \tr{f_2}$ are isomorphic spans. By~\eqref{eq:defcomposedmonad}, composition in $\Fop \bicomp{\Perm} \F$ is defined in terms of $\lambdapb$, thus is by pullback. In the terminology of~\cite{BenabouBicategories}, one can see $\Fop \bicomp{\Perm} \F$ as the \emph{classifying category} of the bicategory $\Span{\F}$, obtained by identifying the isomorphic $1$-cells and forgetting the $2$-cells. \label{ex:composingSemPROP2}
    \item Dually, there exist a distributive law $\lambdapo \: \Fop \bicomp{\Perm} \F \to \F \bicomp{\Perm} \Fop$ defined by pushout in $\F$~\cite{Lack2004a}. The composite PROP $\F \bicomp{\Perm} \Fop$ is the classifying category of $\glssymbol{Cospan}(\F)$,  the bicategory where 1-cells are cospans in $\F$ and composition is by pushout. \label{ex:composingSemPROP3}
        \end{enumerate}
\end{example}

\paragraph{Composing free PROPs} We now turn our attention to the case in which the PROPs $\T_1$ and $\T_2$ involved in the composition $\lambda \: \T_2 \bicomp{\Perm} \T_1 \to \T_1 \bicomp{\Perm} \T_2$ are generated by SMTs, say $(\Sigma_1,E_1)$ and $(\Sigma_2,E_2)$ respectively. It turns out that also $\T_1 \bicomp{\Perm} \T_2$ enjoys a presentation by generators and equations, which we now describe in steps. First, by definition a composable pair $\tr{f \in \T_1}\tr{g \in \T_2}$ consists of a $\Sigma_1$-term $f$ modulo $E_1$ followed by a $\Sigma_2$-term $g$ modulo $E_2$. One can then see $\tr{f}\tr{g}$ as a $\Sigma_1\uplus\Sigma_2$-term $f \poi g$ modulo $E_1 \uplus E_2$. Motivated by this observation, we shall take $\Sigma_1\uplus\Sigma_2$ as the signature for $\T_1 \bicomp{\Perm} \T_2$. As observed in \S~\ref{sec:coproduct}, terms $n \to m$ generated by $\Sigma_1\uplus\Sigma_2$ can be regarded as sequences of composable $\Sigma_1$- and $\Sigma_2$-terms, here represented by blue and red arrows respectively:
$$\xymatrix{n\ar@[red][r] &\ar@[blue][r] &\ar@[red][r] &\ar@[blue][r] & \dots \ar@[red][r] &\ar@[blue][r] &m}$$
Which equations we shall impose to put these sequences in 1-1 correspondence with the arrow $[\tr{\in \T_1}\tr{\in \T_2}]_{\equiv_{\Perm}}$ of $\T_1  \bicomp{\Perm} \T_2$?
The key is to read the graph of $\lambda$ as a set $E_{\lambda}$ of (directed) equations between $\Sigma_1\uplus\Sigma_2$-terms modulo $E_1 \uplus E_2$, calculated as follows:
\begin{itemize}[noitemsep,topsep=2pt,parsep=2pt,partopsep=2pt]
\item suppose that $\lambda$ maps the arrow $[\tr{g \in \T_2}\tr{f \in \T_1}]_{\equiv_{\Perm}}$ of $\T_2  \bicomp{\Perm} \T_1$ to the arrow $[\tr{f' \in \T_1}\tr{g' \in \T_2}]_{\equiv_{\Perm}}$ of $\T_1  \bicomp{\Perm} \T_2$. Then, put the equation $g \poi f \feq f'\poi g'$ in $E_{\lambda}$.
\end{itemize}
We can now use the equations in $E_{\lambda}$ to rewrite any $\Sigma_1\uplus\Sigma_2$-term into one of the shape $f \poi g$ --- in our graphical representation, a term where all red arrows precede any blue arrow:
$$ \def \lambdaarrow {\ar@{=>}[]+<0ex,-2ex>;[d]!<0ex,2ex>^{\lambda}}
 \def \lambdaarrowtwo {\ar@{=>}[]+<0ex,2ex>;[d]!<0ex,5ex>^{\lambda}}
\xymatrix@C=40pt@R=20pt{
n \ar@[red][r] & \ar@[blue][r] \ar@[red][dr] & \lambdaarrow \ar@[red][r] & \ar@[blue][r] \ar@[red][dr] & \lambdaarrow \ar@[red][r] & \dots \ar@[red][dr] \ar@[blue][r] \ar@[red][dr] & \lambdaarrow \ar@[red][r] & \ar@[blue][r] & m \\
& & \ar@[blue][ur] \ar@[red][dr] & \lambdaarrowtwo &  \ar@[blue][ur] \ar@[red][dr] &  \lambdaarrowtwo & \ar@[blue][ur]& &\\
& & & \dots \ar@[blue][ur] \ar@[red][dr] & \lambdaarrowtwo & \dots \ar@[blue][ur]  & & &\\
& & & &\ar@[blue][ur] & & & & &}$$

Compatibility of $\lambda$ with unit and multiplication of the monads $\T_1$ and $\T_2$ (see~\eqref{eq:distrlawequationsUnit}-\eqref{eq:distrlawequationsMult}), that is, identity and composition of $\T_1$ and $\T_2$ as categories, guarantees that any different rewriting reaching the form $f\poi g$ yields a term of the same equivalence class under $E_1 \uplus E_2\uplus E_{\lambda}$.
Therefore, arrows of $\T_1 \bicomp{\Perm} \T_2$ are the same thing as $\Sigma_1\uplus\Sigma_2$-terms modulo the equations $E_1 \uplus E_2\uplus E_{\lambda}$. We fix our conclusion with the following statement.
\begin{proposition}[{\cite[Prop. 4.7]{Lack2004a}}] \label{prop:SMTforComposition} Let $(\Sigma_1,E_1)$ and $(\Sigma_2,E_2)$ be SMTs generating PROPs $\T_1$ and $\T_2$ respectively. Suppose there is a distributive law $\lambda \: \T_2  \bicomp{\Perm} \T_1 \to \T_1  \bicomp{\Perm} \T_2$, yielding a set of equations $E_{\lambda}$ in the way described above. Then the SMT $(\Sigma_1\uplus\Sigma_2, E_1 \uplus E_2\uplus E_{\lambda})$ presents the PROP $\T_1  \bicomp{\Perm} \T_2$.
\end{proposition}

  More suggestively, one can read Proposition~\ref{prop:SMTforComposition} as saying that $\T_1 \bicomp{\Perm} \T_2$ is the quotient of the sum $\T_1 + \T_2$ under the equations $E_{\lambda}$ encoded by $\lambda$. This agrees with the intuition that composing PROPs amounts to expressing, in the form of a distributive law, \emph{compatibility conditions} between the algebraic structures that the PROPs describe.

\begin{example}~ \label{ex:distrlawsyntactic}
\begin{enumerate}[(a)]
\item \label{ex:distrlawsyntactic1} We show how the PROP $\Mon$ of commutative monoids can be factorised as the composite $\gls{PROPMult} \bicomp{\Perm} \gls{PROPUnit}$. Here $\MULT$ is the PROP freely generated by the theory of a commutative multiplication, that is, the SMT with signature $\{\Wmult\}$ and equations \eqref{eq:wmonassoc}-\eqref{eq:wmoncomm}. The PROP $\UNIT$ is freely generated by the theory of the unit, consisting of a signature $\{\Wunit\}$ and no equations. As we did for $\Mon$, we can view the arrows of $\MULT$ and $\UNIT$ as representing the graphs of functions on ordinals. It is easy to verify that arrows of $\MULT$ are in bijective correspondence with surjective functions and arrows of $\UNIT$ with injections, that is $\MULT \cong \Surj$ and $\UNIT \cong \Inj$. By this correspondence, a distributive law $\lambda \: \UNIT \bicomp{\Perm} \MULT \To \MULT \bicomp{\Perm} \UNIT$ is yielded by the one of type $\Inj \bicomp{\Perm} \Surj \to \Surj \bicomp{\Perm} \Inj$ described in Example~\ref{ex:composingSemPROP}\ref{ex:composingSemPROP1}. By Proposition~\ref{prop:SMTforComposition}, $\MULT \bicomp{\Perm} \UNIT$ is presented by the signature $\{\Wmult, \Wunit\}$ and equations \eqref{eq:wmonassoc}-\eqref{eq:wmoncomm} plus the set $E_{\lambda}$ of equations encoded by $\lambda$. 
    In fact, we can easily give a \emph{finite} presentation for $\MULT \bicomp{\Perm} \UNIT$: since $\MULT \bicomp{\Perm} \UNIT \cong \Surj \bicomp{\Perm} \Inj \cong \F \cong \Mon$, the SMT generating $\MULT \bicomp{\Perm} \UNIT$ is just the one of commutative monoids, meaning that from \eqref{eq:wmonassoc}, \eqref{eq:wmoncomm} and \eqref{eq:wmonunitlaw} one can derive all the equations in $E_{\lambda}$. Note that, since \eqref{eq:wmonassoc} and \eqref{eq:wmoncomm} are already part of the SMT for $\MULT$, \eqref{eq:wmonunitlaw} is the one responsible of encoding $\lambda$: it tells how the only generator $\Wunit$ of $\UNIT$ distributes over the only generator $\Wmult$ of $\MULT$. Also, the type of $\lambda$ suggests a \emph{left-to-right} orientation for \eqref{eq:wmonunitlaw}:
    \begin{eqnarray*}
\lower7pt\hbox{$\includegraphics[height=1.3cm]{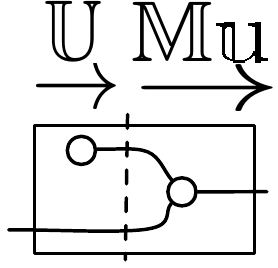}$}
& \Rightarrow &
\lower7pt\hbox{$\includegraphics[height=1.3cm]{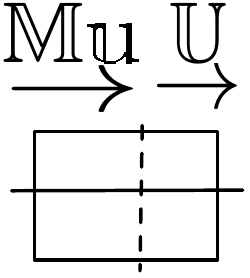}$}.
\end{eqnarray*}
\item \label{ex:distrlawsyntactic2} We now show how composing $\Com$ and $\Mon$ yields the PROP $\B$ of bialgebras (Example~\ref{ex:equationalprops}). Since $\Mon \cong \F$ and $\Com \cong \F^{\op}$, a distributive law $\lambdapb \: \Mon \bicomp{\Perm} \Com \To \Com \bicomp{\Perm} \Mon$ is yielded by the one of type $\F \bicomp{\Perm} \Fop \to \Fop \bicomp{\Perm} \F$ described in Example~\ref{ex:composingSemPROP}\ref{ex:composingSemPROP2}. By Proposition~\ref{prop:SMTforComposition}, the PROP $\Com \bicomp{\Perm} \Mon$ is presented by the generators and equations of $\Com +\Mon$ together with the equations $E_{\lambda}$ encoded by $\lambdapb$. By definition of $\lambdapb$, one can ``read off'' equations in $E_{\lambda}$ from pullback squares in $\F$, for instance:
\[
\xymatrix@R=10pt@C=15pt{
& 1 & && & & 1 \ar[dr]^{\Bcounit} &\\
\ar[ur]^{} 2 & & 0 \ar[ul]_{} & \ar@{|=>}[r] & & \ar[ur]^{\Wmult} 2 \ar[dr]_{\twoBcounit} & & 0 \\
& \ar[ul]^{} 0 \pullbackcorner \ar[ur]_{} & & && & 0 \ar[ur]_{\ZeronetT} & }
\quad
\lower20pt\hbox{$
\text{yields}
\quad
\Wmult ; \Bcounit = \twoBcounit \poi \ZeronetT
$}
\]
\noindent
where the arrows in the pullback diagram are given by initiality of $0$/finality of $1$ and the second diagram is obtained from the pullback by applying the isomorphisms  $\F \cong \Mon$ and
$\Fop \cong \Com$.
The above equation is just \eqref{eq:unitsl}. In fact, \eqref{eq:unitsl}-\eqref{eq:bwbone} yield a \emph{sound and complete} axiomatisation for $E_{\lambdapb}$: each of these four equations corresponds to a pullback square in the above sense (soundness), and the four of them together with the equations in $\MULT + \UNIT$ suffice to derive any equation associated with a pullback in $\F$ (completeness)~\cite{Lack2004a}. 
We can thereby conclude that $\Com \bicomp{\Perm} \Mon$ is isomorphic to the PROP $\B$ of bialgebras. From the point of view of $\B$, this characterisation brings two properties: first, it gives a \emph{decomposition} of any arrow $n \tr{\in \B} m$ as $n \tr{\in \Com}\tr{\in \Mon}m$; second, the isomorphism $\B \cong \Fop \bicomp{\Perm} \F$ can be seen as a \emph{semantics} interpreting any string diagram of $\B$ as a span of functions.

Later in this chapter (\S~\ref{sec:lawvere}) we will give yet another perspective on the characterisation $\B \cong \Com \bicomp{\Perm} \Mon$, reading it as the statement that $\B$ is the \emph{Lawvere theory} for commutative monoids. In that approach, the isomorphism $\B \cong \Com \bicomp{\Perm} \Mon$ will follow as a corollary of a more general result about Lawvere theories as composed PROPs --- Theorem~\ref{Th:LawvereCompositePROP}.
\item \label{ex:distrlawsyntactic3} Instead of pullbacks, one may combine $\Mon$ and $\Com$ via the distributive law defined by pushout in $\F$ of Example~\ref{ex:composingSemPROP}\ref{ex:composingSemPROP3}. The induced distributive law $\lambdapo$ on the freely generated PROPs will have the type $\Com \bicomp{\Perm} \Mon\To \Mon \bicomp{\Perm} \Com$. Just as in the above case, the resulting PROP $\Mon \bicomp{\Perm} \Com$ turns out to be finitely axiomatisable, being presented by the equations~\eqref{eq:BWFrob}-\eqref{eq:BWSep} of special Frobenius algebras~\cite{Lack2004a}. As a consequence, the PROP $\FROB$ introduced in Example~\ref{ex:equationalprops} is isomorphic to $\Mon \bicomp{\Perm} \Com$, each arrow $n \tr{\in \FROB} m$ enjoys a decomposition as $n\tr{\in \Mon}\tr{\in \Com}m$ and a full and faithful interpretation as a cospan of functions. In particular, the factorisation property is widely used in categorical quantum mechanics where is known as the \emph{spider theorem}--- see e.g.~\cite{Coecke2008,CoeckeDuncanZX2011}. 
\end{enumerate}
\end{example}

\subsection{Distributive Laws by Pullback and Pushout}\label{sec:distrLawPullback}

Distributive laws between PROPs (\S \ref{sec:composingprops}) allow for more flexibility than distributive laws between categories (\S \ref{sec:composingCats}): factorisations are not defined on-the-nose but just up-to permutation.

Nonetheless, this situation is still quite restrictive. Suppose $\T$ is a PROP with pullbacks. Ideally, as in Example~\ref{ex:composingSemPROP}\ref{ex:composingSemPROP2}, we would like to generate the PROP of spans in $\T$ via a distributive law $\lambdapb \: \T \bicomp{\Perm} \T^{\op} \to \T^{\op} \bicomp{\Perm} \T$ defined by pullback. Symmetrically, if $\T$ has pushouts, as in Example~\ref{ex:composingSemPROP}\ref{ex:composingSemPROP3} we may want to regard the PROP of cospans in $\T$ as arising by a distributive law $\lambdapo \: \T^{\op} \bicomp{\Perm} \T \to \T \bicomp{\Perm} \T^{\op}$ defined by pushout. The above examples show that $\lambdapb$ and $\lambdapo$ are indeed distributive laws when defined on the PROP $\F$. However, this is a particularly fortunate situation: (co)limits in $\F$ happen to be uniquely defined up-to permutation, as the permutations in this category coincide with the isomorphisms. In general, there is no reason why the isomorphisms in $\T$ should all be permutations, implying that mappings like $\lambdapb$ and $\lambdapo$ may not even be well-defined on $\equiv_{\Perm}$-equivalence classes: given two cospans equal modulo-$\equiv_{\Perm}$, their pullback spans are isomorphic but this iso is not necessarily a permutation, i.e. they can be in a different $\equiv_{\Perm}$-equivalence class.
We propose an approach that allows to define distributive laws by pullback (and pushout) for an arbitrary PROP $\T$ with pullbacks (and pushouts). The key observation is that, when composing PROPs as 1-cells in $\Mod{\Span{\MonC}}$, we are only identifying $\Perm$ as common structure between them. Instead, for distributive laws involving $\T$ and $\T^{\op}$, potentially we could identify more structure that is shared by the two PROPs, as for instance the sub-PROP $\gls{PROPCore}$ (called the \emph{core} of $\T$) whose arrows are the isomorphisms in $\T$. Formally, this will amount to view PROPs $\T$, $\T^{\op}$ not as monads on $\Perm$ but rather on $\PJ$. Their composition in $\Mod{\Span{\MonC}}$ is now defined by tensoring over $\PJ$, implying that objects of $\T \bicomp{\PJ} \T^{\op}$ and of $\T^{\op} \bicomp{\PJ} \T$ are composable pairs of arrows equal up-to $\PJ$, i.e., up-to iso in $\T$. It is then possible to meaningfully define distributive laws by pullback and pushout between these monads.

\begin{lemma}\label{lemma:PROPsMonadsOverIso} Let $\T$ be a PROP and $\PJ$ the core of $\T$. Then $\T$ and $\T^{\op}$ are monads on $\PJ$ in $\Mod{\Span{\MonC}}$.
\end{lemma}
\begin{proof}
The argument follows closely the one given for Corollary~\ref{cor:PROPsAreMonads}. Given a PROP $\T$, the corresponding bimodule has underlying span $\Ob{\T} \tl{} \Ar{\T} \tr{} \Ob{\T}$, with legs given by source and target of arrows in $\T$. Actions $\lambda \: \PJ \poi \T \to \T$ and $\rho \: \T \poi \PJ \to \T$ are defined by composition in $\T$, any arrow in $\PJ$ being also one in $\T$. The monad structure is as follows: the multiplication $\mu \: \T\bicomp{\PJ} \T \to \T$ is defined by composition and the unit $\eta \: \PJ \to \T$ by interpreting an arrow of $\PJ$ as one of $\T$.

The construction of a monad starting from the PROP $\T^{\op}$ is completely analogous.
\end{proof}

\begin{proposition}\label{prop:distrLawPbPo} Let $\T$ be a PROP and $\PJ$ the core of $\T$.
\begin{enumerate}
    \item If $\T$ has pullbacks, there is a distributive law $\lambdapb \: \T \bicomp{\PJ} \T^{\op} \to \T^{\op} \bicomp{\PJ} \T$, defined by pullback, whose resulting PROP $\T^{\op}\bicomp{\PJ} \T$ is the classifying category of $\Span{\T}$.
    \item If $\T$ has pushouts, there is a distributive law $\lambdapo \: \T^{\op} \bicomp{\PJ} \T \to \T \bicomp{\PJ} \T^{\op}$, defined by pushout, whose resulting PROP $\T\bicomp{\PJ} \T^{\op}$ is the classifying category of $\Cospan{\T}$.
\end{enumerate}
\end{proposition}
\begin{proof} We confine ourselves to the first statement, the argument for the second being completely analogous. We regard $\T$ and $\T^{\op}$ as monads on $\PJ$ in $\Mod{\Span{\MonC}}$, as allowed by Lemma~\ref{lemma:PROPsMonadsOverIso}. The underlying span $\T \bicomp{\PJ} \T^{\op}$ has carrier consisting of composable pairs of arrows $\tr{\in \T}\tr{\in \T^{\op}}$, quotiented by the equivalence relation:
\begin{equation}\label{eq:UpToIso}
\tr{f}\tr{g} \ {\equiv}_{\PJ} \ \tr{f'}\tr{g'} \qquad \text{ iff } \qquad \text{there is $\tr{i \in \PJ}$ making }
\vcenter{
            \xymatrix@R=13pt@C=23pt{
                & \ar[dr]^{g} & \\
               \ar[ur]^{f} \ar[r]_{f'} & \ar[u]_-{i} \ar[r]_{g'} &
            }
}
\text{ commute.}
\end{equation}
The 1-cell $\T^{\op} \bicomp{\PJ} \T$ enjoys a symmetric description. We let $\lambdapb \: \T \bicomp{\PJ} \T^{\op} \to \T^{\op} \bicomp{\PJ} \T$ be defined by pullback. This is indeed a well-defined 2-cell: if two pairs of arrows $\tr{f \in \T}\tr{g \in \T^{\op}}$ and $\tr{f'\in \T}\tr{g'\in \T^{\op}}$ are in the same $\equiv_{\PJ}$-class then $\tr{f}\tl{g}$ and $\tr{f'}\tl{g'}$ are isomorphic cospans, meaning that they are mapped on isomorphic pushout spans.


It remains to check that $\lambdapb$ is compatible with the unit and multiplication of the monads $\T$ and $\T^{\op}$:
\begin{itemize}
\item for~\eqref{eq:distrlawequationsMult}, let us focus on the topmost diagram, the argument for the bottommost being analogous.
    \begin{equation}\label{eq:distrlawequationsMultTTop}
    \vcenter{
    \xymatrix@C=40pt@R=20pt{
    \T \bicomp{\PJ} \T^{\op} \bicomp{\PJ} \T^{\op} \ar[d]_{\T \mu^{\T^{\op}}} \ar[r]^{\lambdapb \T^{\op}} & \T^{\op}\bicomp{\PJ}\T \bicomp{\PJ}\T^{\op} \ar[r]^{\T^{\op} \lambdapb \T^{\op}} & \T^{\op} \bicomp{\PJ}\T^{\op}\bicomp{\PJ}\T \ar[d]^{\mu^{\T^{\op}} \T}\\
    \T \bicomp{\PJ} \T^{\op} \ar[rr]^{\lambdapb} && \T^{\op}\bicomp{\PJ}\T
    }
    }
    \end{equation}
    Given an element $[\tr{p\in \T}\tr{q\in \T^{\op}}\tr{q'\in \T^{\op}}]_{\equiv_{\PJ}}$ of $\T \poi \T^{\op} \poi \T^{\op}$, the two paths in~\eqref{eq:distrlawequationsMultTTop} are calculated in $\T$ by pullback ($\lambdapb$) and composition ($\mu^{\T^{\op}}$):
    \[
    \text{top-right path:} \quad
    \vcenter{
            \xymatrix@=15pt{
            & & \ar[dl]_{f_2} \pushoutcorner \ar[dr]^{g} & \\
            & \ar[dl]_{f_1} \pushoutcorner \ar[dr] & &\ar[dl]^{q'} \\
            \ar[dr]_{p} & & \ar[dl]^{q} &\\
            & & &
            }
    } \qquad \qquad
    \text{down-left path:} \quad
    \vcenter{
        \xymatrix@=15pt{
        & & \ar[ddll]_{f'} \pushoutcorner \ar[dr]^{g'} & \\
        & & &\ar[dl]^{q'} \\
        \ar[dr]_{p} & &\ar[dl]^{q} &\\
        & & &
        }
    }
    \]
    By universal property of pullback the resulting spans $\tl{f_1 \poi f_2}\tr{g}$ and $\tl{f'}\tr{g'}$ are isomorphic, meaning that they belong to the same $\equiv_{\PJ}$-equivalence class of $\T \bicomp{\PJ} \T^{\op}$: this makes~\eqref{eq:distrlawequationsMultTTop} commute.
\item 
    For~\eqref{eq:distrlawequationsUnit}, let $\eta^{\T^{\op}} \: \PJ \to \T^{\op}$ be the unit of $\T^{\op}$, given by interpreting an arrow of $\PJ$ as one of $\T^{\op}$. We only check commutativity of the upper triangle in~\eqref{eq:distrlawequationsUnit}, the proof for the lower triangle being symmetric. That diagram instantiates to:
     \[
    \xymatrix@R=15pt@C=15pt{
    \T  \ar[d]_-{\cong} & & & \\
     \PJ \bicomp{\PJ} \T \bicomp{\PJ} \PJ \ar[d]_{\PJ \bicomp{\PJ} \T\, \eta^{\T^{\op}}}  \ar@/^1pc/[drr]^{\eta^{\T^{\op}} \T\bicomp{\PJ} \PJ} & & &\\
    \PJ \bicomp{\PJ} \T \bicomp{\PJ} \T^{\op} \ar[r]^{\PJ \lambdapb} &  \PJ \bicomp{\PJ} \T \bicomp{\PJ} \T^{\op} \ar[r]^{\cong} & \T^{\op} \bicomp{\PJ} \T \bicomp{\PJ} \PJ \ar[r]^-{\cong} & \T^{\op} \bicomp{\PJ} \T.
    }
    \]
    The isomorphisms above are given by $\PJ$ being the identity 1-cell on $\PJ$ (see Remark~\ref{rmk:unitlawMod} for more details). We fix an arrow of $\T$ and check commutativity:
    \[
    \xymatrix@R=15pt@C=15pt{
    \tr{f \in \T} \ar@{|->}[d]_-{\cong} & & & \\
      [\tr{\id \in\PJ}\tr{f \in \T}\tr{\id \in \PJ}]_{\equiv_{\PJ}} \ar@{|->}[d]_{{\scriptstyle \PJ \bicomp{\PJ} \T}\, \eta^{\T^{\op}}} \ar@{|->}@/^1pc/[drr]^{\eta^{\T^{\op}} \T\bicomp{\PJ} \PJ} & & &\\
    [\tr{\id \in\PJ}\tr{f \in \T}\tr{\id \in \T^{\op}}]_{\equiv_{\PJ}} \ar@{|->}[r]^{\PJ \lambdapb} & [\tr{\id \in\PJ}\tr{\id \in \T^{\op}}\tr{f \in \T}]_{\equiv_{\PJ}} \ar@{|->}[r]^-{\cong} & [\tr{\id \in \T^{\op}}\tr{f \in \T}\tr{\id \in\PJ}]_{\equiv_{\PJ}} \ar@{|->}[r]^-{\cong} & [\tr{\id \in \T^{\op}}\tr{f \in \T}]_{\equiv_{\PJ}}.
    }
    \]
    In the diagram, we pick $\tl{\id}\tr{f}$ as pullback of $\tr{f}\tl{\id}$ given by $\lambdapb$. This is harmless: any other choice of a pullback cospan would give the same $\equiv_{\PJ}$-equivalence class.
    \end{itemize}
Arrows of the newly defined PROP $\T^{\op} \bicomp{\PJ} \T$ are equivalence classes $[\tr{\in \T^{\op}}\tr{\in \T}]_{\equiv_{\PJ}}$, that means, they are spans $\tl{}\tr{}$ in $\T$ which are identified whenever there is a span isomorphism between them. Composition is defined by $\lambda$, thus is by pullback. Therefore, $\T^{\op} \bicomp{\PJ} \T$ is the classifying category of $\Span{\T}$.
\end{proof}

Proposition~\ref{prop:distrLawPbPo} will find application in the next chapter, where we shall investigate distributive laws defined by pullback and pushout in the PROP of matrices.

\begin{remark} In~\cite{RosebrRWood_fact} distributive laws of categories in the sense of \S\ref{sec:distrlawCats} are discussed and the authors also investigate the possibility of defining distributive laws by pullback and pushout. With this aim, they propose to relax the definition of distributive law so that diagrams~\eqref{eq:distrlawequationsUnit}-\eqref{eq:distrlawequationsMult} are required to commute only up-to an arrow of a fixed groupoid $\PJ$ (in our case, $\PJ$ is the core of $\T$): this yields a bicategory as the result of composition, which can be turned into a category by quotienting hom-sets by equivalence up-to $\PJ$.

This construction does not immediately generalise to the case of PROPs: differently from categories, distributive laws of PROPs need to be well-defined as mappings between $\equiv_{\Perm}$-equivalence classes, which as explained above is not guaranteed for the case of pullback and pushout. Our approach handles this additional challenge and also does not require any tweak of the definition of distributive law in a bicategory. 
\end{remark}

\subsection{Operations on Distributive laws: Composition, Quotient and Dual}\label{sec:iteratedDistrLaws}

Defining explicitly a distributive law and proving that it is one can be rather challenging. It is thus useful to provide some basic techniques to canonically form new distributive laws from existing ones. The operations that we present are:
\begin{itemize}[itemsep=1pt,topsep=1pt,parsep=1pt,partopsep=1pt]
\item the composition of distributive laws;
\item the quotient of a distributive law by a set of equations;
\item given $\PS \bicomp{\Perm} \T$, the construction of a distributive laws yielding the composite PROP $\T^{\op} \bicomp{\Perm} \PS^{\op}$.
\end{itemize}

\paragraph{Composing distributive laws} For our developments it is useful to generalise PROP composition to the case when there are more than two theories interacting with each other.
The following result, which is easily provable by diagram chasing, is part of a thorough investigation of iterated distributive laws pursued in~\cite{Cheng_IteratedLaws}.

\begin{proposition}[\cite{Cheng_IteratedLaws}] \label{prop:iterateddistrlaw}
Let $\funF$, $\G$, $\FH$ be monads in the same bicategory and suppose there are distributive laws
\begin{gather*}
\lambda \: \FH \poi \funF \to \funF \poi \FH \qquad\qquad
\chi \: \FH \poi \G \to \G \poi \FH \qquad \qquad
\psi \: \G \poi \funF \to \funF \poi \G
\end{gather*}
satisfying the following ``Yang-Baxter'' equation:
\begin{equation}\label{eq:yangbaxter}
\vcenter{
\xymatrix@R=10pt{
& \FH \poi\funF \poi \G  \ar[r]^{\lambda \G} &\funF \poi \FH \poi \G \ar[dr]^{\funF \chi} &\\
\FH \poi \G \poi \funF \ar[ur]^{\FH \psi} \ar[dr]_{\chi \funF} & & &\funF \poi \G \poi \FH  \\
& \G \poi \FH \poi \funF\ar[r]^{\G \lambda} & \G \poi \funF \poi \FH \ar[ur]_{\psi \FH} &\\
}
}
\end{equation}
then
\begin{gather}
\alpha \: \FH \poi \funF \poi \G\xrightarrow{\lambda \G} \funF \poi \FH \poi \G\xrightarrow{\funF \chi} \funF \poi \G \poi \FH \label{eq:threelaws1}\\
\beta \: \G \poi \FH \poi \funF\xrightarrow{\G \lambda} \G \poi \funF\poi \FH\xrightarrow{\psi \FH} \funF \poi \G \poi \FH \label{eq:threelaws2}
\end{gather}
are distributive laws yielding the same monad structure on $\funF \poi \G \poi \FH$.
\end{proposition}
Note that $\alpha$ distributes $\FH$ over $\funF \poi \G$ and $\beta$ distributes $\G \poi \FH$ over $\funF$. Yielding ``the same'' monad structure means that the unit and multiplication of $(\funF \poi \G) \poi \FH$ and of $\funF \poi (\G \poi \FH)$, defined according to $\alpha$ and $\beta$ respectively, are the same 2-cell. 

Supposing that $\funF$, $\G$ and $\FH$ are PROPs, what are the equations arising by their composition? It turns out that $\funF \poi \G \poi \FH$ can be presented as the quotient of $\funF + \G + \FH$ by the equations in $\chi$, $\psi$ and $\lambda$.

\begin{proposition}\label{prop:SMTforIteratedComp} Let $\funF$, $\FH$ and $\G$ be PROPs presented by SMTs $(\Sigma_{\funF}, E_{\funF})$, $(\Sigma_{\FH}, E_{\FH})$ and $(\Sigma_{\G}, E_{\G})$ respectively. Suppose there are distributive laws $\lambda$, $\chi$ and $\psi$ yielding composed PROPs $\funF \bicomp{\Perm} \FH$, $\G \bicomp{\Perm} \FH$ and $\funF \bicomp{\Perm} \G$ respectively and satisfying \eqref{eq:yangbaxter}. Call $E_{\lambda}$, $E_{\chi}$ and $E_{\psi}$ the sets of equations encoding the three laws. Then there exist a composed PROP $\funF \bicomp{\Perm} \G \bicomp{\Perm} \FH$ presented by the SMT with signature $\Sigma_{\funF}\uplus \Sigma_{\FH}\uplus \Sigma_{\G}$ and equations $E_{\funF}\uplus E_{\FH}\uplus E_{\G}\uplus E_{\lambda}\uplus E_{\chi}\uplus E_{\psi}$.
\end{proposition}
\begin{proof}
Consider first the PROP $(\funF \bicomp{\Perm} \G) \bicomp{\Perm} \FH$, defined by the distributive law $\alpha$. Proposition~\ref{prop:SMTforComposition} gives us the following recipe for the SMT presenting $(\funF \bicomp{\Perm} \G) \bicomp{\Perm} \FH$:
\begin{itemize}
\item the signature consits of the signature $\Sigma_{\funF}\uplus\Sigma_{\G}$ of $\funF \bicomp{\Perm} \G$ and the signature $\Sigma_{\FH}$ of $\FH$;
\item the set of equations consists of the equations $E_{\funF}\uplus E_{\G}\uplus E_{\psi}$ of $\funF \bicomp{\Perm} \G$, the equations $E_{\FH}$ of $\FH$ and the equations $E_{\alpha}$ encoded by $\alpha$. By definition of $\alpha$ according to~\eqref{eq:threelaws1}, $E_{\alpha} = E_{\lambda}\uplus E_{\chi}$.
\end{itemize}
We can thereby conclude that $(\funF \bicomp{\Perm} \G) \bicomp{\Perm} \FH$ is presented by the SMT described in the statement. One can reach the same conclusion starting from $\funF \bicomp{\Perm} (\G \bicomp{\Perm} \FH)$: in this case, $E_{\chi}$ will be already part of the SMT of $\G \bicomp{\Perm} \FH$ and the sets $E_{\lambda}$ and $E_{\psi}$ are added when composing $\funF$ and $\G \bicomp{\Perm} \FH$ using~$\beta$.
\end{proof}

\begin{example}[\textbf{Partial functions}] \label{ex:partialfunctions} In Example~\ref{ex:composingSemPROP}\ref{ex:composingSemPROP1} we gave a modular construction for the PROP $\F$ of functions. We now show how the PROP $\gls{PROPPF}$ of \emph{partial} function can be presented modularly using iterated distributive laws.



The leading intuition is to use $\Wmult$, $\Wunit$ and $\Ocounit$ as generators for string diagrams to represent the graph of a partial function. The interpretation extends the one given when explaining the isomorphism $\Mon \cong \F$ at the end of \S~\ref{sec:props}: $\Wmult$ will mean that two elements in the domain are mapped into the same element of the codomain, and $\Wunit$ will mean that a certain element of the codomain is not in the image of the function. The generator $\Ocounit$, which was not part of the theory $\Mon$, indicates partiality: a certain element is not in the domain of the function.

To make this formal, we shall form the composite of three PROPs: the theory $\MULT$ of $\Wmult$ and the theory $\UNIT$ of $\Wunit$, introduced in Example~\ref{ex:distrlawsyntactic}\ref{ex:distrlawsyntactic1}, and the newly introduced theory $\gls{PROPCounit}$. This is the PROP generated by the signature $\{\Ocounit\}$ and no equations: modulo the orange/white colouring, it is just $\UNIT^{\op}$. Following the recipe of Proposition~\ref{prop:iterateddistrlaw}, we now combine these PROPs together via three distributive laws:
\begin{gather*}
\lambda \: \UNIT \bicomp{\Perm} \COUNIT \to \COUNIT \bicomp{\Perm} \UNIT \qquad\qquad
\chi \: \UNIT \bicomp{\Perm} \MULT \to \MULT \bicomp{\Perm} \UNIT \qquad \qquad
\psi \: \MULT \bicomp{\Perm} \COUNIT \to \COUNIT \bicomp{\Perm} \MULT
\end{gather*}
For defining them, we can use the concrete characterisations of $\UNIT$ as $\Inj$, of $\COUNIT$ as $\Inj^{\op}$ and of $\MULT$ as $\Surj$ (Example~\ref{ex:distrlawsyntactic}\ref{ex:distrlawsyntactic1}). We define the distributive law $\chi$ by epi-mono factorisation as in Example~\ref{ex:distrlawsyntactic}\ref{ex:distrlawsyntactic1}; therefore, the resulting PROP $\MULT \bicomp{\Perm} \UNIT$ is $\MULT + \UNIT$ quotiented by \eqref{eq:wmonunitlaw}. Because pullbacks preserve monos, we can define $\lambda$ by pullback in $\F$. Given the extremely simple signature of $\UNIT$ and $\COUNIT$, it is immediate to check that the only equation needed to present $E_{\lambda}$ is
\begin{equation}
\label{eq:wbbone}\tag{P1}
\lower6pt\hbox{$\includegraphics[height=.6cm]{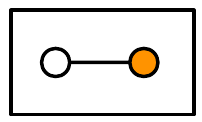}$}
=
\lower5pt\hbox{$\includegraphics[width=16pt]{graffles/idzerocircuit.pdf}$}.
\end{equation}
 Therefore, the composite PROP $\COUNIT \bicomp{\Perm} \UNIT$ is the quotient of $\COUNIT + \UNIT$ by~\eqref{eq:wbbone}. The distributive law $\psi$ can also be defined by pullback, because pullbacks in $\F$ preserve epis (surjections). Also in this case, it is easy to verify that $E_{\psi}$ amounts to the equation
\begin{equation}
\label{eq:wbbialgunit}\tag{P2}
\lower8pt\hbox{$\includegraphics[height=.7cm]{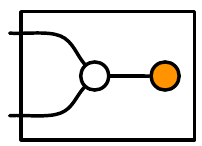}$}
\! = \!
\lower8pt\hbox{$\includegraphics[height=.7cm]{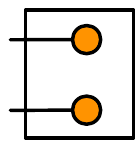}$}.
\end{equation}
The resulting PROP $\MULT \bicomp{\Perm}\COUNIT$ is the quotient of $\MULT+\COUNIT$ by~\eqref{eq:wbbialgunit}. By case analysis on a string diagram in $\UNIT \bicomp{\Perm} \MULT \bicomp{\Perm} \COUNIT$, one can check that $\lambda$, $\chi$ and $\psi$ verify the Yang-Baxter equation~\eqref{eq:yangbaxter}. Thus by Proposition~\ref{prop:iterateddistrlaw} there are distributive laws $\alpha$, $\beta$ yielding the same PROP $\COUNIT \bicomp{\Perm} \MULT \bicomp{\Perm} \UNIT$:
\begin{eqnarray*}
\alpha & \: & \UNIT \bicomp{\Perm}\COUNIT \bicomp{\Perm} \MULT \xrightarrow{\lambda\MULT} \COUNIT\bicomp{\Perm} \UNIT \bicomp{\Perm} \MULT \xrightarrow{\COUNIT \chi} \COUNIT \bicomp{\Perm} \MULT \bicomp{\Perm} \UNIT \\
\beta & \: & \MULT \bicomp{\Perm} \UNIT \bicomp{\Perm} \COUNIT \xrightarrow{\MULT \lambda} \MULT \bicomp{\Perm} \COUNIT \bicomp{\Perm} \UNIT \xrightarrow{\psi \UNIT} \COUNIT \bicomp{\Perm} \MULT \bicomp{\Perm} \UNIT.
\end{eqnarray*}
 It is instructive to give a direct description of $\alpha$: it maps a given cospan $\tr{f \in \Surj}\tr{g \in \Inj}\tl{h \in \Inj}$ to the span $\tl{p \in \Inj}\tr{q \in \Surj}\tr{o \in \Inj}$ obtained by first taking the pullback $\tl{p' \in \Inj}\tr{o \in \Inj}$ of $\tr{g \in \Inj}\tl{h \in \Inj}$ as prescribed by $\lambda$ and then pull back $\tr{f \in \Surj}\tl{p' \in \Inj}$ as an application of $\chi$.
\[
\xymatrix@=13pt{
& \ar[dl]_{p} \pushoutcorner \ar[dr]^{q} & & \\
\ar[dr]_{f} & & \ar[dl]_{p'} \pushoutcorner \ar[dr]^{o} & \\
& \ar[dr]_{g} & & \ar[dl]^{h}\\
& & &
}
\]
By Proposition~\ref{prop:SMTforIteratedComp}, the equations presenting $\COUNIT \bicomp{\Perm} \MULT \bicomp{\Perm} \UNIT$ are those generated by $\lambda$, $\psi$ and $\chi$. Therefore, $\COUNIT \bicomp{\Perm} \MULT \bicomp{\Perm} \UNIT$ is the quotient of $\COUNIT + \MULT + \UNIT$ by \eqref{eq:wmonunitlaw}, \eqref{eq:wbbone} and \eqref{eq:wbbialgunit}. 

We now claim that $\COUNIT \bicomp{\Perm} \MULT \bicomp{\Perm} \UNIT \cong \PF$, that is, the equational theory that we just described presents the PROP of partial functions. To see this, observe that partial functions $n \tr{f \in \PF} m$ are in bijective correspondence with spans $n \tl{i \in \Inj} z \tr{f \in \F}m$: the injection $i$ tells on which elements $\ord{z}$ of $\ord{n}$ the function $f$ is defined. Since $\Inj^{\op} \cong \COUNIT$ and $\F \cong \Mon \cong \MULT \bicomp{\Perm} \UNIT$, this correspondence yields the desired isomorphism $\PF \cong \Inj^{\op} \bicomp{\Perm} \F \cong \COUNIT \bicomp{\Perm} \MULT \bicomp{\Perm} \UNIT$. 

As a concluding remark, let us note that the factorisation property of $\PMon$ allows to interpret any arrow of this PROP as the graph of a partial function, following the intuition given at the beginning of this example. For instance,
$$\includegraphics[height=1cm]{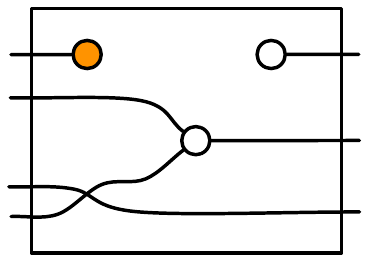}$$
represents the function $\ord{4} \to \ord{3}$ undefined on $1$ and mapping $2,4$ to $2$ and $3$ to $3$. This establishes a bijective correspondence with $\PF$, analogously to the one between $\Mon$ and $\F$.
\end{example}

\paragraph{Quotients of distributive laws} Definition~\ref{def:monadquotient} introduced the notion of \emph{quotient} $\Mquot \: \funF \to \G$ of a monad $\funF$: the idea is that the monad $\G$ is obtained by imposing additional equations on the algebraic theory described by $\funF$. As one may expect, distributive laws are compatible with monad quotients, provided that the law preserves the newly added equations. This folklore result appears in various forms in the literature: \cite{BHKR13} gives it for distributive laws of endofunctors over monads and \cite{BMSZ-daggerJournal,BonchiZanasiLPJournal} for distributive laws of monads. All these references concern distributive laws in $\Cat$. For our purposes, it is useful to state the result for arbitary bicategories.

\begin{proposition}\label{prop:quotientDistrLaw} Suppose that $\lambda \: \funF \poi \FH \to \FH \poi \funF$ is a distributive law in a bicategory $\bicat$, $\Mquot \: \funF \to \G$ a monad quotient and $\lambda' \: \G \poi \FH \to \FH \poi \G$ another 2-cell of $\bicat$ making the following diagram commute.
\begin{equation}\label{eq:quotientDistrLaw}
\vcenter{
\xymatrix{
    \funF \poi \FH \ar[r]^{\Mquot \FH} \ar[d]_{\lambda} & \G \poi \FH \ar[d]^{\lambda'} \\
    \FH \poi \funF \ar[r]^{\FH \Mquot} & \FH \poi \G
}
}
\end{equation}
Then $\lambda'$ is a distributive law of monads.
\end{proposition}
\begin{proof} The diagrams for compatibility of $\lambda'$ with unit and multiplication of $\G$ commute because $\Mquot$ is a monad morphism and \eqref{eq:quotientDistrLaw} commutes. For compatibility of $\lambda'$ with unit and multiplication of $\FH$, one needs to use commutativity of \eqref{eq:quotientDistrLaw} and the fact that $\Mquot$ is epi.
\end{proof}

We remark that Proposition~\ref{prop:quotientDistrLaw} holds also in the version in which one quotients the monad $\FH$ instead of $\funF$. It is now useful to instantiate the result to the case of distributive laws of PROPs.

\begin{proposition}\label{prop:quotientDistrLawsPROP} Let $\T$ be the PROP freely generated by $(\Sigma,E)$ and $E' \supseteq E$ be another set of equations on $\Sigma$-terms. Suppose that there exists a distributive law $\lambda \: \T \bicomp{\Perm} \PS \to \PS \bicomp{\Perm} \T$ such that, if $E'$ implies $c = d$, then $\lambda(\tr{c \in \T}\tr{e \in \PS}) = \lambda(\tr{d \in \T}\tr{e \in \PS})$. Then there exists a distributive law $\lambda' \: \T' \bicomp{\Perm} \PS \to \PS \bicomp{\Perm} \T'$ presented by the same equations as $\lambda$, i.e., $E_{\lambda'} = E_{\lambda}$.
\end{proposition}
\begin{proof}
There is a PROP morphism $\Mquot \: \T \to \T'$ defined by quotienting string diagrams in $\T$ by $E'$. This is a monad quotient in the bicategory $\Mod{\Span{\MonC}}$ where PROPs are monads. We now define another 2-cell $\lambda' \: \T' \bicomp{\Perm} \Com \to \Com\bicomp{\Perm}\T'$ as follows: given $\tr{e \in \T'} \tr{c \in \Com}$, pick any $\tr{d \in \T}$ such that $\Mquot(d) = e$ and let $\tr{c' \in \Com}\tr{d' \in \T}$ be $\lambda(\tr{d \in \T}\tr{c \in \Com})$. Define $\lambda'(\tr{e \in \T'} \tr{c \in \Com})$ as $\tr{c' \in \Com}\tr{\Mquot(d') \in \T'}$. $\lambda'$ is well-defined because, by assumption, if $\Mquot(d_1) = \Mquot(d_2)$ then $E'$ implies that $d_1 = d_2$ and thus $\lambda(\tr{d_1}\tr{c}) = \lambda(\tr{d_2}\tr{c})$.

Now, $\lambda$, $\lambda'$ and $\Mquot$ satisfy the assumptions of Proposition~\ref{prop:quotientDistrLaw}. In particular, \eqref{eq:quotientDistrLaw} commutes by definition of $\lambda'$ in terms of $\lambda$ and $\Mquot$. The conclusion of Proposition~\ref{prop:quotientDistrLaw} guarantees that $\lambda'$ is a distributive law. By construction, $\lambda'$ is presented by the same equations as $\lambda$.
\end{proof}

We will see an application of Proposition~\ref{prop:quotientDistrLawsPROP} in the next section (Lemma~\ref{lemma:LawNoEquations}).

\paragraph{Distributive laws on the opposite PROPs} Let $\PS$ and $\T$ be PROPs. For later reference, we observe that a distributive law $\lambda \: \T \bicomp{\Perm} \PS \to \PS\bicomp{\Perm} \T$ canonically induces one $\lambda' \: \PS^{\op} \bicomp{\Perm} \T^{\op}  \to \T^{\op} \bicomp{\Perm} \PS^{\op}$ defined as:
\[
\xymatrix{
[\tr{\in \PS^{\op}}\tr{\in \T^{\op}}]_{\equiv_{\Perm}} \ar@{|->}@/^3pc/[rrr]^{\lambda'}  \ar@{}[r]|{=} & [\tl{\in \PS}\tl{\in \T}]_{\equiv_{\Perm}} \quad \ar@{|->}[r]^{\lambda} & \quad [\tl{\in \T}\tl{\in \PS}]_{\equiv_{\Perm}} \ar@{}[r]|{=} & [\tr{\in \T^{\op}}\tr{\in \PS^{\op}}]_{\equiv_{\Perm}}
}
\]

\begin{proposition}\label{prop:distrlawOP} $\lambda' \: \PS^{\op} \bicomp{\Perm} \T^{\op}  \to \T^{\op} \bicomp{\Perm} \PS^{\op}$ is a distributive law of PROPs.
\end{proposition}
\begin{proof} The main observation is that the unit and multiplication of $\PS^{\op}$ and of $\T^{\op}$ can be expressed in terms of those of $\PS$ and $\T$ in the same way as $\lambda'$ is defined from $\lambda$. Then,
using the fact that $\lambda$ is a distributive law, it is immediate to check that $\lambda'$ makes diagrams \eqref{eq:distrlawequationsUnit}-\eqref{eq:distrlawequationsMult} commute. \end{proof}

\subsection{Lawvere Theories as Composed PROPs}\label{sec:lawvere}
This section studies a class of distributive laws whose equations present Lawvere theories. First, recall that \emph{Lawvere theories}~\cite{LawvereOriginalPaper,hyland2007category} are a special kind of categories adapted to the study of universal algebra. They are closely related to PROPs: the essential difference is that, wheareas a Lawvere theory is required to be a category with finite products, a PROP may carry any symmetric monoidal structure, not necessarily cartesian.

Just as PROPs, Lawvere theories can be also freely obtained by generators and equations. By analogy with symmetric monoidal theories introduced in \S~\ref{sec:props}, we organise these data as a \emph{cartesian theory}: it simply amounts to the notion of equational theory that one typically finds in abstract algebra, see e.g.~\cite{BurrisSankappanavar1981}.
\begin{definition}\label{def:Lawvere} A \emph{cartesian theory} $(\Sigma,E)$ consists of a signature $\Sigma=\{o_1 \: n_1 \to 1,\dots,o_k \: n_k \to 1\}$ and a set $E$ of equations between \emph{cartesian $\Sigma$-terms}, which are defined as follows:
\begin{itemize}[noitemsep,topsep=0pt,parsep=0pt,partopsep=0pt]
\item for each $i \in \N$, the \emph{variable} $x_i$ is a cartesian term;
\item suppose $o \: n \to 1$ is a generator in $\Sigma$ and $t_1 , \dots , t_n$ are cartesian terms. Then $o(t_1,\dots,t_n)$ is a cartesian term.
\end{itemize}

The \emph{Lawvere theory} $\gls{LawSigmaE}$ freely generated by $(\Sigma, E)$ is the category whose objects are the natural numbers and arrows $n \to m$ are lists $\tpl{t_1,\dots,t_m}$ of cartesian $\Sigma$-terms quotiented by $E$, such that, for each $t_i$, only variables among $x_1,\dots,x_n$ appear in $t_i$. Composition is by substitution:
\[ \left(n\tr{\tpl{t_1,\dots,t_m}}m\right) \poi \left(m \tr{\tpl{s_1,\dots,s_z}}z\right) = n \tr{\tpl{s_1[t_i / x_i \mid 1 \leq i \leq m],\dots,s_z[t_i / x_i \mid 1 \leq i \leq m]}} z\]
where $t[t' / x]$ denotes the cartesian term $t$ in which all occurrences of the variable $x$ have been substituted with $t'$.

$\LwA{\Sigma,E}$ is equipped with a product $\times$ which is defined on objects by addition and on arrows by list concatenation and suitable renaming of variables:
\[ \left(n\tr{\tpl{t_1,\dots,t_m}}m\right) \times  \left(z \tr{\tpl{s_1,\dots,s_l}} l\right) = n+z \tr{\tpl{t_1,\dots,t_m,s_1[x_{i+m} / x_i \mid 1 \leq i \leq l],\dots,s_l[x_{i+m} / x_i \mid 1 \leq i \leq l]}} m+l.\]
\end{definition}
We use notation $\var{t}$ for the list of occurrences of variables appearing (from left to right) in $t$ and, more generally, $\var{t_1,\dots,t_m}$ for the list $\var{t_1}::\dots::\var{t_m}$. Also, $\size{l} \in \N$ denotes the length of a list $l$. We say that a list $\tpl{t_1,\dots,t_m} \: n \to m$ is \emph{linear} if each variable among $x_1,\dots,x_n$ appears exactly once in $\var{t_1,\dots,t_m}$. 

Our first observation is that Lawvere theories are a particular kind of PROP.
\begin{proposition}\label{proposition:LawvereArePROPs} $\LwA{\Sigma,E}$ is a PROP. \end{proposition}
\begin{proof} Let $\times$ act as the monoidal product, $0$ as its unit and define the symmetry $n + m \to m+n$ as the list $\tpl{x_{n+1},\dots,x_{n+m},x_1,\dots,x_n}$. It can be verified that $\LwA{\Sigma,E}$ equipped with this structure satisfies the laws of symmetric strict monoidal categories, thus it is a PROP.\end{proof}

As a side observation, note that the unique PROP morphism $\Perm \to \LwA{\Sigma,E}$ sends $p \: n \to n$ to $\tpl{x_{p^{-1}(1)},\dots,x_{p^{-1}(n)}}$.

\begin{remark}\label{rmk:Law_sigmaToCartesianTerms} In spite of Proposition~\ref{proposition:LawvereArePROPs}, cartesian theories are \emph{not} a subclass of symmetric monoidal theories: in fact, the two concepts are orthogonal. On the one hand, a symmetric monoidal theory $(\Sigma, E)$ is cartesian if  and only if all generators in $\Sigma$ have coarity $1$ and, for all equations $t=s$ in $E$, $t$ and $s$ are $\Sigma$-terms with coarity $1$. Under these conditions, there is a canonical way to interpret any $\Sigma$-term $n \to m$ as a list of $m$ cartesian $\Sigma$-terms on variables $x_1,\dots,x_n$. Below, $o$ ranges over $\Sigma$:
\[\Idnet \: 1 \to 1 \ \mapsto \ \tpl{x_1} \qquad  \symNet \: 2 \to 2 \mapsto \tpl{x_2,x_1} \qquad \OpDiag \: n \to 1 \ \mapsto \  \tpl{o(x_1,\dots,x_n)}\]
The inductive cases are defined using the operations $\poi$ and $\tns$ on lists given in Definition~\ref{def:Lawvere}. Note that $\Sigma$-terms always denote \emph{linear} lists of cartesian terms. This explains why, conversely, not all the cartesian theories are symmetric monoidal: their equations possibly involve non-linear cartesian $\Sigma$-terms, which are not expressible with (symmetric monoidal) $\Sigma$-terms. The subtlety here is that, in a sense, we can still \emph{simulate} a cartesian theory on signature $\Sigma$ with a symmetric monoidal theory, which however will be based on a larger signature $\Sigma'$, recovering the possibility of copying and discarding variables by the use of additional generators. This point will become more clear below, where we will see how copier and discharger, i.e., the cartesian structure, can be mimicked with the use of the PROP $\Com$.
\end{remark}

\begin{example} The SMT $(\Sigma_{\scriptscriptstyle M},E_{\scriptscriptstyle M})$ of commutative monoids is cartesian. It generates the Lawvere theory $\LwA{\Sigma_{\scriptscriptstyle M},E_{\scriptscriptstyle M}}$ whose arrows $n \to m$ are lists $\tpl{t_1,\dots,t_m}$ of elements of the free commutative monoid on $\{x_1,\dots,x_n\}$.
\end{example}

The example of commutative monoids is particularly convenient to sketch our approach to Lawvere theories as composed PROPs. First, note that the Lawvere theory $\LwA{\Sigma_{\scriptscriptstyle M},E_{\scriptscriptstyle M}}$ \emph{includes} the PROP $\Mon$ freely generated by $(\Sigma_{\scriptscriptstyle M},E_{\scriptscriptstyle M})$. Indeed, any string diagram of $\Mon$ can be interpreted as a list of terms following the recipe of Remark~\ref{rmk:Law_sigmaToCartesianTerms}. For instance,
\[ \lower12pt\hbox{$\includegraphics[height=1.2cm]{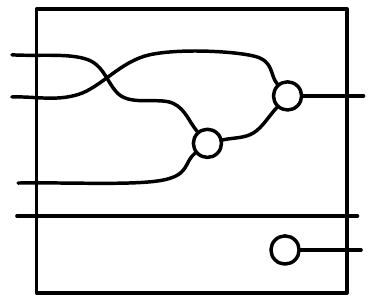}$}\: 4 \to 3 \text{ is interpreted as } \tpl{\Wmult(x_2,\Wmult(x_1,x_3)), x_4, \Wunit} \: 4 \to 3 \]
As we observed above, string diagrams of $\Mon$ can only express \emph{linear} lists.  What makes $\LwA{\Sigma_{\scriptscriptstyle M},E_{\scriptscriptstyle M}}$ more general than $\Mon$ is the ability of \emph{copying} and \emph{discarding} variables, formally explained by the fact that the monoidal product $\LwA{\Sigma_{\scriptscriptstyle M},E_{\scriptscriptstyle M}}$ is the cartesian product~\cite{eilenberg1966closed}. These operations are expressed by arrows
\[ \tpl{x_1,x_1}  \: 1 \to 2 \qquad \qquad \text{ and } \qquad \qquad \tpl{\,} \: 1 \to 0. \]
How can we recover copy and discard in the language of string diagrams for the PROP $\LwA{\Sigma_{\scriptscriptstyle M},E_{\scriptscriptstyle M}}$? A basic observation is that these two operations satisfy the equations of commutative comonoids: discarding is commutative and associative, and copying and then discarding is the same as not doing anything. It thus makes sense to interpreted them as the generators of the PROP $\Com$:
\[ \Bcomult \: 1 \to 2 \qquad \qquad \qquad \qquad \qquad \qquad \Bcounit \: 1 \to 0.\]
Our approach suggests that a copy of $\Mon$ and of $\Com$ ``live'' inside $\LwA{\Sigma_{\scriptscriptstyle M},E_{\scriptscriptstyle M}}$. Also, we claim that these two PROPs provide a complete description of $\LwA{\Sigma_{\scriptscriptstyle M},E_{\scriptscriptstyle M}}$, that means, any arrow of $\LwA{\Sigma_{\scriptscriptstyle M},E_{\scriptscriptstyle M}}$ can be presented diagrammatically by using $\Mon$ and $\Com$. For instance,
\[ \tpl{\Wmult(x_2,\Wmult(x_1,x_4)), x_1, \Wunit} \: 4 \to 3 \text{ corresponds to }\lower16pt\hbox{$\includegraphics[height=1.6cm]{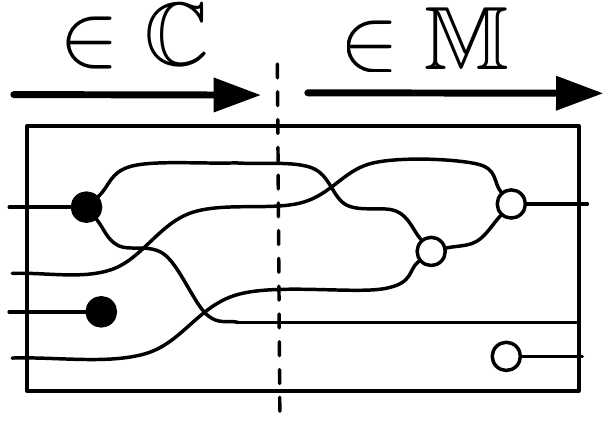}$}\: 4 \to 3  \]
Observe that the diagram is of the factorised form $\tr{\in \Com}\tr{\in \Mon}$. Intuitively, $\Com$ is deputed to model the interplay of variables --- in this case, the fact that $x_1$ is copied and $x_3$ is deleted --- and $\Mon$ describes the syntactic tree of the terms. Of course, to claim that this factorisation is always possible, we need additional equations to model composition of factorised diagrams. For instance:
\begin{eqnarray*}
\tpl{\Wmult(x_1,x_2),x_1} \poi \tpl{x_1,\Wmult(x_1,x_2)} & = &  \tpl{\Wmult(x_1,x_2),\Wmult(\Wmult(x_1,x_2),x_1)}. \\ \\
\qquad\lower8pt\hbox{$\includegraphics[height=.8cm]{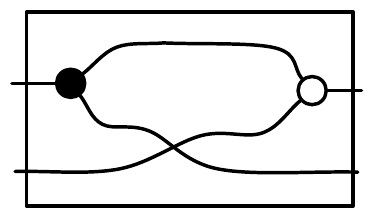}$} \qquad \poi \qquad \lower8pt\hbox{$\includegraphics[height=.8cm]{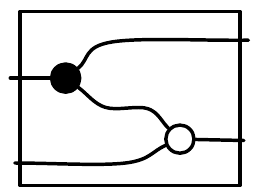}$}\quad \qquad & \eql{?} & \qquad \qquad\qquad
\lower10pt\hbox{$\includegraphics[height=1cm]{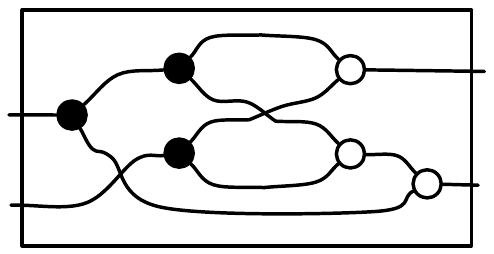}$}
\end{eqnarray*}
The second equality holds if we assume the equation~\eqref{eq:bialg} of the SMT of \emph{bialgebras}.
Thus the example suggests that composition by subsitution in $\LwA{\Sigma_{\scriptscriptstyle M},E_{\scriptscriptstyle M}}$ can be mimicked at the diagrammatic level by allowing the use of bialgebra equations, which as we know from Example~\ref{ex:distrlawsyntactic}\ref{ex:distrlawsyntactic1} present the composite PROP $\Com \bicomp{\Perm} \Mon$. Therefore, the conclusive conjecture of our analysis is that $\LwA{\Sigma_{\scriptscriptstyle M},E_{\scriptscriptstyle M}}$ must be isomorphic to $\Com \bicomp{\Perm} \Mon$ and can be presented by equations~\eqref{eq:wmonassoc}-\eqref{eq:bwbone}.

We now formally develop the above approach. The following is the main result of this section.
\begin{theorem} \label{Th:LawvereCompositePROP} Suppose that $(\Sigma,E)$ is an SMT which is also cartesian and let $\T$ be its freely generated PROP. Then $\LwA{\Sigma,E}$ is the composite PROP $\Com \bicomp{\Perm} \T$. The distributive law $\T \bicomp{\Perm} \Com \to \Com \bicomp{\Perm} \T$ yielding $\LwA{\Sigma,E}$ is presented by equations

\medskip

\noindent\begin{multicols}{2}\noindent
{%
\begin{equation}\label{eq:LawDistrBcounit}
\tag{Lw1}
\lower10pt\hbox{$\includegraphics[height=1cm]{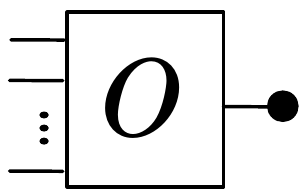}$} = \lower10pt\hbox{$\includegraphics[height=1cm]{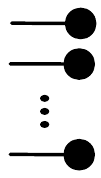}$}
\end{equation}
\vspace{-.8cm}}%
{%
\begin{equation}\label{eq:LawDistrBcomult} \tag{Lw2}
\lower10pt\hbox{$\includegraphics[height=1cm]{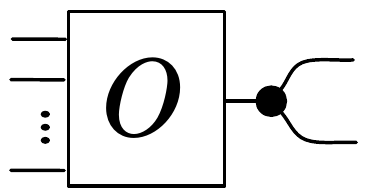}$} = \lower20pt\hbox{$\includegraphics[height=1.7cm]{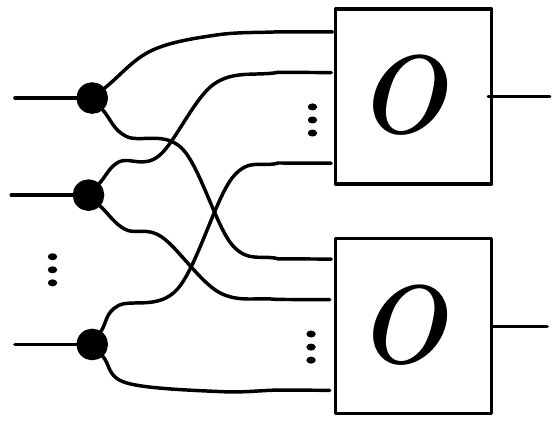}$}
\end{equation}
}%
\end{multicols}
for each $o \in \Sigma$.
\end{theorem}
Before moving to the proof of Theorem~\ref{Th:LawvereCompositePROP}, let us explain its significance, which we believe is two-fold. First, it gives a deeper understanding of the nature of Lawvere theories and how they formally relate to PROPs. In fact, the role played by $\Com$ and equations~\eqref{eq:LawDistrBcounit}-\eqref{eq:LawDistrBcomult} is known in the literature, see~\cite{Burroni1993,Lafont95-equationalReasoningTwoDimDiagrams,GadducciCorradini-termgraph}: what we believe is original of our approach is the formulation in terms of distributive laws, which reveals the provenance of \eqref{eq:LawDistrBcounit}-\eqref{eq:LawDistrBcomult}.

Second, Theorem~\ref{Th:LawvereCompositePROP} provides a technique to define distributive laws for a wide range of PROPs, showing that the result of composition is a well-known categorical notion and giving a finite axiomatisation for it. For instance, for $\T = \Mon$, \eqref{eq:LawDistrBcounit}-\eqref{eq:LawDistrBcomult} become the bialgebra equations~\eqref{eq:unitsl}-\eqref{eq:bwbone} and the result that $\B \cong \Com \bicomp{\Perm}\Mon$ is now an immediate consequence of Theorem~\ref{Th:LawvereCompositePROP}. Our characterisation will also be useful for the developments of the next chapter.

\medskip

The rest of the section is devoted to proving Theorem~\ref{Th:LawvereCompositePROP}. First we observe that, by the following lemma, it actually suffices to check our statement for SMTs with no equations. This reduction has just the purpose of making computations in $\Lw{(\Sigma,E)}$ easier, by working with terms instead of equivalence classes.

\begin{lemma}\label{lemma:LawNoEquations} If Theorem~\ref{Th:LawvereCompositePROP} holds for $E = \emptyset$, then it holds for any cartesian SMT $(\Sigma, E)$.\end{lemma}
\begin{proof} Let $(\Sigma,E)$ be a cartesian SMT and $\T$, $\T_{\scriptscriptstyle E}$ be the PROPs freely generated, respectively, by $(\Sigma,\emptyset)$ and $(\Sigma,E)$. By assumption, Theorem~\ref{Th:LawvereCompositePROP} holds for $(\Sigma,\emptyset)$, giving us a distributive law $\lambda \: \T \bicomp{\Perm} \Com \to \Com \bicomp{\Perm} \T$. The existence of a distributive law $\lambda' \: \T_{\scriptscriptstyle E} \bicomp{\Perm} \Com \to \Com\bicomp{\Perm}\T_{\scriptscriptstyle E}$ with the required properties is guaranteed by Proposition~\ref{prop:quotientDistrLawsPROP}, provided that $\lambda$ preserves the equations of $E$ --- we refer to Appendix~\ref{app:Lawvere} for the proof of this statement.
\end{proof}


 In the sequel, let us abbreviate $\LwA{\Sigma,\emptyset}$ as $\LwA{\Sigma}$. By virtue of Lemma~\ref{lemma:LawNoEquations}, we shall prove Theorem~\ref{Th:LawvereCompositePROP} for $\LwA{\Sigma}$ and by letting $\T$ be the PROP freely generated by $(\Sigma,\emptyset)$. Our strategy to construct the distributive law $\T \bicomp{\Perm} \Com \to \Com \bicomp{\Perm} \T$ is to first interpret string diagrams $\tr{ \in \T}\tr{\in \Com}$ as arrows of $\LwA{\Sigma}$ and then show that any arrow of $\LwA{\Sigma}$ admits a decomposition as $\tr{\in \Com}\tr{\in \T}$. In order to claim that this approach indeed yields a distributive law, we use the following result from~\cite{Lack2004a}.

\begin{proposition}[\cite{Lack2004a}] \label{prop:fact->distrLaw} Let $\PS$ be a PROP and $\T_1$, $\T_2$ be sub-PROPs of $\PS$. Suppose that each arrow $n \tr{f \in \PS} m$ can be factorised as $n \tr{g_1 \in \T_1}\tr{g_2 \in \T_2} m$ uniquely up-to permutation, that is, for any other decomposition $n \tr{h_1 \in \T_1}\tr{h_2 \in \T_2}m$ of $f$, there exists a permutation $\tr{p \in \Perm}$ such that the following diagram commutes.
\[\vcenter{
            \xymatrix@R=13pt@C=23pt{
                & \ar[dr]^{g_2} & \\
               \ar[ur]^{g_1} \ar[r]_{h_1} &\ar[u]_-{p} \ar[r]_{h_2} &
            }
}.\]
Then there exists a distributive law $\lambda \: \T_2 \bicomp{\Perm}\T_1 \to \T_1 \bicomp{\Perm} \T_2$, defined by associating to a composable pair $\tr{f\in \T_2}\tr{g\in \T_1}$ the factorisation of $f \poi g$ in $\PS$, which yields $\PS \cong \T_1 \bicomp{\Perm} \T_2$.
\end{proposition}

We now give some preliminary lemmas that are instrumental for the definition of the factorisation and the proof of the main result. We begin by showing how string diagrams of $\Com$ and $\T$ are formally interpreted as arrows of $\LwA{\Sigma}$.

\begin{lemma}\label{lemma:LawvereSubprops} ~
\begin{itemize}
\item $\Com$ is the sub-PROPs of $\LwA{\Sigma}$ whose arrows are lists of variables. The inclusion of $\Com$ in $\LwA{\Sigma}$ is the PROP morphism $\Com \to \LwA{\Sigma}$ defined on generators of $\Com$ by
    \[ \Bcomult \ \mapsto \ \tpl{x_1,x_1} \: 1 \to 2 \qquad \qquad \Bcounit \ \mapsto \ \tpl{} \: 0 \to 1. \]
\item $\T$ is the sub-PROPs of $\LwA{\Sigma}$ whose arrows are linear lists. The inclusion of $\T$ in $\LwA{\Sigma}$ is the PROP morphism $\T \to \LwA{\Sigma}$ defined on generators of $\T$ by
    \[ \OpDiag \ \mapsto \ \tpl{o(x_1,\dots,x_n)} \: n \to 1 \qquad \qquad (o\: n \to 1) \in \Sigma. \]
\end{itemize}
 \end{lemma}
 \begin{proof} First, it is immediate to verify that lists of variables are closed under composition, monoidal product and include all the symmetries of $\LwA{\Sigma}$: therefore, they form a sub-PROP. The same holds for linear lists.

 We now consider the first statement of the lemma. There is a 1-1 correspondence between arrows $n \tr{f \in \LwA{\Sigma}}m$ that are lists of variables and functions $\ord{m} \to \ord{n}$: the function for $f$ maps $k$, for $1 \leq k \leq m$, to the index $l$ of the variable $x_l$ appearing in position $k$ in $f$. This correspondence yields an isomorphism between the sub-PROP of $\LwA{\Sigma}$ whose arrows are lists of variables and $\Fop$. Composing this isomorphism with the one $\Com \tr{\cong} \Fop$ yields the PROP morphism inductively defined in the statement of the lemma. 

 We now turn to the second statement. By construction, arrows of $\T$ are $\Sigma$-terms: the interpretation of Remark~\ref{rmk:Law_sigmaToCartesianTerms} defines a faithful PROP morphism $\T \to \LwA{\Sigma}$ which coincides with the one of the statement. The image of this PROP morphism are precisely the linear lists in $\LwA{\Sigma}$.  
 \end{proof}


Lemma~\ref{lemma:LawvereSubprops} allows us to use $\LwA{\Sigma}$ as an environment where letting $\Com$ and $\T$ interact. The following statement guarantees the soundness of the interaction described by \eqref{eq:LawDistrBcounit}-\eqref{eq:LawDistrBcomult}.

\begin{lemma}\label{lemma:LawEquationsPreserveSubst} Equations \eqref{eq:LawDistrBcounit} and \eqref{eq:LawDistrBcomult} are sound in $\LwA{\Sigma}$ under the interpretation of $\T$ and $\Com$ as sub-PROPs of $\LwA{\Sigma}$.
\end{lemma}
\begin{proof} We first focus on~\eqref{eq:LawDistrBcounit}. Following the isomorphisms of Lemma~\ref{lemma:LawvereSubprops}, $\OpDiag \in \T[n,1]$ is interpreted as the arrow $\tpl{o(x_1,\dots,x_n)} \in \LwA{\Sigma}[n,1]$ and $\Bcounit \in \Com[0,1]$ as $\tpl{} \in \LwA{\Sigma}[1,0]$. The left-hand side of~\eqref{eq:LawDistrBcounit} is then the composite $\tpl{o(x_1,\dots,x_n)} \poi \tpl{} \in \LwA{\Sigma}[n,0]$, which is equal by definition to $\tpl{} \in \LwA{\Sigma}[n,0]$. Therefore, the left- and right-hand side of~\eqref{eq:LawDistrBcounit} are the same arrow of $\LwA{\Sigma}$.

 It remains to show soundness of~\eqref{eq:LawDistrBcomult}. Following Lemma~\ref{lemma:LawvereSubprops}, the left-hand side $\OpDiag \poi \Bcomult$ is interpreted in $\LwA{\Sigma}$ as the composite $\tpl{o(x_1,\dots,x_n)} \poi \tpl{x_1,x_1}$. The right-hand side is interpreted as $\tpl{x_1,\dots,x_n,x_1,\dots,x_n} \poi \tpl{o(x_1,\dots,x_n), o(x_{n+1},\dots,x_{n+n})}$. By definition, both composites are equal to $\tpl{o(x_1,\dots,x_n), o(x_1,\dots,x_n)}$ in $\LwA{\Sigma}$. Therefore,~\eqref{eq:LawDistrBcomult} is also sound in~$\LwA{\Sigma}$.   \end{proof}

It is useful to observe that \eqref{eq:LawDistrBcounit}-\eqref{eq:LawDistrBcomult} allows us to copy and discard not only the generators but arbitrary string diagrams of $\T$.

\begin{lemma}\label{lemma:LawEquationsDerived} Suppose $\diagD$ is a string diagram of $\T$. Then the following holds in $\T + \Com$ quotiented by \eqref{eq:LawDistrBcounit}-\eqref{eq:LawDistrBcomult}.
\noindent
\begin{multicols}{2}\noindent
\begin{equation}\label{eq:LawDistrBcomultGen} \tag{Law3}
\lower18pt\hbox{$\includegraphics[height=1.5cm]{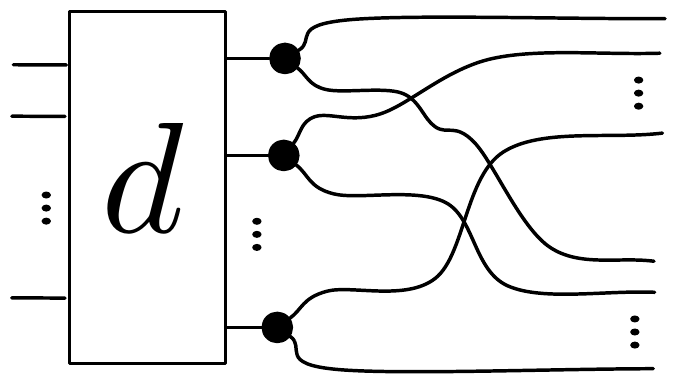}$} = \lower20pt\hbox{$\includegraphics[height=1.7cm]{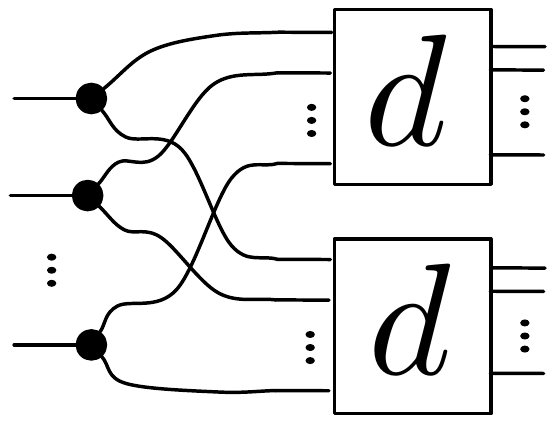}$}
\end{equation}
{
\begin{equation}\label{eq:LawDistrBcounitGen}
\lower7pt\hbox{$\includegraphics[height=.7cm]{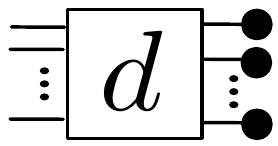}$} =  \lower8pt\hbox{$\includegraphics[height=.8cm]{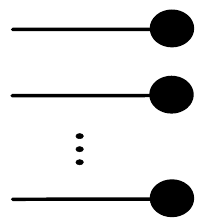}$}
\tag{Law4}
\end{equation}\vspace{-1.2cm}
}
\end{multicols}
\end{lemma}
\begin{proof}
The proof is by induction on $\diagD$. For~\eqref{eq:LawDistrBcomultGen}, the base cases of $\Idnet$ and $\symNet$ follow by \eqref{eq:sliding1} and \eqref{eq:sliding2}. The base case of $\OpDiag$, for $o$ a generator in $\Sigma$, is given by \eqref{eq:LawDistrBcomult}. The inductive cases of composition by $\poi$ and $\tns$ immediately follow by inductive hypothesis. The proof of \eqref{eq:LawDistrBcounitGen} is analogous.
\end{proof}

We can now show the factorisation lemma.
\begin{lemma} \label{lemma:LawvereFactorisation} Any arrow $n\tr{f \in \LwA{\Sigma}}m$ has a factorisation $n \tr{\hat{c}\in \Com}\tr{\hat{d}\in \T}m$ which is unique up-to permutation. \end{lemma}
\begin{proof}
Since the cartesian theory generating $\LwA{\Sigma}$ has no equations, $n \tr{f}m$ is just a list of cartesian $\Sigma$-terms $\tpl{t_1,\dots,t_m}$. The factorisation consists in replacing all variables appearing in $\tpl{t_1,\dots,t_m}$ with fresh ones $x_1,\dots,x_z$, so that no repetition occurs: this gives us the second component of the decomposition as a linear list $z\tr{\hat{d}\in \T}m$. The first component $\hat{c}$ will be the list $n\tr{\var{t_1,\dots,t_m}\in \Com}z$ of variables originally occurring in $f$, so that post-composition with $\hat{d}$ yields $\tpl{t_1,\dots,t_m}$. The reader can consult Appendix~\ref{app:Lawvere} for further details.
\end{proof}

\begin{example} Suppose $\T$ is the free PROP on signature $\Sigma = \{o \: 2\to 1\}$. The arrow $\tpl{ o(x_1,o(x_1,x_3)), x_2} \: 3\to 2$ in $\LwA{\Sigma}$ is factorised as $\tpl{x_1,x_1,x_3,x_2} \poi \tpl{o(x_1,o(x_2,x_3)),x_4}$. Observe that this decomposition is unique only-up to permutation: for instance, another factorisation is {$\tpl{x_1,x_1,x_2,x_3} \poi \tpl{o(x_1,o(x_2,x_4)),x_3}$}  --- in this case, the mediating permutation is $\tpl{x_1,x_2,x_4,x_3}$.
\end{example}

We now have all the ingredients to conclude the proof of our main statement.

\begin{proof}[Proof of Theorem~\ref{Th:LawvereCompositePROP}] Using the conclusion of Lemma~\ref{lemma:LawvereFactorisation}, Proposition~\ref{prop:fact->distrLaw} gives us a distributive law $\lambda \: \T \bicomp{\Perm} \Com \to \Com \bicomp{\Perm} \T$ such that $\LwA{\Sigma} \cong \Com \bicomp{\Perm} \T$. It remains to show that \eqref{eq:LawDistrBcounit}-\eqref{eq:LawDistrBcomult} allow to prove all the equations arising from $\lambda$. By Proposition~\ref{prop:fact->distrLaw}, $\lambda$ maps a composable pair $n\tr{d \in \T}\tr{c\in \Com}m$ to the factorisation $n\tr{c' \in \Com}\tr{d' \in\T}m$ of $d \poi c$ in $\LwA{\Sigma}$, calculated according to Lemma~\ref{lemma:LawvereFactorisation}. The corresponding equation generated by $\lambda$ is $d \poi c = c' \poi d'$, with $d,c,c',d'$ now seen as string diagrams of $\T + \Com$. The equational theory of $\LwA{\Sigma} \cong \Com \bicomp{\Perm} \T$ consists of all the equations arising in this way plus those of $\T + \Com$. What we need to show is that

\starredtext{the string diagrams $d \poi c$ and $c' \poi d'$ are equal modulo the equations of $\T + \Com$ and \eqref{eq:LawDistrBcounit}-\eqref{eq:LawDistrBcomult}.}{$\dag$}

Since our factorisation is unique up-to permutation, it actually suffices to show a weaker statement, namely that

\starredtext{\emph{there exists} a factorisation $n\tr{c'' \in \Com}\tr{d'' \in\T}m$ of $d \poi c$ in $\LwA{\Sigma}$ such that the string diagrams $d \poi c$ and $c'' \poi d''$ are equal modulo the equations of $\T + \Com$ and \eqref{eq:LawDistrBcounit}-\eqref{eq:LawDistrBcomult}.}{$\ddag$}

Statement $(\ddag)$ implies $(\dag)$ because, by uniqueness of the factorisation $c' \poi d'$ up-to permutation, there exists $\tr{p \in \Perm}$ such that $d' = p \poi d''$ and $c'' = c' \poi p$ in $\LwA{\Sigma}$. Since $p$ is an arrow of both sub-PROPs $\T$ and $\Com$, the first equality also holds in $\T$ and the second in $\Com$. This means that $c' \poi d' = c' \poi p \poi d'' = c'' \poi d''$ in $\T + \Com$.

Therefore, we turn to a proof of $(\ddag)$. We describe a procedure to transform the string diagram $\tr{d \in \T}\tr{c\in \Com}$ into the form $\tr{c''\in \Com}\tr{d''\in \T}$ by only using the equations in $\T + \Com$ plus \eqref{eq:LawDistrBcounit}-\eqref{eq:LawDistrBcomult}. Lemmas \ref{lemma:LawvereSubprops} and ~\ref{lemma:LawEquationsPreserveSubst} guarantee that $d \poi c = c''\poi d''$ as arrows of $\LwA{\Sigma}$.
\begin{enumerate}
\item First, there is a preparatory step in which we move all symmetries to the outmost part of the string diagram $d \poi c$, to ease the application of \eqref{eq:LawDistrBcounit}-\eqref{eq:LawDistrBcomult}. By definition, $d$ only contains components of the kind $\OpDiag \: k \to 1$, for $o \in \Sigma$, $\symNet \: 2 \to 2$ and $\Idnet \: 1 \to 1$. We can move all components $\symNet$ to the left of components $\OpDiag$ by repeatedly applying the following instances of~\eqref{eq:sliding2}:
\begin{eqnarray*}
\lower8pt\hbox{$\includegraphics[height=.8cm]{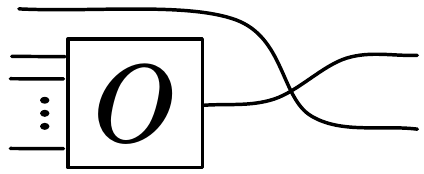}$} \Longrightarrow \lower8pt\hbox{$\includegraphics[height=.8cm]{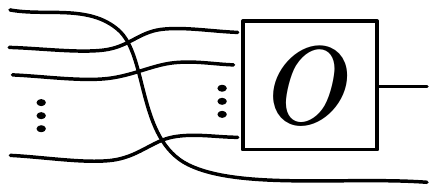}$}
&\qquad \qquad&
\lower8pt\hbox{$\includegraphics[height=.8cm]{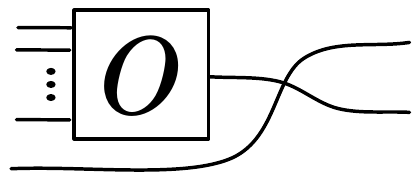}$} \Longrightarrow \lower8pt\hbox{$\includegraphics[height=.8cm]{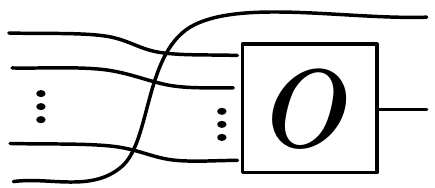}$}
\end{eqnarray*}
    The result is a string diagram $p \poi \bar{d} \poi c'$, where $p$ only contains components $\symNet$ and $\Idnet$ --- i.e., it is a string diagram of $\Perm$ --- and $\bar{d}$ is a string diagram of $\T$ where $\symNet$ does not appear.
\begin{eqnarray*}
\lower10pt\hbox{$\includegraphics[height=1cm]{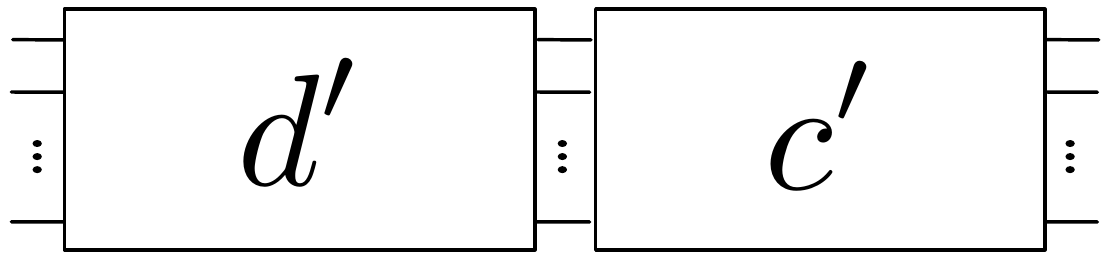}$} \qquad \Longrightarrow \qquad \lower10pt\hbox{$\includegraphics[height=1cm]{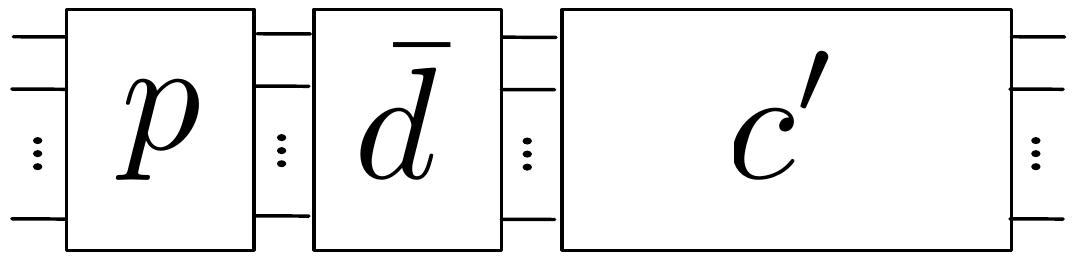}$}
\end{eqnarray*}
    We then perform a symmetric transformation on the string diagram $c$. By definition, $c$ contains components $\Bcomult \: 1 \to 2$, $\Bcounit \: 0 \to 1$, $\symNet \: 2 \to 2$ and $\Idnet \: 1 \to 1$. We can move all components $\symNet$ to the right of any component $\Bcomult$ and $\Bcounit$ by repeatedly applying the following instances of~\eqref{eq:sliding2}:
\begin{eqnarray*}
\lower7pt\hbox{$\includegraphics[height=.7cm]{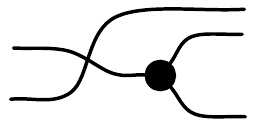}$} \Longrightarrow \lower7pt\hbox{$\includegraphics[height=.7cm]{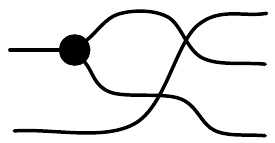}$}
&\qquad \qquad&
\lower5pt\hbox{$\includegraphics[height=.5cm]{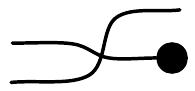}$} \Longrightarrow \lower4pt\hbox{$\includegraphics[height=.4cm]{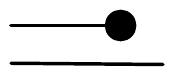}$}
 \\ \\
\lower7pt\hbox{$\includegraphics[height=.7cm]{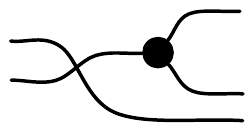}$} \Longrightarrow \lower7pt\hbox{$\includegraphics[height=.7cm]{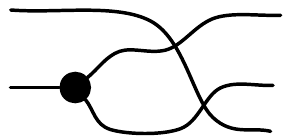}$}
&\qquad \qquad&
\lower5pt\hbox{$\includegraphics[height=.5cm]{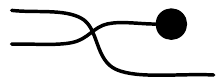}$} \Longrightarrow \lower4pt\hbox{$\includegraphics[height=.4cm]{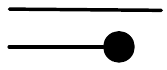}$}
\end{eqnarray*}
    The result is a string diagram $p \poi \bar{d} \poi \bar{c} \poi p'$, where $\bar{c}$ is a string diagram of $\Com$ in which $\symNet$ does not appear and $p'$ is a string diagram of $\Perm$.
\begin{eqnarray*}
\lower10pt\hbox{$\includegraphics[height=1cm]{graffles/LawFinitePresDer2.pdf}$} \qquad \Longrightarrow \qquad \lower10pt\hbox{$\includegraphics[height=1cm]{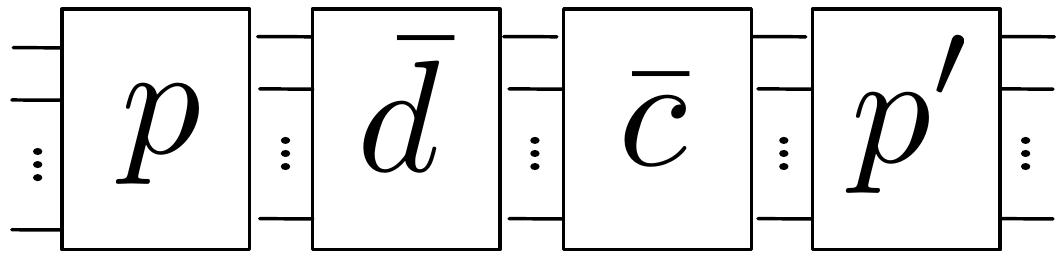}$}
\end{eqnarray*}
\item We now make $\bar{d}$ and $\bar{c}$ interact. First note that, since $\bar{d}$ does not contain $\symNet$ and all generators $o \in \Sigma$ have coarity $1$, $\bar{d}$ must the $\tns$-product $\bar{d_1} \tns \dots \tns \bar{d_z}$ of string diagrams $\bar{d_i} \: k_i \to 1$ of $\T$.
\begin{eqnarray*}
\lower10pt\hbox{$\includegraphics[height=1cm]{graffles/LawFinitePresDer3.pdf}$} \qquad = \qquad \lower20pt\hbox{$\includegraphics[height=2cm]{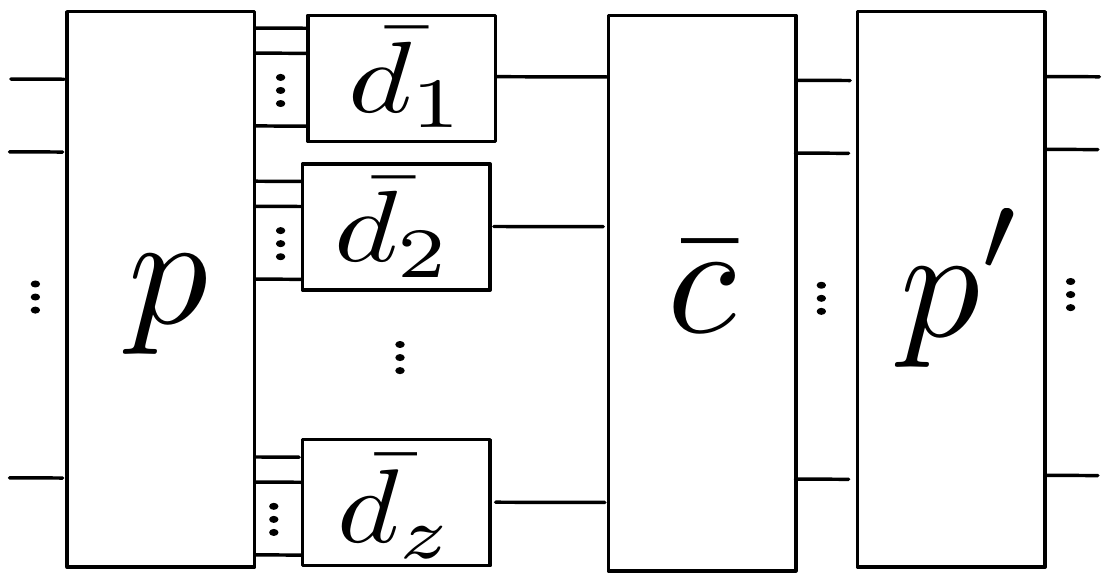}$}
\end{eqnarray*}
    For analogous reasons, $\bar{c}$ is also a $\tns$-product $\bar{c_1} \tns \dots \tns\bar{c_z}$ where, for $1 \leq i \leq z$,
    \begin{eqnarray}\label{eq:LawchoicesforbarCi}
\lower9pt\hbox{$\includegraphics[height=.9cm]{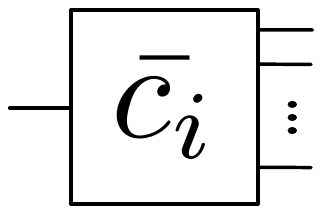}$} \quad = \quad \lower10pt\hbox{$\includegraphics[height=1cm]{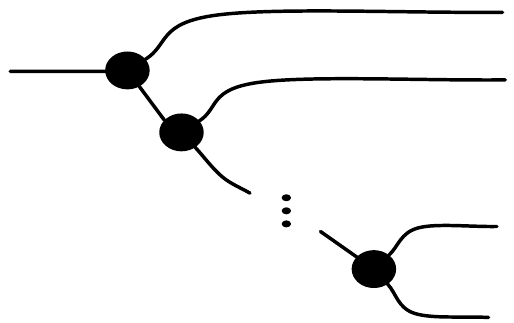}$}
 \qquad\qquad \text{ or } \qquad \qquad
 \lower9pt\hbox{$\includegraphics[height=.9cm]{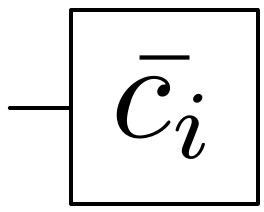}$} \quad = \quad \lower0pt\hbox{$\includegraphics[height=.2cm]{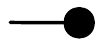}$} \ .
\end{eqnarray}
We thus can present $\bar{c}$ as follows:
\begin{eqnarray*}
\lower20pt\hbox{$\includegraphics[height=2cm]{graffles/LawFinitePresDer4.pdf}$} \qquad = \qquad \lower20pt\hbox{$\includegraphics[height=2cm]{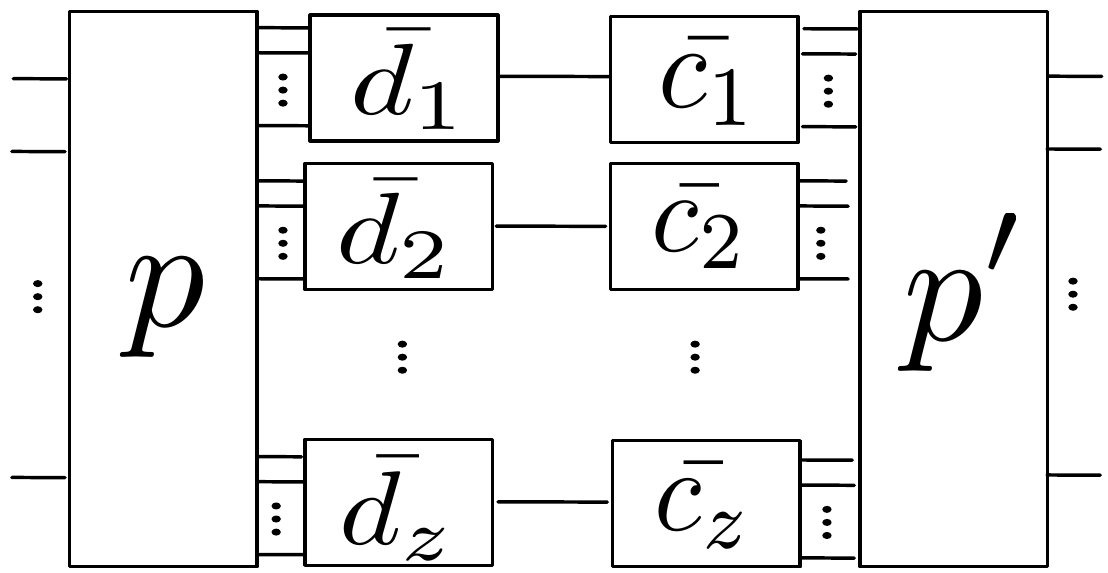}$}.
\end{eqnarray*}
We are now in position to distribute each $\bar{d_i}$ over the corresponding $\bar{c_i}$. Suppose first $\bar{c_i}$ satisfies the left-hand equality in \eqref{eq:LawchoicesforbarCi}. By assumption, all the equations of $\T + \Com$, \eqref{eq:LawDistrBcounit} and \eqref{eq:LawDistrBcomult} hold. Thus, by Lemma~\ref{lemma:LawEquationsDerived}, also \eqref{eq:LawDistrBcomultGen} holds. Starting from $\bar{d_i}\poi\bar{c_i}$, we can iteratively apply \eqref{eq:LawDistrBcomultGen} (and \eqref{eq:sliding1}) to obtain a string diagram of shape $\tr{\in \Com}\tr{\in \T}$:
\begin{eqnarray*}
\lower11pt\hbox{$\includegraphics[height=1.1cm]{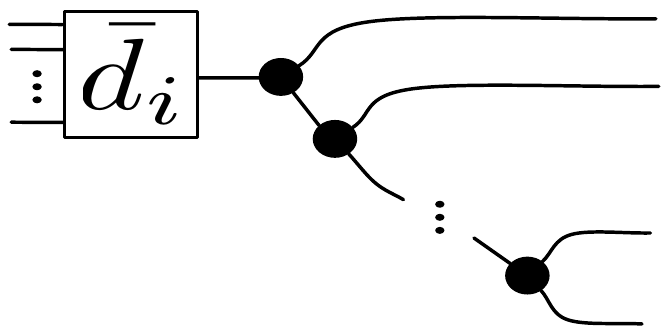}$} \Rightarrow  \lower17pt\hbox{$\includegraphics[height=1.7cm]{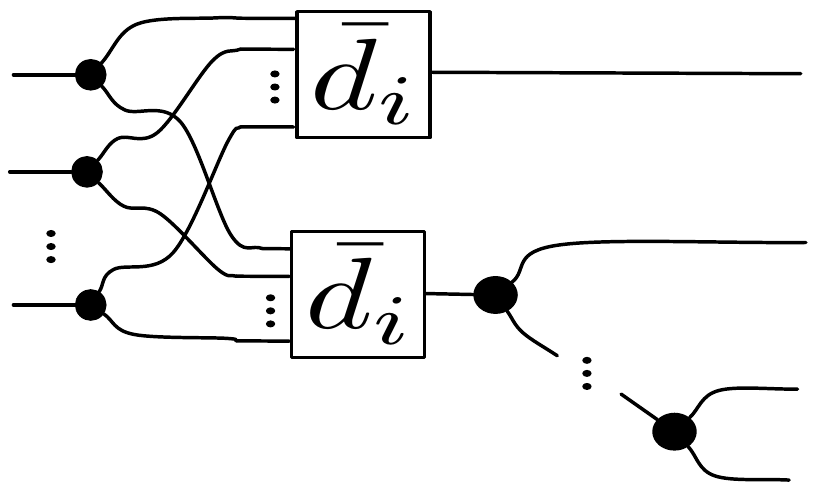}$} \Rightarrow  \lower22pt\hbox{$\includegraphics[height=2.2cm]{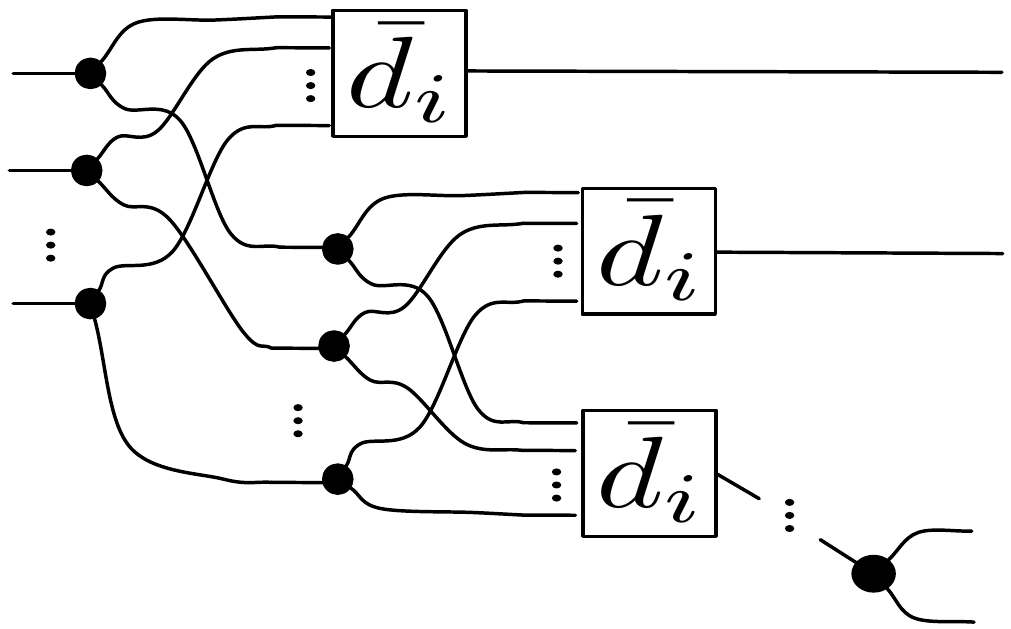}$} \Rightarrow \ \cdots \ \Rightarrow  \lower26pt\hbox{$\includegraphics[height=2.9cm]{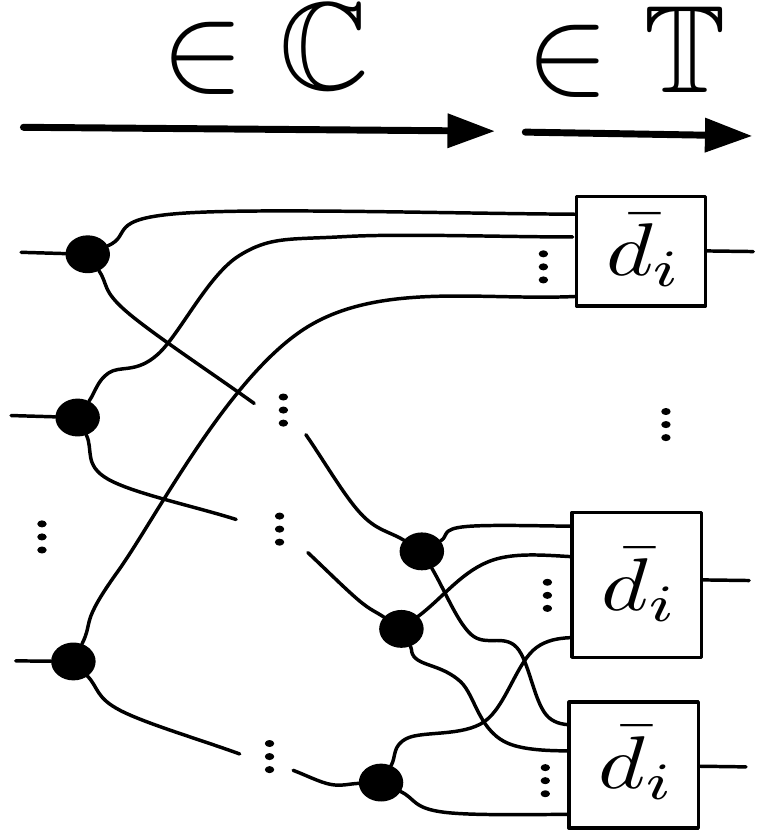}$}.
\end{eqnarray*}
In the remaining case, $\bar{c_i}$ satisfies the right-hand equality in \eqref{eq:LawchoicesforbarCi}. Then, one application of \eqref{eq:LawDistrBcounit} also gives us a string diagram of shape $\tr{\in \Com}\tr{\in \T}$.
\begin{eqnarray*}
\lower8pt\hbox{$\includegraphics[height=.8cm]{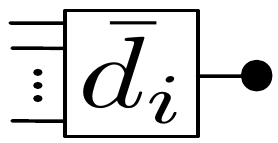}$} \Rightarrow  \lower8pt\hbox{$\includegraphics[height=1.3cm]{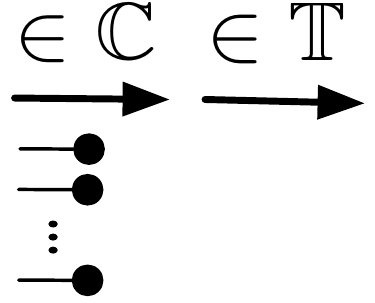}$}
\end{eqnarray*}
Applying the above transformations for each $\bar{d_i}\poi\bar{c_i}$ yields a string diagram of the desired shape $\tr{c''\in \Com}\tr{d''\in \T}$.
\begin{eqnarray*}
\lower20pt\hbox{$\includegraphics[height=2cm]{graffles/LawFinitePresDer5.pdf}$} \qquad \Longrightarrow \qquad \lower20pt\hbox{$\includegraphics[height=2.6cm]{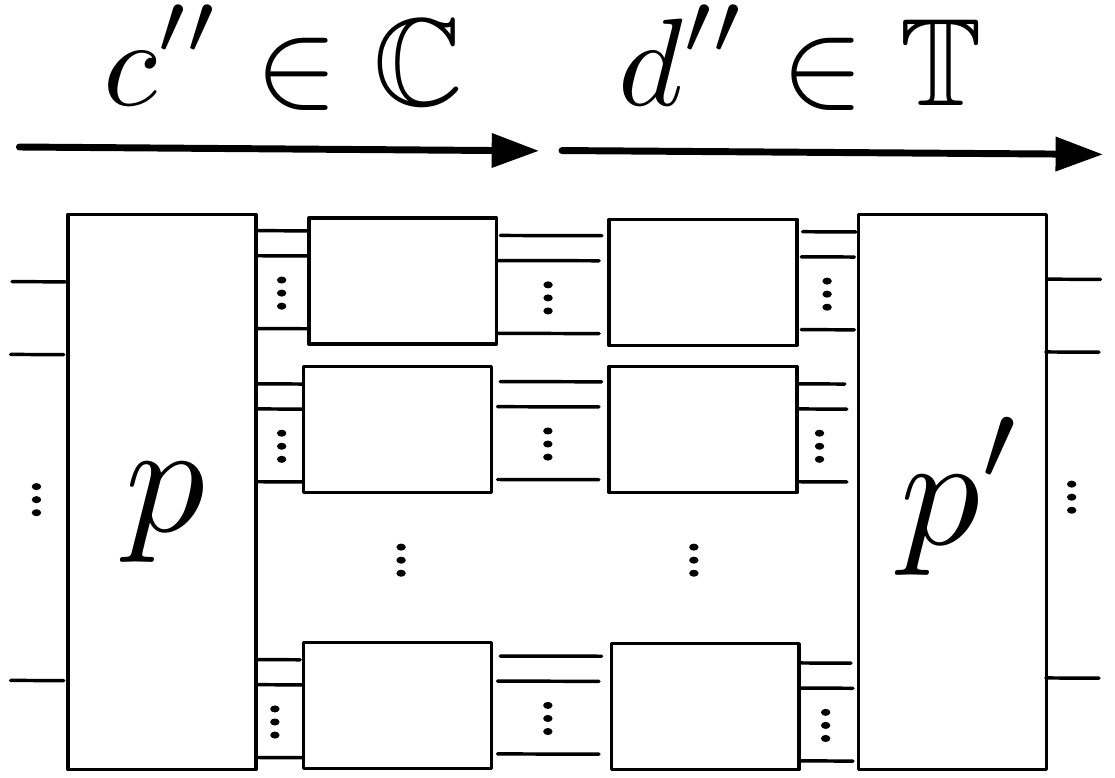}$}.
\end{eqnarray*}
\end{enumerate}
Observe that all the transformations that we described only used equations in $\T + \Com$, \eqref{eq:LawDistrBcounit} and \eqref{eq:LawDistrBcomult}. This concludes the proof of $(\ddag)$ and thus of the main theorem.
\end{proof}
\begin{remark}[Term graphs] The factorised form $\tr{\in \Com}\tr{\in \T}$ of arrows in $\Com \bicomp{\Perm} \T$ represents cartesian $\Sigma$-terms by their syntatic tree --- the $\T$-part --- with the possibility of explictly indicating which variables are shared by sub-terms --- the $\Com$-part. Interestingly, by weakening the structure of $\Com \bicomp{\Perm} \T$, we are able to capture a different, well-known representation for cartesian $\Sigma$-terms, namely \emph{term graphs}, which are acyclic graphs labeled over $\Sigma$. With respect to the standard tree representation, the benefit of term graphs is that the sharing of any common sub-term can be represented explicitly, making them particularly appealing for efficient rewriting algorithms, see e.g.~\cite{TermGraphSurvey} for a survey on the subject. As shown in~\cite{GadducciCorradini-termgraph}, $\Sigma$-term graphs are in 1-1 correspondence with the arrows of the free \emph{gs-monoidal category} generated by $\Sigma$. This categorical construction actually amounts to forming the sum of PROPs $\Com + \T$. Thus the only difference between term graphs and the representation of terms given by $\Com \bicomp{\Perm} \T$ is in the validity of laws~\eqref{eq:LawDistrBcounit}-\eqref{eq:LawDistrBcomult}. Intuitively, these equations do not hold because term graphs are \emph{resource-sensitive}, e.g. they are not invariant with respect to the representation of a shared sub-term as distinct copies of it. 
\end{remark}

\section{Fibered Sum of PROPs} \label{sec:pushout}
So far we have been focusing on two PROP operations: sum and composition. Another useful construction is the \emph{fibered sum} (pushout) of PROPs. We focus here on the case of pushouts taken along a common sub-PROP $\T_1 \rightrightarrows \T_2, \T_3$. The leading intuition is that, while composition $\T_2 \bicomp{\Perm} \T_3$ would quotient the sum $\T_2 + \T_3$ by new equations, the fibered sum only identifies the structure $\T_1$ which is in common between $\T_2$ and $\T_3$. This kind of construction is typical in algebra, from geometric gluing constructions of topological spaces to amalgamated free products of groups.
\begin{proposition}\label{prop:SMTpushout} Suppose that $\T_1$, $\T_2$ and $\T_3$ are PROPs freely generated by $(\Sigma_1,E_1)$, $(\Sigma_2,E_2)$ and $(\Sigma_3,E_3)$ respectively, with $\Sigma_1 \subseteq \Sigma_i$ and $E_1 \subseteq E_i$ for $i \in \{2,3\}$. Let the following be the pushout along the PROP morphisms $\Phi \: \T_1 \to \T_2$ and $\Psi \: \T_1 \to \T_3$ defined by interpreting a $\Sigma_1$-term modulo $E_1$ as a $\Sigma_i$-term modulo $E_i$, for $i \in \{2,3\}$.
$$\xymatrix@=10pt{ \T_1 \ar[r]^{\Phi} \ar[d]_{\Psi} & \T_2 \ar[d]^{\Lambda} \\ \T_3 \ar[r]_{\Gamma} & \T}$$
Then $\T$ is presented by the SMT with signature $\Sigma_2 \uplus (\Sigma_3 \gls{setminus} \Sigma_1)$ and equations $E_2 \uplus (E_3\setminus E_1)$. The PROP morphisms $\Lambda$ and $\Gamma$ are the obvious interpretations of a term of the smaller theory as one of the larger theory.
\end{proposition}
\begin{proof} Pushouts in $\PROP$ can be calculated as in $\Cat$. The case of the span $\T_2 \tl{} \T_1 \tr{} \T_3$ under consideration is simplified by the fact that only categories with the same set of objects $\N$ and identity-on-object functors are involved. This means that the pushout $\T$ has $\N$ as set of objects and $\Lambda$, $\Gamma$ are identity-on-objects functors. Concerning the arrows, $\T$ is the quotient of $\T_2 + \T_3$ by the congruence generated by
\[ \tr{g\in \T_2}\ \sim\ \tr{h\in \T_3} \text{ if and only if there exists} \tr{f \in \T_1} \text{ such that }\Phi(\tr{f})\ =\ \tr{g} \text{ and } \Psi(\tr{f})\ =\ \tr{h}. \]
In terms of the equational description of $\T_1$, $\T_2$ and $\T_3$, this condition means that any $\Sigma_2$-term modulo $E_2$ and $\Sigma_3$-term modulo $E_3$ that are the same $\Sigma_1$-term modulo $E_1$ should be treated as the same arrow in $\T$; that amounts to saying that arrows of $\T$ are generated by the SMT described in the statement. Being $\T$ a quotient of $\T_2 + \T_3$, the cospan $\T_2 \tr{} \T \tl{} \T_3$ is formed by mapping the arrows of $\T_2$, $\T_3$ to the corresponding equivalence class. \end{proof}

In our applications $\T_2$ and $\T_3$ will typically be PROPs modeling, respectively, spans and cospans of the arrows of $\T_1$. Their fibered sum will result in a PROP $\T$ whose arrows model \emph{relations} of some kind. We shall see a detailed example of this construction in \S~\ref{sec:ER} below, where we give a presentation of the PROP of equivalence relations using a pushout construction. \S~\ref{sec:PER} will extend our approach to the case of partial equivalence relations. Fibered sum will also be pivotal for the developments of the next chapter, in which we will push out the PROPs modeling spans and cospans of matrices to obtain the PROP of linear subspaces.

\subsection{Case Study I: Equivalence Relations}\label{sec:ER}

In this section we detail an illustrative example of how fibered sum can be used to characterise the PROP of equivalence relations.

\begin{definition}
Let $\gls{PROPER}$ be the PROP with arrows $n \to m$ the equivalence relations on $\ord{n}\uplus\ord{m}$. Given $e_1 \: n \to z$, $e_2 \: z \to m$, their composite $e_1 \poi e_2$ is defined as follows. First we compose $e_1$ and $e_2$ as relations: 
\begin{eqnarray}\label{eq:defCompRelPrel}
e_1 \relcomp e_2 & \df & \{(v,w) \mid \exists u.\; (v,u)\in e_1 \wedge (u,w)\in e_2\}.
\end{eqnarray}
Then, we consider the equivalence relation freely generated by $e_1 \relcomp e_2$ (that is, the reflexive, symmetric and transitive closure of $e_1 \relcomp e_2$, denoted by $\eqr{e_1 \relcomp e_2}$) and discard any pair containing at least one element from $\ord{z}$. 
\begin{eqnarray} \label{eq:defCompRel}
e_1 \poi e_2 & \df & \{(u,w) \mid (u,w) \in \eqr{e_1 \relcomp e_2} \wedge u,w \in \ord{n}\uplus\ord{m} \}
\end{eqnarray}
The monoidal product $e_1 \tns e_2$ on equivalence relations $e_1, e_2$ is just their disjoint union.
\end{definition}


Our approach in characterising $\ER$ stems from the observation (e.g. mentioned in~\cite{Bruni01somealgebraic}) that any cospan $n \tr{p\in F}\tl{q \in \F}m$ gives rise to an equivalence relation $e$ on $\ord{n}\uplus\ord{m}$ as follows:
\begin{eqnarray}\label{eq:CospanToER}
(v,w) \in e \qquad \text{ iff } \qquad
\begin{cases}
p(v) = q(w) &\mbox{if } v \in \ord{n}, w \in \ord{m} \\
q(v) = p(w) & \mbox{if } v \in \ord{m}, w \in \ord{n} \\
p(v) = p(w) & \mbox{if } v , w \in \ord{n} \\
q(v) = q(w) & \mbox{if } v , w \in \ord{m}.
\end{cases}
\end{eqnarray}
This interpretation is particularly neat if we take the string diagrammatic presentation of a cospan $n \tr{\F}\tl{\F}m$ as an arrow $n \tr{\in \FROB} m$, using the characterisation $\FROB \cong \Mon \bicomp{\Perm} \Com \cong \F \bicomp{\Perm} \Fop$ (Example~\ref{ex:equationalprops} and~\ref{ex:distrlawsyntactic}\ref{ex:distrlawsyntactic3}). For instance, consider the following arrow $5 \tr{\in \FROB} 7$:
\begin{equation}
\lower20pt\hbox{\includegraphics[height=2cm]{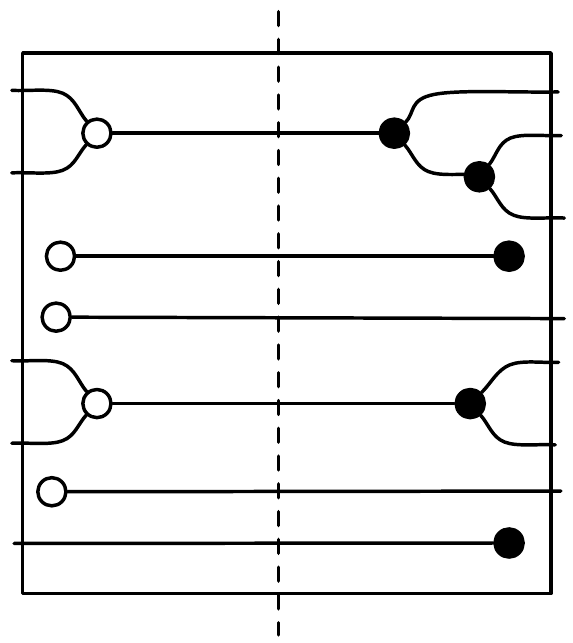}.\label{eq:spider}}
\end{equation}
The dotted line emphasizes that we are visualising this arrow as a string diagram in the factorised form $n \tr{\in \Mon}\tr{\in \Com} m$ yielded by the decomposition $\FROB \cong \Mon \bicomp{\Perm} \Com$. The string diagram~\eqref{eq:spider} defines an equivalence relation $e$ on $\ord{5}\uplus\ord{7}$ by letting $(v,w) \in e$ if and only if the port associated with $v$ and the one associated with $w$ are linked in the graphical representation. For instance, $1,2 \in \ord{5}$ on the left boundary are in the same equivalence class as $1,2,3 \in \ord{7}$ on the right boundary, whereas $5 \in \ord{5}$ and $4 \in \ord{7}$ are the only members of their equivalence class.

Observe that the sub-diagram $\lower5pt\hbox{$\includegraphics[height=.5cm]{graffles/WBbone.pdf}$}$ in \eqref{eq:spider} does not have any influence on the interpretation described above. In fact, the following diagrams represent the same equivalence relations.
\begin{equation}\label{eq:spidersbone}
\lower20pt\hbox{$\includegraphics[height=2cm]{graffles/spider.pdf}$} \qquad \text{and} \qquad \lower20pt\hbox{$\includegraphics[height=2cm]{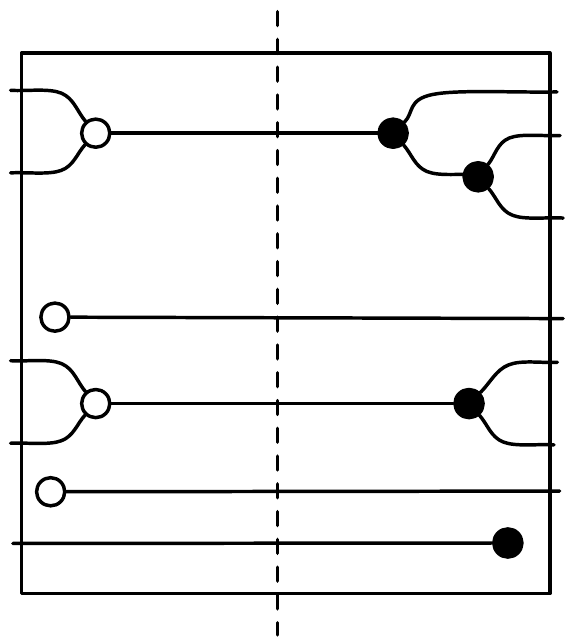}$}
\end{equation}
The above discussion substantiates the claim that $\FROB$ is \emph{almost} the good candidate as a theory of equivalence relations. The additional axiom that we need is \eqref{eq:bwbone} --- first appearing in Example~\ref{ex:equationalprops} --- which tells that any sub-diagram $\lower5pt\hbox{$\includegraphics[height=.5cm]{graffles/WBbone.pdf}$}$ provides redundant information which can be removed. We call $\gls{PROPIFROB}$ (\textbf{i}rredundant \textbf{Fr}obenius algebras) the PROP obtained by quotienting $\FROB$ by~\eqref{eq:bwbone}.

\begin{theorem}\label{th:IFROB=ER} $\IFROB \cong \ER$. \end{theorem}

The rest of this section is devoted to proving Theorem~\ref{th:IFROB=ER}. The argument will essentially rely on showing that $\IFROB$ and $\ER$ are pushouts of isomorphic spans in $\PROP$.

\subsubsection{$\IFROB$ as a Fibered Sum}

We begin with the characterisation of $\IFROB$. Following the above discussion, the idea is that $\IFROB$ will result as the fibered sum of $\FROB$ and the smallest theory containing~\eqref{eq:bwbone}.

For defining such a theory, recall the theories $\UNIT$ and $\COUNIT$, introduced in Example~\ref{ex:distrlawsyntactic}.\ref{ex:distrlawsyntactic1} and~\ref{ex:partialfunctions}, which characterise $\Inj$ and $\Injop$ respectively. Consider the distributive law $\Inj \bicomp{\Perm} \Injop \to \Injop \bicomp{\Perm} \Inj$ given by pullback in $\Inj$. It yields a PROP $\Injop \bicomp{\Perm} \Inj \cong \COUNIT \bicomp{\Perm} \UNIT$ which, according to Proposition~\ref{prop:SMTforComposition}, is presented by generators $\Wunit$ (from $\UNIT$), $\Bcounit$ (from $\COUNIT$) and the equations encoded by $\Inj \bicomp{\Perm} \Injop \to \Injop \bicomp{\Perm} \Inj$. Giving a finite presentation for these equations is particularly simple because each of $\Inj \cong \UNIT$ and $\Injop \cong \COUNIT$ features a single generator, $\Wunit$ and $\Bcounit$ respectively. It suffices to check a single pullback square, namely the one for the cospan $\Wunit \poi \Bcounit = 0\tr{\initVect} 1 \tl{\initVect} 0$: the resulting equation is \eqref{eq:bwbone}. 


We can thereby conclude that $\COUNIT \bicomp{\Perm} \UNIT$ is presented by generators $\Wunit$, $\Bcounit$ and the equation \eqref{eq:bwbone}. We now want to merge $\COUNIT \bicomp{\Perm} \UNIT$ and $\FROB$. The common structure that we need to identify is the one of the PROP $\UNIT + \COUNIT$. There are obvious PROP morphisms $\UNIT + \COUNIT \to \COUNIT \bicomp{\Perm} \UNIT$ and $\UNIT + \COUNIT \to \FROB$ given by quotienting an arrow of $\UNIT + \COUNIT$ by the equations of $\COUNIT \bicomp{\Perm} \UNIT$ and $\FROB$ respectively. We now take the pushout in $\PROP$ along these morphisms.
\begin{equation} \label{eq:pushoutIFROB}
\raise15pt\hbox{$
\xymatrix@C=40pt{
{\UNIT + \COUNIT} \ar[r] \ar[d] & {\COUNIT \bicomp{\Perm}\UNIT} \ar[d] \\
{\FROB} \ar[r] & {\IFROB}
}$}
\end{equation}
By Proposition~\ref{prop:SMTpushout}, $\IFROB$ is the quotient of $\FROB$ by \eqref{eq:bwbone}. Also, the morphisms $\COUNIT \bicomp{\Perm}\UNIT \to \IFROB$ and $\FROB \to \IFROB$ quotient arrows of the source PROP by the equations of $\IFROB$.


\subsubsection{$\ER$ as a Fibered Sum}

We now turn to the pushout construction for $\ER$. The idea is to mimic \eqref{eq:pushoutIFROB} using the concrete representations of $\UNIT$ as $\Inj$ and of $\Mon$ as $\F$, with $\COUNIT$ and $\Com$ characterising their duals $\Inj^{\op}$ and $\F^{\op}$. The resulting diagram is the following, with $\ER$ the candidate concrete representation for $\IFROB$.
\begin{equation} \label{eq:pushoutER}
\raise15pt\hbox{$
\xymatrix@C=40pt{
{\Inj + \Inj^{\op}} \ar[r]^-{[\kappa_1,\,\kappa_2]} \ar[d]_{[\iota_1,\,\iota_2]} & {\Inj^{\op} \bicomp{\Perm} \Inj} \ar[d]^{\SpanToER} \\
{\F \bicomp{\Perm} \Fop} \ar[r]_-{\CospanToER} & {\ER}
}$}
\end{equation}
We define the PROP morphisms in \eqref{eq:pushoutER}:
\begin{itemize}
 \item morphisms $\kappa_1\: \Inj \to \Injop\bicomp{\Perm}\Inj$, $\kappa_2\: \Injop \to \Injop\bicomp{\Perm}\Inj$, $\iota_1\: \Inj \to \F\bicomp{\Perm}\Fop$ and $\iota_2\: \Injop \to \F\bicomp{\Perm}\Fop$ are given by
\[\kappa_1(n\tr{f} m) = (n \tl{\id} n \tr{f} m),\  \kappa_2(n \tr{f} m) = (n \tl{f} m \tr{\id} m ),\]
\[\iota_1(n \tr{f} m) = (n \tr{f} m \tl{\id} m)\text{ and }
\iota_2(n \tr{f} m) = (n \tr{\id} n \tl{f} m).\]
\item the value of $\CospanToER \: \F \bicomp{\Perm} \Fop \to \ER$ on an arrow $n\tr{p \in \F} z \tl{q\in \F} m$ in $\F \bicomp{\Perm} \Fop$ is given according to~\eqref{eq:CospanToER}.
 \item the value of $\SpanToER \: \Injop \bicomp{\Perm} \Inj \to \ER$ on an arrow $n\tl{f \in \Inj} z \tr{g\in \Inj} m$ in $\Injop \bicomp{\Perm} \Inj$ is given by
 \[ (v,w) \in \SpanToER(n\tl{f} z \tr{g} m) \quad \text{iff} \quad \left\{
	\begin{array}{ll}
		 \inv{f}(v) = \inv{g}(w)  & \mbox{ if } v \in \ord{n},w\in \ord{m} \\
		 \inv{f}(w) = \inv{g}(v) & \mbox{ if } w \in \ord{n},v\in \ord{m}  \\
        v=w & \mbox{ otherwhise.}
	\end{array}
\right. \]
\end{itemize}
The definition of the $\kappa$s and $\iota$s are the standard interpretations of an (injective) function as a span/cospan. The mapping $\CospanToER$ implements the idea, described above, of interpreting a cospan as an equivalence relation. For $\SpanToER$, the key observation is that spans of injective functions can also be seen as equivalence relations. Once again, the graphical representation of an arrow of $\Injop \bicomp{\Perm} \Inj$ as a string diagram in $\COUNIT \bicomp{\Perm}\UNIT$ can help visualising this fact. The factorised arrow
\begin{equation}\label{eq:injspan}
\lower13pt\hbox{$\includegraphics[height=1.3cm]{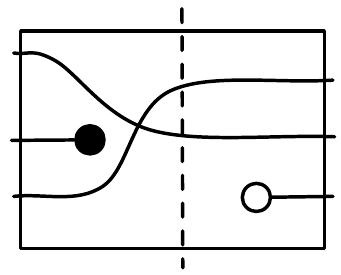}$}
\end{equation}
of $\COUNIT \bicomp{\Perm}\UNIT$ can be interpreted as the equivalence relation associating $1$ on the left boundary with $2$ on the right boundary, $3$ on the left with $1$ on the right and letting $2$ on the left, $3$ on the right be the only representatives of their equivalence class. Note that this interpretation would not work the same way for spans of non-injective functions, as their graphical representation in $\Fop \bicomp{\Perm} \F$ may involve $\Wmult$ and $\Bcomult$ --- more on this in Remark~\ref{rmk:pushoutTrivialTheory}.

One can verify that $\CospanToER$ and $\SpanToER$ are indeed functors (see Appendix~\ref{App:proofsER}). Towards a proof that \eqref{eq:pushoutER} is a pushout, a key step is to understand exactly when two cospans in $\F$ should be identified by $\CospanToER$. Example~\eqref{eq:spidersbone} gives us a lead: two cospans represent the same equivalence relation precisely when their string diagrammatic presentation in $\FROB$ is the same modulo~\eqref{eq:bwbone}. Now, what is the cospan counterpart of applying~\eqref{eq:bwbone}? Recall that this equation arises by a distributive law $\F \bicomp{\Perm} \Fop \to \Fop \bicomp{\Perm} \F$ defined by pullback in $\F$ (Example~\ref{ex:composingSemPROP}\ref{ex:composingSemPROP2}). Thus one could be tempted of saying that the characterising property of $\CospanToER$ should be: two cospans are identified when they have the same pullback.
However, this approach fails because it identifies too many cospans. Here is a counterexample: the leftmost and rightmost cospans below have the same pullback (at the center). However, they represent, respectively, an equivalence relation with one equivalence class and one with two equivalence classes.
\begin{equation}\label{eq:injspan}
\lower13pt\hbox{$\includegraphics[height=1.2cm]{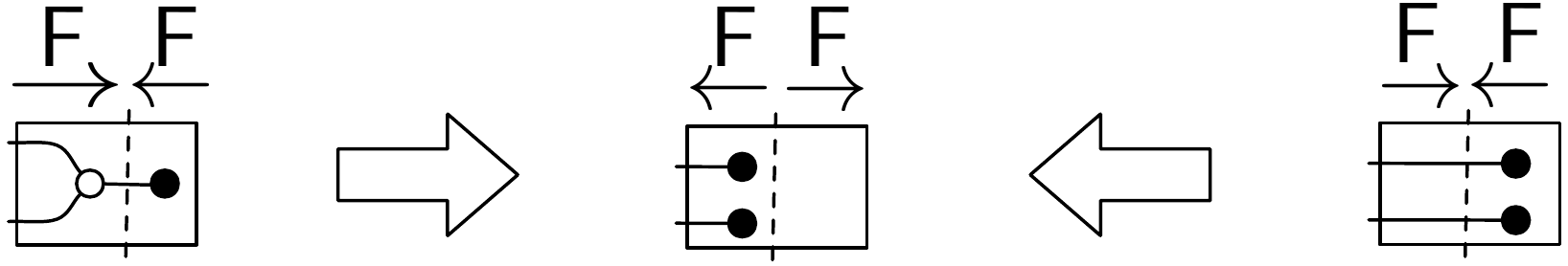}$}
\end{equation}
The correct approach is slightly more subtle. The key observation is that we only need to rewrite $\WBboneText$ as $\ZeronetT$: this does not require to pull back the whole cospan, but just the region where all sub-diagrams of shape $\WBboneText$ lie. To make this formal, recall that $\F$ itself is a composite PROP $\Surj \bicomp{\Perm} \Inj$ (Example~\ref{ex:composingSemPROP}\ref{ex:composingSemPROP1}). We can then factorise a cospan $\tr{\in \F}\tl{\in\F}$ as $\tr{\in \Surj}\tr{\in \Inj}\tl{\in\Inj}\tl{\in\Surj}$. For instance, the leftmost string diagram in~\eqref{eq:spidersbone} is factorised as the leftmost diagram below. On the right below, we pull back the middle cospan $\tr{\in \Inj}\tl{\in\Inj}$. Since $\ZeronetT$ is the pullback of $\WBboneText$, this transformation removes all sub-diagrams of shape $\WBboneText$. The result is the rightmost diagram in~\eqref{eq:spidersbone}.
\begin{equation}\label{eq:spidersboneepimono}
\lower20pt\hbox{$\includegraphics[height=2.4cm]{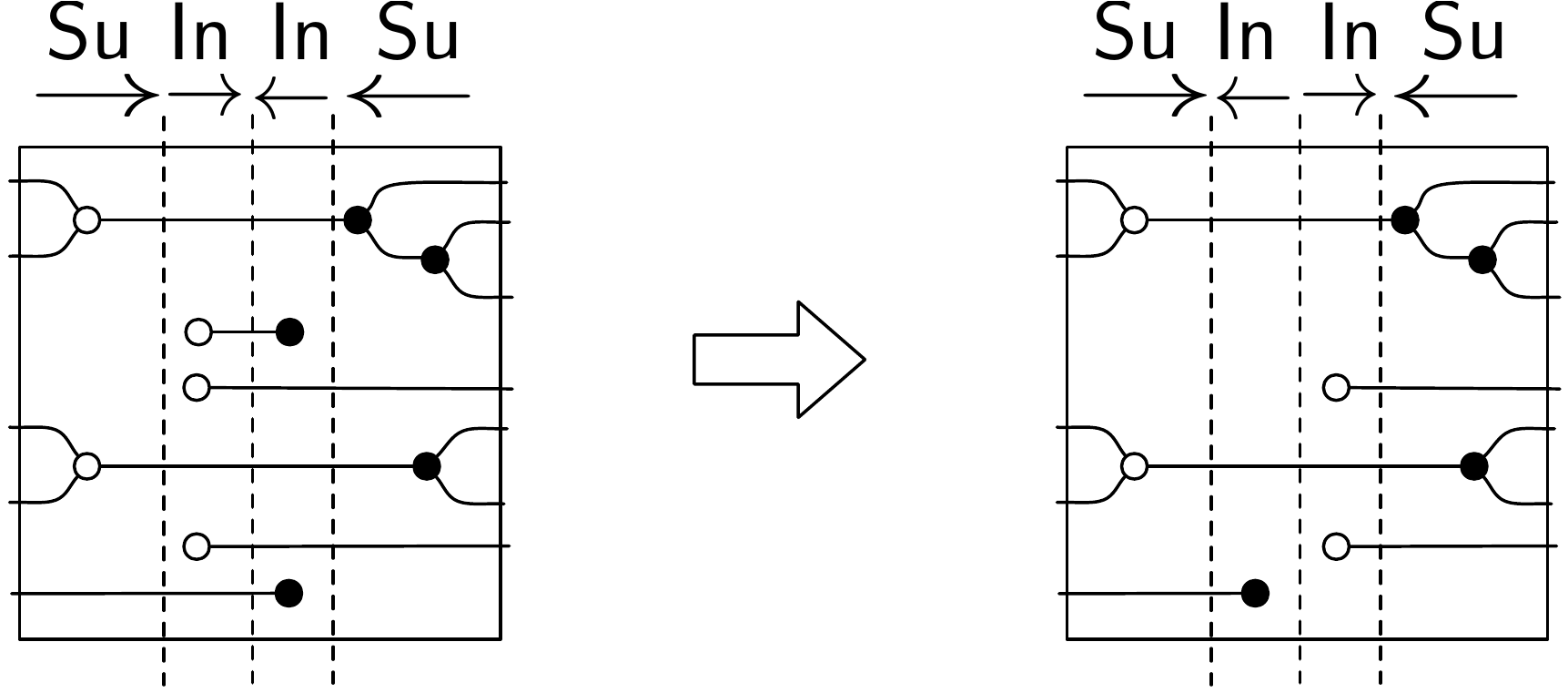}$}
\end{equation}
We fix our observations with the following definition.

\begin{definition}\label{def:equalzeros} We say that two cospans $n \tr{p_1 \in \F} z \tl{q_1 \in \F} m$ and $n \tr{p_2 \in \F} r \tl{q_2 \in \F} m$ are \emph{equal modulo-zeros} if there is an epi-mono factorisation $ \tr{e^1_p \in \Surj}\tr{m^1_p\in\Inj} \tl{m^1_q\in\Inj}\tl{e^1_q \in \Surj}$ of $\tr{p_1} \tl{q_1}$, and one  $ \tr{e^2_p \in \Surj}\tr{m^2_p \in\Inj} \tl{m^2_q\in\Inj}\tl{e^2_q \in \Surj}$ of $ \tr{p_2}  \tl{q_2}$ such that $\tr{m^1_p} \tl{m^1_q}$ and $\tr{m^2_p} \tl{m^2_q}$ have the same pullback and $e^1_p = e^2_p$, $e^1_q =e^2_q$.
\end{definition}

\begin{remark} It may be insightful to give a different view on Definition~\ref{def:equalzeros}. Note that two cospans are equal modulo-zeros precisely when they are in the equivalence relation generated by
\begin{center} $\left(n\tr{p}z\tl{q}m\right) \sim \left(n\tr{p}z\tr{h}z'\tl{h}z\tl{q}m\right)$, where $h$ is an injection. \end{center}
The idea is that $z\tr{h}z'\tl{h}z$ plays a role akin to a repeated use of equation~\eqref{eq:bwbone} in the diagrammatic language: it deflates the codomain of $[p,q] \: n+m \to z$ so as to ``make it surjective'' --- \emph{cf.}~\eqref{eq:spidersbone}. 
\end{remark}

Our proof that \eqref{eq:pushoutER} is a pushout will rely on showing that $\CospanToER$ equalizes two cospans precisely when they are equal modulo-zeros. As a preliminary step, we shall need to establish some properties holding for any $\Gamma$, $\Delta$ and $\mathbb{X}$ making the following diagram commute.
\begin{equation}
\label{eq:arbitraryER}
\raise15pt\hbox{$
\xymatrix@C=40pt{
{\Injplus} \ar[r]^-{[\kappa_1,\,\kappa_2]} \ar[d]_{[\iota_1,\,\iota_2]} & {\Injop \bicomp{\Perm} \Inj} \ar[d]^{\Gamma} \\
{\F\bicomp{\Perm}\Fop} \ar[r]_-{\Delta} & {\mathbb{X}}
}$}
\end{equation}

\begin{lemma}
\label{lemma:arbitraryPROP_ER}
Given a PROP $\mathbb{X}$ and a commutative diagram~\eqref{eq:arbitraryER}, the following hold.
\begin{enumerate}[(i)]
\item If $\tr{p} \tl{q} $ is a cospan in $\Inj$ with pullback (in $\Inj$) $\tl{f} \tr{g}$, then $\Gamma(\tl{f} \tr{g}) = \Delta(\tr{p} \tl{q} )$.
\item If $ \tl{p_1} \tr{q_1} $ and $\tl{p_2} \tr{q_2}$ are cospans in $\Inj$ with the same pullback then $\Delta(\tr{p_1} \tl{q_1} ) = \Delta(\tr{p_2}  \tl{q_2})$.
\item If $ \tr{p_1}  \tl{q_1} $ and $ \tr{p_2}  \tl{q_2} $ are equal modulo-zeros then $\Delta(\tr{p_1} \tl{q_1}) = \Delta(\tr{p_2} \tl{q_2})$.
\item If $ \tl{f} \tr{g} $ is a span in $\Inj$ with pushout (in $\F$) $ \tr{p}  \tl{q} $, then $\Gamma(\tl{f} \tr{g} ) = \Delta(\tr{p} \tl{q})$.
\end{enumerate}
\end{lemma}
\begin{proof}
(i) The statement is given by the following derivation.
\begin{align*}
\Delta(\tr{p}\tl{q}) &= \Delta(\iota_1 p \poi \iota_2 q) \\
&= \Delta\iota_1 p\poi \Delta\iota_2 q \\
&= \Gamma\kappa_1 p \poi \Gamma\kappa_2 q \\
&= \Gamma(\kappa_1 p \poi \kappa_2 q) \\
&= \Gamma(\tl{f}\tr{g})
\end{align*}
(ii) Let $\tl{f} \tr{g}$ be the pullback of both $\tr{p_1}\tl{q_1}$ and $\tr{p_2}\tl{q_2}$. By (i) $\Gamma(\tl{f} \tr{g}) = \Delta(\tr{p_1}\tl{q_1})$ and $\Gamma(\tl{f} \tr{g}) = \Delta(\tr{p_2}\tl{q_2})$. The statement follows.

\noindent (iii) Let $n \tr{e^1_p}\tr{m^1_p} z  \tl{m^1_q}\tl{e^1_q} m$ and $n \tr{e^2_p} \tr{m^2_p} r \tl{m^2_q}\tl{e^2_q} m$ be epi-mono factorisations of ${n \tr{p_1} z \tl{q_1} m}$ and $n \tr{p_2} r \tl{q_2} m$ respectively, with the properties described in Definition \ref{def:equalzeros}. By assumption $e^1_p = e^2_p$, $e^1_q = e^2_q$. Also, $\tr{m^1_p} \tl{m^1_q}$ and $\tr{m^2_p} \tl{m^2_q}$ have the same pullback, call it $\tr{f} \tl{g}$. The statement is given by the following derivation.
\begin{align*}
\Delta(\tr{p_1}\tl{q_1}) &=  \Delta(\tr{e^1_p}\tr{m^1_p} \tl{m^1_q}\tl{e^1_q})  \\
&= \Delta(\tr{e^2_p}\tr{m^1_p} \tl{m^1_q}\tl{e^2_q}) \\
& = \Delta(\tr{e^2_p}\tl{\id})\poi \Delta(\tr{m^1_p}  \tl{m^1_q}) \poi \Delta(\tr{\id}\tl{e^2_q}) \\
&\hspace{-.2cm}\eql{(ii)}\hspace{-.2cm} \Delta(\tr{e^2_p}\tl{\id})\poi \Delta(\tr{m^2_p}  \tl{m^2_q}) \poi \Delta(\tr{\id}\tl{e^2_q}) \\
&= \Delta(\tr{e^2_p}\tr{m^2_p} \tl{m^2_q}\tl{e^2_q}) \\
&= \Delta(\tr{p_2}\tl{q_2}).
\end{align*}
(iv) The statement is given by the following derivation.
\begin{align*}
\Gamma(\tl{f}\tr{g}) &= \Gamma(\kappa_2 f \poi \kappa_1 g) \\
&= \Gamma(\kappa_2 f )\poi \Gamma(\kappa_1 g) \\
&= \Delta(\iota_2 f)\poi \Delta(\iota_1 g) \\
&= \Delta( \iota_2 f \poi \iota_1 g) \\
&= \Delta( \tr{p}\tl{q} ).~\qedhere
\end{align*}
\end{proof}

In particular, all properties of Lemma~\ref{lemma:arbitraryPROP_ER} hold for \eqref{eq:pushoutER}.

\begin{lemma}\label{lemma:cubeERcommutes} \eqref{eq:pushoutER} commutes.\end{lemma}
\begin{proof}
It suffices to show that \eqref{eq:pushoutER} commutes on the two injections into $\Injplus$, that means, for any
$f \: n\to m$ in $\Inj$,
\begin{equation*}
\SpanToER(\tl{\id}\tr{f}) = \CospanToER(\tr{f}\tl{\id}) \qquad \text{ and } \qquad \SpanToER(\tl{f}\tr{\id}) =
\CospanToER(\tr{\id}\tl{f}).\end{equation*}
These are clearly symmetric, so it is enough to check one:
\begin{align*}
\SpanToER(\tl{\id}\tr{f}) &=
\eqr{\{ (v,w) \mid v=\inv{f}(w) \}} \\
& = \eqr{\{ (v,w) \mid f(v)= w  \}} \\
& = \eqr{\{ (v,w) \mid f(v)= w \vee f(v) = f(w) \vee v = w \}} && (f \text{ is injective})\\
& = \CospanToER(\tr{f}\tl{\id}).~\qedhere
\end{align*}
\end{proof}

 Lemma~\ref{lemma:arbitraryPROP_ER} states that a commutative diagram~\eqref{eq:arbitraryER} equalizes all cospans that are equal modulo-zeros. We now verify that, for \eqref{eq:pushoutER}, also the converse statement holds.
\begin{lemma}
\label{lemma:pullbackpsi}
The following are equivalent
\begin{enumerate}[(a)]
\item \label{point:pullbackpsi1}$n \tr{p_1} z  \tl{q_1}m$ and $n\tr{p_2} r \tl{q_2}m$ are equal modulo zeros.
\item \label{point:pullbackpsi2}$\CospanToER(\tr{p_1}\tl{q_1})=
\CospanToER(\tr{p_2}\tl{q_2})$.
\end{enumerate}
\end{lemma}
\begin{proof}
The conclusions of
Lemmas~\ref{lemma:cubeERcommutes} and~\ref{lemma:arbitraryPROP_ER}
give that \ref{point:pullbackpsi1} $\Rightarrow$ \ref{point:pullbackpsi2}. It thus suffices to show that
\ref{point:pullbackpsi2} $\Rightarrow$ \ref{point:pullbackpsi1}. For this purpose, it is useful to first verify the following properties:
\begin{enumerate}[(i)]
\item for all $u,u' \in \ord{n}$, $p_1(u) = p_1(u')$ if and only if $p_2(u) = p_2(u')$
\item for all $v,v' \in \ord{m}$, $q_1(v) = q_1(v')$ if and only if $q_2(v) = q_2(v')$
\item for all $u \in \ord{n}$, $v \in \ord{m}$, $p_1(u) = q_1(v)$ if and only if $p_2(u) = q_2(v)$
\item Let $p_1[\ord{n}]$ be the number of elements of $\ord{n}$ that are in the image of $p_1$, and similarly for $p_2[\ord{n}]$. Then $p_1[\ord{n}] = p_2[\ord{n}]$.
\item $q_1[\ord{n}] = q_2[\ord{n}]$.
\end{enumerate}
For statement (i), observe that, by definition of $\CospanToER$ as in~\eqref{eq:CospanToER}, for any two elements $u, u' \in \ord{n}$ the pair $(u,u')$ is in $\CospanToER(\tr{p_1}\tl{q_1})$ if and only if $p_1(u) = p_1(u')$. Similarly, $(u,u') \in \CospanToER(\tr{p_2}\tl{q_2})$ if and only if $p_2(u) = p_2(u')$. Since by assumption $\CospanToER(\tr{p_1}\tl{q_1}) = \CospanToER(\tr{p_2}\tl{q_2})$, we obtain (i). A symmetric reasoning yields (ii).
The argument for statement (iii) is analogous: for $i \in \{1,2\}$ and $u \in \ord{n}$, $v \in \ord{m}$, by definition of $\CospanToER$, $(u,v) \in \CospanToER(\tr{p_i}\tl{q_i})$ if and only if $p_i(u) = q_i(v)$. Since $\CospanToER(\tr{p_1}\tl{q_1}) = \CospanToER(\tr{p_2}\tl{q_2})$, we obtain (iii).
Statement (iv) is an immediate consequence of (i), and (v) of (ii).

Now, by virtue of properties (i)-(v), it should be clear that we can define epi-mono factorisations $n \tr{e^1_p}\tr{m^1_p} z  \tl{m^1_q}\tl{e^1_q} m$ and $n \tr{e^2_p} \tr{m^2_p} r \tl{m^2_q}\tl{e^2_q} m$ of $n \tr{p_1} z \tl{q_1} m$ and $n \tr{p_2} r \tl{q_2} m$ respectively, with the following properties.
\begin{itemize}
    \item[(vi)] $e^1_p$ and $e^2_p$ are the same function, with source $n$ and target $p_1[\ord{n}] = p_2[\ord{n}]$. Also $e^1_q$ and $e^2_q$ are the same function, with source $m$ and target $q_1[\ord{m}] = q_2[\ord{m}]$.
    \item[(vii)] For all $u\in p_1[\ord{n}]=p_2[\ord{n}]$ and $v \in q_1[\ord{n}] = q_2[\ord{n}]$, $m^1_p(u) = m^1_q(v)$ iff $m^2_p(u) = m^2_q(v)$.
\end{itemize}
It remains to prove that $\tr{m^1_p}  \tl{m^1_q}$ and $\tr{m^2_p} \tl{m^2_q}$ have the same pullback. For this purpose, let the following be
pullback squares in $\Inj$:
\[
\xymatrix{
{h_1} \ar[d]_{g_1} \ar[r]^{f_1} & {q_1[\ord{n}]} \ar[d]^{m^1_q} \\
{p_1[\ord{n}]} \ar[r]_{m^1_p} & {z}
}
\qquad
\xymatrix{
{h_2} \ar[d]_{g_2} \ar[r]^{f_2} & {q_1[\ord{n}]} \ar[d]^{m^2_q} \\
{p_1[\ord{n}]} \ar[r]_{m^2_p} & {r}
}
\]
By the way pullbacks are computed in $\Inj$ (i.e., in $\F$), using (vii) we can conclude that $m^1_p g_2=m^1_q f_2$ and $m^2_p g_1=m^2_q f_1$. By universal property of pullbacks, this implies that the spans $\tl{g_1}\tr{f_1}$ and
$\tl{g_2}\tr{f_2}$ are isomorphic. 
 \end{proof}

It is also useful to record the following observation.
\begin{lemma}\label{lemma:psifull} $\CospanToER \: \F\bicomp{\Perm}\Fop \to \ER$ is full.\end{lemma}
\begin{proof} Let $c_1,\dots,c_k$ be the equivalence classes of an equivalence relation $e$ on $\ord{n}\uplus\ord{m}$. We define a cospan $n \tr{p} k \tl{q} m$ by letting $p$ map $v \in \ord{n}$ to the equivalence class $c_i$ to which $v$ belongs, and symmetrically for $q$ on values $w \in \ord{m}$. It is routine to check that $\CospanToER (\tr{p}\tl{q}) = e$.\end{proof}

We now have all the ingredients to conclude the characterisation of $\ER$.

\begin{proposition}
\label{lemma:pushoutER}
\eqref{eq:pushoutER} is a pushout in $\PROP$.
\end{proposition}
\begin{proof}
Commutativity of \eqref{eq:pushoutER} is given by Lemma~\ref{lemma:cubeERcommutes}, thus it remains to show the universal property. Suppose that we have a commutative diagram of PROP morphisms as in~\eqref{eq:arbitraryER}. It suffices to show that there exists a PROP morphism $\Theta \: \ER\to\mathbb{X}$ with $\Theta\SpanToER = \Gamma$ and $\Theta\CospanToER=\Delta$ -- uniqueness is automatic by fullness of $\CospanToER$ (Lemma~\ref{lemma:psifull}).

Given an equivalence relation $e \: n\to m$, there exist a cospan $\tr{p}\tl{q}$ such that
$\CospanToER(\tr{p}\tl{q})=e$.
We let $\Theta(e)=\Delta(\tr{p}\tl{q})$.
This is well-defined: if $\tr{p'}\tl{q'}$ is another
cospan such that $\CospanToER(\tr{p'}\tl{q'})=e$ then applying the conclusion of Lemma~\ref{lemma:pullbackpsi} gives us that $\tr{p}\tl{q}$ and $\tr{p'}\tl{q'}$ are equal modulo-zeros and thus, by Lemma~\ref{lemma:arbitraryPROP_ER}, $\Delta(\tr{p}\tl{q})=\Delta(\tr{p'}\tl{q'})$.
This  argument also shows that, generally, $\Theta\CospanToER=\Delta$.
Finally, $\Theta$ preserves composition:
\begin{align*}
\Theta(e\poi e') &= \Theta(\CospanToER(\tr{p}\tl{q})\poi\CospanToER(\tr{p'}\tl{q'})) \\
&= \Theta(\CospanToER((\tr{p}\tl{q})\poi (\tr{p'}\tl{q'}))) \\
&= \Delta((\tr{p}\tl{q})\poi (\tr{p'}\tl{q'})) \\
&= \Delta(\tr{p}\tl{q})\poi\Delta(\tr{p'}\tl{q'})\\
&= \Theta(e)\poi\Theta(e').
\end{align*}

We can now also show that
$\Theta\SpanToER = \Gamma$: given a span $\tl{f}\tr{g}$ in $\Inj$, let
$\tr{p}\tl{q}$ be its pushout span in $\F$.
Using the conclusions of Lemma~\ref{lemma:arbitraryPROP_ER}.(iv),
$\Gamma(\tl{f}\tr{g}) = \Delta(\tr{p}\tl{q})
=\Theta\CospanToER(\tr{p}\tl{q})=\Theta\SpanToER(\tl{f}\tr{g})$.\end{proof}

\subsubsection{The Cube for Equivalence Relations}

We are now ready to supply a proof of the main statement of this section. 

\begin{proof}[Proof of Theorem~\ref{th:IFROB=ER}]
We construct the following diagram in $\PROP$.
\begin{equation}\label{eq:cubeER}
\raise30pt\hbox{$
\xymatrix@=5pt{
& {\UNIT + \COUNIT} \ar[dd]_(.3){\cong}|{\hole}
\ar[dl] \ar[rr] & & {\COUNIT \bicomp{\Perm}\UNIT} \ar[dl] \ar[dd]^{\cong} \\
{\FROB} \ar[rr] \ar[dd]_{\cong}  & & {\IFROB} \ar@{.>}[dd] \\
& {\Inj + \Injop} \ar[dl] \ar[rr]|(.57){\hole} & & {\Injop \bicomp{\Perm} \Inj} \ar[dl] \\
{\F \bicomp{\Perm} \Fop} \ar[rr] & & {\ER}
}$}
\end{equation}
The top and the bottom face of the cube are \eqref{eq:pushoutIFROB} and \eqref{eq:pushoutER} respectively. The three non-dotted vertical morphisms are the PROP isomorphisms yielded by $\UNIT \cong \Inj$, $\COUNIT \cong \Injop$ and $\FROB \cong \F \bicomp{\Perm} \Fop$ (see Example~\ref{ex:distrlawsyntactic}\ref{ex:distrlawsyntactic1}-\ref{ex:distrlawsyntactic3}). The two rear faces of the cube commute by definition.
Since the top and the bottom face are pushouts, the dotted vertical morphism $\IFROB \to \ER$, which is given by universal property of the topmost pushout, is also an isomorphism.
\end{proof}

\begin{remark}\label{rmk:ERcomposedPROP} As hinted by the rightmost diagram in \eqref{eq:spidersboneepimono}, one can give an alternative characterisation of $\ER$ as the composite PROP $\Surj \bicomp{\Perm} \Injop \bicomp{\Perm} \Inj \bicomp{\Perm} \Surj$. This will rely on defining the appropriate distributive laws and combine them together using the results presented in \S~\ref{sec:iteratedDistrLaws}. It turns out that the equations necessary to present the laws yielding $\Surj \bicomp{\Perm} \Injop \bicomp{\Perm} \Inj \bicomp{\Perm} \Surj$ are precisely those of $\IFROB$. Then, showing that factorised arrows of $\Surj \bicomp{\Perm} \Injop \bicomp{\Perm} \Inj \bicomp{\Perm} \Surj$ are in bijective correspondence with equivalence relations in $\ER$ completes the proof that $\IFROB \cong \ER$.

We did not show this characterisation in detail for two reasons. First, we find obtaining the isomorphism $\IFROB \cong \ER$ by universal property of pushout more elegant and simple than passing through the construction ``by hand'' of an isomorphism $\Surj \bicomp{\Perm} \Injop \bicomp{\Perm} \Inj \bicomp{\Perm} \Surj \cong \ER$. Second, we wanted to offer an illustrative example of the cube construction~\eqref{eq:cubeER}, as it will be pivotal for the characterisation results of the next chapter.\end{remark}

\begin{remark}\label{rmk:pushoutTrivialTheory} Our construction merges the theory of cospans of functions and of spans of injective functions to form the theory of equivalence relations. One may wonder what happens with a more symmetric approach, namely if we consider spans of arbitrary functions. We saw in Example~\ref{ex:distrlawsyntactic}\ref{ex:distrlawsyntactic2} that the theory of spans in $\F$ is the PROP $\B$ of bialgebras. Thus mimicking the cube construction~\eqref{eq:cubeER} would result in the following diagram
\begin{equation}\label{eq:cubetrivial}
\raise30pt\hbox{$
\xymatrix@=5pt{
& {\Mon + \Com} \ar[dd]_(.3){\cong}|{\hole}
\ar[dl] \ar[rr] & & {\B} \ar[dl] \ar[dd]^{\cong} \\
{\FROB} \ar[rr] \ar[dd]_{\cong}  & & {\T} \ar@{.>}[dd] \\
& {\F + \Fop} \ar[dl] \ar[rr]|(.57){\hole} & & {\Fop \bicomp{\Perm} \F} \ar[dl] \\
{\F \bicomp{\Perm} \Fop} \ar[rr] & & {.}
}$}
\end{equation}
where the top and the bottom face are pushouts. By Proposition~\ref{prop:SMTpushout}, the SMT for $\T$ features equations \eqref{eq:wmonassoc}-\eqref{eq:bwbone}, \eqref{eq:BWFrob} and \eqref{eq:BWSep}. Thus the following holds in $\T$.
\begin{equation*}
\lower6pt\hbox{$\includegraphics[height=.6cm]{graffles/idcircuit.pdf}$}
\eql{\eqref{eq:wmonunitlaw},\eqref{eq:bcomonunitlaw}}
\lower7pt\hbox{$\includegraphics[height=.7cm]{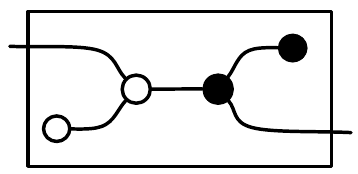}$}
\eql{\eqref{eq:BWFrob}}
\lower9pt\hbox{$\includegraphics[height=.9cm]{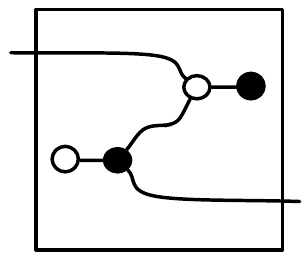}$}
\eql{\eqref{eq:unitsr}}
\lower9pt\hbox{$\includegraphics[height=.9cm]{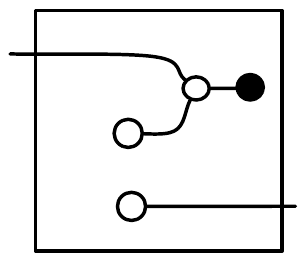}$}
\eql{\eqref{eq:unitsl}}
\lower9pt\hbox{$\includegraphics[height=.9cm]{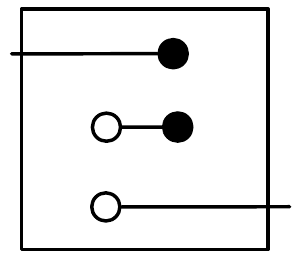}$}
\eql{\eqref{eq:bwbone}}
\lower6pt\hbox{$\includegraphics[height=.6cm]{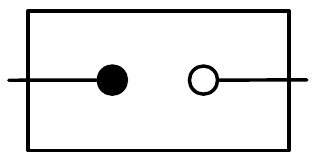}$}
\end{equation*}
The above equation trivializes $\T$, as it implies that any two arrows $n \tr{f \in \T} m$ and $n \tr{g \in \T} m$ are equal. This makes $\T$ (and thus also the pushout object of the bottom face in~\eqref{eq:cubetrivial}) the \emph{terminal object} in $\PROP$: for any other PROP $\PS$ there is a unique morphism which maps any arrow $n \tr{h \in \PS} m$ into the unique arrow with that source and target in $\T$.
\end{remark}

\subsection{Case Study II: Partial Equivalence Relations} \label{sec:PER}

Building on the results of the previous section, we shall now construct the SMT of \emph{partial equivalence relations} (PERs) via pushout of PROPs. Recall that a relation $e \subseteq X \times X$ is a PER if it is symmetric and transitive --- equivalently, $e$ is an equivalence relation on a subset $Y \subseteq X$. PERs are used extensively in the semantics of higher order $\lambda$-calculi (e.g., \cite{JacobsCLTT,Streicher_SemanticsTT}) and, more recently, of quantum computations (e.g., \cite{Jacobs_quantumPER,HasuoH11_GOIQuantumPER}).

\begin{definition}
Let $\gls{PROPPER}$ be the PROP with arrows $n \to m$ partial equivalence relations on $\ord{n}\uplus\ord{m}$. Given $e_1 \: n \to z$, $e_2 \: z \to m$, their composite $e_1 \poi e_2$ is defined by
\begin{equation} \label{eq:compPER}
\begin{aligned}
e_1 \relcomp e_2 \df & \{(v,w) \mid \exists u.\; (v,u)\in e_1 \wedge (u,w)\in e_2\} &\\
e_1 \poi e_2 \df & \{(u,w) \mid u,w \in \ord{n+m} \wedge (u,w) \in \peqr{e_1 \relcomp e_2} \} &\\
& \cup \{(u,u) \mid u \in \ord{n+m} \wedge (u,u) \in e_1 \cup e_2\} &
\end{aligned}
\end{equation}
where $\peqr{e_1 \relcomp e_2}$ denotes the PER generated by $e_1 \relcomp e_2$, that is, its symmetric and transitive closure.
The monoidal product $e_1 \tns e_2$ on partial equivalence relations $e_1, e_2$ is given by disjoint union.
\end{definition}
Observe that composition in $\PER$ is defined in the same way as in $\ER$, the only difference being that here we take the reflexive closure only for those elements of $\ord{n}\uplus\ord{m}$ on which either $e_1$ or $e_2$ are defined.

\subsubsection{Partial Frobenius Algebras}

As we did for equivalence relations, we establish some preliminary intuition on representing PERs with string diagrams. For functions (Example~\ref{ex:partialfunctions}), partiality was captured at the graphical level by incorporating an additional generator $\Ocounit$. The strategy for PERs is analogous: for the elements on which a PER $e$ is defined, the diagrammatic description will be the same given for equivalence relations in~\eqref{eq:spider}; the elements on which $e$ is undefined will instead correspond to ports where we plug in $\Ocounit$ (if on the left) or $\Ounit$ (if on the right).

Therefore, the string diagrammatic theory for PERs will involve $\FROB$ expanded with generators $\Ocounit$, $\Ounit$, subject to suitable compatibility conditions. To make this formal, we introduce the PROP of ``partial'' special Frobenius algebras. Its definition will rely on the PROP $\PMon$ of partial commutative monoids discussed in Example~\ref{ex:partialfunctions}.

\begin{definition} The PROP $\gls{PROPPFROB}$ is defined by letting arrows $n \to m$ be arrows of $\PMon + (\PMon)^{\op}$ quotiented by equations~\eqref{eq:BWFrob}, \eqref{eq:BWSep} and the following.
\begin{multicols}{2}\noindent
 \begin{equation}\tag{PFR1}
\label{eq:ocobone}
\lower7pt\hbox{$\includegraphics[height=.7cm]{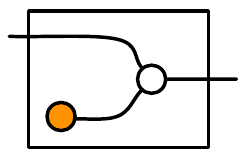}$}
=
\lower6pt\hbox{$\includegraphics[height=.6cm]{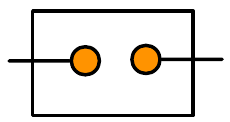}$}
 =
 \lower7pt\hbox{$\includegraphics[height=.7cm]{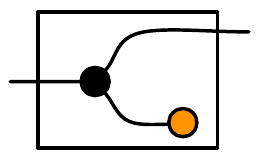}$}
\end{equation}
\begin{equation}\tag{PFR2}
\label{eq:obone}
\lower6pt\hbox{$\includegraphics[height=.6cm]{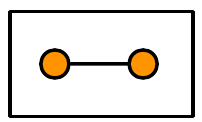}$}
=
\lower5pt\hbox{$\includegraphics[width=16pt]{graffles/idzerocircuit.pdf}$}
\end{equation}
 \end{multicols}
 \end{definition}

 \begin{example} It is instructive at this point to consider the following example~\eqref{eq:compPER1} of a string diagram $4 \to 2$ in $\PFROB$. Diagram~\eqref{eq:compPER2} is just the identity on $2$.
 \begin{multicols}{2}\noindent
 \begin{equation}
\label{eq:compPER1}
\lower9pt\hbox{$\includegraphics[height=.9cm]{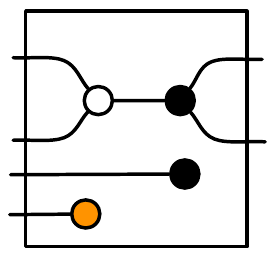}$}
\end{equation}
\begin{equation}
\label{eq:compPER2}
\lower9pt\hbox{$\includegraphics[height=.9cm]{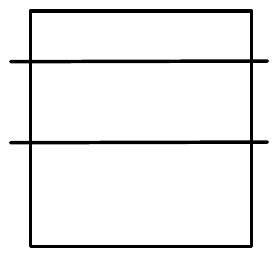}$}
\end{equation}
 \end{multicols}
\eqref{eq:compPER1} represents the PER $e$ on $\ord{4}\uplus{2}$ with equivalence classes $[1 \in \ord{4},2 \in \ord{4}, 1\in \ord{2}, 2\in \ord{2}], [3 \in \ord{4}]$ and undefined on $4 \in \ord{4}$. As~\eqref{eq:compPER2} is the identity, the composite $e \poi \id$ should be equal to~\eqref{eq:compPER1}. It becomes now more clear why the addition of pairs $(u,u)$ in \eqref{eq:compPER} is necessary to get the right notion of composition. In fact, $e \relcomp \id$ does not contain the equivalence class $[3 \in \ord{4}]$: intuitively, this is because we are only taking connectivity of ports into account, thus not distinguishing between $\Bcounit$ and $\Ocounit$ in~\eqref{eq:compPER1}. Including in $e \poi \id$ the set of pairs $\{(u,u) \mid u \in \ord{3}\uplus\ord{2} \wedge (u,u) \in e \cup \id\}$ allows to make this distinction and pick elements whose port is connected with $\Bcounit$.
 \end{example}

 As a partial version of $\FROB$, we expect $\PFROB$ to characterise cospans of partial functions. To phrase this statement, note that $\PF$ is equivalently described as the coslice category $1 / \F$ (that is, the skeletal category of pointed finite sets and functions) and thus has pushouts inherited from $\F$. We can then form the PROP $\PF \bicomp{\Perm} \PF^{\op}$ via a distributive law $\PF^{\op} \bicomp{\Perm}\PF \to \PF \bicomp{\Perm} \PF^{\op}$ defined by pushout, analogously to the case of functions (Example~\ref{ex:composingSemPROP}\ref{ex:composingSemPROP3}).
 \begin{proposition}\label{th:PFROBcomplete} $\PFROB \cong \PF \bicomp{\Perm} \PF^{\op}$. \end{proposition}
 For proving Proposition~\ref{th:PFROBcomplete}, we need to check that the equation associated with any pushout diagram in $\PF$ is provable by the equations of $\PFROB$. The key observation is that we can confine ourselves to just pushouts involving the generators of $\PMon$.

 Before making this formal, we illustrate the idea of the argument with the following example. The leftmost diagram below is a diagram representing a span $\tl{f}\tr{g}$ (left), which we transform into a cospan (right) pushing out $\tl{f}\tr{g}$, only using equations of $\PFROB$.
 \begin{equation*}
\lower9pt\hbox{$\includegraphics[height=1.1cm]{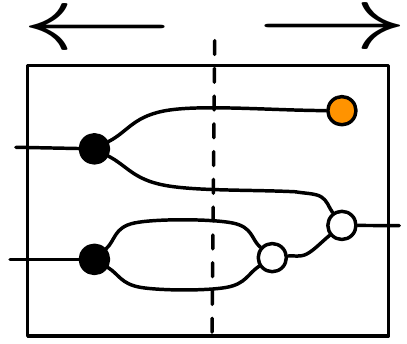}$} \Rightarrow
\lower9pt\hbox{$\includegraphics[height=1.1cm]{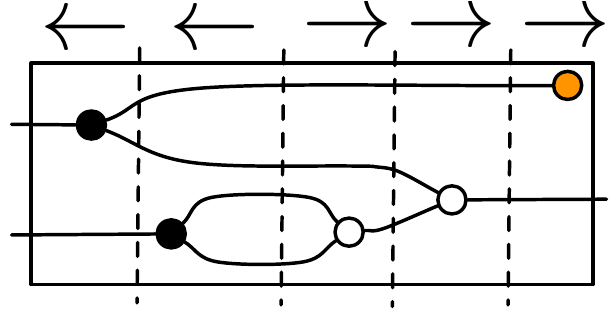}$} \Rightarrow
\lower9pt\hbox{$\includegraphics[height=1.1cm]{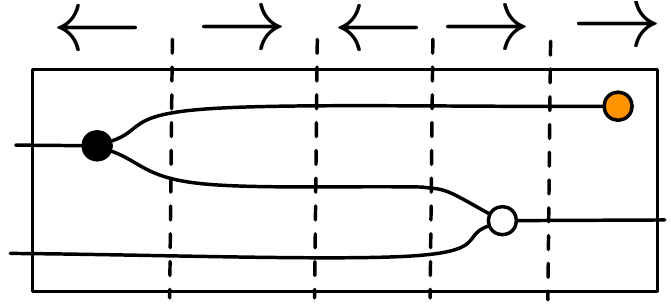}$} \Rightarrow
\lower9pt\hbox{$\includegraphics[height=1.1cm]{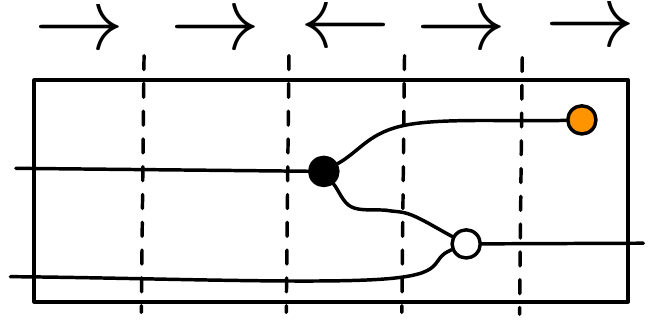}$} \Rightarrow
\lower9pt\hbox{$\includegraphics[height=1.1cm]{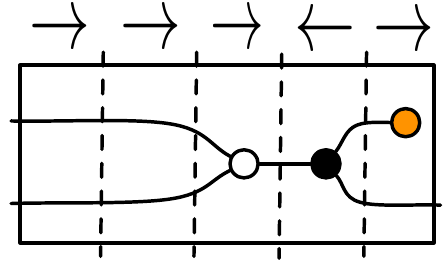}$} \Rightarrow
\lower9pt\hbox{$\includegraphics[height=1.1cm]{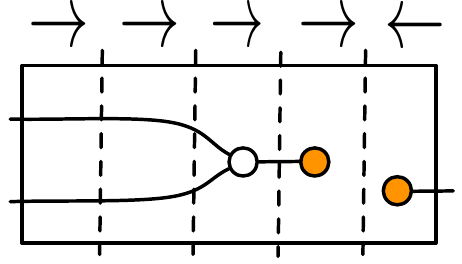}$}
\end{equation*}
  The steps are as follows. First, we expand $\tl{f}$ and $\tr{g}$ as $\tl{f_1}\tl{f_2}$ and $\tr{g_1}\tr{g_2}\tr{g_3}$ respectively, in such a way that each $f_i$ and $g_i$ contains at most one generator of $\PF$ and $\PF^{\op}$. In the following steps, we proceed pushing out spans $\tl{f_i}\tr{g_j}$ whenever possible: graphically, this amounts to apply valid equations of $\PFROB$ of a very simple kind, namely those describing the interaction of a single generator of $\PF^{\op}$ with one of $\PF$ (here we also consider instances of $\Idnet$ and $\symNet$ as generators). Note that pushing out spans of this form always gives back a cospan $\tr{p}\tl{q}$ with $p$, $q$ containing at most one generator, meaning that the procedure can be applied again until no more spans appear. The resulting diagram (the rightmost above) is the pushout of the leftmost one by pasting properties of pushouts. Therefore, we just proved that the equation
   \begin{equation*}
\lower9pt\hbox{$\includegraphics[height=1.1cm]{graffles/PFROB_ex1.pdf}$}
=\lower9pt\hbox{$\includegraphics[height=1.1cm]{graffles/PFROB_ex6.pdf}$}
\end{equation*}
 arising by the distributive law $\PF^{\op} \bicomp{\Perm} \PF \to \PF \bicomp{\Perm} \PF^{\op}$ is provable in $\PFROB$.

 We now formalise the argument sketched above. Let us call \emph{atom} any diagram of $\PMon$ of shape $\tr{f}\tr{b}\tr{g}$, where $f$ and $g$ consist of components $\Idnet$ and $\symNet$ composed together via $\tns$ or $\poi$, and $b$ is either $\id_0$ or a generator of $\PMon$. The following lemma establishes that $\PFROB$ is complete for pushouts involving atoms.
 \begin{lemma}\label{lemma:atomicCircuits} Let $\tl{f} \tr{g}$ be a span in $\PF$ where $f$ and $g$ are in the image (under the isomorphism $\PMon \cong \PF$) of atoms and suppose that the following is a pushout square.
  \begin{equation}\label{eq:atomicCircuitspushout}
   \xymatrix@=10pt{ & \ar[dl]_{f} r \ar[dr]^{g} & \\
 \ar[dr]_{p} m & & \ar[dl]^{q} n \\
 & z & }
 \end{equation}
 Then:
\begin{enumerate}
  \item $p$ and $q$ are also in the image of atoms;
   \item the associated equation is provable in $\PFROB$.
  \end{enumerate}
 \end{lemma}
 \begin{proof}
 The two points are proved by case analysis on all the possible choices of generators of $\PMon$ and $(\PMon)^{\op}$.
 \end{proof}
 We can now provide a proof of Proposition~\ref{th:PFROBcomplete}.

 \begin{proof}[Proof of Proposition~\ref{th:PFROBcomplete}] Let the following be a pushout square in $\PF$.
 \begin{equation}\label{eq:poPFROBcompl}
  \vcenter{\xymatrix@=10pt{ & \ar[dl]_{f} r \ar[dr]^{g} & \\
 \ar[dr]_{p} m & & \ar[dl]^{q} n \\
 & z & }}
 \end{equation}
Pick expansions $f = f_1 \poi \dots \poi f_k$ and $g = g_1 \poi \dots \poi g_j$, with each $f_i$ and $g_i$ in the image of an atom. We can calculate the pushout above by tiling pushouts of atoms as follows:
  \begin{equation}\label{eq:poatomicexpansion}
  \vcenter{ \xymatrix@R=2pt{
&&&\ar[dl]_{f_1} z \ar[dr]^{g_1} &&& \\
&& \ar[dl]_{f_2} \ar[dr] && \ar[dl] \ar[dr]^{g_2} &&\\
& \dots\ar[dl]_{f_k} \ar[dr] &&\ar[dr] \ar[dl] &&\dots \ar[dl] \ar[dr]^{g_j} &\\
\ar[dr] &&&\dots &&& \ar[dl] \\
&\ar[dr] &&\ar[dl]\dots \ar[dr]&& \ar[dl] &\\
&&\dots\ar[dr] && \dots\ar[dl] &&\\
&&&&&&
}
}
 \end{equation}
 Point 1 of Lemma~\ref{lemma:atomicCircuits} guarantees that each inner square only involves arrows in the image of some atom and Point 2 ensures that all the associated equations are provable in $\PFROB$. It follows that also the equation associated with the outer pushout~\eqref{eq:poatomicexpansion} is provable.

In order to derive from this that also the equation associated with~\eqref{eq:poPFROBcompl} is provable, it suffices to show that the pushout cospan in~\eqref{eq:poatomicexpansion} (call it $\tr{p'}\tl{q'}$) and the pushout cospan $\tr{p}\tl{q}$ in~\eqref{eq:poPFROBcompl} are provably equal once interpreted as string diagrams of $\PFROB$.
 Because $\tr{p}\tl{q}$ and $\tr{p'}\tl{q'}$ push out the same span, there is an isomorphism $\tr{i\in \PF}$ (in fact, a permutation) such that
\[\vcenter{
            \xymatrix@R=13pt@C=23pt{
                & & \\
               \ar[ur]^{p} \ar[r]_{p'}&\ar[u]_-{i} & \ar[ul]_{q}\ar[l]^{q'}
            }
}\]
commutes. It then follows that in $\PF + \PF^{\op}$
\[ \tr{p}\tl{q} \ =\  \tr{p'}\tr{i}\tl{q}  \ =\  \tr{p'}\tr{i}\tl{i}\tl{q'} \ = \ \tr{p'}\tr{i}\tr{i^{-1}}\tl{q'} \ = \ \tr{p'}\tl{q'}.\]
Because of the isomorphism $\PMon \cong \PF$, string diagrams associated with $\tr{p'}\tl{q'}$ and $\tr{p}\tl{q}$ are provably equal in $\PMon + (\PMon)^{\op}$ and thus also in its quotient $\PFROB$.
  \end{proof}

\begin{remark}\label{rmk:pullbacksnoatom} It is worth mentioning that the completeness argument for $\PFROB$ can be used as well to show completeness of the equations of $\FROB$ with respect to $\F \bicomp{\Perm} \Fop$ (Ex.~\ref{ex:distrlawsyntactic}\ref{ex:distrlawsyntactic3}). However, an analogous way of reasoning fails for showing the dual characterisation $\B \cong \Fop \bicomp{\Perm} \F$ (Ex.~\ref{ex:distrlawsyntactic}\ref{ex:distrlawsyntactic2}). The problem lies in proving point 1 of Lemma~\ref{lemma:atomicCircuits}, with~\eqref{eq:atomicCircuitspushout} now a pullback square. A counterexample is the pullback associated with \eqref{eq:bialg}: that equations shows a cospan of atoms (the left-hand side) whose pullback span (the right-hand side) does not consist of atoms.
\end{remark}

Now that we have an understanding of Frobenius algebras with partiality, the last ingredient that we need to model PERs is the axiom~\eqref{eq:bwbone}. Let us call $\gls{PROPIPFROB}$ (\textbf{i}rredundant \textbf{p}artial \textbf{Fr}obenius algebras) the PROP obtained by quotienting $\PFROB$ by~\eqref{eq:bwbone}. The rest of this section will be devoted to proving that the equations of $\IPFROB$ present the PROP of PERs.

\begin{theorem}\label{th:PIFROB=PER} $\IPFROB \cong \PER$.\end{theorem}

Our argument proceeds analogously to the case of equivalence relations discussed in \S~\ref{sec:ER}. The next two subsections describe $\IPFROB$ and $\PER$ as the fibered sum of isomorphic spans in $\PROP$.

\subsubsection{$\IPFROB$ as a Fibered Sum}

In showing Theorem~\ref{th:PIFROB=PER} it will be useful to modularly reconstruct $\IPFROB$, in the following way.

\begin{proposition} The following diagram is a pushout in $\PROP$, with PROP morphisms the standard interpretations of a diagram of the smaller theory as one of the larger theory.
\begin{equation}\label{eq:topPER}
\vcenter{
\xymatrix@C=15pt{
& \ar[dl]_{} \UNIT + \COUNIT \ \ar[rr]^{} && \ar[dl] \ \FROB \ \ar[rr] && \ar[dl] \ \PFROB \ \\
\COUNIT \bicomp{\Perm}\UNIT \ar[rr] && \IFROB \ar[rr] && \IPFROB &
}
}
\end{equation}
\end{proposition}
\begin{proof}
The two inner diagrams are pushouts by Proposition~\ref{prop:SMTpushout}, thus diagram~\eqref{eq:topPER} is also a pushout.
\end{proof}

\subsubsection{$\PER$ as a Fibered Sum}

We now turn to the construction of a pushout diagram for $\PER$. It will mimick construction~\eqref{eq:topPER} using the ``semantic'' characterisations of PROPs $\UNIT$, $\COUNIT$, $\FROB$, $\IFROB$ and $\PFROB$ developed in the previous sections.
\begin{equation}\label{eq:bottomPER}
\vcenter{
\xymatrix@=15pt{
& \ar[dl]_{[\kappa_1,\kappa_2]} \Inj + \Injop \ar[rr]^{[\iota_1,\iota_2]} && \ar[dl]_{\CospanToER} \F\bicomp{\Perm}\Fop \ar[rr]^{\CospanToPCospan} && \ar[dl]_{\CospanToER'} \PF\bicomp{\Perm}\PF^{\op}\\
\Injop \bicomp{\Perm} \Inj \ar[rr]_{\SpanToER} && \ER \ar[rr]_{\ERToPER} && \PER &
}
}
\end{equation}
PROP morphisms in the leftmost inner square of~\ref{eq:bottomPER} have been defined in \S~\ref{sec:ER}. We need to specify $\CospanToPCospan$, $\ERToPER$ and $\CospanToER'$:
\begin{itemize}
    \item For $\CospanToPCospan$, we recall the existence of an adjunction between $\F$ and $\PF$:
     \[
\xymatrix@R=5pt{
\F \ar@(ur,ul)[rr]^{F} &\bot & \ar@(dl,dr)[ll]^{U} \PF.
}
 \]
$F \: \F \to \PF$ is the obvious embedding of functions into partial functions. $U \: \PF \to \F$ maps $\ord{n}$ to $\ord{n+1}$ and $f \: \ord{n} \to \ord{m}$ to the function $\ord{n+1} \to \ord{m+1}$ acting as the identity on $\star \in \ord{1}$, and mapping $u \in \ord{n}$ into $f(u)$, if defined, and into $\star \in \ord{1}$ otherwise. Now, $\CospanToPCospan$ is defined as the embedding of $\F \bicomp{\Perm} \Fop$ into $\PF \bicomp{\Perm} \PF^{\op}$ induced by the functor $F \: \F \to \PF$. Note that this assignment is indeed functorial because left adjoints preserve colimits and composition of cospans is by pushout.
\item Similarly to the case of $\F$ and $\PF$, we let $\ERToPER$ be the obvious embedding of $\ER$ into $\PER$. This assignment is functorial because composition in $\PER$ behaves as composition in $\ER$ on partial equivalence relations that are totally defined.
\item The PROP morphism $\CospanToER' \: \PF\bicomp{\Perm}\PF^{\op} \to \PER$ is the extension of $\CospanToER \: \F\bicomp{\Perm}\Fop \to \ER$ to partial functions, defined by the same clause~\eqref{eq:CospanToER}. Note that the generality of $\PER$ is necessary: the value $e$ of $\CospanToER'$ on a cospan $\tr{p}\tl{q}$ in $\PF$ is possibly not a reflexive relation, since $p$ and $q$ may be undefined on some elements of $\ord{n}$ and $\ord{m}$. One can check that $\CospanToER'$ is a functorial assignment: the argument is essentially the same as for $\CospanToER$, because composition in $\PF\bicomp{\Perm}\PF^{\op}$ and $\PER$ is defined analogously to the one in $\F\bicomp{\Perm}\Fop$ and $\ER$.
\end{itemize}

We now show the following result.

\begin{proposition}\label{prop:pushoutPER} \eqref{eq:bottomPER} is a pushout in $\PROP$. \end{proposition}

For this purpose, it suffices to prove that the rightmost square (which clearly commutes by definition of $\CospanToER$, $\CospanToER'$, $\ERToPER$ and $\CospanToPCospan$) is a pushout. The following is the key lemma.

\begin{lemma}\label{lemma:PERpushout} Let $e \in \ER[n,m]$ and $\tr{p}\tl{q} \in \PF\bicomp{\Perm}\PF^{\op}$. The following are equivalent.
\begin{itemize}
\item[(i)] $\ERToPER(e) = \CospanToER'(\tr{p}\tl{q})$.
\item[(ii)] There are cospans $\tr{p_1}\tl{q_1}$, $\dots$, $\tr{p_k}\tl{q_k}$ in $\F\bicomp{\Perm}\Fop[n,m]$ such that 
    \begin{align*}
    e &= \CospanToER(\tr{p_1}\tl{q_1}) \\
    \CospanToPCospan(\tr{p_1}\tl{q_1}) &= \CospanToPCospan(\tr{p_2}\tl{q_2}) \\
    \CospanToER(\tr{p_2}\tl{q_2}) &= \CospanToER(\tr{p_3}\tl{q_3})\\
    \dots & \dots \\
    \CospanToPCospan(\tr{p_k}\tl{q_k}) &= \tr{p}\tl{q}.
    \end{align*}
\end{itemize}
\end{lemma}
\begin{proof}
First we observe that, because $\CospanToPCospan$ is an embedding, $\CospanToPCospan(\tr{p_i}\tl{q_i}) = \CospanToPCospan(\tr{p_{i+1}}\tl{q_{i+1}})$ implies $\tr{p_i}\tl{q_i} = \tr{p_{i+1}}\tl{q_{i+1}}$. It follows that (ii) is equivalent to the statement that (iii) there exist $\tr{p'}\tl{q'} \in \F\bicomp{\Perm}\Fop[n,m]$ such that $e = \CospanToER(\tr{p'}\tl{q'})$ and $\CospanToPCospan(\tr{p'}\tl{q'}) = \tr{p}\tl{q}$.

It is very easy to show that (iii) implies (i):
\begin{align*}
\ERToPER(e) & = \ERToPER \CospanToER(\tr{p'}\tl{q'})&& \text{(iii)} \\
 & = \CospanToER'\CospanToPCospan(\tr{p'}\tl{q'}) && \text{commutativity of \eqref{eq:bottomPER}} \\
 & = \CospanToER'(\tr{p}\tl{q}) && \text{(iii)}
\end{align*}

For the converse direction, suppose that we can show (*) the existence of $\tr{p'}\tl{q'} \in \F\bicomp{\Perm}\Fop[n,m]$ such that $\tr{p'}\tl{q'} = \CospanToPCospan(\tr{p}\tl{q})$. Then the following derivation gives statement (iii):
\begin{align*}
\ERToPER(e) & = \CospanToER'(\tr{p}\tl{q}) && \text{(i)} \\
 & = \CospanToER'\CospanToPCospan(\tr{p'}\tl{q'}) && \text{(*)} \\
 & = \ERToPER\CospanToER(\tr{p'}\tl{q'}) && \text{commutativity of \eqref{eq:bottomPER}}
\end{align*}
Indeed, because $\ERToPER$ is an embedding, the derivation above implies that $e = \CospanToER(\tr{p'}\tl{q'})$. Therefore it suffices to show (*). For this purpose, we just need to prove that both $n\tr{p \in \PF}z$ and $m\tr{q \in \PF}z$ are \emph{total} functions. Let $u$ be an element of $n$: since $\CospanToER'(\tr{p}\tl{q}) = \ERToPER(e)$ and $\ERToPER$ embeds equivalence relations into PERs, then $\CospanToER'(\tr{p}\tl{q})$ is in fact an equivalence relation, meaning that $u$ belongs to some equivalence class of the partition induced by $\CospanToER'(\tr{p}\tl{q})$. It follows by definition of $\CospanToER'$ that $p \: n \to z$ is defined on $u$. With a similar argument, one can show that $q \: m \to z$ is defined on all elements of $m$ and thus both $p$ and $q$ are total functions. This implies that $\tr{p}\tl{q}$ is in the image of the embedding $\CospanToPCospan$.
\end{proof}

\begin{proof}[Proof of Proposition~\ref{prop:pushoutPER}] The leftmost square in~\eqref{eq:bottomPER} has already been proven to be a pushout (Lemma~\ref{lemma:pushoutER}), thus it suffices to show that the rightmost square is also a pushout. With this aim, recall that pushouts in $\PROP$ can be calculated as in $\Cat$. In particular, \eqref{eq:bottomPER} involves categories all with the same objects and identity-on-objects functors. This means that the pushout object is the quotient of $\ER$ and $\F\bicomp{\Perm}\Fop$ along the equivalence relation generated by
\begin{equation}\label{eq:quotientPushoutPER} \{(e,\tr{p}\tl{q}) \mid \text{there is} \tr{p'}\tl{q'} \text{ such that } \CospanToER(\tr{p'}\tl{q'}) = e \text{ and } \CospanToPCospan(\tr{p'}\tl{q'}) =\ \tr{p}\tl{q} \}.
\end{equation}
 Lemma \ref{lemma:PERpushout} proves that $\CospanToER'$ and $\ERToPER$ map $n\tr{e \in \ER}m$ and $n\tr{p}\tl{q}m$ to the same arrow exactly when they are in the equivalence relation described above. This means that $\PER$ indeed quotients by~\eqref{eq:quotientPushoutPER} and thus is the desired pushout object.\end{proof}

\subsubsection{The Cube for Partial Equivalence Relations}

We are now able to conclude the characterisation of partial equivalence relations.

\begin{proof}[Proof of Theorem~\ref{th:PIFROB=PER}] We construct the following diagram in $\PROP$.
\begin{equation*}
\xymatrix@R=10pt{
& \ar[dl]_{} \UNIT + \COUNIT\ar[rr]^{} \ar[dd]_(.3){\cong}|{\hole} && \ar[dl] \FROB \ar[rr] \ar[dd]_(.3){\cong}|{\hole} && \ar[dl] \PFROB \ar[dd]^{\cong}\\
\COUNIT \bicomp{\Perm}\UNIT \ar[rr]  \ar[dd]_{\cong} && \IFROB \ar[rr] \ar[dd]_(.3){\cong} && \IPFROB \ar@{.>}[dd] & \\
& \ar[dl] \Inj + \Injop \ar[rr]|(.49){\hole} && \ar[dl] \F\bicomp{\Perm}\Fop \ar[rr]|(.48){\hole} && \ar[dl] \PF\bicomp{\Perm}\PF^{\op}\\
\Injop \bicomp{\Perm} \Inj \ar[rr] && \ER \ar[rr] && \PER &
}
\end{equation*}
The top and bottom face are pushout diagrams~\eqref{eq:topPER} and \eqref{eq:bottomPER} respectively. The vertical edges are PROP isomorphisms. Also, all the rear faces commute by definition. Therefore, the arrow $\IPFROB \to \PER$ given by universal property of the topmost pushout is also a PROP isomorphism.
\end{proof}

\begin{remark}Differently from the case of equivalence relations (Remark~\ref{rmk:ERcomposedPROP}), it is not clear how a modular story of partial equivalence relations could be reconstructed only using PROP composition. Intuitively, to present $\IPFROB$ as a composed PROP, we would need in particular equations to push components $\Ounit$ to the extreme left of any diagram. This should also involve a distributive law $\lambda \: \UNIT\bicomp{\Perm} \MULT \to \MULT \bicomp{\Perm}\UNIT$. Consider now the value of $\lambda$ on \lower6pt\hbox{$\includegraphics[height=.6cm]{graffles/OunitWmult.pdf}$}: because $\Ounit$ now indicates partiality, it should \emph{not} be defined as $\Idnet$ --- following~\eqref{eq:wmonunitlaw} --- but rather as \lower5pt\hbox{$\includegraphics[height=.5cm]{graffles/Ocobone.pdf}$} --- following~\eqref{eq:ocobone}. However, \lower5pt\hbox{$\includegraphics[height=.5cm]{graffles/Ocobone.pdf}$} is not of the correct type $\UNIT \bicomp{\Perm} \MULT$.
\end{remark}   
\chapter{Interacting Hopf Algebras}\label{chapter:hopf}

\section{Overview}\label{sec:overview}

\begin{figure}
\hspace{-.5cm}\includegraphics[width=16cm]{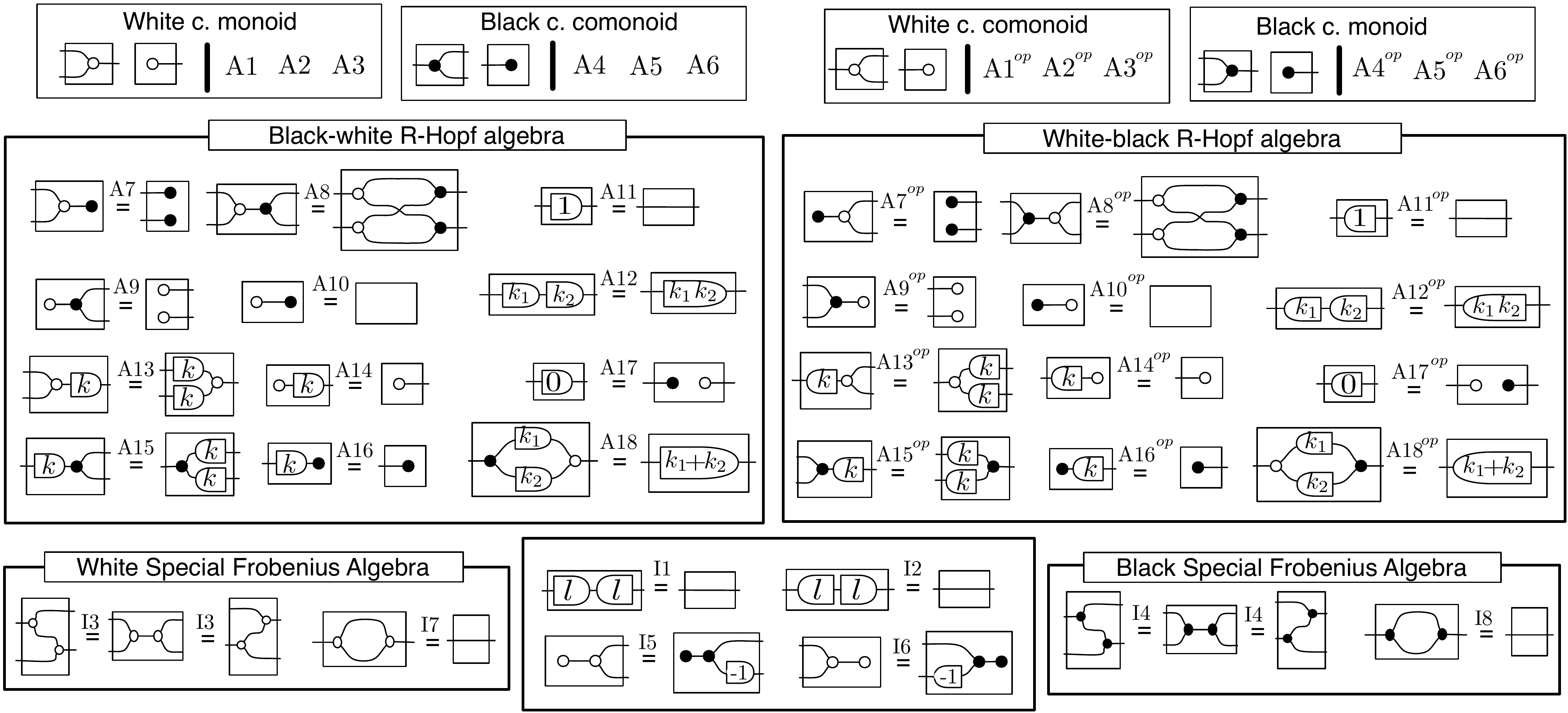}
 \caption{Equational theory of $\IBR$.}\label{fig:INTRO-axiomsIBR}
\end{figure}

The equational theory of bialgebras --- see Ex.~\ref{ex:equationalprops} --- can be parametrised over any commutative semiring $\PID$ \footnote{Recall that a semiring, sometimes called ``rig'', is a ring without the requirement that additive inverses exist.}. That amounts to adding to the bialgebra equations laws expressing the semiring structure of $\PID$ and how the bialgebra distributes over it: the resulting axiomatisation consists of equations~\eqref{eq:wmonassoc}-\eqref{eq:scalarsum} displayed in Figure~\ref{fig:INTRO-axiomsIBR}, where elements $k$, $k_1$ and $k_2$ range over $\PID$.
An easy observation is that, whenever $\PID$ is a ring, this set of equations actually defines an \emph{Hopf algebra}, with $-1 \in \PID$ playing the role of the antipode. The following is the derivation of the Hopf law.
\begin{align}\label{eq:HopfLawDer}
\lower8pt\hbox{$\includegraphics[height=.8cm]{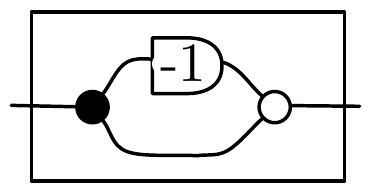}$}
=
\lower8pt\hbox{$\includegraphics[height=.8cm]{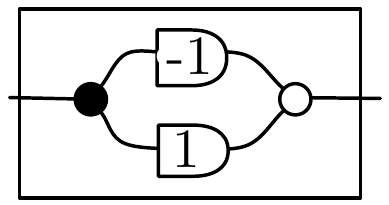}$}
=
\lower5pt\hbox{$\includegraphics[height=.5cm]{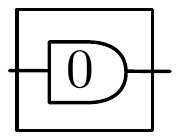}$}
=
\lower5pt\hbox{$\includegraphics[height=.5cm]{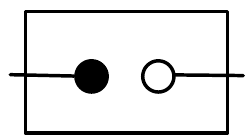}$}
=
\lower5pt\hbox{$\includegraphics[height=.5cm]{graffles/Hopf3.pdf}$}
=
\lower8pt\hbox{$\includegraphics[height=.8cm]{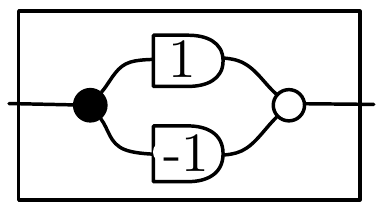}$}
=
\lower8pt\hbox{$\includegraphics[height=.8cm]{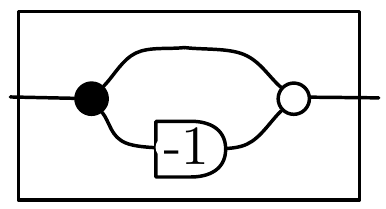}$}
\end{align}
We call the theory of \emph{$\PID$-Hopf algebras} the PROP $\ABR$ freely obtained from generators $$\Wmult \quad \Wunit  \quad \Bcomult  \quad \Bcounit  \quad \scalar \qquad \qquad k \in \PID$$ and equations~\eqref{eq:wmonassoc}-\eqref{eq:scalarsum}. An important property of $\PID$-Hopf algebras is that they characterise $\PID$-matrices. That means, there is an isomorphism of PROPs
\begin{equation}\label{eq:IntroABR=Vect} \ABR \cong \VectR \end{equation}
between $\ABR$ and the PROP $\VectR$ whose $n \to m$ arrows are $m\times n$ matrices over $\PID$. Because of this characterisation, $\PID$-Hopf algebras are a ubiquitous presence in computer science and the subject of many recent works:
 \begin{itemize}[noitemsep,topsep=0pt,parsep=0pt,partopsep=0pt]
  \item $\HA{\N}$ is isomorphic to the PROP $\B$ of bialgebras (Example~\ref{ex:equationalprops}). The isomorphism maps $\zeroscalar$ to $\spaceXaxisTiny$, $\unitscalar$ to $\idcircuit$  and, for $n \geq 2$, $\circuitn$ to the diagram $c_n \poi m_n$, where $c_n \: 1 \to n$ consists of $n-1$ copies of $\Bcomult$ and $m_n \: n \to 1$ consists of $m-1$ copies of $\Wmult$, arranged in the evident way. $\N$-matrices are the same as spans in $\F$ and \eqref{eq:IntroABR=Vect} is the result, appearing in many recent works~\cite{HylandPower_sketches,Lafont2003,Pirashvili_bialgebras,Lack2004a,Mimram_structureFirstOrderCausality} and reported in Example~\ref{ex:distrlawsyntactic}\ref{ex:distrlawsyntactic2}, that bialgebras characterise the PROP of spans.
 \item Spans, characterised by $\HA{\N}$, can be also described as multirelations. One obtains relations by considering matrices over the boolean semiring $\Bool$, that is the one given by the lattice on $\{0,1\}$ where one interprets $\vee$ as multiplication. Similarly to the case of multirelations, also the theory $\HA{\Bool}$ for relations has a finite presentation, being axiomatisable as $\B$ quotiented by \eqref{eq:BWSep}, which is the instance of \eqref{eq:scalarsum} expressing the law $1+1=1$.
 \item Another interesting case --- e.g. relevant for categorical quantum mechanics~\cite{CoeckeDuncanZX2011,BialgAreFrob14} --- is the one of integers modulo $2$. The PROP $\HA{\Z_{\scriptscriptstyle 2}}$ of $\Z_{\scriptscriptstyle 2}$-Hopf algebras is isomorphic to $\B$ quotiented by the following equation, which is the istance of \eqref{eq:scalarsum} expressing $1+1=0$:
     \begin{equation*}
\lower6pt\hbox{$\includegraphics[height=.7cm]{graffles/BWSep.pdf}$} =  \lower4pt\hbox{$\includegraphics[height=15pt]{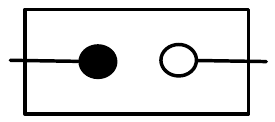}$}.
\end{equation*}
 \item As we shall prove in \S\ref{sec:instances}, also the PROP $\HA{\Z}$ of $\Z$-Hopf algebras enjoys a finite presentation: the only required generator of the kind $\scalar$ is the antipode $\scalarminusone$. In fact, $\Z$-Hopf algebras coincide with ordinary Hopf algebras as they are commonly defined in the literature --- see e.g.~\cite{Fiore2013,HopfRef1,HopfRef2}. Statement \eqref{eq:IntroABR=Vect} yields the folklore result that Hopf algebras characterise integer matrices.
 \item For $\PID$ any field, $\PID$-Hopf algebras are studied in rewriting theory, see e.g.~\cite{Rannou_phdthesis,Lafont2003}.
 Actually, to the best of our knowledge, \cite{Lafont2003} is the earliest reference for \eqref{eq:IntroABR=Vect}, since the argument for a field is the same as for a commutative semiring.
  \item The general case of a commutative semiring was recently considered in \cite{wadsley2015props}, where the authors prove \eqref{eq:IntroABR=Vect} with a different approach than \cite{Lafont2003}. We will show \eqref{eq:IntroABR=Vect} with yet a different argument, relying on the techniques developed in Chapter~\ref{sec:background}. Our construction will reveal the provenance of the axioms of $\PID$-Hopf algebras from distributive laws of PROPs.
\end{itemize}

\medskip
Now let $\PID$ be a principal ideal domain (PID). The main contribution of this chapter is the study of the interaction of $\ABR$ and its opposite PROP $\ABRop$. Our ultimate goal is to show the following characterisation result
\begin{equation}\label{eq:IntroIBR=SVR} \IBR \cong \SVR .\end{equation}
Here, $\SVR$ is the PROP of linear subspaces over the \emph{field of fractions} $\frPID$ of $\PID$ (which exists because $\PID$ is a PID): an arrow $n \to m$ of $\SVR$ is a subspace of $\frPID^{n}\times \frPID^{m}$ and composition is relational. $\IBR$ is the PROP of the theory of \emph{interacting $\PID$-Hopf algebras}: it is freely constructed by the generators of $\ABR + \ABR^{\op}$
\[ \Wmult \qquad \Wunit \qquad \Bcomult \qquad \Bcounit \qquad \scalar \qquad \Wcomult \qquad \Wcounit \qquad \Bmult \qquad \Bunit \qquad \coscalar \qquad \qquad k \in \PID \]
and the equations in Figure~\ref{fig:INTRO-axiomsIBR}, with $k$ ranging over $\PID$ and $l$ over $\PID\setminus \{0\}$. It may be useful to anticipate what meaning the isomorphism in~\eqref{eq:IntroIBR=SVR} assigns to the generators. $\Wmult$ is interpreted as \emph{sum}, i.e. it maps to the subspace $2 \to 1$ of pairs $({\scriptsize \left(\begin{array}{c} \!\!\! k \!\!\! \\ \!\!\! l \!\!\! \end{array}\right)},(k+l))$. $\Bcomult$ is interpreted as \emph{copy}, i.e. it maps to the subspace $1 \to 2$ of pairs $((k),{\scriptsize \left(\begin{array}{c} \!\!\! k \!\!\! \\ \!\!\! k \!\!\! \end{array}\right)})$. The subspaces for $\Wunit$ and $\Bcounit$ are $\{(0)\}$ and $\{(k) \mid k \in \PID\}$: they act as the neutral elements for sum and copy respectively. $\scalar$ is interpreted as multiplication by $k$, i.e. it maps to the subspace $1 \to 1$ of pairs $((l),(kl))$. Finally, the meaning of the last five generators above, from $\ABRop$, is the converse of the relation assigned to their counterpart in $\ABR$, e.g. $\coscalar$ is mapped to the subspace of pairs $((kl),(l))$.

Let us now focus on the axiomatisation of $\IBR$, given in Figure~\ref{fig:INTRO-axiomsIBR}. It has a very symmetric character: in a nutshell, monoid-comonoid pairs of different colour interact according to the laws of Hopf $\PID$-algebras, while pairs of the same colour yield \emph{special Frobenius algebras} --- \emph{cf.} Example~\ref{ex:equationalprops}. The remaining axioms of $\IBR$ state that each diagram $\scalar$ has $\coscalar$ as a formal inverse and that there are ``cup'' and ``cap'' structures, expressible either with white or black components.

The proof of \eqref{eq:IntroIBR=SVR} relies on a \emph{modular} construction of $\IBR$ based on the techniques introduced in Chapter~\ref{sec:background}. The key ingredient are two distributive laws between $\ABR$ and $\ABRop$, defined respectively by pullback and by pushout. Recall that in Chapter~\ref{sec:background} we considered the same kind of interaction between $\Mon$ and $\Mon^{\op} \cong \Com$: in that case, we used $\F$ as the ambient category where to compute the pullbacks and pushouts defining the distributive laws --- see Example~\ref{ex:distrlawsyntactic}\ref{ex:distrlawsyntactic2}-\ref{ex:distrlawsyntactic3}.
 For $\ABR$, we can take advantage of~\eqref{eq:IntroABR=Vect} and compute them in $\VectR$: the existence of pullbacks and pushouts is guaranteed by the assumption that $\PID$ is a principal ideal domain. The desired distributive laws will then have the following type, where $\PJ$ is the core of $\VectR$ (\emph{cf.}~\S~\ref{sec:distrLawPullback}):
\begin{eqnarray*}
\lambdapb \: \VectR \bicomp{\PJ} \VectRop \to \VectRop \bicomp{\PJ} \VectR & \qquad \qquad &
\lambdapo \: \VectRop \bicomp{\PJ} \VectR \to \VectR \bicomp{\PJ} \VectRop.
\end{eqnarray*}
As usual, the challenge is to provide a \emph{complete} set of axioms presenting the equations generated by $\lambdapb$ and by $\lambdapo$. In the case of distributive laws computed in $\F$, the complete axiomatisations were given by the PROP of bialgebras, for $\Fop \bicomp{\Perm} \F$, and by the PROP of special Frobenius algebras, for $\F \bicomp{\Perm} \Fop$ --- see Ex.~\ref{ex:distrlawsyntactic}. We shall characterise $\lambdapb$ and $\lambdapo$ with PROPs $\IBRw$ and $\IBRb$ respectively:
\begin{eqnarray*}
\IBRw \cong  \VectRop \bicomp{\PJ} \VectR & \qquad \qquad &  \IBRb \cong \VectR \bicomp{\PJ} \VectRop.
\end{eqnarray*}
The chosen name hints at the fact that $\IBRw$ and $\IBRb$ have the same generators as $\IBR$ and a slightly weaker equational theory: the axiomatisation for $\IBRw$ presents an asymmetry towards the white structure and the axiomatisation for $\IBRb$ one towards the black structure.

Once we have equational presentations for spans and cospans of matrices, there is an elegant construction yielding the desired characterisation~\eqref{eq:IntroIBR=SVR} for subspaces: it is given by the following ``cube'' diagram in $\PROP$.
\begin{equation}\label{eq:cube}
\tag{\mancube}
\raise30pt\hbox{$
\xymatrix@=5pt{
& {\ABR + \ABRop} \ar[dd]_(.3){\cong}|{\hole}
\ar[dl] \ar[rr] & & {\IBRw} \ar[dl] \ar[dd]^{\cong} \\
{\IBRb} \ar[rr] \ar[dd]_{\cong}  & & {\IBR} \ar@{.>}[dd]^(.3){\cong} \\
& {\VectR+ \VectRop} \ar[dl] \ar[rr]|(.51){\hole}& & {\VectRop \bicomp{\PJ} \VectR} \ar[dl] \\
{\VectR \bicomp{\PJ} \VectRop} \ar[rr] & & {\SVR}
}$}
\end{equation}
The top face is a pushout, that means, $\IBR$ is the result of merging the equational theories of $\IBRw$ and $\IBRb$. As we will prove, the bottom face of the cube is also a pushout, describing the linear algebraic nature of our theories. Commutativity of the diagram makes the isomorphism $\IBR \cong \SVR$ arise by universal property of the topmost pushout.

\medskip

Before proceeding with the formal developments, we would like to stress the importance of the \emph{modular} account of $\IBR$. Identifying the two distributive laws $\lambdapb$ and $\lambdapo$ not only enables the formulation of the isomorphism $\IBR \cong \SVR$ as a universal arrow, but also reveals the provenance of the axioms of $\IBR$: the Hopf algebras are the \emph{building blocks} of the theory, whereas the Frobenius algebras are \emph{derivative}, as they arise by the interaction of the Hopf algebras according to laws $\lambdapb$ and $\lambdapo$. Another remarkable consequence of our analysis is that $\IBR$ inherits the two factorisation properties induced by $\lambdapb$ and $\lambdapo$, that means, every string diagram of $\IBR$ enjoys a decomposition as a \emph{span} and as a \emph{cospan} of $\ABR$-diagrams.
This property allows for intriguing observations about the nature of $\IBR$. For instance, it is instrumental in showing how linear algebraic manipulations in $\SVR$ translate to equational reasoning in $\IBR$: we will see in \S~\ref{sec:graphicallinearlagebra} that the span form expresses the representation of a subspace in terms of a \emph{basis}, whereas the cospan form can be read as the representation of a subspace as a \emph{system of equations}; the fact that each diagram of $\IBR$ has both decompositions means that these encodings are all equivalent --- moreover, they are described \emph{uniformly} with the graphical syntax. Somewhat surprisingly, the two factorisations are relevant also for the reading of string diagrams as executable circuits that we will pursue in Chapter~\ref{chapter:SFG}: it turns out that being in span or cospan form has implications for the operational semantics of diagrams --- for instance, spans prevent \emph{deadlocks}.

\medskip

As a concluding note, let us mention that, just as $\ABR$, also $\IBR$ is a ubiquitous presence appearing in various research threads. In categorical quantum mechanics, $\IH{\Z_2}$ is known as the phase-free fragment of the ZX-calculus~\cite{CoeckeDuncanZX2011,BialgAreFrob14}. In concurrency theory, $\IH{\Z_2}$ is closely related to the calculus of stateless connectors~\cite{Bruni2006} and the Petri calculus~\cite{Soboci'nski2010,Sobocinski2013a}. Another important instance is $\IBpoly$~\cite{Bonchi2014b,Bonchi2015,Baez2014,BaezErbele-CategoriesInControl}, which we shall investigate extensively in Chapter~\ref{chapter:SFG} as an environment for modeling stream processing circuits. The last example is $\IH{\Z}$, which by~\eqref{eq:IntroIBR=SVR} gives a presentation by generators and relations for $\Q$-subspaces: we study $\IH{\Z}$ at the end of this chapter.
\paragraph{Synopsis} Our exposition is organised as follows.
\begin{itemize}
\item \S~\ref{sec:theorymatr} introduces the theory of $\PID$-Hopf Algebras and gives a novel, modular proof of the fact that it characterises the PROP of $\PID$-matrices (Proposition~\ref{prop:ab=vect}). Although the characterisation holds for any commutative semiring, for the sake of uniformity we assume from the start that $\PID$ is a principal ideal domain: the additional properties of $\PID$ will become relevant in \S~\ref{sec:ibrw}.
\item \S~\ref{sec:ibrw} discusses the theories $\IBRw$ and $\IBRb$. First, in \S~\ref{sec:cc} we introduce $\IBRw$ and describe its compact closed structure (Proposition~\ref{prop:snakecc}). \S~\ref{sec:completeness} is the heart of the chapter: it contains the technical developments showing that the equations of $\IBRw$ present a distributive law defined by pullback in $\VectR$. This makes us able to prove that $\IBRw$ characterises spans of $\PID$-matrices (Theorem~\ref{th:Span=IBw}). Finally, in \S~\ref{sec:IBRbCospan} we also introduce the theory $\IBRb$ and show that it characterises cospans of $\PID$-matrices (Theorem~\ref{th:IBRb=Cospan}).
\item \S~\ref{sec:cubetop} presents the theory $\IBR$ and constructs the cube~\eqref{eq:cube}. \S\ref{sec:cubebottom} shows that the bottom face is a pushout (Theorem~\ref{th:bottomfacePushout}). \S\ref{sec:cubeback} shows commutativity of the rear faces: this passes through an inductive presentation of the isomorphisms $\IBRw \cong \VectRop \bicomp{\PJ} \VectR$ (Proposition~\ref{prop:semanticsIBRwIso}) and $\IBRb \cong \VectR \bicomp{\PJ} \VectRop$ (Proposition~\ref{prop:indPresIsoIBRb}). In \S~\ref{sec:cuberebuilt} we construct the isomorphism $\IBR \cong \SVR$ (Theorem~\ref{th:IBR=SVR}) and prove the span and cospan factorisation property of $\IBR$ (Theorem~\ref{Th:factIBR}). An inductive presentation of the isomorphism $\IBR \cong \SVR$, sketched above, is given in Definition~\ref{def:semIBRInd}. \S\ref{sec:graphicallinearlagebra} proposes some illustrative example of how linear algebraic manipulations can be carried out graphically in $\IBR$: we discuss matrices and their kernels (Example~\ref{ex:graphicalLA-matrices}), different encodings of subspaces (Example~\ref{ex:subspacesgraphical}) and prove graphically a lemma about injective matrices (Proposition~\ref{prop:graphicalinjectivematrix}).
\item \S~\ref{sec:instances} shows an example of our cube construction: the theory $\IH{\Z}$ of interacting Hopf algebras for rational subspaces.
\end{itemize}

\section{Hopf Algebras: the Theory of Matrices}\label{sec:theorymatr}

Recall that a \emph{principal ideal domain} (PID) is a commutative ring with no zero divisors in which every ideal is principal, i.e., can be generated by a single element. We refer to~\cite[\S 23]{HandbookLinearAlgebra} for an overview of linear algebra for PIDs. Throughout the chapter we fix a principal ideal domain $\PID$.

In this section we give a presentation by generators and relations of the PROP of $\PID$-matrices.
\begin{definition} The PROP $\VectR$ of $\PID$-matrices is defined as follows:
\begin{itemize}[noitemsep,topsep=0pt,parsep=0pt,partopsep=0pt]
\item arrows $n \to m$ are $m \times n$-matrices over $\PID$.
\item Composition $n \tr{A} z \tr{B} m$ is matrix multiplication $BA \: n \to m$.
\item The monoidal product $A \tns B$ is the matrix $\tiny{\left(%
                \begin{array}{cc}
                 \!\!\! A \!\! & \!\! 0 \!\!\!\\
                 \!\!\! 0 \! \!& \!\! B\!\!\!
                \end{array}\right)}$.
\item The symmetries are the rearrangements of the rows of the identity matrix. For instance, the symmetry $1+2 \to 2+1$ is the matrix $\tiny{\left(%
                \begin{array}{ccc}
                 \!\!\! 0 \!\! &\!\!\! 1 \!\! & \!\! 0 \!\!\!\\
                 \!\!\! 0 \!\! &\!\!\! 0 \! \!& \!\! 1\!\!\! \\
                 \!\!\! 1 \!\! &\!\!\! 0 \!\! & \!\! 0 \!\!\!
                \end{array}\right)}$.
\end{itemize}
\end{definition}

Before the technical developments, let us anticipate the leading idea of a diagrammatic presentation of $\PID$-matrices. An $m \times n$-matrix $M$ will correspond to a string diagram with $n$ ports on the left (representing columns) and $m$ ports on the right (representing rows). We draw a link from port $i$ on the left to port $j$ on the right whenever $M_{ji}$ has a non-zero value $k \in \PID$, in which case we weight the link with a $1 \to 1$ diagram $\scalar$. If $k = 1$, we can also omit drawing $\scalar$ on the link. Additional generators are $\Bcomult$, $\Bcounit$ and $\Wmult$, $\Wunit$, giving the branching that permits zero or more column to be connected to zero or more rows. An example is given below:
\begin{eqnarray}\label{eq:matrixform}
M = {\scriptsize \left(%
\begin{array}{ccc}
  k_1 & 0 & 0 \\
  1 & 0 & 0 \\
  k_2 & 1 & 0 \\
  0 & 0 & 0
\end{array}\right)}
 & \qquad \qquad &
\lower25pt\hbox{$\includegraphics[width=70pt]{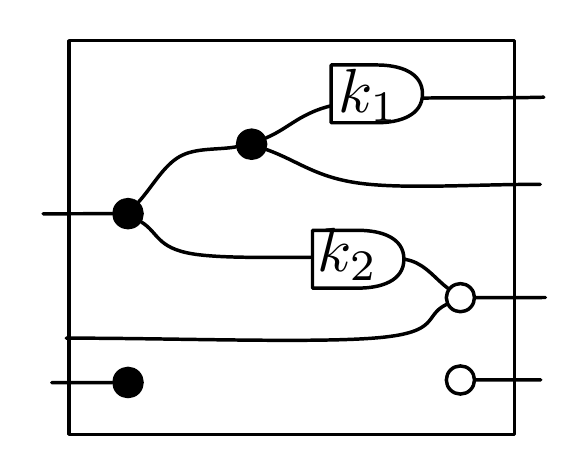}$}
\end{eqnarray}
%

The diagrammatic theory for $\Mat{\PID}$ will be constructed in a modular fashion, using the technique of composing PROPs developed in Chapter~\ref{sec:background}. The above sketch suggests that three PROPs will be involved: $\bcom$ (for $\Bcomult$ and $\Bcounit$), $\wmon$ (for $\Wmult$ and $\Wunit$) --- introduced in Ex.~\ref{ex:equationalprops} --- and a newly introduced PROP $\PROPR$, for diagrams of shape $\scalar$.
Formally, we let $\PROPR$ be the PROP generated by the signature consisting of $\scalar$ for each $k \in \PID$ and the following equations, where $k_1,k_2$ range over~$\PID$. 
\begin{multicols}{2}\noindent
\begin{equation}
\label{eq:unitscalar}\tag{A11}
\lower6pt\hbox{$\includegraphics[height=.6cm]{graffles/unitscalar.pdf}$}
=
\lower6pt\hbox{$\includegraphics[height=.6cm]{graffles/idcircuit.pdf}$}
\end{equation}
\begin{equation}
\label{eq:scalarmult}
\tag{A12}
\lower6pt\hbox{$\includegraphics[height=.7cm]{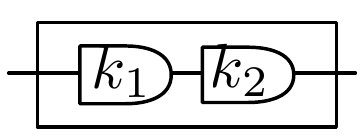}$}
=
\lower5pt\hbox{$\includegraphics[height=.7cm]{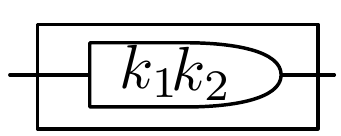}$}
\end{equation}
\end{multicols}
Now, to understand how our building blocks $\wmon$, $\bcom$ and $\PROPR$ should be composed together, observe the factorisation of the string diagram in~\eqref{eq:matrixform}:
\[\includegraphics[width=70pt]{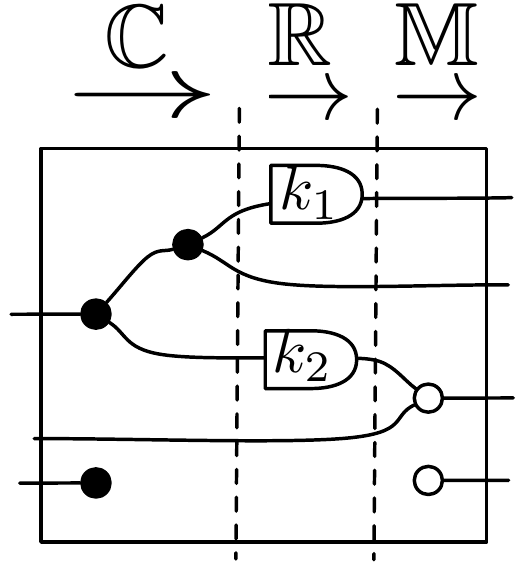}.\]
This suggests to base the theory of matrices on the composite $\bcom \bicomp{\Perm} \PROPR \bicomp{\Perm} \wmon$, which we shall construct using iterated distributive laws of PROPs (\S~\ref{sec:iteratedDistrLaws}). For this purpose, we describe first how $\PROPR$ interacts with $\bcom$ and $\wmon$. 
\begin{lemma}~
\label{lemma:threelaws}
\begin{enumerate}[(a)]
  \item There is a distributive law $\chi \: \wmon \bicomp{\Perm} \PROPR \To \PROPR \bicomp{\Perm} \wmon$ yielding a PROP $\PROPR \bicomp{\Perm} \wmon$ presented by the equations of $\PROPR + \wmon$ and, for all $k \in \PID$:
       \begin{multicols}{2}\noindent
\begin{equation}
\label{eq:scalarwmult}
\tag{A13}
\lower11pt\hbox{$\includegraphics[height=.9cm]{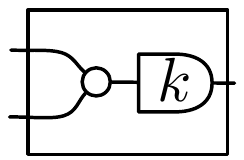}$}
=
\lower11pt\hbox{$\includegraphics[height=.9cm]{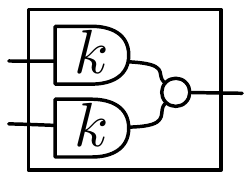}$}
\end{equation}
\begin{equation}
\label{eq:scalarwunit}
\tag{A14}
\lower7pt\hbox{$\includegraphics[height=.7cm]{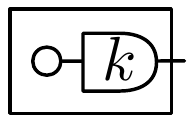}$}
=
\lower7pt\hbox{$\includegraphics[height=.7cm]{graffles/Wunit.pdf}$}
\end{equation}
\end{multicols}
  \item There is a distributive law $\psi \: \PROPR \bicomp{\Perm} \bcom \To \bcom\bicomp{\Perm}\PROPR$ yielding a PROP $\bcom\bicomp{\Perm}\PROPR$ presented by the equations of $\bcom + \PROPR$ and, for all $k \in \PID$:
\begin{multicols}{2}
\noindent
\begin{equation}
\label{eq:scalarbcomult}
\tag{A15}
\lower10pt\hbox{$\includegraphics[height=.9cm]{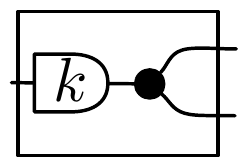}$}
=
\lower10pt\hbox{$\includegraphics[height=.9cm]{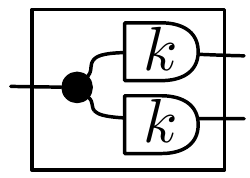}$}
\end{equation}
\begin{equation}
\label{eq:scalarbcounit}
\tag{A16}
\lower7pt\hbox{$\includegraphics[height=.7cm]{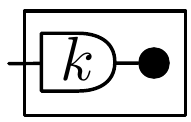}$}
=
\lower7pt\hbox{$\includegraphics[height=.7cm]{graffles/Bcounit.pdf}$}
\end{equation}
\end{multicols}
\end{enumerate}
\end{lemma}
\begin{proof}
First observe that $\PROPR$ is self-dual, the iso $\PROPR \cong \PROPR^{\op}$ being the unique PROP morphism that maps a generator $\scalar$ to $\coscalar$. Since we also have $\Mon \cong \Com^{\op}$, by Proposition~\ref{prop:distrlawOP} statement $(b)$ implies $(a)$. The proof is concluded by observing that $(b)$ is just an instance of Theorem~\ref{Th:LawvereCompositePROP}. \end{proof}

The distributive laws of Lemma~\ref{lemma:threelaws} are now combined with $\lambda \: \wmon \bicomp{\Perm} \bcom \To \bcom \bicomp{\Perm} \wmon$ introduced in Ex.~\ref{ex:distrlawsyntactic}.\ref{ex:distrlawsyntactic1} to construct the composite PROP $\bcom \bicomp{\Perm} \PROPR \bicomp{\Perm} \wmon$.
\begin{proposition} There is a distributive law yielding a PROP $\bcom \bicomp{\Perm} \PROPR \bicomp{\Perm} \wmon$ presented by the signature $\{\Wmult,\Wunit,\Bcounit,\Bcomult,\scalar \mid k \in \PID\}$ and equations \eqref{eq:wmonassoc}-\eqref{eq:scalarbcounit}.
\end{proposition}
\begin{proof} By Proposition~\ref{prop:iterateddistrlaw}, to form the PROP $\bcom \bicomp{\Perm} \PROPR \bicomp{\Perm} \wmon$ it suffices to check that the three distributive laws $\lambda$, $\sigma$ and $\psi$ satisfy the Yang-Baxter compatibility condition. This is given by commutativity of the following diagram, which can be easily verified by case analysis on the string diagrams $\tr{\in \Mon}\tr{\in \PROPR}\tr{\in \Com}$.
\[\xymatrix@R=10pt{
& \wmon \bicomp{\Perm}\bcom \bicomp{\Perm}\PROPR \ar[r]^{\lambda \PROPR} & \bcom \bicomp{\Perm}\wmon \bicomp{\Perm}\PROPR \ar[dr]^{\bcom \chi} &\\
\wmon \bicomp{\Perm}\PROPR \bicomp{\Perm}\bcom \ar[ur]^{\wmon \psi} \ar[dr]_{\chi {\bcom}} & & & \bcom \bicomp{\Perm}\PROPR \bicomp{\Perm}\wmon\\
& \PROPR \bicomp{\Perm}\wmon \bicomp{\Perm}\bcom \ar[r]^{\PROPR \lambda} & \PROPR \bicomp{\Perm} \bcom \bicomp{\Perm}\wmon \ar[ur]_{\psi {\wmon}} &\\
}\]
The PROP $\bcom \bicomp{\Perm} \PROPR \bicomp{\Perm} \wmon$ is equivalently defined by distributive laws $\lambda \PROPR \poi \bcom \chi$ and $\PROPR \lambda \poi \psi \wmon$. Proposition~\ref{prop:SMTforIteratedComp} guarantees that $\bcom \bicomp{\Perm} \PROPR \bicomp{\Perm} \wmon$ is presented by generators and equations as described in the statement.
\end{proof}

The PROP $\PROPR$ only accounts for the multiplicative part of $\PID$. In order to describe also its additive part, and thus faithfully capture $\PID$-matrices, we need to quotient $\bcom \bicomp{\Perm} \PROPR \bicomp{\Perm} \wmon$ by two more equations.

\begin{definition}\label{def:HA} The PROP $\ABR$ ($\PID$-\textbf{H}opf \textbf{A}lgebras) is defined as the quotient of $\bcom \bicomp{\Perm} \PROPR \bicomp{\Perm} \wmon$ by the following equations, for all $k_1,k_2 \in \PID$:
\begin{multicols}{2}
\noindent
\begin{equation}
\label{eq:zeroscalar}
\tag{A17}
\lower6pt\hbox{$\includegraphics[height=.6cm]{graffles/zeroscalar.pdf}$}
=
\!\!
\lower11pt\hbox{$\includegraphics[height=1cm]{graffles/zeroscalar2.pdf}$}
\end{equation}
\begin{equation}
\label{eq:scalarsum}
\tag{A18}
\lower12pt\hbox{$\includegraphics[height=1cm]{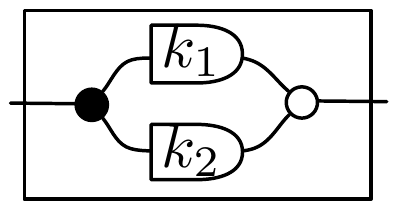}$}
=
\lower7pt\hbox{$\includegraphics[height=.7cm]{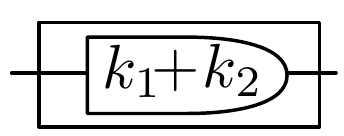}$}
\end{equation}
\end{multicols}
\end{definition}
\begin{remark}\label{rmk:hopf} The name ``$\PID$-Hopf algebra'' is justified by the case of the initial principal ideal domain $\Z$. Indeed, as observed in \S~\ref{sec:overview}, $\HA{\Z}$ encodes Hopf algebras as they are traditionally defined in the literature (see e.g.~\cite{Fiore2013,HopfRef1,HopfRef2}). That means, $\HA{\Z}$ can be presented by the bialgebra equations~\eqref{eq:wmonassoc}-\eqref{eq:bwbone} and the instances of \eqref{eq:unitscalar}-\eqref{eq:scalarsum} where $k$ is either $0$, $1$ or $-1$. In particular, $\scalarminusone$ plays the role of the \emph{antipode}, for which reason hereafter we fix notation $\antipode \df\! \scalarminusone$. The usual \emph{Hopf law} is formulated as follows:
  \begin{equation} \label{eq:hopf}
  \tag{Hopf}
  \lower11pt\hbox{$\includegraphics[height=1.1cm]{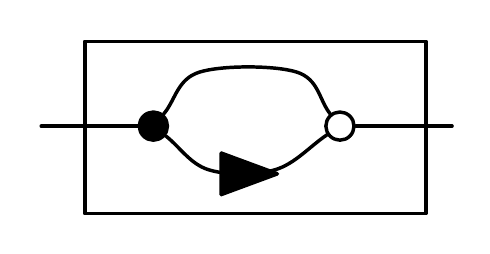}$}
  = \lower10pt\hbox{$\includegraphics[height=1cm]{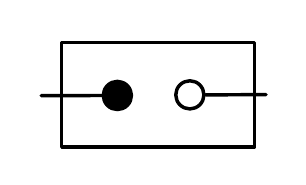}$}
  = \lower11pt\hbox{$\includegraphics[height=1.1cm]{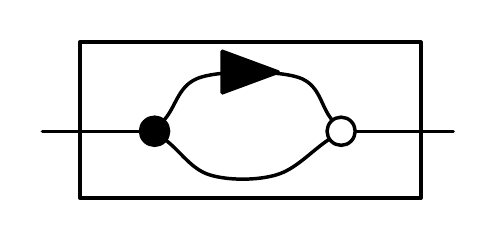}$}
   \end{equation}
   The derivation of \eqref{eq:hopf} in $\HA{\PID}$, for any PID $\PID$, is reported in \eqref{eq:HopfLawDer}. We will return on the analysis of $\HA{\Z}$ in~\S\ref{sec:instances}.
  \end{remark}

We now recover the correspondence in~\eqref{eq:matrixform} by describing a class of string diagrams of $\ABR$ for which one can easily read off the associated matrices.

\begin{definition}\label{def:matrixform}
A string diagram $\tr{\in \ABR}$ is in \emph{matrix form} if it is of shape $\tr{\in \bcom}\tr{\in \PROPR}\tr{\in \wmon}$, any port on the left boundary has exactly one connection with any port on the right boundary and any such connection passes through exactly one scalar $\scalar$. We say that there is a \emph{$k$-path from $i$ to $j$} if $k$ is the scalar on the path from the $i$th port on the left to the $j$th port on the right, assuming a top-down enumeration.
\end{definition}
When drawing matrix forms, for the sake of readability it will be often convenient to massage the above definition as follows: we typically omit to draw the scalar $k =1$, by virtue of~\eqref{eq:unitscalar}, and omit the scalar $k=0$, by~\eqref{eq:zeroscalar}, leaving the ports in question disconnected.

\begin{example}\rm\label{ex:matrixform} The string diagram $c$ in~\eqref{eq:matrixform} is in matrix form. The matrix on the left can be derived from $c$ as follows: put $M_{ij}=k$ exactly when there is a $k$-path from port $j$ to port $i$ in $c$.
\end{example}

\noindent We will often write $\circuitAdots$ for the string diagram, in matrix form, encoding the matrix $A$. We now extend the matrix interpretation to arbitrary string diagrams of $\ABR$. 

\begin{definition}\label{def:sem}
The PROP morphism $\sem{\ABR} \: \ABR \to \VectR$ is defined inductively on terms of the SMT of~$\ABR$:
\begin{eqnarray}
\begin{aligned}\label{eq:defsemABR1}
                    \Wunit \mapsto\  \initVect &&
                    \Bcounit\mapsto\ \finVect
                     &&
                   \Wmult \mapsto   {\scriptsize\left(%
                \begin{array}{cc}
                \!\!\!  1 \! &\!\! 1 \!\!\!
                \end{array}\right)}
                &&
                 \Bcomult \mapsto\ \tiny{\left(%
                \begin{array}{c}
                 \!\! 1 \!\!\\
                 \!\! 1\!\!
                \end{array}\right)}
                &&
                \scalar \mapsto\ {\scriptsize  \left(%
                \begin{array}{c}
                 \!\!\! k\!\!\!
                \end{array}\right)}
            \end{aligned}
\end{eqnarray}
\vspace{-.4cm}
\begin{eqnarray*}
            \begin{aligned}
            \IdnetT \mapsto {\scriptsize \left(%
		\begin{array}{c}
		  \!\! 1\!\!	
                 \end{array}\right)}
            &&
\symNetT \mapsto {\scriptsize \left(%
		\begin{array}{c c}
		  \!\! 0 & 1 \!\!\\
		  \!\! 1 & 0 \!\!	
                 \end{array}\right)}
                 &&
                   s\tns t  \mapsto  \sem{\ABR}(s) \tns\sem{\AB}(t)
                   &&
                  s \poi t \mapsto\ \sem{\ABR}(s) \poi \sem{\AB}(t)
                  &&
                 \end{aligned}
\end{eqnarray*}
                 where $\initVect \: 0 \to 1$ and $\finVect \: 1 \to 0$ are the arrows given by, respectively, initiality and finality of $0$ in $\VectR$. 
                 It readily follows that $\sem{\ABR}$ is well-defined, as it respects the equations of $\ABR$. 
\end{definition}
It is immediate to check that the above assignment is indeed functorial and respects the symmetric monoidal structure. In fact, $\sem{\ABR}$ could also have been defined as the unique PROP morphism mapping the generators of $\ABR$ as in \eqref{eq:defsemABR1}. In the rest of this thesis, we will introduce several inductively defined maps and, to be concise, we will usually adopt this second formulation.

One can easily verify that, for string diagrams in matrix form, the above inductive definition coincides with the intuitive picture given in Example~\ref{ex:matrixform}. We now show the main result of this section: $\ABR$ is a presentation by generators and equations of $\Mat{\PID}$, with $\sem{\ABR}$ the witnessing isomorphism.
\begin{proposition}\label{prop:ab=vect} $\sem{\ABR} \: \ABR \to \VectR$ is an isomorphism of PROPs.
\end{proposition}
Since the two categories have the same objects, it suffices to prove that $\sem{\ABR}$ is full and faithful.
Fullness is immediate: given a matrix $M$, it is clear how to generalise the correspondence described in Example~\ref{ex:matrixform} to construct a string diagram $c$ in matrix form encoding $M$. One can then verify that $\sem{\ABR}(c) = M$.

To show faithfulness, the idea is to reduce to the case of matrix forms. This is possible because, as shown in the following lemma, any string diagram of $\ABR$ can be put in matrix form: the key point is that $\ABR$ has been constructed modularly, using the technique of composing PROPs, and thus it enjoys a suitable factorisation property.

\begin{lemma}\label{lem:path1}
For all string diagrams $c \in \ABR[n,m]$, there exists $d \in \ABR[n,m]$ in matrix form such that $c = d$.
\end{lemma}
\begin{proof}
Because $\ABR$ is a quotient of the composite PROP $\bcom \bicomp{\Perm} \PROPR \bicomp{\Perm} \wmon$, it inherits its factorisation property: this means that $\tr{c \in \ABR}$ can be decomposed as $\tr{\in \bcom}\tr{\in \PROPR}\tr{\in \wmon}$. Then, by using associativity of $\Wmult$ and $\Bcomult$, \eqref{eq:sliding2} and \eqref{eq:zeroscalar}-\eqref{eq:scalarsum}, we make any port on the left boundary have exactly one connection with any port on the right and, by applying \eqref{eq:scalarmult}, we make any such connection pass through exactly one scalar $\scalar$. The resulting string diagram is in matrix form.
\end{proof}

The next step is to show faithfulness for matrix forms.

\begin{lemma} \label{lem:path2} Let $c$ and $d$ be string diagrams $n \to m$ in matrix form such that $\sem{\ABR}(c) = \sem{\ABR}(d)$. Then $c = d$ in $\ABR$.\end{lemma}
\begin{proof} The cases when $n=0$ or $m=0$ are easy: for any $n,m\in \N$
there is a unique string diagram of type $0\to m$ in matrix form, and similarly,
a unique matrix form diagram of type $n\to 0$, displayed below.
\[
\lower8pt\hbox{$\includegraphics[height=.8cm]{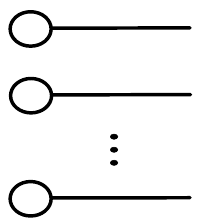}$}  \ \ \: 0 \to m \qquad\qquad\qquad
\lower8pt\hbox{$\includegraphics[height=.8cm]{graffles/manyBcounits.pdf}$} \ \ \: n \to 0
\]
We then argue by induction on $n$ and $m$ bigger than $0$. The base case is $m=n=1$: this means that $c$ and $d$ denote a $1 \times 1$ matrix $(k)$ via $\sem{\ABR}$. As they are in matrix form, both $c$ and $d$ must be the string diagram $\scalar$ and the statement clearly holds.

Consider now the case $c,d \: n+1 \to m$. Let $M$ the $m \times (n+1)$ matrix denoted by $c$ and $d$ and let $k_1 = M_{1,m+1}, k_2 = M_{2,m+1}, \dots , k_m = M_{n,m+1}$ be the entries in the last column of $M$.
By definition of matrix form, there is a $k_i$-path from port $n+1$ on the left boundary of $c$ to port $i$ on the right boundary, for each $i \in \{1,\dots,m\}$. This means that, using laws~\eqref{eq:sliding1}-\eqref{eq:sliding2} of symmetric strict monoidal categories, we can transform $c$ into the string diagram:
\[
c_1 \quad \df \quad \lower25pt\hbox{$\includegraphics[height=2.2cm]{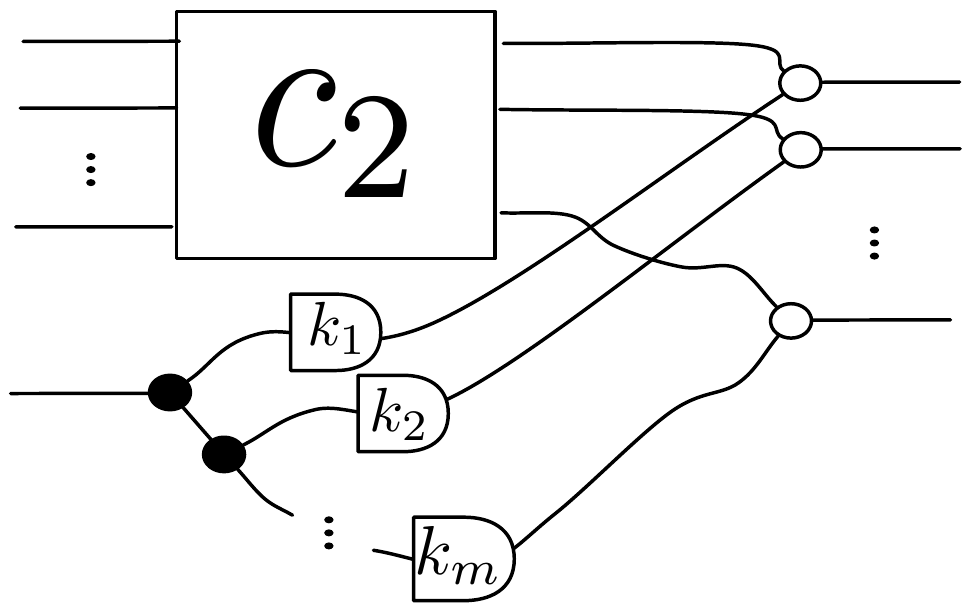}$}
\]
The idea is that $c_1$ has the same connectivity as $c$, but we pulled out the part of the diagram that encodes column $n+1$. Because $d$ is also in matrix form and denotes $M$, then port $n+1$ on the left boundary of $d$ must have a $k_i$-path to port $i$ on the right boundary, for each $i \in \{1,\dots,m\}$. We can thus transform $d$ analogously:
\[
d_1 \quad \df \quad \lower25pt\hbox{$\includegraphics[height=2.2cm]{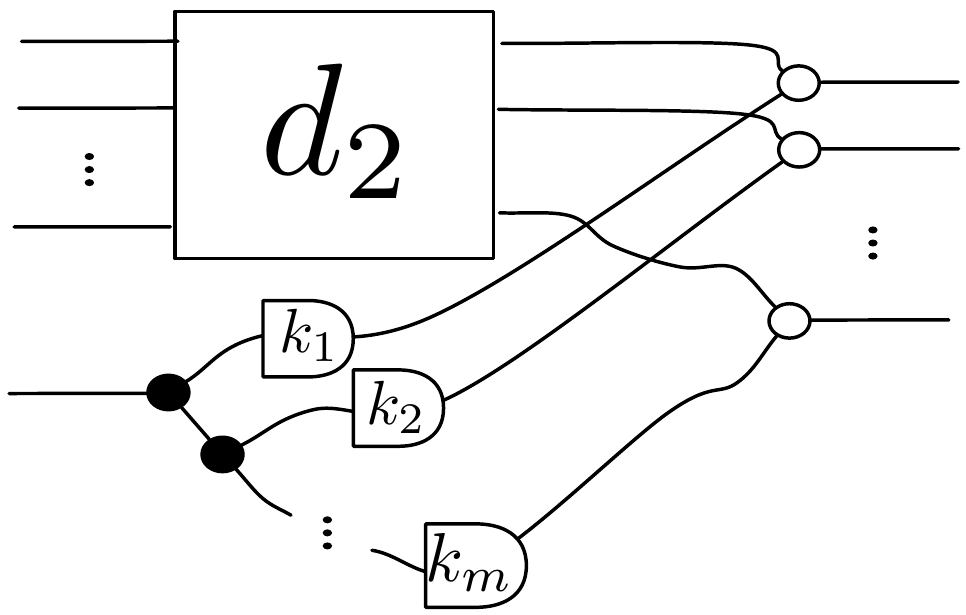}$}
\]
Now, string diagrams $c_2, d_2 \: n \to m$ are in matrix form and denote via $\sem{\ABR}$ the matrix consisting of the first $n$ columns of $M$. Thus, by inductive hypothesis, $c_2 = d_2$ in $\ABR$. It is then immediate that also $c_1 = d_1$. Because $c_1$ and $d_1$ were derived from $c$ and $d$ by only applying equations of $\ABR$, it follows that $c=d$.

The remaining inductive case to consider is the one in which $c$ and $d$ have type $n \to m+1$. We can reason in a dual manner. First, let $N$ be the $(m+1) \times n$ matrix denoted by $c$ and $d$ and let $l_1 = N_{m+1,1}, l_2 = N_{m+1,2}, \dots , l_n = N_{m+1,n}$ be the entries in the last row of $N$. Symmetrically to the other inductive case, we can pull out the part of diagrams $c$ and $d$ that describes row $m+1$:
\[
c_1' \df \lower25pt\hbox{$\includegraphics[height=2.2cm]{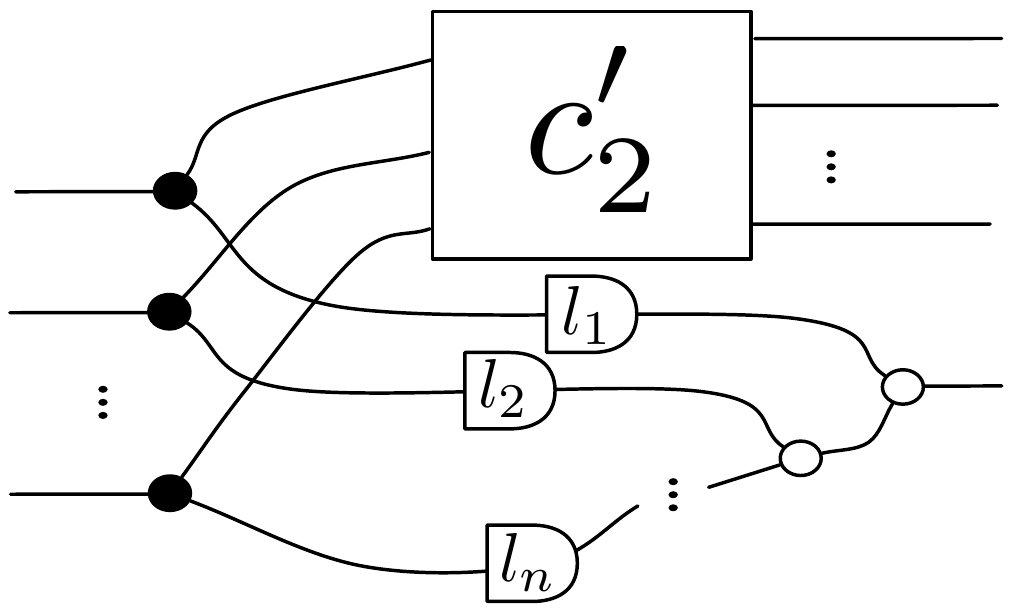}$} \qquad \qquad d_1' \df \lower25pt\hbox{$\includegraphics[height=2.2cm]{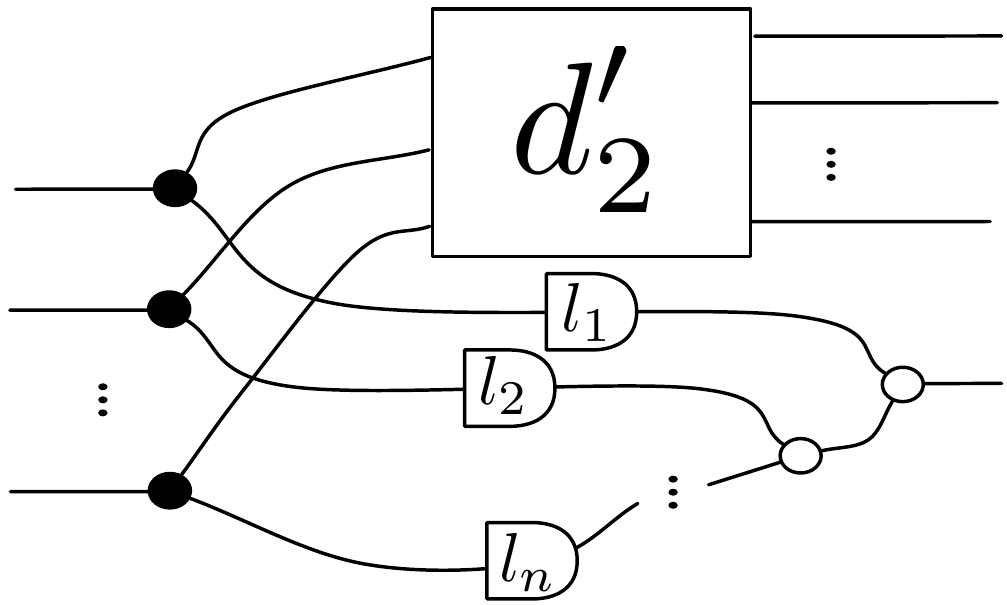}$}
\]
String diagrams $c_1'$ and $d_1'$ can be obtained from $c$ and $d$ by applying laws~\eqref{eq:sliding1}-\eqref{eq:sliding2} of SMCs. Diagrams $c_2'$ and $d_2'$ are in matrix form and denote the matrix consisting of the first $m$ rows of $N$. Thus by inductive hypothesis $c_2' = d_2'$ and we conclude that also $c = d$.
\end{proof}

\begin{proof}[Proof of Proposition~\ref{prop:ab=vect}]
We complete the proof by showing the faithfulness of $\sem{\ABR}$.
Suppose that $\sem{\ABR}(c_1)=\sem{\ABR}(c_2)$. Then, by the conclusion of Lemma~\ref{lem:path1}
there exist $d_1, d_2$ in matrix form with $d_1=c_1$ and $d_2=c_2$. It follows that $\sem{\ABR}(d_1) = \sem{\ABR}(d_2)$ and, by the conclusion of Lemma~\ref{lem:path2}, $d_1=d_2$. Therefore, $c_1=d_1=d_2=c_2$.
\end{proof}

%

In the rest of the chapter we shall investigate the composition of $\ABR$ with its opposite PROP $\ABRop$. Following the same convention used in Chapter~\ref{sec:background}, we shall represent the string diagrams of $\ABRop$ as those of $\ABR$ reflected about the $y$-axis. That means, $\ABRop$ is presented by generators
$$ \Bunit \quad \Bmult \quad \Wcounit \quad \Wcomult \quad \coscalar \quad\quad
\quad k \in \PID$$
and equations \eqref{eq:wmonassoc}$^{\op}$-\eqref{eq:scalarsum}$^{\op}$ of Figure~\ref{fig:axiomsIH}. The duality between $\ABR$ and $\ABRop$ is witnessed by the contravariant identity PROP morphism $\contrid{\cdot} \: \ABR \to \ABRop$ that reflects a diagram about the $y$-axis. It can be inductively defined as follows:
  \begin{multicols}{5}
\noindent
      \begin{eqnarray*}
     \Bcounit \mapsto \Bunit
    \end{eqnarray*}
  \begin{eqnarray*}
     \Wunit \mapsto \Wcounit
    \end{eqnarray*}
     \begin{eqnarray*}
     \Wmult \mapsto \Wcomult
    \end{eqnarray*}
  \begin{eqnarray*}
   \Bcomult \mapsto \Bmult
   \end{eqnarray*}
           \begin{eqnarray*}
     \scalar \mapsto \coscalar
    \end{eqnarray*}
    \end{multicols}
    \smallskip
     \begin{multicols}{2}
     \noindent
     \begin{equation*}    \lower10pt\hbox{$\includegraphics[height=.9cm]{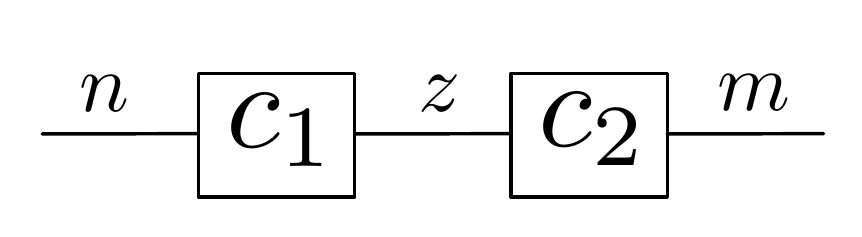}$}
  \!\mapsto\! \lower8pt\hbox{$\includegraphics[height=.7cm]{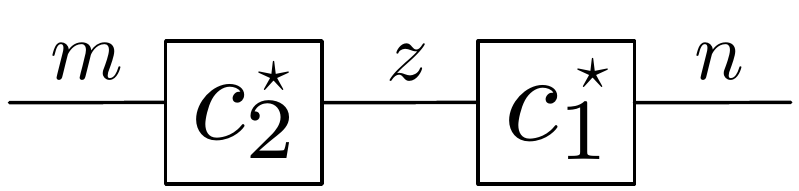}$}
   \end{equation*}
        \begin{equation*}
         \lower20pt\hbox{$\includegraphics[height=1.5cm]{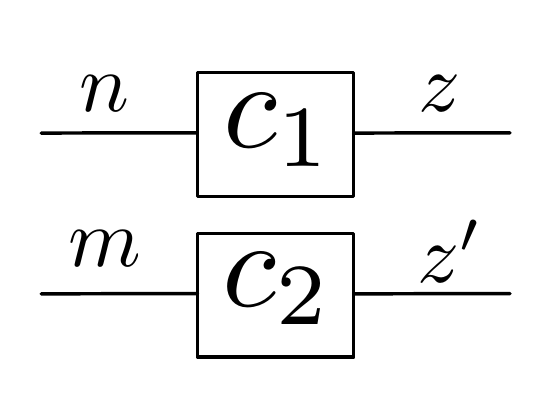}$}
   \!\mapsto\!\lower18pt\hbox{$\includegraphics[height=1.3cm]{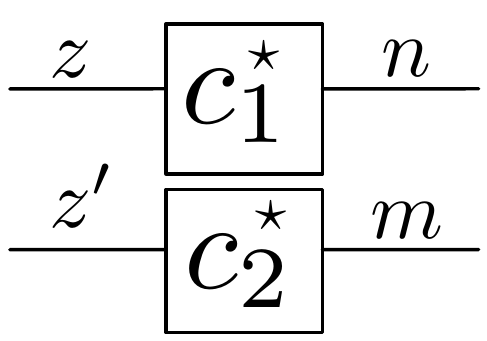}$}
   \end{equation*}
 \end{multicols}
 Just as diagrams of $\ABR$, also diagrams of $\ABRop$ can be interpreted as matrices. The interpretation is given by the  isomorphism $\sem{\ABR}^{\op} \: \ABRop \to \VectRop$ induced by the isomorphism $\sem{\ABR} \: \ABR \to \VectR$. To understand the action of $\sem{\ABR}^{\op}$, consider arrows in $\VectRop[n,m]$ as matrices in $\VectR[m,n]$. This means that, since $\sem{\ABR}$ maps $\Bcomult$ to ${\tiny \matrixOneOne} \in \VectR[1,2]$, then $\sem{\ABR}^{\op}$ maps $\Bmult$ to ${\tiny \matrixOneOne} \in \VectRop[2,1]$. Therefore, to associate matrices with string diagrams of $\ABRop$, one should intuitively follow the same procedure of Example~\ref{ex:matrixform}, but now reading the diagram from right to left, i.e., columns are ports on the right boundary and rows are ports on the left boundary.

 Throughout the thesis, we shall draw \lower5pt\hbox{$\includegraphics[height=.5cm]{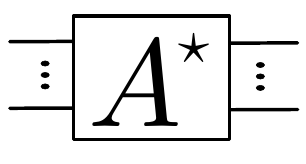}$} for the image under $\contrid{\cdot}$ of the $\ABR$-circuit $\circuitAdots$ in matrix form representing the matrix $A$.

 \begin{remark}[Transpose] \label{rmk:graphicaltranspose}
 The PROP $\VectR$ is self-dual: the operation of matrix transpose yields the desired isomorphism $\trasp{\cdot} \: \VectR \to \VectRop$. Note that the transpose has nothing to do with the transformation $\contrid{\cdot} \: \ABR \to \ABRop$ defined above. Instead, the diagrammatic counterpart of $\trasp{\cdot}$ is the (covariant) PROP isomorphism $\ABR \to \ABRop$ that swaps the black/white colouring of a string diagram and the orientation of components $\scalar$, for $k \in \PID$. For instance, the transpose of the matrix $M$ in \eqref{eq:matrixform}, below left, is represented (via $\sem{\ABR}^{\op}$) by the diagram $c \in \ABRop[3,4]$ on the right: one can obtain $c$ from the $\ABR$-diagram in~\eqref{eq:matrixform} by applying the operation that we just described. 
 \begin{eqnarray*}
\trasp{M} = {\scriptsize \left(%
\begin{array}{cccc}
  k_1 & 1 & k_2 & 0 \\
  0 & 0 & 1 & 0\\
  0 & 0 & 0 & 0
\end{array}\right)}
 & \qquad \qquad &
{\sem{\ABR}^{\op}}^{-1}(\trasp{M}) = c = \lower25pt\hbox{$\includegraphics[width=70pt]{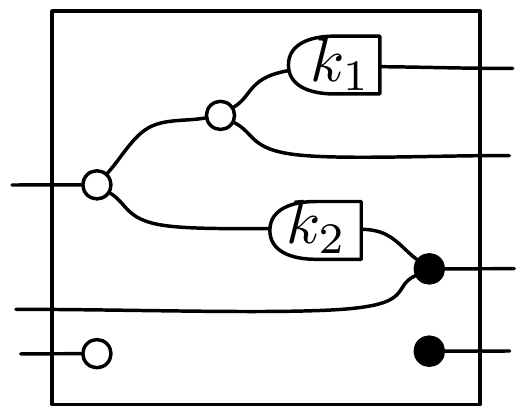}$}
\end{eqnarray*}

 We will investigate the diagrammatic representation of transposing a matrix more thoroughly in~\S~\ref{sec:IBRbCospan}. This is just a glimpse of our graphical perspective on linear algebra, of which we shall give more examples throughout the chapter (see e.g. Example~\ref{ex:HNF}, Proposition~\ref{prop:graphicalinjectivematrix} and \S~\ref{sec:graphicallinearlagebra}).
 \end{remark}

\section{Interacting Hopf Algebras I: Spans and Cospans of Matrices}\label{sec:ibrw}

In this section we commence the exploration of several theories that arise from
composing $\ABR$ with $\ABRop$, which is the main focus and contribution of this chapter.
Collectively, we refer to them as \emph{interacting Hopf algebras}.

We first introduce $\IBRw$ --- the superscript $\SupSpan$ represents the fact that $\IBRw$ will be shown to be the theory of \emph{spans} of $\PID$-matrices. In \S~\ref{sec:IBRbCospan} we introduce $\IBRb$, which will be shown to be the theory of \emph{cospans} of $\PID$-matrices. 



\begin{definition}\label{def:IBRw} The PROP $\IBRw$ is the quotient of $\ABR + \ABRop$ by the following additional equations, where $l$ is any non-zero element and $k$ any element of $\PID$.
\begin{multicols}{2}\noindent
 \begin{equation}
\label{eq:lcm}
\tag{W1}
\lower7pt\hbox{$\includegraphics[height=.7cm]{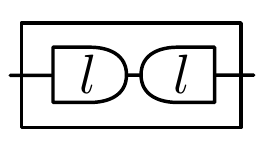}$}
=
\lower6pt\hbox{$\includegraphics[height=.6cm]{graffles/idcircuit.pdf}$}
\end{equation}
\begin{equation}
\label{eq:wbone}
\tag{W2}
\lower4pt\hbox{$\includegraphics[height=.5cm]{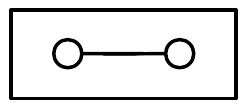}$}
=  \lower4pt\hbox{$\includegraphics[height=.5cm]{graffles/idzerocircuit.pdf}$}
\end{equation}
\end{multicols}
 \begin{multicols}{2}\noindent
\begin{equation}
\label{eq:WFrob}
\tag{W3}
\lower12pt\hbox{$\includegraphics[height=1.2cm]{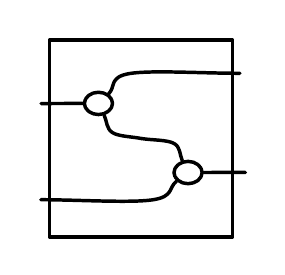}$}
\!\!
=
\!\!
\lower10pt\hbox{$\includegraphics[height=1cm]{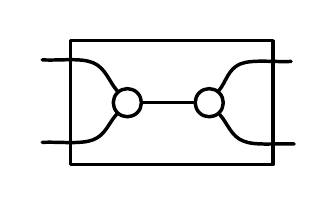}$}
\!\!
=
\!\!
\lower12pt\hbox{$\includegraphics[height=1.2cm]{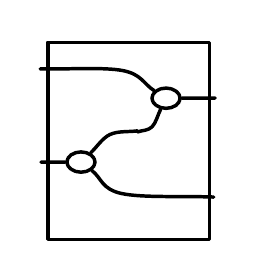}$}
\end{equation}
\begin{equation}
\label{eq:BFrob}
\tag{W4}
\lower12pt\hbox{$\includegraphics[height=1.2cm]{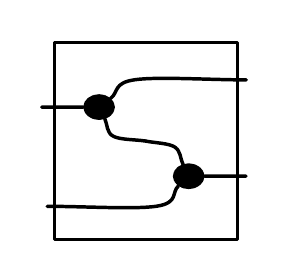}$}
\!\!
=
\!\!
\lower10pt\hbox{$\includegraphics[height=1cm]{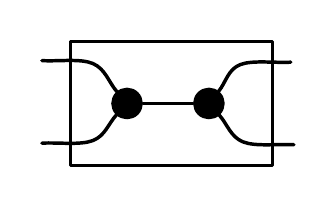}$}
\!\!
=
\!\!
\lower12pt\hbox{$\includegraphics[height=1.2cm]{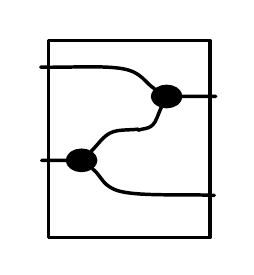}$}
\end{equation}
\end{multicols}
 \begin{multicols}{2}\noindent
\begin{equation}
\label{eq:lcc}
\tag{W5}
\lower12pt\hbox{$\includegraphics[height=1cm]{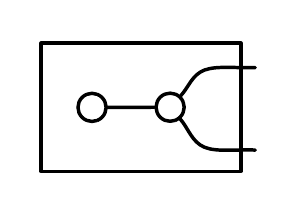}$}
\!\!
=
\!\!
\lower12pt\hbox{$\includegraphics[height=1cm]{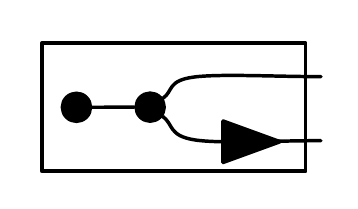}$}
\end{equation}
\begin{equation}
\label{eq:rcc}
\tag{W6}
\lower12pt\hbox{$\includegraphics[height=1cm]{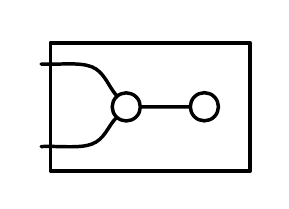}$}
\!\!
=
\!\!
\lower12pt\hbox{$\includegraphics[height=1cm]{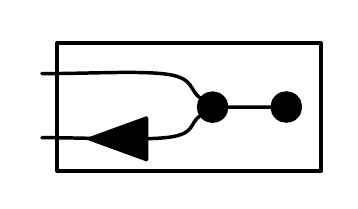}$}
\end{equation}
\end{multicols}
\begin{multicols}{2}\noindent
\begin{equation}
\label{eq:BccscalarAxiomOne}
\tag{W7}
\lower9pt\hbox{$\includegraphics[height=.8cm]{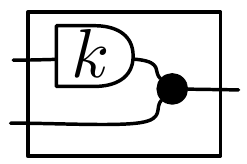}$}
=
\lower9pt\hbox{$\includegraphics[height=.8cm]{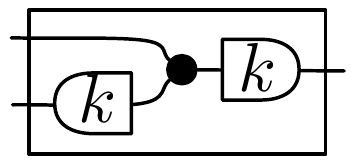}$}
\end{equation}
\begin{equation}
\label{eq:BccscalarAxiomTwo}
\tag{W8}
\lower9pt\hbox{$\includegraphics[height=.8cm]{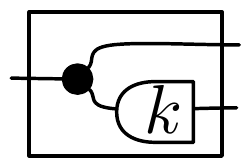}$}
=
\lower9pt\hbox{$\includegraphics[height=.8cm]{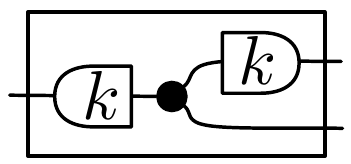}$}
\end{equation}
\end{multicols}
\end{definition}

We fix notation $\HAToIHw \: \ABR \to \IBRw$ and $\HAopToIHw \: \ABRop \to \IBRw$ for the PROP morphisms interpreting string diagrams of $\ABR$ and $\ABRop$, respectively, as string diagrams of $\IBRw$. Syntactically speaking, the generators of $\ABR$ together with those of $\ABRop$ are also the generators of $\IBRw$ and therefore we will often abuse notation by confusing $c$ in $\ABR$ with $\HAToIHw(c)$ in $\IBRw$, and the same for $\ABRop$.

\medskip

The following are some of the derived laws of $\IBRw$, where $k$ is any element and $l$ any non-zero element of $\PID$. We refer to Appendix~\ref{AppDerLaws} for their equational proof. In~\eqref{eq:wunitcancelbcomult} below and in the sequel, we shall use the shorthand notation $\lower7pt\hbox{$\includegraphics[height=18pt]{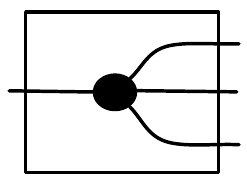}$}$ for the comultiplication $\lower7pt\hbox{$\includegraphics[height=18pt]{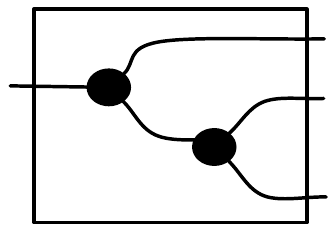}$}$ from $1$ to $3$, and more generally $\lower7pt\hbox{$\includegraphics[height=18pt]{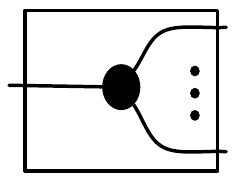}$}$ for the one from $1$ to an arbitrary $n$. This convention is harmless by~\eqref{eq:bcomonassoc}. We will adopt an analogous notation for multiplications $\Wmult$ of arity bigger than $2$.

\begin{multicols}{2}\noindent
\begin{equation}
\label{eq:lccb}
\tag{D1}
\lower7pt\hbox{$\includegraphics[height=.7cm]{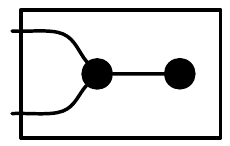}$}
\!
=
\!
\lower7pt\hbox{$\includegraphics[height=.7cm]{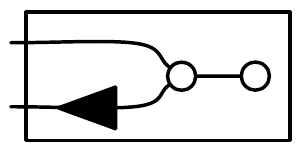}$}
\end{equation}
\begin{equation}
\label{eq:rccb}
\tag{D2}
\lower7pt\hbox{$\includegraphics[height=.7cm]{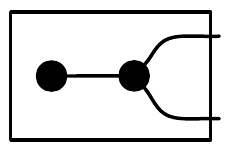}$}
\!
=
\!
\lower7pt\hbox{$\includegraphics[height=.7cm]{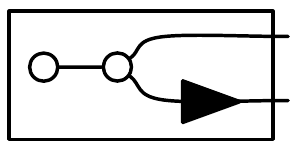}$}
\end{equation}
\end{multicols}
\begin{multicols}{2}\noindent
\begin{equation}
\label{eq:uniqueantipode}
\lower3pt\hbox{$
\tag{D3}
\lower5.5pt\hbox{$\includegraphics[height=.55cm]{graffles/antipode.pdf}$}
=
\lower5.5pt\hbox{$\includegraphics[height=.55cm]{graffles/antipodeop.pdf}$}
$}
\end{equation}
\begin{equation}\label{eq:QFrob}
\tag{D4}
\lower12pt\hbox{$\includegraphics[height=.9cm]{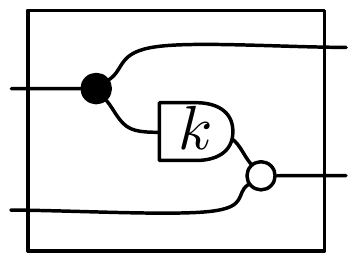}$} =
\lower12pt\hbox{$\includegraphics[height=.9cm]{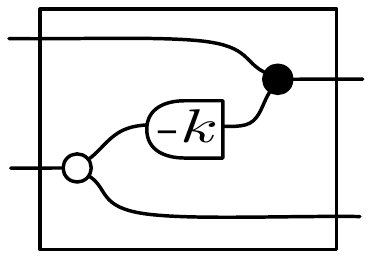}$}
\end{equation}
\end{multicols}
\begin{multicols}{2}\noindent
\begin{equation}
\label{eq:coscalarwunit}
\tag{D5}
\lower8pt\hbox{$\includegraphics[height=.7cm]{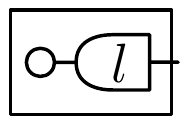}$}
=
\lower6pt\hbox{$\includegraphics[height=.6cm]{graffles/Wunit.pdf}$}
\end{equation}
\begin{equation}
\label{eq:scalarwcounit}
\tag{D6}
\lower7pt\hbox{$\includegraphics[height=.7cm]{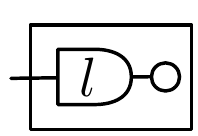}$}
 =
\lower6pt\hbox{$\includegraphics[height=.6cm]{graffles/Wcounit.pdf}$}
\end{equation}
\end{multicols}
\begin{multicols}{2}\noindent
\begin{align}
\label{eq:coscalarbcomult}
\tag{D7}
\lower10pt\hbox{$\includegraphics[height=.9cm]{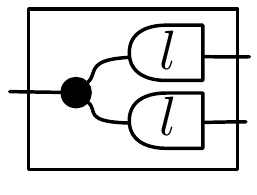}$}
=
\lower8pt\hbox{$\includegraphics[height=.8cm]{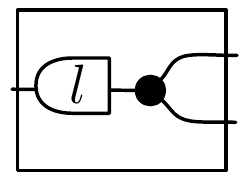}$}
\end{align}
\begin{equation}
\label{eq:scalarbmult}
\tag{D8}
\lower8pt\hbox{$\includegraphics[height=.8cm]{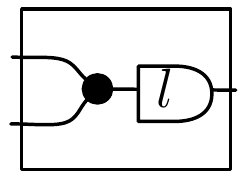}$}
=
\lower7.5pt\hbox{$\includegraphics[height=.9cm]{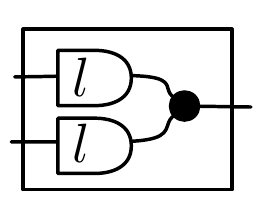}$}
\end{equation}
\end{multicols}
\noindent
\begin{equation}\label{eq:papillon}
\tag{D9}
\lower12pt\hbox{$\includegraphics[width=2.6cm]{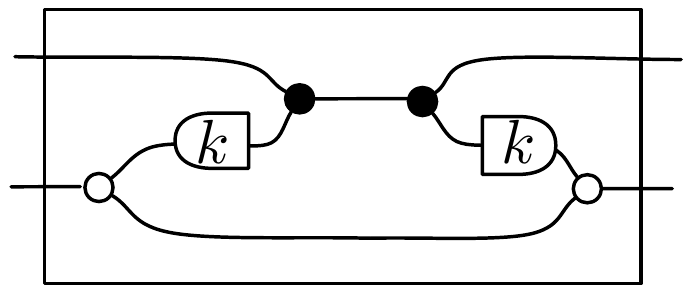}$}
=
\lower12pt\hbox{$\includegraphics[height=1.1cm]{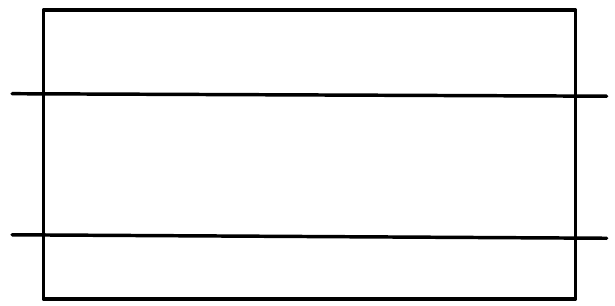}$}
 =
\lower12pt\hbox{$\includegraphics[width=2.4cm]{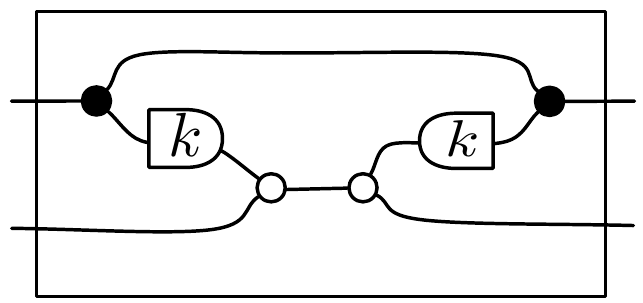}$}
\end{equation}
\begin{multicols}{2}\noindent
\begin{align}
\label{eq:wunitcancelbcomult}
\tag{D10}
\lower9pt\hbox{$\includegraphics[height=.9cm]{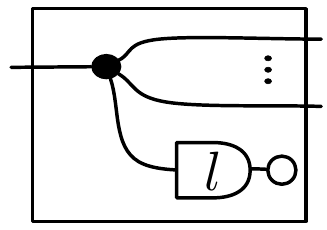}$}
=
\lower7pt\hbox{$\includegraphics[height=.7cm]{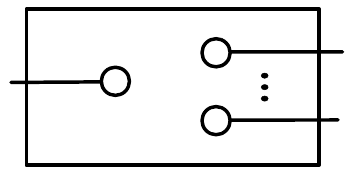}$}
\end{align}
\begin{align}
\label{eq:Bsep}
\tag{D11}
\lower6pt\hbox{$\includegraphics[height=.7cm]{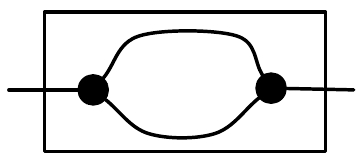}$}
=
\lower6pt\hbox{$\includegraphics[height=.7cm]{graffles/idcircuit.pdf}$}
\end{align}
\end{multicols}
Equation~\eqref{eq:uniqueantipode} states that the antipodes of $\ABR$ and $\ABRop$ coincide in $\IBRw$, which allows us to use the same notation $\antipodesquare$ for the two of them. Also observe that, because of \eqref{eq:BFrob} and \eqref{eq:Bsep}, the black structure in $\IBRw$ forms a special Frobenius algebra --- \emph{cf.} Ex.\ref{ex:distrlawsyntactic}\ref{ex:distrlawsyntactic2}. The white structure, by \eqref{eq:WFrob}, also forms a Frobenius algebra that however is not special, that is, $\WSep\!\! \neq \idcircuit$. The situation is dual (the white Frobenius algebra is special, the black is not) for the theory $\IBRb$ that we will investigate in \S~\ref{sec:IBRbCospan}. Interestingly, besides the two monochromatic Frobenius algebras, a tweak~\eqref{eq:QFrob} of the Frobenius law also holds between the white and the black structure.

Moreover, we remark that \eqref{eq:wbone} can be actually derived by the other axioms of $\IBRw$ plus \eqref{eq:Bsep} --- see Appendix~\ref{AppDerLawsIH} --- meaning that one could axiomatise $\IBRw$ with \eqref{eq:Bsep} in place of \eqref{eq:wbone}.

\subsection{Compact Closed Structure of $\IBRw$}\label{sec:cc}


An important property of the PROP $\IBRw$ is that it enjoys a (self-dual) \emph{compact closed structure}~\cite{kelly1980compactclosed}. Before going into the formal details, we illustrate the general idea. Defining a compact closed structure on a monoidal category $\catC$ requires to associate to any object $x$ a \emph{dual} object $\coc{x}$ and then define maps $x \tns \coc{x} \to I$ and $I \to x \tns \coc{x}$, where $I$ is the unit of $\tns$ in $\catC$. These maps are typically called \emph{cups} and \emph{caps} respectively, because of their graphical representation (below). For the case of $\IBRw$, $I$ will be $0$ and we let each object $n$ be its own dual. Cups and caps will be then special diagrams of type $n+n \to 0$ and $0 \to n+n$ respectively:
 \begin{equation*}
   \lower12pt\hbox{$\includegraphics[height=1.2cm]{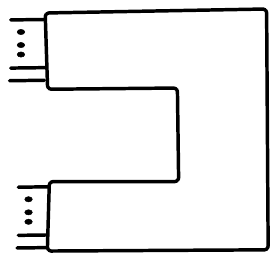}$}
  \qquad \qquad \qquad \qquad \qquad \lower12pt\hbox{$\includegraphics[height=1.2cm]{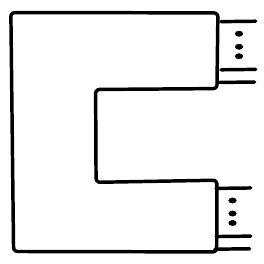}$}.
 \end{equation*}
Caps and cups behave as bent identities: they need to obey the following yanking axiom:
  \begin{equation}\label{eq:yankingInt}
  \lower16pt\hbox{$\includegraphics[height=1.6cm]{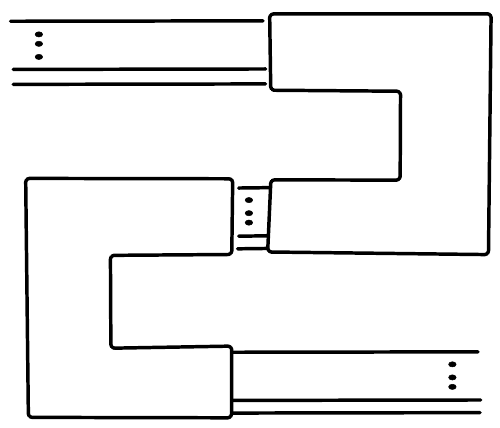}$}
  \ =\  \lower3pt\hbox{$\includegraphics[height=.3cm]{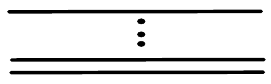}$}
  \ =\  \lower16pt\hbox{$\includegraphics[height=1.6cm]{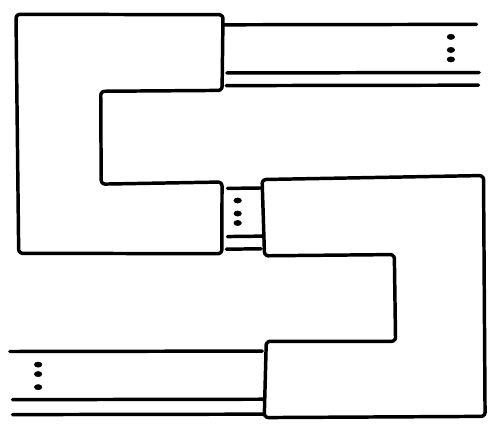}$}.
   \end{equation}
  A useful consequence is that any diagram $c$ of the appropriate type can be moved along them:
\begin{multicols}{2}
\noindent
  \begin{equation*}
  \lower12pt\hbox{$\includegraphics[height=1.2cm]{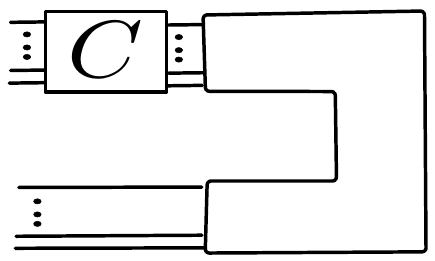}$}
  \ =\  \lower12pt\hbox{$\includegraphics[height=1.2cm]{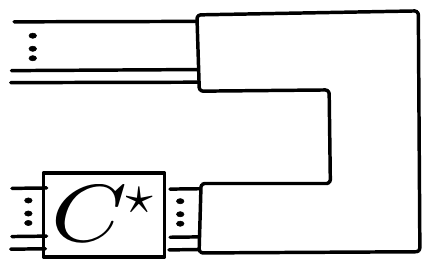}$}
   \end{equation*}
  \begin{equation}\label{eq:SlidingCCInt}
  \lower12pt\hbox{$\includegraphics[height=1.2cm]{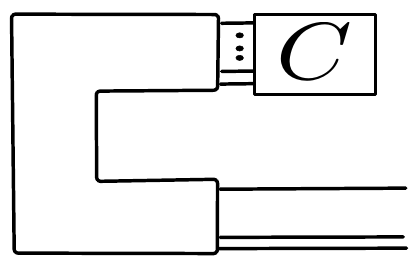}$}
  \ =\  \lower12pt\hbox{$\includegraphics[height=1.2cm]{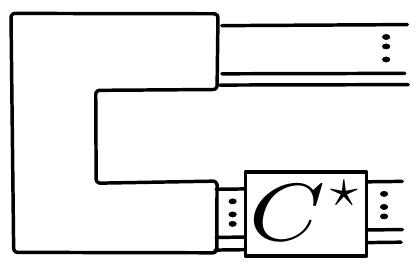}$}
   \end{equation}
 \end{multicols}
This movement turns $c$ into $\coc{c}$, which one should think of as $c$ reflected about the $y$-axis --- see Proposition~\ref{prop:star=refl} below.

 \medskip
 We now give the formal definition of the compact closed structure of $\IBRw$ and explore its properties. For constructing cups and caps we shall use components of the kind $\lccB$ and $\rccB$. For each $n$, we formally define cup and cap as diagrams $\eta_n \: 0 \to n + n$ and $\epsilon_n \: n+n \to 0$ respectively, given by the following induction on $n$:
\begin{multicols}{2}
\noindent
 \begin{eqnarray*}
  \hspace{-2cm}
   \lower18pt\hbox{$
   \alpha_0 \: 2 \to 2 \df  
 \lower5pt\hbox{$\includegraphics[height=.6cm]{graffles/symmetryalt.pdf}$}
   $}
 \end{eqnarray*}
  \begin{eqnarray*}
   \hspace{-2cm}
   \alpha_{n+1} \: 2(n+1) \to 2(n+1) \df \lower22pt\hbox{$\includegraphics[height=1.8cm]{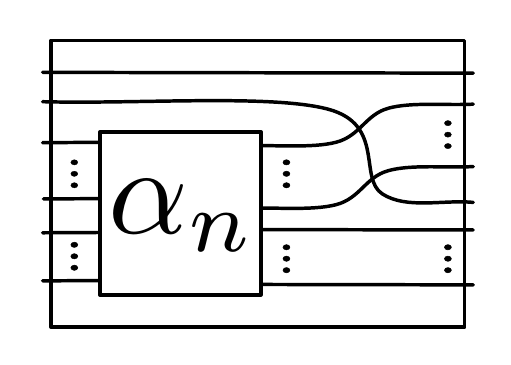}$}
 \end{eqnarray*}
 \end{multicols}
 \begin{multicols}{2}
\noindent
  \begin{eqnarray*}
  \hspace{-1.6cm}
       \lower12pt\hbox{$
   \eta_0 \: 0 \to 0 \df \lower5pt\hbox{$\includegraphics[height=.6cm]{graffles/idzerocircuit.pdf}$}
   $}
   \end{eqnarray*}
  \begin{eqnarray*}
  \hspace{-2cm}
   \eta_{n+1} \: 0 \to 2(n+1)   \df  \lower22pt\hbox{$\includegraphics[height=1.6cm]{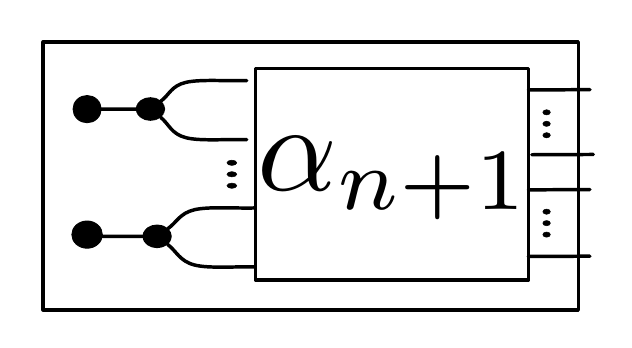}$}
   \end{eqnarray*}
   \end{multicols}
 \begin{multicols}{2}
\noindent
  \begin{eqnarray*}
  \hspace{-2cm}
     \lower12pt\hbox{$
   \beta_0 \: 2 \to 2  \df  
    \lower5pt\hbox{$\includegraphics[height=.6cm]{graffles/symmetryalt.pdf}$}
   $}
    \end{eqnarray*}
  \begin{eqnarray*}
  \hspace{-2cm}
   \beta_{n+1} \: 2(n+1) \to 2(n+1)  \df \lower22pt\hbox{$\includegraphics[height=1.6cm]{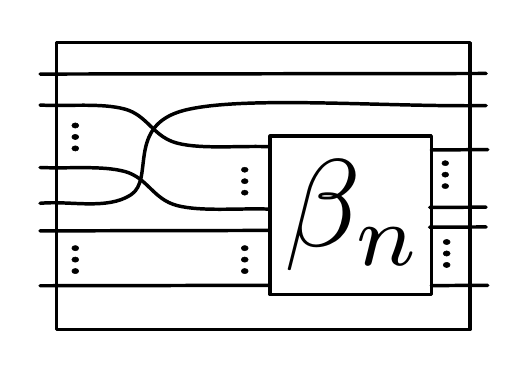}$} \end{eqnarray*}
   \end{multicols}
 \begin{multicols}{2}
\noindent
  \begin{eqnarray*}
  \hspace{-1.5cm}
       \lower18pt\hbox{$
   \epsilon_0 \: 0 \to 0   \df  \lower5pt\hbox{$\includegraphics[height=.6cm]{graffles/idzerocircuit.pdf}$}
   $}
    \end{eqnarray*}
  \begin{eqnarray*}
  \hspace{-2cm}
   \epsilon_{n+1} \: 2(n+1) \to 0   \df  \lower22pt\hbox{$\includegraphics[height=1.8cm]{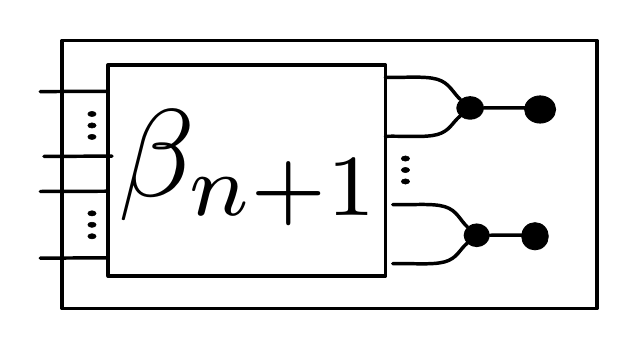}$}
   \end{eqnarray*}
 \end{multicols}

\noindent For a more concrete grip on that definition, we show the first values of $\eta_n$:
  \begin{multicols}{3}
\noindent
  \begin{eqnarray*}
       \lower20pt\hbox{$
   \eta_1 = \lower11pt\hbox{$\includegraphics[height=.9cm]{graffles/lccl.pdf}$}
   $}
    \end{eqnarray*}
  \begin{eqnarray*}
   \eta_2 =  \lower20pt\hbox{$\includegraphics[height=1.7cm]{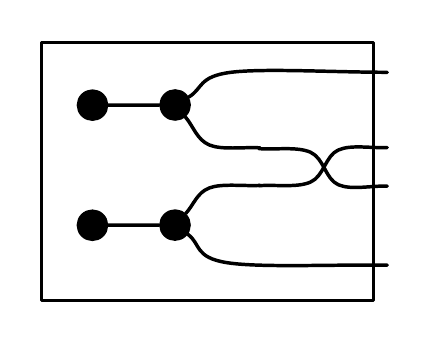}$}
   \end{eqnarray*}
     \begin{eqnarray*}
   \eta_3 =  \lower22pt\hbox{$\includegraphics[height=1.8cm]{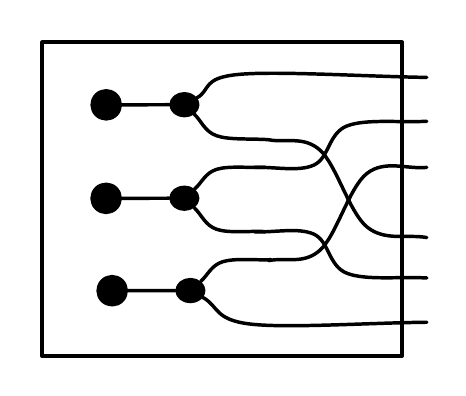}$}
   \end{eqnarray*}
 \end{multicols}

 For the sequel, we fix notation \lccn for $\eta_n$ and \rccn for $\epsilon_n$. Also, we let \idncircuit be notation for $\id_n$. Similarly, we write \wcounitn (respectively, \bcounitn) for the $n$-fold $\tns$-product of \Wcounit\ (respectively, \Bcounit) and $\nscalar$ for the $n$-fold $\tns$-product of $\scalar$.

 The next statement proves that the $\eta$s and the $\epsilon$s satisfy the yanking axiom~\eqref{eq:yankingInt} and thus yield a compact closed structure.

\begin{proposition}\label{prop:snakecc} $\IBRw$ is self-dual compact closed with structure given by $\eta_n$ and $\epsilon_n$ for each object $n \in \IBRw$. \end{proposition}
\begin{proof} It suffices to verify the following equality, for each $n \in \IBRw$.
  \begin{equation}
  \label{eq:gensnake}
  \tag{CC1}
  \lower18pt\hbox{$\includegraphics[height=1.4cm]{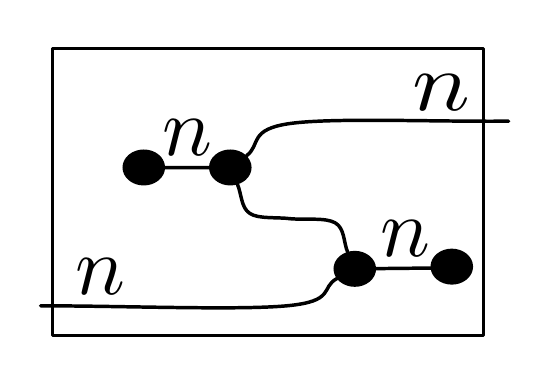}$}
  = \lower12pt\hbox{$\includegraphics[height=1cm]{graffles/idncircuit.pdf}$}
  = \lower18pt\hbox{$\includegraphics[height=1.4cm]{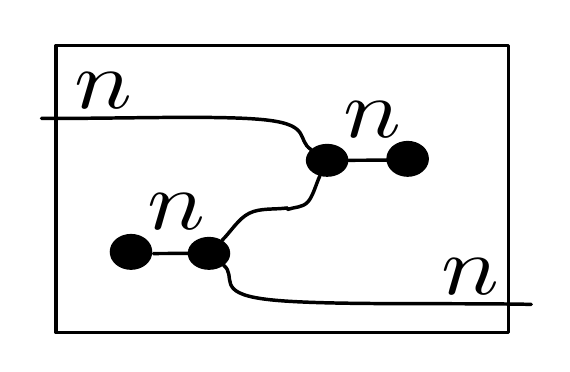}$}
   \end{equation}
   The details of this derivation in $\IBRw$ can be found in Appendix~\ref{AppCC}.\end{proof}

We now formalise and prove the sliding equations~\eqref{eq:SlidingCCInt}. As observed in \cite[Rmk.~2.1]{Selinger07DaggerCC}, the compact closed structure allows to define a contravariant PROP morphism $\coc{\cdot} \: \IBRw \to \IBRw$ as:
\begin{align}\label{eq:defstar}
 \lower12pt\hbox{$\includegraphics[height=1.2cm]{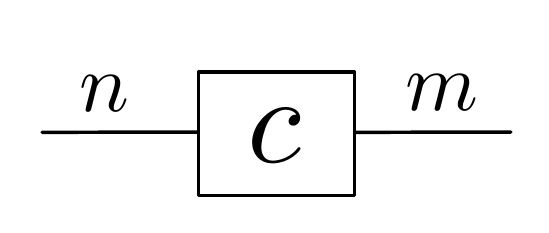}$} \mapsto \lower12pt\hbox{$\includegraphics[height=1.2cm]{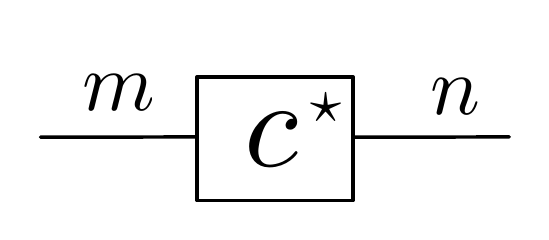}$} \df \lower18pt\hbox{$\includegraphics[height=1.7cm]{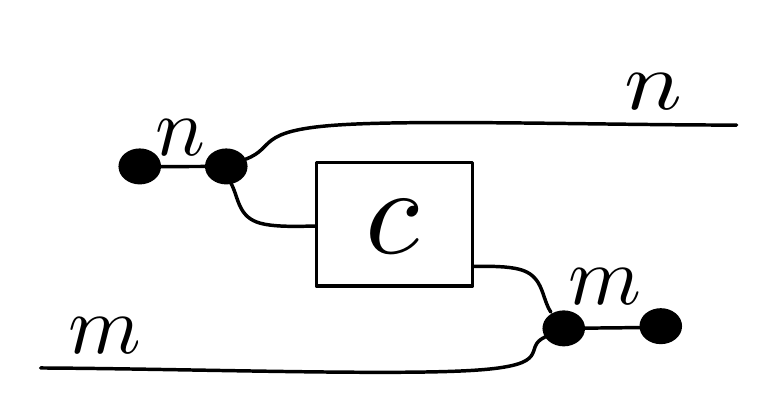}$}.
 \end{align}
 \begin{corollary} For any string diagram $c \: n \to m$ of $\IBRw$,
 \begin{multicols}{2}\noindent
  \begin{equation}
  \hspace{-.1cm}
  \label{eq:ccsliding}
  \tag{CC2}
  \lower10pt\hbox{$\includegraphics[height=1.1cm]{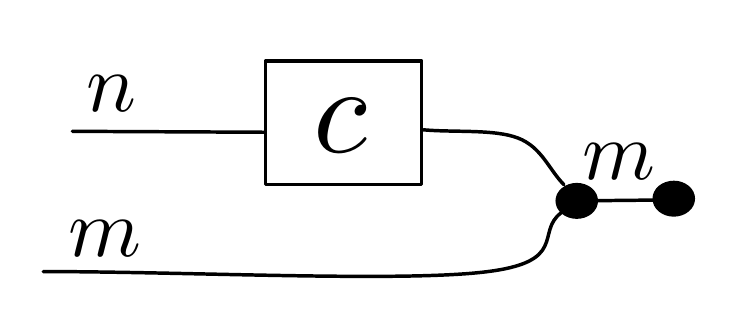}$}
   \! =  \! \lower13pt\hbox{$\includegraphics[height=1.1cm]{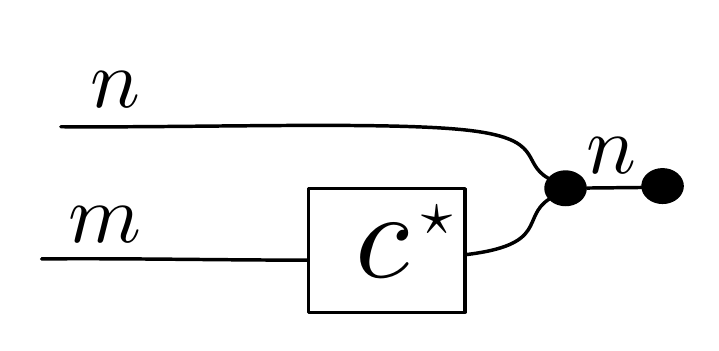}$}
   \end{equation}
     \begin{equation}
  \label{eq:ccsliding2}
  \tag{CC3}
  \lower8pt\hbox{$\includegraphics[height=1.1cm]{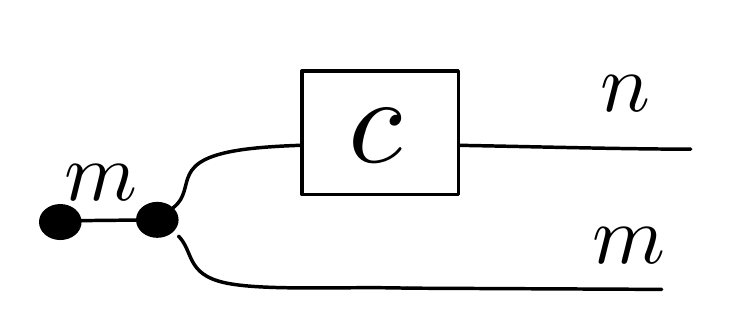}$}
  \! = \! \lower10pt\hbox{$\includegraphics[height=1.1cm]{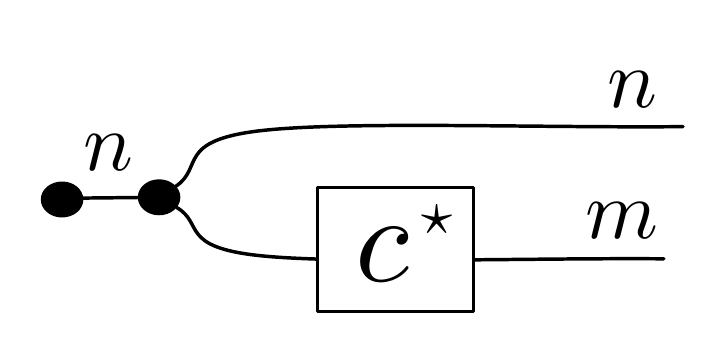}$}
   \end{equation}
   \end{multicols}
 \end{corollary}
 \begin{proof} The following is the derivation of \eqref{eq:ccsliding} in $\IBRw$. The one of \eqref{eq:ccsliding2} is analogous.
  \begin{equation*}
  \lower12pt\hbox{$\includegraphics[height=1.3cm]{graffles/ccslidingr.pdf}$}
  \eql{Def. $\coc{c}$} \lower22pt\hbox{$\includegraphics[height=1.7cm]{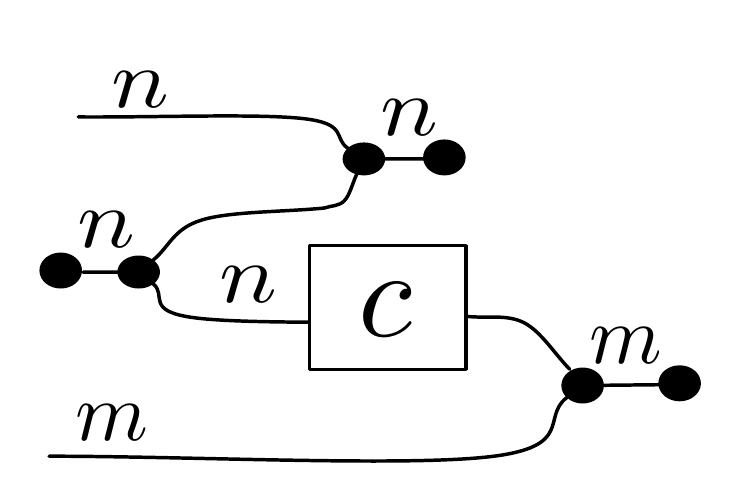}$}
  \eql{\eqref{eq:gensnake}} \lower12pt\hbox{$\includegraphics[height=1.3cm]{graffles/ccslidingl.pdf}$}
   \end{equation*}
 \end{proof}

The following proposition ensures that the notation $\coc{\cdot}$ used above actually does not conflict with the one used for the contravariant PROP morphism $\ABR \to \ABRop$ defined in \S~\ref{sec:theorymatr} --- that means, $\coc{\HAToIHw(c)} = \HAToIHw(\coc{c})$, where $\HAToIHw \: \ABR \to \IBRw$ is the inclusion of $\ABR$ in $\IBRw$, see Definition~\ref{def:IBRw}. First, let $\refl{\cdot} \: \IBRw \to \IBRw$ be the contravariant PROP morphism given inductively as follows:
  \begin{multicols}{4}
\noindent
      \begin{eqnarray*}
     \Bcounit \mapsto \Bunit
    \end{eqnarray*}
   \begin{eqnarray*}
     \Bunit \mapsto \Bcounit
    \end{eqnarray*}
  \begin{eqnarray*}
     \Wunit \mapsto \Wcounit
    \end{eqnarray*}
  \begin{eqnarray*}
   \Wcounit \mapsto \Wunit
   \end{eqnarray*}
    \end{multicols}
    \smallskip
     \begin{multicols}{4}
     \noindent
     \begin{eqnarray*}
     \Wmult \mapsto \Wcomult
    \end{eqnarray*}
  \begin{eqnarray*}
   \Wcomult \mapsto \Wmult
   \end{eqnarray*}
        \begin{eqnarray*}
     \Bmult \mapsto \Bcomult
    \end{eqnarray*}
  \begin{eqnarray*}
   \Bcomult \mapsto \Bmult
   \end{eqnarray*}
      \end{multicols}
      \smallskip
     \begin{multicols}{2}
     \noindent
           \begin{eqnarray*}
     \scalar \mapsto \coscalar
    \end{eqnarray*}
  \begin{eqnarray*}
   \coscalar \mapsto \scalar
   \end{eqnarray*}
    \end{multicols}
    \smallskip
     \begin{multicols}{2}
     \noindent
     \begin{equation*}    \lower10pt\hbox{$\includegraphics[height=.9cm]{graffles/reflcompl.pdf}$}
  \!\mapsto\! \lower10pt\hbox{$\includegraphics[height=.9cm]{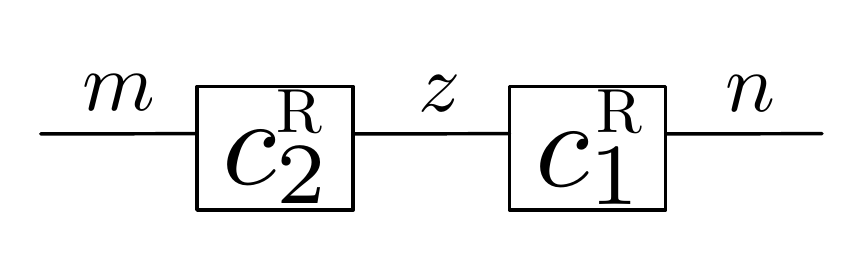}$}
   \end{equation*}
        \begin{equation*}
         \lower20pt\hbox{$\includegraphics[height=1.5cm]{graffles/refltnsl.pdf}$}
   \!\mapsto\!\lower20pt\hbox{$\includegraphics[height=1.5cm]{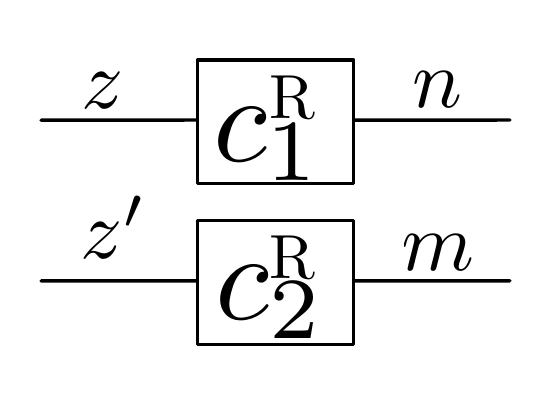}$}
   \end{equation*}
 \end{multicols}

\begin{proposition}\label{prop:star=refl} $\coc{c} = \refl{c}$ for all string diagrams $c \: n \to m$ of $\IBRw$.
\end{proposition}
\begin{proof} The proof is by induction on $c$. See Appendix~\ref{AppCC} for the details. \end{proof}

By virtue of Proposition~\ref{prop:star=refl}, one can think of \eqref{eq:ccsliding}-\eqref{eq:ccsliding2} as saying that any diagram $c$ can be moved along cups and caps, with the effect of flipping it horizontally: this agrees with the intuitive picture provided in~\eqref{eq:SlidingCCInt}.

%


\subsection{$\IBRw$: the theory of Spans of $\PID$-matrices}\label{sec:completeness}

In this section we prove that $\IBRw$ characterises spans of matrices. Towards this result, we need to show that the PROP $\SpanMat$ of spans of $\PID$-matrices --- here $\PJ$ is the core of $\Mat{\PID}$, i.e., the PROP of invertible $\PID$-matrices --- can be actually constructed, that means, $\Mat{\PID}$ has pullbacks. For this purpose, we shall first illustrate how finite limits and colimits work in $\Mat{\PID}$.

\subsubsection{(Co)limits of $\PID$-matrices}\label{sec:colimitsMat}

For computing (co)limits in $\Mat{\PID}$ it is useful to fix the following categories:
\begin{itemize}[noitemsep,topsep=0pt,parsep=0pt,partopsep=0pt]
  \item the abelian category $\RMod{\PID}$ of $\PID$-modules and linear maps;
  \item its full subcategory $\FRMod{\PID}$ consisting of the finitely-generated free $\PID$-modules and linear maps between them.
\end{itemize}
There is an equivalence of categories between $\FRMod{\PID}$ and $\Mat{\PID}$: a finitely-generated free $\PID$-module in $\FRMod{\PID}$, say of dimension $n$, is isomorphic to $\PID^n$ and thus we can associate it with the object $n$ in $\VectR$. A linear map $f \: V \to W$ in $\FRMod{\PID}$ is represented by a matrix $M \: n \to m$, where $V \cong \PID^n$ and $W \cong \PID^m$.

The equivalence with $\FRMod{\PID}$ guarantees various properties of $\Mat{\PID}$. An important one is the existence of biproducts: the object $n+m$ is both a coproduct and a product of $n$ and $m$. It is convenient to fix notation for the matrices given by the two universal properties:
\[
\xymatrix@C=40pt{
& z  & \\
n \ar[ur]^{A} &\ar[l] n + m \ar[r]\ar@{-->}[u]_<<<{(A \mid B)} & m \ar[ul]_{B}
} \qquad \qquad
\xymatrix@C=40pt{
& \ar[dr]^{D} r \ar[dl]_{C} \ar@{-->}[d]^{(\frac{C}{D})} & \\
n \ar[r] & n + m & m \ar[l]
}
\]
The notation reflects the way in which these matrices are constructed, by putting $A$ and $B$ side-by-side and $C$ above $D$.

We now turn to the question of pullbacks in $\Mat{\PID}$. In the following, we use the notation $\Ker{A}$ for the matrix representing the kernel of $A \in \Mat{\PID}[n,z]$ --- or, more precisely, the categorical kernel of $A$ in the abelian category $\RMod{\PID}$, obtainable via the pullback:
\begin{equation*}
\vcenter{
\xymatrix@R=8pt@C=10pt{
&\ar[dl]_{\Ker{A}} \PID^r \pushoutcorner \ar[dr]^{\finVect} & \\
\PID^n \ar[dr]_{A} & & \PID^0 \ar[dl]^{\initVect}\\
& \PID^z  & }
}
\end{equation*}
 Our assumption that $\PID$ is a PID becomes now crucial to show that $\Mat{\PID}$ has pullbacks, which can be obtained by kernel computation in $\RMod{\PID}$.
\begin{proposition}\label{prop:PullbacksMatR} Provided that $\PID$ is a principal ideal domain, $\Mat{\PID}$ has pullbacks. \end{proposition}
\begin{proof}
A cospan $n \tr{A} z \tl{B} m$ in $\VectR \simeq \FRMod{\PID}$ can be interpreted in $\RMod{\PID}$ as a cospan $\PID^n \tr{A} \PID^z \tl{B} \PID^m$. Since $\RMod{\PID}$ is an abelian category, one can form its pullback starting from the kernel $\KerAB \: V \to \PID^n \tns \PID^m$. Here $V$ is the submodule of $\PID^n \tns \PID^m$ of pairs $(\vv,\ww)$ such that $A\vv \minusB \ww = \zerov$, i.e. $A \vv = B \ww$. The fact that $\PID$ is a PID guarantees that every submodule of a free module is itself free (see e.g.~{\cite[\S 23]{HandbookLinearAlgebra}}). Therefore, $V$ is a free module, meaning that $V \cong \PID^r$ for some natural number $r \leq n+m$. We can then express the pullback of $A$ and $B$ in $\RMod{\PID}$ as follows:
\[\xymatrix@=40pt{
&\ar[dl]_{C} \PID^r \ar[d]|{\Ker{A \, \mid \, {\scalebox{0.55}[1.0]{\( - \)}\! B}}} \ar[dr]^{D} & \\
\PID^n \ar[dr]_{A}  & \PID^n \tns \PID^m \ar[d]|{(A \, \mid \, {\scalebox{0.55}[1.0]{\( - \)}\! B})} & \PID^m \ar[dl]^{B} \\
& \PID^z  & }\]
Since $\KerAB$ ranges over $\PID^n \tns \PID^m$, it is of shape $(\frac{C}{D}) \: \PID^r \to \PID^n \tns \PID^m$ and postcomposition with the product projections $\pi_1 \: n \tns m \to n$ and $\pi_2 \: n \tns m \to m$ yields matrices $C$ and $D$ as in the diagram above. As a consequence, the following is a pullback square in $\VectR$:
\[\xymatrix@C=10pt@R=10pt{
&\ar[dl]_{C} r \pushoutcorner \ar[dr]^{D} & \\
n \ar[dr]_{A}  & & m \ar[dl]^{B} \\
& z  & }.\]
\end{proof}

\begin{remark} It is worth mentioning that the same reasoning does not apply to pushouts: they exist in $\VectR$ for purely formal reasons, being the category self-dual, but cannot be generally calculated as in $\RMod{\PID}$. This asymmetry actually plays a role in our developments: we shall return to it in Remark~\ref{rmk:PushoutsMatrices}.
\end{remark}

By Proposition~\ref{prop:distrLawPbPo} and~\ref{prop:PullbacksMatR} there is a distributive law $\lambdapb \: \CospanMat \to \SpanMat$ defined by pullback in $\Mat{\PID}$, with $\PJ$ being the core of $\Mat{\PID}$. The PROP $\SpanMat$ has arrows the spans of $\PID$-matrices, made equal whenever there is a span isomorphism between them. Composition is by pullback in $\Mat{\PID}$.

\subsubsection{Soundness and Completeness of $\IBRw$}

The rest of the section will be devoted to showing that the theory $\IBRw$ of interacting Hopf algebras is a presentation by generators and relations of $\SpanMat$.

\begin{theorem}\label{th:Span=IBw} $\IBRw \cong \SpanMat$. \end{theorem}

Towards a proof of Theorem~\ref{th:Span=IBw}, first observe that, by Proposition~\ref{prop:ab=vect}, $\SpanMat \cong \ABRop \bicomp{\PJ} \ABR$. This correspondence gives us a recipe for the SMT presenting $\SpanMat$: by Proposition~\ref{prop:SMTforComposition}, it will have generators those of $\ABR + \ABRop$, and equations those of $\ABR + \ABRop$ plus the set $E_{\lambdapb}$ of equations induced by the distributive law $\lambdapb \: \CospanMat \to \SpanMat$.

Now, by definition $\IBRw$ is generated by the same signature as $\ABR + \ABRop$. For the equations, we need to show:
\begin{description}[noitemsep,topsep=2pt,parsep=2pt,partopsep=2pt]
\item[Soundness] all the equations of $\IBRw$ are provable from the equations of $\ABR + \ABRop$ plus $E_{\lambdapb}$;
\item[Completeness] all the equations of $\ABR + \ABRop$ and in $E_{\lambdapb}$ are provable in $\IBRw$.
\end{description}

\begin{remark}\label{rmk:readoffequationsMatrix} For a better grip on the two statements, let us first describe the set $E_{\lambdapb}$. By definition of $\lambdapb$, it will consists of those equations that can be read off by pullback squares in $\Mat{\PID}$. The reading process is completely analogous to the one illustrated for pullbacks in $\F$ in Example~\ref{ex:distrlawsyntactic}\ref{ex:distrlawsyntactic1}. We pick any pair $\tr{p \in \ABR}\tr{\coc{q} \in \ABRop}$ in the source of $\lambdapb$. This yields a cospan $\tr{p}\tl{q}$ in $\ABR$. We compute its pullback in $\Mat{\PID}$:
\begin{equation}\label{eq:exReadOffEquationMat}
\vcenter{
\xymatrix@R=15pt@C=15pt{
&\ar[dl]_{\sem{\ABR}(f)}  \pushoutcorner \ar[dr]^{\sem{\ABR}(g)} & \\
 \ar[dr]_{\sem{\ABR}(p)} & &  \ar[dl]^{\sem{\ABR}(q)}\\
&   & }
}
\end{equation}
The resulting span $\tl{f}\tr{g}$ in $\ABR$ yields a pair $\tr{\coc{f} \in \ABRop}\tr{g \in \ABR}$, which $\lambdapb$ associates to $\tr{p \in \ABR}\tr{\coc{q} \in \ABRop}$. Thus the equation $p \poi \coc{q} \feq \coc{f} \poi g$ will be part of $E_{\lambdapb}$.
\end{remark}

Now, turning to the soundness claim, by definition of $\IBRw$ it suffices to show the statement for~\eqref{eq:lcm}-\eqref{eq:BccscalarAxiomTwo}. We observe that each of those equations is of the shape $p \poi \coc{q} \feq \coc{f}\poi g$, with $\tl{f}\tr{g}$ the pullback of $\tr{p}\tl{q}$ in $\ABR$. Following our discussion in Remark~\ref{rmk:readoffequationsMatrix}, this implies that~\eqref{eq:lcm}-\eqref{eq:BccscalarAxiomTwo} are all in $E_{\lambdapb}$, thus proving soundness of $\IBRw$.

\begin{example} It is useful to give some representative example of the soundness argument.
\begin{itemize}
\item Consider the two pullback squares in $\ABR$ on the left below. The pullback computation is carried out in $\VectR$ as on the right below.
\begin{eqnarray}\label{eq:exWFrob1}
\vcenter{
    \xymatrix@R=12pt@C=15pt{
    & \ar[dl]_{\WFrobL} 3 \pushoutcorner \ar[dr]^{\WFrobR} & \\
    2 \ar[dr]_{\Wmult} && 2 \ar[dl]^{\Wmult} \\
    & 1 &
    }
}
    & \qquad \raise30pt\hbox{${\xymatrix{\\ \ar@{|->}[r]^{\sem{\ABR}} &}}$}\qquad &
\vcenter{
    \xymatrix@R=12pt@C=15pt{
    & \ar[dl]_{\tiny{\left(%
                   \begin{array}{ccc}
                     \!\! 1 \!\! & \!\! 1 \!\! & \!\! 0\!\!\\
                     \!\! 0\!\! &\!\! 0\!\! & \!\! 1 \!\!
                    \end{array}\right)}} 3 \pushoutcorner \ar[dr]^{\tiny{\left(%
                   \begin{array}{ccc}
                     \!\! 1 \!\! & \!\! 0 \!\! & \!\! 0\!\!\\
                     \!\! 0\!\! &\!\! 1\!\! & \!\! 1 \!\!
                    \end{array}\right)}} & \\
    2 \ar[dr]_{\tiny{\left(%
                     \begin{array}{cc}
                    \!\!  1 \!\! &\! 1 \!\!
                    \end{array}\right)}} && 2 \ar[dl]^{\tiny{\left(%
                     \begin{array}{cc}
                    \!\!  1 \!\! &\! 1 \!\!
                    \end{array}\right)}} \\
    & 1 &
    }
}
 \\ \label{eq:exWFrob2}
\vcenter{
    \xymatrix@R=12pt@C=15pt{
    & \ar[dl]_{\WFrobR} 3 \pushoutcorner \ar[dr]^{\WFrobL} & \\
    2 \ar[dr]_{\Wmult} && 2 \ar[dl]^{\Wmult} \\
    & 1 &
    }
}
    & \qquad \raise30pt\hbox{$\xymatrix{\\ \ar@{|->}[r]^{\sem{\ABR}} &}$}\qquad &
\vcenter{
    \xymatrix@R=12pt@C=15pt{
    & \ar[dl]_{\tiny{\left(%
                   \begin{array}{ccc}
                     \!\! 1 \!\! & \!\! 0 \!\! & \!\! 0\!\!\\
                     \!\! 0\!\! &\!\! 1\!\! & \!\! 1 \!\!
                    \end{array}\right)}} 3 \pushoutcorner \ar[dr]^{\tiny{\left(%
                   \begin{array}{ccc}
                     \!\! 1 \!\! & \!\! 1 \!\! & \!\! 0\!\!\\
                     \!\! 0\!\! &\!\! 0\!\! & \!\! 1 \!\!
                    \end{array}\right)}} & \\
    2 \ar[dr]_{\tiny{\left(%
                     \begin{array}{cc}
                    \!\!  1 \!\! &\! 1 \!\!
                    \end{array}\right)}} && 2 \ar[dl]^{\tiny{\left(%
                     \begin{array}{cc}
                    \!\!  1 \!\! &\! 1 \!\!
                    \end{array}\right)}} \\
    & 1 &
    }
}
\end{eqnarray}
 The pullbacks in $\ABR$ above describe the two equations of axiom~\eqref{eq:WFrob}. The left-hand equation arises from~\eqref{eq:exWFrob1}:
 \begin{center}
 the cospan $2\tr{\Wmult}1\tl{\Wmult}2$ yields the string diagram $\lower6pt\hbox{$\includegraphics[height=.7cm]{graffles/WFrobX.pdf}$}$ in~\eqref{eq:WFrob}\\
 its pullback span $2 \tl{\WFrobL} 3 \tr{\WFrobR} 2$ yields the string diagram $\lower8pt\hbox{$\includegraphics[height=.9cm]{graffles/WFrobS.pdf}$}$ in~\eqref{eq:WFrob}
 \end{center}
 because, as described in Remark~\ref{rmk:readoffequationsMatrix},  the right leg of the cospan and the left leg of the span stand for arrows of $\ABRop$, thus they should be read reflected about the $y$-axis, e.g. $\coc{\Wmult} = \Wcomult$. Similarly, the right-hand equation in~\eqref{eq:WFrob} arises from the second pullback~\eqref{eq:exWFrob2}.

    Note that~\eqref{eq:exWFrob1}-\eqref{eq:exWFrob2} represent two pullbacks for the same cospan. The resulting spans are indeed isomorphic, the isomorphism being the invertible matrix
    \begin{equation*}
    {\tiny \begin{pmatrix}
     1 \!&\! 1 \!&\! 0 \\
    0 \!&\! -1 \!&\! 0 \\
    0 \!&\! 1 \!&\! 1
    \end{pmatrix} }.
    \end{equation*}
\item The axiom~\eqref{eq:BFrob} corresponds to the two pullback square in $\ABR$ on the left below. The pullback computation is carried out in $\VectR$ as on the right below.
\begin{eqnarray*}
\xymatrix@R=12pt@C=15pt{
& \ar[dl]_{\Bcomult} 1 \pushoutcorner \ar[dr]^{\Bcomult} & \\
2 \ar[dr]_{\BFrobL} && 2 \ar[dl]^{\BFrobR} \\
& 3 &
}
& \qquad \xymatrix{\\ \ar@{|->}[r]^{\sem{\ABR}} &}\qquad &
\xymatrix@R=12pt@C=15pt{
& \ar[dl]_{\tiny{\left(%
                 \begin{array}{c}
                \!\!  1 \!\! \\
                \! 1 \!\!
                \end{array}\right)}} 1 \pushoutcorner \ar[dr]^{\tiny{\left(%
                 \begin{array}{c}
                \!\!  1 \!\! \\
                \! 1 \!\!
                \end{array}\right)}} & \\
2 \ar[dr]_{\tiny{\left(%
               \begin{array}{cc}
                 \!\! 1 \!\! & \!\! 0 \!\!  \\
                 \!\! 1\!\! &\!\! 0\!\! \\
                 \!\! 0\!\! & \!\! 1 \!\!
                \end{array}\right)}} && 2 \ar[dl]^{\tiny{\left(%
               \begin{array}{cc}
                 \!\! 1 \!\! & \!\! 0 \!\!  \\
                 \!\! 0\!\! &\!\! 1\!\! \\
                 \!\! 0\!\! & \!\! 1 \!\!
                \end{array}\right)}} \\
& 3 &
} \\
\xymatrix@R=12pt@C=15pt{
& \ar[dl]_{\Bcomult} 1 \pushoutcorner \ar[dr]^{\Bcomult} & \\
2 \ar[dr]_{\BFrobR} && 2 \ar[dl]^{\BFrobL} \\
& 3 &
}
& \qquad \xymatrix{\\ \ar@{|->}[r]^{\sem{\ABR}} &}\qquad &
\xymatrix@R=12pt@C=15pt{
& \ar[dl]_{\tiny{\left(%
                 \begin{array}{c}
                \!\!  1 \!\! \\
                \! 1 \!\!
                \end{array}\right)}} 1 \pushoutcorner \ar[dr]^{\tiny{\left(%
                 \begin{array}{c}
                \!\!  1 \!\! \\
                \! 1 \!\!
                \end{array}\right)}} & \\
2 \ar[dr]_{\tiny{\left(%
               \begin{array}{cc}
                 \!\! 1 \!\! & \!\! 0 \!\!  \\
                 \!\! 0\!\! &\!\! 1\!\! \\
                 \!\! 0\!\! & \!\! 1 \!\!
                \end{array}\right)}} && 2 \ar[dl]^{\tiny{\left(%
               \begin{array}{cc}
                 \!\! 1 \!\! & \!\! 0 \!\!  \\
                 \!\! 1\!\! &\!\! 0\!\! \\
                 \!\! 0\!\! & \!\! 1 \!\!
                \end{array}\right)}} \\
& 3 &
}
\end{eqnarray*}
    Exhibiting the associated pullbacks allows us to notice that the black Frobenius equation arises in a different way from the white one: the directionality given by the distributive law $\lambdapb$ is $\cdot \Rightarrow \cdot \Leftarrow \cdot$ in~\eqref{eq:BFrob} and $\cdot \Leftarrow \cdot \Rightarrow \cdot$ in~\eqref{eq:WFrob}.
\item The axiom~\eqref{eq:lcc} corresponds to the pullback square in $\ABR$ on the left below. The pullback computation is carried out in $\VectR$ as on the right below.
\begin{eqnarray*}
\xymatrix@R=12pt@C=15pt{
& \ar[dl]_{\Bcounit} 1 \pushoutcorner \ar[dr]^{\BcomultSingleAntipode} & \\
0 \ar[dr]_{\Wunit} && 2 \ar[dl]^{\Wmult} \\
& 1 &
}
& \qquad \xymatrix{\\ \ar@{|->}[r]^{\sem{\ABR}} &}\qquad &
\xymatrix@R=12pt@C=15pt{
& \ar[dl]_{\finVect} 1 \pushoutcorner \ar[dr]^{\tiny{\left(%
               \begin{array}{c}
                 \!\! 1 \!\!\\
                 \!\! -1\!\!
                \end{array}\right)}} & \\
0 \ar[dr]_{\initVect} && 2 \ar[dl]^{\tiny{\left(%
                 \begin{array}{cc}
                \!\!  1 \!\! &\! 1 \!\!
                \end{array}\right)}} \\
& 1 &
}
\end{eqnarray*}
\end{itemize}
\end{example}

While soundness does not pose particular problems, the real challenge is showing completeness:
 we need to verify that any pullback in $\VectR$ yields an equation which is provable in $\IBRw$.

\begin{proposition} \label{prop:IBwComplete} For any pullback square in $\VectR$ (as on the left), the corresponding diagrammatic equation (on the right) is derivable in $\IBRw$.
\[\xymatrix@R=10pt@C=10pt{
&\ar[dl]_{C} r \pushoutcorner \ar[dr]^{D} & \\
n \ar[dr]_{A} & & m \ar[dl]^{B}\\
& z  & }
\qquad \qquad
\lower20pt\hbox{
\lower8pt\hbox{$\includegraphics[height=.8cm]{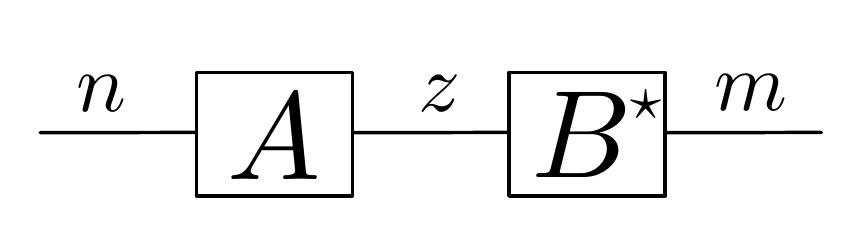}$} =
\lower8pt\hbox{$\includegraphics[height=.8cm]{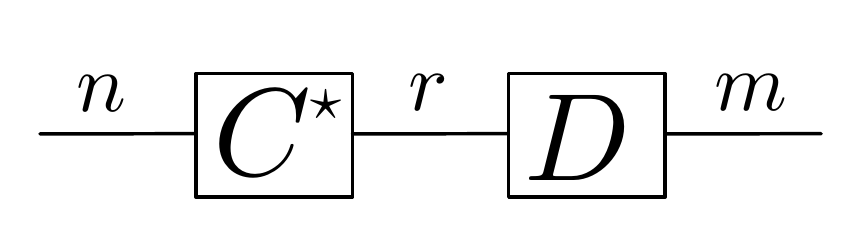}$}
}
\]
\end{proposition}

Note that, for the same reasons as for pullbacks in $\F$ (see Remark~\ref{rmk:pullbacksnoatom}), we cannot prove the completeness statement by relying on a notion of atomic diagram, as we did for cospans in $\PF$ (Proposition~\ref{th:PFROBcomplete}). We shall instead give an argument based on the way pullbacks are canonically computed in $\Mat{\PID}$: this is developed in the remaining of this section.

By the above discussion, the proof of Theorem~\ref{th:Span=IBw} is completed by showing Proposition~\ref{prop:IBwComplete}. 

\subsubsection{Invertible Matrices, Graphically}

A key step towards a proof of Proposition~\ref{prop:IBwComplete} is to understand how to canonically represent pullback diagrams graphically, in terms of string diagrams. Computing pullbacks involves standard algebraic concepts, amongst which invertible matrices. For this reason, we now prove some basic properties of the diagrammatic representation of invertible matrices.

\begin{lemma}\label{lemma:invertiblestar} For $U \in \VectR[n,n]$ invertible, the following holds in $\IBRw$:
\begin{equation}\label{eq:invertiblestar}
\lower11pt\hbox{$\includegraphics[height=1cm]{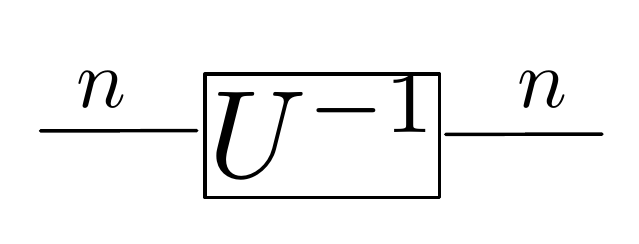}$} =
\lower11pt\hbox{$\includegraphics[height=1cm]{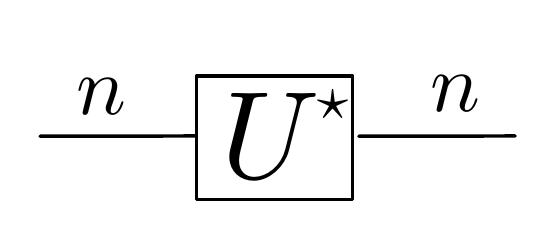}$}
\end{equation}
\end{lemma}
\begin{proof} Recall that an invertible $n \times n$ $\PID$-matrix is one obtainable from the identity $n \times n$ matrix by application of elementary row operations. Thus we can prove our statement by induction on the number of applied operations.

The base case is the one in which no row operation is applied and thus $U = \id_n$. Then we have the following equality in $\IBRw$, yielding~\eqref{eq:invertiblestar}.
 \[
\lower11pt\hbox{$\includegraphics[height=1cm]{graffles/circuitUminusone.pdf}$} =
\lower4pt\hbox{$\includegraphics[height=.8cm]{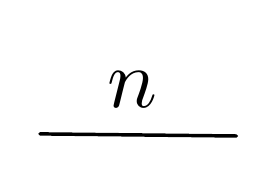}$} =
\lower11pt\hbox{$\includegraphics[height=1cm]{graffles/circuitUstar.pdf}$}
\]

Inductively, suppose that $U$ is obtained by swapping two rows of an invertible matrix $V$. We can assume without loss of generality that the two rows are one immediately above the other, with $j$ the number of rows above them and $k$ the number of rows below, where $n = j+2+k$. In diagrammatic terms, this means that
\[
\lower11pt\hbox{$\includegraphics[height=1cm]{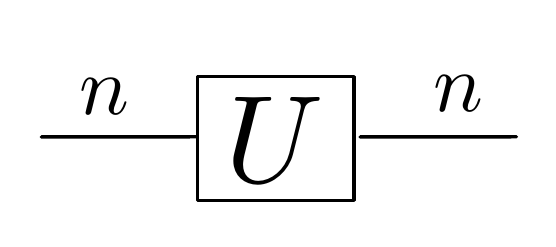}$} =
\lower11pt\hbox{$\includegraphics[height=1cm]{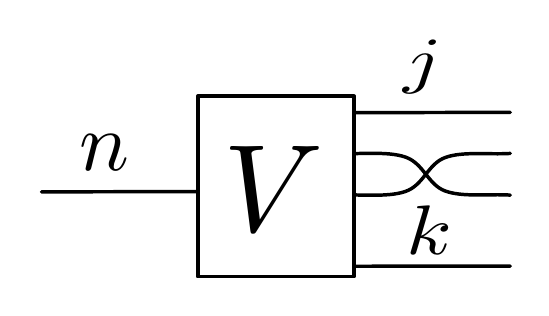}$}
\]
In order to show~\eqref{eq:invertiblestar}, it suffices to prove that the string diagram representing $\coc{U}$ is the inverse of $U$, that is, $U\poi \coc{U} = \id_n = \coc{U}\poi U$. This is given by the following derivations.
\begin{eqnarray*}
\lower14pt\hbox{$\includegraphics[height=1.3cm]{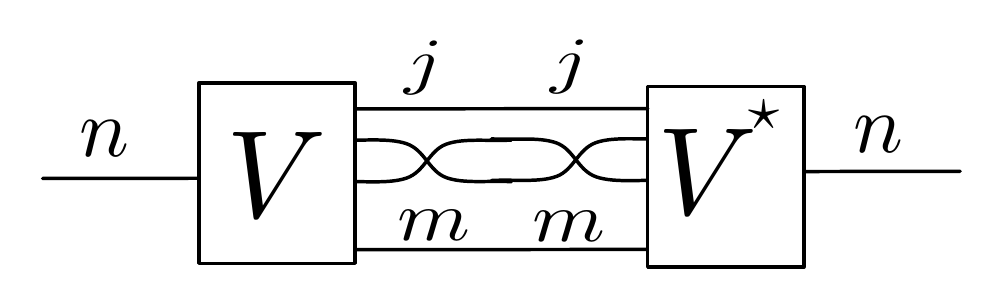}$}
\eql{\eqref{eq:SymIso}}
\lower11pt\hbox{$\includegraphics[height=1cm]{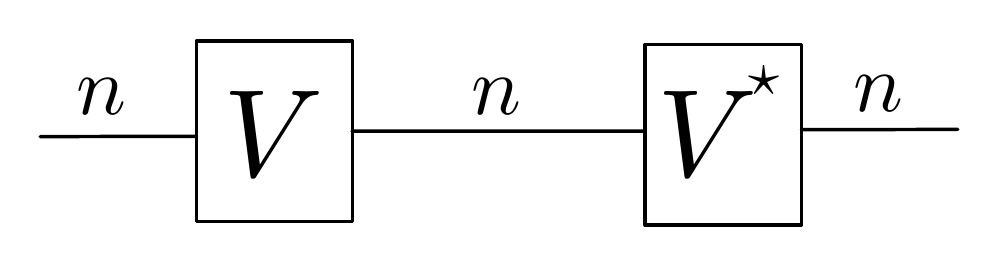}$}
\eql{Ind. hyp.}
\lower4pt\hbox{$\includegraphics[height=.8cm]{graffles/circuitidn.pdf}$} \\
\lower14pt\hbox{$\includegraphics[height=1.3cm]{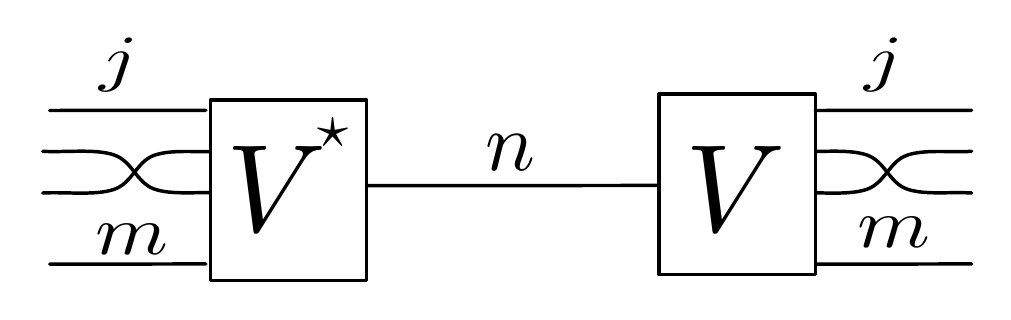}$}
\eql{Ind. hyp.}
\lower12pt\hbox{$\includegraphics[height=1.2cm]{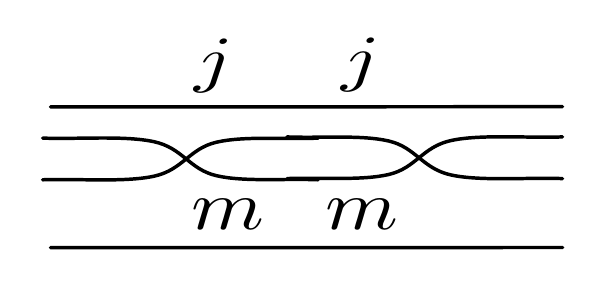}$}
\eql{\eqref{eq:SymIso}}
\lower4pt\hbox{$\includegraphics[height=.8cm]{graffles/circuitidn.pdf}$}
 \end{eqnarray*}

 The next inductive case that we consider is the one of row sum: $U$ has been obtained from an invertible matrix $V$ by replacing a row $r$ with the sum $r + k r'$, for some other row $r'$ and element $k \in \PID$.
 As above, we may assume that $r$ and $r'$ are adjacent rows of $V$. By this description, we can represent $U$ as the following string diagram, where $j+2+m = n$:
\[
\lower11pt\hbox{$\includegraphics[height=1cm]{graffles/circuitU.pdf}$} =
\lower18pt\hbox{$\includegraphics[height=1.5cm]{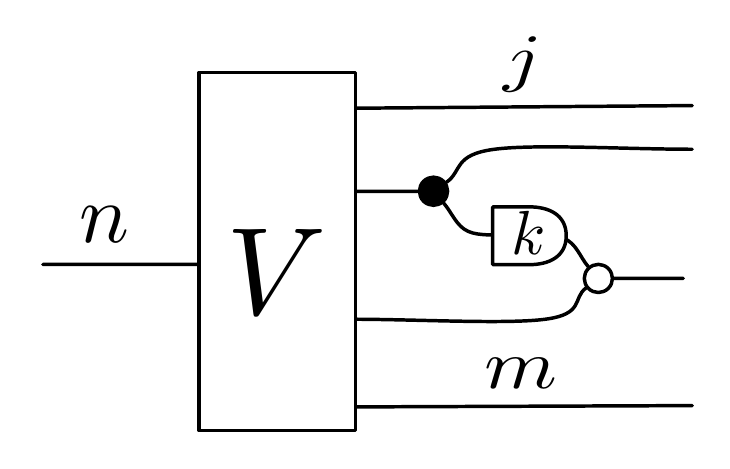}$}
\]
The following two derivations prove that $\coc{U}$ is the inverse of $U$:
\begin{gather*}
\lower18pt\hbox{$\includegraphics[height=1.9cm]{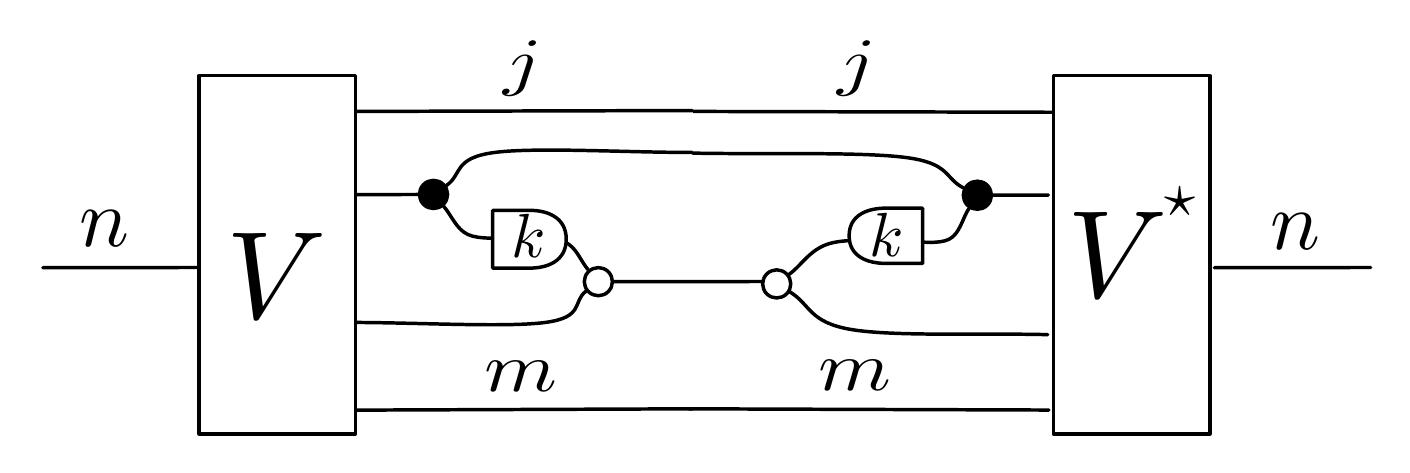}$}
\eql{\eqref{eq:papillon}}
\lower11pt\hbox{$\includegraphics[height=1.2cm]{graffles/circuitVVstar.pdf}$}
\eql{Ind. hyp.}
\lower4pt\hbox{$\includegraphics[height=.8cm]{graffles/circuitidn.pdf}$} \\
\lower18pt\hbox{$\includegraphics[height=2.1cm]{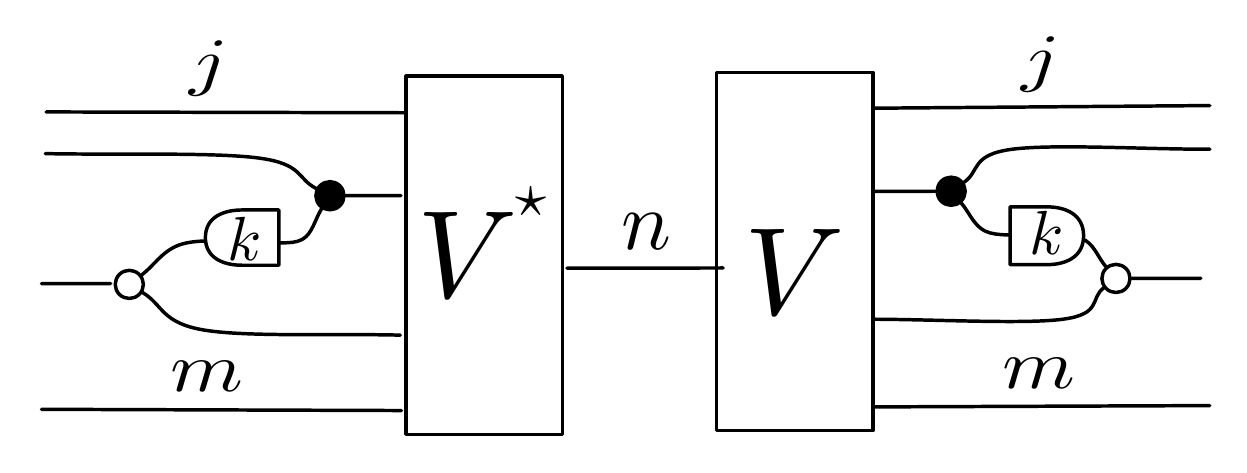}$}
\eql{Ind. hyp.}
\lower16pt\hbox{$\includegraphics[height=2cm]{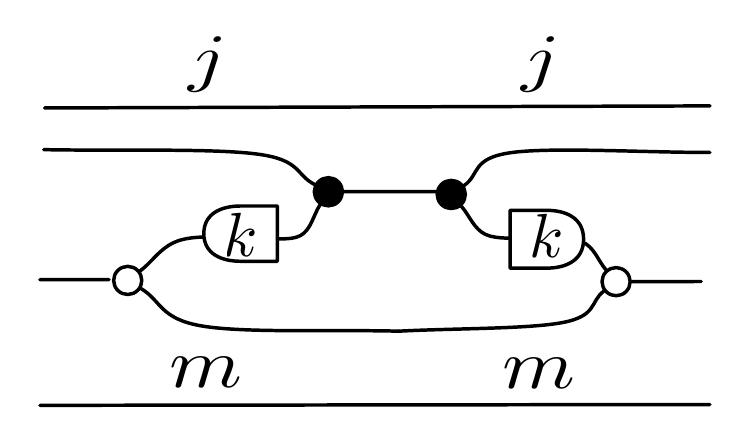}$}
\eql{\eqref{eq:papillon}}
\lower4pt\hbox{$\includegraphics[height=.8cm]{graffles/circuitidn.pdf}$}
 \end{gather*}
 Finally, we have the inductive case in which $U$ is obtained by $V$ via multiplication of a row by a invertible element $i \in \PID$. Let us write $i^{-1} \in \PID$ for the multiplicative inverse of $i$. We can represent $U$ as the following string diagram, where $z+1+m=n$:
 \[
\lower11pt\hbox{$\includegraphics[height=1cm]{graffles/circuitU.pdf}$} =
\lower14pt\hbox{$\includegraphics[height=1.3cm]{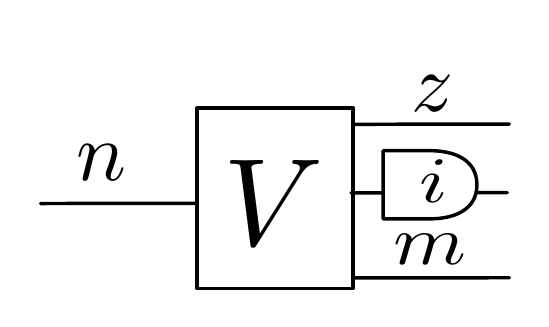}$}
\]
The desired equalities are derivable in $\IBRw$ as follows.
 \begin{eqnarray*}
\lower13pt\hbox{$\includegraphics[height=1.6cm]{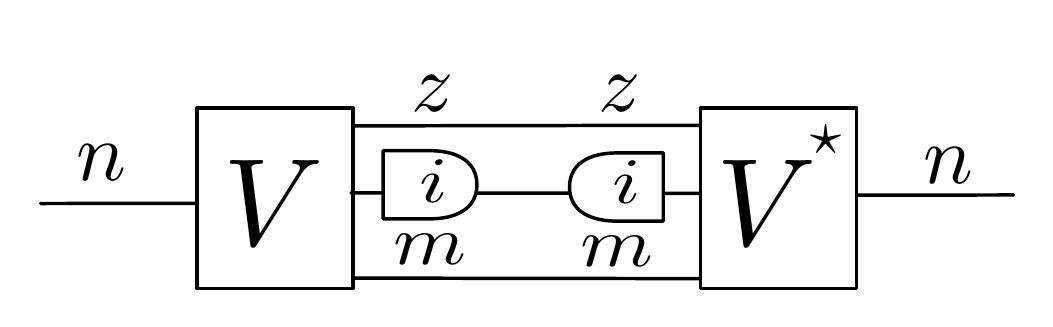}$}
\eql{\eqref{eq:lcm}}
\lower11pt\hbox{$\includegraphics[height=1cm]{graffles/circuitVVstar.pdf}$}
\eql{Ind. hyp.}
\lower4pt\hbox{$\includegraphics[height=.8cm]{graffles/circuitidn.pdf}$}
\end{eqnarray*}
\begin{eqnarray*}
\lower13pt\hbox{$\includegraphics[height=1.6cm]{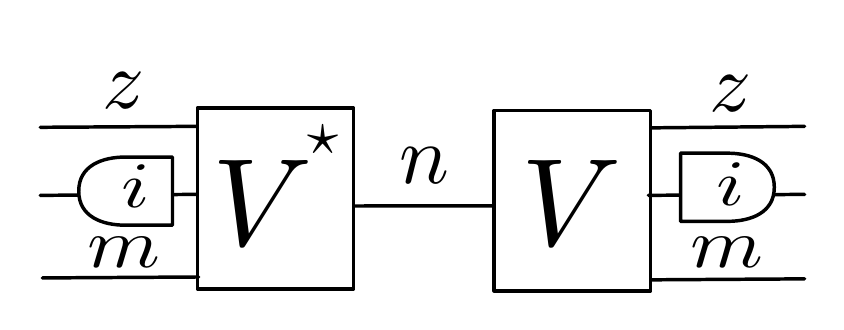}$}
\!\!\eql{Ind. hyp.}\!\!
\lower12pt\hbox{$\includegraphics[height=1.5cm]{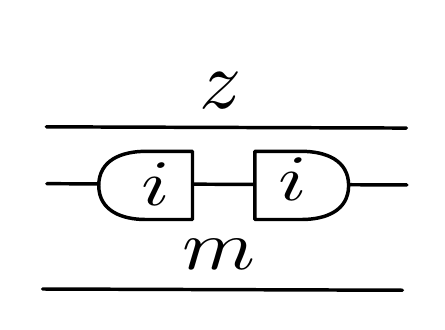}$}
& \!\!\eql{\eqref{eq:lcm}}\!\! &
\lower12pt\hbox{$\includegraphics[height=1.7cm]{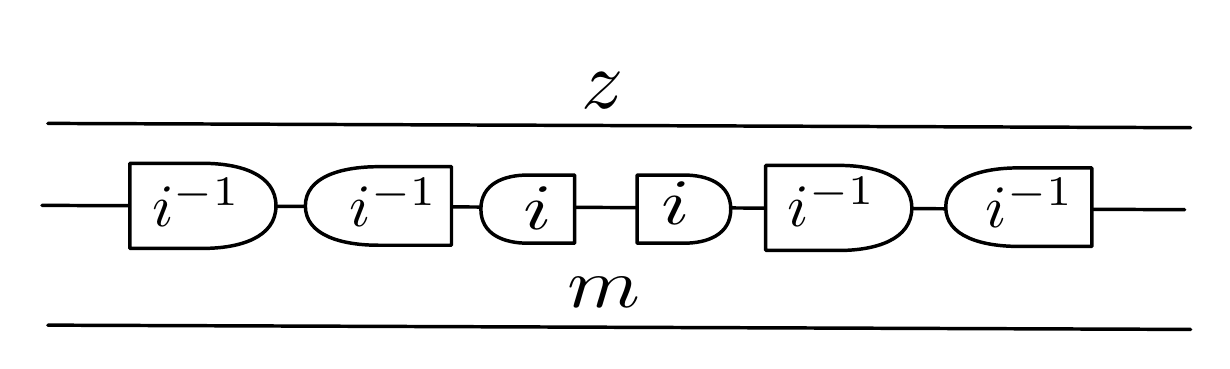}$}\\
& \!\!\eql{\eqref{eq:scalarmult}}\!\! &
\lower12pt\hbox{$\includegraphics[height=1.5cm]{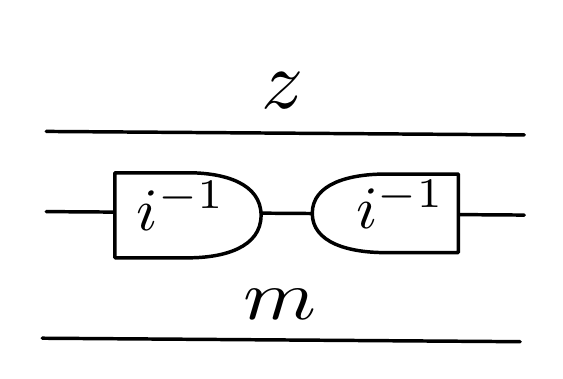}$}
\!\!\eql{\eqref{eq:lcm}}\!\!
\lower4pt\hbox{$\includegraphics[height=.8cm]{graffles/circuitidn.pdf}$}
 \end{eqnarray*}
\end{proof}

The next lemma guarantees that spans which are identified in $\SpanMat$ are not distinguished by the equational theory of $\IBRw$. Recall that $n \tl{A} z \tr{B} m$ and $n \tl{C} z \tr{D} m$ are identified as arrows of $\SpanMat$ precisely when they are isomorphic spans. Since isomorphisms in $\VectR$ are precisely the invertible $\PID$-matrices, this conditions amounts to the existence of an invertible matrix $U \in \VectR[z,z]$ such that the following diagram commutes: 
\begin{eqnarray}\label{diag:isospan}
\vcenter{
\xymatrix@C=15pt@R=15pt{ && \ar[drr]^{B} z  \ar[dll]_{A}&& \\
 n && \ar[ll]^{C} z \ar[u]^>>>>U \ar[rr]_{D}&& m }
 }
\end{eqnarray}
\begin{lemma}\label{lemma:mirror} Let $A,B,C,D,U$ be as in \eqref{diag:isospan}. Then the following equation holds in $\IBRw$:
\[
\lower11pt\hbox{$\includegraphics[height=1cm]{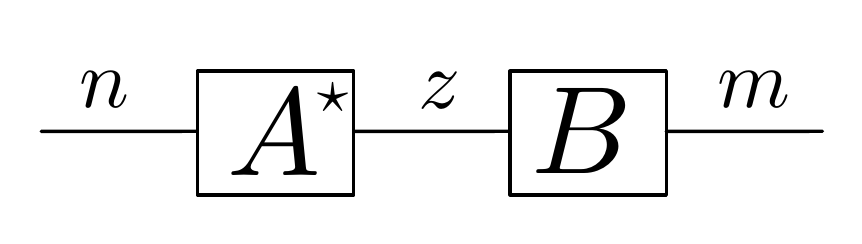}$} = \lower11pt\hbox{$\includegraphics[height=1cm]{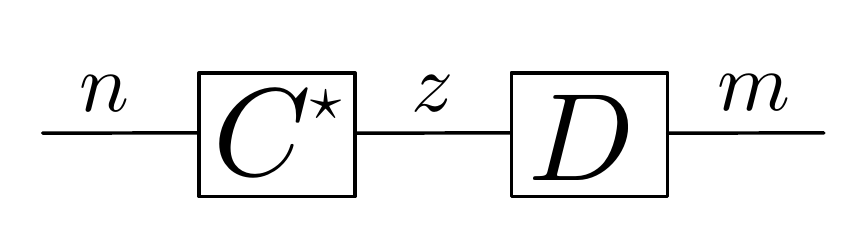}$}
\]
\end{lemma}

\begin{remark}
We emphasize the importance of Lemma~\ref{lemma:mirror} for a proof of our completeness statement: instead of showing that \emph{any} equation arising from a pullback square in $\VectR$ is provable in $\IBRw$, we can now confine ourselves to a canonical choice of pullback for any cospan in $\VectR$. All the equations arising from other pullback diagrams for the same cospan will be provable in $\IBRw$ by Lemma~\ref{lemma:mirror}. Our canonical choice will involve kernel computation, see Lemma~\ref{lemma:pbKernel} below.
\end{remark}

\begin{proof}[Proof of Lemma~\ref{lemma:mirror}] Since $\ABR \cong \VectR$, commutativity of \eqref{diag:isospan} yields the following equations in $\ABR$:
        \begin{equation}\label{eq:mirror1}
         \lower9pt\hbox{$\includegraphics[height=1cm]{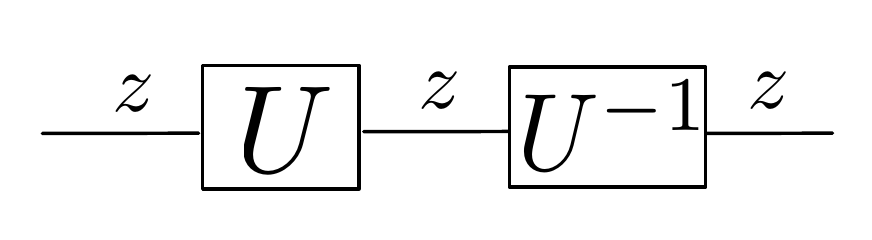}$}
   =\lower4pt\hbox{$\includegraphics[height=1cm]{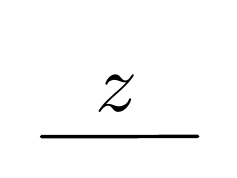}$}
   =          \lower9pt\hbox{$\includegraphics[height=1cm]{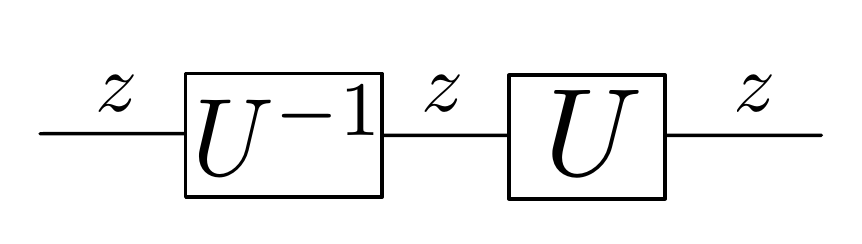}$}
   \end{equation}
        \begin{multicols}{2}
     \noindent
        \begin{equation}\label{eq:mirror2}    \lower7pt\hbox{$\includegraphics[height=.8cm]{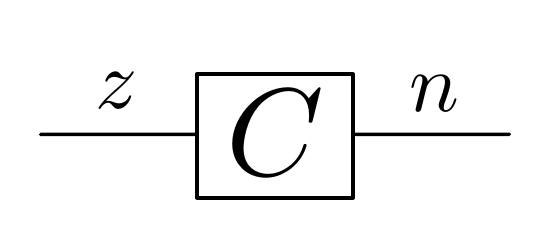}$}
  \!=\! \lower7pt\hbox{$\includegraphics[height=.8cm]{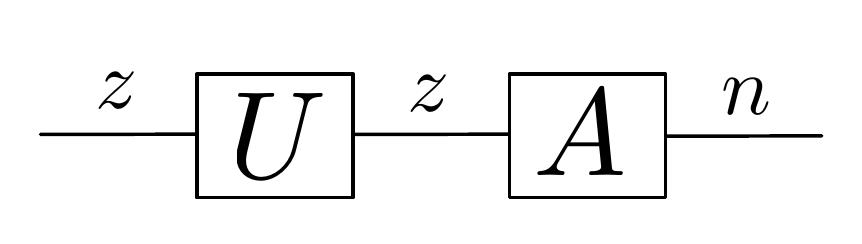}$}
   \end{equation}
        \begin{equation}   \label{eq:mirror3} \lower7pt\hbox{$\includegraphics[height=.8cm]{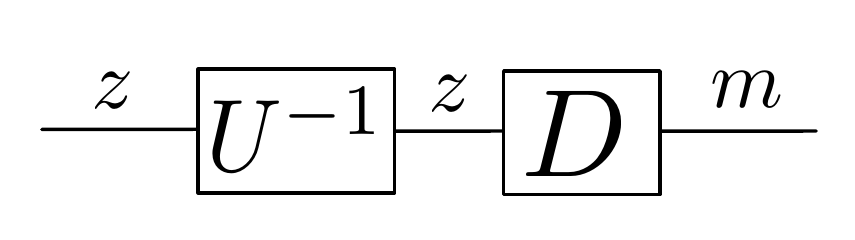}$}
  \!=\! \lower7pt\hbox{$\includegraphics[height=.8cm]{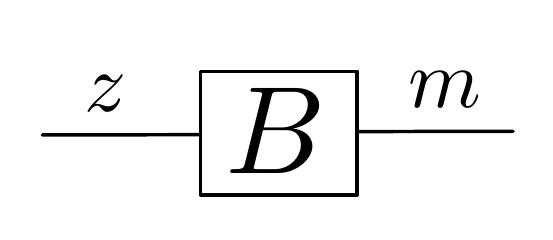}$}
   \end{equation}
 \end{multicols}
Since $\ABR$ is a sub-theory of $\IBRw$, these equations are also valid in $\IBRw$. The statement of the lemma is then given by the following derivation.
\begin{eqnarray*}
\lower11pt\hbox{$\includegraphics[height=1cm]{graffles/circuitCstarDz.pdf}$}
& \eql{\eqref{eq:mirror1}} & \lower11pt\hbox{$\includegraphics[height=1cm]{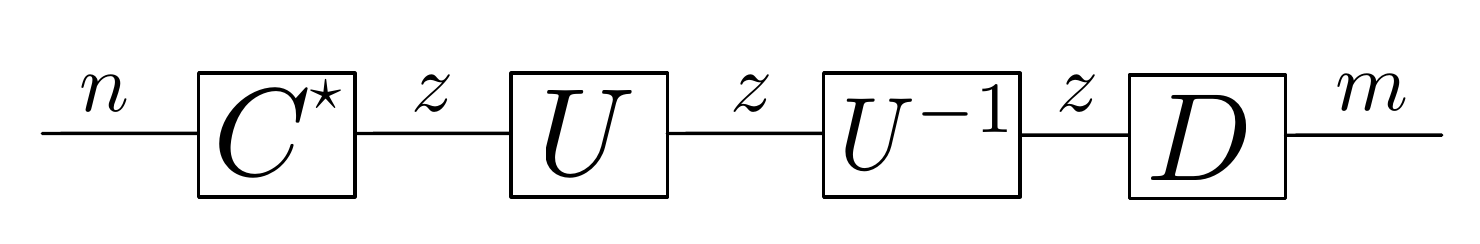}$} \\
& \eql{\eqref{eq:mirror3}} & \lower11pt\hbox{$\includegraphics[height=1cm]{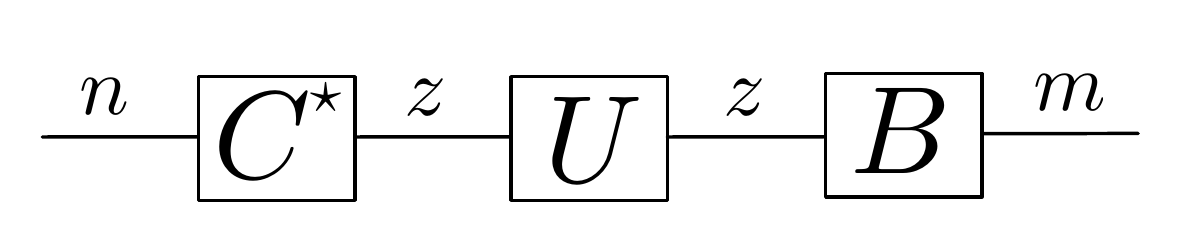}$} \\
& \eql{\eqref{eq:mirror2}} &
\lower13pt\hbox{$\includegraphics[height=1.2cm]{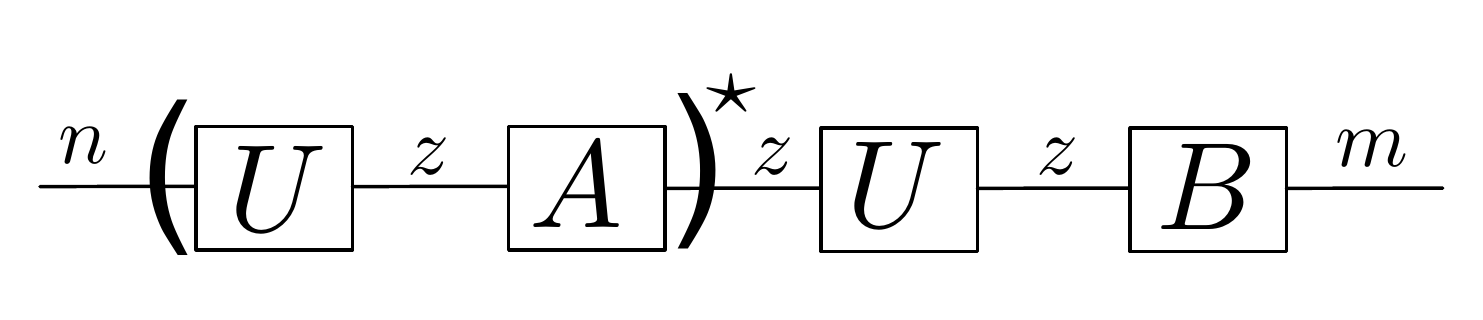}$} \\
& \eql{Def. $\coc{(\cdot)}$} &
\lower11pt\hbox{$\includegraphics[height=1cm]{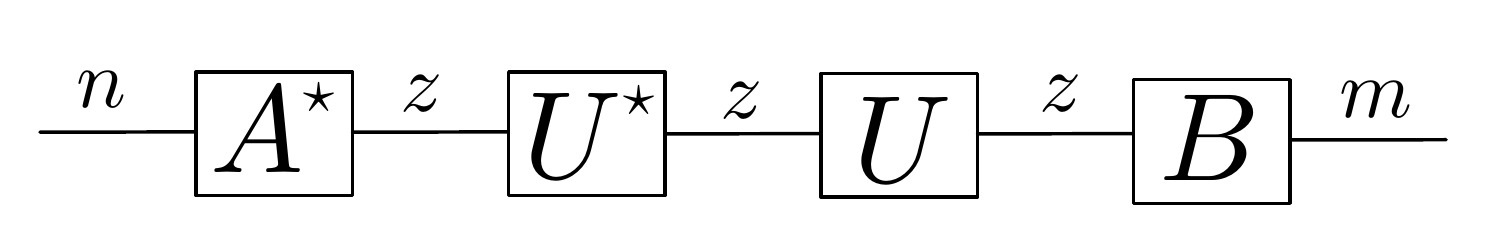}$} \\
& \eql{Lemma~\ref{lemma:invertiblestar}} &
\lower11pt\hbox{$\includegraphics[height=1cm]{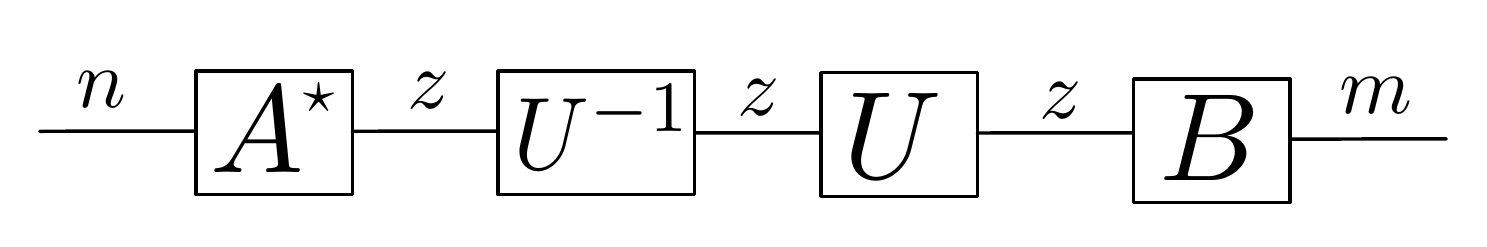}$} \\
& \eql{\eqref{eq:mirror1}} &
\lower11pt\hbox{$\includegraphics[height=1cm]{graffles/circuitAstarB.pdf}$}
\end{eqnarray*}
\end{proof}

The next lemma is another important ingredient in the proof of Proposition~\ref{prop:IBwComplete}: it allows us to reduce, in the graphical theory, the computation of pullbacks to the computation of kernels.

\begin{lemma}\label{lemma:pbKernel} Given a pullback square in $\VectR$ as on the left, the equation on the right holds in $\IBRw$:
\[\xymatrix@R=10pt@C=10pt{
&\ar[dl]_{C} r \pushoutcorner \ar[dr]^{D} & \\
n \ar[dr]_{A} & & m \ar[dl]^{B}\\
& z  & }
\qquad \qquad
\lower20pt\hbox{
\lower12pt\hbox{$\includegraphics[height=1cm]{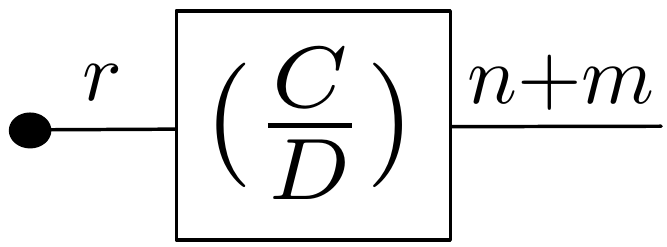}$} =
\lower8pt\hbox{$\includegraphics[height=.8cm]{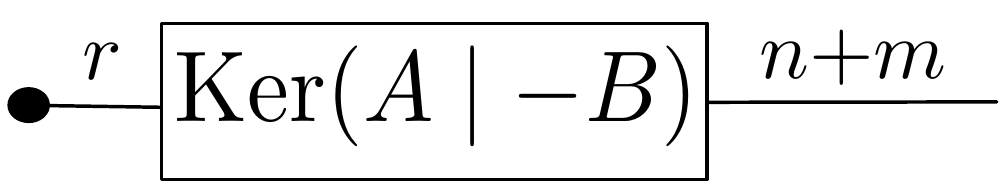}$}
}
\]
\end{lemma}
\begin{proof}
As shown in the proof of Proposition~\ref{prop:PullbacksMatR}, the pullback of $n \tr{A} z \tl{B} m$ can be canonically obtained as follows. First, compute the kernel $\KerAB \: r \to n+m$ of the matrix $(A \mid \minusB)$. The universal property of $n+m$ (see \S~\ref{sec:colimitsMat}) yields $C' \: r \to n$ and $D' \: r \to m$ such that ${\KerAB} = (\frac{C'}{D'})$: then, the span $n \tl{C'} r \tr{D'} m$ pulls back $n \tr{A} z \tl{B} m$.

Since by assumption $n \tl{C} r \tr{D} m$ pulls back the same cospan, then $ \tl{C'} \tr{D'}$
 and $\tl{C} \tr{D}$ are isomorphic spans. Using the conclusion of Lemma~\ref{lemma:mirror}, we infer that
\begin{equation}\label{eq:mirrorAppliedToKernel} \tag{$\nabla$}
\lower8pt\hbox{$\includegraphics[height=.8cm]{graffles/circuitCstarDr.pdf}$} =
 \lower6pt\hbox{$\includegraphics[height=.6cm]{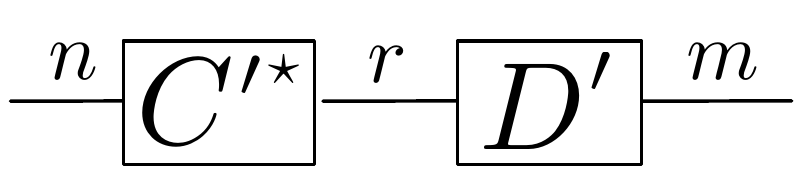}$}
\end{equation}
 from which follows that
\begin{equation}\label{eq:mirrorAppliedToKernel2} \tag{$\triangle$}
\lower10pt\hbox{$\includegraphics[height=1cm]{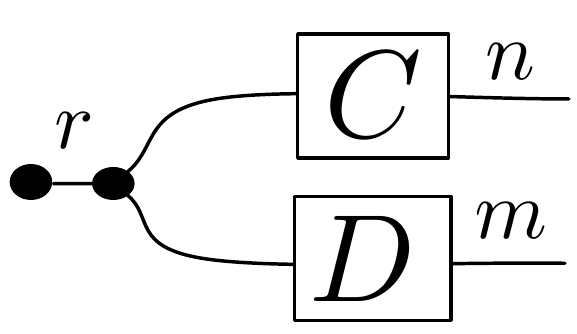}$} \eql{\eqref{eq:ccsliding2}}
\lower10pt\hbox{$\includegraphics[height=1cm]{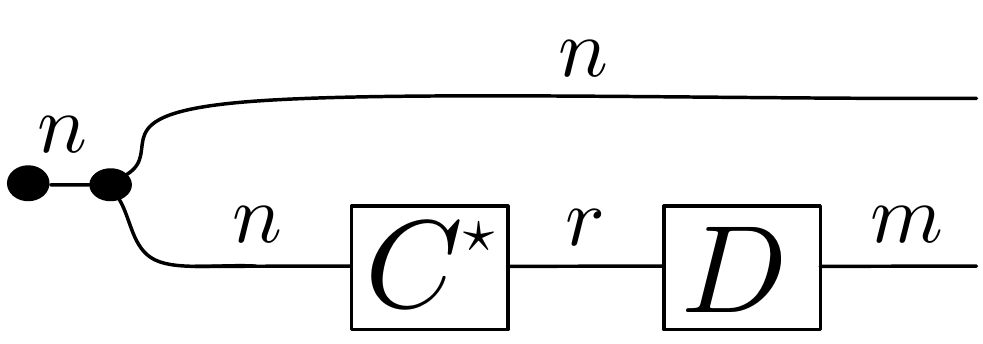}$} \eql{\eqref{eq:mirrorAppliedToKernel}}
\lower10pt\hbox{$\includegraphics[height=1cm]{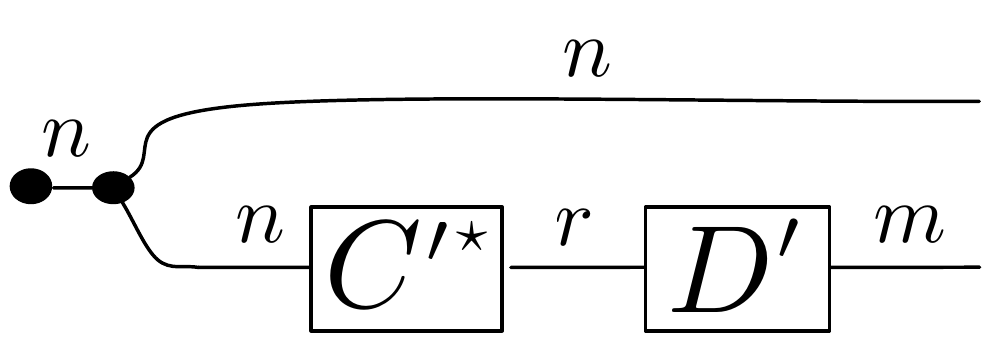}$} \eql{\eqref{eq:ccsliding2}}
\lower10pt\hbox{$\includegraphics[height=1cm]{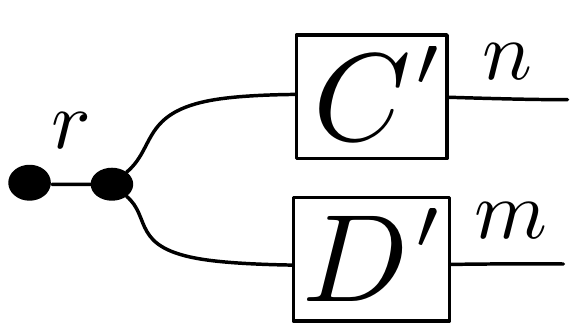}$}.
\end{equation}
 We can now conclude the proof of our statement:
 \[
\lower8pt\hbox{$\includegraphics[height=.8cm]{graffles/BcounitCoverD.pdf}$} \!\eql{Def. $\left(\frac{C}{D}\right)$}\!
\lower10pt\hbox{$\includegraphics[height=1cm]{graffles/ccCD.pdf}$} \!\eql{\eqref{eq:mirrorAppliedToKernel2}}\!
\lower10pt\hbox{$\includegraphics[height=1cm]{graffles/ccCprimeDprime.pdf}$} \!\eql{Def. $\left(\frac{C'}{D'}\right)$}\!
\lower8pt\hbox{$\includegraphics[height=.8cm]{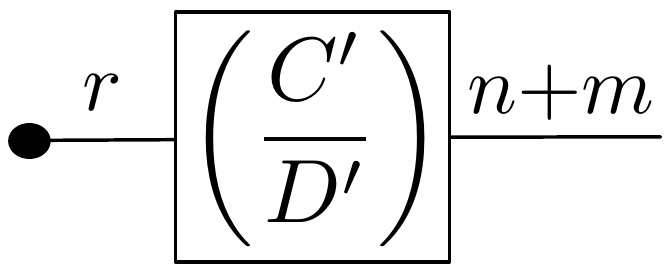}$} \!\eql{}\!
\lower6pt\hbox{$\includegraphics[height=.6cm]{graffles/BcounitkernelAminusB.pdf}$}.
\]
\end{proof} 

\subsubsection{Computing Kernels in $\IBRw$}

 Lemma~\ref{lemma:pbKernel} suggests that kernels play a pivotal role in proving completeness of $\IBRw$. We shall now describe how the kernel computation of a matrix can be formulated
 within the equational theory of $\IBRw$. This is the heart of the completeness argument and the main technical ingredient in the proof of Proposition~\ref{prop:IBwComplete}, which will be given at the end of the section.

 We first recall some linear algebra that will be used in our argument.
 \begin{definition}\label{Def:HNF} An $m \times n$ matrix $A$ is said to be in \emph{Hermite normal form} (HNF) if there is a natural number $r \leq n$ and a strictly increasing function $f \: [r+1,n] \to [1,m]$ mapping column $i$ to a row $f(i)$, such that:
 \begin{enumerate}
   \item the first $r$ columns of $A$ have all entries with value $0$;
   \item for all columns $i$ with $r+1 \leq i \leq n$, $A_{f(i),i} \neq 0$ and
   \item for all $j \gr f(i)$, $A_{j,i} = 0$.
 \end{enumerate}
\end{definition}

\begin{example}\label{ex:HNFmatrix} The following are examples of integer matrices in HNF.
\begin{eqnarray*}
{\scriptsize
\left(
  \begin{array}{cccc}
    0 & 0 & 2 & -1\\
    0 & 4 & 1 & -3\\
    0 & 0 & 1 & 0\\
    0 & 0 & 0 & 0\\
    0 & 0 & 0 & 3
  \end{array}
\right)
\qquad \qquad
\left(
  \begin{array}{cccc}
    0 & 0 & 2 & -1\\
    0 & 0 & -1 & 0\\
    0 & 0 & 0 & 4\\
    0 & 0 & 0 & 0\\
  \end{array}
\right)
\qquad \qquad
\left(
  \begin{array}{ccc}
     0 & -2 & 0\\
     4 & 3 & 3\\
     0 & 0 & 0\\
     0 & 1 & 0\\
     0 & 0 & 3
  \end{array}
\right)
}
\end{eqnarray*}
\end{example}
The Hermit normal form plays for matrices over a PID a role akin to reduced echelon form for matrices over a field. For a thorough account of HNF, see e.g. \cite[\S 2.4.2]{Cohen:1995ComputAlg} and \cite{Conti:1990:HermitePID}. In the following we list some useful properties.
\begin{lemma}\label{lemma:HNFtriangular} Suppose that $A$ is an $m \times n$ matrix in HNF and fix a column $i \leq n$. Then $A_{f(i),j} = 0$ for all columns $j \ls i$.
\end{lemma}
\begin{proof} If $j \leq r$ then $A_{f(i),j} = 0$ by property 1 of HNF. Otherwise, fix $j$ such that $r \ls j \ls i$. Since $f$ is strictly increasing, $f(i) \gr f(j)$. Then by property 3 of HNF, $A_{f(i),j} = 0$.
\end{proof}

A proof of the next statement can be found in \cite[Th. 2.4.3]{Cohen:1995ComputAlg}. Recall that a matrix $A$ is \emph{column-equivalent} to a matrix $B$ if $B = AU$ for some invertible matrix $U$ (which encodes the sequence of column operations --- column sum, swap and multiplication --- allowing to obtain $B$ from $A$).
\begin{proposition} Every $\PID$-matrix is column-equivalent to a matrix in HNF. \end{proposition}

Finding $B = AU$ in HNF gives a means to compute the kernel of $A$ as follows.

\begin{proposition}\label{prop:kernelfromHNF} Let $A$ be an $m \times n$ $\PID$-matrix, $U$ an invertible $n \times n$ $\PID$-matrix such that $B = AU$ is in HNF and $r \leq n$ the number of initial $0$-columns of $B$. Then the first $r$ columns of $U$ form a basis for the kernel of $A$.
\end{proposition}

\begin{example}\label{ex:HNFmatrix2}
The integer matrix $A$, below left, can be turned into the matrix $B = AU$ in HNF, below right, with $U$ an invertible matrix. Since $B$ has only one initial $0$-column, the kernel space for $A$ is generated by the first column of $U$ (in bold).
\begin{equation*}
 \overbrace{ \scriptsize
 \begin{pmatrix}
    3 & -1 & 0 & 0\\
    -2 & 1 & 4 & -4\\
    -1 & 0 & 0 & 0\\
    0 & 0 & 0 & 0\\
    -15 & 3 & 0 & 0
  \end{pmatrix}
  }^{\text{A}} \ \text{column-reduces to} \
\overbrace{ \scriptsize
 \begin{pmatrix}
    0 & 0 & 2 & -1\\
    0 & 4 & 1 & -3\\
    0 & 0 & 1 & 0\\
    0 & 0 & 0 & 0\\
    0 & 0 & 0 & 3
  \end{pmatrix}
  }^{\text{B}} =
  \overbrace{ \scriptsize
 \begin{pmatrix}
    3 & -1 & 0 & 0\\
    -2 & 1 & 4 & -4\\
    -1 & 0 & 0 & 0\\
    0 & 0 & 0 & 0\\
    -15 & 3 & 0 & 0
  \end{pmatrix}
  }^{\text{A}} \cdot
  \overbrace{ \scriptsize
 \begin{pmatrix}
    \textbf{0} & 0 & -1 & 0\\
    \textbf{0} & 0 & -5 & 1\\
    \textbf{1} & 3 & -2 & 0\\
    \textbf{1} & 2 & -3 & 1
  \end{pmatrix}
  }^{\text{U}}.
 \end{equation*}
\end{example}

\begin{proof}[Proof of Proposition~\ref{prop:kernelfromHNF}] A proof can be found for the PID of integers in \cite[Prop. 2.4.9]{Cohen:1995ComputAlg}, which we reformulate here for an arbitrary PID $\PID$. We include the details because the next result will be essentially a graphical rendition of the argument.

 For $i \leq r$, let $\uu_i$ be the $i$-th column of $U$. By definition $A \uu_i = B_i$, which is a $0$-vector because $i \leq r$. Thus all first $r$ columns of $U$ are elements of the kernel of $A$. Conversely, let $\xx$ be a vector such that $A\xx = 0$. Then $A\xx = A U U^{-1} \xx = B U^{-1} \xx$ because $U$ is invertible. Let $y_1, \dots, y_n$ be the coordinates of $\yy \df U^{-1}\xx$. For each $i$ in $[r+1,n]$, we show that $y_i = 0$, by backward induction on $i$. This unfolds as a kind of ``chain reaction'':
\begin{itemize}[itemsep=.5pt,topsep=3pt, partopsep=0pt]
  \item[(I)] \label{pt:chainreaction1} if $i =n$, let $f(n)$ be as in Definition~\ref{Def:HNF}. Since $B\yy = 0$, the $f(n)$-th coordinate of $B\yy$ is
      \begin{equation}\label{eq:propHNF1} \tag{$\bigtriangleup$}
          B_{f(n),1}y_1 + \dots + B_{f(n),n}y_n = 0.
      \end{equation}
      By Lemma~\ref{lemma:HNFtriangular}, $B_{f(n),1}, \dots, B_{f(n),{n-1}}$ are all equal to $0$, meaning by \eqref{eq:propHNF1} that $B_{f(n),n}y_n = 0$. By property 2 of HNF, $B_{f(n),n} \neq 0$ and thus, since $\PID$ has no non-zero divisors, $y_n = 0$.
  \item[(II)] \label{pt:chainreaction2} For $i$ with $r \ls i \ls n$, the $f(i)$-th coordinate of $B\yy$ is $B_{f(i),1}y_1 + \dots + B_{f(i),n}y_n = 0$ and by induction hypothesis $y_j = 0$ for all $j$ such that $i \ls j \leq n$. By Lemma~\ref{lemma:HNFtriangular}, $B_{f(i),1}, \dots, B_{f(i),{i-1}}$ are all equal to $0$, which means, analogously to the base case, that $B_{f(i),i}y_i = 0$ and since $B_{f(i),i}$ then $y_i = 0$.
  \item[(III)] Thus we proved that the coordinates $y_{r+1},\dots,y_n$ of $\yy$ are equal to $0$. Instead the first $r$ coordinates of $\yy$ can be arbitrary, because the $j$-th row of $B\yy$, for $j \leq r$, is give by $B_{j,1}y_1 + \dots + B_{j,n}y_n = 0$ and we know that, by property 1 of HNF, the entries $B_{j,1},\dots, B_{j,n}$ have value $0$.
\end{itemize}
Therefore the kernel of $B$ is generated by the first $r$ canonical basis vectors $\vv_1, \dots \vv_r$ of $\PID^n$. Since $B = AU$, then $U\vv_1, \dots, U\vv_r$ form a basis for the kernel of $A$. But those are just the first $r$ columns of $U$: hence we have proven the statement of the theorem.
\end{proof}

We now recast the core of Proposition~\ref{prop:kernelfromHNF} ``in purely graphical terms''. For an instance of the construction used in the proof, see Example~\ref{ex:HNF}.

\begin{lemma}\label{lemma:BHNFequalKernel} Let $B$ be an $m \times n$ $\PID$-matrix in HNF and $r$ the number of initial $0$-columns of $B$. Then the following holds in $\IBRw$:
$$\lower11pt\hbox{$\includegraphics[height=1cm]{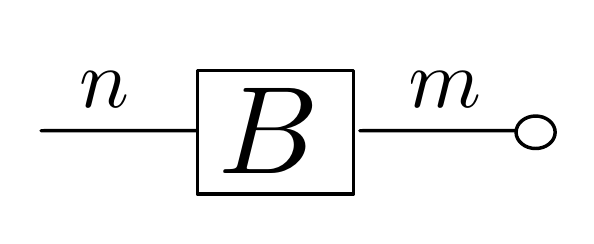}$} = \lower11pt\hbox{$\includegraphics[height=1.2cm]{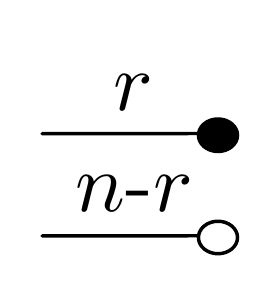}$}$$
\end{lemma}
\begin{proof}
The idea of the proof is to show that the kernel computation described in the proof of Proposition~\ref{prop:kernelfromHNF} can be carried out directly on string diagrams using the equational theory of $\IBRw$. Since $B$ is in HNF, the corresponding string diagram (in matrix form) can be assumed of a particular shape, that we depict below right. 

\noindent \begin{minipage}[t]{.60\textwidth}
Following the factorisation property of matrix form (Def.~\ref{def:matrixform}), we can partition the diagram for $B$ into sub-diagrams from $\Com$, $\PROPR$ and $\Mon$. In particular, the sub-diagram $P$ is only made of basic components from $\PROPR$, of the kind $\symNet$, $\Idnet$ and $\scalar$. Now, by property 1 of HNF, the first $r$ columns of $B$ only have $0$ entries, meaning that the topmost $r$ ports on the left boundary are not connected to the right boundary. Also, by Lemma~\ref{lemma:HNFtriangular} we know that the $f(n)$-th row of $B$ (where $f \: [r+1,n] \to [1,m]$ is as in Definition~\ref{Def:HNF}) has only one non-zero value $k \in \PID$, at position $B_{f(n),n}$.  In diagrammatic terms, this allows us to assume that the $f(n)$-th port on the right boundary only connects to the $n$-th and last port on the left boundary.
\end{minipage}
\begin{minipage}[t]{.35\textwidth}
\vspace{-.8cm}$$\includegraphics[height=5.5cm]{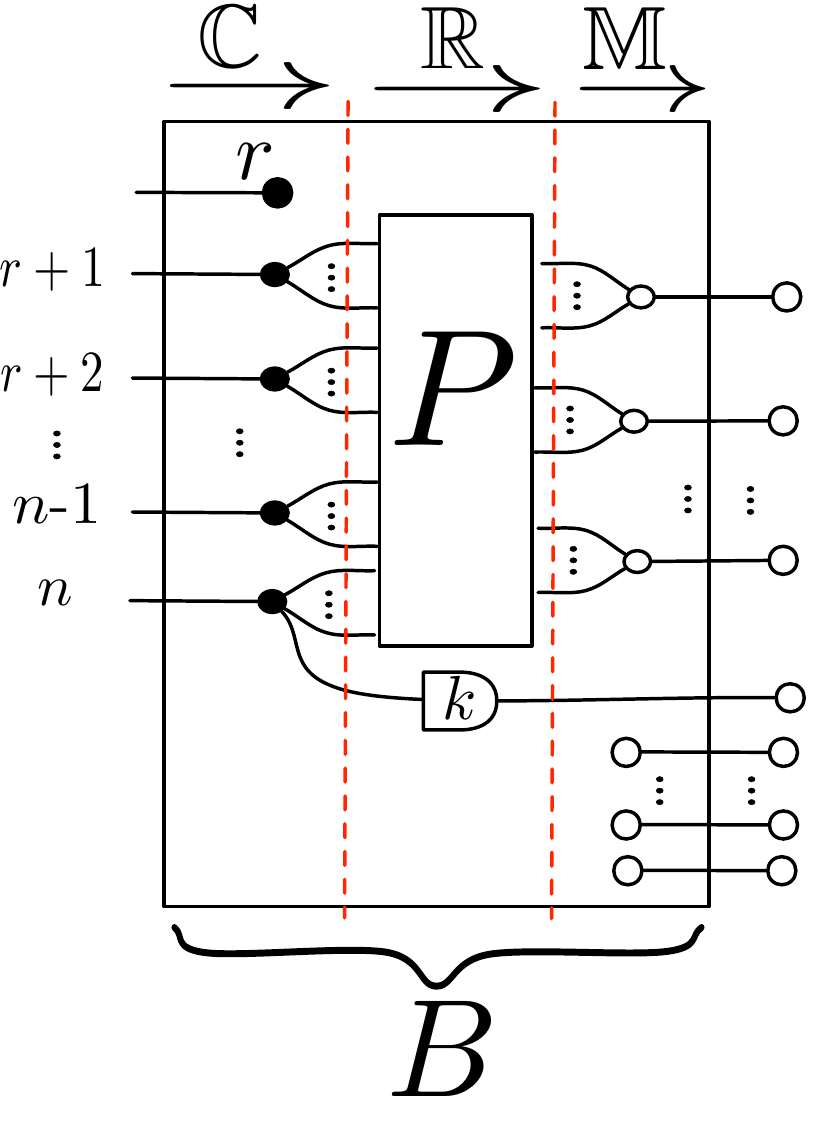}$$
\end{minipage}

\noindent As yet another consequence of the definition of HNF, we know that, for each $i$ with $m \geq i \gr f(n)$, row $i$ only has $0$ entries, allowing us to represent all the rows below $f(n)$ in the string diagram above as ports on the right boundary not connected to any port on the left. Once we plug counits on the right of the string diagram representing $B$, we trigger the chain reaction described in the proof of Proposition~\ref{prop:kernelfromHNF}, which we now reproduce in diagrammatic terms. By backward induction on $i$ with $n \geq i \gr r$, we construct string diagrams $B_n, \dots , B_{r+1}$ such that:
\[
\lower11pt\hbox{$\includegraphics[height=1cm]{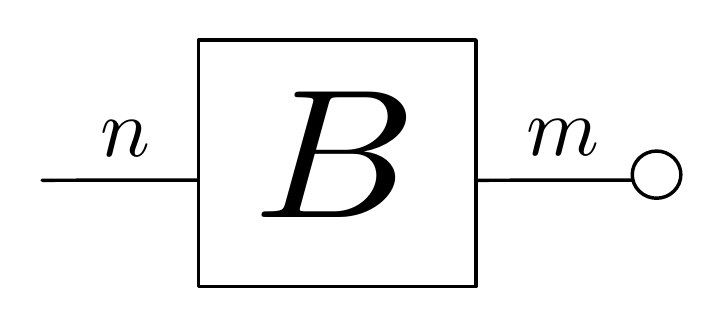}$} =
\lower11pt\hbox{$\includegraphics[height=1cm]{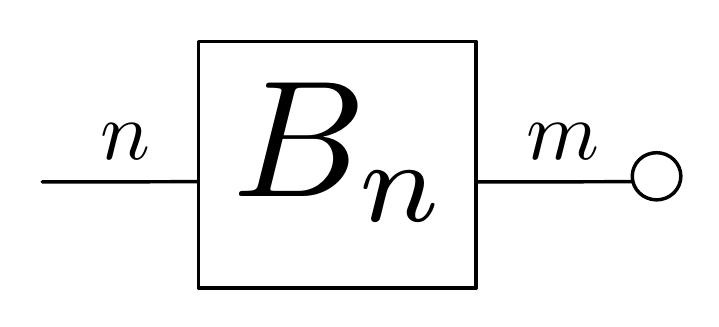}$}
= \dots =
\lower11pt\hbox{$\includegraphics[height=1cm]{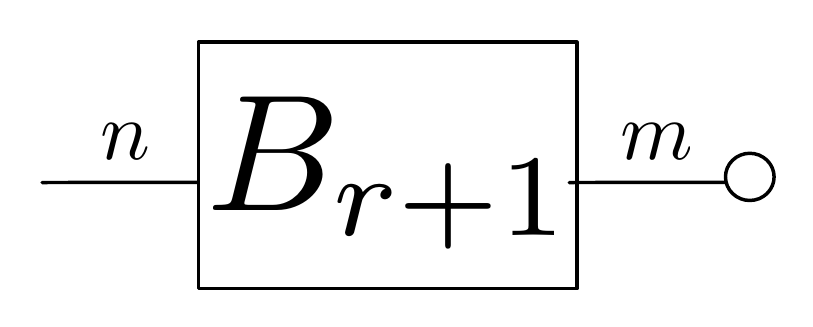}$} =
\lower11pt\hbox{$\includegraphics[height=1.2cm]{graffles/circuitCounitsrn-r.pdf}$}
\]
Clearly, this suffices to prove the main statement.
\begin{itemize}
\item[(I)] For the base case, suppose $i = n$. Since $k \neq 0$, we can use the derived law \eqref{eq:wunitcancelbcomult} of $\IBRw$ to ``disconnect'' the $n$-th port on the left from any port on the right. We define $B_n$ in terms of the resulting string diagram.
\begin{eqnarray*}
\lower81pt\hbox{$\includegraphics[height=5cm]{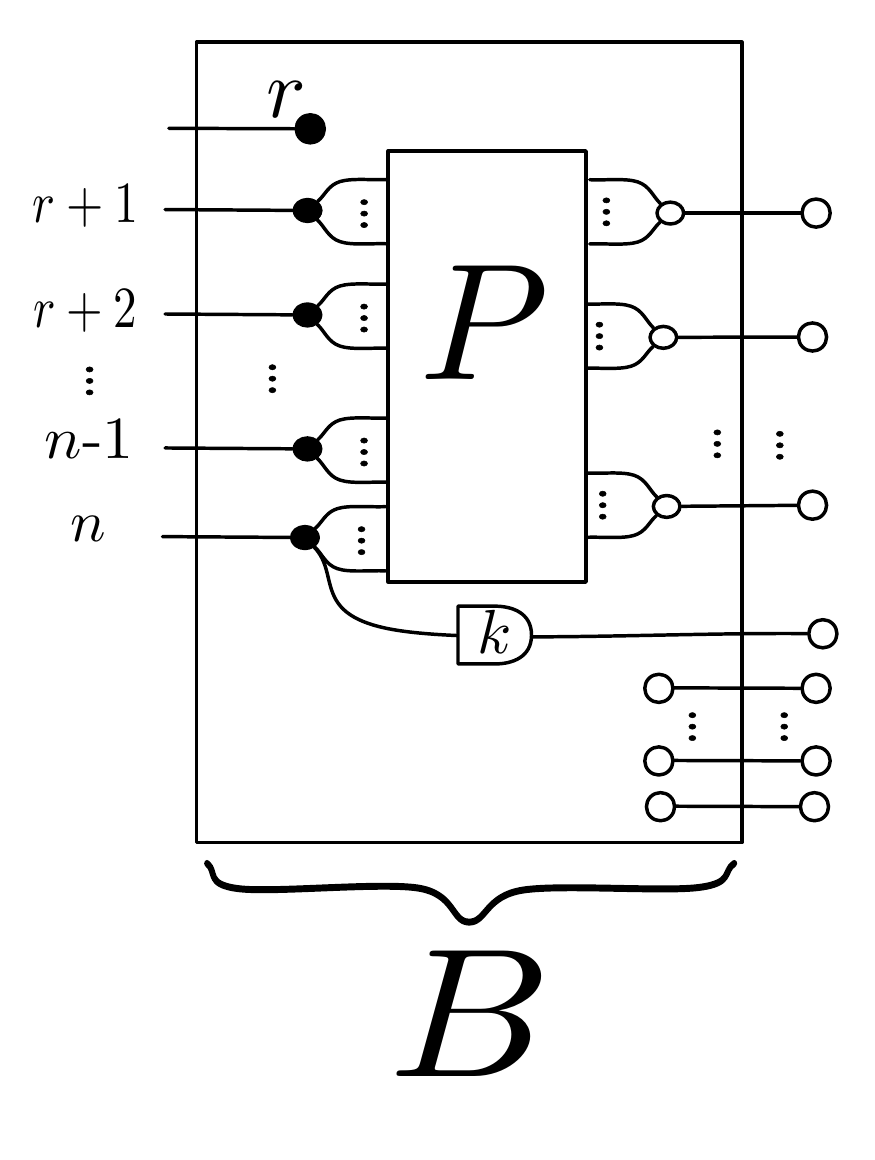}$} \ \eql{\eqref{eq:wunitcancelbcomult}} \
\lower84pt\hbox{$\includegraphics[height=5cm]{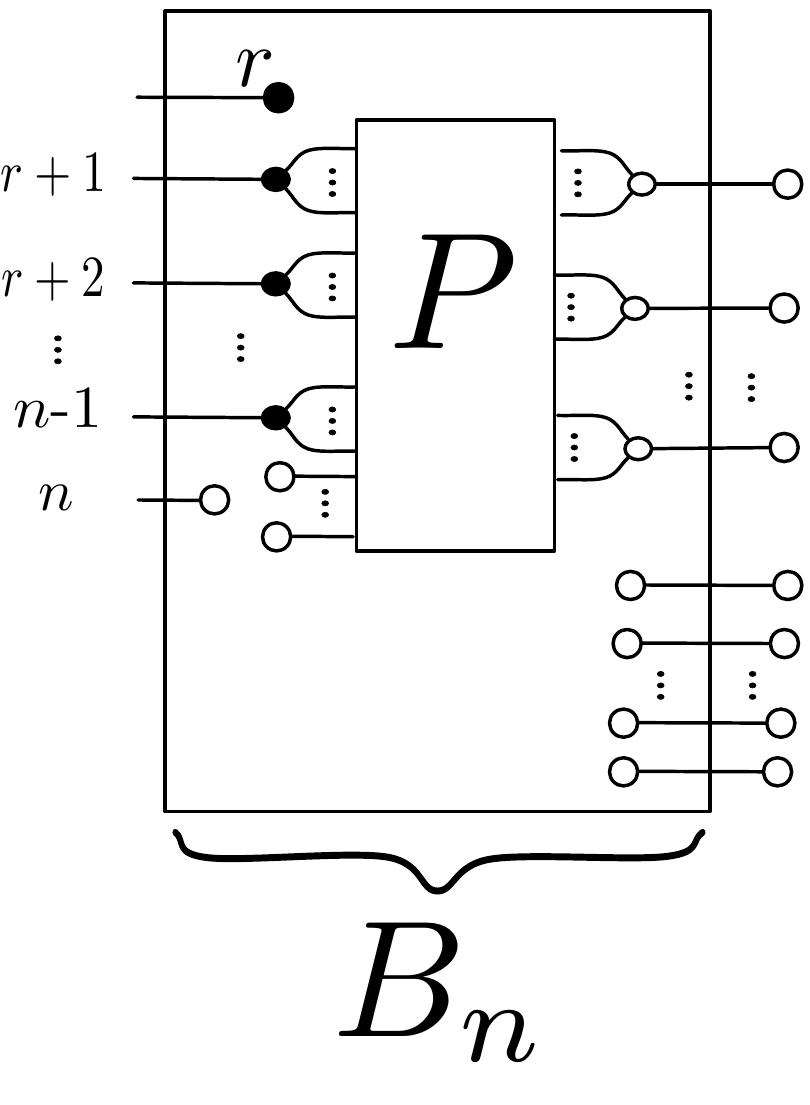}$} \dfop
\lower11pt\hbox{$\includegraphics[height=1.2cm]{graffles/circuitBncounits.pdf}$}
\end{eqnarray*}
We assign the name $P_n$ to the string diagram $P$ depicted above and proceed with the inductive step of $i$ with $n \gr i \gr r$.

\item[(II)]
\noindent \begin{minipage}[t]{.60\textwidth}
The inductive construction
gives us a string diagram $B_{i+1}$ as on the right. The $i$-th port on the left boundary corresponds to column $i$ in $B$ and thus it is assigned a row $f(i)$. This corresponds to the $f(i)$-th port on the right boundary of the string diagram representing $B_{i+1}$. By Lemma~\ref{lemma:HNFtriangular}, such a port has no connections with ports $1,\dots,i-1$ on the left boundary. Moreover, by inductive hypothesis it also has no connections with ports $i+1,\dots,n$ on the left boundary. Therefore port $f(i)$ on the right connects only to port $i$ on the left. These connections are part of the string diagram $P_{i+1}$ --- which by inductive construction only contains $\symNet$, $\Idnet$ and $\scalar$.
\end{minipage}
\begin{minipage}[t]{.35\textwidth}
\vspace{-.5cm}\[
\includegraphics[height=5cm]{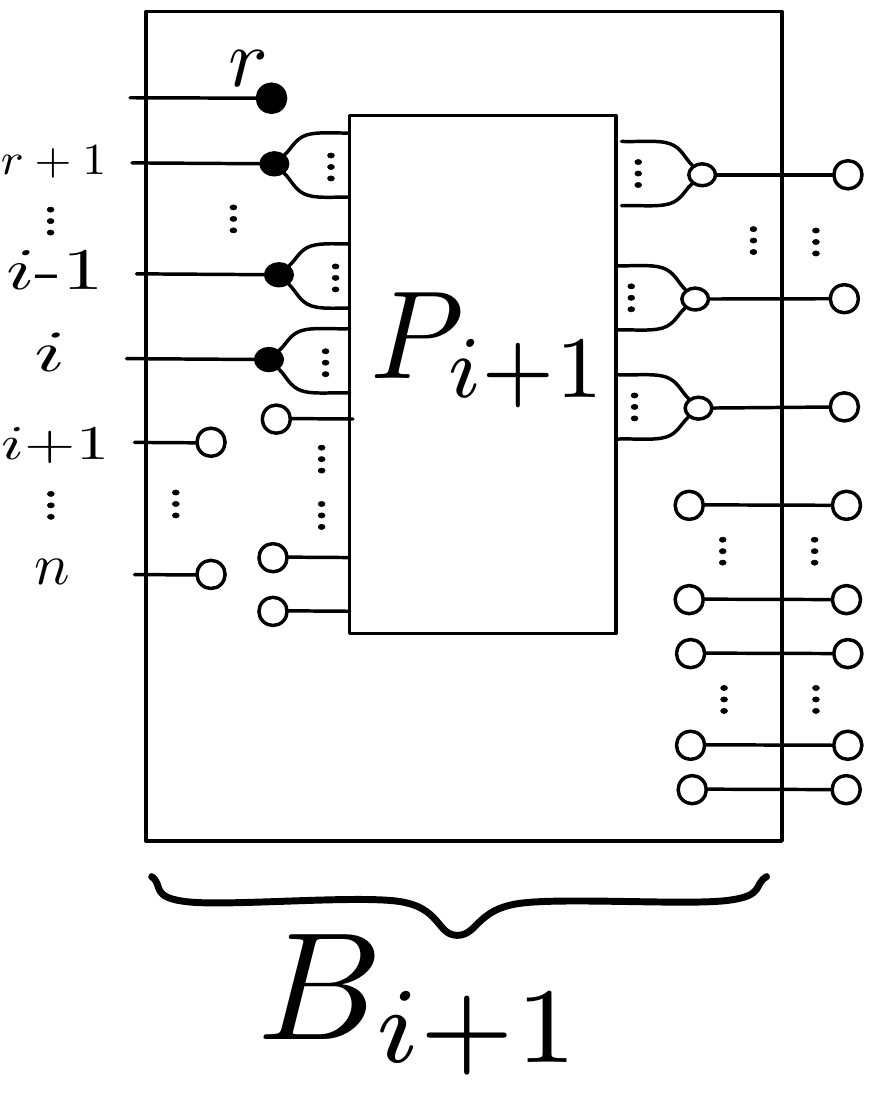}
\]
\end{minipage}
 \noindent It should then be clear that we can use~\eqref{eq:sliding1}-\eqref{eq:sliding2} to ``move port $f(i)$ towards the left side of the diagram'', isolating its connections from the others in $P_{i+1}$.
The resulting string diagram is depicted below, where $P_i$ results from the rearrangement of $P_{i+1}$ in order to allow the move of port $f(i)$ towards the left side of the diagram.
\begin{equation} \label{eq:circuitBn+1}
\lower70pt\hbox{$\includegraphics[height=5cm]{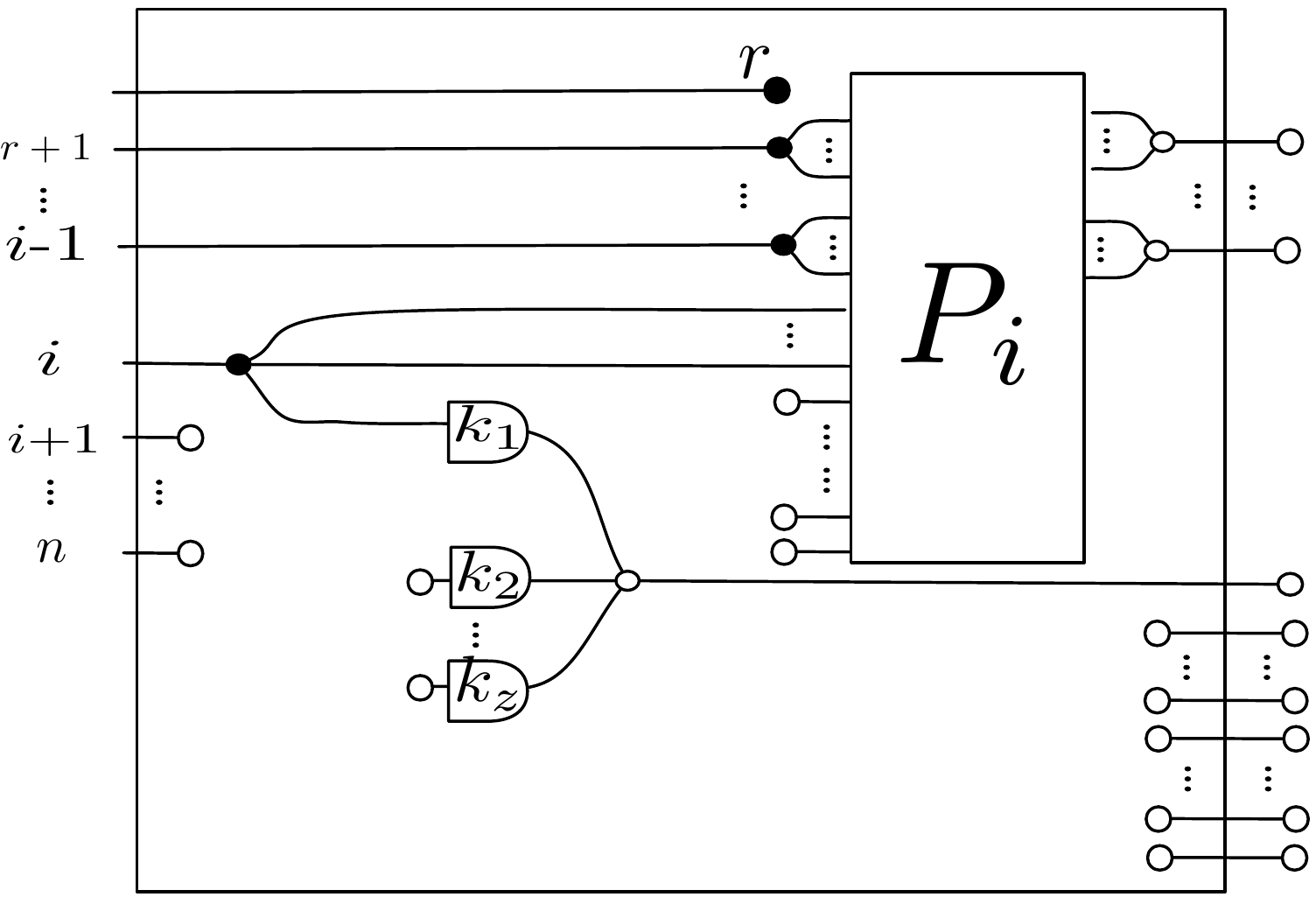}$}
 \end{equation}
We now focus on the sub-diagram depicting the connection of port $i$ on the left with (former) port $f(i)$. In the derivation below, \eqref{eq:wunitcancelbcomult} can be applied because $k_1 = B_{f(i),i} \neq 0$.
\begin{equation*}
\lower35pt\hbox{$\includegraphics[height=2.5cm]{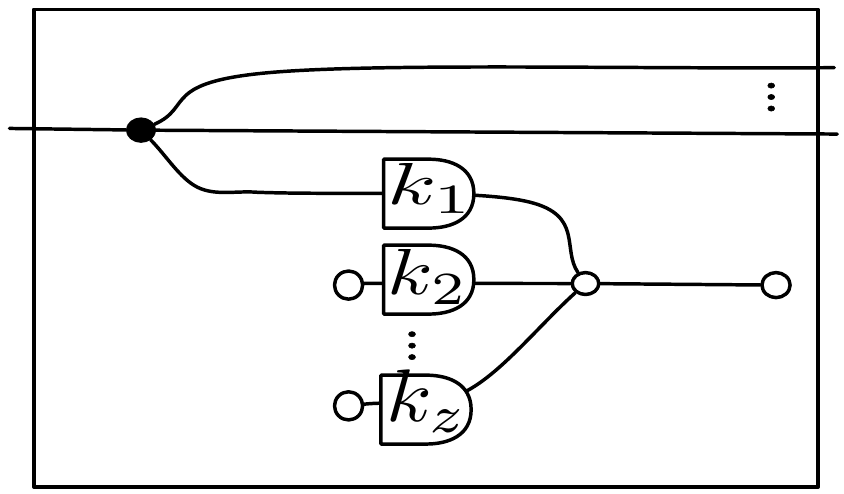}$}
 \eql{\eqref{eq:scalarwunit},\eqref{eq:wmonunitlaw}}
\lower21pt\hbox{$\includegraphics[height=1.6cm]{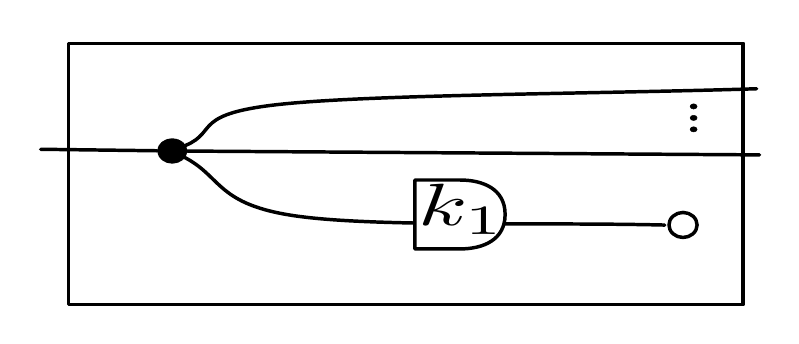}$}
 \eql{\eqref{eq:wunitcancelbcomult}}
\lower15pt\hbox{$\includegraphics[height=1.1cm]{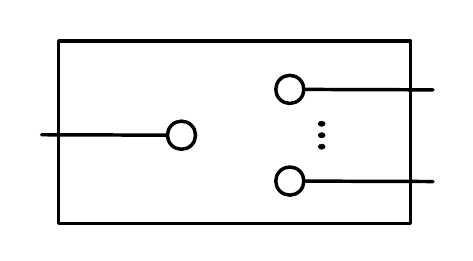}$}
 \end{equation*}
Thus \eqref{eq:circuitBn+1} is equal to the string diagram below left, from which we define $B_{i}$.
\[
\lower62pt\hbox{$\includegraphics[height=4.4cm]{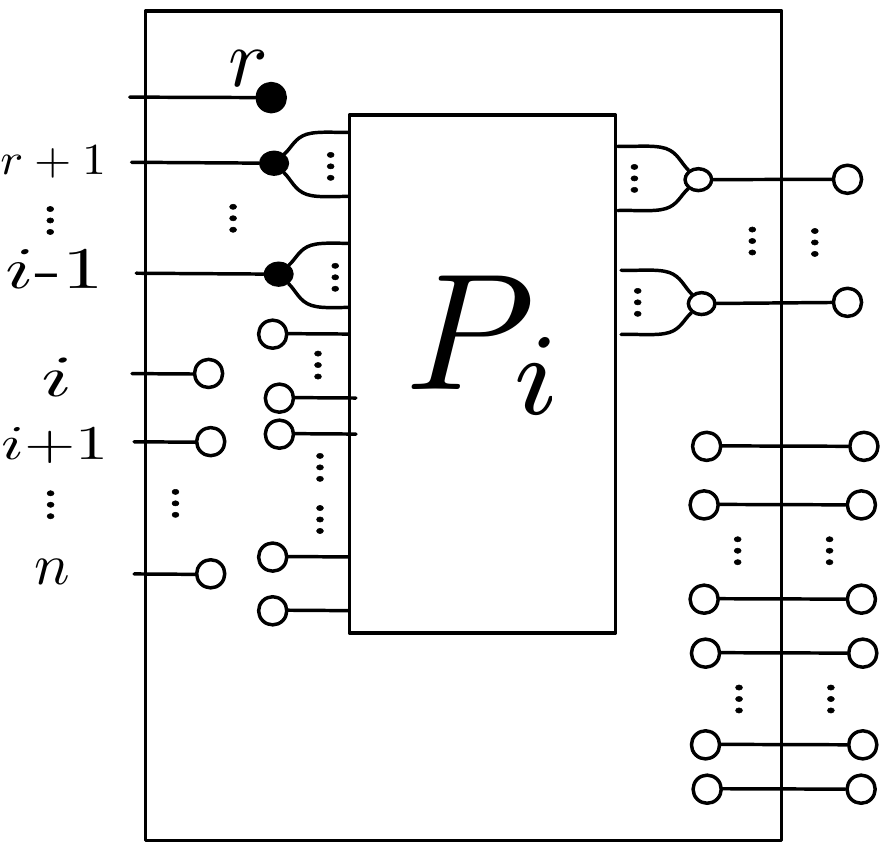}$}
\quad \dfop \quad
\lower15pt\hbox{$\includegraphics[height=1.5cm]{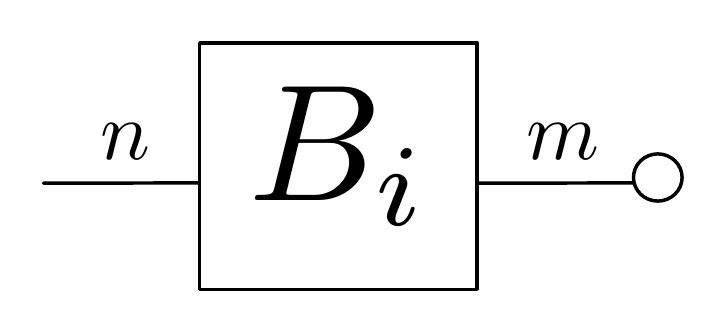}$}
\]
\item[(III)] Finally, at step $r+1$, our inductive construction produces a string diagram as on the left below. Through our inductive construction we have disconnected all ports $i$ on the left and all ports $f(i)$ on the right: in $P_{r+1}$, components $\scalar$ can only appear in correspondence of rows not in the image of $f$ (if any). We can easily remove also this last piece of information:
\begin{align*}
\lower80pt\hbox{$\includegraphics[height=5cm]{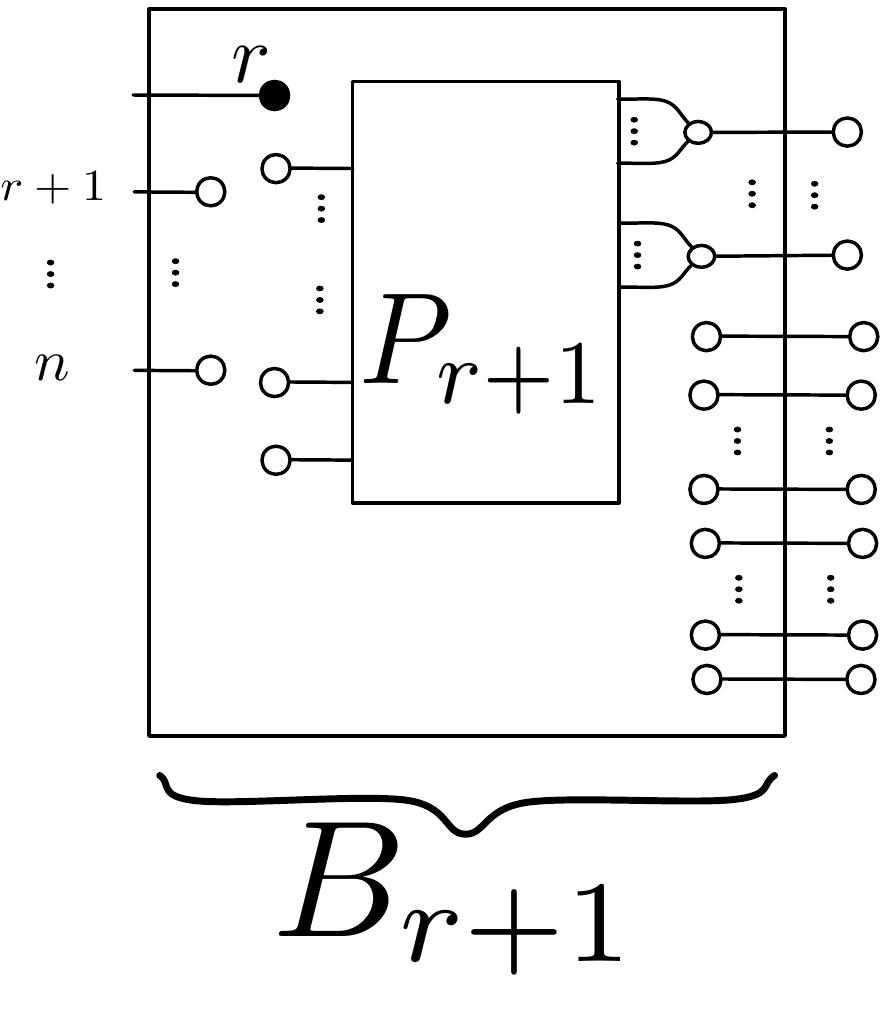}$}\ \ \  & \eql{\eqref{eq:scalarwunit},\eqref{eq:wmonunitlaw}}   & 
\lower42pt\hbox{$\includegraphics[height=3.6cm]{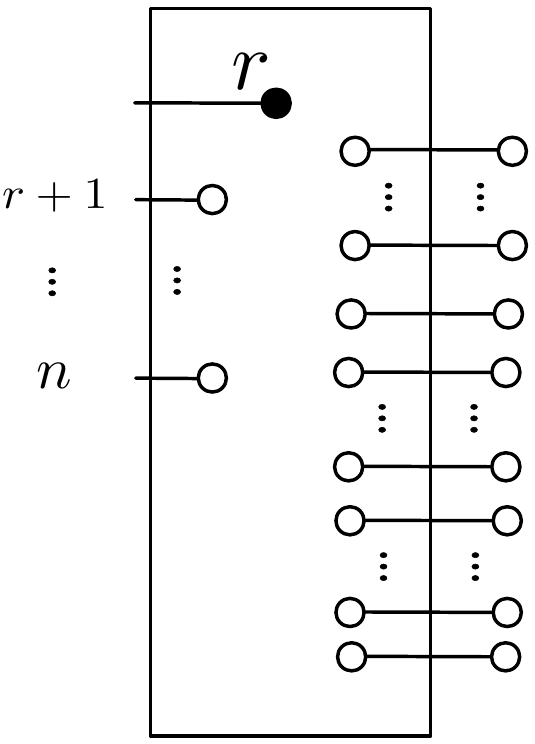}$}
\ \ \ & \eql{\eqref{eq:wbone}} &
\lower15pt\hbox{$\includegraphics[height=1.6cm]{graffles/circuitCounitsrn-r.pdf}$}.
\end{align*}
For the first equality, observe that by inductive construction $P_{r+1}$ is only made of basic components of the kind $\symNet$, $\Idnet$ and $\scalar$ : the white units plugged on the left boundary of $P_{r+1}$ cancel $\symNet$ by application of~\eqref{eq:sliding2} and cancel $\scalar$ by~\eqref{eq:scalarwunit}. The second equality holds by repeated application of~\eqref{eq:wbone}.
\end{itemize}
\end{proof}


 \begin{example}\label{ex:HNF} We show the construction of Lemma~\ref{lemma:BHNFequalKernel} on a string diagram in matrix form that represents the matrix $B$ in HNF of Example~\ref{ex:HNFmatrix}.

\begin{center}
\includegraphics[height=7cm]{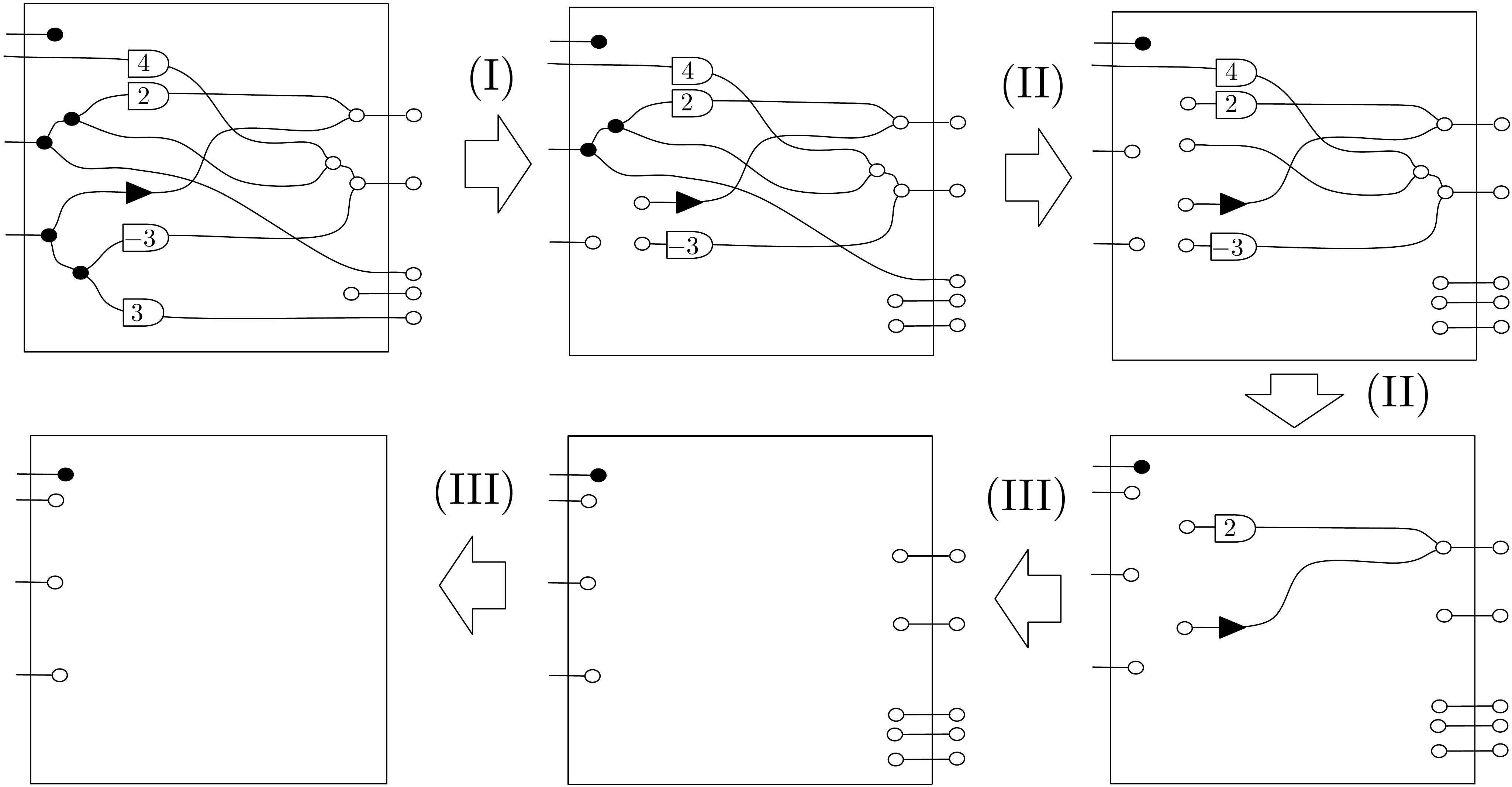}
\end{center}
\end{example}

Given $A \in \VectR[n,m]$ and $r \leq n$, let the \emph{$r$-restriction of $A$} be the matrix $\restr{A}{r} \in \VectR[r,m]$ consisting of the first $r$ columns of $A$. It is useful to make the following observation.

\begin{lemma}\label{lemma:invertiblerestriction} Let $U \in \VectR[n,m]$ be a matrix and fix $r \leq n$. Then the following holds in $\IBRw$:
$$\lower11pt\hbox{$\includegraphics[height=.9cm]{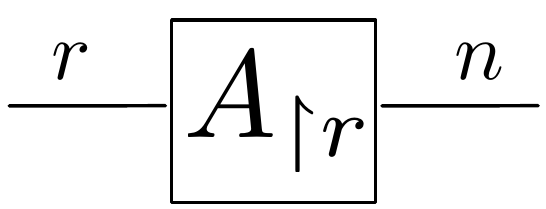}$} =
\lower9pt\hbox{$\includegraphics[height=.9cm]{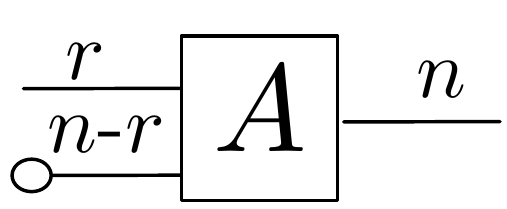}$}$$
\end{lemma}
\begin{proof} One can compute that $\sem{\ABR}(\lower4pt\hbox{\includegraphics[height=.7cm]{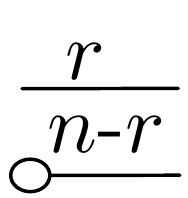}})\poi A = \restr{A}{r}$. Then the statement holds because $\sem{\ABR} \: \ABR \to \VectR$ is an isomorphism. 
\end{proof}

We now have all the ingredients to state the soundness of kernel computation for an arbitrary $\PID$-matrix of $\VectR$.

\begin{proposition}\label{prop:matrixequalkenrel}
Let $A \in \VectR[n,m]$ be a $\PID$-matrix. Then the equation below left, which corresponds to the pullback on the right, is valid in $\IBRw$.
\begin{eqnarray}\label{eq:kernelplback}
\lower9pt\hbox{$\includegraphics[height=.9cm]{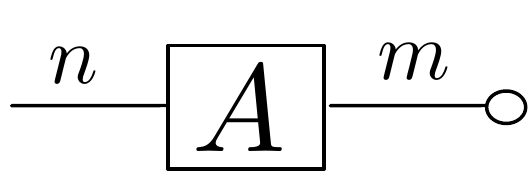}$} \  =
\lower13pt\hbox{$\includegraphics[height=1.2cm]{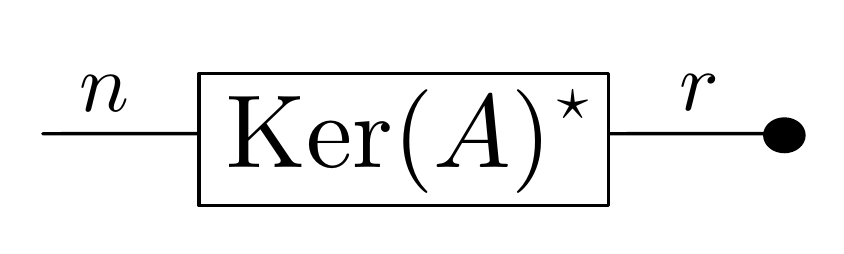}$} & \qquad &
\vcenter{
\xymatrix@R=8pt@C=10pt{
&\ar[dl]_{\Ker{A}} r \pushoutcorner \ar[dr]^{\finVect} & \\
n \ar[dr]_{A} & & 0 \ar[dl]^{\initVect}\\
& m  & }
}
\end{eqnarray}
\end{proposition}
\begin{proof}Let $B = AU$ be the HNF of $A$ for some invertible matrix $U \: n \to n$. We derive in $\IBRw$:
\begin{eqnarray*}
\lower10pt\hbox{$\includegraphics[height=.9cm]{graffles/circuitAwcounits.pdf}$} &\eql{}&
\lower9pt\hbox{$\includegraphics[height=.8cm]{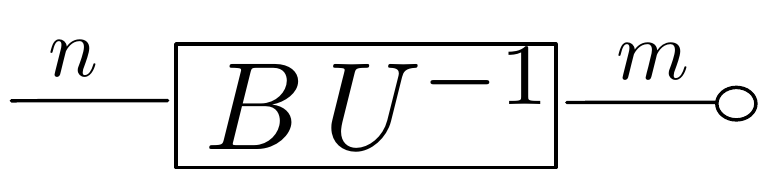}$}
\\
&\eql{}&
\lower11pt\hbox{$\includegraphics[height=1cm]{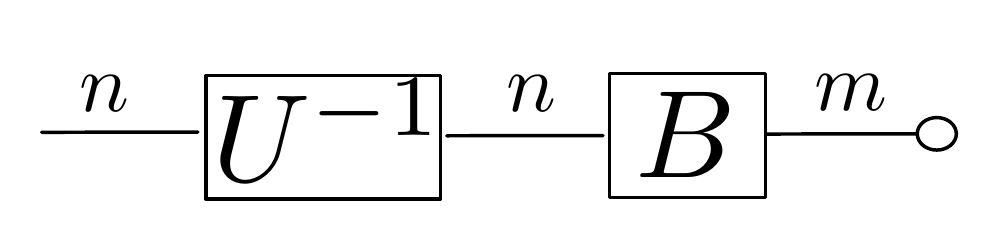}$}
\\
&\eql{Lemma \ref{lemma:invertiblestar}}&
\lower11pt\hbox{$\includegraphics[height=1cm]{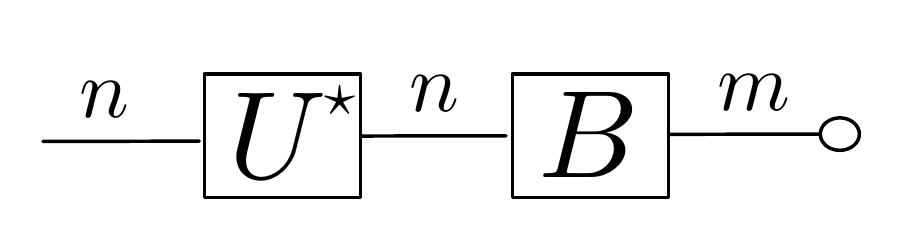}$} \\
&\eql{Lemma \ref{lemma:BHNFequalKernel}}&
\lower13pt\hbox{$\includegraphics[height=1.2cm]{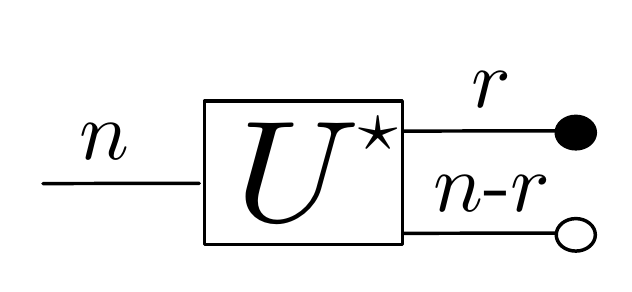}$} \\
&\eql{Prop. \ref{prop:star=refl}}&
\lower15pt\hbox{$\includegraphics[height=1.4cm]{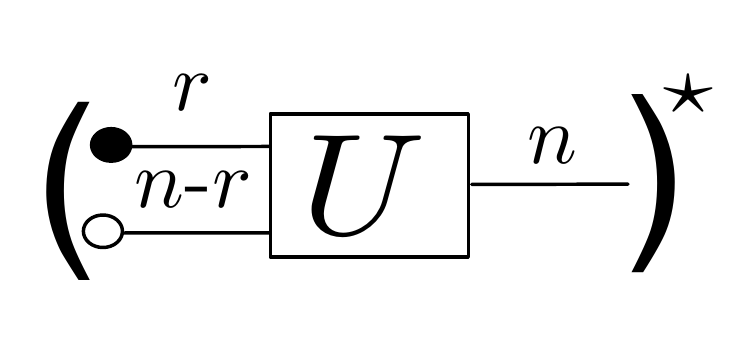}$} \\
&\eql{Lemma \ref{lemma:invertiblerestriction}}&
\lower15pt\hbox{$\includegraphics[height=1.4cm]{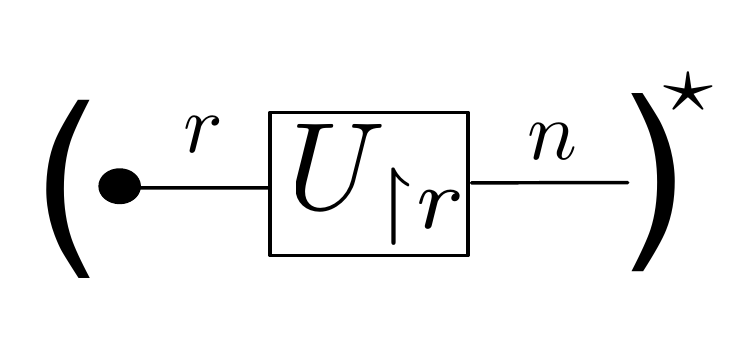}$} \\
&\eql{Prop. \ref{prop:star=refl}}&
\lower12pt\hbox{$\includegraphics[height=1.1cm]{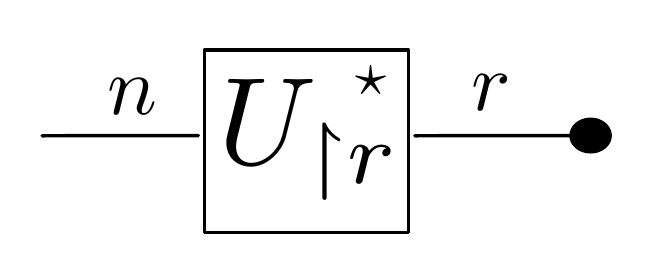}$}.
 \end{eqnarray*}
By Proposition~\ref{prop:kernelfromHNF}, the columns of the matrix $\restr{U}{r} \: r \to n$ yield a basis for the kernel of $A$. Thus $n\tl{\restr{U}{r}} r \tr{\finVect} 0$ is also a pullback span in \eqref{eq:kernelplback} and, since $\sem{\ABR}(\finVect \: r \to 0) = \circuitrbcounits$, we know by Lemma~\ref{lemma:mirror} that
\begin{equation*}
\lower12pt\hbox{$\includegraphics[height=1.1cm]{graffles/circuitUrestrrcounits.pdf}$}
=
\lower12pt\hbox{$\includegraphics[height=1.1cm]{graffles/circuitKerAstar.pdf}$} \end{equation*}
which concludes the proof of our statement.
\end{proof} 


We now have all the ingredients to provide a proof of our completeness statement, from which the characterization result of Theorem~\ref{th:Span=IBw} follows.

\begin{proof}[Proof of Proposition~\ref{prop:IBwComplete}]

Let $A,B,C,D$ be as in the statement of Proposition~\ref{prop:IBwComplete} and consider the following derivation in $\IBRw$:
\begin{eqnarray}\label{eq:dercompl}
\nonumber \lower20pt\hbox{$\includegraphics[height=1.7cm]{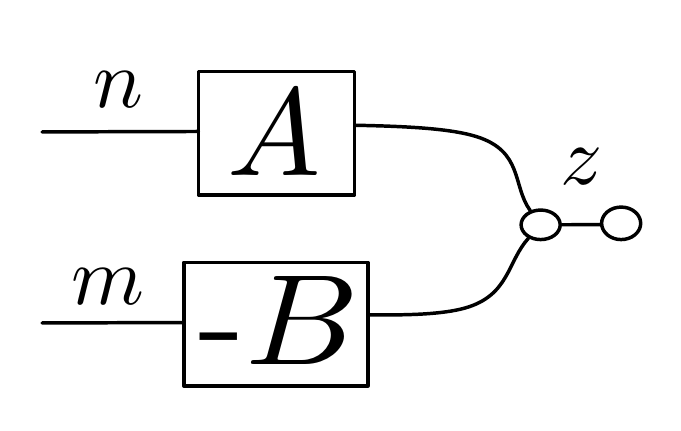}$}
&\eql{Def. $A|\minus B$}&
\lower10pt\hbox{$\includegraphics[height=1cm]{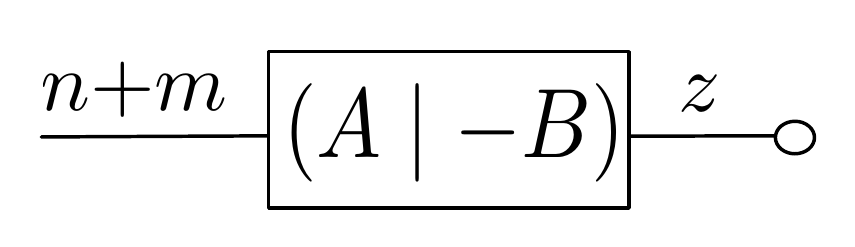}$} \\ \nonumber
&\eql{Prop.~\ref{prop:matrixequalkenrel}}&
\lower10pt\hbox{$\includegraphics[height=1cm]{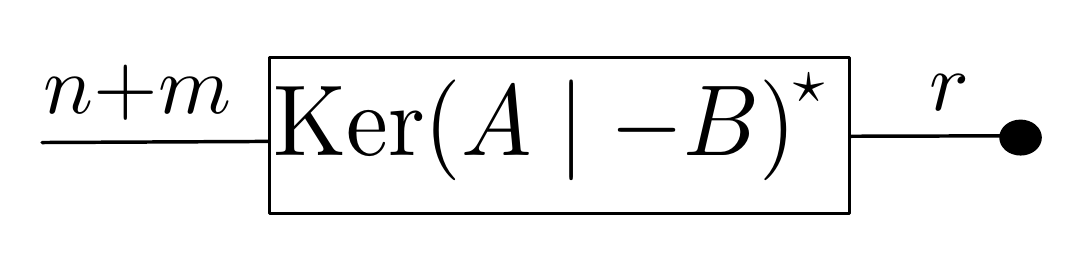}$} \\
&\eql{Lemma~\ref{lemma:pbKernel}}&
\lower15pt\hbox{$\includegraphics[height=1.2cm]{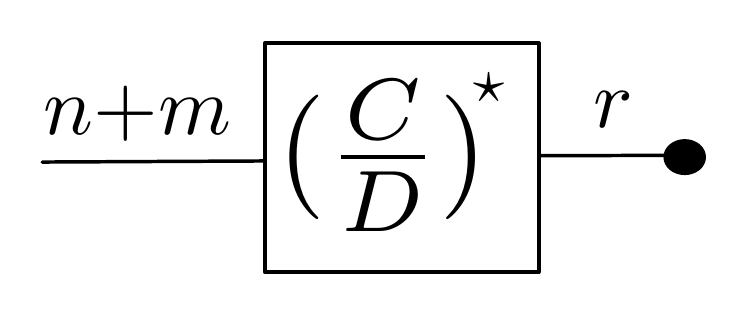}$} \\ \nonumber
&\eql{Def. $\left(\frac{C}{D}\right)$}&
\lower20pt\hbox{$\includegraphics[height=1.7cm]{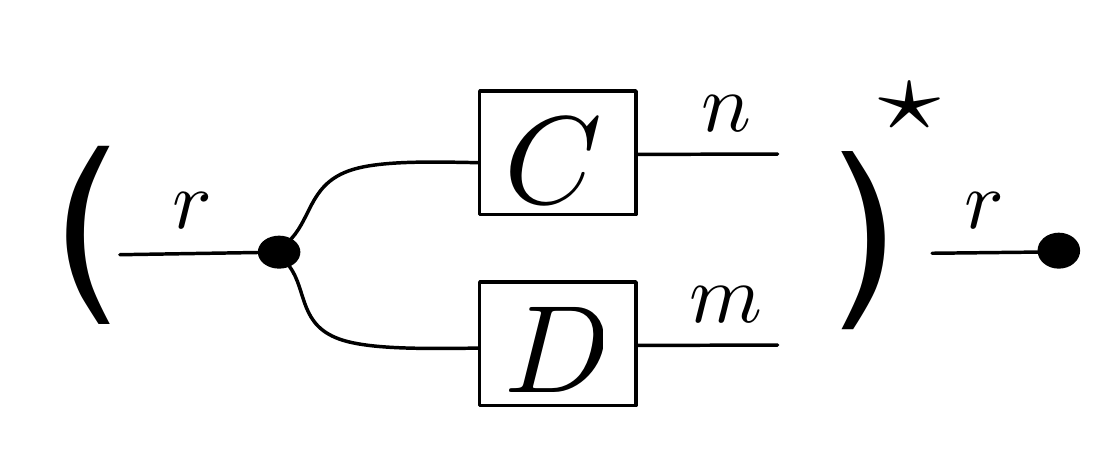}$} \\ \nonumber
&\eql{Prop. \ref{prop:star=refl}}&
\lower19pt\hbox{$\includegraphics[height=1.6cm]{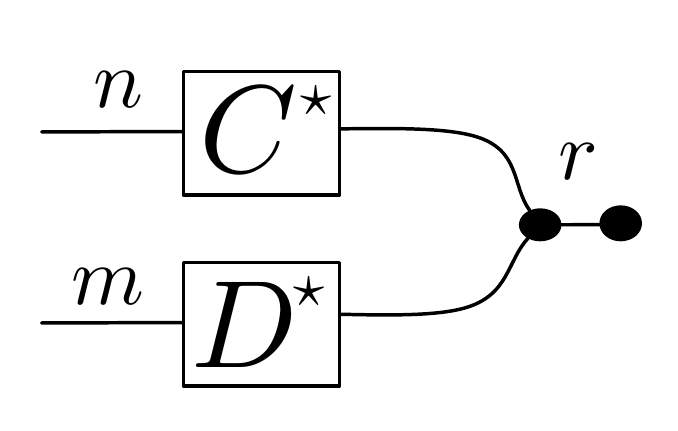}$}
 \end{eqnarray}
The proof is concluded by the following derivation, yielding the desired equation in $\IBRw$.
\begin{eqnarray*}
\lower10pt\hbox{$\includegraphics[height=1cm]{graffles/circuitABstar.pdf}$}
&\eql{\eqref{eq:defstar}}&
\lower18pt\hbox{$\includegraphics[height=2cm]{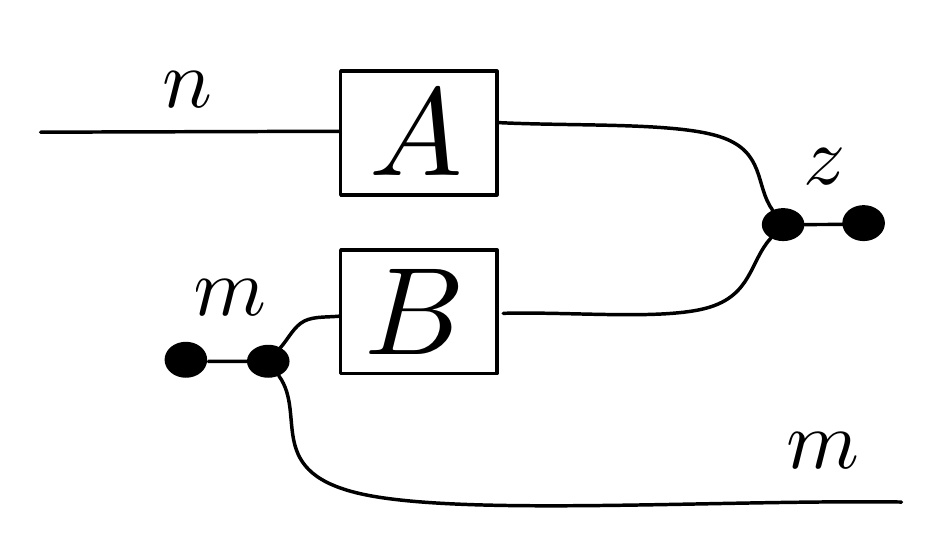}$} \\
&\eql{\eqref{eq:lccb}}&
\lower18pt\hbox{$\includegraphics[height=2cm]{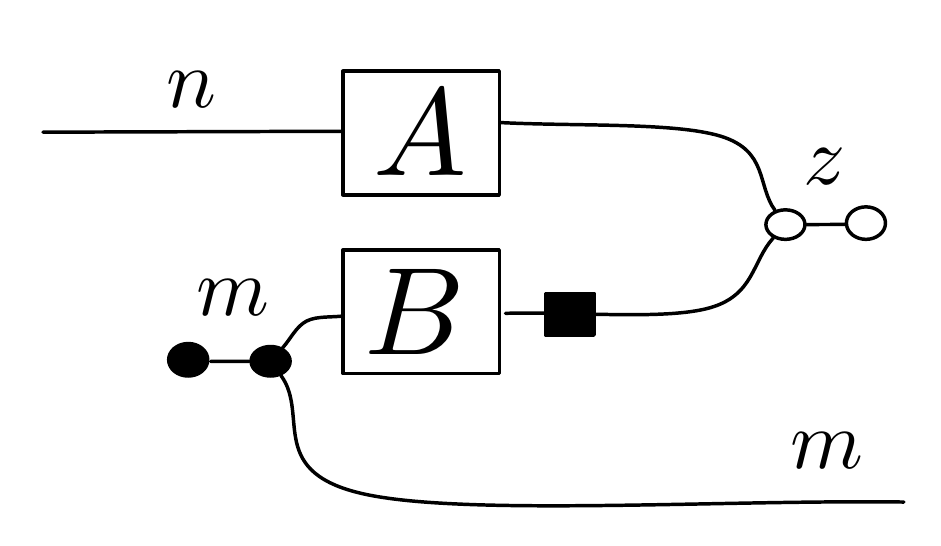}$} \\
&\eql{\eqref{eq:scalarwmult},\eqref{eq:scalarmult}}&
\lower18pt\hbox{$\includegraphics[height=2cm]{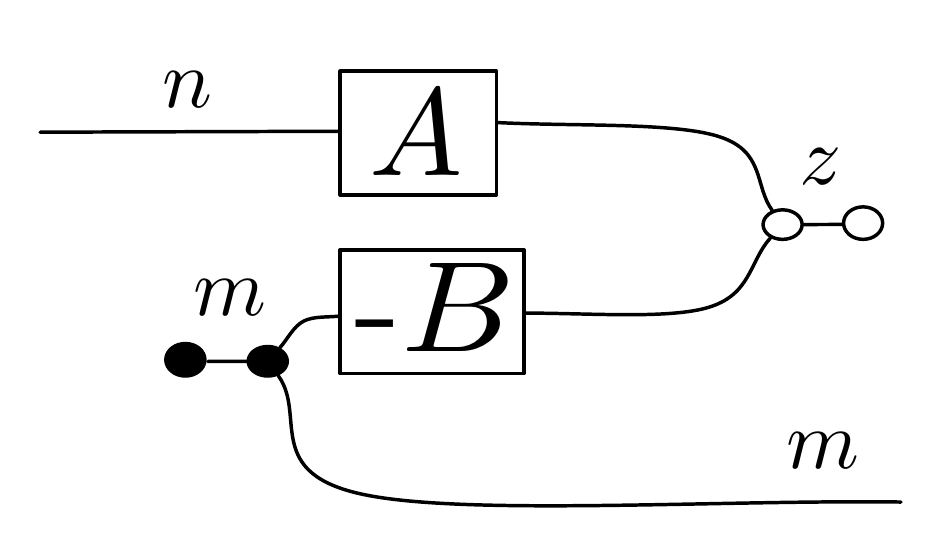}$} \\
&\eql{\eqref{eq:dercompl}}&
\lower18pt\hbox{$\includegraphics[height=2cm]{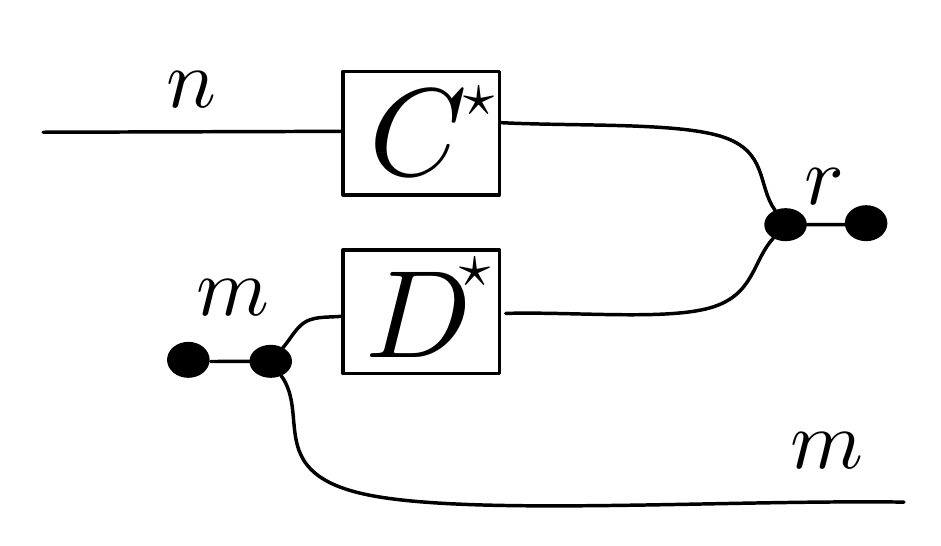}$} \\
&\eql{\eqref{eq:defstar}}&
\lower10pt\hbox{$\includegraphics[height=1cm]{graffles/circuitCstarDr.pdf}$}
 \end{eqnarray*}
We detail the various derivation steps. First, we can ``bend'' our string diagram using the compact-closed structure $\coc{\cdot}$. Then we iteratively apply~\eqref{eq:lccb} to turn the rightmost part of the compact-closed structure from black into white. This produces $z$ copies of the antipode \antipodesquare. The third equality is given by iteratively applying~\eqref{eq:scalarwmult} to push the antipodes in front of each scalar in diagram $B$, and then multiply all those scalars by the antipode value $-1$ using~\eqref{eq:scalarmult}. As a result, we obtain (a string diagram representing) the matrix $\minusB$. Then we can conclude using derivation~\eqref{eq:dercompl}.\end{proof}

This concludes the proof of Theorem~\ref{th:Span=IBw}. As an immediate consequence, we obtain the following factorisation property. Recall that $\HAopToIHw \: \ABRop \to \IBRw$ and $\HAToIHw \: \ABR \to \IBRw$ are the evident inclusions of theories, introduced after Definition~\ref{def:IBRw} above. Also, as expected, an \emph{invertible} arrow $h \: n \to m$ of $\ABR$ is one for which there exists $h^{-1} \: m \to n$ such that $h \poi h^{-1} = \id_n$ and $h^{-1} \poi h = \id_m$ --- equivalently, $\ABR \tr{\cong} \VectR$ maps $h$ to an invertible matrix.

\begin{corollary}\label{cor:factorisationIBRw} For any string diagram $c$ of $\IBRw$ there exist string diagrams $c_1$ of $\ABRop$ and $c_2$ of $\ABR$ such that $c = \HAopToIHw(c_1) \poi \HAToIHw(c_2)$.

Also, this factorisation is unique up-to isomorphism of spans. That means, for any other two string diagrams $d_1$ of $\ABRop$ and $d_2$ of $\ABR$ such that $c = \HAopToIHw(d_1) \poi \HAToIHw(d_2)$, there exists $h$ invertible making the following diagram commute in $\ABR$.
\[
\xymatrix@C=50pt{
& \ar[dl]_<<<<<<<<{\coc{c_1}} \ar[dr]^>>>>>>>>>>>>>{c_2}& \\
& \ar[l]_<<<<<<<<{\coc{d_1}} \ar[u]^{h}  \ar[r]^>>>>>>>>>>>>>{d_2} &
}
\]
\end{corollary}

The second part of Corollary~\ref{cor:factorisationIBRw} stems from the observation that an arrow of $\VectRop \bicomp{\PJ} \VectR$ is not a span, but rather an isomorphism class of spans --- see~\S~\ref{sec:distrLawPullback}. Thus the characterisation $\IBRw \cong \VectRop \bicomp{\PJ} \VectR$ of Theorem~\ref{th:Span=IBw} maps a string diagram to a unique span only up-to span isomorphism. This isomorphism --- an arrow of the core $\PJ$ of $\VectR$ --- is an invertible matrix, which corresponds through $\VectR \cong \ABR$ to the invertible arrow $h$ of $\ABR$ appearing in Corollary~\ref{cor:factorisationIBRw}.

\subsection{$\IBRb$: the Theory of Cospans of $\PID$-matrices}\label{sec:IBRbCospan}

In this section we provide a presentation by generators and equations of $\CospanMat$, the PROP of cospans of $\PID$-matrices. The key point is that the matrix transpose operation on $\VectR$ (see Rmk.~\ref{rmk:graphicaltranspose}) can be used to map spans to cospans. Our strategy will be to understand the transpose in graphical terms: starting from the string diagrammatic theory $\IBRw$ of $\SpanMat$, this will give ``for free'' also the theory of $\CospanMat$. First, we introduce the PROP claimed to characterise $\CospanMat$.

\begin{definition}\label{def:IBRb} The PROP $\IBRb$ is the quotient of $\ABR + \ABRop$ by the following equations, for $k$ any element and $l$ any non-zero element of $\PID$.
\begin{multicols}{2}\noindent
 \begin{equation}
\label{eq:lcmop}
\tag{B1}
\lower7.5pt\hbox{$\includegraphics[height=.75cm]{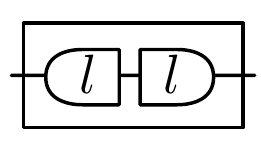}$}
=
\lower6pt\hbox{$\includegraphics[height=.6cm]{graffles/idcircuit.pdf}$}
\end{equation}
\begin{equation}
\label{eq:bbone}
\tag{B2}
\lower6pt\hbox{$\includegraphics[height=.6cm]{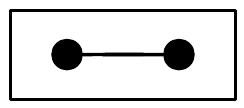}$}
=  \lower4pt\hbox{$\includegraphics[height=.5cm]{graffles/idzerocircuit.pdf}$}
\end{equation}
\end{multicols}
 \begin{multicols}{2}\noindent
\begin{equation}\tag{B3}
\lower15pt\hbox{$\includegraphics[height=1.3cm]{graffles/BFrobS.pdf}$}
\!\!
=
\!\!
\lower11pt\hbox{$\includegraphics[height=1cm]{graffles/BFrobX.pdf}$}
\!\!
=
\!\!
\lower15pt\hbox{$\includegraphics[height=1.3cm]{graffles/BFrobZ.pdf}$}
\end{equation}
\begin{equation}\tag{B4}
\lower15pt\hbox{$\includegraphics[height=1.3cm]{graffles/WFrobS.pdf}$}
\!\!
=
\!\!
\lower11pt\hbox{$\includegraphics[height=1cm]{graffles/WFrobX.pdf}$}
\!\!
=
\!\!
\lower15pt\hbox{$\includegraphics[height=1.3cm]{graffles/WFrobZ.pdf}$}
\end{equation}
\end{multicols}
 \begin{multicols}{2}\noindent
\begin{equation}\tag{B5}
\lower7pt\hbox{$\includegraphics[height=.7cm]{graffles/bccr.pdf}$}
=
\!\!
\lower11pt\hbox{$\includegraphics[height=1cm]{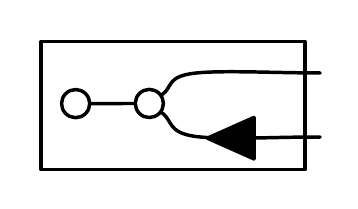}$}
\end{equation}
\begin{equation}\tag{B6}
\lower7pt\hbox{$\includegraphics[height=.7cm]{graffles/bccl.pdf}$}
=
\!\!
\lower11pt\hbox{$\includegraphics[height=1cm]{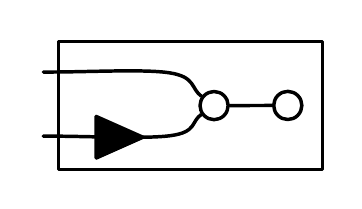}$}
\end{equation}
\end{multicols}
\begin{multicols}{2}\noindent
\begin{equation}\tag{B7}
\label{eq:WcccoscalarAxiomOne}
\lower9pt\hbox{$\includegraphics[height=.8cm]{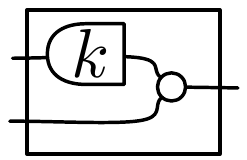}$}
=
\lower9pt\hbox{$\includegraphics[height=.8cm]{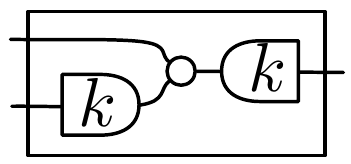}$}
\end{equation}
\begin{equation}\tag{B8}
\label{eq:WcccoscalarAxiomTwo}
\lower9pt\hbox{$\includegraphics[height=.8cm]{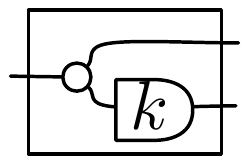}$}
=
\lower9pt\hbox{$\includegraphics[height=.8cm]{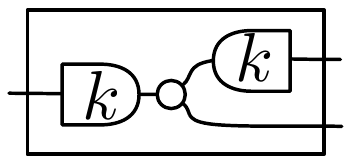}$}
\end{equation}
\end{multicols}
\end{definition}

Similarly to the case of $\IBRw$, we write $\tau_1 \: \ABR \to \IBRb$ and $\tau_2 \: \ABRop \to \IBRb$ for the PROP morphisms interpreting string diagrams of $\ABR$ and $\ABRop$, respectively, as string diagrams of $\IBRb$.

\medskip

The equations of $\IBRb$ are the \emph{photographic negative} of the ones of $\IBRw$, that is, they are the same modulo swapping the black and white colors (and the orientation of diagrams representing elements of $\PID$). More formally, we inductively define a PROP morphism $\pn \: \IBRb \to \IBRw$ by the following mapping.
  \begin{multicols}{4}
\noindent
      \begin{eqnarray*}
     \Bcounit \mapsto \Wcounit
    \end{eqnarray*}
   \begin{eqnarray*}
  \hspace{-.7cm}   \Bunit \mapsto \Wunit
    \end{eqnarray*}
  \begin{eqnarray*}
   \hspace{-.7cm}  \Wunit \mapsto \Bunit
    \end{eqnarray*}
  \begin{eqnarray*}
   \Wcounit \mapsto \Bcounit
   \end{eqnarray*}
    \end{multicols}
     \begin{multicols}{4}
     \noindent
     \begin{eqnarray*}
     \Wmult \mapsto \Bmult
    \end{eqnarray*}
  \begin{eqnarray*}
  \hspace{-.7cm} \Wcomult \mapsto \Bcomult
   \end{eqnarray*}
        \begin{eqnarray*}
  \hspace{-.7cm}   \Bmult \mapsto \Wmult
    \end{eqnarray*}
  \begin{eqnarray*}
   \Bcomult \mapsto \Wcomult
   \end{eqnarray*}
      \end{multicols}
     \begin{multicols}{4}
     \noindent
           \begin{eqnarray*}
     \scalar \mapsto \coscalar
    \end{eqnarray*}
  \begin{eqnarray*}
 \hspace{-.7cm}  \coscalar \mapsto \scalar
   \end{eqnarray*}
     \begin{equation*}
     \lower2pt\hbox{\hspace{-1.4cm}$c\poi c' \mapsto \pn(c)\poi\pn(c')$}
   \end{equation*}
     \begin{equation*}
     \lower2pt\hbox{\hspace{-.7cm}$c\tns c' \mapsto \pn(c)\tns\pn(c')$}
   \end{equation*}
 \end{multicols}

\noindent The next lemma confirms that $\pn$ is well-defined.

\begin{lemma}\label{lemma:negativerespectequations} For all string diagrams $c,c'$ of $\IBRb$,  $c = c'$ in $\IBRb$ iff $\pn(c) = \pn(c')$ in $\IBRw$.
\end{lemma}
\begin{proof} Observe that the equations~\eqref{eq:lcmop}-\eqref{eq:WcccoscalarAxiomTwo} presenting $\IBRb$ are the image under $\pn$ of the equations~\eqref{eq:lcm}-\eqref{eq:BccscalarAxiomTwo} presenting $\IBRw$, that is, (W$i$) $c \feq c'$ corresponds to (B$i$) $\pn(c) \feq \pn(c')$. Thus the statement is also true for all the derived laws of the two theories.  \end{proof}

\begin{lemma} $\pn$ is an isomorphism of PROPs. \end{lemma}
\begin{proof} Fullness of $\pn$ is easily verified by induction on $c \in \IBRw$ and faithfulness follows by the ``only if'' direction of Lemma~\ref{lemma:negativerespectequations}.
 \end{proof}

We now specify the matrix counterpart of $\pn$. As observed in Remark~\ref{rmk:graphicaltranspose}, the operation of taking the transpose of a matrix yields a PROP isomorphism $\trasp{\cdot} \: \VectR \cong \VectRop$. 
This also induces a PROP morphism $\tra \: \SpanMat \to \CospanMat$ mapping $n \tl{A} z \tr{B} m$ into $n \tr{A^T} z \tl{B^T} m$. To see that this assignment is functorial, observe that pushouts in $\VectR$ --- giving composition in $\CospanMat$ --- can be calculated by transposing pullbacks of transposed matrices. Because $\trasp{\cdot}$ is an isomorphism, also $\tra$ is an isomorphism.

We can now obtain an isomorphism between $\IBRb$ and $\CospanMat$ as:
\begin{equation}\label{eq:IsoCospan}  \xymatrix{ \IBRb \ar[r]^{\pn}&  \IBRw \ar[r]^-{\cong} & \SpanMat \ar[r]^{\tra} & \CospanMat} \text{.}  \end{equation}

\begin{theorem}\label{th:IBRb=Cospan}
$\IBRb \cong \CospanMat$.
\end{theorem}

\section{Interacting Hopf Algebras II: Linear Subspaces}\label{sec:cubetop}

In this section we characterise by generators and equations the PROP of subspaces over the field of fractions of $\PID$.
\begin{definition}\label{def:sv} The PROP $\SVR$ of subspaces over a field $\frPID$ is defined as follows:
\begin{itemize}
\item arrows $n\to m$ are subspaces of $\frPID^n \times \frPID^m$, considered as a $\frPID$-vector space.
\item Composition is relational: given $V \: n \to z$, $W \: z \to m$,
\[
(\xx,\zz)\in V\poi W \quad\Leftrightarrow\quad \exists \yy.\; (\xx,\yy)\in V \wedge (\yy,\zz)\in W.
\]
\item The monoidal product is given by direct sum.
\item The symmetry $n + m \to m+n$ is the subspace $\{(\matrixTwoROneC{\xx}{\yy},\matrixTwoROneC{\yy}{\xx}) \mid \xx \in \field^n \wedge \yy \in \field^m\}$.
\end{itemize}
\end{definition}
\begin{convention}\label{conv:subspaceLinRel} By the way composition is defined in $\SVR$, it is legitimate to regard an arrow $V \: n \to m$ in $\SVR$ as a \emph{relation} between $\field^n$ and $\field^m$, which we call \emph{linear} to emphasise that it is also a sub-vector space. We shall use the terms subspace and linear relation interchangeably. The latter will exclusively appear in Chapter~\ref{chapter:SFG}, as we find it particularly convenient for the operational reading of diagrams.
%
\end{convention}

Recall that, for $\PID$ a principal ideal domain, its \emph{field of fractions} $\frPID$ is canonically constructed by letting elements of $\frPID$ be fractions $\frac{k_1}{k_2}$, where $k_1 , k_2 \in \PID$, $k_2 \neq 0$ and $\frac{k_1}{k_2}$ represents an equivalence class of the relation $(k_1,k_2) \sim (k_3,k_4)$ on pairs of elements of $\PID$ defined by \[(k_1,k_2) \sim (k_3,k_4) \text{ if } k_2, k_4 \neq 0 \text{ and }k_1 k_4 = k_3 k_2 \text{.}\]
Throughout this section, we reserve the notation $\frPID$ for the field of fractions of $\PID$. The equational presentation for $\SVR$ will be given by the following PROP.

\begin{definition}\label{def:IBR} The PROP $\IBR$ is the quotient of $\ABR + \ABRop$ by the following equations, for $l$ any non-zero element of $\PID$.
  \begin{multicols}{2}\noindent
 \begin{equation}\tag{I1} \label{eq:lcmIH}
\lower8pt\hbox{$\includegraphics[height=.75cm]{graffles/lcml_l.pdf}$}
=
\lower6pt\hbox{$\includegraphics[height=.6cm]{graffles/idcircuit.pdf}$}
\end{equation}
 \begin{equation}\tag{I2} \label{eq:lcmopIH}
\lower8pt\hbox{$\includegraphics[height=.75cm]{graffles/lcmopl_l.pdf}$}
=
\lower6pt\hbox{$\includegraphics[height=.6cm]{graffles/idcircuit.pdf}$}
\end{equation}
\end{multicols}
 \begin{multicols}{2}\noindent
\begin{equation}\tag{I3}\label{eq:WFrobIBR}
\lower15pt\hbox{$\includegraphics[height=1.3cm]{graffles/WFrobS.pdf}$}
\!\!
=
\!\!
\lower11pt\hbox{$\includegraphics[height=1cm]{graffles/WFrobX.pdf}$}
\!\!
=
\!\!
\lower15pt\hbox{$\includegraphics[height=1.3cm]{graffles/WFrobZ.pdf}$}
\end{equation}
\begin{equation}\tag{I4}\label{eq:BFrobIBR}
\lower15pt\hbox{$\includegraphics[height=1.3cm]{graffles/BFrobS.pdf}$}
\!\!
=
\!\!
\lower11pt\hbox{$\includegraphics[height=1cm]{graffles/BFrobX.pdf}$}
\!\!
=
\!\!
\lower15pt\hbox{$\includegraphics[height=1.3cm]{graffles/BFrobZ.pdf}$}
\end{equation}
\end{multicols}
 \begin{multicols}{2}\noindent
\begin{equation}\tag{I5}\label{eq:lccIBR}
\lower12pt\hbox{$\includegraphics[height=1cm]{graffles/lccr.pdf}$}
\!\!
=
\!\!
\lower12pt\hbox{$\includegraphics[height=1cm]{graffles/blackcceta.pdf}$}
\end{equation}
\begin{equation}\tag{I6}\label{eq:rccIBR}
\lower12pt\hbox{$\includegraphics[height=1cm]{graffles/rccl.pdf}$}
\!\!
=
\!\!
\lower12pt\hbox{$\includegraphics[height=1cm]{graffles/blackccepsilon.pdf}$}
\end{equation}
\end{multicols}
\begin{multicols}{2}\noindent
\begin{equation}\tag{I7}\label{eq:WSepIBR}
\lower7pt\hbox{$\includegraphics[height=.7cm]{graffles/WSep.pdf}$}
=
\lower6pt\hbox{$\includegraphics[height=.6cm]{graffles/idcircuit.pdf}$}
\end{equation}
\begin{equation}\tag{I8}\label{eq:BSepIBR}
\lower7pt\hbox{$\includegraphics[height=.7cm]{graffles/BSep.pdf}$}
=
\lower6pt\hbox{$\includegraphics[height=.6cm]{graffles/idcircuit.pdf}$}
\end{equation}
\end{multicols}
\end{definition}

\begin{remark} In the case in which the PID under consideration is actually a field, we can replace \eqref{eq:lcmIH} and \eqref{eq:lcmopIH} by the following condition, for all $l \neq 0$: 
\begin{equation}\label{eq:flipscalar}\tag{Inv}
\lower8pt\hbox{$\includegraphics[height=.7cm]{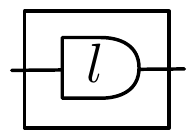}$} = \lower8pt\hbox{$\includegraphics[height=.7cm]{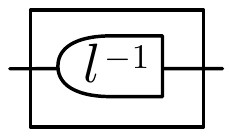}$}.
\end{equation}
\end{remark}

The idea is that the equations of $\IBR$ are obtained by merging the equational theories of $\IBRw$ and $\IBRb$, while identifying the common structure $\ABR + \ABR^{\op}$. This modular perspective is captured by expressing $\IBR$ as a fibered sum (\emph{cf.} \S\ref{sec:pushout}) of $\IBRw$ and $\IBRb$. First, consider the following diagram in $\PROP$.
 \begin{equation}\label{eq:topface}
 \tag{Top}
\vcenter{
\xymatrix@=20pt{ & \ar[dl]_{[\HAToIHw,\HAopToIHw]} \ABR+\ABRop \ar[dd]|{[\varphi_1,\varphi_2]} \ar[dr]^{[\tau_1,\tau_2]} & \\
 \IBRw \ar[dr]_{\Theta} & & \ar[dl]^{\Lambda} \IBRb \\
 & \IBR &
 }
 }
 \end{equation}
The PROP morphism $\Theta$ quotients $\IBRw$ by the equations of $\IBRb$ and $\Lambda$ quotients $\IBRb$ by the ones of $\IBRw$.  For later use, it is useful to also fix the PROP morphism $[\varphi_1,\varphi_2] \: \ABR + \ABRop \to \IBR$ defined by commutativity of the diagram.

\begin{proposition}\label{prop:topfacepushout} \eqref{eq:topface} is a pushout in $\PROP$.
\end{proposition}
\begin{proof}
Proposition~\ref{prop:SMTpushout} gives us the conditions on the signature and equations of $\IBR$ to prove the statement. For the signature, there is nothing to prove, because all the four PROPs in~\eqref{eq:topface} have the same. For the equations, we need to check the following: an equation is provable in $\IBR$ if and only if it is provable from the union of the axioms of $\IBRb$ and $\IBRw$. For the left to right direction, note that: equations~\eqref{eq:WFrobIBR}-\eqref{eq:rccIBR} are provable in both theories; \eqref{eq:lcmIH} and \eqref{eq:BSepIBR} are provable in $\IBRw$; \eqref{eq:lcmopIH} and \eqref{eq:WSepIBR} are provable in $\IBRb$. For the converse direction, observe that the only axioms of the two theories that are missing in the presentation of $\IBR$ are \eqref{eq:wbone}, \eqref{eq:BccscalarAxiomOne}, \eqref{eq:BccscalarAxiomTwo} (from $\IBRw$), \eqref{eq:bbone}, \eqref{eq:WcccoscalarAxiomOne} and \eqref{eq:WcccoscalarAxiomTwo} (from $\IBRb$). They are all provable in $\IBR$: we refer to Appendix~\ref{AppDerLawsIH} for the derivations.
\end{proof}

It is worth emphasizing the following consequence of the modular account of $\IBR$, that will play a fundamental role in our developments.
\begin{theorem}[Factorisation of $\IBR$] \label{Th:factIBR} Let $c \in \IBR[n,m]$ be a string diagram.
 \begin{enumerate}[(i)]
   \item There exist $c_1$ in $\ABRop$ and $c_2$ in $\ABR$ such that $c = \varphi_2(c_1) \poi\varphi_1(c_2)$.
   \item  There exist $c_3$ in $\ABR$ and $c_4$ in $\ABRop$ such that $c = \varphi_1(c_3) \poi \varphi_2(c_4)$.
 \end{enumerate}
\end{theorem}
\begin{proof} The main observation is that, since $\IBR$ is both a quotient of $\IBRw$ and of $\IBRb$, it inherits their factorisation property. The first statement of Theorem~\ref{Th:factIBR} follows by the factorisation result for $\IBRw$ (Corollary \ref{cor:factorisationIBRw}). Since $\IBRb$ has been shown to be isomorphic to $\ABR \bicomp{\Perm} \ABRop$, then the dual of Corollary \ref{cor:factorisationIBRw} holds for $\IBRb$. This yields the second statement of Theorem~\ref{Th:factIBR}.
\end{proof}

Theorem~\ref{Th:factIBR} states that any string diagram of $\IBR$ can be put in the form of a span and of a cospan of $\ABR$-diagrams. When interpreting diagrams of $\IBR$ as subspaces, we shall see that these two canonical forms correspond to well-known ways of representing a subspace, namely as a basis and as a system of linear equations (Example~\ref{ex:subspacesgraphical}).

\begin{remark}
Note that, differently from the case of $\IBRw$ (Corollary~\ref{cor:factorisationIBRw}), the span factorisation for $\IBR$ is \emph{not} unique up-to isomorphism in $\VectR$. Consider for instance $\Wcomult \poi \Wmult$ and $\idcircuit \poi \idcircuit$: they correspond to non-isomorphic spans in $\VectR$ --- and, indeed, different string diagrams of $\IBRw$ --- but are identified in $\IBR$ by the axiom~\eqref{eq:WSepIBR}. The same observation holds for the cospan factorisation: \eqref{eq:BSepIBR} shows non-isomorphic cospans that are identified as diagrams of $\IBR$.
\end{remark}

We now state the main result of this section.

\begin{theorem}\label{th:IBR=SVR} 
$\IBR \cong \SVR$.
\end{theorem}
The proof will be reminiscent of the argument for equivalence relations (Theorem~\ref{th:IFROB=ER}) and partial equivalence relations (Theorem~\ref{th:PIFROB=PER}). We shall construct the cube~\eqref{eq:cube} introduced in~\S~\ref{sec:overview}. We already showed that the top face~\eqref{eq:topface} is a pushout. We next prove that the bottom face is also a pushout (\S~\ref{sec:cubebottom}). Then, in \S~\ref{sec:cubeback}, we show commutativity of the rear faces, whose vertical arrows are isomorphisms.  The isomorphism $\IBR \to \SVR$ will then be given by universal properties of the top and bottom faces (\S~\ref{sec:cuberebuilt}).

\subsection{The Cube: Bottom Face}\label{sec:cubebottom}

In this section we show that the following diagram, which is the bottom face of the cube \eqref{eq:cube}, is a pushout in $\PROP$. The argument is similar to the one used for showing that $\ER$ is a fibered sum (see~\S\ref{sec:ER}).
\begin{equation}\tag{Bot}
\label{eq:bottomface}
\raise15pt\hbox{$
\xymatrix@C=40pt{
{\VectR + \VectRop} \ar[r]^-{[\MatToSpan,\,\MatopToSpan]} \ar[d]_{[\MatToCospan,\,\MatopToCospan]} & {\SpanMat} \ar[d]^{\Phi} \\
{\CospanMat} \ar[r]_-{\Psi} & {\SVR}
}$}
\end{equation}
In the diagram above, we define
\[\MatToSpan(A \: n\to m) = (n \tl{\id} n \tr{A} n),\  \MatopToSpan(A \: n\to m) = (n \tl{A} m \tr{\id} m ),\]
\[\MatToCospan(A \: n\to m) = (n \tr{A} m \tl{\id} m)\text{ and }
\MatopToCospan(A \: n\to m) = (n \tr{\id} n \tl{A} m).\]
For the definition of $\Phi$, we let $\Phi(n\tl{A} z \tr{B} m)$ be the subspace
\[
\{\,(\xx,\yy)\ |\ \xx\in \frPID^n,\, \yy\in \frPID^m,\, \exists \zz\in \frPID^z.\; A\zz=\xx \wedge B\zz = \yy \,\}.
\]
$\Psi(n \tr{A} z \tl{B} m)$ is defined to be the subspace \[
\{\, (\xx,\yy) \ | \ \xx\in \frPID^n,\, \yy\in \frPID^m,\ A\xx = B\yy\, \}.
\]
In the sequel we verify that $\Phi$ and $\Psi$ are indeed functorial assignments. This requires some preliminary work.
Recall SMCs $\RMod{\PID}$ and $\FRMod{\PID}$ introduced in \S\ref{sec:completeness}. We define $\RMod{\frPID}$ and $\FRMod{\frPID}$ analogously, as the SMCs of finite-dimensional $\frPID$-vector spaces and of finite-dimensional free $\frPID$-vectors spaces respectively (note that, of course, $\FRMod{\frPID} \cong \RMod{\frPID}$).
There is an obvious PROP morphism $I \: \VectR \to \Mat{\frPID}$ interpreting a matrix with entries in $\PID$ as one with entries in $\frPID$. Similarly, we have an inclusion $J \: \FRMod{\PID} \to \FRMod{\frPID}$. This yields the following commutative diagram, where $\simeq$ stands for equivalence.
\[
\xymatrix{
{\Mat{\PID}} \ar[d]_I \ar[r]^-{\simeq} & {\FRMod{\PID}} \ar[d]^J \\
{\Mat{\frPID}} \ar[r]_-{\simeq} & {\FRMod{\frPID}}
}
\]

\begin{lemma} \label{lem:pullbacksinmat} $I \: \Mat{\PID}\to \Mat{\frPID}$ preserves pullbacks and pushouts.
\end{lemma}
\begin{proof} Because both $\Mat{\PID}$ and $\Mat{\frPID}$ are self-dual (Rmk.~\ref{rmk:graphicaltranspose}) and $I \: \Mat{\PID}\to\Mat{\frPID}$ respects duality, it preserves pullbacks iff it preserves pushouts. It is thus enough to show that it preserves pullbacks. This can be easily be proved directly as follows. Suppose that the diagram
\begin{equation}
\label{eq:pbinmatr}
\tag{$\star$}
\vcenter{
\xymatrix@=15pt{
{r} \ar[d]_A \ar[r]^B& {m}\ar[d]^D\\
{n} \ar[r]_C  & {z}
}
}
\end{equation}
is a pullback in $\Mat{\PID}$. We need to show that it is also a pullback in $\Mat{\frPID}$.
Suppose that, for some $P \: q\to n$, $Q \: q\to m$ in $\Mat{\frPID}$ we have that
$P \poi C = Q \poi D$ in $\Mat{\frPID}$, that is, by definition of composition, $CP$ and $DQ$ are the same matrix. Since $\PID$ is a PID we can find least common multiples:
thus let $d$ be a common multiple of all the denominators that appear in $P$ and $Q$.
Then $dP \: q \to n$, $dQ \: q\to m$ are in $\Mat{\PID}$ and we have
$dp \poi C = C(dP)=d(CP)=d(DQ)=D(dQ) = dQ \poi D$. Since \eqref{eq:pbinmatr} is a pullback in $\Mat{\PID}$, there exists a unique $H \: q\to r$ with $H \poi A = AH =dP$ and $H \poi B = BH = dQ$. This means that we have found a mediating arrow, $H/d \: q \to r$, in $\Mat{\frPID}$ since $H/d \poi A = A(H/d)=AH/d=dP/d=P$ and similarly $H/d \poi B =Q$. Uniqueness in $\Mat{\frPID}$ can also be translated in a straightforward way to uniqueness in $\Mat{\PID}$. Basically if $H'$ is another mediating morphism and $d'$ is the least common multiple of denominators in $H'$ then we must have $d'(H/d)=d'H'$ because of the universal property in $\Mat{\PID}$. Dividing both sides by $d'$ yields the required equality.
\end{proof}

We are now able to show that
\begin{lemma}
$\Phi\: \SpanMat \to \SVR$ is a PROP morphism.	
\end{lemma}
\begin{proof}
We must verify that $\Phi$ preserves composition. In the diagram below let the centre square be a pullback diagram in $\VectR$.
\[
\xymatrix{
& & {r} \ar[dl]_{F_2'} \ar[dr]^{G_1'} \\
& {z_1} \ar[dl]_{F_1} \ar[dr]^{G_1} & & {z_2} \ar[dl]_{F_2} \ar[dr]^{G_2} \\
{n} & & {z} & & {m}
}
\]
By definition of composition in $\SpanMat$, $(\tl{F1} \tr{G_1}) \poi (\tl{F_2}\tr{G_2})
=\ \tl{F_1F_2'}\tr{G_2G_1'}$.

Now, by definition, if $(\xx,\zz)\in \Phi(\tl{F_1F_2'}\tr{G_2G_1'})$
then there exist $\ww$ with $\xx= F_1F_2'\ww$ and $\zz=G_2G_1'\ww$.
Therefore
$(\xx,\zz)\in {\Phi(\tl{F1} \tr{G_1})}\poi
{\Phi(\tl{F_2}\tr{G_2})}$ by commutativity of the square.

Conversely, if $(\xx,\zz)\in \Phi(\tl{F_1} \tr{G_1})\poi
\Phi(\tl{F_2}\tr{G_2})$ then
for some $\yy$ we must have
$(\xx,\yy)\in \Phi(\tl{F_1} \tr{G_1})$ and
$(\yy,\zz)\in \Phi(\tl{F_2}\tr{G_2})$.
Thus there exists $\uu$ with
$\xx=F_1\uu$ and $\yy=G_1\uu$
and there exists $\vv$ with $\yy=F_2\vv$
and $\zz=G_2\vv$.
By Lemma~\ref{lem:pullbacksinmat}, the square is also a pullback in $\Mat{\frPID}$ and then it translates to a pullback diagram in $\VectSpFr$. It follows the existence of $\ww$ with $F_2'\ww=\uu$
and $G_1'\ww=\vv$: thus $(\xx,\zz)\in \Phi((\tl{F1} \tr{G_1})\poi (\tl{F_2}\tr{G_2}))$.
This completes the proof.
\end{proof}

The proof that $\Psi$ is also a functor will rely on the following lemma.

\begin{lemma}\label{lemma:pushoutinvect}
Let the following be a pushout diagram in $\FRMod{\frPID}$.
\[
\xymatrix@=15pt{
 {U} \ar[d]_f \ar[r]^g &  {W} \ar[d]^q\\
{V}\ar[r]_p & {T}
}
\]	
Suppose that there exist $\vv\in V$, $\ww\in W$ such that $p\vv=q\ww$.
Then there exists $\uu\in U$ with $f\uu=\vv$ and $g\uu=\ww$.
\end{lemma}
\begin{proof}
Pushouts in $\FRMod{\frPID} \cong \RMod{\frPID}$ can be constructed by quotienting the vector space $V+W$ by the subspace generated by $\{\,(f\uu,g\uu)\,|\,\uu\in U\,\}$. Thus,
if $p(\vv) = q(\ww)$ then there exists a chain $\uu_1,\uu_2,\dots,\uu_k$
with $f(\uu_1)=\vv$, $g(\uu_1)=g(\uu_2)$, $f(\uu_2)=f(\uu_3)$, \dots, $f(\uu_{k-1})=f(\uu_{k-1})$ and $g(\uu_k)=\ww$.
If $k=1$ then we are finished. Otherwise, to construct an inductive argument we need to consider a chain $\uu_1, \uu_2, \uu_3$ with $f(\uu_1)=\vv$,
$g(\uu_1)=g(\uu_2)$, $f(\uu_2)=f(\uu_3)$ and $g(\uu_3)=\ww$.
Now $f(\uu_1-\uu_2+\uu_3)=f(\uu_1)-f(\uu_2)+f(\uu_3)=\vv$ and
$g(\uu_1-\uu_2+\uu_3)=g(\uu_1)-g(\uu_2)+g(\uu_3)=\ww$, so we have reduced the size of the chain to one.
\end{proof}

\begin{lemma}
$\Psi\: \CospanMat \to \SVR$ is a PROP morphism.
\end{lemma}
\begin{proof}
We must verify that $\Psi$ preserves composition. Let the square in the diagram below be a pushout in $\Mat{\PID}$. By definition of composition in $\CospanMat$
we have $(\tr{P_1}\tl{Q_1})\poi(\tr{P_2}\tl{Q_2})
=\ \tr{R_1P_1}\tl{R_2Q_2}$.
\[
\xymatrix{
{n} \ar[dr]_{P_1} & & {z} \ar[dl]^{Q_1} \ar[dr]_{P_2} & & {m} \ar[dl]^{Q_2} \\
& {z_1} \ar[dr]_{R_1} & & {z_2}\ar[dl]^{R_2} \\
& & {r}
}
\]
Consider $(\xx,\zz)\in \Psi(\tr{R_1P_1} \tl{R_2Q_2} )$.
Then $R_1P_1 \xx  = R_2 Q_2 \zz = \yy \in \frPID^r$.
Since the pushout diagram maps to a pushout diagram in $\VectSpFr$, we can use the conclusions of Lemma~\ref{lemma:pushoutinvect} to obtain $\yy\in \frPID^{z}$ such that $Q_1\yy=P_1\xx$ and $P_2\yy=Q_2\zz$. In other words, we have
$(\xx,\yy)\in\Psi(\tr{P_1}\tl{Q_1})$
and
$(\yy,\zz)\in\Psi(\tr{P_2}\tl{Q_1})$, meaning that $(\xx,\zz) \in \Psi(\tr{P_1}\tl{Q_1}) \poi\Psi(\tr{P_2}\tl{Q_1})$.

Conversely if $(\xx,\zz)\in\Psi(\tr{P_1}\tl{Q_1})\poi\Psi(\tr{P_2}\tl{Q_2})$ then $\exists \yy\in\frPID^{z}$ such that $(\xx,\yy)\in\Psi(\tr{P_1}\tl{Q_1})$
and $(\yy,\zz)\in\Psi(\tr{P_2}\tl{Q_2})$. It follows that $R_1P_1\xx=R_1Q_1\yy=R_2P_2\yy=R_2Q_2\zz$ and thus $(\xx,\zz)\in \Psi(\tr{R_1P_1} \tl{R_2Q_2} )$ as required.
\end{proof}

\begin{remark}\label{rmk:PushoutsMatrices} The proof of Lemma~\ref{lemma:pushoutinvect} relies on the fact that, for $\frPID$ a field, pushouts in $\FRMod{\frPID}$ coincide with those in $\RMod{\frPID}$. It would not work for an arbitrary PID $\PID$: $\FRMod{\PID}$ has pushouts for purely formal reasons, because it has pullbacks and is self-dual. However, differently from pullbacks (for which one can use that submodules of a free $\PID$-module are free, see Proposition~\ref{prop:PullbacksMatR}), pushouts are not calculated as in $\RMod{\PID}$. This asymmetry is the reason why proving functoriality of $\Psi$ requires more work than for $\Phi$.
\end{remark}

We now proceed in steps identifying the properties that make~\eqref{eq:bottomface} a pushout. As for the pushout characterisation of $\ER$ (\emph{cf.} Lemma~\ref{lemma:arbitraryPROP_ER}), the first step is to fix some useful conditions satisfied by any commutative diagram
\begin{equation}
\label{eq:arbitrary}
\raise15pt\hbox{$
\xymatrix@C=40pt{
{\VectR + \VectRop} \ar[r]^-{[\MatToSpan,\,\MatopToSpan]} \ar[d]_{[\MatToCospan,\,\MatopToCospan]} & {\SpanMat} \ar[d]^{\Gamma} \\
{\CospanMat} \ar[r]_-{\Delta} & {\mathbb{X}}.
}$}
\end{equation}

\begin{lemma}
\label{lemma:arbitraryPROP}
Let $\mathbb{X}$ be an arbitrary PROP, $\Delta$ and $\Psi$ PROP morphisms making \eqref{eq:arbitrary} commute. Consider the following diagram in $\Mat\PID$:
\begin{equation}
\label{eq:square}
\raise10pt\hbox{$
\xymatrix{
{} \ar[r]^{G} \ar[d]_{F} & {} \ar[d]^Q \\
{} \ar[r]_{P} & {}
}$}
\end{equation}
\begin{enumerate}[(i)]
\item if \eqref{eq:square} is a pushout diagram then
$\Gamma(\tl{F}\tr{G})=\Delta(\tr{P}\tl{Q})$.
\item if \eqref{eq:square} is a pullback diagram then
$\Gamma(\tl{F}\tr{G})=\Delta(\tr{P}\tl{Q})$.
\item if $\tl{F_1}\tr{G_1}$ and
$\tl{F_2}\tr{G_2}$ have the same pushout cospan
in $\Mat{\PID}$ then
$\Gamma(\tl{F_1}\tr{G_1})=\Gamma(\tl{F_2}\tr{G_2})$.
\item if $\tr{P_1}\tl{Q_1}$ and
$\tr{P_2}\tl{Q_2}$ have
the same pullback span in $\Mat{\PID}$ then
${\Delta(\tr{P_1}\tl{Q_1})}=
\Delta(\tr{P_2}\tl{Q_2})$.
\end{enumerate}
\end{lemma}
\begin{proof}~
\begin{enumerate}[(i)]
\item Suppose that $\tr{P}\tl{Q}$ is the cospan obtained by pushing out
$\tl{F}\tr{G}$ in $\Mat\PID$. Then
\begin{align*}
\Gamma(\tl{F}\tr{G}) &= \Gamma(\MatopToSpan F \poi \MatToSpan G) \\
&= \Gamma(\MatopToSpan F )\poi \Gamma(\MatToSpan G) \\
&= \Delta(\MatopToCospan F)\poi \Delta(\MatToCospan G) \\
&= \Delta( \MatopToCospan F \poi \MatToCospan G) \\
&= \Delta( \tr{P}\tl{Q} ).
\end{align*}
\item Suppose that $\tl{F}\tr{G}$ is the span obtained by pulling back
$\tr{P}\tl{Q}$. Then, reasoning in a similar way to (i), we get $\Delta(\tr{P}\tl{Q})= \Gamma(\tl{F}\tr{G})$.
\item Suppose that $\tr{P}\tl{Q}$ is the cospan obtained by pushing out
$\tl{F_1}\tr{G_1}$ and $\tl{F_2}\tr{G_2}$.
Using (i) we get $\Gamma(\tl{F_1}\tr{G_1})=\Delta(\tr{P}\tl{Q})=\Gamma(\tl{F_2}\tr{G_2})$.
\item The proof of (iv) is similar and uses (ii).
\end{enumerate}
\end{proof}

We now verify some properties of~\eqref{eq:bottomface}.
\begin{lemma}
\label{lemma:bottomfacecommutes}
\eqref{eq:bottomface} commutes.
\end{lemma}
\begin{proof}
It suffices to show that it commutes on the two injections into
$\VectR + \VectRop$. This means that we have to show, for any
$A \: n\to m$ in $\Mat{\PID}$, that
\[
\Phi(\tl{\id}\tr{A}) = \Psi(\tr{A}\tl{\id})
\]
and
\[
\Phi(\tl{A}\tr{\id}) =
\Psi(\tr{\id}\tl{A}).
\]
These are clearly symmetric, so it is enough to check one.
But this follows directly from the definition of $\Phi$ and $\Psi$:
\[
\Phi(\tl{\id}\tr{A}) =
\{\, (\xx,\yy) \,|\, A \xx = \yy \,\} = \Psi(\tr{A}\tl{\id})
\]
\end{proof}

%
%

\begin{lemma}
\label{lemma:pullbackpsi}
The following are equivalent:
\begin{enumerate}[(i)]
\item $n \tr{P_1} z_1 \tl{Q_1} m$ and $n \tr{P_2} z_2 \tl{Q_2}m$ have the same pullback in $\Mat{\PID}$.
\item $\Psi(\tr{P_1}\tl{Q_1})=
\Psi(\tr{P_2}\tl{Q_2})$.
\end{enumerate}
\end{lemma}
\begin{proof}
The conclusions of
Lemmas~\ref{lemma:bottomfacecommutes} and~\ref{lemma:arbitraryPROP}
give that (i) $\Rightarrow$ (ii). It thus suffices to show that
(ii) $\Rightarrow$ (i). Indeed, suppose that
$\Psi(\tr{P_1}\tl{Q_1})=
\Psi(\tr{P_2}\tl{Q_2})$. In particular
on elements $\xx\in\PID^n$, $\yy\in\PID^m$ we have $(\dag)$
$P_1\xx=Q_1\yy$ if and only if $P_2\xx = Q_2\yy$. Compute the following
pullbacks in $\Mat{\PID}$:
\[
\xymatrix{
{r_1} \ar[d]_{G_1} \ar[r]^{F_1} & {m} \ar[d]^{Q_1} \\
{n} \ar[r]_{P_1} & {z_1}
}
\qquad
\xymatrix{
{r_2} \ar[d]_{G_2} \ar[r]^{F_2} & {m} \ar[d]^{Q_2} \\
{n} \ar[r]_{P_2} & {z_2}
}
\]
By $(\dag)$ we can conclude that $P_1G_2=Q_1F_2$ and $P_2G_1=Q_2F_1$. This, using the universal property of pullbacks, implies that the spans $\tl{G_1}\tr{F_1}$ and
$\tl{G_2}\tr{F_2}$ are isomorphic.
\end{proof}

\begin{lemma}
\label{lemma:pushoutphi}
The following are equivalent:
\begin{enumerate}[(i)]
\item $n \tl{F_1} z_1 \tr{G_1} m$ and
$n \tl{F_2} z_2 \tr{G_2} m$ have the same pushout
in $\Mat{\PID}$
\item $\Phi(\tl{F_1}\tr{G_1})=
\Phi(\tl{F_2}\tr{G_2})$.
\end{enumerate}
\end{lemma}
\begin{proof}
The conclusions of
Lemmas~\ref{lemma:bottomfacecommutes} and~\ref{lemma:arbitraryPROP} again
give us that (i) $\Rightarrow$ (ii). It thus suffices to show that
(ii) $\Rightarrow$ (i).
Assume $\Phi(\tl{F_1}\tr{G_1})=
\Phi(\tl{F_2}\tr{G_2})$.
Compute the following
pushouts in $\Mat{\PID}$:
\[
\xymatrix{
{z_1} \ar[d]_{G_1} \ar[r]^{F_1} & {n} \ar[d]^{Q_1} \\
{m} \ar[r]_{P_1} & {r_1}
}
\qquad
\xymatrix{
{z_2} \ar[d]_{G_2} \ar[r]^{F_2} & {n} \ar[d]^{Q_2} \\
{m} \ar[r]_{P_2} & {r_2}
}
\]
By the conclusion of Lemma~\ref{lemma:arbitraryPROP}, we have
$\Psi(\tr{P_1}\tl{Q_1})=
\Psi(\tr{P_2}\tl{Q_2})$.
Applying the conclusion of Lemma~\ref{lemma:pullbackpsi},
$\tr{P_1}\tl{Q_1}$ and
$\tr{P_2}\tl{Q_2}$ have the same pullback span.
Call this span $\tl{A}\tr{B}$. Then both
$\tr{P_1}\tl{Q_1}$ and
$\tr{P_2}\tl{Q_2}$ are the pushout cospan of
$\tl{A}\tr{B}$, thus they must be isomorphic.
\end{proof}

\begin{lemma}
\label{lemma:phipsifull}
$\Phi \: \SpanMat \to\SVR$ and $\Psi \: \CospanMat\to\SVR$ are both full.	
\end{lemma}
\begin{proof}
Take any subspace
$S \: n\to m$ in $\SVR$. Picking any finite basis (say, of size $r$) for this subspace and multiplying out fractions gives us a finite set of elements in $\PID^{n+m}$.
In the obvious way, this yields
\[
n \tl{S_1} r \tr{S_2} m
\]	
in $\SpanMat$ with $\Phi(\tl{S_1}\tr{S_2})=S$. Thus $\Phi$ is full.
Let $\tr{R_1}\tl{R_2}$ be the cospan obtained from pushing out
$\tl{S_1}\tr{S_2}$ in $\Mat{\PID}$. By the conclusion of Lemma~\ref{lemma:arbitraryPROP}, $\Psi(\tr{R_1}\tl{R_2})=\Phi(\tl{S_1}\tr{S_2})=S$, which shows that $\Psi$ is full. \end{proof}

We now have all the ingredients to conclude our characterisation of $\SVR$.
\begin{theorem}\label{th:bottomfacePushout}
\eqref{eq:bottomface} is a pushout in $\PROP$.
\end{theorem}
\begin{proof}
Suppose that we have a commutative diagram of PROP morphisms as in~\eqref{eq:arbitrary}.
By the conclusions of Lemma~\ref{lemma:phipsifull} it suffices to show that there exists a PROP morphism $\Theta \: \SVR\to\mathbb{X}$ with $\Theta\Phi = \Gamma$ and $\Theta\Psi=\Delta$ -- uniqueness is automatic by fullness of $\Phi$ (or of $\Psi$).

Given a subspace $S \: n\to m$, by Lemma~\ref{lemma:phipsifull} there
exists a span $\tl{S_1}\tr{S_2}$ with
$\Phi(\tl{S_1}\tr{S_2})=S$.
We let $\Theta(S)=\Gamma(\tl{S_1}\tr{S_2})$.
This is well-defined: if $\tl{S_1'}\tr{S_2'}$ is another
span with $\Phi(\tl{S_1'}\tr{S_2'})=S$ then applying the conclusions of Lemma~\ref{lemma:pushoutphi} gives us that $\tl{S_1}\tr{S_2}$ and $\tl{S_1'}\tr{S_2'}$ have the same pushout in $\Mat{\PID}$. Now the conclusions of Lemma~\ref{lemma:arbitraryPROP} give us that $\Gamma(\tl{S_1}\tr{S_2})=\Gamma(\tl{S_1'}\tr{S_2'})$.
This  argument also shows that, generally, $\Theta\Phi=\Gamma$.
Finally, $\Theta$ preserves composition:
\begin{align*}
\Theta(R\poi S) &= \Theta(\Phi(\tl{R_1}\tr{R_2})\poi\Phi(\tl{S_1}\tr{S_2})) \\
&= \Theta(\Phi((\tl{R_1}\tr{R_2})\poi (\tl{S_1}\tr{S_2}))) \\
&= \Gamma((\tl{R_1}\tr{R_2})\poi (\tl{S_1}\tr{S_2})) \\
&= \Gamma(\tl{R_1}\tr{R_2})\poi\Gamma(\tl{S_1}\tr{S_2})\\
&= \Theta(R)\poi\Theta(S).
\end{align*}

It is also easy to show that
$\Theta\Psi = \Delta$: given a cospan $\tr{F}\tl{G}$ let
$\tl{P}\tr{Q}$ be its pullback span in $\Mat\PID$.
Using the conclusions of Lemma~\ref{lemma:arbitraryPROP},
${\Delta(\tr{F}\tl{G})} = \Gamma(\tl{P}\tr{Q})
=\Theta\Phi(\tl{P}\tr{Q})=\Theta\Psi(\tr{F}\tl{G})$.
\end{proof}

\subsection{The Cube: Rear Faces}\label{sec:cubeback}

To complete the proof of Theorem~\ref{th:IBR=SVR}, it remains to show that the rear faces of the cube~\eqref{eq:cube} commute.
\begin{equation}\label{backwardcube}
\tag{Rear}
\raise20pt
\hbox{$
\xymatrix@C=35pt{
{\IBRb} \ar[d]_{\sem{\IBRb}} & {\ABR + \ABRop} \ar[l]_{[\tau_1,\tau_2]} \ar[d]|{\sem{\ABR} + \sem{\ABR}^{\op} } \ar[r]^-{[\HAToIHw,\HAopToIHw]} & {\IBRw} \ar[d]^{\sem{\IBRw}} \\
\CospanMat & \VectR + \VectRop \ar[r]_{[\MatToSpan,\MatopToSpan]} \ar[l]^{[\MatToCospan,\MatopToCospan]} & \SpanMat
}$}
\end{equation}
For this purpose, it is useful to give an explicit description of the isomorphisms $\IBRw \to \SpanMat$ and $\IBRb \to \CospanMat$, whose existence has been shown in \S~\ref{sec:completeness}-\ref{sec:IBRbCospan}, in the same inductive way as $\sem{\ABR}$ is defined.

The two isomorphisms are indicated in~\eqref{backwardcube} with $\sem{\IBRb}$ and $\sem{\IBRw}$ respectively.
For the definition of $\MatToCospan$, $\MatopToCospan$, $\MatToSpan$ and $\MatopToCospan$ see the beginning of \S~\ref{sec:cubebottom}.
The PROP morphisms $\HAToIHw \: \ABR \to \IBRw$, $\HAopToIHw \: \ABRop \to \IBRw$ and $\tau_1 \: \ABR \to \IBRb$, $\tau_2 \: \ABRop \to \IBRb$ have been introduced after Definition~\ref{def:IBRw} and \ref{def:IBRb} respectively.

\subsubsection{An Inductive Presentation of $\sem{\IBRw}$}

Recall that, by definition, the set of generators of $\IBRw$ is the union of the set of generators of $\ABR$ and of $\ABRop$. This allows us to make the following inductive definition. 
\begin{definition} Let $\sem{\IBRw} \: \IBRw \to \SpanMat$ be the PROP morphism defined by the following mapping of the generators of $\IBRw$.
\begin{align*}
   c  \mapsto \left\{
	\begin{array}{ll}
        \MatToSpan(\sem{\AB}(c')) & \text{ if } c = \HAToIHw(c') \text{ and $c'$ is a generator of }\ABR \\
        \MatopToSpan( \sem{\AB}^{\op} (c'))  & \text{ if } c = \HAopToIHw(c') \text{ and $c'$ is a generator of }\ABRop. \\
        \end{array}
\right.
 \end{align*}
 \end{definition}
  The mapping is well-defined as all the equations of $\IBRw$ are sound w.r.t. $\sem{\IBRw}$. It is clear by definition that $\sem{\IBRw}$ makes the rightmost square in \eqref{backwardcube} commute. It remains to show the following result.

\begin{proposition}\label{prop:semanticsIBRwIso}  $\sem{\IBRw}$ is an isomorphism of PROPs.
\end{proposition}


\begin{proof}
For fullness, let $n \tl{A} z \tr{B} m$ be an arrow in $\SpanMat$. By fullness of $\sem{\ABR}$ there exist $c_1 \in \ABR[z,n]$ and $c_2 \in \ABR[z,m]$ such that $\sem{\ABR}(c_1) = A$ and $\sem{\ABR}(c_2) = B$. The following derivation shows that $n \tl{A} z \tr{B} m$ is targeted by $\HAopToIHw(\contrid{c_1}) \poi \HAToIHw(c_2) \in \IBR[n,m]$.
\begin{align*}
\sem{\IBRw}(\HAopToIHw(\contrid{c_1}) \poi \HAToIHw(c_2)) &= \sem{\IBRw}(\HAopToIHw(\contrid{c_1})) \poi \sem{\IBRw}(\HAToIHw(c_2)) \\
&=
\MatopToSpan(\sem{\ABR}^{op}(\contrid{c_1})) \poi \MatopToSpan(\sem{\ABR}(c_2)) \\
&=
\MatopToSpan(A \: n \to z) \poi \MatopToSpan(B \: z \to m) \\
&=
(n \tl{A}z \tr{\id} z) \poi (z \tl{\id}z \tr{B} m) \\
&=
n \tl{A}z \tr{B} m.
 \end{align*}
 It remains to show faithfulness. For this purpose, fix $c \in \IBRw[n,m]$ and $c' \in \IBRw[n,m]$ and suppose that $\sem{\IBRw}(c)  = \sem{\IBRw}(c')$. By Corollary~\ref{cor:factorisationIBRw} it follows that
\begin{eqnarray*}
\sem{\IBRw}(c) = n \tl{\sem{\ABR}(\contrid{c_1})} z \tr{\sem{\ABR}(c_2)} m \\
\sem{\IBRw}(c') = n \tl{\sem{\ABR}(\contrid{c_1'})} z' \tr{\sem{\ABR}(c_2')} m
\end{eqnarray*}
 where $c_1, c_1'$ are in $\ABRop$, $c_2,c_2'$ in $\ABR$ and $c = \HAopToIHw(c_1)\poi \HAToIHw(c_2)$, $c' = \HAopToIHw(c_1') \poi \HAToIHw(c_2')$. Since $\sem{\IBRw}(c)  = \sem{\IBRw}(c')$ are the same arrow of $\SpanMat$, then they are isomorphic spans: thus there is an invertible matrix $U \in \VectR[z',z]$ making the following diagram commute (and witnessing that $z = z'$).
\begin{eqnarray*}
\vcenter{
\xymatrix@C=25pt@R=15pt{ && \ar[drr]^{\sem{\ABR}(c_2)} z  \ar[dll]_{\sem{\ABR}(\contrid{c_1})}&& \\
 n && \ar[ll]^{\sem{\ABR}(\contrid{c_1'})} z' \ar[u]^>>>>U \ar[rr]_{\sem{\ABR}(c_2)}&& m }
 }
\end{eqnarray*}
Then by Lemma~\ref{lemma:mirror} we have that $c$ and $c'$ are equal as arrows of $\IBRw$.

 \end{proof}

%

\subsubsection{An Inductive Presentation of $\sem{\IBRb}$}
Similarly to what we did for $\IBRw$, we give an inductive presentation for the iso ${\IBRb \cong \CospanMat}$.
\begin{definition}Let $\sem{\IBRb} \: \IBRb \to \CospanMat$ be the PROP morphism defined on the generators of $\IBRb$ as follows:
\begin{align*}
   c  \mapsto \left\{
	\begin{array}{ll}
        \MatToCospan(\sem{\AB}(c')) & \text{ if } c = \tau_1(c') \text{ and $c'$ is a generator of }\ABR \\
        \MatopToCospan( \sem{\AB}^{\op} (c'))  & \text{ if } c = \tau_2(c') \text{ and $c'$ is a generator of }\ABRop. \\
        \end{array}
\right.
 \end{align*}
 \end{definition}
 The mapping is well-defined as all the equations of $\IBRb$ are sound w.r.t. $\sem{\IBRb}$. Also, $\sem{\IBRb}$ clearly makes the leftmost part of~\eqref{backwardcube} commute.

\begin{proposition}\label{prop:indPresIsoIBRb} $\sem{\IBRb}$ is an isomorphism of PROPs. \end{proposition}
\begin{proof} Following~\eqref{eq:IsoCospan}, it suffices to show that $\sem{\IBRb} = \pn \poi \sem{\IBRw} \poi \tra$ --- see~\S\ref{sec:IBRbCospan} for the definition of $\pn \: \IBRb \to \IBRw$ and $\tra \: \SpanMat \to \CospanMat$. The equality can be easily verified by induction on diagrams of $\IBRb$: for instance, $\sem{\IBRb}$ maps $\Bmult \: 2 \to 1$ to $2 \tr{\id} 2 \tl{\tiny \matrixOneOne} 1$ and $\pn \poi \sem{\IBRw} \poi \tra$ maps $\Bmult$ first to $\Wmult$, then to $2 \tl{\id} 2 \tr{\tiny \matrixOneOneFlat} 1$ and finally to $2 \tr{\id} 2 \tl{\tiny \matrixOneOne} 1$.  \end{proof}

%
%

\subsection{The Cube Rebuilt}\label{sec:cuberebuilt}

The results of the previous two sections conclude the proof of Theorem~\ref{th:IBR=SVR}. We are now in position to patch together all the faces of the cube \eqref{eq:cube}. This will also give us an inductive presentation of the isomorphism $\sem{\IBR} \: \IBR \to \SVR$.
\begin{equation}
\label{eq:cuberebuilt}
\vcenter{
    \xymatrix@=35pt{
    & {\ABR + \ABRop} \ar[dr]^{[\varphi_1,\varphi_2]} \ar[dd]_(.7){\sem{\ABR}+\sem{\ABR}^{\op}}|{\hole}
    \ar[dl]_{[\tau_1,\tau_2]} \ar[rr]^{[\HAToIHw,\HAopToIHw]} & & {\IBRw} \ar[dl]_{\Theta} \ar[dd]^{\sem{\IBRb}} \\
    {\IBRb} \ar[rr]^(.7){\Lambda} \ar[dd]_{\sem{\IBRb}}  & & {\IBR} \ar@{.>}[dd]^(.3){\sem{\IBR}} \\
    & {\VectR+ \VectRop}  \ar[dr]^{[\psi_1,\psi_2]}  \ar[dl]_{[\MatToCospan,\MatopToCospan]} \ar[rr]^(.4){[\MatToSpan,\MatopToSpan]}|(.51){\hole} & & \SpanMat \ar[dl]_{\Phi} \\
    {\CospanMat} \ar[rr]^{\Psi} & & {\SVR}
    }
}
\end{equation}
Above we draw the PROP morphism $[\psi_1,\psi_2] \: \VectR + \VectRop \to \SVR$ defined by commutativity of the bottom face. Commutativity of all the faces yields commutativity of the ``section'':
 \begin{equation}\label{eq:cubesection}
 \tag{Sec}
 \vcenter{
 \xymatrix@=15pt{
 & {\ABR + \ABRop} \ar[dr]^{[\varphi_1,\varphi_2]} \ar[dd]_{\sem{\ABR}+\sem{\ABR}^{\op}}
& &  \\
& & {\IBR} \ar@{.>}[dd]^{\sem{\IBR}} \\
& {\VectR+ \VectRop}  \ar[dr]_{[\psi_1,\psi_2]} & & \\
 & & {\SVR}
 }
 }
 \end{equation}

 \subsubsection{An Inductive Presentation of $\sem{\IBR}$}

Diagram~\eqref{eq:cubesection} provides us a recipe for an inductive presentation of $\sem{\IBR}$, similarly to what we previously did for $\sem{\IBRb}$ and $\sem{\IBRw}$:
\begin{align*}
   c  \mapsto \left\{
	\begin{array}{ll}
        \psi_1(\sem{\ABR}(c')) & \text{ if } c = \varphi_1(c') \text{ and $c'$ is a generator of }\ABR \\
        \psi_2( \sem{\ABR}^{\op} (c'))  & \text{ if } c = \varphi_2(c') \text{ and $c'$ is a generator of }\ABRop. \\
        \end{array}
\right.
 \end{align*}

By using commutativity of~\eqref{eq:cuberebuilt}, one can actually give a more direct description of the behaviour of $\sem{\IBR}$.
 In the definition below, the notation $[(\xx_1,\yy_1), \dots, (\xx_z,\yy_z)]$ for an arrow in $\SVR[n,m]$ indicates the subspace of $\frPID^n \times \frPID^m$ spanned by pairs $(\xx_1,\yy_1), \dots, (\xx_z,\yy_z)$ of vectors, where each $\xx_i$ is in $\PID^n$ and each $\yy_i$ is in $\PID^m$. Also, $\matrixNull$ denotes the unique element of the space of dimension $0$.

\begin{definition}\label{def:semIBRInd} The following is an assignment to each generator $c \: n \to m$ of $\IBR$ of a subspace in $\SVR[n,m]$. For the generators which are also in the signature of $\ABR$:
\begin{multicols}{2}\noindent
\begin{equation*}
\Bcomult \longmapsto [(%
\small{\matrixOne},\tiny{\left(\begin{array}{c}
                \!\!  1 \!\!\\
                 \!\! 1 \!\!
                \end{array}\right)})]
\end{equation*}
\begin{equation*}
\Wmult  \longmapsto  [\tiny{(\left(%
                \begin{array}{c}
                \!\!\!  0 \!\!\!\\
                \!\!\!  1 \!\!\!
                \end{array}\right)},\small{\matrixOne}),(\tiny{\left(%
                \begin{array}{c}
                \!\!\!  1 \!\!\!\\
                \!\!\!  0 \!\!\!
                \end{array}\right)},\small{\matrixOne})]
\end{equation*}
\end{multicols}
\smallskip
\begin{multicols}{3}\noindent
\begin{equation*}
\Bcounit \longmapsto [(\small{\matrixOne,\matrixNull})]
\end{equation*}
\begin{equation*}
\Wunit \longmapsto \{(\small{\matrixNull,\matrixZero})\}
\end{equation*}
\begin{equation*}
\scalar \longmapsto [(\small{\matrixOne},\small{\matrixK})]
\end{equation*}
\end{multicols}
The remaining generators of $\IBR$ are those coming from the signature of $\ABRop$.
\begin{multicols}{2}\noindent
\begin{equation*}
\Bmult \longmapsto [((%
\tiny{\left(\begin{array}{c}
                \!\!  1 \!\!\\
                 \!\! 1 \!\!
                \end{array}\right)}),\small{\matrixOne})]
\end{equation*}
\begin{equation*}
\Wcomult  \longmapsto  [(\small{\matrixOne},\tiny{(\left(%
                \begin{array}{c}
                \!\!\!  0 \!\!\!\\
                \!\!\!  1 \!\!\!
                \end{array}\right))}),(\small{\matrixOne},\tiny{\left(%
                \begin{array}{c}
                \!\!\!  1 \!\!\!\\
                \!\!\!  0 \!\!\!
                \end{array}\right)})]
\end{equation*}
\end{multicols}
\smallskip
\begin{multicols}{3}\noindent
\begin{equation*}
\Bunit \longmapsto [(\small{\matrixNull,\matrixOne})]
\end{equation*}
\begin{equation*}
\Wcounit \longmapsto \{(\small{\matrixZero,\matrixNull})\}
\end{equation*}
\begin{equation*}
\coscalar \longmapsto [(\small{\matrixK},\small{\matrixOne})].
\end{equation*}
\end{multicols}
\end{definition}
Note the symmetry between generators $c$ of $\ABR$ and generators $\coc{c}$ of $\ABRop$: in Definition~\ref{def:semIBRInd}, $\coc{c}$ is being mapped to the inverse of the subspace assigned to $c$. We now verify that this inductive assignment provides an equivalent definition of the PROP morphism $\sem{\IBR}$.

\begin{proposition}\label{prop:SemIBRInductive} $\sem{\IBR}\: \IBR \to \SVR$ in \eqref{eq:cubesection} is equivalently described as the PROP morphism inductively defined on the generators of $\IBR$ according to Definition~\ref{def:semIBRInd}.
\end{proposition}
\begin{proof}
One can verify that the PROP morphism $\IBR \to \SVR$ determined by Definition~\ref{def:semIBRInd} makes the front faces of~\eqref{eq:cuberebuilt} commute. As $\sem{\IBR}$ is by construction the unique PROP morphism with this property, the two must coincide.
\end{proof}

\subsection{A Return on Graphical Linear Algebra}\label{sec:graphicallinearlagebra}

In \S~\ref{sec:overview} we claimed that $\IBR$ offers a different perspective on linear algebra, where computations with matrices, subspaces and other encodings are replaced by purely diagrammatic equational reasoning. We already gave a glimpse of this approach in the cube construction~\eqref{eq:cube}, when we used a graphical rendition of Gaussian elimination (proof of Lemma~\ref{lemma:BHNFequalKernel}) and of the transpose operation (\S~\ref{sec:IBRbCospan}).

Now that the theory $\IBR$ and its subspace interpretation $\sem{\IBR} \: \IBR \to \SVR$ have been formally introduced, we have all the ingredients for more examples. Our illustration should also give to the reader a more concrete grip on the inductive definition of $\sem{\IBR}$.

\begin{example}[Matrices and kernels] \label{ex:graphicalLA-matrices} 
We begin with some simple observations concerning the diagrammatic representation of matrices and their kernel. For matrices, we already described (Def.~\ref{def:matrixform}) a canonical representation as string diagrams of $\ABR$. As expected, this correspondence lifts to $\IBR$: following Definition~\ref{def:semIBRInd}, one can verify that, for any $m \times n$ matrix $A$, $\sem{\IBR} \: \IBR \to \SVR$ maps the corresponding string diagram $\circuitA \in \IBR[n,m]$ in matrix form to the subspace generated by $A$:
\[ \sem{\IBR} (\circuitA) = \{(\xx, A \xx) \mid \xx \in \PID^n\} \in \SV[n,m]. \]
In $\ABRop$, $A$ is represented by the $m\to n$ string diagram $\Astar$. The PROP morphism $\sem{\IBR}$ behaves on $\Astar \in \IBR[m,n]$ as on generators of $\ABRop$, mapping it to the inverse of the relation $\sem{\IBR}(\circuitA)$:
\[ \sem{\IBR} (\Astar) = \{(A \xx, \xx) \mid \xx \in \PID^n\} \in \SV[m,n]. \]
Combining diagrams from $\ABR$ and $\ABRop$, it is possible to describe in $\IBR$ the \emph{kernel (null) space} of $A$: one simply needs to post-compose the matrix form for $A$ with white counits, thus obtaining $\kernelA \in \IBR[n,0]$. This can be verified by computing the semantics of $\kernelA$:
\begin{eqnarray*}
 \sem{\IBR}(\kernelA)& = & \sem{\IBR}(\circuitA) \poi \sem{\IBR}(\lower1pt\hbox{$\includegraphics[height=.4cm]{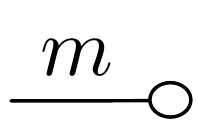}$}) \\
 & = & \{(\xx, A \xx) \mid \xx \in \PID^n\} \poi \{ (\zerov, \matrixNull) \} \\
 & = & \{(\xx , \matrixNull) \mid A \xx = \zerov\}
\end{eqnarray*}
where $\zerov$ here indicates the $0$-vector of length $m$.
\end{example}

\begin{example}[Subspaces] \label{ex:subspacesgraphical} By Theorem~\ref{th:IBR=SVR}, the theory $\IBR$ associates to any subspace $S \in \SVR[n,m]$ a class of equivalent string diagrammatic representation. In particular, our modular construction disclosed two factorisation properties of $\IBR$ (Theorem~\ref{Th:factIBR}), telling that among these equivalent representations there is one which is a span of string diagrams of $\ABR$ and another one which is a cospan of string diagrams of $\ABR$. We now show that the two canonical forms correspond to well-known encodings of subspaces.

First, let us focus on the cospan form: it witnesses the representation of a subspace as an homogeneous system of linear equations. We illustrate this with a simple example:
\begin{equation}\label{eq:systemequations}
\begin{aligned}
k_1 x_1 + k_2 x_2 + k_3 x_3 & = & 0 \  \\
l_1 x_1 + l_2 x_2 + l_3 x_3 & = & 0.
\end{aligned}
\end{equation}
The above system of equations yields a subspace $S$ of $\PID^3$. There is a simple recipe to compute a string diagram $c_1 \in \IBR[3,0]$ in cospan form representing $S$, starting from~\eqref{eq:systemequations}.
\begin{equation*}
   c_1 \quad \df \quad \lower45pt\hbox{$\includegraphics[height=4.5cm]{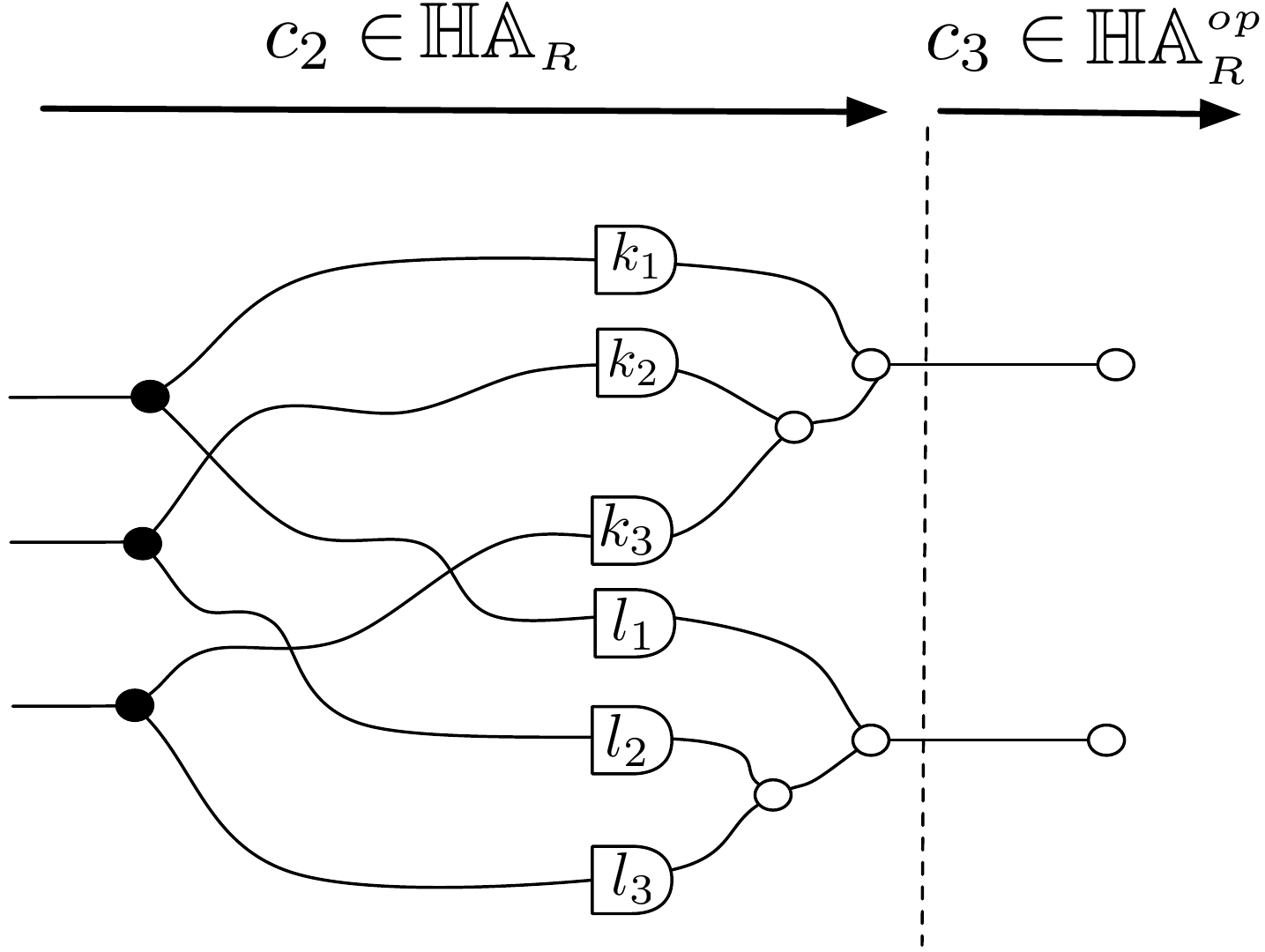}$}.
 \end{equation*}
The dotted line separates the components describing the left and right side of the equations. Ports represent variables $x_1, x_2, x_3$. The branching given by $\Bcomult$ connects them to components $\Wmult$ which are responsible of modeling sum in~\eqref{eq:systemequations}. Composition with $\Wcounit$ means that each sum must be equal to $0$. This description is confirmed when computing $\sem{\IBR}(c_1)$ according to Definition~\ref{def:semIBRInd}:
\begin{eqnarray*}
 \sem{\IBR}(c_1) &=& \sem{\IBR}(c_2)\poi\sem{\IBR}(c_3) \\
&=& \{(\vv,\ww) \mid
{\tiny\begin{pmatrix}
 k_1 & k_2 & k_3 \\
 l_1 & l_2 & l_3
\end{pmatrix}} \vv = \ww \} \poi \{ ({\tiny\begin{pmatrix} 0 \\ 0 \end{pmatrix}}, \matrixNull)\} \\
&=& \{(\vv,\matrixNull) \mid
{\tiny\begin{pmatrix}
 k_1 & k_2 & k_3 \\
 l_1 & l_2 & l_3
\end{pmatrix}} \vv = {\tiny\begin{pmatrix} 0 \\ 0 \end{pmatrix}} \} \\
& = & S\in \SVR[3,0].
\end{eqnarray*}
Drawning on the observations of Example~\ref{ex:graphicalLA-matrices}, one can also describe $c_1$ as the diagrammatic representation of the kernel space of the matrix encoded by $c_2$.

It is worth remarking that $c_1$ is not the only way of representing $S$ as a string diagram in cospan form. By definition of $\SVR$, one can model $S \subseteq \PID^3$ as an arrow of any type $n \to m$ in $\SVR$ such that $n+m = 3$. Different choices of $n$ and $m$ intuitively correspond to move variables in \eqref{eq:systemequations} on the other side of each equation. For instance, we choose $n=2$ and $m=1$ and let $x_1$ be the variable moved on the right side:
\begin{equation}\label{eq:systemequationstwo}
\begin{aligned}
k_2 x_2 + k_3 x_3 & = & -k_1 x_1 \\
l_2 x_2 + l_3 x_3 & = & -l_1 x_1.
\end{aligned}
\end{equation}

The transformation of~\eqref{eq:systemequations} into \eqref{eq:systemequationstwo} has a neat string diagrammatic description: it corresponds to ``rewiring'' the first port on the left in $c_1$ to a port on the other boundary, by pre-composing with $\rccB$. We can then use equational reasoning in $\IBR$ to obtain a string diagram $c_4 \in \IBR[2,1]$ in cospan form modeling~\eqref{eq:systemequationstwo}:
\begin{equation*}
   \includegraphics[height=9cm,width=15cm]{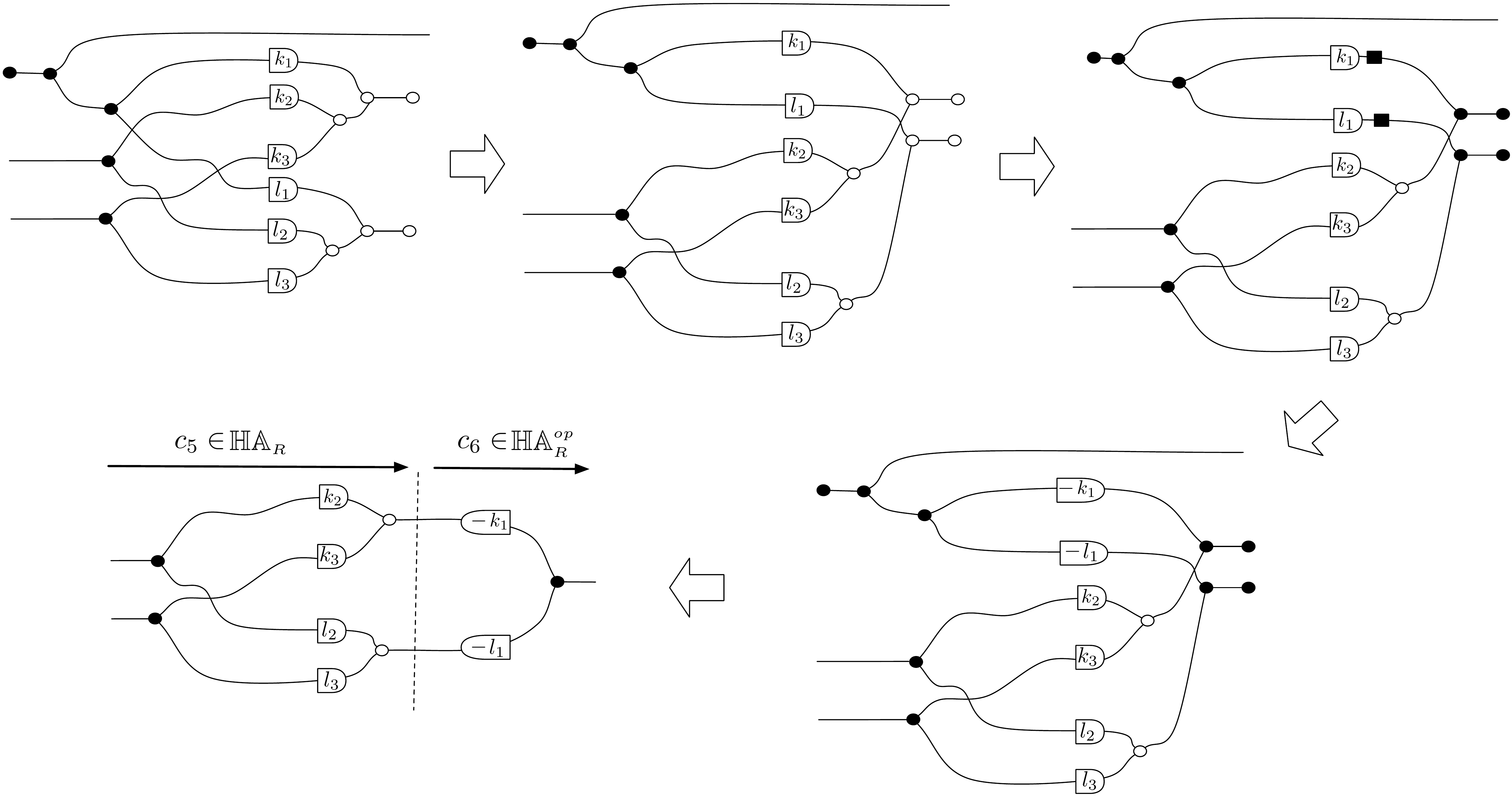}.
 \end{equation*}
In the above derivation, the first step pulls out the part of the diagram of interest using~\eqref{eq:sliding1}-\eqref{eq:sliding2}; then, \eqref{eq:rcc}, \eqref{eq:scalarmult} and the conclusion of Proposition~\ref{prop:star=refl} are applied in the next steps. Computing with $\sem{\IBR}$ confirms that $c_4$ represents the subspace $S$:
\begin{eqnarray*}
 \sem{\IBR}(c_4) &=& \sem{\IBR}(c_5)\poi\sem{\IBR}(c_6) \\
&=& \{(\yy,\zz) \mid
{\tiny\begin{pmatrix}
 k_2 \!& k_3 \\
 l_2 \!& l_3
\end{pmatrix}} \yy = \zz \} \poi \{ (\zz,\xx) \mid\zz = {\tiny\begin{pmatrix}
 -k_1 \\  -l_1
\end{pmatrix}}\xx\} \\
&=& \{(\yy,\xx) \mid
{\tiny\begin{pmatrix}
 k_2 \!\! \!& k_3 \\
 l_2 \!\! \!& l_3
\end{pmatrix}} \yy = {\tiny\begin{pmatrix}
 -k_1 \\  -l_1
\end{pmatrix}} \xx \} \\
&=& \{(\yy,\xx) \mid
{\tiny\begin{pmatrix}
 k_1 \!\! \!& k_2 \!\! \!& k_3 \\
 l_1 \!\! \!& l_2 \!\! \!& l_3
\end{pmatrix}} {\tiny\begin{pmatrix}
 \xx \\
 \yy
\end{pmatrix}} =  {\tiny\begin{pmatrix} 0 \\ 0 \end{pmatrix}}\} \\
& = & S\in \SVR[2,1].
\end{eqnarray*}
We now move to analysing the span form for string diagrams of $\IBR$. While cospans give a canonical representation of systems of equations, with spans we can encode \emph{bases} for a subspace. Suppose $S' \subseteq \PID^4$ is the subspace spanned by vectors
\[ \vv = {\tiny\left(%
				\begin{array}{c}
				  \!\!\!k_1\!\!\! \\
				  \!\!\!k_2\!\!\! \\
				  \!\!\!k_3\!\!\!\\
				  \!\!\!k_4\!\!\!
				\end{array}\right)} \qquad \qquad
   \ww = {\tiny\left(%
				\begin{array}{c}
				  \!\!\!l_1\!\!\! \\
				  \!\!\!l_2\!\!\! \\
				  \!\!\!l_3\!\!\!\\
				  \!\!\!l_4\!\!\!
				\end{array}\right)}.
\]
On the base of this information, we can construct a string diagram $c_7 \in \IBR[2,2]$ in span form:
\begin{equation*}
   c_7 \quad \df \quad \lower35pt\hbox{$\includegraphics[height=3.5cm]{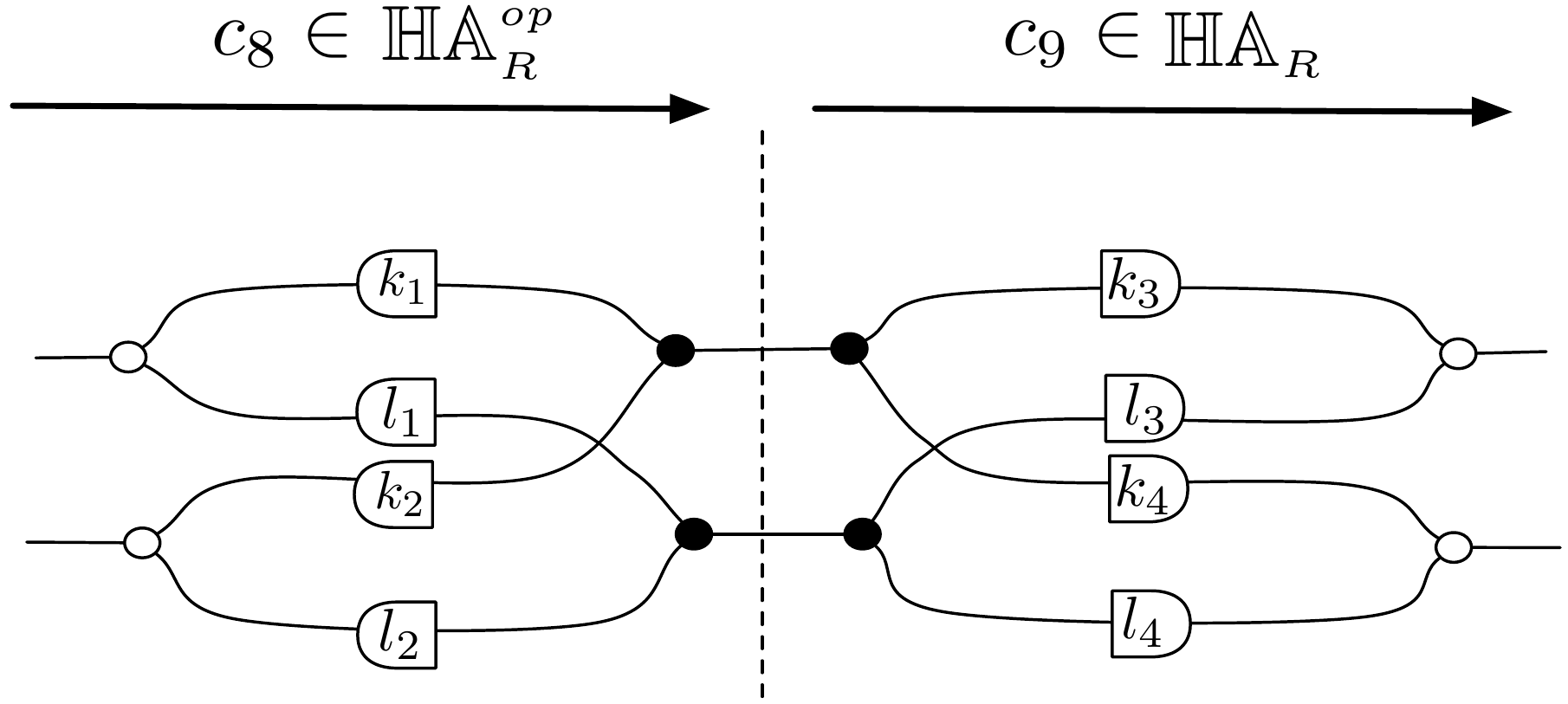}$}.
 \end{equation*}
The intuition is that the two ports on middle boundary of $c_7$ represent $\vv$ and $\ww$ respectively. On the outer boundaries, there are as many ports as the dimension of $\PID^4$. To see that $c_5$ models $S'$, note that \begin{eqnarray*}
\sem{\IBR}(c_8)
=
\{ (\xx,\matrixTwoROneC{r_1}{r_2}) \mid \xx =
{\tiny\begin{pmatrix}  k_1\!\!\! & \!\!\! l_1\!\!\! \\ \!\!\!k_2\!\!\! & \!\!\! l_2 \end{pmatrix}} \matrixTwoROneC{r_1}{r_2}
\}
=
\{ (\xx,\matrixTwoROneC{r_1}{r_2}) \mid \xx = r_1 \matrixTwoROneC{k_1}{k_2} + r_2 \matrixTwoROneC{l_1}{l_2}
\}
\\
 \sem{\IBR}(c_9) = \{ (\matrixTwoROneC{r_1}{r_2},\zz) \mid \zz =
{\tiny\begin{pmatrix} k_3\!\!\! & \!\!\! l_3\!\!\! \\ \!\!\!k_4\!\!\! & \!\!\! l_4 \end{pmatrix}} \matrixTwoROneC{r_1}{r_2}
\}
=
\{ (\matrixTwoROneC{r_1}{r_2},\zz) \mid \zz = r_1 \matrixTwoROneC{k_3}{k_4} + r_2 \matrixTwoROneC{l_3}{l_4}
\}.
\end{eqnarray*}
Therefore,
$$\sem{\IBR}(c_7) = \sem{\IBR}(c_8)\poi\sem{\IBR}(c_9) = \{ (\xx,\zz) \mid \text{ there is }\matrixTwoROneC{r_1}{r_2} \text{ such that }\matrixTwoROneC{\xx}{\zz} = r_1 \vv + r_2 \vv\} = S'.$$
Albeit our examples only showed how to pass from a system of equations/a basis to a string diagram in cospan/span form, the converse is also possible: the computation is made easy by the fact that the sub-diagrams in $\ABR$ and $\ABRop$ of any string diagram in cospan/span form can be always put in matrix form, by virtue of Lemma~\ref{lem:path1}.

These observations allow us to see the factorisation result of Theorem~\ref{Th:factIBR} under a new light: it amounts to the familiar result that any subspace (any string diagram of $\IBR$) is equivalently presented as (is equal in $\IBR$ to) a system of linear equations (a cospan in $\ABR$) and as a basis (a span in $\ABR$). The advantage of our description is that all these nonhomogeneous encodings of the same entity are now uniformly described by \emph{a unique} string diagrammatic syntax.
\end{example}

In concluding our excursus on graphical linear algebra, we shall prove a simple fact about kernels of matrices using equational reasoning in $\IBRw$.
\begin{proposition}\label{prop:graphicalinjectivematrix} A matrix is injective if and only if its kernel is the empty space. \end{proposition}
For the proof of Proposition~\ref{prop:graphicalinjectivematrix}, it is useful first to formulate its statement in purely diagrammatic terms.

\begin{lemma}\label{lemmainjmatrices1}A matrix is injective if and only if $\IBR$ proves
\begin{equation}\label{eq:lemmainjmatrices1}
\lower5pt\hbox{$\includegraphics[height=.6cm]{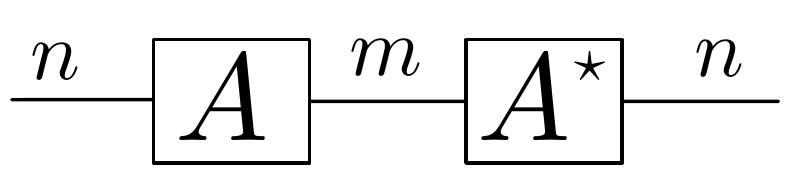}$}
=
\lower0pt\hbox{$\includegraphics[height=.4cm]{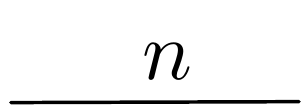}$}.
\end{equation}
\end{lemma}
\begin{proof}
For the left to right direction, suppose $A$ is an injective matrix. Thus it is a monomorphism in $\VectR$, meaning that
\begin{equation}\label{eq:monopullbackid}
\vcenter{
    \xymatrix@=8pt{
    & \ar[dl]_{\id}  \ar[dr]^{\id} & \\
    \ar[dr]_{A} && \ar[dl]^{A} \\
    &&
    }
}
\end{equation}
is a pullback square in $\VectR$. As observed in Remark~\ref{rmk:readoffequationsMatrix}, this implies that the equation \eqref{eq:lemmainjmatrices1} arises from the distribuitve law $\lambdapb \: \VectR \bicomp{\PJ} \VectRop \to \VectRop \bicomp{\PJ} \VectR$, and thus it is part of the equational theory of $\IBRw \cong \VectRop \bicomp{\PJ} \VectR$. Since $\IBR$ is a quotient of $\IBRw$, it also proves \eqref{eq:lemmainjmatrices1}.

For the converse direction, suppose that \eqref{eq:lemmainjmatrices1} holds true in $\IBR$. Then we have that
\begin{eqnarray*}
\{(\xx , A \xx) \mid \xx \in \PID^n \} \poi \{(A \xx ,\xx) \mid \xx \in \PID^n\} & = &  \sem{\IBR}(\lower5pt\hbox{$\includegraphics[height=.6cm]{graffles/AAstar.pdf}$}) \\
& = & \sem{\IBR}(\lower0pt\hbox{$\includegraphics[height=.4cm]{graffles/idnNoFrame.pdf}$}) \\
& = & \{ (\xx,\xx) \mid \xx \in \PID^n \}.
\end{eqnarray*}
Suppose now $\vv$ and $\ww$ are vectors such that $A \vv = A \ww$. We have $(\vv, A \vv) \in \{(\xx , A \xx) \mid \xx \in \PID^n \}$ and $(A \ww, \ww) \in \{(A \xx ,\xx) \mid \xx \in \PID^n\}$. Thus, using that $A \vv = A \ww$, we have $(\vv, \ww) \in \{(\xx , A \xx) \mid \xx \in \PID^n \} \poi \{(A \xx ,\xx) \mid \xx \in \PID^n\}$. By the above derivation, this implies $(\vv , \ww) \in \{ (\xx,\xx) \mid \xx \in \PID^n \}$ and thus $\vv = \ww$. This proves injectivity of $A$.
\end{proof}

\begin{lemma}\label{lemmainjmatrices2} The kernel of a matrix $A$ is the empty space if and only if $\IBR$ proves
\begin{equation}\label{eq:lemmainjmatrices2}
\lower6pt\hbox{$\includegraphics[height=.7cm]{graffles/kernelA.pdf}$}
=
\lower2pt\hbox{$\includegraphics[height=.6cm]{graffles/emptyspace.pdf}$}.
\end{equation}
\end{lemma}
\begin{proof} The fact that the kernel of $A$ is equal to the empty space can be expressed in $\SVR$ by saying
\begin{equation}\label{eq:lemmainjmatrices2mezzo}
 \{ \big(\xx, \matrixNull \big) \mid A\xx = \zerov\} \quad = \quad \{ \big( \zerov, \matrixNull \big) \} \qquad \in \SVR[n,0]
 \end{equation}
 where $\zerov$ stands for a $0$-vector of the appropriate length --- $m$ in the rightmost and $n$ in the leftmost occurrence.

Following Example~\ref{ex:graphicalLA-matrices}, $\sem{\IBR} \: \IBR \to \SVR$ gives us a recipe for representing these two subspaces as string diagrams:
\[ \{ \big(\xx, \matrixNull \big) \mid A\xx = \zerov\} = \sem{\IBR}(\kernelA ) \qquad \qquad \{ \big( \zerov, \matrixNull \big) \} = \sem{\IBR}(\emptyspace) \]
Because $\sem{\IBR}$ is faithful, \eqref{eq:lemmainjmatrices2mezzo} holds in $\SVR$ if and only if \eqref{eq:lemmainjmatrices2} is provable in $\IBR$.
 \end{proof}

We will also make use of the following simple observations.

\begin{lemma}\label{lemmainjmatrices3} For any string diagram $c \in \ABR[n,m]$, $\IBR$ proves
\begin{multicols}{3}\noindent
\begin{equation*}
\lower10pt\hbox{$
\lower6pt\hbox{$\includegraphics[height=.7cm]{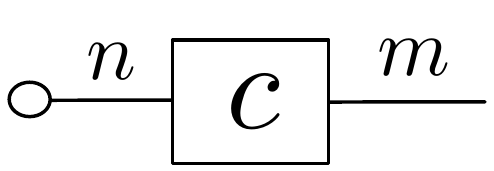}$}
=
\lower3pt\hbox{$\includegraphics[height=.6cm]{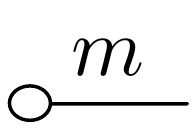}$}
$}
\end{equation*}
\begin{equation*}
\lower8pt\hbox{$\includegraphics[height=1cm]{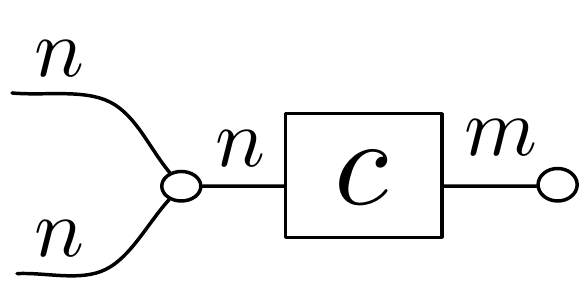}$}
=
\lower10pt\hbox{$\includegraphics[height=1cm]{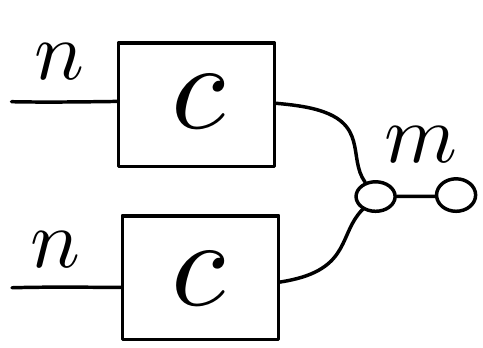}$}
\end{equation*}
\begin{equation*}
\lower11pt\hbox{$
\lower6pt\hbox{$\includegraphics[height=.6cm]{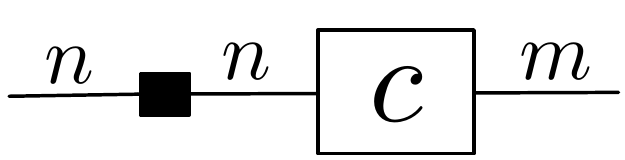}$}=
\lower6pt\hbox{$\includegraphics[height=.6cm]{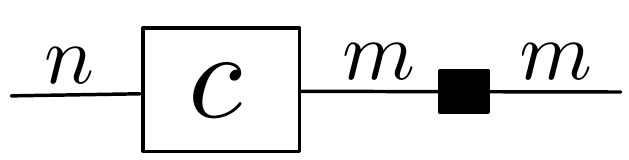}$}
$}
\end{equation*}
\end{multicols}
\end{lemma}
\begin{proof} 
%
The three statements are proved by a simple induction on $c$. \end{proof}

We now have all the ingredients to give a diagrammatic proof of Proposition~\ref{prop:graphicalinjectivematrix} within $\IBR$.

\begin{proof}[Proof of Proposition~\ref{prop:graphicalinjectivematrix}] Let us fix a matrix $A \in \VectR[n,m]$.
 By Lemma~\ref{lemmainjmatrices1} and \ref{lemmainjmatrices2}, the statement of the proposition reduces to proving that
 \begin{eqnarray}\label{eq:matrixinjkernel_reformulation}
\lower5pt\hbox{$\includegraphics[height=.6cm]{graffles/AAstar.pdf}$}
=
\lower0pt\hbox{$\includegraphics[height=.4cm]{graffles/idnNoFrame.pdf}$}
 & \Longleftrightarrow &
\lower5pt\hbox{$\includegraphics[height=.6cm]{graffles/kernelA.pdf}$}
=
\lower0pt\hbox{$\includegraphics[height=.4cm]{graffles/emptyspace.pdf}$}.
 \end{eqnarray}
  For the left to right direction,
 \begin{eqnarray*}
\lower5pt\hbox{$\includegraphics[height=.6cm]{graffles/AAstar.pdf}$}
=
\lower0pt\hbox{$\includegraphics[height=.4cm]{graffles/idnNoFrame.pdf}$}
 & \Longrightarrow &
\lower5pt\hbox{$\includegraphics[height=.6cm]{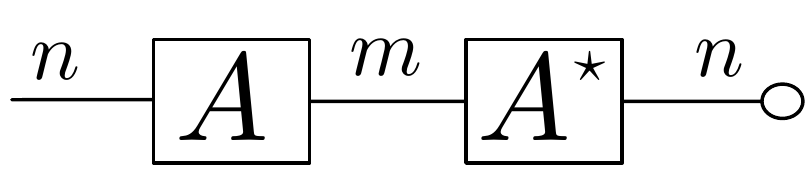}$}
=
\lower0pt\hbox{$\includegraphics[height=.4cm]{graffles/emptyspace.pdf}$} \\
 & \Longrightarrow &
\lower5pt\hbox{$\includegraphics[height=.6cm]{graffles/kernelA.pdf}$}
=
\lower0pt\hbox{$\includegraphics[height=.4cm]{graffles/emptyspace.pdf}$}
 \end{eqnarray*}
 where, for the last implication, we apply the dual of the first equation in Lemma~\ref{lemmainjmatrices3} to cancel $\Astar$. For the converse direction:
 \begin{eqnarray*}
\lower5pt\hbox{$\includegraphics[height=.6cm]{graffles/AAstar.pdf}$}
&\eql{Prop.\ref{prop:star=refl}}&
\lower15pt\hbox{$\includegraphics[height=1.6cm]{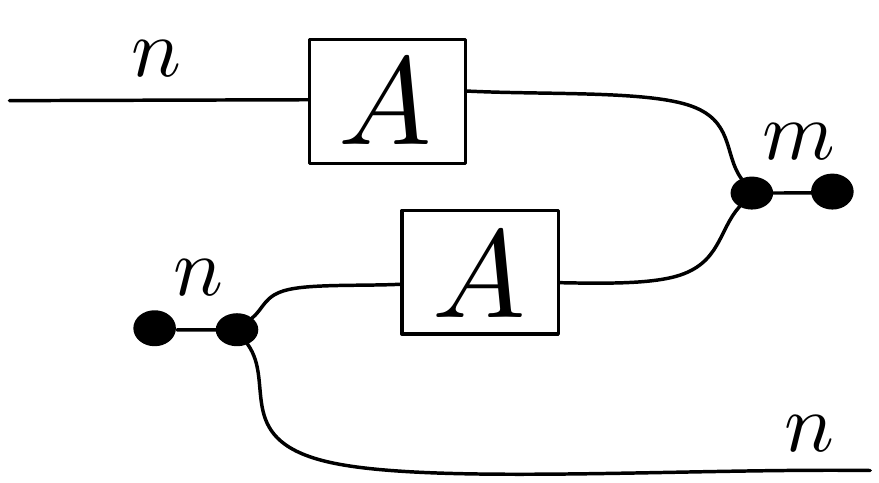}$} \\
 & \eql{\eqref{eq:lccb}} &
\lower15pt\hbox{$\includegraphics[height=1.6cm]{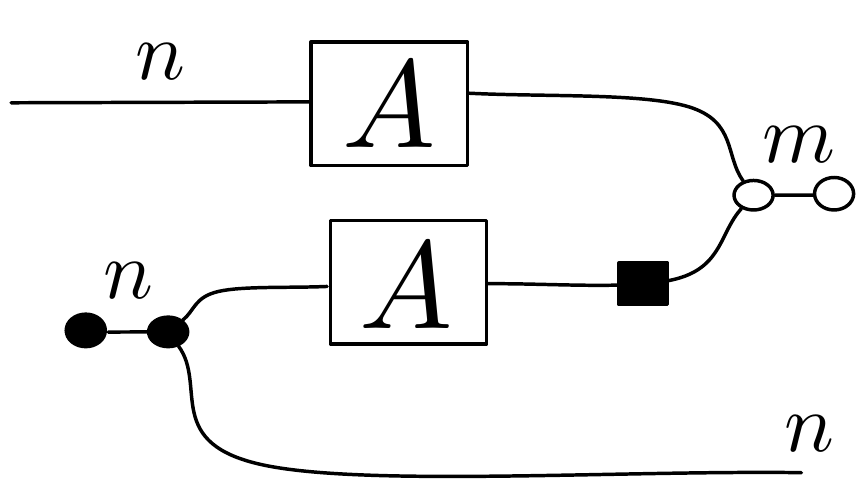}$} \\
 & \eql{Lemma~\ref{lemmainjmatrices3}} &
\lower15pt\hbox{$\includegraphics[height=1.6cm]{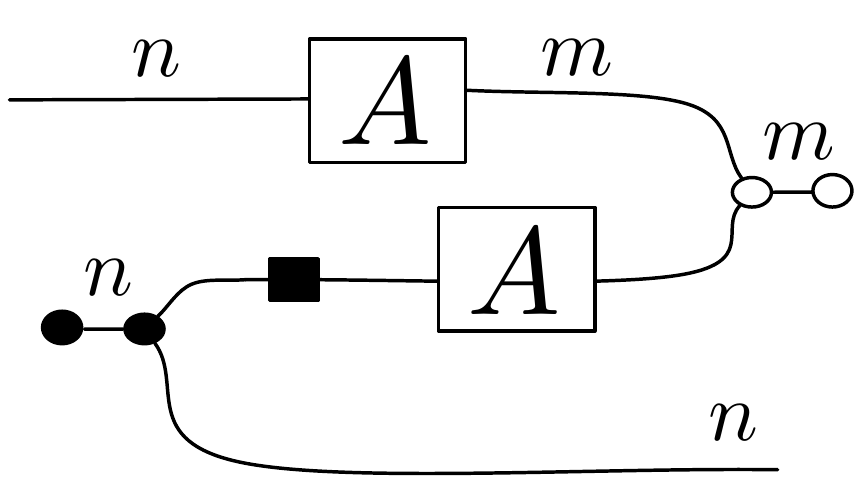}$} \\
 & \eql{Lemma~\ref{lemmainjmatrices3}} &
\lower15pt\hbox{$\includegraphics[height=1.6cm]{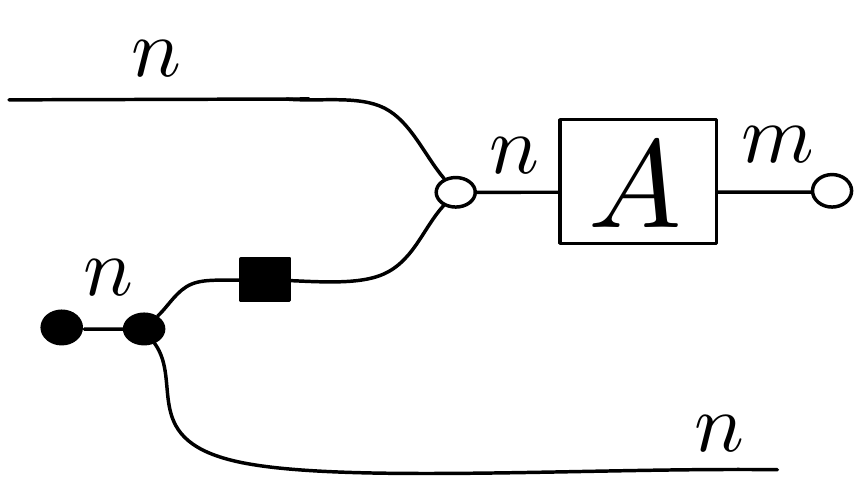}$} \\
 & \eql{Assumption} &
\lower15pt\hbox{$\includegraphics[height=1.6cm]{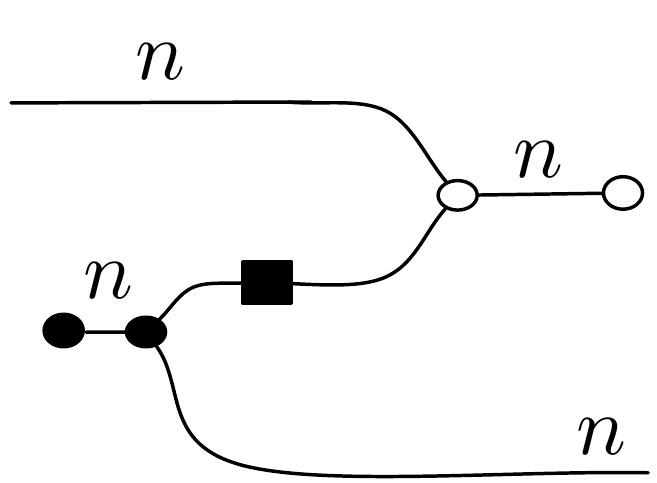}$} \\
 & \eql{\eqref{eq:lccb}} &
\lower13pt\hbox{$\includegraphics[height=1.6cm]{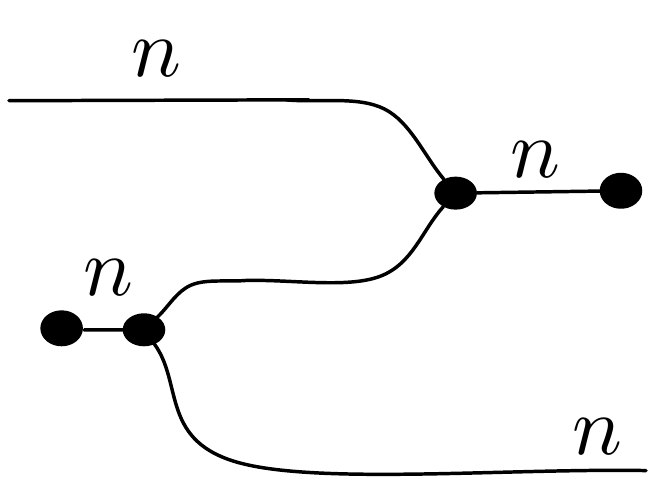}$} \\
 & \eql{\eqref{eq:gensnake}} &
\raise2pt\hbox{$\includegraphics[height=.4cm]{graffles/idnNoFrame.pdf}$}.
 \end{eqnarray*}
\end{proof}

\section{Example: Interacting Hopf Algebras for Rational Subspaces}\label{sec:instances}

In concluding this chapter, we exhibit a simple, yet important, example of our cube construction~\eqref{eq:cube}: the axiomatisation $\IH{\Z}$ for the PROP of rational subspaces. A more thorough example, concerning polynomials and formal power series, will be the subject of the next chapter. We begin by describing the sub-theory of $\IH{\Z}$ modeling integer matrices.

\paragraph{The theory of integer matrices} By Proposition~\ref{prop:ab=vect}, the PROP $\Mat{\Z}$ of integer matrices is presented by the axioms \eqref{eq:wmonassoc}-\eqref{eq:scalarsum} of $\HA{\Z}$. In fact, a finite axiomatisation is possible: let $\HA{}$ be the PROP freely generated by the SMT with signature $\{\antipode ,\Bcounit , \Bcomult , \Wunit , \Wmult \}$ and equations:
 \begin{multicols}{3}\noindent
\begin{equation*}
\lower9pt\hbox{$\includegraphics[height=.9cm]{graffles/Wunitlaw.pdf}$}
\!\!\!
=\!
\lower5pt\hbox{$\includegraphics[height=.6cm]{graffles/idcircuit.pdf}$}
\end{equation*}
\begin{equation*}
\lower5pt\hbox{$\includegraphics[height=.6cm]{graffles/Wmult.pdf}$}
\!
=
\!\!\!\!
\lower11pt\hbox{$\includegraphics[height=1cm]{graffles/Wcomm.pdf}$}
\end{equation*}
\begin{equation*}
\lower12pt\hbox{$\includegraphics[height=1cm]{graffles/Wassocl.pdf}$}
\!\!\!
=
\!\!\!
\lower12pt\hbox{$\includegraphics[height=1cm]{graffles/Wassocr.pdf}$}
\end{equation*}
\end{multicols}
\begin{multicols}{3}\noindent
\begin{equation*}
\lower9pt\hbox{$\includegraphics[height=.9cm]{graffles/Bcounitlaw.pdf}$}
\!\!\!
=
\!
\lower5pt\hbox{$\includegraphics[height=.6cm]{graffles/idcircuit.pdf}$}
\end{equation*}
\begin{equation*}
\lower5pt\hbox{$\includegraphics[height=.6cm]{graffles/Bcomult.pdf}$}
\!
=
\!\!\!
\lower11pt\hbox{$\includegraphics[height=1cm]{graffles/Bcomm.pdf}$}
\end{equation*}
\begin{equation*}
\lower11pt\hbox{$\includegraphics[height=1cm]{graffles/Bcoassocl.pdf}$}
\!\!\!
=
\!\!\!
\lower11pt\hbox{$\includegraphics[height=1cm]{graffles/Bcoassocr.pdf}$}
\end{equation*}
\end{multicols}
\begin{multicols}{4}
\noindent
\begin{equation*}
\lower3pt\hbox{$
\lower5pt\hbox{$\includegraphics[height=.6cm]{graffles/lunitsl.pdf}$}
=
\lower5pt\hbox{$\includegraphics[height=.6cm]{graffles/lunitsr.pdf}$}
$}
\end{equation*}\noindent\begin{equation*}
\lower3pt\hbox{$
\lower5pt\hbox{$\includegraphics[height=.6cm]{graffles/runitsl.pdf}$}
=
\lower5pt\hbox{$\includegraphics[height=.6cm]{graffles/runitsr.pdf}$}
$}
\end{equation*}
\begin{equation*}
\lower5pt\hbox{$\includegraphics[height=.6cm]{graffles/bialgl.pdf}$}
=
\lower10pt\hbox{$\includegraphics[height=.9cm]{graffles/bialgr.pdf}$}
\end{equation*}
\begin{equation*}
\lower4pt\hbox{$
\lower4pt\hbox{$\includegraphics[height=.5cm]{graffles/unitsl.pdf}$}
=\id_0\phantom{\quad}
$}
\end{equation*}
\end{multicols}
\begin{multicols}{4}\noindent
\begin{equation*}
\lower11pt\hbox{$\includegraphics[height=1cm]{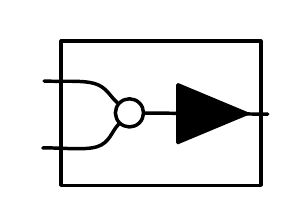}$}
\!\!
=
\!\!
\lower11pt\hbox{$\includegraphics[height=1cm]{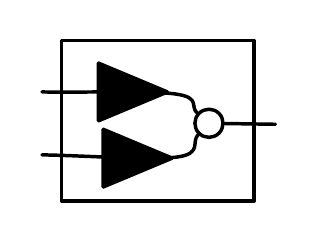}$}
\end{equation*}
\begin{equation*}
\lower10pt\hbox{$\includegraphics[height=.9cm]{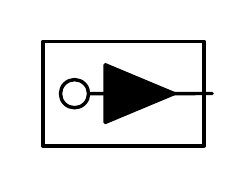}$}
\!\!
=
\!\!
\lower10pt\hbox{$\includegraphics[height=.9cm]{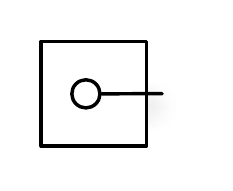}$}
\end{equation*}
\begin{equation*}
\lower11pt\hbox{$\includegraphics[height=1cm]{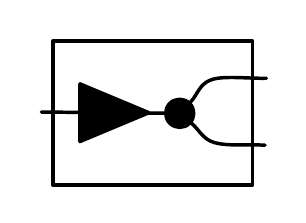}$}
\!\!
=
\!\!
\lower11pt\hbox{$\includegraphics[height=1cm]{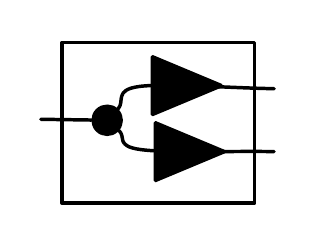}$}
\end{equation*}
\begin{equation*}
\lower10pt\hbox{$\includegraphics[height=.9cm]{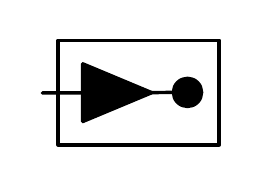}$}
\!\!
=
\!\!
\lower10pt\hbox{$\includegraphics[height=.9cm]{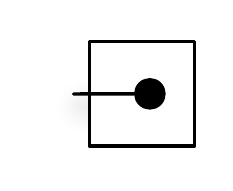}$}
\end{equation*}
\end{multicols}\vspace{-.3cm}
\begin{equation*}
  \lower11pt\hbox{$\includegraphics[height=1cm]{graffles/hopfr.pdf}$}
  = \lower8pt\hbox{$\includegraphics[height=.8cm]{graffles/hopfc.pdf}$}
   \end{equation*}

\begin{proposition}\label{prop:finitePresHAZ} $\HA{} \cong \HA{\Z}$.
\end{proposition}
\begin{proof}
Let $\alpha \: \HA{\Z} \to \HA{}$ be the PROP morphism defined on generators as follows. It is the identity on $\Bcounit$, $\Bcomult$, $\Wunit$ and $\Wmult$. For $k \in \Z$, $\alpha(\scalar)$ is given by:
\begin{align*}
\lower4pt\hbox{$\includegraphics[height=15pt]{graffles/zeroscalar.pdf}$} \mapsto  \lower4pt\hbox{$\includegraphics[height=15pt]{graffles/finitereprZ_zeror.pdf}$} %
 \ \ \ \ &&
  \lower10pt\hbox{$\includegraphics[height=25pt]{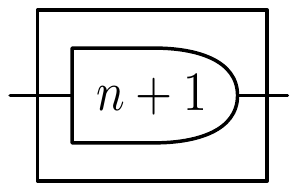}$} \mapsto \lower15pt\hbox{$\includegraphics[height=35pt]{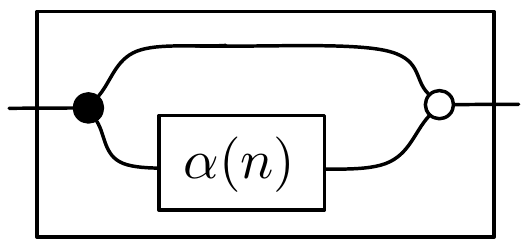}$} \ \ \ \  &&
  \lower10pt\hbox{$\includegraphics[height=25pt]{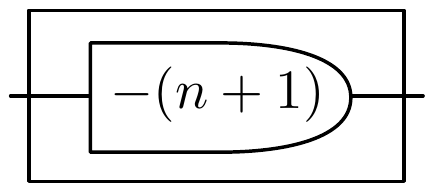}$} \mapsto \lower10pt\hbox{$\includegraphics[height=25pt]{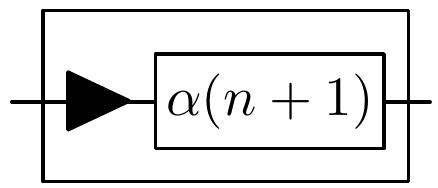}$}
\end{align*}
An inductive argument confirms that $\alpha$ is well-defined, in the sense that it preserves equality of diagrams in $\HA{\Z}$. Fullness is clear by construction. For faithfulness, just observe that all axioms of $\HA{}$ are also axioms of $\HA{\Z}$.
\end{proof}

A pleasant example of graphical reasoning in $\HA{}$ is the derivation showing that the antipode $\antipode$ is involutive:
$$\includegraphics[height=2.2cm]{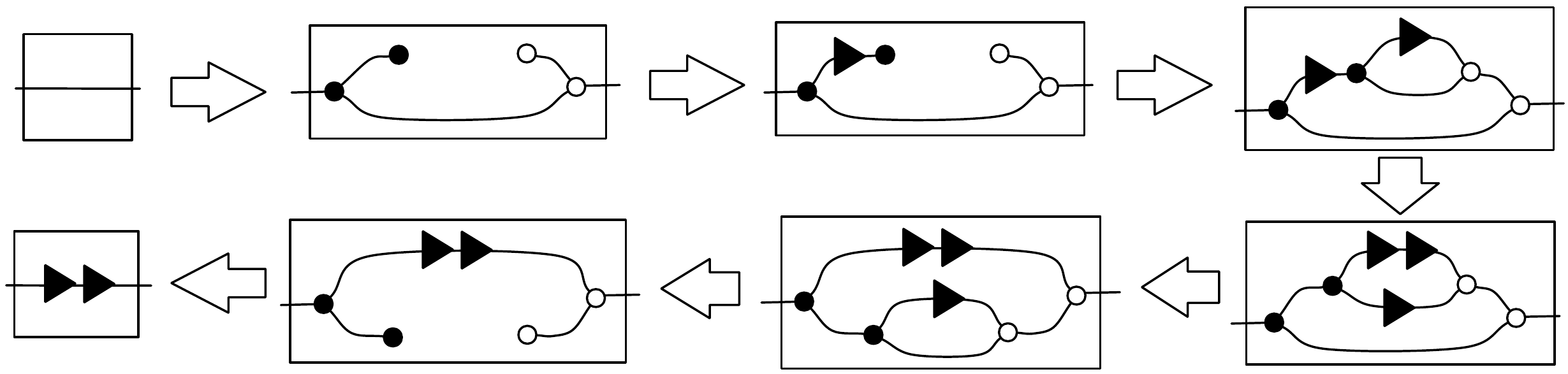}$$

\paragraph{The theory of rational subspaces} By Theorem~\ref{th:IBR=SVR}, $\IH{\Z}$ is isomorphic to the PROP $\SVH{\Q}$ of subspaces over the field $\Q$ of rational numbers. In view of Proposition~\ref{prop:finitePresHAZ}, we can give an alternative presentation of $\IH{\Z}$ based on the finite signature of $\HA{}+\HA{}^{\op}$: in axioms \eqref{eq:lcmIH}-\eqref{eq:lcmopIH}, $\scalar$ and $\coscalar$ become notational conventions for $\alpha(\scalar)$ and $\alpha^{\op}(\coscalar)$, respectively.

Note that, differently from the case of $\HA{}$, $\IH{\Z}$ has an \emph{infinite} presentation, because axiom schemas \eqref{eq:lcmIH}-\eqref{eq:lcmopIH} range over non-zero integers. An interesting observation, suggested to the author independently by Pawel Sobocinski and Peter Selinger, is that actually $\IH{\Z}$ \emph{cannot} be presented by a finite number of axioms.

\begin{proposition} $\IH{\Z}$ is not finitely axiomatisable.
\end{proposition}
\begin{proof} For $n \geq 0$, let $A_n$ be the set consisting of all the equations of $\HA{} + \HA{}^{\op}$,~\eqref{eq:WFrobIBR}-\eqref{eq:BSepIBR} and the instances of \eqref{eq:lcmIH}-\eqref{eq:lcmopIH} with $1 \leq l \leq n$. Then each $A_n$ is finite. Also, let $A$ be the union of all the $A_n$s. We verify that $A$ is a sound and complete axiomatisation for $\IH{\Z}$. It suffices to show \eqref{eq:lcmIH}-\eqref{eq:lcmopIH} for the case $l = -1$
\begin{equation*}
\lower6pt\hbox{$\includegraphics[height=.6cm]{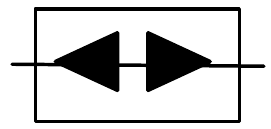}$}
 =
\lower9pt\hbox{$\includegraphics[height=.9cm]{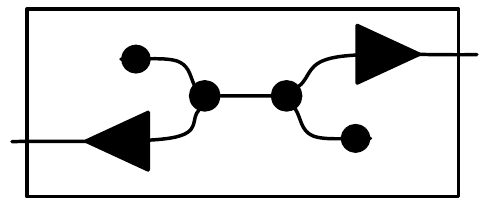}$}
 = 
\lower11pt\hbox{$\includegraphics[height=1.1cm]{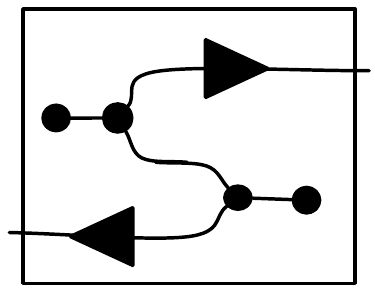}$}
 = 
\lower11pt\hbox{$\includegraphics[height=1.1cm]{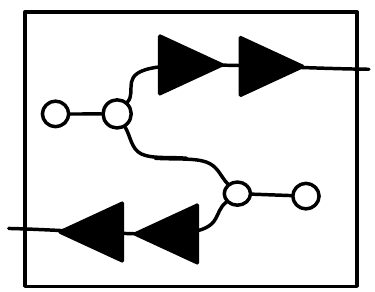}$}
= 
\lower11pt\hbox{$\includegraphics[height=1.1cm]{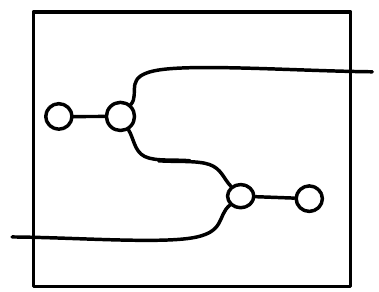}$}
=
\lower9pt\hbox{$\includegraphics[height=.9cm]{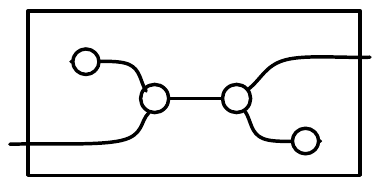}$}
 = 
\lower6pt\hbox{$\includegraphics[height=.6cm]{graffles/idcircuit.pdf}$}
\end{equation*}
\begin{equation*}
\lower6pt\hbox{$\includegraphics[height=.6cm]{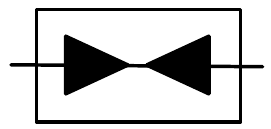}$}
 =
\lower6pt\hbox{$\includegraphics[height=.6cm]{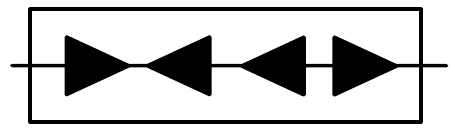}$}
 =
\lower6pt\hbox{$\includegraphics[height=.6cm]{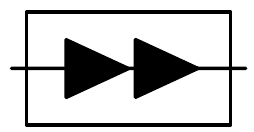}$}
 =
\lower6pt\hbox{$\includegraphics[height=.6cm]{graffles/idcircuit.pdf}$}
\end{equation*}
and for the case of a negative integer $l = -k$, where $k > 0$:
\begin{equation*}
\lower7pt\hbox{$\includegraphics[height=.7cm]{graffles/lcmopl_l.pdf}$}
 =
\lower6pt\hbox{$\includegraphics[height=.6cm]{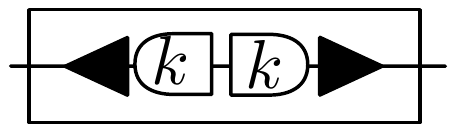}$}
 =
\lower6pt\hbox{$\includegraphics[height=.6cm]{graffles/lcmopAntipode.pdf}$}
 =
\lower6pt\hbox{$\includegraphics[height=.6cm]{graffles/idcircuit.pdf}$}
 =
\lower6pt\hbox{$\includegraphics[height=.6cm]{graffles/lcmAntipode.pdf}$}
 =
\lower6pt\hbox{$\includegraphics[height=.6cm]{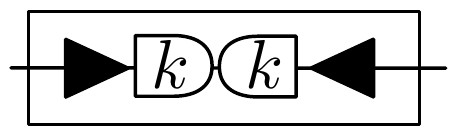}$}
=
\lower7pt\hbox{$\includegraphics[height=.7cm]{graffles/lcml_l.pdf}$}.
\end{equation*}
.

We now show that no single $A_n$ is complete. For a fixed $n$, let $p > n$ be a prime and consider the field $\Z_{\scriptscriptstyle p}$ of integers modulo $p$. Then all the equations of $A_n$ are sound for $\IH{\Z_{\scriptscriptstyle p}}$. However, for $l = p$,~\eqref{eq:lcmIH} is not sound in $\IH{\Z_{\scriptscriptstyle p}}$, whence it is also not a consequence of $A_n$. Therefore, no $A_n$ is complete for $\IH{\Z}$.

We can now use a compactness argument to conclude that there is no finite axiomatisation of $\IH{\Z}$. Suppose that $\widetilde{A}$ were such a finite set of equations. Since $A$ is complete, each equation $c = d$ in $\widetilde{A}$ follows from $A$. Now, the derivation of $c=d$ from $A$ uses only a finite number of equations in $A$, and since $\widetilde{A}$ contains finitely many equations, there must be a finite subset $A'$ of $A$ such that each equation of $\widetilde{A}$ is already derivable from $A'$. Then, by construction of $A$, we have that $A' \subseteq A_n$ for some $n$. This means that $A_n$ also implies all of $\widetilde{A}$ and thus is a complete axiomatisation for $\IH{\Z}$, a contradiction.
\end{proof}

\medskip

For a glimpse of the graphical reasoning in $\IH{\Z}$, we now give a combinatorial diagrammatic description of the subspaces of the 2-dimensional rational space (where $k_1,k_2$ are non-zero integers):
\begin{align}\label{eq:subspacesof2}
\spaceFull && \spaceZero && \spaceXaxis && \spaceYaxis && \spacekonektwo.
\end{align}
The diagram $\spaceFull$ denotes (via $\sem{\IH{\Z}}$) the full space $[ \tiny{\left(\begin{array}{c}
                \!\! 1 \!\! \\
                \!\! 0 \!\!
                \end{array}\right),\left(\begin{array}{c}
                \!\! 0 \!\! \\
                \!\! 1 \!\!
                \end{array}\right)}]$ and $\spaceZero$ the $0$-dimensional subspace $\{\tiny{\left(\begin{array}{c}
                \!\! 0 \!\! \\
                \!\! 0 \!\!
                \end{array}\right)}\}$. The remaining subspaces, all of dimension $1$, are conventionally represented as lines through the origin on the $2$-dimensional cartesian coordinate system. Three kinds of string diagrams suffice to represent all of them: $\spaceXaxis$ denotes the $x$-axis; $\spaceYaxis$ denotes the $y$-axis; for $k_1, k_2 \neq 0$, $\spacekonektwo$ denotes the line with slope $\frac{k_2}{k_1}$.

Conversely, using the modular structure of $\IH{\Z}$ we can check that the above combinatorial analysis~\eqref{eq:subspacesof2} covers all the $1 \to 1$ diagrams:
\begin{proposition} for all $c \in \IH{\Z}[1,1]$, $c$ is equal in $\IH{\Z}$ to a diagram in \eqref{eq:subspacesof2}.
\end{proposition}
\begin{proof}
By Proposition~\ref{Th:factIBR}, $c$ can be factorised as a cospan $1 \tr{c_1 \in \HA{}} n \tl{c_2 \in \HA{}} 1$. If we now take the pullback in $\HA{}$
\begin{equation}\label{eq:pb11}
\vcenter{
\xymatrix@=15pt{
r \ar[d]_{d_1} \ar[r]^{d_2} & 1 \ar[d]^{c_2}
\\ 1 \ar[r]_{c_1}& n}
}
\end{equation}
then $r$ must be either $0$, $1$ or $2$. We check that in all the cases $\varphi_2(\contrid{d_1})\poi\varphi_1(d_2)$ is equal to a diagram in~\eqref{eq:subspacesof2}.
\begin{itemize}
\item If $r=2$, then a pullback span of $\tr{c_1}\tl{c_2}$ is the one given by the projections $\lower4pt\hbox{$\includegraphics[height=15pt]{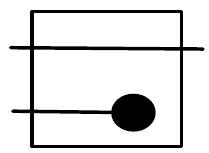}$} \: 2 \to 1$ and $\lower4pt\hbox{$\includegraphics[height=15pt]{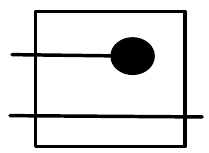}$} \: 2 \to 1$. Since $\tl{d_1}\tr{d_2}$ also pulls back $\tr{c_1}\tl{c_2}$, by Lemma~\eqref{lemma:mirror}, $\varphi_2(d_1)\poi\varphi_1(d_2)$ is equal in $\IH{\Z}$ to the composite $\spaceFull$ of $\contrid{\lower4pt\hbox{$\includegraphics[height=15pt]{graffles/proj1.pdf}$}} = \lower4pt\hbox{$\includegraphics[height=15pt]{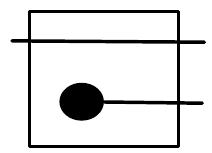}$}$ and $\lower4pt\hbox{$\includegraphics[height=15pt]{graffles/proj2.pdf}$}$.
\item If $r=1$, then $d_1$ and $d_2$ are $1\times 1$ matrices $(k_1)$, $(k_2)$ respectively. If $k_1=k_2=0$, then $\varphi_2(\contrid{d_1})\poi\varphi_1(d_2) = \spaceYaxis \poi \spaceXaxis = \spaceZero$ by \eqref{eq:zeroscalar}, \eqref{eq:zeroscalar}$^{\op}$, \eqref{eq:wbone}. If $k_1 = 0$ and $k_2 \neq 0$, then $d_1 = \spaceXaxis$  by \eqref{eq:zeroscalar} and $\varphi_2(\contrid{d_1})\poi\varphi_1(d_2) = \spaceYaxis \poi \scalarktwo= \spaceYaxis$ by \eqref{eq:scalarbcounit}$^{\op}$, \eqref{eq:lcmopIH}. Symmetrically, if $k_2 = 0$  and $k_1 \neq 0$ then $\varphi_2(\contrid{d_1})\poi\varphi_1(d_2) = \spaceXaxis$ by  by \eqref{eq:zeroscalar}$^{\op}$, \eqref{eq:scalarbcounit}, \eqref{eq:lcmopIH}. In case both $k_1$ and $k_2$ are different from $0$, $\varphi_2(d_1)\poi\varphi_1(d_2) = \spacekonektwo$.
\item Otherwise, $r$ is $0$. By initiality, $d_1 = d_2 = \Wunit$, meaning that $\varphi_2(\contrid{d_1})\poi\varphi_1(d_2) = \spaceZero$.
\end{itemize}
By Proposition~\ref{prop:IBwComplete}, the equation associated with the pullback \eqref{eq:pb11} holds in $\IH{\Z}^{\scriptscriptstyle w}$, meaning that $\varphi_1(c_1)\poi\varphi_2(\contrid{c_2}) = \varphi_2(\contrid{d_1})\poi\varphi_1(d_2)$ in $\IH{\Z}$. Since $c = \varphi_1(c_1)\poi\varphi_2(\contrid{c_2})$ by construction, this concludes the proof of the statement. \end{proof}

Notice that $\IH{\Z}[1,1]$ contains within its structure all of rational arithmetic:
$0$ can be identified with \spaceXaxis, and $\frac{k_2}{k_1}$, for $k_1\neq 0$, with \spacekonektwo. Multiplication
$\cdot\: \IH{\Z}[1,1]\times\IH{\Z}[1,1]\to \IH{\Z}[1,1]$
is composition $x\cdot y = x\poi y$, addition $+\: \IH{\Z}[1,1]\times\IH{\Z}[1,1]\to \IH{\Z}[1,1]$
is defined \[x+y = \Bcomult \poi (x \oplus y) \poi \Wmult.\]
Multiplication is associative
but not commutative in general: of course, it \emph{is} commutative when restricted to diagrams representing rationals.
Associativity and commutativity of addition follow from associativity and commutativity in $\bcom$ and $\wmon$.  
\chapter{The Calculus of Signal Flow Diagrams}\label{chapter:SFG}

\epigraphfontsize{\small\itshape}
\setlength\epigraphwidth{8cm}
\setlength\epigraphrule{0pt}
\epigraphfontsize{\small\itshape}
\epigraph{[T]he reason why physics has ceased to look for causes is that in fact there are no such things. The law of causality, I believe, like much that passes muster among philosophers, is a relic of a bygone age, surviving, like the monarchy, only because it is erroneously supposed to do no harm.}{--- \textup{Bertrand Russell}, On the Notion of Cause \textup{(1913)}.}

\section{Overview}\label{sec:intro}

Feedback and related notions such as self-reference and recursion are at the core of several disciplines, including computer science, engineering and control theory. In control, \emph{linear dynamical systems} are amongst the most extensively studied and well-understood classes of systems with feedback. They are signal transducers with two standard interpretations: \emph{discrete}, where---roughly speaking---signals come one after the other in the form of a stream, and \emph{continuous}, where signals are typically well-behaved real-valued functions.

From the earliest days, diagrams played a central role in motivating the subject matter.  Graphical representations were not merely intuitive, but also closely resembled physical manifestations (implementations) of linear dynamic systems, such as electrical circuits. While differing in levels of formality and minor technical details, the various notions share the same set of fundamental features---and for this reason we will group them all under the umbrella of \emph{signal flow graphs}. These features are: (i) the ability to copy, (ii) to add and (iii) to amplify signals, (iv) the ability to \emph{delay} a signal (in the discrete, stream-based interpretation) or to \emph{differentiate}/\emph{integrate} a signal (in the continuous interpretation), (v) the possibility of feedback loops and (vi) the concept of directed signal flow. 
Notably, while features (i)-(v) are usually present in physical manifestations, (vi) seems to have been included to facilitate human understanding as well as to avoid ``nonsensical'' diagrams where the intended signal flow seems to be incompatible or paradoxical. Of course, physical electrical wires do not insist on a particular orientation of electron flow; both are possible and the actual flow direction depends on the context.

In this chapter we introduce a string diagrammatic theory of signal flow graphs, which we call the \emph{signal flow calculus}. The syntax is based on the following operations, sequential ($\poi$)  and parallel ($\tns$) composition. Terms generated by the syntax are referred to as \emph{circuits}.
\vspace{-.2cm}
$$\begin{array}{lcrlcr}
\multicolumn{3}{c}{\stackrel{\FCeq}{\overbrace{\phantom{a}\qquad\qquad\qquad}}} &
\multicolumn{3}{c}{\stackrel{\FCopeq}{\overbrace{\phantom{a}\qquad\qquad\qquad}}} \\
\Bcomult\!\! & \!\!\circuitX\!\!\!\! & \Wmult\!\!  & \Wcomult\!\! & \!\!\circuitXop\!\!\!\! & \Bmult\\
\, \Bcounit\!\! & \scalar\!\! & \!\!\Wunit\!  &\,\Wcounit\!\!\! & \coscalar\!\! & \!\!\Bunit\,\\
\multicolumn{6}{c}{\stackrel{\underbrace{\phantom{aaaa}\qquad\qquad\qquad\qquad\qquad\qquad}}{\CD}}
\end{array}$$
We concentrate on the discrete interpretation; thus circuits are given a \emph{stream semantics}\footnote{The continuous interpretation will be the subject of a future work --- see Chapter~\ref{chapter:conclusion}.}. The intuition is that wires carry elements of a field $\field$ that enter and exit through boundary ports.  In particular, for circuits built from components in the leftmost three columns, which we refer to as being in $\FC$, the signal enters from the left and exits from the right boundary. Computation is synchronous, and at each iteration fresh elements are processed from input streams on the left and emitted as elements of output streams on the right. The basic components $\Bcomult$, $\Wmult$, $\scalar$ ($k\in\field$) and $\circuitX$ realise features (i)-(iv). That means, $\Bcomult$ \emph{duplicates} the input signal, $\Wmult$ \emph{sums} the two input signals and $\scalar$ \emph{multiplies} the signal by a scalar $k\in \field$. $\circuitX$ is a \emph{delay}: when a sequence of signals $k_0, k_1, k_2,\dots$ arrives on the left, it outputs the sequence $0,k_0,k_1\dots$ It can thus be thought as a synchronous one-place buffer initialised with $0$. The remaining components, $\Bcounit$ and $\Wunit$, are the units of $\Bcomult$ and $\Wmult$ respectively: $\Bcounit$ accepts any signal and discards it, while $\Wunit$ constantly outputs the signal $0$.

Each operation in $\FC$ has a symmetric counterpart in $\FCop$. Whereas in $\FC$ the signal flows left-to-right, in $\FCop$ it flows right-to-left: for instance, $\Bmult$ releases on the left two copies of the signal received on the right. The signal flow calculus $\CD$ consists of all the circuits obtained by freely combining circuits of $\FC$ and $\FCop$. Note that there is no primitive for recursion, but feedbacks are constructible at the level of $\CD$ by the use of ``bent identity wires'' $\Bunit \poi\Bcomult$ and $\Bmult\poi \Bcounit$. This allows us to recover the traditional notion of signal flow graph~\cite{mason1953feedback,Lahti} as the sub-class $\SFGform$ of $\CD$ formed by closing $\FC$ under feedbacks passing at least one delay. For instance,
\begin{equation*}
\includegraphics[height=1.4cm]{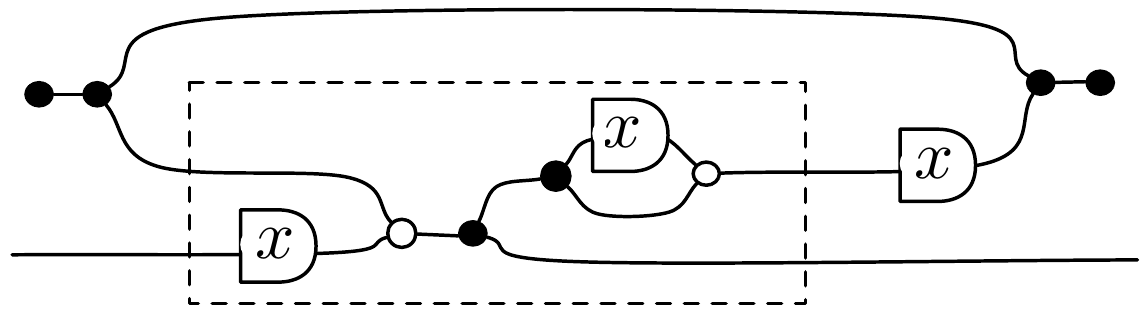}
\end{equation*}
is the circuit expressing the Fibonacci signal flow graph, that means, an input $1,0,0,\dots$ on the left produces the Fibonacci sequence $1,1,2,3,5,\dots$ on the right. It is in $\SFGform$ because it consists of a $\FC$-circuit (in the dotted square) equipped with a feedback passing through a delay. We will examine the Fibonacci circuit more in detail in Example~\ref{ex_fibonacci}. The intended execution behaviour of circuits in $\CD$ is formalised through the rules of a \emph{structural operational semantics} --- Fig.~\ref{fig:operationalSemantics} --- based on the intuitive explanation given above.

 We also give a \emph{denotational semantics} to the signal flow calculus. Circuits in $\FC$ and $\FCop$ are associated with linear \emph{functions} over streams, because we can always identify one boundary (left for $\FC$, right for $\FCop$) as the input boundary. This is no longer possible for circuits of $\CD$, which are built out by mixing the operations of $\FC$ and $\FCop$: for their denotational semantics, we need the generality of linear \emph{relations} over streams. We must also use an extended notion of streams, \emph{Laurent series}, typical in algebraic approaches~\cite{Barnabei19983} to signal processing---roughly speaking, these streams are allowed to start in the past. 



Our approach to the denotational semantics is based on the realisation that circuits of $\FC$ can be also interpreted as string diagrams of the theory $\ABpoly$ of $\poly$-Hopf algebras, where $\poly$ is the ring of polynomials with coefficients from $\field$ and unknown $x$. For $\FCop$, we use the dual theory $\ABpolyop$. For the whole $\CD$, the interpretation is in the theory $\IBpoly$ of Interacting $\poly$-Hopf algebras. In Chapter~\ref{chapter:hopf} we saw that $\ABpoly$ characterises $\poly$-matrices, whereas $\IBpoly$ characterises linear relations --- i.e., subspaces, \emph{cf.} Convention~\ref{conv:subspaceLinRel} --- over the field $\frpoly$ of \emph{fractions} of polynomials. Then, the passage to the stream semantics simply consists in interpreting polynomials and their fractions as streams --- see Table~\ref{tableStreams}. Since the characterisations for $\ABpoly$ and $\IBpoly$ are given both with an universal property and an inductive definition, we are able to present also the stream semantics of $\FC$ and $\CD$ in the two fashions. Figure~\ref{fig:roadmap} summarises the construction of the denotational semantics for the signal flow calculus.
\bgroup
\def\arraystretch{1.5}%
\begin{table}[t] \centering
\begin{tabular}{|l|l|l|}
\hline
 $\poly$ & the ring of polynomials & $\sum_0^n k_i x^i$ for some natural $n$\\ \hline
 $\frpoly$ & the field of fractions of polynomials & $\frac{p}{q}$ for $p,q\in \poly$ with $q\neq 0$ \\ \hline
 $\ratio$ & the ring of rationals & $\frac{\sum_0^n k_i x^i}{\sum_0^m l_j x^j}$ with $l_0 \neq 0$\\ \hline
 $\fps$ & the ring of formal power series & $\sum_0^{\infty} k_i x^i$\\ \hline
 $\laur$ & the field of Laurent series & $\sum_{d}^{\infty} k_i x^i$ for some interger $d$\\ \hline
 \end{tabular}
 \caption{Rings and fields depending on a field $\field$ ($k_i$ and $l_j$ range over $\field$).}\label{tableStreams}
\end{table}
\egroup
\begin{figure}[t]
$$ \xymatrix@R=15pt{
 \FCop \ar@{^{(}->}[d] \ar@{->>}[r] & \ABpolyop \ar@{^{(}->}[d] \ar[r]^{\cong} & \Matpolyop \ar@{^{(}->}[r] \ar@{^{(}->}[d] &  \Matfpsop \ar@{^{(}->}[d] \\
 \CD 
 \ar@{->>}[r] & \IBpoly \ar[r]^{\cong} & \SVpoly \ar@{^{(}->}[r] &  \SVfps \\
\FC \ar@{^{(}->}[u] \ar@{->>}[r] & \ABpoly \ar@{^{(}->}[u] \ar[r]^{\cong} & \Matpoly \ar@{^{(}->}[r] \ar@{^{(}->}[u] &  \Matfps \ar@{^{(}->}[u] } \qquad \qquad
\xymatrix@R=15pt{
 \CD  \ar@{->>}[r] & \IBpoly \ar[r]^-{\cong} & \SVpoly  \\
 \SFGform \ar@{^{(}->}[u] \ar@{->>}[r] & \SFG \ar@{^{(}->}[u] \ar[r]^-{\cong} & \Matratio  \ar@{^{(}->}[u] \\
\FC \ar@{^{(}->}[u] \ar@{->>}[r] & \ABpoly \ar@{^{(}->}[u] \ar[r]^-{\cong} & \Matpoly  \ar@{^{(}->}[u]
}$$
\caption{\label{fig:roadmap}A technical roadmap of results. Rings $\poly$, $\frpoly$, $\ratio$, $\fps$ and $\laur$ are described in Table~\ref{tableStreams}. 
 As in Chapter~\ref{chapter:hopf}, for a ring $\PID$ and field $\field$, $\Mat{\PID}$ and $\Subspace{\field}$ denote the PROP of $\PID$-matrices and $\field$-linear relations, respectively. In the diagrams, the double-headed arrows are the interpretation of syntax within an algebraic theory (i.e.\ quotienting w.r.t.\ a set of equations); the tailed arrows are embeddings and the arrows labeled with $\cong$ are isomorphisms.
The middle row in the left diagram is the factorization of the stream semantics $\CD \to \SVfps$. The diagram on the right shows the status of the class of orthodox signal flow diagrams $\SFGform$, where $\SFG$ is the quotient of $\SFGform$ by the equations of $\IBpoly$.} 
\end{figure}

Remarkably, $\ABpoly$ and $\IBpoly$ do not only give us a mean to define the stream semantics, but also a \emph{sound and complete axiomatisation} for it. Indeed, it immediately follows by our construction that the equations of $\ABpoly$ suffice to prove denotational equivalence of any two circuits of $\FC$, and the same for $\IBpoly$ and $\CD$.

A first application of completeness is the use of equational reasoning to recast, under a new light, fundamental results of signal flow diagrams theory. A well-known theorem (see e.g.~\cite{Lahti}) states that circuits which we identify as the sub-class $\SFGform$ represent exactly those behaviours expressible by matrices with entries from $\ratio$, the ring of \emph{rationals}: those fractions of polynomials where the constant term in the denominator is non-zero --- \emph{cf.} Table~\ref{tableStreams}.
We give a novel, graphical proof of this result (Theorem~\ref{th:SFGcharactRationals}). Differently from traditional approaches, our formulation features a formal syntax $\SFGform$ for diagrams and a complete set of axioms $\IBpoly$ for semantic equivalence, motivating the appellative of \emph{Kleene's theorem} for the characterisation of rational matrices.

Another well-known fact in signal flow diagrams theory is a normal form: every circuit is equivalent to one where all delays occur in the feedbacks. We give a concise graphical proof of this result (Proposition \ref{prop:normalformSFG}) based on the observation that feedbacks ``guarded'' by delays are a trace in the categorical sense~\cite{Selinger2009}.

%

\bigskip

The second part of the chapter focuses on comparing the operational perspective with the established denotational model. This question turns out to be quite subtle. 
In a sense, the denotational semantics is too abstract: finite computations that reach \emph{deadlocks} are ignored. Such deadlocks can arise for instance when circuits of $\FC$ are composed with the those of $\FCop$ and, intuitively, the signal flows from the left and right toward the middle. For an example, consider the circuit below on the left.
\begin{equation}\label{problems}
\circuitXcospan \hspace{3cm} \circuitXspan
\end{equation}
In a first step, the signals arriving from left and right are stored in the two buffers. Then, the stored values are compared in the middle of the circuit: if they do not agree then the computation gets stuck.
The circuit on the right features another problem, which we call \emph{initialisation}. Intuitively, the flow goes from the middle toward left and right. All its computations are forced to start by emitting on the left and on the right the value $0$ which is initially stored in the two buffers. The two circuits are denotationally equivalent --- and equal in $\IBpoly$, but their operational behaviour can be obviously distinguished: the leftmost does not have initialisation and the rightmost cannot deadlock.

Deadlock and initialisation are dual problems at the heart of the mismatch of operational and denotational semantics.
We show that circuits in \emph{cospan form}, namely circuits built from a circuit of $\FC$ followed by one of $\FCop$ --- like the left-hand circuit in \eqref{problems}, are free from initialisation. On the other hand, circuits in \emph{span form}, i.e., those built from a circuit of $\FCop$ followed by one of $\FC$ --- like the right-hand in \eqref{problems}, are free from deadlock.
This is interesting because our modular account of interacting Hopf algebras in Chapter~\ref{chapter:hopf} implies that any circuit is equivalent in $\IBpoly$ to both one in cospan and one in span form (Theorem~\ref{Th:factIBR}). The duality of deadlock and initialisation helps us in proving a \emph{full abstraction} theorem: for those circuits that are free from both deadlock and initialisation, the operational and the denotational semantics agree (Corollary~\ref{cor:fullabstractInitDeadFree}).

Our analysis spotlights circuits, like those in~\eqref{problems}, in which the mix of $\FC$- and $\FCop$-components makes impossible to coherently determine flow direction through the wires. For these diagrams, the operational semantics is not meant to describe the running of a state-machine, like if we had an input/output partition of ports, but rather describes a notion of ``equilibrium'' between boundaries. It becomes then questionable whether our approach presents a model which is computationally sensible: in absence of flow, how can we claim that the behaviour denoted by a circuit with deadlock/initialisations is really implemented in the signal flow calculus? 
Our answer to this question is a \emph{realisability} result (Theorem~\ref{thm:realisability}): every circuit diagram $c$ in $\CD$ can be transformed using the equations of $\IBpoly$ into at least one suitably \emph{rewired} circuit $d$ in $\SFGform$. Rewired circuits in $\SFGform$ are deadlock and initialisation free, meaning that their operational semantics is fully abstract; we can determine flow directionality in the wires of $d$ and compute its operational semantics as the step-by-step evolution of a state-machine. Therefore, $d$ can be really thought as an executable circuit, which properly \emph{realises} the behaviour denoted by $c$.

The realisability theorem is the culmination of our approach. Keeping the direction of signal flow out of definitions enabled us to propose a compositional model and disclose the algebraic landscape $\IBpoly$ underlying the signal flow calculus. Realisability ensures that our departure from the orthodoxy has no real tradeoff: whenever one is interested in a proper operational understanding of a diagram in $\CD$, there is a procedure that, using equational reasoning in $\IBpoly$, transforms it into an executable machine expressing the same stream transducer. On a different perspective, this result tells that the signal flow calculus is not more expressive than orthodox signal flow diagrams; viewed as transducers, they define the same class.




The chapter is concluded by showing how flow directionality can be formally treated as a \emph{derivative} notion of our theory. We define the \emph{directed} signal flow calculus, where the wires appearing in diagrams have an explicit orientation, and give an interpretation $\E$ of the directed calculus into $\CD$, which ``forgets'' flow directionality. The purpose is two-fold. First, using realisability, it allows us to observe that any behaviour denoted by a circuit of $\CD$ can be properly implemented, modulo $\E$, by some directed circuit. Second, using full abstraction, we are able to show that any two directed circuits which, under $\E$, are provable equal in $\IBpoly$, have the same operational behaviour: therefore, there is no harm in reasoning about signal flow graphs without explicit indications of signal flow.

The conclusion that we draw from our analysis is a re-evaluation of \emph{causality} as central ingredient for the theory of signal flow graphs. In 1953 Mason \cite{mason1953feedback} wrote: ``flow graphs differ from electrical network graphs in that their branches are directed. In accounting for branch directions it is necessary to take an entirely different line of approach from that adopted in electrical network topology.'' Instead, our results suggest that, like for electrical circuits, also for signal flow graphs directionality is \emph{not} a primitive notion as originally advocated by Mason.



\paragraph{Synopsis} Our exposition will be organised as follows.
\begin{itemize}
\item \S~\ref{sec:SFcalculus} introduces the syntax and the operational semantics of the signal flow calculus.
\item \S~\ref{sec:polysem} connects $\FC$ and $\CD$ to PROPs $\ABpoly$ and $\IBpoly$ respectively. This yields an interpretation for $\FC$ in terms of polynomial matrices (\S~\ref{sec:polysemFC}) and one for $\CD$ in terms of linear relations over fractions of polynomials (\S~\ref{sec:polysemCD}). 
\item \S~\ref{sec:stream} introduces the denotational stream semantics. Extending the polynomial interpretations of the previous section, we first consider circuits of $\FC$ (\S\ref{sec:streamHA}) and then generalise to circuits of $\CD$ (\S\ref{ss:streamIB}). We show soundness and completeness of $\IBpoly$ for denotational equivalence (Corollary~\ref{prop:isoIHsoundcomplete}).
\item \S~\ref{sec:SFG} shows that, up-to equality in $\IBpoly$, circuits of $\SFGform$ characterise the rational behaviours (Theorem~\ref{th:SFGcharactRationals}). \S~\ref{sec:trace} proves that, using the trace structure of $\IBpoly$, circuits of $\SFGform$ can be put in a normal form where delays only appear in the feedbacks (Proposition~\ref{prop:normalformSFG}).
\item \S~\ref{sec:fullabstract} compares denotational and operational equivalence. In \S\ref{sec:dualitydeadlockinit} we analyse the phenomena of deadlock and initialisation and give syntactic characterisations for circuits without these design flaws (Proposition~\ref{thm:spandeadlock} and~\ref{thm:cospaninit}). In~\ref{sec:trfls} we prove full abstraction for deadlock and initialisation free circuits (Corollary~\ref{cor:fullabstractInitDeadFree}).
\item \S~\ref{sec:realisability} proves the realisability theorem (Theorem~\ref{thm:realisability}, Corollary~\ref{cor:rea}) and investigates some of its consequences.
\item \S~\ref{sec:types} introduces a directed version of the signal flow calculus and shows how $\IBpoly$ can be used to reason about directed circuits.
\end{itemize}  

\section{Syntax and Operational Semantics}\label{sec:SFcalculus}


In this section we define the syntax and the structural operational semantics of a simple process calculus, to which we shall refer to as the \emph{signal flow calculus}.

\paragraph{Syntax} Throughout this chapter we fix an arbitrary field $\field$. The syntax, given below, does not feature  binding nor primitives for recursion, while $k$ ranges over $\field$. As we shall see, the indeterminate $x$ plays a formal role akin to that in the algebra of polynomials.
\smallskip
\begin{align}
c \bnfEq &  \BcounitT \bnfSep \BcomultT \bnfSep \scalarT \bnfSep \circuitXT   \bnfSep \WmultT \bnfSep \WunitT \bnfSep \label{eq:SFcalculusSyntax1} \\
& \,\BunitT \!\bnfSep \BmultT \bnfSep \scalaropT \bnfSep \!\circuitXopT \bnfSep  \!\WcomultT \bnfSep \WcounitT \bnfSep \label{eq:SFcalculusSyntax2}\\ 
& \,\ZeronetT \bnfSep \IdnetT \bnfSep \!\!\symNetT   \bnfSep c\poi c \bnfSep c \tns c \label{eq:SFcalculusSyntax3}
\end{align}

\begin{figure}[t]
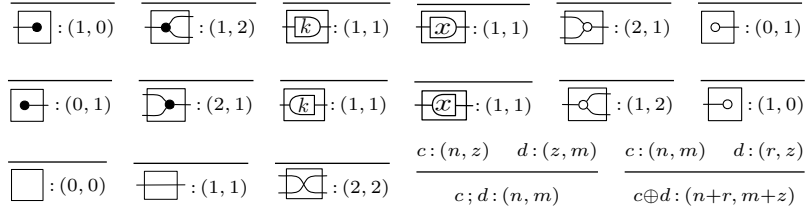

\[
\reductionRule{}{ \typeJudgment{}{\BcounitT}{\sort{1}{0}} }\quad
\reductionRule{}{ \typeJudgment{}{\BcomultT}{\sort{1}{2}} }\quad
\reductionRule{}{ \typeJudgment{}{\scalarT}{\sort{1}{1}} }\quad
\reductionRule{}{ \typeJudgment{}{\circuitXT}{\sort{1}{1}} }\quad
\reductionRule{}{ \typeJudgment{}{\WmultT}{\sort{2}{1}} }\quad
\reductionRule{}{ \typeJudgment{}{\WunitT}{\sort{0}{1}} }
\]
\[
\reductionRule{}{ \typeJudgment{}{\BunitT}{\sort{0}{1}} } \quad
\reductionRule{}{ \typeJudgment{}{\BmultT}{\sort{2}{1}} }\quad
\reductionRule{}{ \typeJudgment{}{\scalaropT}{\sort{1}{1}} }\quad
\reductionRule{}{ \typeJudgment{}{\circuitXopT}{\sort{1}{1}} }\quad
\reductionRule{}{ \typeJudgment{}{\WcomultT}{\sort{1}{2}} }\quad
\reductionRule{}{ \typeJudgment{}{\WcounitT}{\sort{1}{0}} }
\]
\[
\reductionRule{}{ \typeJudgment{}{\ZeronetT}{\sort{0}{0}} }\quad
\reductionRule{}{ \typeJudgment{}{\IdnetT}{\sort{1}{1}} }\quad
\reductionRule{}{ \typeJudgment{}{\symNetT}{\sort{2}{2}} }\quad
\reductionRule{ \typeJudgment{}{c}{\sort{n}{z}} \quad \typeJudgment{}{d}{\sort{z}{m}} }
{ \typeJudgment{}{c\poi d}{\sort{n}{m}} }\quad
\reductionRule{ \typeJudgment{}{c}{\sort{n}{m}} \quad \typeJudgment{}{d}{\sort{r}{z}} }
{ \typeJudgment{}{c \tns d}{\sort{n+r}{m+z}} }
\]
\caption{Sort inference rules.\label{fig:sortInferenceRules}}
\end{figure}

A \emph{sort} is a pair $\sort{n}{m}$, with $n,m\in \N$. We shall consider only terms that are sortable, according to the rules of Fig.~\ref{fig:sortInferenceRules}. A simple inductive argument confirms uniqueness of sorting: if $\typeJudgment{}{c}{\sort{n}{m}}$ and $\typeJudgment{}{c}{\sort{n'}{m'}}$ then $n=n'$ and $m=m'$.
We shall refer to sortable terms as \emph{circuits} since, intuitively, a term $\typeJudgment{}{c}{\sort{n}{m}}$ represents a circuit with $n$ ports on the left and $m$ ports on the right. Following the convention established in Remark~\ref{rmk:hopf}, we use notation $\antipode$ and $\antipodeop$ for $\scalarminusone$ and $\scalarminusoneop$ respectively. The reader may already notice a close relationship between circuits and the string diagrams of the theory of Interacting Hopf algebras studied in Chapter~\ref{chapter:hopf}: we shall explore this link in the subsequent sections.

\begin{figure}[t]
\[
\BcomultT \dtrans{\,k\,}{k \labelSep\labelSep k} \BcomultT\quad\quad
\BcounitT \dtrans{k}{} \BcounitT \quad\quad
\scalarT  \dtrans{\,\,l\,\,}{kl} \scalarT \quad \quad
\circuitXT ^{\labelSep l}\! \dtrans{k}{l} \circuitXT ^{\labelSep k} \quad \quad
\, \WmultT \dtrans{k \labelSep\labelSep l }{k+l} \WmultT \quad\quad
\WunitT \dtrans{\phantom{b}}{0} \WunitT
\]
\[
\BmultT\! \dtrans{k \labelSep\labelSep k}{k} \!\BmultT \quad\quad
\BunitT \dtrans{\phantom{k}}{k} {\BunitT} \quad\quad
\scalaropT  \dtrans{kl}{\,l\,} \scalaropT \quad \quad
\circuitXopT ^{\labelSep l}\! \dtrans{l}{k} \circuitXopT ^{\labelSep k} \quad \quad
\WcomultT\! \dtrans{ k+l}{k \labelSep\labelSep l }\! \WcomultT\quad\quad
\WcounitT \dtrans{0}{} \WcounitT
\]
\[
\IdnetT  \dtrans{k}{k} \IdnetT \quad\quad
\symNetT \dtrans{k \labelSep\labelSep l}{l \labelSep\labelSep k} \symNetT \quad\quad
\derivationRule{s\dtrans{\uu}{\vv} s' \quad t\dtrans{\vv}{\ww} t'}
{s\poi t \dtrans {\uu}{\ww} s' \poi t'}{} \quad \quad\quad
\derivationRule{s\dtrans{\uu_1}{\vv_1} s'\quad t\dtrans{\uu_2}{\vv_2} t'}
{s\tns t \dtrans{\uu_1 \labelSep \uu_2}{\vv_1 \labelSep \vv_2} s'\tns t'}{}
\]
\caption{Structural rules for operational semantics, with $k,l$ ranging over $\field$ and $\uu,\vv,\ww$ vectors of elements of $\field$ of the appropriate length\label{fig:operationalSemantics}.}
\end{figure}

\begin{remark}\label{remark:intuition}
Recalling the intuition established in~\S\ref{sec:intro}, we can consider circuits built up of the components in row \eqref{eq:SFcalculusSyntax1} as taking signals --- values in $\field$ --- from the left boundary to the right: thus $\Bcomult$ is a \emph{copier}, duplicating the signal arriving on the left; $\Bcounit$ accepts any signal on the left and discards it, producing nothing on the right; $\Wmult$ is an \emph{adder} that takes two signals on the left and emits their sum on the right, and $\Wunit$ constantly emits the signal $0$ on the right;  $\scalar$ is an \emph{amplifier}, multiplying the signal on the left by the scalar $k\in \field$. Finally, $\circuitX$ is a \emph{delay}, a synchronous one-place buffer initialised with $0$.

The terms of row \eqref{eq:SFcalculusSyntax2} are those of row \eqref{eq:SFcalculusSyntax1} reflected about the $y$-axis. Their behaviour is symmetric --- indeed, here it can be helpful to think of signals as flowing from right to left.
In row \eqref{eq:SFcalculusSyntax3}, $\symNetT$ is a \emph{twist}, swapping two signals, $\ZeronetT$ is the empty circuit and $\IdnetT$ is the \emph{identity} wire: the signals on the left and on the right ports are equal. Terms are combined with two binary operators: sequential ($\poi$) and parallel ($\tns$) composition. The intended behaviour of terms is now formalised through an operational semantics.
\end{remark}

\paragraph{Operational semantics} The operational semantics is a transition system with circuits augmented by states, where each delay component ($\circuitXT$ and $\circuitXopT$) is assigned some value $k \in \field$. Thus  \emph{states} are obtained by replacing the delays in the syntax specification with \emph{registers} $\circuitXT^{\labelSep k}$ and $\circuitXopT^{\labelSep k}$ for each $k \in \field$. We shall conveniently refer to a state for a circuit $c$ as a \emph{$c$-state}. As for circuits, we only consider sortable states, 
which are defined by adding \begin{center}$\typeJudgment{}{\circuitXT^{\labelSep k}}{\sort{1}{1}} $ and $\typeJudgment{}{\circuitXopT^{\labelSep k}}{\sort{1}{1}}$.\end{center}
to the rules in Fig.~\ref{fig:sortInferenceRules}.

Structural inference rules for operational semantics are given in Fig.~\ref{fig:operationalSemantics} where we use strings of length $n$ to represent vectors in $\field^n$. So, the empty string stands for $\matrixNull$, the only vector of $\field^0$, and $\vv=k_1\dots k_n$ for the column vector ${\tiny\left(%
				\begin{array}{c}
				  \!\!\!k_1\!\!\! \\
				  \!\!\!\vdots\!\!\! \\
				  \!\!\!k_n\!\!\!
				\end{array}\right)}$ in $\field^n$. Also, we use notation $\zerov$ for the string representing the vector where each value is $0$. When not explicit, the length of $\zerov$ will be evident from the context.

If state $s \mathrel{:} \sort{n}{m}$ is the source of a transition $\dtrans{\vv}{\ww} t$ then $t$ is also a state with sort $\sort{n}{m}$ and
$\vv$ and $\ww$ are strings representing vectors of $\field^n$ and $\field^m$, respectively.
Intuitively, $s \dtrans{\vv}{\ww} t$ means that $s$ can become $t$ whenever the signals on the $n$ ports on the left agree with $\vv$ and the signals on the $m$ ports on the right agree with $\ww$. Each circuit $c$ then yields a transition system with a chosen \emph{initial state} $s_0$ of $c$, obtained by replacing the delays $\circuitXT$ and $\circuitXopT$ in $c$ with registers $\circuitXT^{\labelSep 0}$ and $\circuitXopT^{\labelSep 0}$ containing $0$.

A \emph{computation} of a circuit $c$, is a (possibly infinite) path
$s_0 \dtrans{\vv_0}{\ww_0} s_1 \dtrans{\vv_1}{\ww_1} \dots$ in the
transition system of $c$, starting from its initial state $s_0$.
When $c$ has sort $\sort{n}{m}$, each $\vv_i$ is a string $k_{i1} \dots k_{in}$ representing a $\field$-vector of length $n$ and $\ww_i$ is a string $l_{i1} \dots l_{im}$ representing a $\field$-vector of length $m$. The \emph{trace} of a computation $s_0 \dtrans{\vv_0}{\ww_0} s_1 \dtrans{\vv_1}{\ww_1} \dots$ is then a pair of vectors ${\tiny \left(%
				\begin{array}{c}
				  \!\!\!\alpha_1\!\!\! \\
				  \!\!\!\vdots\!\!\! \\
				  \!\!\!\alpha_n\!\!\!
				\end{array}\right), \left(%
				\begin{array}{c}
				  \!\!\!\beta_1\!\!\! \\
				  \!\!\!\vdots\!\!\! \\
				  \!\!\!\beta_m\!\!\!
				\end{array}\right)}$ of sequences $\alpha_j=k_{0j} k_{1j} \dots$ and $\beta_j=l_{0j} l_{1j} \dots $. Occasionally we will use the notation $(\vlist{\alpha},\vlist{\beta})$ for such a pair. 
Moreover, we write $\alpha_j(i)$ and $\beta_j(i)$ for the $i$-th elements of $\alpha_j$ and $\beta_j$.

Note that, in a computation of length $z$, all $\alpha_j, \beta_j$ have length $z$, while for an infinite computation all $\alpha_j,\beta_j$ are infinite.
In the former case, we say that a trace is \emph{finite}, in the latter that it is \emph{infinite}.
We use $\gls{ftrC}$ to denote the set of all finite traces of $c$ and $\gls{itrC}$ for the set of all infinite ones. Properties of traces will be studied in more depth in \S\ref{sec:trfls}.

\begin{example}\label{exm:opsem}
Consider the two circuits below.
\begin{equation}\label{eq:exopsem}
{\lower4pt\hbox{$\includegraphics[height=.7cm]{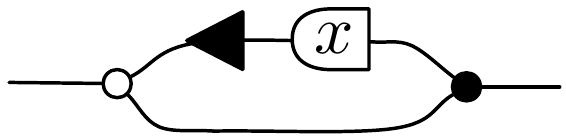}$}}
\qquad\qquad
{\lower10pt\hbox{$\includegraphics[height=1cm]{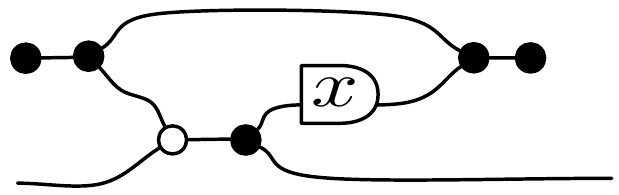}$}}
\end{equation}
The first is a graphical representation of the term
\[
c_1 = (\WcomultT \poi ((\antipodeop \poi \circuitXopT)\tns \IdnetT)) \poi \BmultT
\]
the second of the term
\begin{multline*}
c_2 = ( (\BunitT\poi \BcomultT) \tns \IdnetT) \poi
(\IdnetT \tns (\WmultT \poi \BcomultT)) \\
\poi (((\IdnetT \tns \circuitXT) \tns \IdnetT)
\poi ( (\BmultT \poi \BcounitT) \tns \IdnetT))
\end{multline*}
where we adopted the convention of depicting $\poi$ as horizontal and $\tns$ as vertical juxtaposition of diagrams --- we shall comment below on the adequacy of the graphical representation.
Note that, according to our intuition, in the left-hand circuit the signal flows from right to left, while the right-hand, the signal flows from left to right. Indeed, the terms $\BunitT\poi \BcomultT$
and $\BmultT \poi \BcounitT$ serve as bent identity wires which allow us to  form a feedback loop --- this idea, which was also raised when shaping the compact closed structure of $\IBR$ in \S\ref{sec:cc}, can be now made formal using the the rules in Fig.~\ref{fig:operationalSemantics}:
\[
\derivationRule{\raise3pt\hbox{$\BunitT\dtrans{}{k} \BunitT \quad \BcomultT\dtrans{\, k \,}{k \labelSep k} \BcomultT$}}
{\rccB \dtrans {}{k \labelSep k} \rccB}{} \quad \quad
\derivationRule{\raise3pt\hbox{$\BmultT\dtrans{k \labelSep k}{\, k \,} \BmultT \quad\BcounitT\dtrans{k}{} \BcounitT$}}
{\lccB \dtrans {k \labelSep k}{} \lccB}{} \quad \quad k \in \field.
\]
We now describe the operational behaviour of the diagrams in~\eqref{eq:exopsem}. Let $c_1[k]$ and $c_2[k]$ represent the states of $c_1$ and $c_2$, with $k$ denoting the value at the register. The rules of Fig.~\ref{fig:operationalSemantics} yield the computation
\[
c_i[0] \dtrans{1}{1} c_i[1] \dtrans{0}{1} c_i[1] \dtrans{0}{1} c_i[1] \cdots
\]
for $i\in 0,1$,
which yields the trace $(1000\dots),(1111\dots)$. In fact, as we shall show via a sound and complete axiomatisation, despite of the signal intuitively flowing in different directions, the two circuits in~\eqref{eq:exopsem} have the same observable behaviour.

A slightly more involved example is given below.
\[
\includegraphics[height=1cm]{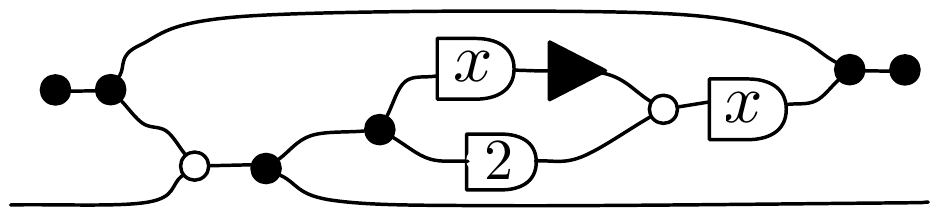}
\]
We leave the reader to write down a term that is represented by the diagram above: call it $c_3$ and let $c_3[k_1,k_2]$ represent the state where the two registers, reading from from left to right, have values $k_1$ and $k_2$. Then, the operational semantics allows us to derive the following computation
\[
c_3[0,0] \dtrans{1}{1}
c_3[1,2] \dtrans{0}{2}
c_3[2,3] \dtrans{0}{3}
c_3[3,4] \dtrans{0}{4} \cdots
\]
that yields the trace $(1000\dots,1234\dots)$.
\end{example}

\paragraph{From syntax to string diagrams}\label{sec:circuitdiag}
In Example~\ref{exm:opsem} we used the graphical language of \emph{string diagrams} to represent syntactic terms. As explained in \S~\ref{sec:props}, we can view these diagrams as the arrows  of a symmetric strict monoidal category, which we call $\CD$. Objects of $\CD$ are the natural numbers. Arrows $n \to m$ are called \emph{circuit diagrams}: they are circuit terms $c \: (n,m)$ quotiented by the laws of symmetric strict monoidal categories (Fig.~\ref{fig:axSMC}). The identity $0 \to 0$ is $\ZeronetT$, the identity $1 \to 1$ is $\IdnetT$ and the symmetry $1+1 \to 1+1$ is $\symNetT$: all the other identities and symmetries arise by their composition. As in Example~\ref{exm:opsem}, composition of arrows is meant to be represented graphically by horizontal (for $\poi$) and vertical (for $\tns$) juxtaposition of circuit diagrams. Laws of symmetric strict monoidal categories guarantee that this representation is sound (see \S~\ref{sec:props}).

We can use the theory of PROPs to give a more succinct definition of $\CD$. 
\begin{definition}\label{def:cicruitdiag}
The PROP $\gls{PROPCD}$ of \emph{circuit diagrams} is freely generated by the signature consisting of generators in \eqref{eq:SFcalculusSyntax1}-\eqref{eq:SFcalculusSyntax2} and no equations.
\end{definition}
Note that arrows of $\CD$ are in fact constructed using all the basic components of the signal flow calculus: the generators in~\eqref{eq:SFcalculusSyntax3} do not appear in Definition~\ref{def:cicruitdiag} because they are built-in building blocks of any freely generated PROP, just as the laws of symmetric strict monoidal categories.

Clearly, any circuit can be graphically rendered as a circuit diagram, as we did in Example~\ref{exm:opsem}, but the syntax carries more information than the diagrammatic notation (e.g. associativity). From the point of view of operational behaviour, however, this extra information is irrelevant. More precisely, one can easily check that, for any axiom $c \feq d$ of symmetric strict monoidal categories in Fig.~\ref{fig:axSMC}, with variables $t_1, t_2, t_3, t_4$ now standing for circuits of the appropriate sort, the circuits $c$ and $d$ yield isomorphic transition systems\footnote{By isomorphism of transition systems we mean a bijection between state-paces that preserves and reflects transitions and initial state.}. As the equations of Fig.~\ref{fig:axSMC} are the only ones valid in $\CD$, this observation justifies our use of the graphical notation and makes harmless to reason up to the laws of symmetric strict monoidal categories. For this reason, in the rest of the chapter we shall refer to circuits diagrams and state diagrams just as \emph{circuits} and \emph{states}, purposefully blurring the line between diagrams and traditional syntax.

For our developments it is useful to identify two sub-categories of $\CD$: $\gls{PROPCDleftright}$ has as arrows only those circuits in $\CD$ that are built from the components in rows \eqref{eq:SFcalculusSyntax1} and \eqref{eq:SFcalculusSyntax3} and $\gls{PROPCDrightleft}$ only those circuits built from the components in rows \eqref{eq:SFcalculusSyntax2} and \eqref{eq:SFcalculusSyntax3}. Equivalently, one can define $\FC$ and $\FCop$ as the PROPs freely generated by the signature~\eqref{eq:SFcalculusSyntax1} and~\eqref{eq:SFcalculusSyntax2} respectively. The notation recalls the intuition that for circuits in $\FC$, signal flow is from left to right, and in $\FCop$ from right to left. Formally, observe that $\FCop$ is the opposite category of $\FC$: any circuit of $\FCop$ can be seen as one of $\FC$ reflected about the $y$-axis. Also, note that $\CD$ is the \emph{sum} $\FC + \FCop$ --- \emph{cf.} \S~\ref{sec:coproduct}.

We say that $\tr{c \in \CD}$ is in \emph{cospan form} if it is of shape $\tr{c_1 \in \FC}\tr{c_2 \in \FCop}$. Dually, it is in \emph{span form} if it is of shape $\tr{c_1 \in \FCop}\tr{c_2 \in \FC}$.




\paragraph{Feedback and signal flow diagrams} Beyond $\FC$ and $\FCop$, we identify another class of circuits of $\CD$ that adhere closely to the orthodox notion of \emph{signal flow diagram} (see \emph{e.g.} \cite{mason1953feedback}), albeit without directed wires. Here, the signal can flow from left to right, as in $\FC$, but with the possibility of \emph{feedbacks}, provided that these pass through at least one delay. This amounts to defining, for all $n$, $m$, a map $\Tr{}(\cdot) \: \CD[n+1,m+1] \to \CD[n,m]$ taking $c \: n+1 \to m+1$ to the $n$-to-$m$ circuit below:
\begin{equation}\label{eq:onefeedbackcircuit}
\lower12pt\hbox{$\includegraphics[height=1cm]{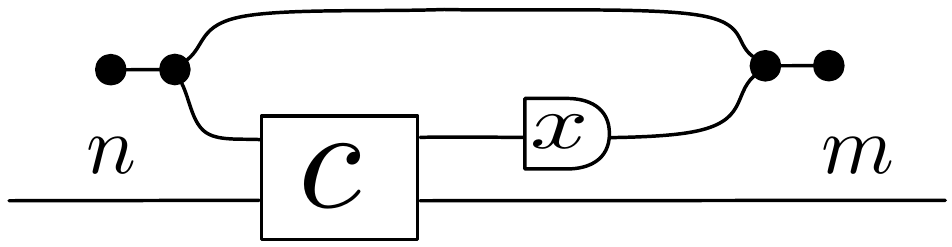}$}
\end{equation}
where $\idzcircuit$ is shorthand notation for the circuit $\id_z$ --- see \S\ref{sec:cc}. Intuitively, $\Tr{}(\cdot)$ equips $c$ with a feedback loop carrying the signal from its topmost right to its topmost left port. 

Circuits of this shape form a PROP $\gls{PROPSFGform}$, defined as the sub-category of $\CD$ inductively given as follows:
\begin{enumerate}[(i)]
\item if $c \in \FC[n,m]$, then $c\in\SFGform[n,m]$
\item if $c \in \SFGform[n+1,m+1]$, then $\Tr{}(c)\in\SFGform[n,m]$
\item if $c_1 \in \SFGform[n,z]$ and $c_2\in \SFGform[z,m]$, then $c_1 \poi c_2\in\SFGform[n,m]$
\item if $c_1 \in \SFGform[n,m]$ and $c_2\in \SFGform[r,z]$, then $c_1 \tns c_2 \in\SFGform[n+r,m+z]$.
\end{enumerate}
Equivalently, $\SFGform$ is the smallest sub-PROP of $\CD$ that contains $\FC$ and is closed under the $\Tr{}(\cdot)$ operation.
For instance, the right-hand circuit in \eqref{eq:exopsem} is in $\SFGform$, whereas the left-hand is in $\FCop$.

\begin{remark}\label{rmk:signalflowOpSem} The rules of Figure~\ref{fig:operationalSemantics} describe the step-by-step evolution of state machines without relying on a fixed flow orientation. This operational semantics is \emph{not} meant to be executable for all circuits: the rule for sequential composition implicitly quantifies existentially on the middle value $\vv$, resulting in potentially unbounded non-determinism. However, for circuits where flow directionality can be assigned, like the class $\SFGform$ above, existential quantification becomes deterministic subject to a choice of inputs to the circuit at each step of evaluation.
We will see in \S~\ref{sec:realisability} that any circuit can be transformed into this form, where the valid transformations are those allowed by the equational theory of Interacting Hopf algebras --- \emph{cf.} Remark~\ref{rmk:executableOpsem}.
\end{remark}

Remark~\ref{rmk:signalflowOpSem} emphasises that $\SFGform$ is a convenient setting for operational considerations on diagrams --- we will indeed pursue this perspective later in the chapter, see e.g. Proposition~\ref{prop:SFGfree}. Nonetheless, $\CD$ is still preferable for investigating the mathematical properties of circuits. Indeed, $\SFGform$ is not generated by any SMT, as $\Tr{}(\cdot)$ cannot be expressed as an operation with an arity and coarity. On the other hand, not only is $\CD$ generated by an SMT, but it is based on the same signature as the theory of interacting Hopf algebras (Def.~\ref{def:IBR}). This observation will enable us to develop a rich mathematical theory for $\CD$ based on the results of the previous chapter.

We shall return to the comparison between $\CD$ and $\SFGform$ in \S\ref{sec:realisability}.

\section{Denotational Semantics I: Polynomials}\label{sec:polysem}

In this section we commence our investigation of the denotational semantics of the signal flow calculus. We shall first give a semantics to $\FC$ and $\FCop$ in terms of polynomial matrices (\S~\ref{sec:polysemFC}) and to $\CD$ in terms of linear relations over the field of fractions of polynomials (\S~\ref{sec:polysemCD}). We will take advantage of the theory developed in the previous chapter to give a sound and complete axiomatisation for both these interpretations.


\subsection{Polynomial Semantics of $\FC$}\label{sec:polysemFC}

\begin{figure}[t]
%
%
%
%
%

\begin{multicols}{3}\noindent
\begin{align}
\notag
\lower10pt\hbox{$\includegraphics[height=1cm]{graffles/Wassocl.pdf}$}
\!\!\eqrule{A1}
\!\!
\lower10pt\hbox{$\includegraphics[height=1cm]{graffles/Wassocr.pdf}$}
\end{align}
\begin{align}
\notag
\lower5pt\hbox{$\includegraphics[height=.6cm]{graffles/Wmult.pdf}$}
\eqrule{A2}
\!\!\!
\lower11pt\hbox{$\includegraphics[height=1cm]{graffles/Wcomm.pdf}$}
\end{align}
\begin{align}
\notag
\lower10pt\hbox{$\includegraphics[height=.9cm]{graffles/Wunitlaw.pdf}$}
\!\!
\eqrule{A3}
\lower5pt\hbox{$\includegraphics[height=.6cm]{graffles/idcircuit.pdf}$}
\end{align}
\end{multicols}
\begin{multicols}{3}\noindent
\begin{align}
\notag
\lower10pt\hbox{$\includegraphics[height=1cm]{graffles/Bcoassocl.pdf}$}
\!\!\eqrule{A4}\!\!
\lower10pt\hbox{$\includegraphics[height=1cm]{graffles/Bcoassocr.pdf}$}
\end{align}
\begin{align}
\notag
\lower5pt\hbox{$\includegraphics[height=.6cm]{graffles/Bcomult.pdf}$}
\eqrule{A5}
\!\!\!
\lower11pt\hbox{$\includegraphics[height=1cm]{graffles/Bcomm.pdf}$}
\end{align}
\begin{align}
\notag
\lower10pt\hbox{$\includegraphics[height=.9cm]{graffles/Bcounitlaw.pdf}$}
\!\!
\eqrule{A6}
\lower6pt\hbox{$\includegraphics[height=.6cm]{graffles/idcircuit.pdf}$}
\end{align}
\end{multicols}
\begin{multicols}{4}
\noindent
\begin{align}
\notag
\lower6pt\hbox{$
\lower6pt\hbox{$\includegraphics[height=.6cm]{graffles/lunitsl.pdf}$}
\eqrule{A7}
\lower6pt\hbox{$\includegraphics[height=.6cm]{graffles/lunitsr.pdf}$}
$}
\end{align}\noindent\begin{align}
\notag
\lower6pt\hbox{$
\lower6pt\hbox{$\includegraphics[height=.6cm]{graffles/runitsl.pdf}$}
\eqrule{A8}
\lower6pt\hbox{$\includegraphics[height=.6cm]{graffles/runitsr.pdf}$}
$}
\end{align}
\begin{align}
\notag
\lower5pt\hbox{$\includegraphics[height=.6cm]{graffles/bialgl.pdf}$}
\eqrule{A9}
\lower9pt\hbox{$\includegraphics[height=.9cm]{graffles/bialgr.pdf}$}
\end{align}
\begin{align}
\notag
\lower7pt\hbox{$
\lower4pt\hbox{$\includegraphics[height=.5cm]{graffles/unitsl.pdf}$}
\eqrule{A10}
\lower4pt\hbox{$\includegraphics[height=.5cm]{graffles/idzerocircuit.pdf}$}
$}
\end{align}
\end{multicols}
\begin{multicols}{4}\noindent
\begin{align}
\notag
\lower8pt\hbox{$
\lower4pt\hbox{$\includegraphics[height=.5cm]{graffles/unitscalar.pdf}$}
\eqrule{A11}
\lower4pt\hbox{$\includegraphics[height=.5cm]{graffles/idcircuit.pdf}$}
$}
\end{align}
\begin{align}
\notag
\lower8pt\hbox{$
\lower5pt\hbox{$\includegraphics[height=.6cm]{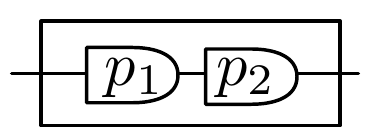}$}
\eqrule{A12}
\lower5pt\hbox{$\includegraphics[height=.6cm]{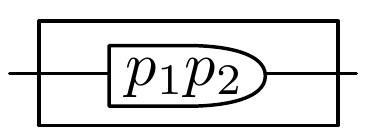}$}
$}
\end{align}
\begin{align}
\notag
\lower6pt\hbox{$
\lower8pt\hbox{$\includegraphics[height=.7cm]{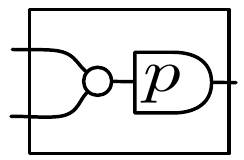}$}
\eqrule{A13}
\lower8pt\hbox{$\includegraphics[height=.7cm]{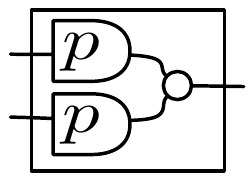}$}
$}
\end{align}
\begin{align}
\notag
\lower6pt\hbox{$
\lower6pt\hbox{$\includegraphics[height=.6cm]{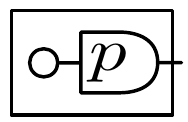}$}
\eqrule{A14}
\lower6pt\hbox{$\includegraphics[height=.6cm]{graffles/Wunit.pdf}$}
$}
\end{align}
\end{multicols}
\begin{multicols}{4}
\noindent
\begin{align}
\notag
\lower18pt\hbox{$\includegraphics[height=.7cm]{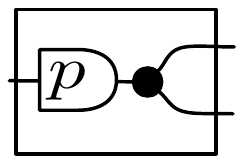}$}
\lower12pt\hbox{$\eqrule{A15}$}
\lower18pt\hbox{$\includegraphics[height=.7cm]{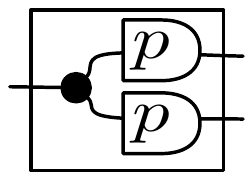}$}
\end{align}
\begin{align}
\notag
\lower12pt\hbox{$
\lower6pt\hbox{$\includegraphics[height=.6cm]{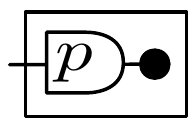}$}
\eqrule{A16}
\lower6pt\hbox{$\includegraphics[height=.6cm]{graffles/Bcounit.pdf}$}
$}
\end{align}
\begin{align}
\notag
\lower12pt\hbox{$
\lower4pt\hbox{$\includegraphics[height=.5cm]{graffles/zeroscalar.pdf}$}
\eqrule{A17}\!\!
\lower6pt\hbox{$\includegraphics[height=.7cm]{graffles/zeroscalar2.pdf}$}
$}
\end{align}
\begin{align}
\notag
\lower8pt\hbox{$\includegraphics[height=.9cm,width=1.3cm]{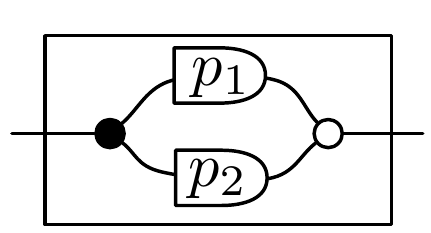}$}
\eqrule{A18}
\lower6pt\hbox{$\includegraphics[height=.6cm,width=1.1cm]{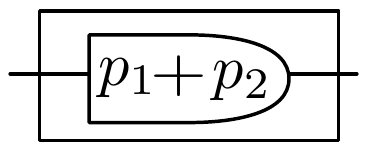}$}
\end{align}
\end{multicols}
\caption{Axioms of $\HA{\poly}$, describing the interaction of the generators in \eqref{eq:SFcalculusSyntax1}.}\label{fig:axiomsHA}
\end{figure}

In \S\ref{sec:theorymatr} we introduced the string diagrammatic theory for matrices over any principal ideal domain. It is now instrumental to instantiate our approach to the ring $\gls{poly}$ of polynomials with unknown $x$ and values over $\field$. By Proposition~\ref{prop:ab=vect}, the PROP $\Matpoly$ of $\poly$-matrices is presented by the generators and equations of the PROP $\HA{\poly}$ of $\poly$-Hopf algebras (Def.~\ref{def:HA}). We report the equational theory of $\HA{\poly}$ in Figure~\ref{fig:axiomsHA}. 

\begin{figure}[t]
%
%
%
%
%
\begin{multicols}{3}\noindent
\begin{align}
\notag
\lower8pt\hbox{$\includegraphics[height=.7cm]{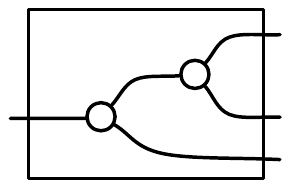}$}
\eqrule{\rev{A1}}
\lower8pt\hbox{$\includegraphics[height=.7cm]{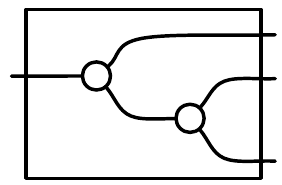}$}
\end{align}
\begin{align}
\notag
\lower5pt\hbox{$\includegraphics[height=.6cm]{graffles/Wcomult.pdf}$}
\eqrule{\rev{A2}}
\!\!\!
\lower6pt\hbox{$\includegraphics[height=.7cm]{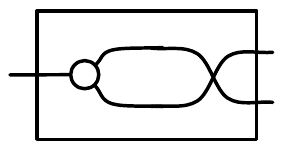}$}
\end{align}
\begin{align}
\notag
\lower6pt\hbox{$\includegraphics[height=.65cm]{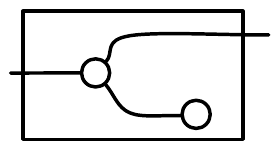}$}
\!\!
\eqrule{\rev{A3}}
\lower5pt\hbox{$\includegraphics[height=.5cm]{graffles/idcircuit.pdf}$}
\end{align}
\end{multicols}
\begin{multicols}{3}\noindent
\begin{align}
\notag
\lower8pt\hbox{$\includegraphics[height=.7cm]{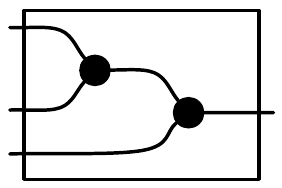}$}
\eqrule{\rev{A4}}
\lower8pt\hbox{$\includegraphics[height=.7cm]{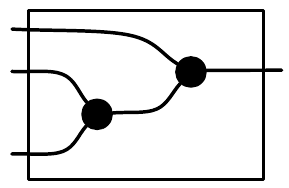}$}
\end{align}
\begin{align}
\notag
\lower5pt\hbox{$\includegraphics[height=.6cm]{graffles/Bmult.pdf}$}
\eqrule{\rev{A5}}
\!\!\!
\lower6pt\hbox{$\includegraphics[height=.7cm]{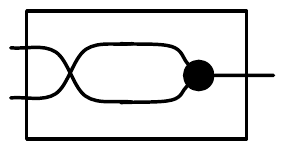}$}
\end{align}
\begin{align}
\notag
\lower6pt\hbox{$\includegraphics[height=.65cm]{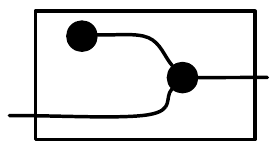}$}
\!\!
\eqrule{\rev{A6}}
\lower5pt\hbox{$\includegraphics[height=.6cm]{graffles/idcircuit.pdf}$}
\end{align}
\end{multicols}
\begin{multicols}{4}
\noindent
\begin{align}
\notag
\lower6pt\hbox{$
\lower6pt\hbox{$\includegraphics[height=.6cm]{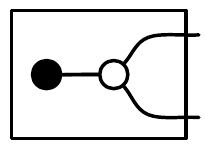}$}
\eqrule{\rev{A7}}
\lower6pt\hbox{$\includegraphics[height=.6cm]{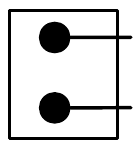}$}
$}
\end{align}\noindent\begin{align}
\notag
\lower6pt\hbox{$
\lower6pt\hbox{$\includegraphics[height=.6cm]{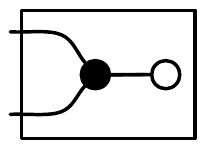}$}
\eqrule{\rev{A8}}
\lower6pt\hbox{$\includegraphics[height=.6cm]{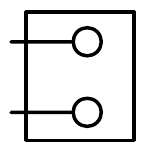}$}
$}
\end{align}
\begin{align}
\notag
\lower5pt\hbox{$\includegraphics[height=.6cm]{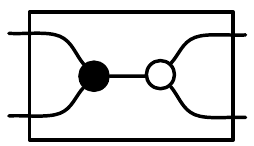}$}
\eqrule{\rev{A9}}
\lower9pt\hbox{$\includegraphics[height=.9cm]{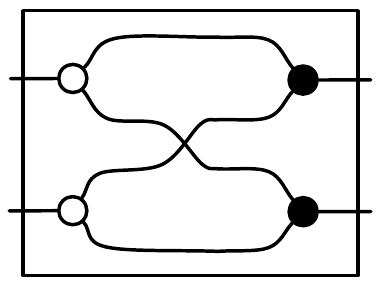}$}
\end{align}
\begin{align}
\notag
\lower7pt\hbox{$
\lower5pt\hbox{$\includegraphics[height=.5cm]{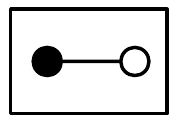}$}
\eqrule{\rev{A10}}
\lower5pt\hbox{$\includegraphics[height=.5cm]{graffles/idzerocircuit.pdf}$}
$}
\end{align}
\end{multicols}
\begin{multicols}{4}\noindent
\begin{align}
\notag
\lower8pt\hbox{$
\lower5pt\hbox{$\includegraphics[height=.5cm]{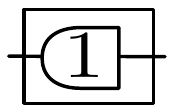}$}
\eqrule{\rev{A11}}
\lower4pt\hbox{$\includegraphics[height=.5cm]{graffles/idcircuit.pdf}$}
$}
\end{align}
\begin{align}
\notag
\lower8pt\hbox{$
\lower7pt\hbox{$\includegraphics[height=.6cm]{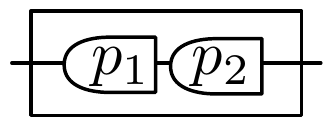}$}
\eqrule{\rev{A12}}
\lower7pt\hbox{$\includegraphics[height=.6cm]{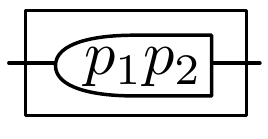}$}
$}
\end{align}
\begin{align}
\notag
\lower6pt\hbox{$
\lower8pt\hbox{$\includegraphics[height=.8cm]{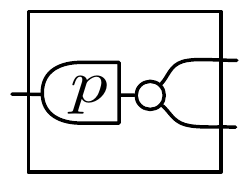}$}
\eqrule{\rev{A13}}
\lower8pt\hbox{$\includegraphics[height=.75cm]{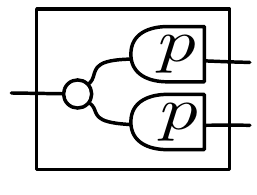}$}
$}
\end{align}
\begin{align}
\notag
\lower6pt\hbox{$
\lower6pt\hbox{$\includegraphics[height=.6cm]{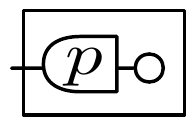}$}
\eqrule{\rev{A14}}
\lower6pt\hbox{$\includegraphics[height=.6cm]{graffles/Wcounit.pdf}$}
$}
\end{align}
\end{multicols}
\begin{multicols}{4}
\noindent
\begin{align}
\notag
\lower18pt\hbox{$\includegraphics[height=.75cm]{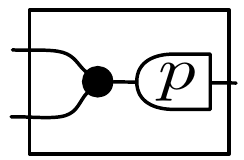}$}
\lower10pt\hbox{$\eqrule{\rev{A15}}$}
\lower19pt\hbox{$\includegraphics[height=.8cm]{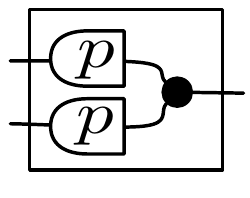}$}
\end{align}
\begin{align}
\notag
\lower12pt\hbox{$
\lower6pt\hbox{$\includegraphics[height=.6cm]{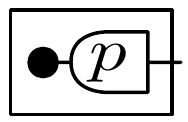}$}
\eqrule{\rev{A16}}
\lower6pt\hbox{$\includegraphics[height=.6cm]{graffles/Bunit.pdf}$}
$}
\end{align}
\begin{align}
\notag
\lower12pt\hbox{$
\lower4pt\hbox{$\includegraphics[height=.5cm]{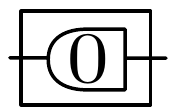}$}
\eqrule{\rev{A17}}\!\!
\lower4pt\hbox{$\includegraphics[height=.5cm]{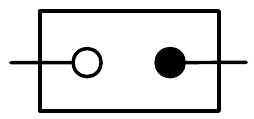}$}
$}
\end{align}
\begin{align}
\notag
\lower18pt\hbox{$\includegraphics[height=.9cm,width=1.3cm]{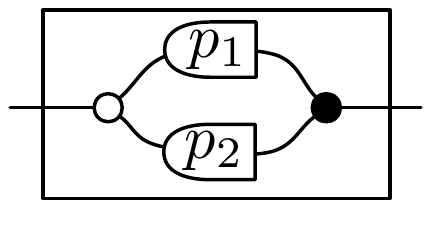}$}
\lower8pt\hbox{$\eqrule{\rev{A18}}$}
\lower14pt\hbox{$\includegraphics[height=.6cm,width=1.1cm]{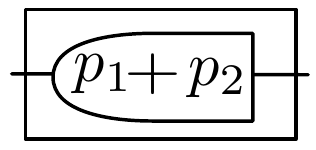}$}
\end{align}
\end{multicols}
\caption{Axioms of $\HA{\poly}^{\op}$, describing the interaction of the generators in \eqref{eq:SFcalculusSyntax2}.}\label{fig:axiomsHAop}
\end{figure}
\begin{figure}
%
%
%
 %
 %
 %
 \begin{multicols}{4}\noindent
\begin{align}
\notag
\lower5pt\hbox{$
\lower5pt\hbox{$\includegraphics[height=.55cm]{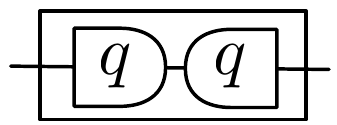}$}
\eqrule{I1}
\lower5pt\hbox{$\includegraphics[height=.55cm]{graffles/idcircuit.pdf}$}
$}
\end{align}
\begin{align}
\notag
\lower5pt\hbox{$
\lower5pt\hbox{$\includegraphics[height=.55cm]{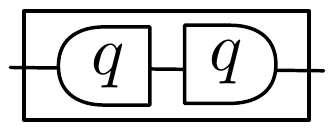}$}
\eqrule{I2}
\lower5pt\hbox{$\includegraphics[height=.55cm]{graffles/idcircuit.pdf}$}
$}
\end{align}
 \begin{align}
\notag
\raise6pt\hbox{$\includegraphics[height=1.1cm]{graffles/WFrobS.pdf}$}
\!\!
\raise19pt\hbox{$\eqrule{I3}$}
\!\!
\raise12pt\hbox{$\includegraphics[height=.8cm]{graffles/WFrobX.pdf}$}
\!\!
\raise19pt\hbox{$\eqrule{I3}$}
\!\!
\raise6pt\hbox{$\includegraphics[height=1.1cm]{graffles/WFrobZ.pdf}$}
\end{align}
\begin{align}
\notag
\includegraphics[height=1.1cm]{graffles/BFrobS.pdf}
\!\!
\raise14pt\hbox{$\eqrule{I4}$}
\!\!
\raise4pt\hbox{$\includegraphics[height=.9cm]{graffles/BFrobX.pdf}$}
\!\!
\raise14pt\hbox{$\eqrule{I4}$}
\!\!
\includegraphics[height=1.1cm]{graffles/BFrobZ.pdf}
\end{align}
\end{multicols}\noindent \vspace{-1cm}
 \begin{multicols}{4}\noindent
\begin{align}
\notag
\lower8pt\hbox{$\includegraphics[height=.8cm]{graffles/lccr.pdf}$}
\!\!
\eqrule{I5}
\lower6pt\hbox{$\includegraphics[height=.6cm]{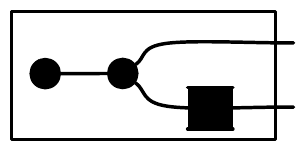}$}
\end{align}
\begin{align}
\notag
\lower8pt\hbox{$\includegraphics[height=.8cm]{graffles/rccl.pdf}$}
\!\!
\eqrule{I6}
\lower6pt\hbox{$\includegraphics[height=.6cm]{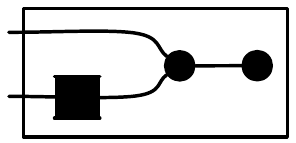}$}
\end{align}
\begin{align}
\notag
\lower5pt\hbox{$
\lower7pt\hbox{$\includegraphics[height=.7cm]{graffles/WSep.pdf}$} \eqrule{I7}
\lower6pt\hbox{$\includegraphics[height=.6cm]{graffles/idcircuit.pdf}$}
$}
\end{align}
\begin{align}
\notag
\lower5pt\hbox{$
\lower6pt\hbox{$\includegraphics[height=.6cm]{graffles/BSep.pdf}$} \eqrule{I8}
\lower5pt\hbox{$\includegraphics[height=.5cm]{graffles/idcircuit.pdf}$}
$}
\end{align}
\end{multicols}
\caption{Axioms of $\IBpoly$, describing the interaction of generators in \eqref{eq:SFcalculusSyntax1} with those in \eqref{eq:SFcalculusSyntax2}.}\label{fig:axiomsIH}
\end{figure}

Note that any string diagram of $\FC$ can be interpreted as one of $\HA{\poly}$: indeed, all the generators of $\FC$ as in~\eqref{eq:SFcalculusSyntax1} are also generators of $\HA{\poly}$. This gives us a PROP morphism $\FCtoHA \: \FC \to \HA{\poly}$, which can be composed with the iso $\HA{\poly} \tr{\cong} \Matpoly$ to obtain an interpretation $\gls{dsemHAO} \: \FC \to \Matpoly$ of circuits as polynomial matrices. Following Definition~\ref{def:sem}, we can present $\dsemHAO$ inductively as the PROP morphism mapping the generators in \eqref{eq:SFcalculusSyntax1} as follows:
\begin{eqnarray}\label{eq:defsemFCpart1}
    \begin{array}{rcl}
\Bcomult &\longmapsto& \tiny{\left(%
                \begin{array}{c}
                 \!\! 1 \!\!\\
                 \!\! 1\!\!
                \end{array}\right)} \\
\Bcounit &\longmapsto&  \initVect  \\
\end{array}
&
    \begin{array}{rcl}
\Wmult  &\longmapsto&   {\scriptsize\left(%
                \begin{array}{cc}
                \!\!\!  1 \! &\!\! 1 \!\!\!
                \end{array}\right)} \\
\Wunit &\longmapsto& \finVect \\
\end{array}
&
    \begin{array}{rcl}
\scalarT &\longmapsto& {\scriptsize\left(\begin{array}{c}
                 \!\!\! k\!\!\!
                \end{array}\right)} \\
\circuitXT &\longmapsto& {\scriptsize\left( \begin{array}{c}
                 \!\!\! x\!\!\!
                \end{array}\right)}%
\end{array}
\end{eqnarray}
\noindent where $\initVect \: 0 \to 1$ and $\finVect \: 1 \to 0$ are given by initiality and finality of $0$ in $\Matpoly$.

Proposition~\ref{prop:ab=vect} yields as a corollary \emph{soundness and completeness} of the semantics $\dsemHAO$.
\begin{corollary}\label{prop:soundcompleteFC}
For all circuits $c,d$ in $\FC$, $\dsemHA{c}=\dsemHA{d}$ iff $c \eqHA d$.
\end{corollary}

In Corollary~\ref{prop:soundcompleteFC}, we use notation $c \eqHA d$ to mean that $\FCtoHA(c) = \FCtoHA(d)$ in $\HA{\poly}$. Hereafter, we shall use the same convention in analogous situations when circuits are compared within a certain equational theory --- also, we do not bother to specify the subscript $\poly$ when the context is unambiguous.

\begin{remark}
There is a slight mismatch between the signatures of $\FC$ and $\HA{\poly}$: whereas $\FC$ only has a generators $\scalarT$ for each scalar $k \in \field$ and one $\circuitXT$ for the unknown $x$, $\HA{\poly}$ has a generator for all the polynomials of $\poly$. Nonetheless, we can still represent any polynomial $p = k_0 + k_1 x + k_2 x^2 + \dots + k_n x^n$ in $\FC$, as the circuit
$$\includegraphics[width=120pt]{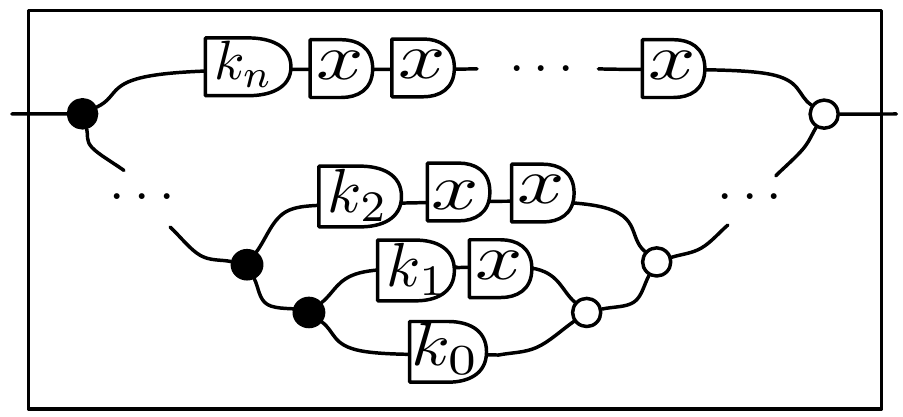}$$
which, under the interpretation $\FCtoHA \: \FC \to \HA{\poly}$, is equal to the circuit $\scalarp$, as expected. This observation also implies that $\FCtoHA$ is a full functor, whose action should be simply thought as quotienting the circuit syntax by the equations of $\poly$-Hopf algebras.
\end{remark}

Corollary~\ref{prop:soundcompleteFC} can be conveniently exploited also for circuits in $\FCop$. Indeed,  $\dsemHAO \: \FC \to \Matpoly$ induces the PROP morphism $\dsemHAO^{\m{op}}$ on the opposite categories, that we hereafter denote by  $\gls{dsemHAOop} \: \FCop \to \Matpoly^{\m{op}}$. The sound and complete axioms for this interpretation are those of $\HA{\poly}^{\op}$, reported in Figure~\ref{fig:axiomsHAop}.


\begin{corollary}\label{prop:soundcompleteFCop}
For all circuits $c,d$ in $\FCop$, $\dsemHAop{c}=\dsemHAop{d}$ iff $c \eqHAop d$.
\end{corollary}

\subsection{Polynomial Semantics of $\CD$}\label{sec:polysemCD}

The semantics of $\FC$ was given in terms of matrices, i.e., linear functions. This approach is coherent with the intuition, explained in Remark~\ref{remark:intuition}, that the signal in circuits of $\FC$ flows from left to right: left ports are inputs and right ports are outputs.

However, these traditional \emph{mores} fail in $\CD$---indeed, in this larger class only some circuits have a functional interpretation. For a counterexample, consider the circuit $\lccB \: 2 \to 0$. As observed in Example~\ref{exm:opsem}, the operational reading is that
$\lccB$ is a bent wire, whose behaviour is relational: its ports are neither inputs nor outputs in any traditional sense. This observation justifies the choice of the PROP $\SVpoly$ as the semantic domain for $\CD$. Here $\gls{frpoly}$ is the field of fractions of $\poly$ and arrows $n \to m$ in $\SVpoly$ are subspaces of $\frpoly^n \times \frpoly^m$, which we conveniently regards as linear \emph{relations} between $\frpoly^n$ and $\frpoly^m$ --- \emph{cf.} Convention~\ref{conv:subspaceLinRel}.

To give an interpretation $\CD \to \SVpoly$ we will adopt the same strategy as for the functional case. Whereas the axiomatisation for $\FC$ was given by the theory $\HA{\poly}$ of $\poly$-Hopf algebras, for $\CD = \FC + \FCop$ we shall use the theory $\IH{\poly}$ of \emph{interacting} $\poly$-Hopf algebras (Def.~\ref{def:IBR}). By definition, $\IH{\poly}$ features the equations of $\HA{\poly}$ (Fig.~\ref{fig:axiomsHA}), the equations of $\HA{\poly}^{\op}$ (Fig.~\ref{fig:axiomsHAop}) plus the equations describing the interaction between $\HA{\poly}$ and $\HA{\poly}^{\op}$ (Fig.~\ref{fig:axiomsIH}, where $q \neq 0$).

Now, note that all the generators of $\CD$ are also generators of $\IH{\poly}$: this yields a PROP morphism $\CDtoIH \: \CD \to \IH{\poly}$. Theorem~\ref{th:IBR=SVR} gives us an isomorphism $\sem{\IBR} \: \IH{\poly} \xrightarrow{\cong} \SVpoly$, which can be precomposed with $\CDtoIH$ to obtain the semantic map $\dsemO \: \CD \to \SVpoly$. Following the description of $\sem{\IBR}$ given by Definition \ref{def:semIBRInd}, we can inductively present $\dsemO$ as the PROP morphism mapping the generators in \eqref{eq:SFcalculusSyntax1} as
\begin{eqnarray*}
 \begin{array}{rcl}
\Bcomult & \!\!\longmapsto\!\! & [(%
                \matrixOnebis,\matrixOneOne)] \\
  \Bcounit & \!\!\longmapsto\!\! & [(\matrixOnebis,\matrixNull)]
 \end{array}
&
 \begin{array}{rcl}
\Wmult & \!\!\longmapsto\!\! & [(\matrixZeroOne,\matrixOnebis),(\matrixOneZero,\matrixOnebis)]\\
\Wunit & \!\!\longmapsto\!\! & \{(\matrixNull,\matrixZerobis)\}
 \end{array}
&
 \begin{array}{rcl}
\scalarT & \!\!\longmapsto\!\! & [(\matrixOnebis,k)] \\
\circuitXT & \!\!\longmapsto\!\! & [(\matrixOnebis,x)] \\
 \end{array}
\end{eqnarray*}
and defined symmetrically on the generators in \eqref{eq:SFcalculusSyntax2}
    \begin{eqnarray*}
 \begin{array}{rcl}
\Bmult & \!\!\longmapsto\!\! & [(%
                \matrixOneOne,\matrixOnebis)] \\
  \Bunit & \!\!\longmapsto\!\! & [(\matrixNull,\matrixOnebis)]
 \end{array}
&
 \begin{array}{rcl}
\Wcomult & \!\!\longmapsto\!\! & [(\matrixOnebis,\matrixZeroOne),(\matrixOnebis,\matrixOneZero)]\\
\Wcounit & \!\!\longmapsto\!\! & \{(\matrixZerobis,\matrixNull)\}
 \end{array}
&
 \begin{array}{rcl}
\scalaropT & \!\!\longmapsto\!\! & [(k,\matrixOnebis)] \\
\circuitXopT & \!\!\longmapsto\!\! & [(x,\matrixOnebis)]. \\
 \end{array}
\end{eqnarray*}
  where $[(\vv_1,\uu_1),\dots,(\vv_n,\uu_n)]$ indicates the space spanned by pairs of vectors $(\vv_1,\uu_1),\dots,(\vv_n,\uu_n)$.
Theorem~\ref{th:IBR=SVR} yields soundness and completeness of $\dsemO$.
\begin{corollary}\label{prop:isoIHsoundcomplete}
For all circuits $c,d$ in $\CD$, $\dsem{c}=\dsem{d}$ iff $c \eqIH d$.
\end{corollary}

It is useful for later reference to conclude with the following observation.
\begin{remark}\label{rmk:matrixformFractions} As explained in \S~\ref{sec:theorymatr}, there is a canonical way of representing any polynomial matrix $M \in \Matpoly[n,m]$ as a circuit $c \in \FC[n,m]$, which we called in \emph{matrix form} --- see Definition~\ref{def:matrixform}. By working in the wider class $\CD$, we can extend this representation to matrices over $\frpoly$.
Consider the following example:
\[
N = \left(
\begin{array}{ccc}
 \!\! p_{1}/q_{1} \!& p_{4}/q_{4} \!& p_{7}/q_{7}  \!\!\\
 \!\! p_{2}/q_{2} \!& p_{5}/q_{5} \!& p_{8}/q_{8}  \!\!\\
 \!\! p_{3}/q_{3} \!& p_{6}/q_{6} \!& p_{9}/q_{9}  \!\!
\end{array}\right)
\quad
\qquad \qquad
\quad
d = \lower50pt\hbox{$\includegraphics[width=130pt]{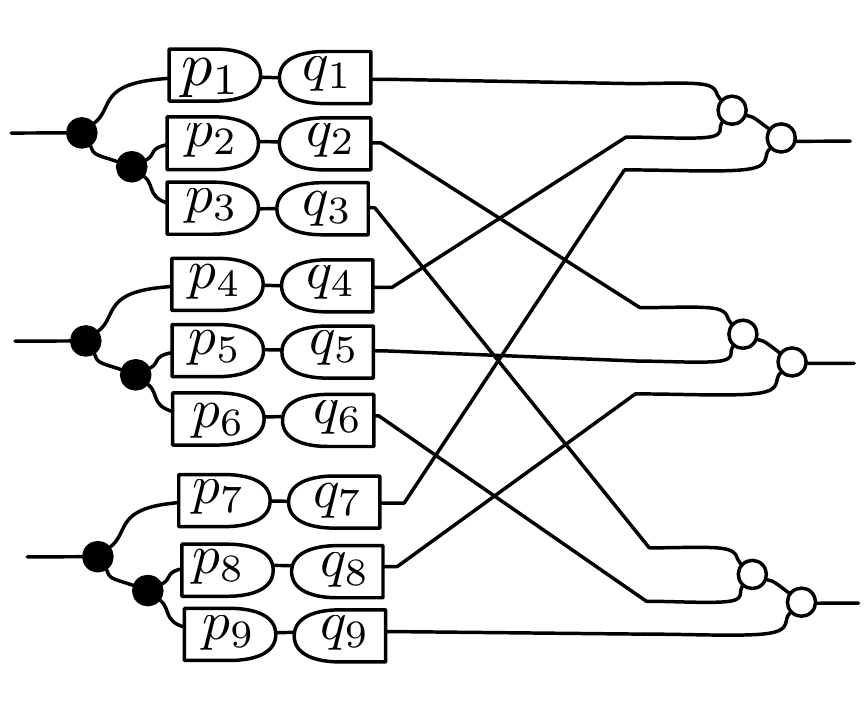}$}
\]
The circuit $d$ encodes $N$ in the following way: for each boundary of $d$, we assume a top-bottom enumerations of the ports, starting from $1$. Then $N_{ij} = \frac{p}{q}$ if and only if, reading the circuit from the left to the right, one finds a path connecting the $j^{th}$ port on the left to the $i^{th}$ port on the right passing through a circuit $\rationalcircuit$. Intuitively, the ports on the left represent columns, the ones on the right rows, and the links between them carry the values in the matrix. One can compute
$$\dsem{d} = \{(\vv,N \vv) \mid \vv \in \frpoly^3\} = [({\ee}_i, N{\ee}_i)]_{i\leq3}$$
where $\{{\ee}_i \ | \ i\leq 3\}$ is the standard basis of $\frpoly^3$.
\end{remark}

\vspace{-0.1cm}
\section{Denotational Semantics II: Streams}\label{sec:stream}
The polynomial semantics is characterised by very strong properties. For instance, not only it is sound and complete, but it is also full: any object in the semantic domain has a counterpart in the syntax. As we shall see, this richness turns out to be very convenient to prove theorems about the signal flow calculus. 
However, it does not quite capture the view of diagrams as signal processing circuits (Remark~\ref{remark:intuition}). In order to model infinite sequences of signals we need something more general than fractions of polynomials, namely \emph{streams}. This perspective motivates the extension of the polynomial semantics to a stream semantics, which we develop in this section.

\medskip

For this purpose, we first need to recall some useful notions.
A \emph{formal Laurent series} (fls) is a function $\sigma\colon \Z \to \field$ for which there exists $i\in \Z$ such that $\sigma(j)=0$ for all $j < i$. The \emph{degree} of $\sigma$ is the smallest $d\in \Z$ such that $\sigma(d)\neq 0$. We shall often write $\sigma$ as $\dots,\sigma(-1),\underline{\sigma(0)}, \sigma(1), \dots$ with position $0$ underlined, or as formal sum $\LSum{d}\sigma(i)x^i$. Using the latter notation, the sum and product of two fls $\sigma=\LSum{d}\sigma(i)x^i$ and $\tau=\LSum{e}\tau(i)x^i$ are given by:
\begin{equation}\label{eq:conv}
\sigma+\tau=\LSum{min(d,e)}\big(\sigma(i)+\tau(i)\big)x^i \qquad \sigma \cdot \tau = \LSum{d+e}\Big( \sum_{k+j=i} \sigma(j)\cdot \tau(k) \Big)x^i \end{equation}
The units for $+$ and $\cdot$ are $\dots 0,\underline{0},0\dots$ and $\dots 0,\underline{1},0\dots$.
Fls form a field $\gls{laur}$, where the inverse $\sigma^{-1}$ of the fls $\sigma$ with degree $d$ is given as follows.
\begin{equation}\label{eq:inverse}
\sigma^{-1}(i) = \,
\begin{cases} 0 & \text{ if } i < -d \\
\sigma(d)^{-1}  &\text{ if } i=-d \\
\frac{\sum_{i=1}^{n} \big( \sigma(d+i) \cdot \sigma^{-1}(-d+n-i)\big)}{-\sigma(d)} & \text{ if } i=-s+n \text{ for } n>0
\end{cases}
\end{equation}
A \emph{formal power series} (fps) is a fls with degree $d\geq 0$. By \eqref{eq:conv}, fps are closed under $+$ and $\cdot$, but not under inverse: it is immediate by \eqref{eq:inverse} that $\sigma^{-1}$ is a fps iff $\sigma$ has degree $d=0$. Therefore fps form a ring --- actually, a principal ideal domain --- which we denote by $\gls{fps}$.

We shall refer to both fps and fls as \emph{streams}. Indeed, fls are sequences with an infinite future, but a finite past.
Analogously to how a polynomial $p$ can be seen as a fraction $\frac{p}{1}$,
an fps $\sigma$ can be interpreted as the fls $\dots ,0,\underline{\sigma(0)}, \sigma(1), \sigma(2), \dots$. A polynomial $p_0+p_1x+\dots +p_nx^n$ can also be regarded as the fps $\LSum{0}p_ix^i$ with $p_i=0$ for all $i>n$. Similarly, polynomial
fractions can be regarded as fls:
we define $\embfrpolyfls\colon \frpoly \to \laur$ as the unique field morphism mapping $k\in \field$
to  $\dots0,\underline{k},0\dots$ and the indeterminate $x$ to $\dots ,0,\underline{0},1,0, \dots$

Differently from polynomials, fractions can denote streams with possibly infinitely many non-zero values. For instance, \eqref{eq:conv} and \eqref{eq:inverse} imply that $\frac{x}{1-x-x^2}$ is the Fibonacci series $\dots,0,\underline{0},1,1,2,3,\dots$. Moreover, while polynomials can be interpreted as fps, fractions need the full generality of fls: $\frac{1}{x}$ denotes $\dots0,0,1,\underline{0},0,\dots$ These translations are ring homomorphisms and are illustrated by the commutative diagram below.

\begin{equation}\label{eq:ringhoms}
\vcenter{
\renewcommand{\labelstyle}{\textstyle}
\xymatrix@R=1pt@C=5pt{
\fps \ar@{^{(}->}[rrrr] &  & & & \laur\\
\\
\\
& \ratio \ar@{_{(}->}[luuu] \ar@{^{(}->}[rrrd]
\\
\poly \ar@{^{(}->}[ur] \ar@{_{(}->}[rrrr] \ar@{^{(}->}[uuuu]^{\embpolyfps} & & & & \frpoly \ar@{_{(}->}[uuuu]_{\embfrpolyfls}
}
}
\end{equation}
At the center, $\gls{ratio}$ is the ring of \emph{rationals}, i.e., fractions  $\frac{k_0+k_1x+k_2x^2 \dots + k_nx^n}{l_0+l_1x+l_2x^2 \dots + l_nx^n}$ where $l_0\neq 0$. Differently from fractions, rationals denote only fps---in other words---bona fide streams that do not start ``in the past''. Indeed, since $l_0\neq 0$, the inverse of $l_0+l_1x+l_2x^2 \dots + l_nx^n$ is, by \eqref{eq:inverse}, a fps.
The streams denoted by $\ratio$ are known in literature as \emph{rational streams}~\cite{Berstel}.

Hereafter, we shall often use polynomials and fractions to denote the corresponding streams. Also, we will use PROPs to reason about these algebras:
 \begin{itemize}[noitemsep,topsep=0pt,parsep=0pt,partopsep=0pt]
 \item the PROP $\Matfps$ of $\fps$-matrices --- see Definition~\ref{def:HA};
 \item the PROP $\Matratio$ of $\ratio$-matrices --- see Definition~\ref{def:HA};
 \item the PROP $\SVfps$ of linear relations over $\laur$ --- see Definition~\ref{def:sv}.
 \end{itemize}

\subsection{Stream Semantics of $\FC$}\label{sec:streamHA}
The polynomial semantics $\dsemHAO\colon \FC \to \Matpoly$ gives us already a mean to regard the circuits in $\FC$ as \emph{stream} transformers.
Indeed, the interpretation $\embpolyfps \: \poly \to \fps$ of a polynomial as an fps --- see~\eqref{eq:ringhoms} --- can be extended pointwise to a faithful PROP morphism $\Matpoly \to \Matfps$, for which we conveniently use the same notation $\embpolyfps$. We can then define the \emph{stream semantics} $\gls{strsemHAO}$ of $\FC$ as the composite $\dsemHAO \poi \embpolyfps \: \FC \to \Matfps$ --- it maps a circuit $c\in \FC[n,m]$ to an $m \times n$ matrix $\strsemHA{c}$ over $\fps$.

\begin{remark} The matrix $\strsemHA{c}$ can be described as a stream transformer with the same kind of behaviour expressed by the operational rules of Fig.~\ref{fig:operationalSemantics} --- \emph{cf}. Remark~\ref{remark:intuition}. 
For instance, $\strsemHA{\circuitX}= \left(  x \right)$ maps every stream $\sigma\in \fps$ into the stream $\sigma \cdot x$ which, by \eqref{eq:conv}, is $\underline{0},\sigma(0), \sigma(1), \sigma(2), \dots$ Thus $\circuitX$ behaves as a \emph{delay}. For $k \in \field$, $\strsemHA{\scalar}= \left(  k \right)$ maps $\sigma$ to $\sigma \cdot k = \underline{k\sigma(0)}, k\sigma(1), k\sigma(2), \dots$ Therefore $\scalar$ is an \emph{amplifier}. For the remaining operations: $\strsemHA{\Wmult} = {\small \left(\begin{array}{cc}\!\!\! 1\!\!\!&\!\!\! 1 \!\!\!\end{array}\right)}$ maps $\matrixTwoROneC{\sigma}{\tau}$ to $\sigma + \tau$, thus $\Wmult$ behaves as an \emph{adder}. Its unit $\Wunit$ emits the constant stream $\underline{0},0,0 \dots$. Finally, $\Bcomult$ behaves as a \emph{copier}, because $\strsemHA{\Bcomult} = \matrixTwoROneC{1}{1}$ maps $\sigma$ to $\matrixTwoROneC{\sigma}{\sigma}$, and its counit $\Bcounit$ is the transformer taking any stream as input and producing no output. We shall return on the comparison between stream semantics and operational behaviour more thoroughly in \S~\ref{sec:fullabstract}.
\end{remark}

The interpretation $\strsemHA{\cdot}$ coincides with the semantics developed by J.~Rutten in \cite[\S~4.1]{DBLP:journals/tcs/Rutten05}. Our approach has the advantage of making the circuits representation formal and allowing for equational reasoning, as shown for instance in Example~\ref{ex:partoseq} below.
Indeed, since $\embpolyfps\colon \Matpoly \to \Matfps$ is faithful, the axiomatization of $\ABpoly$ is sound and complete also for $\strsemHAO$.

\begin{example}\label{ex:partoseq}Consider the following derivation in the equational theory of $\ABpoly$, where \eqref{eq:scalarbcomult} is used at each step.
\begin{equation*}
\lower16pt\hbox{$\includegraphics[height=1.3cm,width=11.2cm]{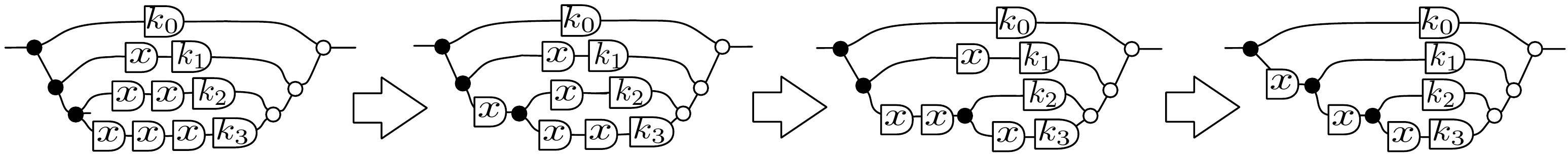}$}
\end{equation*}
Any of the circuits above has stream semantics given by the matrix $\left(p\right) \in \Matfps[1,1]$, where $p = \underline{k_0}, k_1, k_2, k_3, 0 \dots$.
Along the lines of \cite[Prop. 4.12]{DBLP:journals/tcs/Rutten05}, one can think of the derivation above as a procedure that reduces the total number of delays $\circuitX$ appearing in the implementation of  the stream function $f \: \sigma \mapsto \sigma \cdot p$.
\end{example}

\subsection{Stream semantics of $\CD$}\label{ss:streamIB}

In \S\ref{sec:polysemCD} we gave a semantics to $\CD$ in terms of linear relations over fractions of polynomials. We now extend this semantics to linear relations over streams. While formal power series are enough to provide a stream semantics to $\FC$, for the whole of $\CD$ one needs the full generality of Laurent series since, as we have discussed above, not all fractions of polynomials (e.g. $\frac{1}{x}$) denote fps.

 The \emph{stream semantics} of $\CD$ is the unique PROP morphism $\gls{strsemO} \colon \CD \to \SVfps$ mapping the generators in~\eqref{eq:SFcalculusSyntax1} as follows:
\begin{eqnarray*}\noindent
\begin{array}{rcl}\Bcomult &\mapsto& \{ %
				( \sigma ,{\scriptsize\left(%
				\begin{array}{c}
				  \!\!\sigma\!\! \\
				  \!\!\sigma\!\!
				\end{array}\right)} ) \ | \ \sigma \in \laur \}\\
\Bcounit & \mapsto & \{ (\sigma, \matrixNull ) \ |\ \sigma\in\laur \}\\
\scalarT & \mapsto & \{ (\sigma, k \cdot \sigma)\ |\ \sigma\in\laur \}

\end{array}%
&
\begin{array}{rcl}
\Wmult & \mapsto & \{( {\scriptsize\left(%
				\begin{array}{c}
				  \!\!\sigma\!\! \\
				  \!\!\tau\!\!
				\end{array}\right)}, \sigma+\tau) \ | \ \sigma,\tau \in \laur \}\\
\Wunit & \mapsto & \{( \matrixNull ,0) \}\\
\circuitXT & \mapsto & \{ (\sigma, x \cdot \sigma)\ |\ \sigma\in\laur \}
\end{array}
\end{eqnarray*}
and symmetrically for the components in~\eqref{eq:SFcalculusSyntax2}:
\begin{eqnarray*}\noindent
\begin{array}{rcl}\Bmult &\mapsto& \{ %
				({\scriptsize\left(%
				\begin{array}{c}
				  \!\!\sigma\!\! \\
				  \!\!\sigma\!\!
				\end{array}\right)},\sigma ) \ | \ \sigma \in \laur \}\\
\Bunit & \mapsto & \{ (\matrixNull,\sigma ) \ |\ \sigma\in\laur \}\\
\scalaropT & \mapsto & \{ (k \cdot \sigma,\sigma)\ |\ \sigma\in\laur \}

\end{array}%
&
\begin{array}{rcl}
\Wcomult & \mapsto & \{(\sigma+\tau, {\scriptsize\left(%
				\begin{array}{c}
				  \!\!\sigma\!\! \\
				  \!\!\tau\!\!
				\end{array}\right)}) \ | \ \sigma,\tau \in \laur \}\\
\Wcounit & \mapsto & \{(0,\matrixNull) \}\\
\circuitXopT & \mapsto & \{ (x \cdot \sigma,\sigma)\ |\ \sigma\in\laur \}
\end{array}
\end{eqnarray*}
where $0$, $x$ and $k$ here denote streams. The reader may notice a close resemblance between the definition of $\strsemO{}$ and of $\dsemO{}$ --- we shall see later that indeed $\strsemO{}$ factorises through $\dsemO{}$ (Proposition~\ref{lemma:streamSemUnivProperty}).

\begin{example}\label{exm:cospan} Consider the circuit $\circuitXcospan$. We have that
\begin{eqnarray*}	
	\strsem{\circuitXcospan} &\!\!\!\!=\!\!\!\!&
  \strsem{\circuitXT} \poi \strsem{\circuitXopT} \\
&\!\!\!\!=\!\!\!\!& \{ (\sigma , \sigma \cdot x) \mid \sigma \in \laur \} \poi \{ (\sigma \cdot x , \sigma ) \mid \sigma \in \laur \} \\
&\!\!\!\!=\!\!\!\!& \{ (\sigma , \sigma ) \mid \sigma \in \laur \}
\end{eqnarray*}
which is equal to $\strsem{\IdnetT}$.
Note that any fls $\dots \sigma(-1),\underline{\sigma(0)},\sigma(1)\dots$ on the left of $\circuitXT$ is related to $\dots \underline{\sigma(-1)},\sigma(0),\sigma(1)\dots$ on its right and this is in turn related to $\dots \sigma(-1),\underline{\sigma(0)},\sigma(1)\dots$ on the right of $\circuitXop$. The circuit $\circuitXop$ is thus the inverse of $\circuitXT$: while $\circuitXT$ \emph{delays} $\sigma$, $\circuitXopT$ \emph{accelerates} it.

Similarly, consider a circuit $\circuitkkop$. Its semantics is the composite of $\strsem{\scalarT}$ (pairing $\sigma$ with $\sigma \cdot k$) and $\strsem{\scalaropT}$ (pairing $k \cdot \sigma$ with $\sigma$): if $k \neq 0$, we can see it as first multiplying and then dividing $\sigma$ by $k$. Thus for $k \neq 0$ $\circuitkkop$ and $\IdnetT$ have the same denotation.
\end{example}

\begin{example}
In Example \ref{exm:opsem}, we presented the circuit $c_2$ as the composition of four sequential chunks. Their stream semantics is displayed below.
$$\begin{array}{rcl}
\strsem{(\BunitT\poi \BcomultT) \tns \IdnetT}{} & = & \{ (\sigma_1, {\scriptsize\left(%
				\begin{array}{c}
				  \!\!\tau_1\!\! \\
				  \!\!\tau_1\!\! \\
				  \!\!\sigma_1\!\!
				\end{array}\right)}) \ |\ \sigma_1,\tau_1 \in  \laur \}
\\
\strsem{\IdnetT \tns (\WmultT \poi \BcomultT)}{ } & = & \{( {\scriptsize\left(%
				\begin{array}{c}
				  \!\!\tau_2\!\! \\
				  \!\!\sigma_2\!\! \\
				  \!\!\rho_2\!\!
				\end{array}\right)},
				{\scriptsize\left(%
				\begin{array}{c}
				  \!\!\tau_2\!\! \\
				  \!\!\sigma_2+\rho_2\!\! \\
				  \!\!\sigma_2+\rho_2\!\!
				\end{array}\right)}) \ |\ \sigma_2,\tau_2,\rho_2 \in  \laur  \}
\\
\strsem{ (\IdnetT \tns \circuitXT) \tns \IdnetT} { } & = & \{( {\scriptsize\left(%
				\begin{array}{c}
				  \!\!\tau_3\!\! \\
				  \!\!\sigma_3\!\! \\
				  \!\!\rho_3\!\!
				\end{array}\right)},
				{\scriptsize\left(%
				\begin{array}{c}
				  \!\!\tau_3\!\! \\
				  \!\!x \cdot \sigma_3 \!\! \\
				  \!\!\rho_3\!\!
				\end{array}\right)}) \ |\ \sigma_3,\tau_3,\rho_3 \in  \laur  \}
\\				
\strsem{ (\BmultT \poi \BcounitT) \tns \IdnetT}{} & = &  \{ ({\scriptsize\left(%
				\begin{array}{c}
				  \!\!\tau_4\!\! \\
				  \!\!\tau_4\!\! \\
				  \!\!\sigma_4\!\!
				\end{array}\right)}, \sigma_4) \ |\ \sigma_4,\tau_4 \in  \laur \}
\end{array}$$
The composition in $\SVfps$ of the four linear relations above is
$$\{(\sigma_1, \sigma_4)  \ |\ \text{there exist } {\scriptsize \sigma_2,\sigma_3,\tau_1,\dots,\tau_4, \rho_2,\rho_3} \text{ s.t. }
{\scriptsize \begin{cases}
\tau_1=\tau_2=\sigma_2 =\tau_3=\tau_4,   \\
\sigma_2+\rho_2=\sigma_3, \; x \cdot \sigma_3 = \tau_4 \\
\sigma_1 = \rho_2, \; \sigma_2+\rho_2 =\rho_3 = \sigma_4
\end{cases}}\}
$$
By simple algebraic manipulations one can check that the above systems of equations has a unique solution given by $\sigma_4=\frac{1}{1-x}\sigma_1$.		
Since $\strsemO{}$ is a PROP morphism and $c_2$ is the composition of the four chunks above, we obtain  				
$$\strsem{c_2}{ } = \{(\sigma_1, \frac{1}{1-x} \cdot \sigma_1)  \ |\ \sigma_1 \in \laur\}.$$
This relation contains all pairs of streams that can occur on the left and on the right ports of $c_2$. For instance if $\underline{1},0,0 \dots$ is on the left,  $\underline{1},1,1 \dots$ is on the right.

For the other circuit of Example \ref{exm:opsem}, namely $c_1$, it is immediate to see that
$$\strsem{c_1}{ } = \{((1-x)\cdot \sigma_1, \sigma_1)  \ |\ \sigma_1 \in \laur\}$$
which is clearly the same subspace as $\strsem{c_2}{ }$. In Example \ref{ex:twoimplementations}, we will prove the semantic equivalence of the two circuits  by means of the equational theory of $\IBpoly$.  This is always possible since, as stated by the following theorem, the axiomatization of $\IBpoly$ is sound and complete with respect to $\strsemO{}$.
\end{example}



\begin{theorem}[Soundness and completeness] \label{cor:completeness}
 For all $c,d$ in $\CD$, $c\eqIH d$ iff $\strsem{c} = \strsem{d}$.
\end{theorem}
\noindent Theorem~\ref{cor:completeness} can be proved elegantly via a cube construction, similarly to the ones used for characterising (partial) equivalence relations (\S~\ref{sec:pushout}) and linear relations (\S~\ref{sec:cubebottom}). The diagram of interest in $\PROP$ is the following:
\begin{equation}\label{eq:cube2}
\hspace{-.95cm}{
\tag{${\text{\mancube}\above 0pt \text{\mancube}}$}
\vcenter{$
\qquad \qquad
\xymatrix@C=15pt@R=15pt{
& {\Matpoly + \Matpolyop} \ar@{^{(}->}[dd]_(.3){ }|{\hole}
\ar[dl]_(.7){[\MatToCospan,\MatopToCospan]} \ar[rr]^{[\MatToSpan,\MatopToSpan]} & & {\SpanMatpoly} \ar[dl]|{\Phi} \ar@{^{(}->}[dd]^{\Theta} \\
{\CospanMatpoly} \ar[rr]|{\Psi} \ar@{^{(}->}[dd]_{\Upsilon}  & & {\SVpoly} \ar@{.>}[dd]^(0.65){[\embfrpolyfls]} \\
& {\Matfps+ \Matfpsop} \ar[dl]_(0.6){[\MatToCospan',\MatopToCospan']} \ar[rr]|(.49){\hole} \ar@{}[rr]^(0.37){[\MatToSpan',\MatopToSpan']} & & {\SpanMatfps} \ar[dl]|{\Phi'} \\
{\CospanMatfps} \ar[rr]|{\Psi'} & & {\SVfps}
}$}
}
\end{equation}
Note that the top face is in fact the bottom face of the cube~\eqref{eq:cube} for linear relations, instantiated to the case $\PID = \poly$. Thus constructing~\eqref{eq:cube2} can be seen as adding another  ``floor'' below~\eqref{eq:cube}. The argument for Theorem~\ref{cor:completeness} will essentially rely on showing that $\strsemO$ is the composite of the polynomial semantics $\dsemO$ with the universal map $[\embfrpolyfls]$  in \eqref{eq:cube2}. 
%
To this aim, we first distill the components of \eqref{eq:cube2}.
\paragraph{Top face} Since the top face  of \eqref{eq:cube2} is the bottom face of \eqref{eq:cube}, we simply need to instantiate the various definitions to the case $\PID = \poly$.
The map $[\MatToSpan,\MatopToSpan]$ arises from:
\[
\begin{array}{cc}
\MatToSpan \colon \Matpoly \to \SpanMatpoly & \MatopToSpan \colon \Matpolyop \to \SpanMatpoly \\
A \: n\to m \longmapsto (n \tl{\id} n \tr{A} m) & B \: n\to m \longmapsto (n \tl{B} m \tr{\id} m)
\end{array}
\] and, similarly, $[\MatToCospan,\MatopToCospan]$ is the pairing of
\[
\begin{array}{cc}
\MatToCospan \colon \Matpoly \to \CospanMatpoly & \MatopToCospan \colon \Matpolyop \to \CospanMatpoly \\
A \: n\to m \longmapsto (n \tr{A} m \tl{\id} m) & B \: n\to m \longmapsto (n \tr{\id} n \tl{B} m)\text{.}
\end{array}
\]

%
The morphism $\Phi$ maps $n\tl{V} z \tr{W} m$ to the linear relation
\[
\{\,(\uu,\vv)\ |\ \uu\in \frpoly^n,\, \vv\in \frpoly^m,\, \exists \ww\in \frpoly^z.\; \embpolyfrpoly(V)\ww=\uu \wedge \embpolyfrpoly(W)\ww = \vv \,\}\]
where $\delta \: \Matpoly \to \Matfrpoly$ is the obvious embedding, and $\Psi$ acts as follows:
\[
{n \tr{V} z \tl{W} m} \quad \longmapsto \quad
\{\, (\uu,\vv) \ | \ \uu\in \frpoly^n,\, \vv\in \frpoly^m,\ \embpolyfrpoly(V)\uu = \embpolyfrpoly(W)\vv\, \}\text{.}
\]
Theorem~\ref{th:bottomfacePushout} ensures that these maps make the top face a pushout diagram in $\PROP$.

\paragraph{Bottom face} The morphisms of the bottom face, $[\MatToSpan',\MatopToSpan']$, $[\MatToCospan',\MatopToCospan']$, $\Phi'$ and $\Psi'$, are defined analogously.
%
Since $\fps$ is a PID and $\laur$ is its field of fractions, by Theorem~\ref{th:bottomfacePushout}, the bottom face is also a pushout in $\PROP$.

\paragraph{Vertical edges}
The rear morphism follows from the embedding $\embpolyfps \colon \Matpoly \to \Matfps$ described in \S\ref{sec:streamHA}. $\Theta$ maps a span $n\tl{V} z \tr{W} m$ to $n\tl{\wstream{V}} z \tr{\wstream{W}} m$. To verify that this is a morphism of PROPs, one needs to check that it preserves composition in $\Matpoly$:
\begin{lemma}\label{lemma:pullbackpreservation}
 $\embpolyfps \colon \Matpoly \to \Matfps$ preserves pullbacks.
\end{lemma}
\begin{proof} See Appendix~\ref{app:appendixStream}.\end{proof}
Similarly, the leftmost morphism $\Upsilon$ maps $n\tr{V} z \tl{W} m$ to $n\tr{\wstream{V}} z \tl{\wstream{W}} m$. Since $\Matpoly$ and $\Matfps$ are both self-dual (see Rmk.~\ref{rmk:graphicaltranspose}), it follows by Lemma~\ref{lemma:pullbackpreservation} that $\embpolyfps$ also preserves pushouts and, therefore, $\Upsilon$ is a morphism of PROPs.

By definition, the rear faces commute. As a consequence, there exists $[\embfrpolyfls]\colon \SVpoly \to \SVfps$ given by the universal property of the top face of \eqref{eq:cube2}. 
To give a concrete description of $[\embfrpolyfls]$, observe that $\embfrpolyfls \colon \frpoly \to \laur$ in \eqref{eq:ringhoms} can be pointwise extended to matrices and sets of vectors. For a subspace $H$ in $\SVpoly$, let $[\lrn{H}]$ be the space in $ \SVfps$ generated by the set of vectors $\lrn{H}$. 
\begin{lemma}\label{lemma:concreteDescrUnivPropFLS}
 The morphism $[\embfrpolyfls]\colon \SVpoly \to \SVfps$ maps $H$ in $\SVpoly$ to $[\wlrn{H}]$. 
\end{lemma}
\begin{proof} See Appendix~\ref{app:appendixStream}.\end{proof}

\begin{proposition}\label{lemma:streamSemUnivProperty}
$\strsemO = \dsemO\poi [\embfrpolyfls]$.
\end{proposition}
\begin{proof} Clear from definitions of $\strsemO$ and $\dsemO$, and Lemma~\ref{lemma:concreteDescrUnivPropFLS}. \end{proof}

By construction, the morphism $\dsemO\poi [\embfrpolyfls]$ has the desired properties allowing to infer soundness and completeness of $\IBpoly$ with respect to the stream semantics.
\begin{proof}[Proof of Theorem~\ref{cor:completeness}] Let $c$ and $d$ be in $\CD$. By Proposition~\ref{lemma:streamSemUnivProperty}, $\strsem{c} = \strsem{d}$ iff $[\lrnwide{\dsem{c}}] = [\lrnwide{\dsem{d}}]$.
Now, $[\embfrpolyfls]$ is given by the universal property in \eqref{eq:cube2}: since the other vertical edges of \eqref{eq:cube2} are faithful, also $[\embfrpolyfls]$ is faithful. Then $[\lrnwide{\dsem{c}}] = [\lrnwide{\dsem{d}}]$ iff $\dsem{c} = \dsem{d}$ and, therefore, iff $c \eqIH d$ by Corollary~\ref{prop:isoIHsoundcomplete}.
\end{proof}

\section{A Kleene's Theorem for Signal Flow Diagrams}\label{sec:SFG}

The denotational semantics for $\CD$, developed in the previous sections, clearly also gives a denotation for circuits in its sub-PROP $\SFGform$, in terms of \emph{relations} between streams. However, as outlined in \S\ref{sec:SFcalculus}, we expect that signal flow graphs express \emph{functional} behaviours. In this section we shall show that this is the case: our main result is that circuits in $\SFGform$, up-to equality in $\IBpoly$, characterise precisely the \emph{rational} behaviours in $\SVpoly$, i.e., the functional subspaces given by $\ratio$-matrices.

Kleene's theorem~\cite{Kleene56} is the statement that any finite automaton can be characterised by a regular expression and that, conversely, any regular expression can be realised by such automaton. By analogy, one can call Kleene's theorem the syntactic characterisation of rational behaviours for various models of computations: recent results of this kind have been given for instance for Mealy machines~\cite{Silva_KleeneQuantitative}, weighted automata and probabilistic systems~\cite{Silva_KleeneMealy}. Typically, the challenging task in this kind of approach is to provide a \emph{complete axiomatisation} for the equivalence of syntactic expressions: for regular languages, such axiomatisation has been shown only in 1990 by Kozen~\cite{Kozen94acompleteness}.


In the spirit of Kleene's result, we characterise the class of rational behaviours --- the functional subspaces generated by $\ratio$-matrices --- by means of a syntax --- $\SFGform$. In our case, the quest for a complete set of axioms does not require any additional effort: the equations of $\IBpoly$ suffice to prove semantic equivalence of $\SFGform$-circuits.

We remark that the correspondence between signal-flow diagrams and rational matrices is a folklore result in control theory (see e.g.~\cite{DBLP:journals/tcs/Rutten05,Lahti}): the novelty of our approach is in giving a rigorous (string diagrammatic) syntax --- where notions of ``input'', ``output'' and direction of flow are entirely derivative --- and in providing a complete axiomatisation. Interestingly, \cite{Milius_streamaxiom} also presents a Kleene's theorem for $\ratio$-matrices, but the proposed (non-diagrammatic) syntax and the complete set of axioms are of a rather different flavour and only capture $1 \times 1$ matrices. Also, we emphasize that our approach, based on the generators of $\CD$, does not feature any primitive for recursion, which appear instead in \cite{Milius_streamaxiom} as well as in the aforementioned works~\cite{Kleene56,Silva_KleeneQuantitative,Silva_KleeneMealy}.

The following is one direction of the correspondence between $\SFGform$ and $\ratio$-matrices.
\begin{proposition}\label{prop:circuittomatrix}
Suppose that $c\in \SFGform[n,m]$. Then $\dsem{c}$ is the subspace $[({\ee}_i, A{\ee}_i)]_{i\leq n}$ for some $A\in \Matratio [n,m]$, where $\{{\ee}_i \ | \ i\leq n\}$ is the standard basis of $\frpoly^n$.
\end{proposition}
\begin{proof}See Appendix~\ref{appendix:SFG}.\end{proof}

Note that the converse does not hold, in the sense that there are functional subspaces given by rational matrices which are in the image of circuits \emph{not} in $\SFGform$. In order to strengthen our correspondence to an isomorphism with $\Matratio$, we are going to show that all such circuits are provably equivalent in $\IBpoly$ to one in $\SFGform$. The following example illustrates an instance of our general result.

 \begin{example}\label{ex_fibonacci} The rational $\frac{x}{1-x-x^2}$ denoting the Fibonacci sequence can be succinctly represented as the circuit $\circuitX\poi\circuitfibr$, which is not in $\SFGform$. Indeed, composing
  $\dsem{\circuitX} = [(1,x)]$ with $\dsem{\circuitfibr} =[(1-x-x^2,1)]$ yields the $\frpoly$-subspace $[(1,\frac{x}{1-x-x^2})]$.
In terms of streams,  $\strsem{\circuitX\poi\circuitfibr}{}$ is the\ $\laur$-subspace $[(
\quad \underline{1}, 0, 0, \dots \, \quad , \, \quad \underline{0}, 1, 1, 2, 3, 5, \dots \quad)]$.

 The derivation in the equational theory of $\IBpoly$ below shows how we can ``implement'' the Fibonacci specification $\circuitX\poi\circuitfibr$, by transforming it into a string diagram of $\SFGform$. Since $\antipode \eqIH \antipodeop$, hereafter we use notation $\antipodesquare$ for both these circuits, analogously to what we did in Chapter~\ref{chapter:hopf}.
 \begin{equation*}
\includegraphics[height=2.5cm]{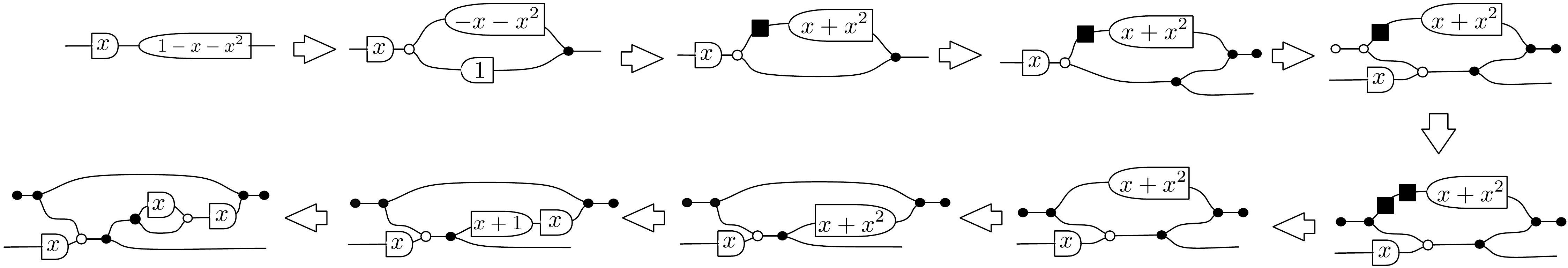}
\end{equation*}
In the above derivation, the strategy is to unfold $\circuitfibr$ (with \eqref{eq:unitscalar}$^{op}$, \eqref{eq:scalarmult}$^{op}$ and \eqref{eq:scalarsum}$^{op}$) and use the Frobenius axioms~\eqref{eq:WFrob}-\eqref{eq:BFrob} to deform the circuit to obtain the feedback loop. Then the sub-circuit representing $x+x^2$ is moved along the bent wire $\lccB$ using \eqref{eq:ccsliding}. 

At the end of \S\ref{sec:types}, we will explain formally in which sense the final circuit of the derivation can be thought as the implementation of the first one. At an intuitive level, this can be explained in terms of flows: in the first circuit it is not possible to assign a direction to the flow, while in the last one signal flows from left to right.
Using the operational rules of Fig.~\ref{fig:operationalSemantics}, the reader can verify that inputting the stream $\underline{1}, 0, 0, \dots$ on the left of the final circuit yields the Fibonacci sequence
 $\underline{0}, 1, 1, 2, 3, 5, \dots$
as output on the right.
 \end{example}

In view of Example~\eqref{ex_fibonacci}, we shall work with $\SFGform$ modulo $\IBpoly$. Since morphisms of PROPs are identity-on-objects, we can simply take the image of $\SFGform$ in $\IBpoly$.
\begin{definition}\label{defn:sfgImage} $\SFG$ is the sub-PROP of $\IBpoly$ given by the image of $\SFGform \to \CD \to \IBpoly$.
\end{definition}
One can think of $\SFG$ as the PROP whose arrows are the circuits of $\CD$ that are equivalent in $\IBpoly$ to one of $\SFGform$, subject to the equations of $\IBpoly$. 
We can now state the main theorem of this section.

\begin{theorem}\label{th:SFGcharactRationals}
$\SFG \cong \Matratio$.
\end{theorem}

The direction from circuits to matrices of Theorem~\ref{th:SFGcharactRationals} is already given by Proposition~\ref{prop:circuittomatrix}. The following statement takes care of the converse.

\begin{proposition}\label{prop:matrixtocircuit}
Suppose that $A\in \Matratio [n,m]$. Then for any $c\in\CD[n,m]$ such that $\dsem{c} = [({\ee}_i, A{\ee}_i)]_{i\leq n}$ there exists a circuit $c'\in \SFGform[m,n]$ such that $c \eqIH c'$.
\end{proposition}
\begin{proof}
Consider a rational of the form $\frac{1}{k+xp}$ in $\ratio$, with $k\neq 0$ and $p\in\poly$. This can be seen as a $1\times 1$ matrix of $\Matratio$, yielding the subspace $[(1,\frac{1}{k+xp})] \: 1 \to 1$. The following derivation shows that a circuit of $\CD$, whose semantics (via $\dsem{\cdot}$) is the subspace $[(1,\frac{1}{k+xp})]$, is equal to one in $\SFGform$.
The sequence of applied laws of $\IBpoly$ is: \eqref{eq:scalarsum}$^{\op}$,  \eqref{eq:bcomonunitlaw} + \eqref{eq:bcomonunitlaw}$^{\op}$, \eqref{eq:QFrob} + \eqref{eq:BFrob}, \eqref{eq:ccsliding}, \eqref{eq:flipscalar}, \eqref{eq:scalarmult}. 
\[
\includegraphics[height=2.7cm]{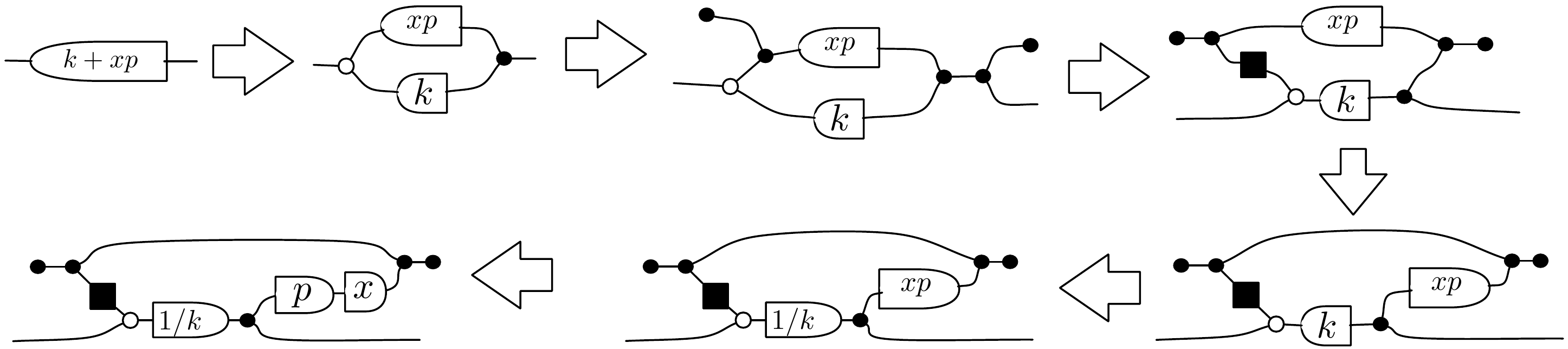}
\]
Now, fix a matrix $A\in \Matratio [n,m]$ and the associated subspace $[({\ee}_i, A{\ee}_i)]_{i\leq n}$. Let $d \in \CD[n,m]$ be the circuit in matrix form constructed as in Remark~\ref{rmk:matrixformFractions}, whose $\dsem{\cdot}$-semantics is the subspace $[({\ee}_i, A{\ee}_i)]_{i\leq n}$: each entry $q$ of the matrix $A$ --- that is, a (rational) fraction $q = p_1/p_2 \in \ratio$ --- is encoded as a component
$\scalarpone \poi \scalarptwoop$ of $d$. By the observation above, we can put any such circuit $\scalarptwoop$ in the form of a circuit of $\SFGform$. Therefore, $d$ is equal in $\IBpoly$ to a circuit $c$ where all components are in $\SFGform$ and, since $\SFGform$ is closed under $\tns$ and $\poi$, then also $c$ is a circuit of $\SFGform$.
\end{proof}

We can now prove our characterisation result.

\begin{proof}[Proof of Theorem \ref{th:SFGcharactRationals}] There is an obvious embedding $\Matratio \to \SVpoly$ mapping $A\in \Matratio [n,m]$ into the subspace $[({\ee}_i, A{\ee}_i)]_{i\leq n}$: the idea is to show that $\SFG$ characterises its image. To do this, define $F \: \SFG \to \Matratio$ as follows. By definition, an arrow $f$ of $\SFG$ is an $\IBpoly$-equivalence class containing a circuit $c$ of $\SFGform$. By Proposition~\ref{prop:circuittomatrix}, $\dsem{c}=[({\ee}_i, A{\ee}_i)]_{i\leq m}$ for some $A$ in $\Matratio$. We let $F$ map $f$ to $A$: Corollary~\ref{prop:isoIHsoundcomplete} guarantees that $F$ is well-defined and faithful. To see that $F$ is full, let $A$ be a matrix in $\Matratio$. Because $\dsem{\cdot}$ is full on $\SVpoly$, there is a circuit $c$ in $\CD$ such that $\dsem{c} = [({\ee}_i, A{\ee}_i)]_{i\leq m}$. By Proposition~\ref{prop:matrixtocircuit}, there is also $d$ in $\SFGform$ such that $c \eqIH d$ and $\dsem{c}=\dsem{d}$.
 We conclude that $F$ is full and faithful and thus an isomorphism. 
\end{proof}

As a consequence of Theorem~\ref{th:SFGcharactRationals}, it is worth mentioning that the restriction of the stream semantics $\strsemO \: \CD \to \SVfps$ to circuits in $\SFGform$ does not actually require the full generality of Laurent series: as shown by~\eqref{eq:ringhoms}, the stream representation of rationals does not need a ``finite past.'' By soundness of $\strsemO$ w.r.t. the equations of $\IBpoly$, we can extend this observation to all the circuits of $\CD$ that are equal to one in $\SFGform$.

\begin{corollary}\label{cor:SFratioNoFinitePast} Suppose $c \in \CD[n,m]$ is a circuit such that $c \eqIH d$ for some $d \in \SFGform[n,m]$. Then $\strsem{c} \in \SVfps[n,m]$ is spanned by vectors over $\fps$. \end{corollary}
\begin{proof} By Theorem~\ref{th:SFGcharactRationals}, $\dsem{d}  = [({\ee}_i, A{\ee}_i)]_{i\leq n}$ for some $A \in \Matratio[n,m]$. Then, by Proposition~\ref{lemma:streamSemUnivProperty},
$\strsem{d} = [\lrnwide{({\ee}_i, A{\ee}_i)}]_{i\leq n}$. Now, each element ${\ee}_i$ of the standard basis of $\frpoly^n$ can be clearly seen also as an element of $\fps^n$. Moreover, as described by diagram~\eqref{eq:ringhoms}, $\ratio$ embeds in $\fps$, meaning that $A{\ee}_i$ can be also seen as an element of $\fps^m$. Therefore, $\strsem{d}$ is spanned by $\fps$-vectors. By Theorem~\ref{cor:completeness}, $c \eqIH d$ implies that $\strsem{c} = \strsem{d}$.  This concludes the proof of the statement. 
\end{proof}



\subsection{Trace Canonical Form for Circuits of $\SFGform$}\label{sec:trace}

In this section we give another demonstration of how $\IBpoly$ can be used to rephrase traditional results of the theory of signal flow graphs in a purely diagrammatic way. We shall show that circuits of $\SFGform$ can always be put, using the equations of $\IBpoly$, into a convenient shape: a core given by a circuit $c$ of $\FC$ without delays, and an exterior part given by a ``bundle'' of feedback loops. We formally introduce this notion below.

\begin{definition} For $n,m,z \in \N$, $c \in \CD[z+n,z+m]$, the \emph{$z$-feedback} $\Tr{z}(c) \in \CD[n,m]$ is the circuit below, for which we use the indicated shorthand notation:
\begin{eqnarray*}
\lower12pt\hbox{$\includegraphics[height=1.3cm]{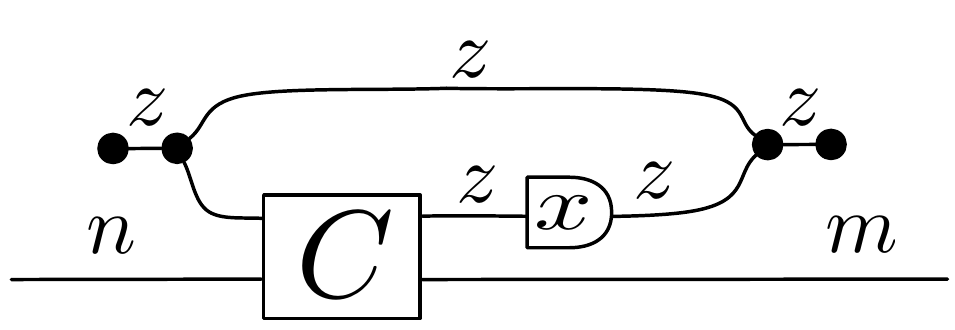}$}
 & \quad \dfop \quad & \lower12pt\hbox{$\includegraphics[height=1.3cm]{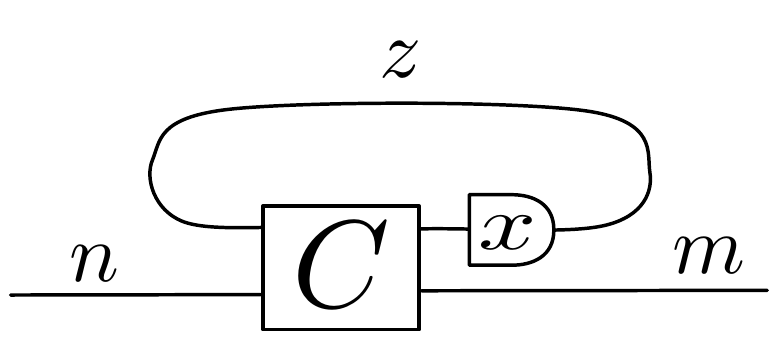}$}
\end{eqnarray*}
\end{definition}
Note that $\Tr{1}(\cdot)$ coincides with the assignment $\gls{TrO}$ given in \eqref{eq:onefeedbackcircuit}, thus notation does not conflict.

\begin{proposition}[Trace form for $\SFGform$]\label{prop:normalformSFG}
Let $\FC \mathsf{\setminus x}$ be the sub-PROP of $\FC$ whose circuits do not contain any delay $\circuitXT$.
For every circuit $d \in \SFGform[n,m]$, there are $z \in \N$ and $c \: z+n \to z+m$ of $\FC\mathsf{\setminus x}$ such that $d\eqIH \Tr{z}(c)$.
%
\end{proposition}

The existence of this form is a folklore result in the theory of signal flow diagrams. Here we provide a novel proof that consists of showing that $\Tr{z}(\cdot)_{z\in\N}$ is a \emph{right trace}~\cite[\S 5.1]{Selinger2009} on the category $\IBpoly$.

\begin{remark}
Note that $\Tr{z}(\cdot)_{z\in\N}$ is \emph{not} the canonical trace induced by the compact closed structure of $\IBpoly$ (\S~\ref{sec:cc}), but a version ``guarded'' by a register. This makes our approach different from other works using traces to model recursion like~\cite{Hasegawa97recursionfrom,stefanescu2000network}: the traces considered there are not guarded by a register and indeed satisfy the ``yanking'' law~\cite{JoyalTracedMonoidalCats}, which is instead false for $\Tr{z}(\cdot)_{z\in\N}$:
    \begin{eqnarray*}
 \lower10pt\hbox{$\includegraphics[height=1.2cm]{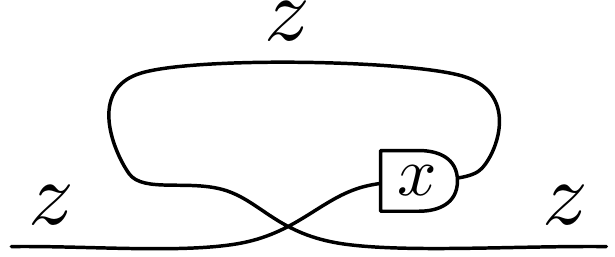}$}
& \neq &
 \lower3pt\hbox{$\includegraphics[height=.3cm]{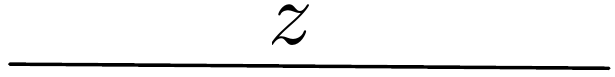}$}.
\end{eqnarray*}
\end{remark}

\medskip


\begin{proposition}\label{prop:trace}The family $\Tr{z}(\cdot)_{z\in\N}$ is a right trace on $\IBpoly$.
\end{proposition}
The axioms of (right) traced categories, as presented in \cite[\S 5.1]{Selinger2009}, are:
\begin{enumerate}
    \item Tightening:
    \begin{eqnarray}\label{eq:traceaxiomTightening}
 \lower9pt\hbox{$\includegraphics[height=1.4cm]{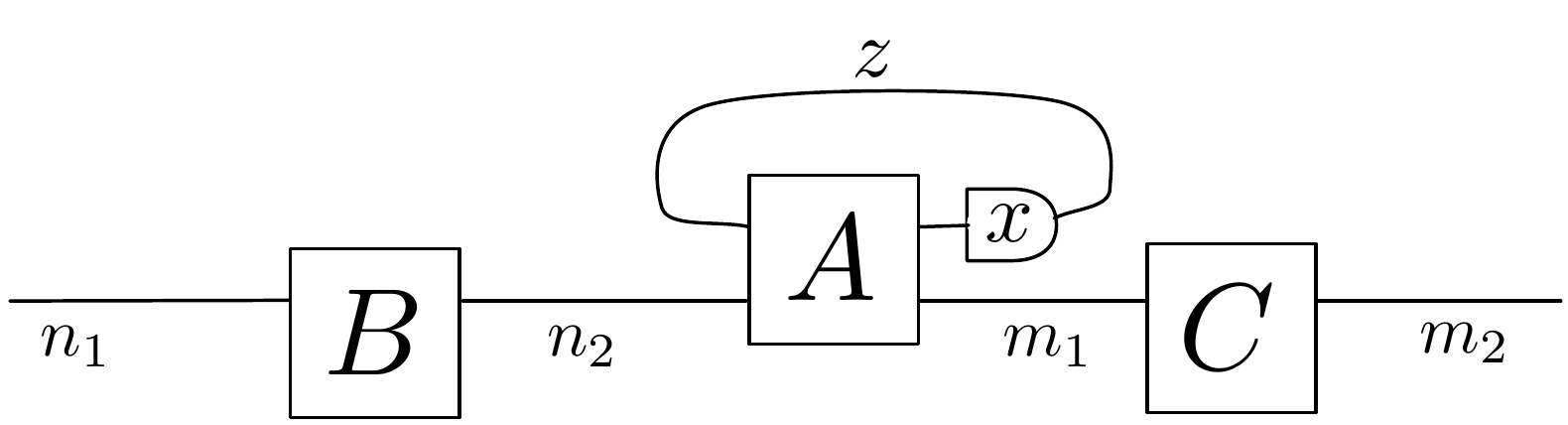}$}
&  \!\!\!\!\!\! \ = \!\!\!\!\!\! \ &
 \lower9pt\hbox{$\includegraphics[height=1.4cm]{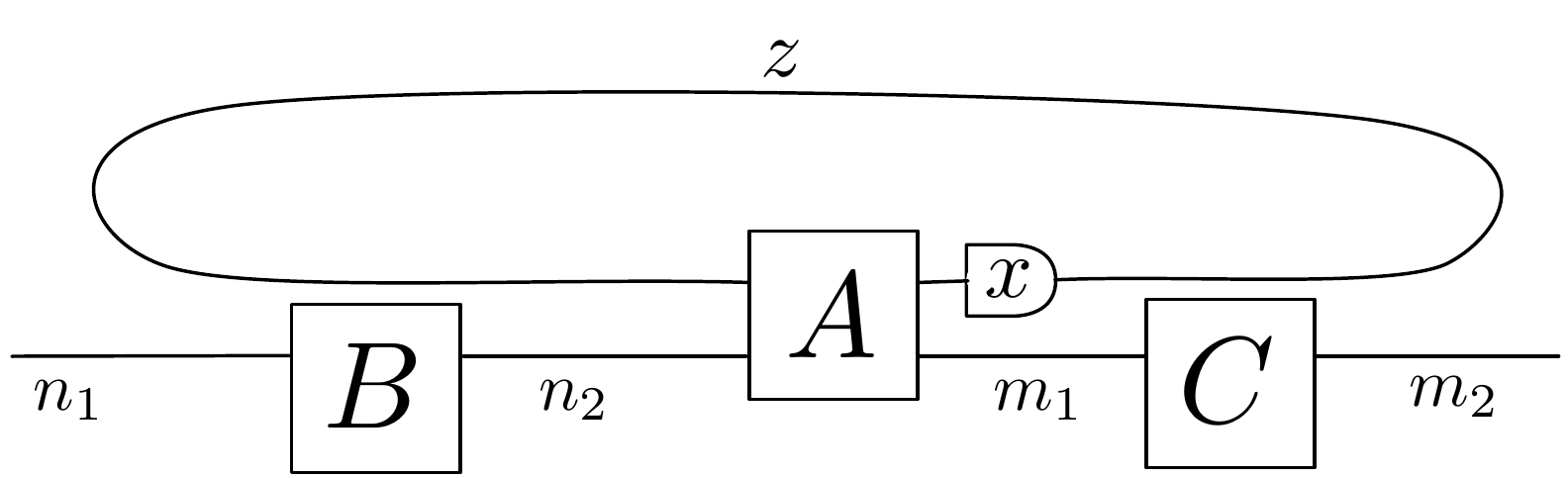}$}\qquad
\end{eqnarray}
    \item Sliding:
    \begin{eqnarray}\label{eq:traceaxiomTightening}
 \lower9pt\hbox{$\includegraphics[height=1.4cm]{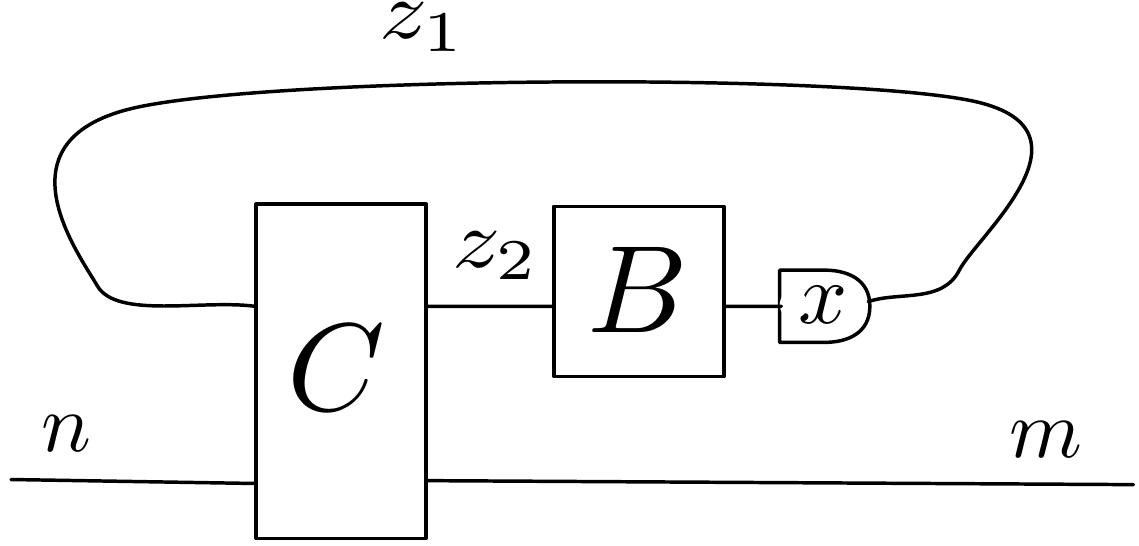}$}
& \ = \ &
 \lower9pt\hbox{$\includegraphics[height=1.4cm]{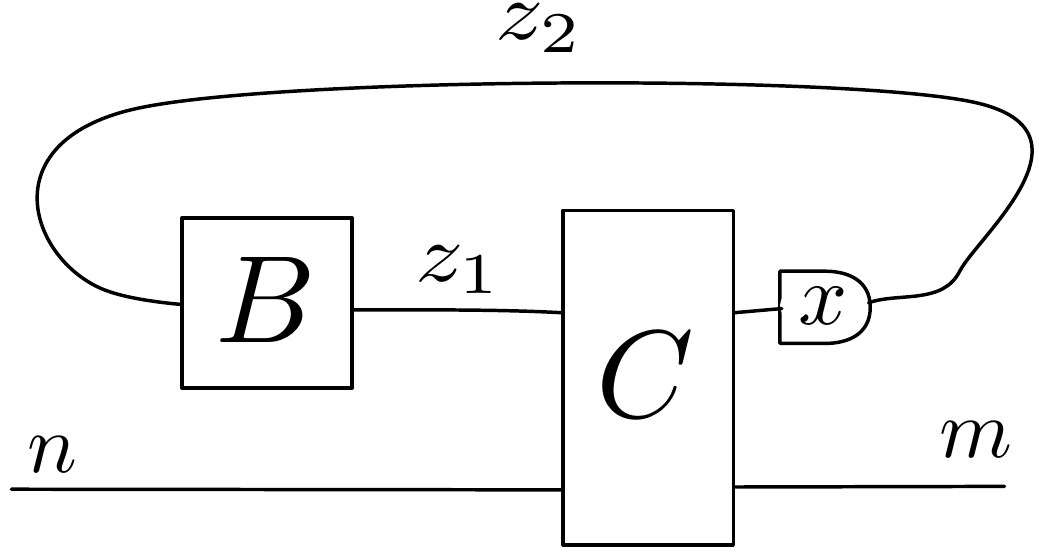}$}
\end{eqnarray}
    \item Vanishing:
    \begin{multicols}{2}\noindent
    \begin{eqnarray}\label{eq:traceaxiomVanishing1}
 \lower9pt\hbox{$\includegraphics[height=1.1cm]{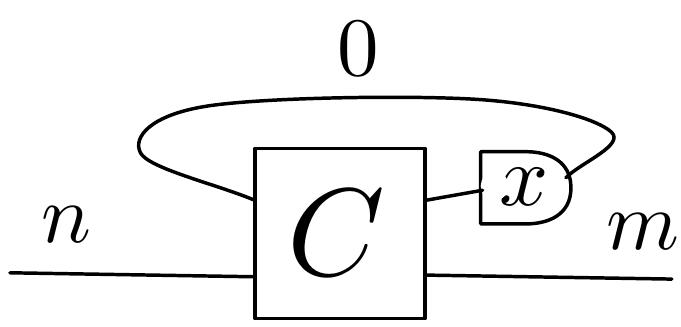}$}
& \ \!\!\!\!\!\!=\!\!\!\!\!\! \ &
 \lower9pt\hbox{$\includegraphics[height=.6cm]{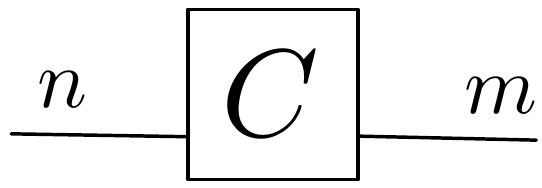}$} \quad\quad\quad
\end{eqnarray}
    \begin{eqnarray}\label{eq:traceaxiomVanishing2}
 \lower14pt\hbox{$\includegraphics[height=1.2cm]{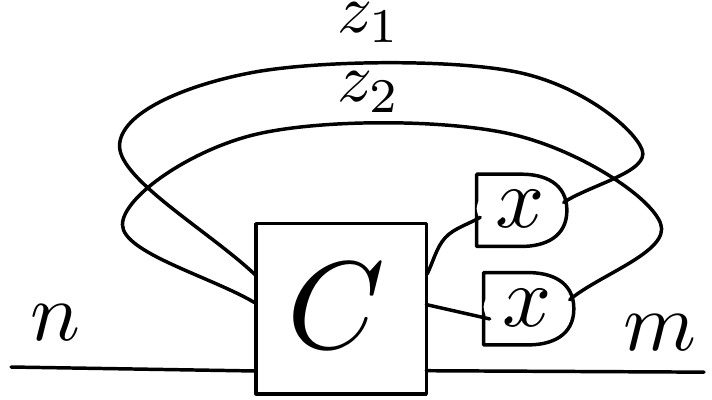}$}
& \ \!\!\!\!\!\!= \!\!\!\!\!\! \ &
 \lower14pt\hbox{$\includegraphics[height=1cm]{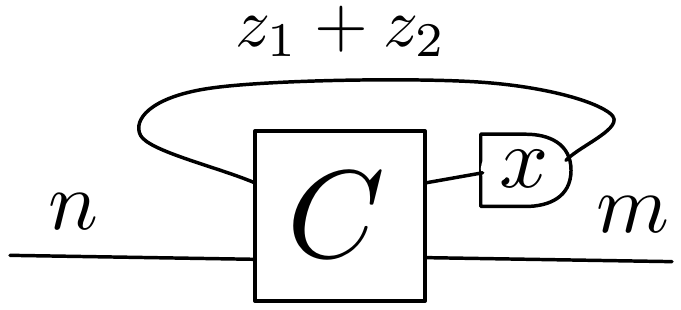}$}
 \quad\quad\quad
\end{eqnarray}
\end{multicols}
    \item Strength:
    \begin{eqnarray}\label{eq:traceaxiomStrength}
 \lower18pt\hbox{$\includegraphics[height=1.8cm]{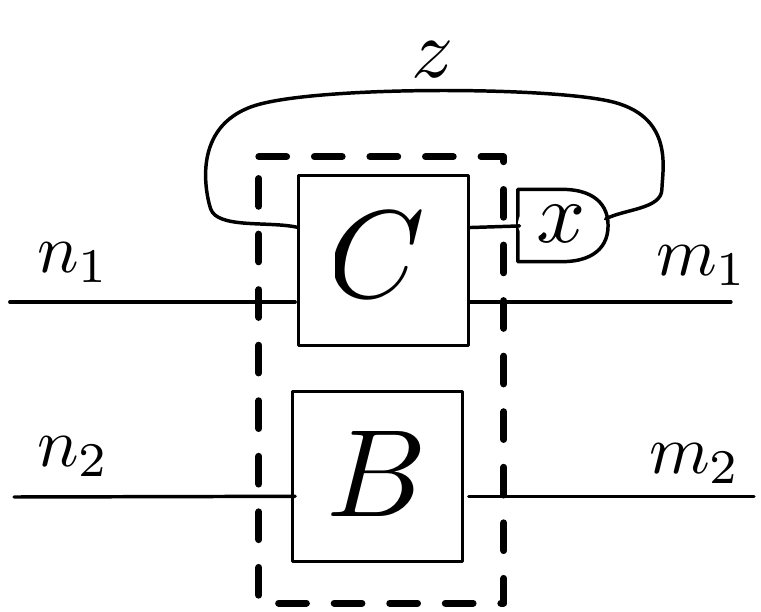}$}
& \ = \ &
 \lower18pt\hbox{$\includegraphics[height=1.8cm]{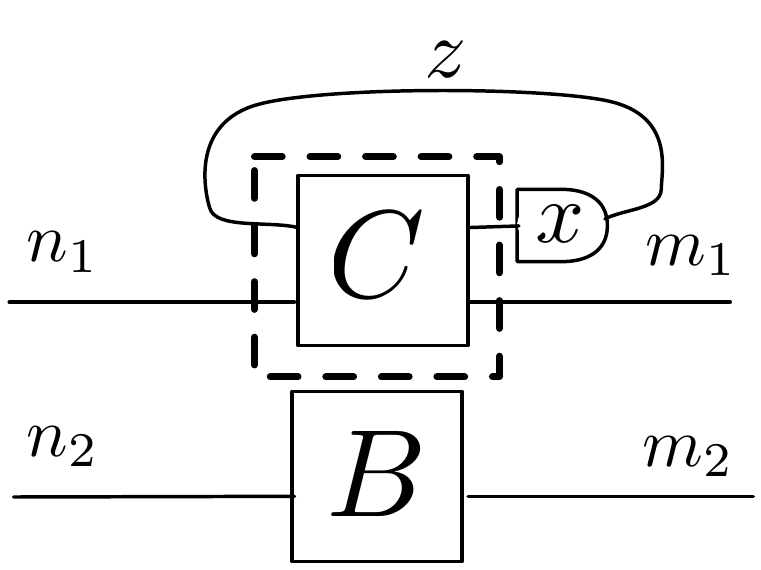}$}
\end{eqnarray}
\end{enumerate}
Tightening and strength hold for our definition of trace simply by laws of symmetric strict monoidal categories. Therefore we focus on sliding and vanishing.

For sliding, it is first useful to record the following lemma.

 \begin{lemma}\label{lemma:xsliding} For any $n,m\in \N$ and circuit $c \in \CD [n,m]$,
 \begin{equation}
   \label{eq:xsliding}
  \lower7pt\hbox{$\includegraphics[height=.7cm]{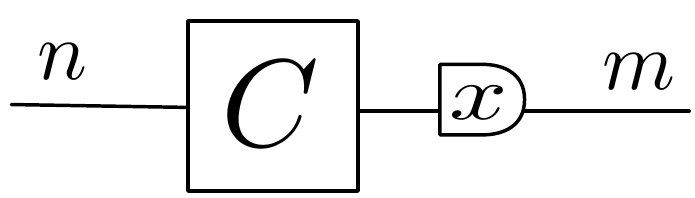}$}
   \quad \eqIH  \quad \lower7pt\hbox{$\includegraphics[height=.7cm]{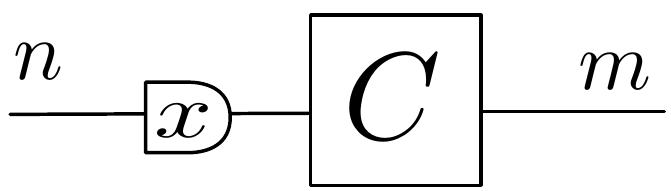}$}
   \end{equation}
 \end{lemma}
 \begin{proof}
The proof is by induction on $c$. For the components in \eqref{eq:SFcalculusSyntax1}, the statement is given for $\Wunit$, $\Wmult$, $\Bcounit$, $\Bcomult$, $\scalar$ and $\circuitX$ by \eqref{eq:scalarwunit}, \eqref{eq:scalarwmult}, \eqref{eq:scalarbcounit}, \eqref{eq:scalarbcomult}, \eqref{eq:scalarmult} and \eqref{eq:scalarmult} respectively. The derivations for $\circuitXop$ and $\coscalar$, with $k \neq 0$, are:
     \begin{eqnarray*}
 \lower5pt\hbox{$\includegraphics[height=.5cm]{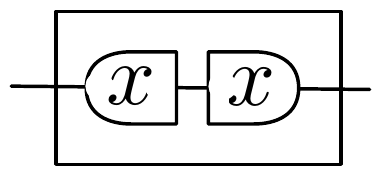}$}
\! \!\!\eql{\eqref{eq:lcmopIH}} \!\!\!
 \lower5pt\hbox{$\includegraphics[height=.5cm]{graffles/idcircuit.pdf}$}
\!  \!\!\eql{\eqref{eq:lcmIH}} \!\!\!
 \lower5pt\hbox{$\includegraphics[height=.5cm]{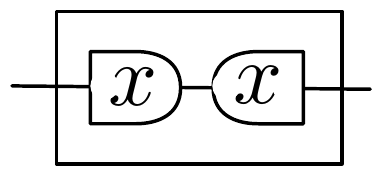}$}
 \qquad
 \lower5pt\hbox{$\includegraphics[height=.5cm]{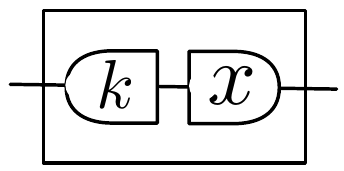}$}
 \! \!\!\eql{\eqref{eq:flipscalar}} \!\!\!
 \lower5pt\hbox{$\includegraphics[height=.5cm]{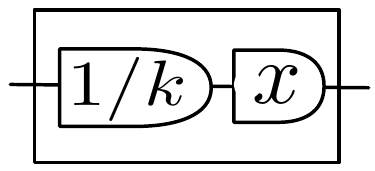}$}
 \! \!\!\eql{\eqref{eq:scalarmult}} \!\!\!
 \lower5pt\hbox{$\includegraphics[height=.5cm]{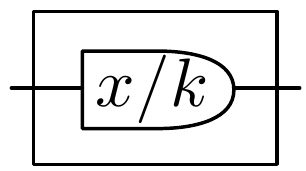}$}
 \! \!\!\eql{\eqref{eq:scalarmult}} \!\!\!
 \lower5pt\hbox{$\includegraphics[height=.5cm]{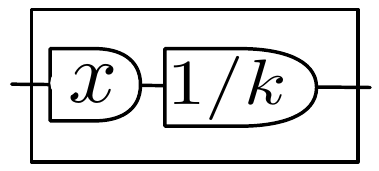}$}
\!  \!\!\eql{\eqref{eq:flipscalar}} \!\!\!
 \lower5pt\hbox{$\includegraphics[height=.5cm]{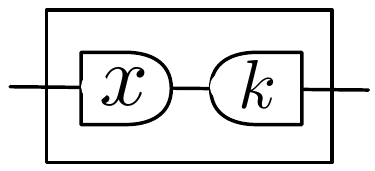}$}
 \end{eqnarray*}
Similarly, one can verify the statement for the remaining cases in \eqref{eq:SFcalculusSyntax2} and \eqref{eq:SFcalculusSyntax3}. The inductive cases of parallel ($\tns$) and sequential ($\poi$) composition of circuits are handled by simply applying the induction hypothesis.
 \end{proof}

\paragraph{Sliding} The following derivation yields the sliding equation:
    \begin{eqnarray*}
 \lower10pt\hbox{$\includegraphics[height=1.4cm]{graffles/traceSlidingl.pdf}$}
 \eql{Lemma~\ref{lemma:xsliding}}
 \lower10pt\hbox{$\includegraphics[height=1.4cm]{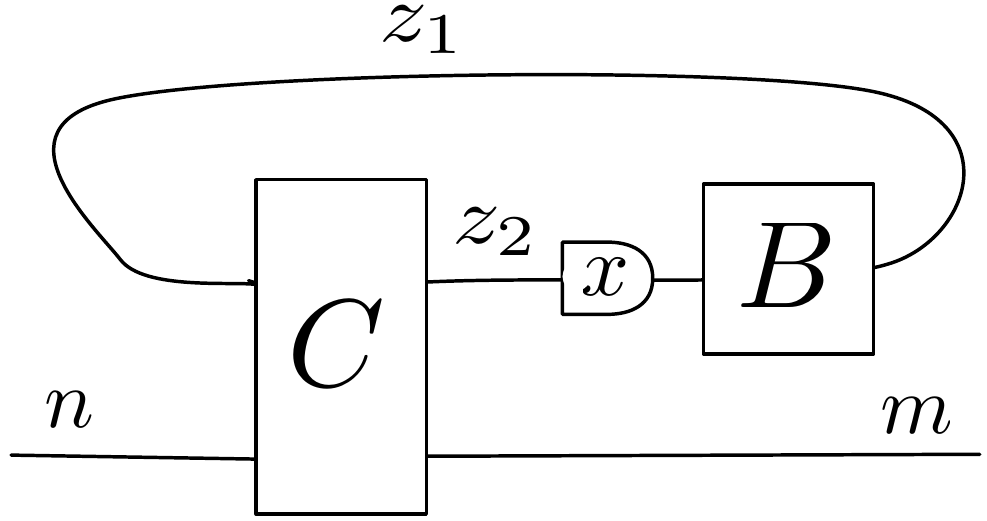}$}
 \eql{\eqref{eq:ccsliding}}
 \lower10pt\hbox{$\includegraphics[height=1.4cm]{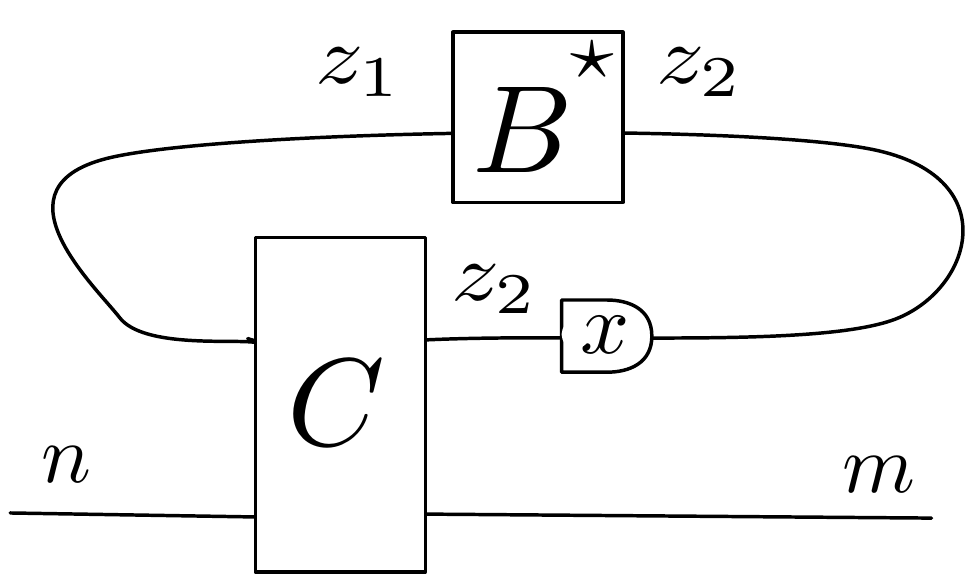}$}
 \eql{\eqref{eq:ccsliding2}}
 \lower10pt\hbox{$\includegraphics[height=1.4cm]{graffles/traceSlidingr.pdf}$}
\end{eqnarray*}
For the two last steps, observe that $\lower4pt\hbox{$\includegraphics[height=.5cm]{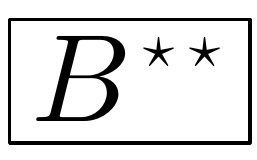}$} = \lower4pt\hbox{$\includegraphics[height=.5cm]{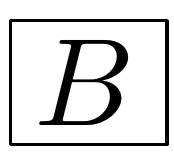}$}$ by definition of $\coc{(\cdot)}$ as in~\eqref{eq:defstar} and by \eqref{eq:gensnake}.

 \paragraph{Vanishing} Concerning vanishing, \eqref{eq:traceaxiomVanishing1} holds because, by definition, $\lccn$, $\rccn$ and $\ncircuitX$ are all equal to $\id_0$ for $n=0$. It remains to check \eqref{eq:traceaxiomVanishing2}. We provide the proof for $z_1,z_2 = 1$. The general case is handled (by induction) by the obvious generalisation of the same argument.

 For this purpose, it will be useful to first introduce the following two equations, holding in $\CD$ by \eqref{eq:sliding2}.
 \begin{multicols}{2}\noindent
    \begin{eqnarray}\label{Bcclsym}
 \lower7pt\hbox{$\includegraphics[height=.7cm]{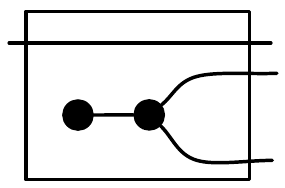}$}
 =
 \lower9pt\hbox{$\includegraphics[height=.8cm]{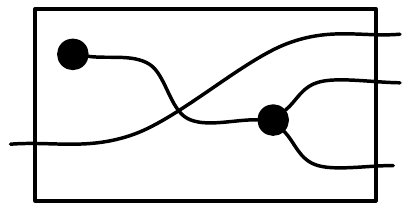}$}
 =
 \lower9pt\hbox{$\includegraphics[height=.8cm]{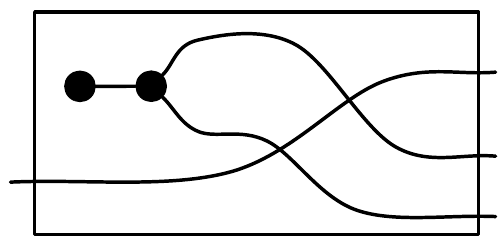}$}\qquad
\end{eqnarray}
    \begin{eqnarray}\label{Bccrsym}
 \lower7pt\hbox{$\includegraphics[height=.7cm]{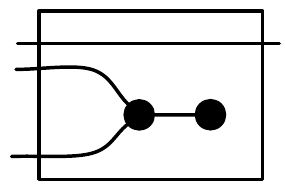}$}
 =
 \lower9pt\hbox{$\includegraphics[height=.8cm]{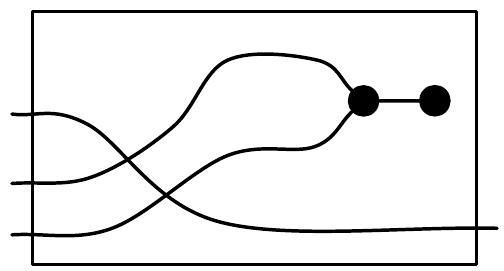}$}
\end{eqnarray}
\end{multicols}
By definition, the first circuit below is $\Tr{1}\Tr{1}(c)$ and the last is $\Tr{2}(c)$. The first step applies \eqref{Bcclsym} and \eqref{Bccrsym}, the second and the third follow by laws~\eqref{eq:sliding1}-\eqref{eq:SymIso} of SMCs.
\begin{equation*}
 \lower14pt\hbox{$\includegraphics[height=3.5cm]{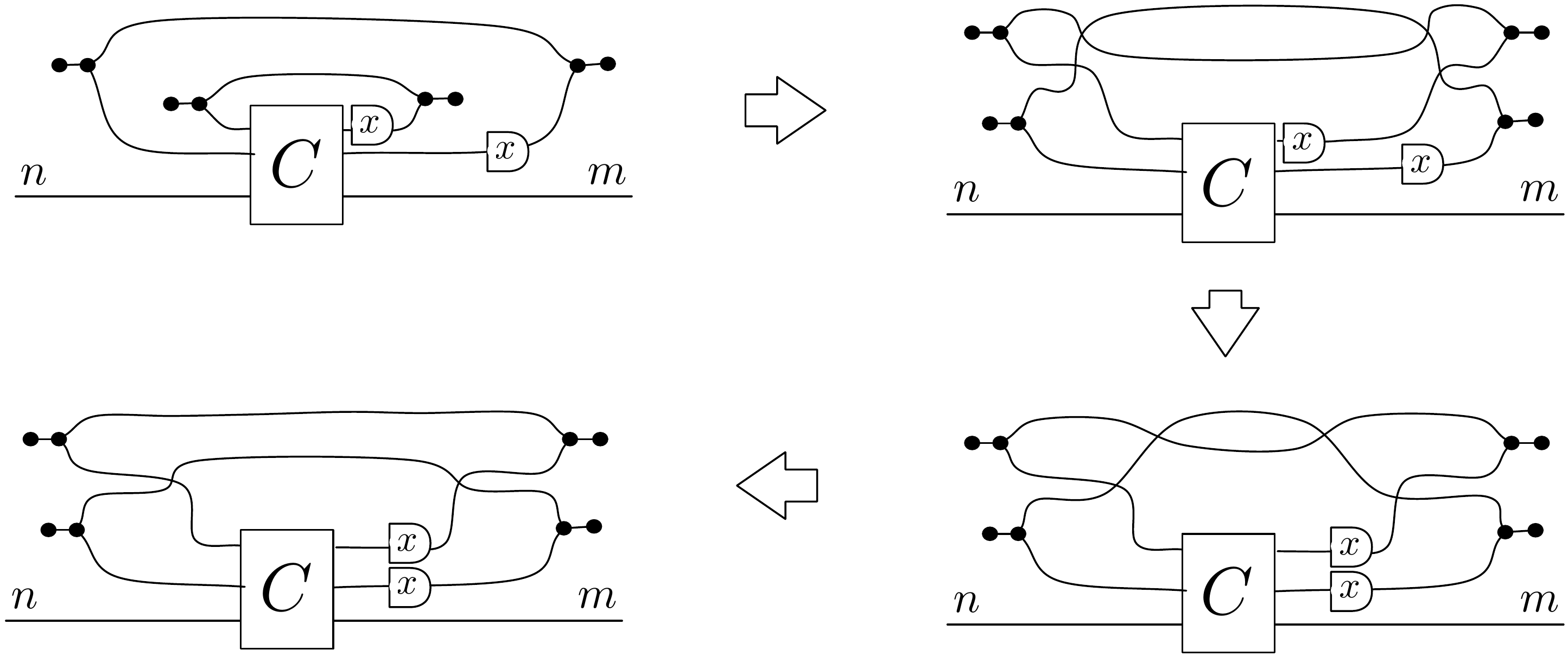}$}
\end{equation*}

This concludes the proof of Proposition~\ref{prop:trace}. We can now give the argument for Proposition~\ref{prop:normalformSFG}.

\begin{proof}[Proof of Proposition~\ref{prop:normalformSFG} ]
The proof goes by induction on a circuit $d$ of $\SFGform$.
If $d$ is a component in \eqref{eq:SFcalculusSyntax1} different from $\circuitXT$, $d = \Tr{0}(d)$. For $\circuitXT$, it is easy to check that $\circuitXT \eqIH \Tr{1}(\symNetT)$.
The second clause of the inductive definition of $\SFGform$ is the case in which $d = \Tr{1}(c)$ for some circuit $c$ of $\SFGform$. By induction hypothesis $c\eqIH \Tr{z}(c')$ for some $c'$ in $\FC\mathsf{\setminus x}$ and thus, by \eqref{eq:traceaxiomVanishing2}, $d \eqIH \Tr{1+z}(c')$. The remaining two cases are the ones in which $d$ is given by sequential or parallel composition of circuits of $\SFGform$---which are, by induction hypothesis, of the form described in the statement. The proof uses the properties of the trace:
\begin{eqnarray*}\label{der:proofcat}
\lower14pt\hbox{\includegraphics[height=1.7cm]{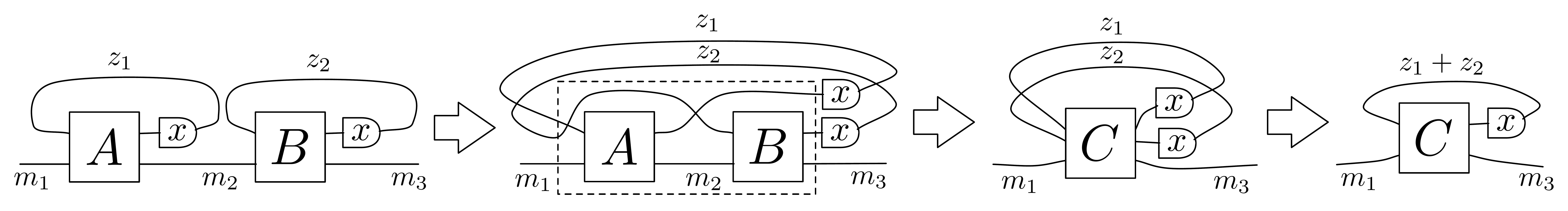}}
\end{eqnarray*}
\begin{eqnarray*}\label{eq:profften}
\lower22pt\hbox{\includegraphics[height=2.2cm]{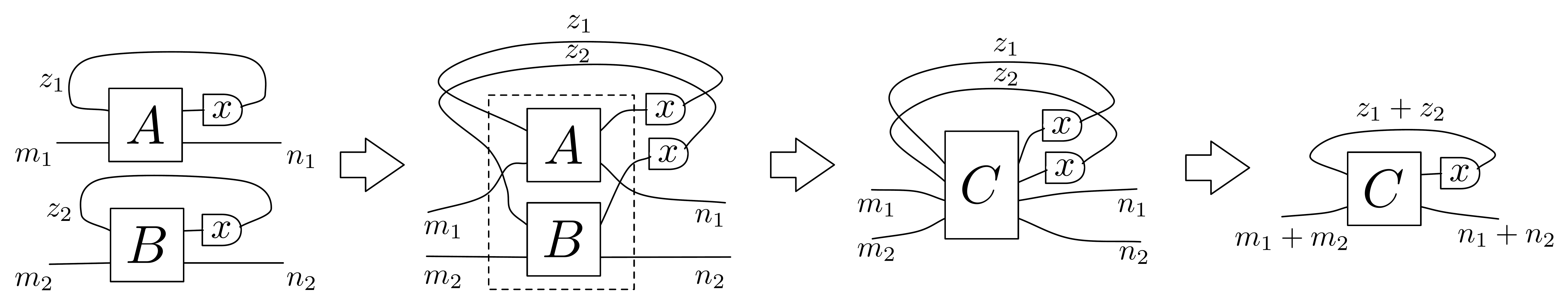}}
\end{eqnarray*}
\end{proof}

\section{Full Abstraction}\label{sec:fullabstract}

So far, the operational and the denotational semantics of the signal flow calculus have been tied only on an intuitive level. In this section we tackle the question of \emph{full abstraction}, that is, how operational and denotational equivalence are formally related. 

\subsection{The Duality of Deadlock and Initialisation}\label{sec:dualitydeadlockinit}

To this aim, an elementary observation is that the denotational semantics is apparently too coarse, since it abstracts away from the finite behaviours that might arise during the executions of the circuits. For example, consider $\circuitXcospan$ and $\IdnetT$: as we have shown in Example~\ref{exm:cospan}, they have the same denotational semantics, namely the set of all pairs $\pair{\sigma}{\sigma}$ of fls. However, some computations of the former circuit can reach a \emph{deadlock}. For instance, we can make a transition from the initial state with labels $k \neq l$, but there are no further possible transitions from the resulting state:
$$\cospanXTzero\dtrans {k}{l} \cospanXTkl \centernot\rightarrow$$
Such failures are not taken into account by the denotational semantics: intuitively, it only considers the successful computations, which are the ones yielding infinite traces. If we restrict to these situations, then $\circuitXcospan$ behaves exactly as the identity circuit $\IdnetT$.
%
Here are more examples of circuits that may reach a deadlocked state. In the first two, the problem is that we may store non-zero values in the registers, whereas in the last two the parallel registers may contain different values in the same state of the computation.
\begin{equation}\label{eq:exdeadlock}
\lower6pt\hbox{$
\includegraphics[height=.5cm,width=.7cm]{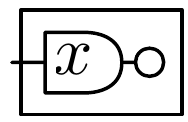} \qquad
\includegraphics[height=.5cm,width=.7cm]{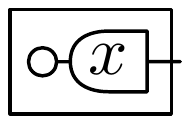} \qquad
\lower4pt\hbox{$\includegraphics[height=.9cm]{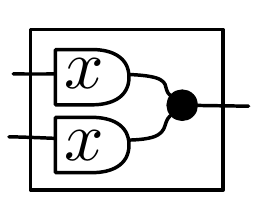}$} \qquad
\lower6pt\hbox{$\includegraphics[height=.9cm]{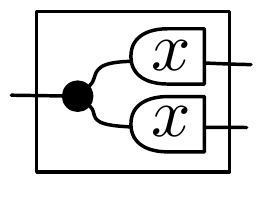}$}
$}
\end{equation}
Our diagnosis is that problematic circuits are those in which internal components (in particular, the delays) have a conflicting design. Note that all the above examples are in \emph{cospan form}, that is, they are of shape $c = c_1 \poi c_2$ with $c_1$ a circuit of $\FC$ and $c_2$ one of $\FCop$. Intuitively, the signal in $c$ is flowing from the left/right boundaries towards the middle, that is, the boundary shared by circuits $c_1$ and $c_2$.

According to our analysis, we can avoid deadlocks by considering instead circuits $d$ in \emph{span form}, i.e. $d = d_1 \poi d_2$ for $d_1$ in $\FCop$ and $d_2$ in $\FC$. Circuits of this shape cannot deadlock since, intuitively, the signal is flowing from the middle boundary towards the left (transmitted by $d_1$) and the right (transmitted by $d_2$).

In order to formalise our observations, first we say that a circuit is \emph{deadlock free} when none of its computations can reach a deadlock - namely, a state from which no transition is derivable. Then we have the following result.
\begin{theorem}\label{thm:spandeadlock}
Circuits of $\CD$ in span form are deadlock free.
\end{theorem}

The proof of Theorem \ref{thm:spandeadlock} is based on the following two lemmas. They formalise the intuition that, in circuits of $\FC$ and $\FCop$, the signal only flows in one direction, thus any input on one boundary (left for $\FC$, right for $\FCop$) produces an output on the other boundary.

\begin{lemma}\label{lemma:spandeadlockFC} Let $c \in \FC[n,m]$ be a circuit. For each $c$-state $s$ and vector $\vv$ of size $n$ there is a $c$-state $t$ and $\uu$ such that $s \dtrans{\vv}{\uu} t$.\end{lemma}
\begin{proof} We proceed by induction on $c$. The base cases and the inductive case of $c = c_1 \tns c_2$ are immediate. For the inductive case $c = c_1 \poi c_2$, let $s$ and $\vv$ be as in the statement. Then $s = s_1 \poi s_2$, where $s_1$ is a $c_1$-state and $s_2$ is a $c_2$-state. By induction hypothesis on $c_1$, there is $\ww$ and $s_1'$ such that $s_1 \dtrans{\vv}{\ww} s_1'$. We can now trigger also the inductive hypothesis on $c_2$ to obtain $\uu$ and $s_2'$ such that $s_2 \dtrans{\ww}{\uu} s_2'$. It follows that $s_1\poi s_2 \dtrans{\vv}{\uu} s_1' \poi s_2'$.
\end{proof}

With a dual argument, we can also show:

\begin{lemma}\label{lemma:spandeadlockFCop} Let $c \in \FCop[n,m]$ be a circuit. For each $c$-state $s$ and vector $\uu$ of size $m$ there is a $c$-state $s$ and $\vv$ such that $s \dtrans{\vv}{\uu} t$.\end{lemma}

\begin{proof}[Proof of Theorem~\ref{thm:spandeadlock}]
Let $c \in \CD[n,m]$ be in span form, that is, $c = c_1 \poi c_2$ with $c_1 \in \FCop[n,z]$ and $c_2 \in \FC[z,m]$ for some $z$.
We need to show that $c_1 \poi c_2$ is deadlock free. For this purpose, let $s$ be a $c_1 \poi c_2$-state. This means that $s =  s_1 \poi s_2$, with $s_1$ a $c_1$-state and $s_2$ a $c_2$-state.
Pick now any vector $\ww$ of size $z$. By Lemma~\ref{lemma:spandeadlockFC} there is $\uu$ and $t_2$ such that $s_2 \dtrans{\ww}{\uu} t_2$. Also, by Lemma~\ref{lemma:spandeadlockFCop}, we get $t_1$ and $\vv$ such that $s_1 \dtrans{\vv}{\ww} t_1$. It follows the existence of a transition $s_1 \poi s_2 \dtrans{\vv}{\uu} t_1 \poi t_2$, meaning that we can always avoid a deadlock situation in computations of $c$.
\end{proof}
Paired with the factorisation Theorem~\ref{Th:factIBR}, Theorem~\ref{thm:spandeadlock} asserts that, for each circuit of $\CD$, there exists an equivalent one in $\IBpoly$ that is deadlock free. This could give us some hope of reconciling the operational and the denotational semantics but, unfortunately, also for some circuits in span form they do not agree: for instance, $\circuitXspan$ and $\IdnetT$ have the same denotational semantics, but all the computations of the former are forced to start with $\dtrans {0}{0}$. Indeed
$$\spanXTzero\dtrans {0}{0} \spanXTk \dtrans {k_1}{k_1} \spanXTktwo \dtrans {k_2}{k_2} \spanXTkthree \dtrans {k_3}{k_3} \dots .$$
Note that after the first transition $\circuitXspan$ behaves exactly as $\IdnetT$: in some sense, the former circuit exhibits a proper behaviour only after an initialisation step. To make this formal, we say that a circuit $c$ is \emph{initialisation free} if, whenever $s_0\dtrans{\zerov}{\zerov}s_1$, then $s_1=s_0$, where $s_0$ is the initial state of $c$. Other basic circuits that suffer from initialisation are displayed below.
\begin{equation}\label{eq:exinit}
\lower11pt\hbox{$
    \raise6pt\hbox{$\includegraphics[height=.5cm,width=.7cm]{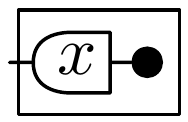}$} \qquad
    \raise6pt\hbox{$\includegraphics[height=.5cm,width=.7cm]{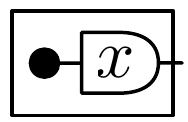}$} \qquad
    \includegraphics[height=.9cm]{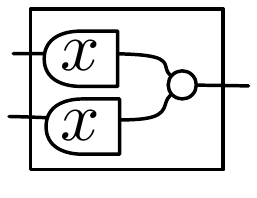} \qquad
    \raise3pt\hbox{$\includegraphics[height=.9cm]{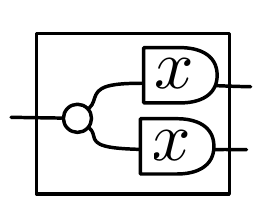}$}
$}
\end{equation}
All problematic circuits above are in span form, meaning that they can be decomposed into $c_1 \poi c_2$, with $c_1$ in $\FCop$ and $c_2$ in $\FC$. The intuition is that any delay in $c_1$ and $c_2$ sends the signal from the common middle boundary towards the outer boundaries, thus requiring a step in which the default value $0$ of each delay is emitted before behaving properly. According to this analysis, such a situation is avoided when we can see all the delays as pointing towards the middle of the circuit. This leads to the following statement.
\begin{theorem}\label{thm:cospaninit}
 Circuits of $\CD$ in cospan form are initialisation free.
\end{theorem}

Also for the proof of Theorem~\ref{thm:cospaninit} it is convenient to record two lemmas. They state a strong form of initialisation freedom for circuits of $\FC$ and $\FCop$.


\begin{lemma}\label{lemma:cospaninitFC} Let $c \in \FC[n,m]$ be a circuit, $s$ its initial state and $t$ another $c$-state. Then $s \dtrans{\zerov}{\uu} t$ implies that $s = t$ and $\uu = \zerov$.
\end{lemma}
\begin{proof} We reason by induction on $c$. 
If $c = \circuitXT$, then its initial state is $\circuitXT^{\labelSep 0}$ and $\circuitXT^{\labelSep 0} \dtrans{0}{l} t$ implies by definition $l = 0$ and $t = \circuitXT^{\labelSep 0}$. The remaining base cases for $c$ and the inductive case $c = c_1 \tns c_2$ are easily verified. We then focus on the case $c = c_1 \poi c_2$. By assumption we are given with a transition $s_1 \poi s_2 \dtrans{\zerov}{\uu} t_1 \poi t_2$. This means that there is $\ww$ such that $s_1 \dtrans{\zerov}{\ww} t_1 $ and $s_2 \dtrans{\ww}{\uu} t_2$. By inductive hypothesis on the transition from $s_1$, $\ww = \zerov$ and $t_1 = s_1$. We can now trigger the inductive hypothesis on the transition from $s_2$ to conclude that $t_2 = s_2$ and $\uu= \zerov$. Then also $s_1 \poi s_2 = t_1 \poi t_2$ and the statement follows.
\end{proof}

\begin{lemma} \label{lemma:cospaninitFCop} Let $c \in \FCop[n,m]$ be a circuit, $s$ its initial state and $t$ another $c$-state. Then $s \dtrans{\vv}{\zerov} t$ implies that $s = t$ and $\vv = \zerov$.
\end{lemma}

\begin{proof} The argument is dual to the one proving Lemma~\ref{lemma:cospaninitFC}.
\end{proof}

\begin{proof}[Proof of Theorem~\ref{thm:cospaninit}]
Let $c \in \CD[n,m]$ be in cospan form, that is, $c = c_1 \poi c_2$ with $c_1 \in \FC[n,z]$ and $c_2 \in \FCop[z,m]$ for some $z$.
We need to show that $c_1 \poi c_2$ is initialisation free. For this purpose, let $s = s_1 \poi s_2$ be the initial state of $c$ and $t = t_1 \poi t_2$ another $c$-state such that $s_1 \poi s_2 \dtrans{\zerov}{\zerov} t_1 \poi t_2$. Then by definition there is $\vv$ such that $s_1 \dtrans{\zerov}{\vv} t_1$ and $s_2 \dtrans{\vv}{\zerov}t_2$. Clearly, since $s_1 \poi s_2$ is the initial $c_1 \poi c_2$-state then $s_1$ is the initial $c_1$-state and $s_2$ the initial $c_2$-state. We can now apply Lemma~\ref{lemma:cospaninitFC} and \ref{lemma:cospaninitFCop} to obtain that $t_1 =s _1$ and $t_2 = s_2$. It follows that $s_1 \poi s_2 = t_1 \poi t_2$, meaning that $c_1 \poi c_2$ is initialisation free.
\end{proof}
Theorems \ref{thm:spandeadlock} and \ref{thm:cospaninit} suggest a duality between deadlock and initialisation, expressible in terms of span and cospan decompositions of circuits of $\CD$. In fact, this is reflected also by the displayed problematic circuits: the ones in \eqref{eq:exdeadlock} are dual to the ones in \eqref{eq:exinit} in a precise sense, namely by changing the black/white colouring and the direction of delays --- this is the ``photographic negative'' transformation $\pn$ described in \S\ref{sec:IBRbCospan}. Moreover, according to our analysis circuits with deadlocks have \emph{more} behaviours (traces) than those prescribed by the denotational semantics, while circuits with initialisation have \emph{fewer} behaviours. We would then expect circuits which are both deadlock and initialisation free to yield exactly the right amount of behaviour: this will be the content of the next section, leading to the full abstraction result (Corollary \ref{cor:fullabstractInitDeadFree}).

\subsection{Reconciling Observation and Denotation}\label{sec:trfls}

Given a circuit $c$ of $\CD$, we define its \emph{observable behaviour} $\osem{c}$ as the pair $\pair{\ftr{c}}{\itr{c}}$ of its finite and infinite traces.
Like the denotational semantics $\strsemO$, also the observable behaviour $\gls{osemO}$ can be expressed in a compositional way, as a PROP morphism from $\CD$ to a certain target PROP that we are going to define below.

In order to do that, we first observe that finite and infinite traces can be equivalently described in terms of polynomials and of fps respectively.
Indeed, in a trace $(\vlist{\alpha},\vlist{\beta})$ of length $z$, each sequence $\alpha_j = k_0  k_1 \dots k_z$ and $\beta_j = l_0  l_1 \dots l_z$ can be encoded as polynomials $k_0 x + k_1 x + \dots + k_z x^z$ and $l_0 x + l_1 x + \dots + l_z x^z$ respectively. Similarly, in an infinite trace $(\vlist{\alpha},\vlist{\beta})$, each stream $\alpha_j = k_0  k_1 \dots $ and $\beta_j = l_0  l_1 \dots$ defines fps $\Sigma^{\infty}_{i=0} k_i x^i$ and $\Sigma^{\infty}_{i=0} j_i x^i$ respectively. We can then see $\ftr{c}$ as a relation between vectors of polynomials and $\itr{c}$ as a relation between vectors of fps.

 On the base of this observation, we take the PROP $\Relpoly \times \Relfps$ as target of $\osemO$, where:
 \begin{itemize}
 \item $\Relpoly$ is the PROP of \emph{$\poly$-relations}, whose arrows $n \to m$ are the relations between $\poly^n$ and $\poly^m$, i.e., the subsets of $\poly^n \times \poly^m$. The monoidal product is disjoint union of relations. Composition, identity and symmetries are defined as for linear relations~(Def.~\ref{def:sv}).
 \item $\Relfps$ is the PROP of \emph{$\fps$-relations}, defined analogously to $\Relpoly$.
 \end{itemize}

 \begin{remark} For the sake of the operational behaviour we are not interested in retaining the information about linearity, whence the choice of plain relations instead of modules. Nonetheless, it is worth remarking that, for any circuit $c \: n \to m$, the set of infinite traces $\itr{c}$ \emph{is} a submodule of $\fps^n \times \fps^m$. 
 \end{remark}

 By definition, arrows of $\Relpoly \times \Relfps$ from $n$ to $m$ are pairs $\pair{f}{g}$ with $f\in \Relpoly[n,m]$ and $g\in \Relfps[n,m]$. The following statement --- whose proof is in Appendix~\ref{app:appendixStream} --- guarantees that $\osemO$ is compositional.

\begin{proposition}\label{prop:OsemFunctor} $\osemO \colon \CD \to \Relpoly \times \Relfps$ is a morphism of PROPs. \end{proposition}

Our aim is now to build a bridge between the observable behaviour and the denotational stream semantics. For this purpose, we first illustrate how to relate infinite traces (i.e., pairs of vectors of fps) and pairs of vectors of formal Laurent series. We say that a trace $\Big(\vlist{\alpha},\vlist{\beta}\Big)\in \fps^n\times\fps^m$ \emph{generates} $\Big(\vlist{\sigma},\vlist{\tau}\Big)\in \laur^n \times \laur^m$  if there exist an \emph{instant} $z \in \Z$ such that

\begin{itemize}
\item[(i)] $\alpha_j$, that means, $\alpha_j(i) = \sigma_j(i +z)$ and  $\beta_h(i) = \tau_h(i +z)$ for all $i \in \N$ and $h$, $j$ with $1 \leq j \leq n$, $1 \leq h \leq m$;
\item[(ii)] $z$ is smaller or equal than any degree of $\sigma_1\dots \sigma_n$, $\tau_1\dots \tau_m$ \footnote{An equivalent, more concise way of expressing conditions $(i)$ and $(ii)$ is by saying that $\sigma_j = \alpha_j \cdot x^z$ and $\tau_j = \beta_j \cdot x^z$, where fps $\alpha_j$ and $\beta_j$ are also regarded as fls. Indeed, multiplication by $x^z$ shifts all the scalars in a fls forward (if $z > 0$) or backward (if $z < 0$) by $z$ positions.}.
    \end{itemize}

To see the intuition behind this notion, recall from \S~\ref{sec:stream} that, whereas fps give a way of encoding streams, fls encode streams ``with a (finite) past''. If we see it as a translation from $(\vlist{\alpha},\vlist{\beta})$ to $(\vlist{\sigma},\vlist{\tau})$, our correspondence takes the fps in $(\vlist{\alpha},\vlist{\beta})$ and fix for them all a common ``present'' moment. For example, the trace $(k_0k_1k_2\dots,\, l_0l_1l_2\dots)\in \fps\times\fps$ generates infinitely many pairs of fls, among which we have:
\begin{center}
$(\dots0\underline{0}k_1k_2k_3\dots,\, \dots0\underline{0}l_1l_2l_3\dots)$ \\
$(\dots00\underline{k_1}k_2k_3\dots,\, \dots00\underline{l_1}l_2l_3\dots)$ \\
$(\dots00k_1\underline{k_2}k_3\dots,\, \dots00l_1\underline{l_2}l_3\dots)$
\end{center}
with choice of present moment $(0,0)$, $(k_1,l_1)$ and $(k_2,l_2)$ respectively. The instant $z \in \Z$ will be $1$ for the first, $0$ for the second and $-1$ for the third pair above.

Conversely, we can start from all the fls in $(\vlist{\sigma},\vlist{\tau})$ and forget about their present moment to obtain streams. The requirement that the instant $z$ is chosen not bigger than any degree implies that only $0$s are removed in the process, that is, there is no information loss. For instance, $(\dots00k_1\underline{k_2}k_3\dots,\, \dots00\underline{l_1}l_2l_3\dots)\in \laur\times\laur$ is generated by infinitely many traces, including:

\begin{center}
 $(k_1k_2k_3\dots, \, 0l_1l_2l_3\dots)$\\
 $(0k_1k_2k_3\dots, \, 00l_1l_2l_3\dots)$\\
 $(00k_1k_2k_3\dots, \, 000l_1l_2l_3\dots).$
 \end{center}
 The instant $z$ is chosen to be $-1$ for the first, $-2$ for the second and $-3$ for the third pair above.

 To complete the picture, we relate finite and infinite traces. A trace $(\vlist{\alpha},\vlist{\beta})\in \poly^n \times \poly^m$ of length $z$ is a \emph{prefix} of an infinite trace $(\vlist{\gamma},\vlist{\delta})\in \fps^n \times \fps^m$ iff $\alpha_j(i)=\gamma_j(i)$ and $\beta_h(i)=\delta_h(i)$ for all $0 \leq i\leq z$, $1 \leq j\leq n$ and $1 \leq h\leq m$.

We are now ready to define a correspondence between observable behaviour and relations of $\laur$-vector spaces.

\begin{definition}\label{def:satforget}
Let $\pair{f}{g}$ be an arrow of $\Relpoly \times \Relfps$.
We define $\Sat\pair{f}{g}$ as the following subset of $\laur^n \times \laur^m$:
\begin{eqnarray*}
\left\{\Big(\vlist{\sigma},\vlist{\tau}\Big) \right.
& \Big| & \text{ there exist a trace } \Big(\vlist{\alpha},\vlist{\beta}\Big) \in g \\
&& \left. \text{ generating }\Big(\vlist{\sigma},\vlist{\tau}\Big)\right\}
\end{eqnarray*}
In the converse direction, given a subset of $\laur^n \times \laur^m$,
we define $\Forget(S)$ as the pair $\pair{f}{g} \in \Relpoly \times \Relfps[n,m]$ where $g \in \Relfps[n,m]$ is given as
\begin{eqnarray*}
\left\{\Big(\vlist{\alpha},\vlist{\beta}\Big) \right.
& \Big| & \text{ there exist a pair } \Big(\vlist{\sigma},\vlist{\tau}\Big) \in S \\
&& \left.  \text{ generated by }\Big(\vlist{\alpha},\vlist{\beta}\Big) \right\}
\end{eqnarray*}
and $f \in \Relpoly[n,m]$ is the set of all prefixes of the traces in $g$.
\end{definition}
Intuitively, the action of $\Sat$ on $\pair{f}{g}$ is to forget the first component and generate all the vectors of fls generated by vectors of fps in $g$.
Conversely, we can describe $\Forget$ as abstracting away the choice of the present for all fls and represent them as fps. This gives the second element in the target pair $\pair{f}{g}$: the first is irrelevant, since it only consists of the prefixes of traces in $g$ (in particular, we do not generate any deadlock trace).
To see how $\Sat$ and $\Forget$ work more precisely, we shall consider the following example. 

\begin{example}\label{ex:SatForget} Recall that, for $c = \circuitXspan$, the set $\itr{c}$ consists of all infinite traces of the form $(0 k_0 k_1 k_2 \dots, 0 k_0 k_1 k_2  \dots)$, that is, $c$ behaves as the identity after one initialisation step. Since this circuit is deadlock free, the set $\ftr{c}$ contains all and only those finite traces which are prefixes of some infinite trace in $\itr{c}$.

One can check that $\Sat\osem{c}$ is the identity relation $\{(\sigma,\sigma) \mid \sigma \in \laur\}$ (which is actually equal to the denotational semantics $\strsem{c}$, see~\S\ref{ss:streamIB}). For instance, the pair of fls \begin{equation}\label{eq:exfls}(\dots00\underline{k_0}k_1k_2k_3\dots,\, \dots00\underline{k_0}k_1k_2k_3\dots)\end{equation} is generated by $(0 k_0 k_1 k_2 \dots, 0 k_0 k_1 k_2  \dots)$. If we then apply $\Forget$ to $\Sat\osem{c}$, we obtain a pair $( f,g)\in \Relpoly \times \Relfps[1,1]$. Note that $\ftr{c}$ and $\itr{c}$ are strictly included in $f$ and $g$ respectively. Indeed, provided that $k_0 \neq 0$, $(k_0 k_1 k_2 \dots, k_0 k_1 k_2  \dots) \notin \itr{c}$ but it belongs to $g$ since it generates the pair \eqref{eq:exfls}. Also, any finite prefix of $(k_0 k_1 k_2 \dots, k_0 k_1 k_2  \dots)$ is in $f$ but not in $\ftr{c}$.

We now focus on circuit $d = \circuitXcospan$. The set $\itr{d}$ consists of all infinite traces of the form $(k_0 k_1 k_2 \dots, k_0 k_1 k_2  \dots)$ and the set $\ftr{d}$ consists of either prefixes of $\itr{d}$ or traces leading to deadlocks having the form $(k_0 k_1 \dots k_u l, k_0 k_1  \dots k_u l')$ with $l\neq l'$. It is easy to check that these traces are lost when applying $\Forget\Sat$ to $\osem{d}$, while no infinite trace is added or removed.
\end{example}

The example above suggests that the composite mapping $\Forget\Sat$ \emph{enlarges} the set of observable behaviours for circuits with initialisation (e.g.\ $\circuitXspan$) and, dually, \emph{restricts} it for circuits with deadlocks (e.g.\ $\circuitXcospan$). The next statement illustrates the extent of these observations.
\begin{theorem}\label{thm:deadinit}
Let $c$ be a circuit of $\CD$. Then the following hold:\footnote{The inclusions are meant to be component-wise, e.g. $\itr{c} \subseteq g$ and $\ftr{c} \subseteq f$ where $\pair{\ftr{c}}{\itr{c}} = \osem{c}$ and $\pair{f}{g} = \Sat\Forget\strsem{c}$.}
\begin{enumerate}[(a)]
 \item $\Sat\osem{c} = \strsem{c}$. \label{point:Posem=Dsem}
 \item if $c$ is deadlock free, then $\osem{c} \subseteq\Forget\Sat\osem{c}$. \label{point:deadIncl}
 \item if $c$ is initialisation free, then $\osem{c} \supseteq\Forget\Sat\osem{c}$. \label{point:initIncl}
\end{enumerate}
\end{theorem}
\begin{proof} See Appendix~\ref{app:appendixStream}. \end{proof}
Statement \ref{point:Posem=Dsem} above is instrumental in showing full abstraction (Corollary \ref{cor:fullabstractInitDeadFree}), but is also of independent interest. Indeed it allows to immediately derive that
\begin{corollary}\label{cor:satosem=dsem}
 For any two circuits $c,d\in \CD$, $\strsem{c}=\strsem{d}$ if and only if $\Sat\osem{c}=\Sat\osem{d}$.
\end{corollary}
In a sense, Corollary \ref{cor:satosem=dsem} tells us under which conditions an external observer cannot distinguish circuits that have the same denotation. This is the case whenever $\Sat\osem{c}=\Sat\osem{d}$, that is, the observation of $c$ and $d$ can be only made ``up-to $\Sat$''. Intuitively, this amounts to imposing the following two conditions, stemming from the definition of $\Sat$. First, we prevent the observation of finite behaviour --- because $\Sat$ disregards $\ftr{c}$ and $\ftr{d}$. This means that we cannot detect deadlock and, for instance, $\circuitXcospan$ and $\IdnetT$ become indistinguishable. Second, we prevent an external agent from choosing \emph{when} to begin the observation. For instance, take $c = \circuitXspan$. By observing the pair $(\dots 00 k_0\underline{k_1} k_2 \dots, \dots 00 k_0\underline{k_1} k_2 \dots)$ in $\Sat\osem{c}$, in principle we are not able to judge whether it has been generated by a trace $(k_0 k_1 k_2 \dots, k_0 k_1 k_2 \dots)$ or $(0 k_0 k_1 k_2 \dots, 0 k_0 k_1 k_2 \dots)$ in $\osem{c}$: the definition of $\Sat$ allows for both (and infinitely many other) options. Since our view is restricted to $\Sat\osem{c}$, we  cannot tell if the observation of the actual stream starts with $0$ or $k_0$. Therefore, from that viewpoint $\circuitXspan$ and $\IdnetT$ are indistinguishable. 

In general, when observations can be made without the restrictions of $\Sat$, one \emph{can} distinguish amongst circuits that are denotationally equivalent, as explained in Example~\ref{ex:SatForget}.
Statements \ref{point:deadIncl} and \ref{point:initIncl} in Theorem \ref{thm:deadinit} allow us to derive that observations and denotations do coincide for the class of well-behaved circuits that do not suffer from deadlocks and initialisation steps.

\begin{corollary}[Full Abstraction]\label{cor:fullabstractInitDeadFree}
 For any two circuits $c$ and $d$ of $\CD$ that are deadlock and initialisation free,
 \begin{center}$\strsem{c}=\strsem{d}$ if and only if $\osem{c}=\osem{d}$.\end{center}
\end{corollary}
\begin{proof}
The statement is given by the following two chains of implications.
\[\osem{c}=\osem{d} \quad \Impl{} \quad \Sat\osem{c}=\Sat\osem{d} \Impl{Th. \ref{thm:deadinit}.\ref{point:Posem=Dsem}} \strsem{c}=\strsem{d}.\]
\begin{eqnarray*}
\strsem{c}=\strsem{d} &\Impl{Th. \ref{thm:deadinit}.\ref{point:Posem=Dsem}} \Sat\osem{c}=\Sat\osem{d}
\quad \Impl{} \quad \Forget\Sat\osem{c}=\Forget\Sat\osem{d}
 \Impl{Th. \ref{thm:deadinit}.\ref{point:deadIncl},\ref{point:initIncl}}& \osem{c}=\osem{d} .
\end{eqnarray*}
\end{proof}
There is an important class of circuits for which the operational semantics is fully abstract.
\begin{proposition}\label{prop:SFGfree}
Every circuit in $\SFGform$ is deadlock and initialisation free.
\end{proposition}
\begin{proof}
The leading intuition is that the rules for inductively constructing circuits of $\SFGform$ do not allow for the conflicting design originating the phenomena of initialisation and deadlock. In particular, adding feedbacks to a circuit of $\SFGform$ preserves initialisation and deadlock freedom. Coming to the formal details, the proof goes by induction on a circuit $c \: n \to m$ in $\SFGform$. In fact, we will need to strengthen our inductive hypothesis by showing the following two properties that imply initialisation and deadlock freedom:
\begin{itemize}
\item[(a)] Let $t_0$ be the initial state of $c$ and $t_1$ another $c$-state. Then $t_0 \dtrans{\zerov}{\uu} t_1$ implies that $t_0 = t_1$ and $\uu = \zerov$. 
\item[(b)] For each $c$-state $t$ and vector $\vv$ of length $n$ there is a $c$-state $t_1$ and $\uu$ such that $t \dtrans{\vv}{\uu} t_1$.
\end{itemize}
For the base case of the induction, suppose that $c$ is in $\FC$. Then the two properties follow by Lemmas~\ref{lemma:spandeadlockFC} and~\ref{lemma:cospaninitFC}. 

For the inductive case, we consider the sequential composition. Suppose that $c=c'\poi c''$ and let $t_0'$ and $t_0''$ be initial states of $c'$ and $c''$ respectively. By definition of initial state, $t_0=t_0'\poi t_0''$. Then $t_0'\poi t_0'' \dtrans{\zerov}{\uu} t_1$ is possible only if
\begin{center}$t_0' \dtrans{\zerov}{\vv} t_1'$ and $t_0'' \dtrans{\vv}{\uu} t_1''$ with $t_1=t_1'\poi t_1''$.\end{center}  We can apply the inductive hypothesis to the first transition to have that $t_1'=t_0'$ and $\vv=\zerov$. Thanks to the latter we can now apply the inductive hypothesis to the second transition to get $t_1''=t_0''$ and $\uu=\zerov$. To conclude the proof of (a), it is sufficient to note that $t_1=t_1'\poi t_1''=t_0'\poi t_0''=t_0$. For (b), observe that any $c$-state $t$ is equal to the composition of a $c'$-state $t'$ and a $c''$-state $t''$. By using the inductive hypothesis on $t'$, we have that for any $\vv$, $t' \dtrans{\vv}{\ww} t_1'$. By the inductive hypothesis on $t''$, we have that $t'' \dtrans{\ww}{\uu} t_1'$. To conclude the proof of (b), observe that $t=t'\poi t'' \dtrans{\vv}{\uu}t_1'\poi t_2'$.

The inductive proof for the tensor is completely analogous, while the one for the trace operator is more challenging. We give it below.

Suppose that $c$ is $\Tr{}(d)$ for some $d$ in $\SFGform$.
 \begin{equation*}
   \includegraphics[height=.8cm]{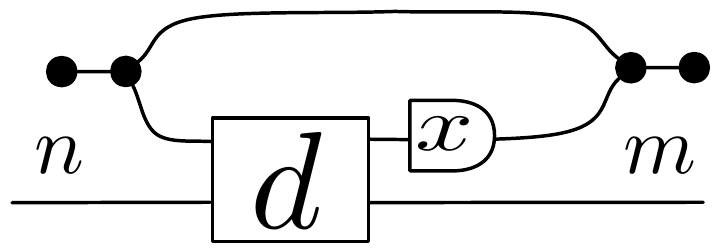}
 \end{equation*}
By construction, the initial state $t_0$ of $c$ can be depicted as follows, where $s_0$ is the initial state of $d$.
  \begin{equation*}
  \includegraphics[height=.8cm]{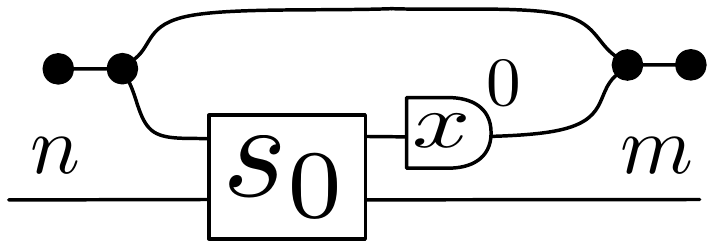}
 \end{equation*}
 Now suppose that we can make the following transition from $t_0$.
  \begin{equation*}
  \lower9pt\hbox{$\includegraphics[height=.8cm]{graffles/SFGdeadinitfree2.pdf}$}
  \qquad
  \dtransF{\scriptsize\left(\begin{array}{c}
                  				  \!\!\!0\!\!\! \\
                				  \!\!\!\vdots\!\!\! \\
                				  \!\!\!0\!\!\!
                				\end{array}\right)}
                {\scriptsize\left(\begin{array}{c}
        				  \!\!\!l_1\!\!\! \\
        				  \!\!\!\vdots\!\!\! \\
        				  \!\!\!l_{m}\!\!\!
                				\end{array}\right)}
        \qquad
    \lower9pt\hbox{$\includegraphics[height=.8cm]{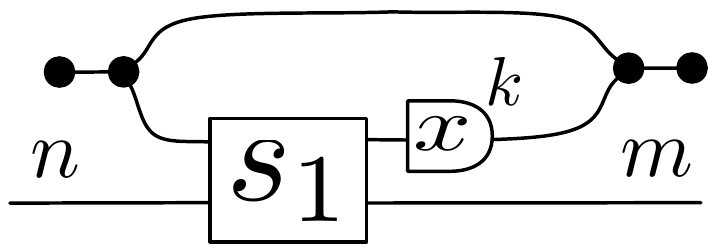}$}
 \end{equation*}
It follows that we can also make the following transition from $s_0$.
  \begin{equation*}
  \lower6pt\hbox{$\includegraphics[height=.5cm]{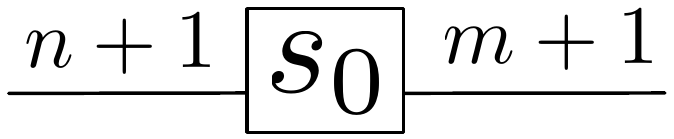}$} \qquad
  \dtransF{\scriptsize\left(\begin{array}{c}
                				  \!\!\!0\!\!\! \\
                				  \!\!\!\vdots\!\!\! \\
                				  \!\!\!0\!\!\!
                				\end{array}\right)}
                {\scriptsize\left(\begin{array}{c}
        				  \!\!\!k\!\!\! \\
        				  \!\!\!l_1\!\!\! \\
        				  \!\!\!\vdots\!\!\! \\
        				  \!\!\!l_{m+1}\!\!\!
                				\end{array}\right)}
    \qquad
    \lower6pt\hbox{$\includegraphics[height=.5cm]{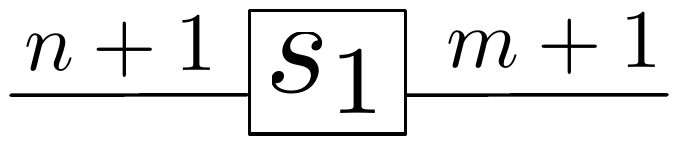}$}
 \end{equation*}
By inductive hypothesis, property (a) holds for $d$. This means that $k, l_1,\dots,l_m$ are all equal to $0$ and $s_0 = s_1$. It follows that $t_0 = t_1$ and property (a) holds for $c$.

\noindent For property (b), let us fix ${\scriptsize \left(\begin{array}{c}
				  \!\!\!k_1\!\!\! \\
				  \!\!\!\vdots\!\!\! \\
				  \!\!\!k_n\!\!\!
				\end{array}\right)}$ and a generic $c$-state as follows.
  \begin{equation*}
  \includegraphics[height=.8cm]{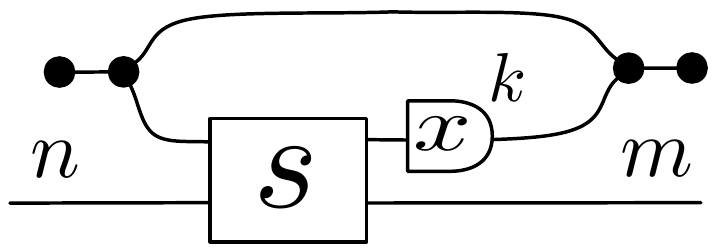}
 \end{equation*}
By construction $s$ is a $d$-state and by inductive hypothesis property (b) holds for $d$. This implies the existence of a transition to a $d$-state $s'$:
  \begin{equation*}
    \lower6pt\hbox{$\includegraphics[height=.5cm]{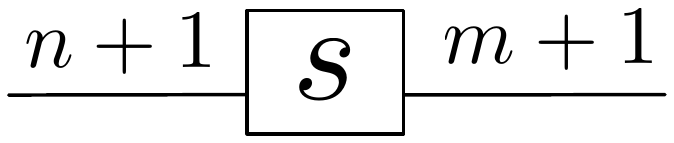}$}
              \qquad
  \dtransF{\scriptsize\left(\begin{array}{c}
                				  \!\!\!k\!\!\! \\
                  				  \!\!\!k_1\!\!\! \\
                				  \!\!\!\vdots\!\!\! \\
                				  \!\!\!k_n\!\!\!
                				\end{array}\right)}
                {\scriptsize\left(\begin{array}{c}
        				  \!\!\!l_1\!\!\! \\
        				  \!\!\!\vdots\!\!\! \\
        				  \!\!\!l_{m+1}\!\!\!
                				\end{array}\right)}
                        \qquad
    \lower6pt\hbox{$\includegraphics[height=.5cm]{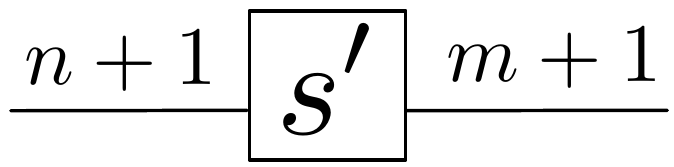}$}
 \end{equation*}
It then follows the existence of the following transition to a $c$-state:
  \begin{equation*}
    \lower9pt\hbox{$\includegraphics[height=.8cm]{graffles/SFGdeadinitfree6.pdf}$} \qquad
  \dtransF{\scriptsize\left(\begin{array}{c}
                  				  \!\!\!k_1\!\!\! \\
                				  \!\!\!\vdots\!\!\! \\
                				  \!\!\!k_n\!\!\!
                				\end{array}\right)}
                {\scriptsize\left(\begin{array}{c}
        				  \!\!\!l_2\!\!\! \\
        				  \!\!\!\vdots\!\!\! \\
        				  \!\!\!l_{m+1}\!\!\!
                				\end{array}\right)}
                                  \qquad
                  \lower9pt\hbox{$\includegraphics[height=.8cm]{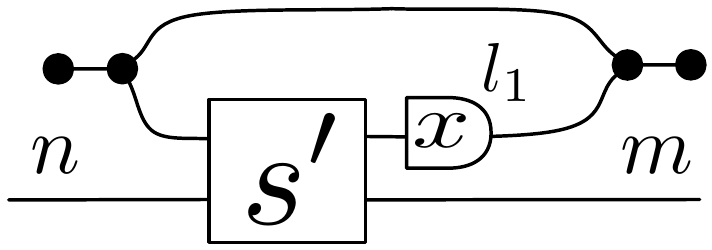}$}
 \end{equation*}
This proves property (b) for the circuit $c$.
\end{proof}
Proposition~\ref{prop:SFGfree} shows that the syntactic restrictions on circuits of $\SFGform$ guarantee deadlock and initialisation freedom. Concerning arbitrary circuits of $\CD$, so far we have seen the span form, preventing deadlocks, and the cospan form, avoiding initialisation steps. In the next section we shall show how to transform, within the equational theory of $\IBpoly$, any circuit $c$ of $\CD$ into one $d$ of $\SFGform$ which, by Proposition~\ref{prop:SFGfree}, features both properties.
%
%

\section{Realisability}
\label{sec:realisability}

The previous sections showed that $\SFGform$ is particularly well-behaved from different points of view:
\begin{enumerate}[(a)]
\item modulo the equational theory of $\IBpoly$, the syntax $\SFGform$ captures the rational behaviours in $\SVpoly$ (Theorem~\ref{th:SFGcharactRationals}). \label{pt:SFGbenefits}
\item any circuit equivalent to one in $\SFGform$ can be given a stream semantics without requiring the full generality of Laurent series (Corollary~\ref{cor:SFratioNoFinitePast}).
\item most importantly, the circuits represented by diagrams of $\SFGform$ are safe with respect to the design flaws identified in \S~\ref{sec:fullabstract}, such as deadlocks and initialisation steps (Proposition~\ref{prop:SFGfree}). The operational semantics for this class of circuits can be really thought as describing the step-by-step execution of a state-machine (\emph{cf.} Remark~\ref{rmk:signalflowOpSem}). Also, operational and denotational equivalence coincide (Corollary~\ref{cor:fullabstractInitDeadFree}).
\end{enumerate}
In light of these results, what can we say about
circuits in $\CD$ that are \emph{not} equivalent to any circuit in $\SFGform$? Do they define a more expressive family of stream transformers under the semantics $\strsemO$? Are they inherently ill-behaved when executed as state-machines?

In this section we demonstrate that the answer to the last questions is \emph{no}: in fact, within the equational theory of $\IBpoly$, $\CD$ is nothing else but a ``jumbled up'' version of $\SFGform$. More precisely, while every
circuit in $\SFGform$ has inputs on the left and outputs on the right, for every circuit in $\CD$ there is a way of partitioning its left and right ports into ``inputs'' and ``outputs'', in the sense that appropriate rewiring yields an $\IBpoly$-equal circuit
in $\SFGform$. The main
result of this section is the \emph{realisability} theorem (Theorem~\ref{thm:realisability}) which guarantees that such an input-output partition
exists, i.e., every circuit in $\CD$ is a rewired circuit in $\SFGform$~\footnote{Note that such a partition is not unique, and this fact corresponds to the physical intuition that in some circuits there is more than one way of orienting flow --- see Examples~\ref{ex:twoimplementations} and~\ref{ex:implementations}.}.


%
%


\smallskip
We begin by giving a precise definition of what we mean by ``jumbling up'' the wires of a circuit.
First, for each $n,m \in \N$, we define circuits $\alpha_n \colon n \to 1+1+n$ and $\beta_m \colon 1+1+m\to m$
in $\CD$ as illustrated below.
            \begin{equation*}
            \alpha_n \df \lower10pt\hbox{$\includegraphics[height=1cm]{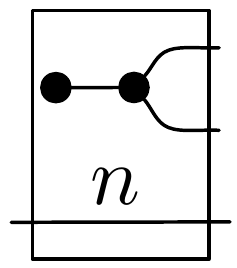}$} \qquad \qquad \beta_m \df \lower10pt\hbox{$\includegraphics[height=1cm]{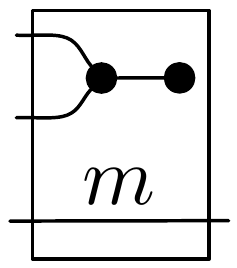}$}
            \end{equation*}
Next, we define the families of operators $\Rwl{n,m} \colon \CD[n+1, m] \to \CD[n, 1+m]$ and $\Rwr{n,m} \colon \CD[n, 1+m] \to \CD[1+n, m]$ as follows: for any circuit $c\in \CD[n+1, m]$,
\[
\Rwl{n,m}(c) \df \alpha_n \poi (id_1 \tns c)
\quad
\left(\lower8pt\hbox{$\includegraphics[height=.7cm]{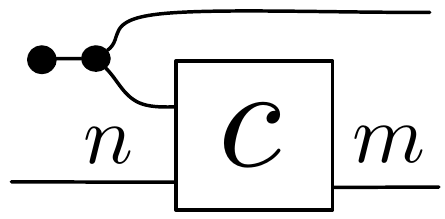}$}\right).
\]
and, for any circuit $d\in \CD[n, m+1]$
\[
\Rwr{n,m}(d) \df (id_1 \tns d) \poi \beta_m.
\quad
\left(\lower8pt\hbox{$\includegraphics[height=.7cm]{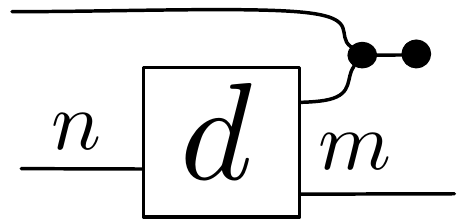}$}\right)
\]
\begin{remark}\label{rmk:rewiringadjunction}
When considered as operations on $\IBpoly$, $ \Rwl{n,m}$ and $\Rwr{n,m}$ enjoy some interesting properties.
Let $1+- \colon \IBpoly \to \IBpoly$ be the functor acting on objects as $k \mapsto 1+k$ and on arrows as $f \mapsto \id_1 \oplus f$.
This functor is self-adjoint: the unit and the counit are the $\alpha_n$ and $\beta_m$ defined as above.
The fact that $\IBpoly$ is a SMC implies naturality of $\alpha$ and $\beta$. They satisfy the triangle equalities by~\eqref{eq:gensnake}: 
\begin{equation*}
\lower14pt\hbox{$\includegraphics[height=1.3cm]{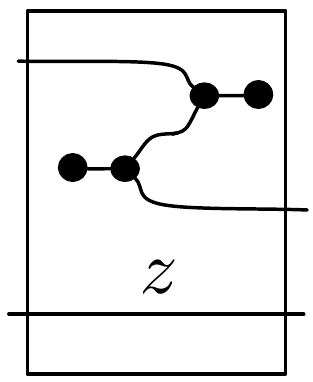}$}
\eqIH
\lower8pt\hbox{$\includegraphics[height=.7cm]{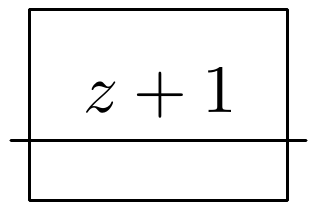}$}
\eqIH
\lower14pt\hbox{$\includegraphics[height=1.3cm]{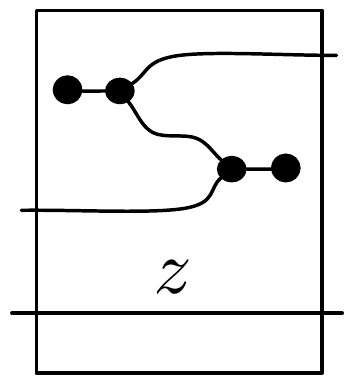}$}
\end{equation*}
The induced isomorphisms are
$$
\Rwl{n,m} \: \IBpoly[n+1, m] \to \IBpoly[n, 1+m] \qquad \qquad \Rwr{n,m} \: \IBpoly[n, 1+m] \to \IBpoly[1+n, m]
$$
defined as above.
We can see $\Rwl{n,m}$ intuitively as ``rewiring'' the first port on the left to the right of the circuit. The fact that $\Rwl{n,m}$ and $\Rwr{n,m}$ are isomorphisms means, of course, that no information is lost -- all such circuits can be ``rewired'' back to their original form.
\end{remark}

\begin{definition}\label{def:rewiring}
A circuit $c_2\in \CD[n_2, m_2]$ is a \emph{rewiring} of $c_1\in \CD[n_1, m_1]$ when $c_2$ can be obtained from $c_1$ by a combination of the following operations:
\begin{enumerate}[(i)]
\item application of $\Rwl{n,m}$, for some $n$ and $m$,
\item application of $\Rwr{n,m}$, for some $n$ and $m$,
\item post-composition with a permutation,
\item pre-composition with a permutation.
\end{enumerate}
\end{definition}

Permutations are needed to rewire an arbitrary---i.e.\ not merely the first---port on each of the boundaries. For instance, they allow to rewire the second port on the right as the third on the left in the circuit $c \: 2 \to 2$ below:
            \begin{center}
            \includegraphics[height=.8cm]{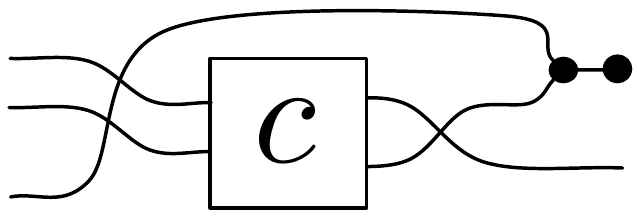}
            \end{center}

\begin{example}\label{ex:rewiring} We illustrate some rewirings of the leftmost circuit below.
\begin{equation*}
\includegraphics[height=.6cm]{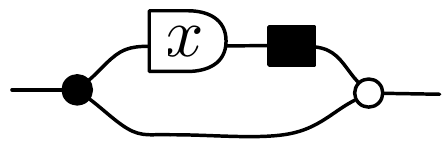} \qquad
\includegraphics[height=.7cm]{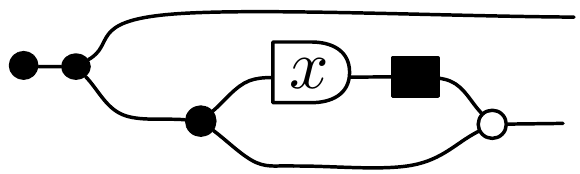} \qquad
\includegraphics[height=.7cm]{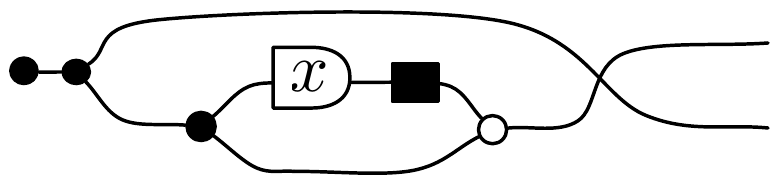} \qquad
\includegraphics[height=.9cm]{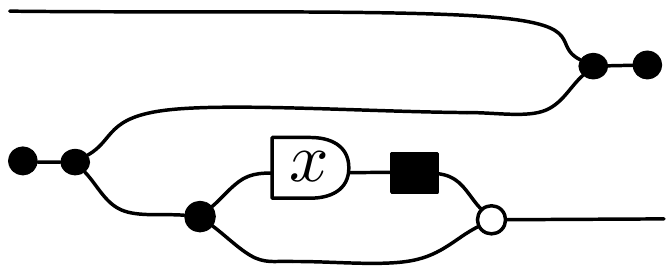} \qquad
\includegraphics[height=.8cm]{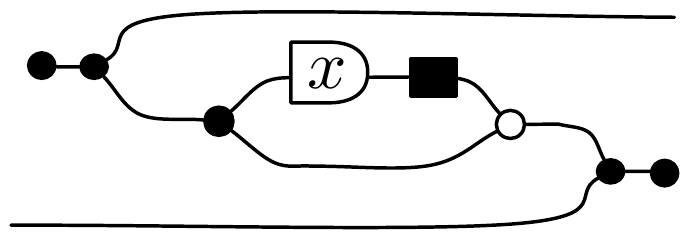} \qquad
\end{equation*}
\end{example}

%

In light of Remark~\ref{rmk:rewiringadjunction},
 ``is a rewiring of'' is an equivalence relation on the circuits of $\CD$ under the equational theory of $\IBpoly$: we
shall say that circuits $c$ and $d$ are \emph{rewiring-equivalent} when $c \eqIH d'$ for some rewiring $d'$ of $d$. For instance, all the circuits of Example~\ref{ex:rewiring} are rewiring-equivalent.


An interesting observation is that, at the semantics level, a rewiring $c_1 \: n_1 \to m_1$ of $c_2 \: n_2 \to m_2$ yields an isomorphisms between subspaces $\dsem{c_1} \in \SVpoly[n_1,n_1]$ and $\dsem{c_2} \in \SVpoly[n_2,m_2]$, with possibly different source and target but such that $n_1+m_1 = n_2+m_2$. We shall see at the end of the section how this property can give us information about the input/output ports of circuits.

\begin{lemma}\label{lemma:rewiringspaces}
If $c_1 \: n_1 \to n_1$ is a rewiring of $c_2 \: n_2 \to m_2$ in $\CD$, then $\dsem{c_1}$ and $\dsem{c_2}$ are isomorphic as subspaces of $\poly^{n_1+m_1} = \poly^{n_2+m_2}$.
\end{lemma}
\begin{proof}
It is enough to observe how $\Rwl{n,m}$, $\Rwr{n,m}$ and permutations
affect the denoted subspaces: \begin{enumerate}[(i)]
\item $\Rwl{n,m}$ induces an isomorphism between $\dsem{c} \in \SVpoly[n+1,m]$ and $\dsem{\Rwl{n,m}(c)} \in \SVpoly[n,m+1]$ defined by
\[
(\begin{footnotesize}\left(\begin{array}{c} \!\!q\!\! \\ \!\!\mathbf{v}\!\!\end{array}\right)\end{footnotesize},\mathbf{w})
\mapsto
(\mathbf{v},\begin{footnotesize}\left(\begin{array}{c} \!\!q\!\! \\ \!\!\mathbf{w}\!\!\end{array}\right)\end{footnotesize}).
\]
\item $\Rwr{n,m}$ induces an isomorphism between $\dsem{c} \in \SVpoly[n,m+1]$ and $\dsem{\Rwl{n,m}(c)}  \in \SVpoly[n+1,m]$ defined by
\[
(\mathbf{v},\begin{footnotesize}\left(\begin{array}{c} \!\!q\!\! \\ \!\!\mathbf{w}\!\!\end{array}\right)\end{footnotesize})
\mapsto
(\begin{footnotesize}\left(\begin{array}{c} \!\! q \!\!\\ \!\!\mathbf{v}\!\!\end{array}\right)\end{footnotesize},\mathbf{w}).
\]
\item post-composition with a permutation $p$ induces an isomorphism
$(\vv,\ww)\mapsto(\vv,\mathbf{w'})$ with $\mathbf{w'}$ obtained from $\ww$
by rearranging its rows according to $p$.
\item pre-composition with a permutation $p$ induces an isomorphism
$(\vv,\ww)\mapsto(\mathbf{v'},\ww)$ with $\mathbf{v'}$ obtained from $\vv$
by rearranging its rows according to $p^{-1}$.
\end{enumerate}

\end{proof}

For our purposes, it is important to also observe that rewiring does not affect the execution properties of circuits studied in \S\ref{sec:fullabstract}.

\begin{lemma}\label{lemma:rew}
Rewiring preserves deadlock and initialisation freedom.
\end{lemma}
\begin{proof}
 The argument goes by induction on the structure of the rewiring. It is immediate to see that the property of being deadlock and initialisation free is preserved by $\Rwl{n,m}$, $\Rwr{n,m}$ and (pre-post-)composition with permutations.
\end{proof}

We are now able to state the main result of this section.
\begin{theorem}\label{thm:realisability}
Every circuit in $\CD$ is rewiring-equivalent to some circuit in~$\SFGform$.
\end{theorem}


Before delving into the proof of Theorem~\ref{thm:realisability}, we illustrate an instance of its statement.
\begin{example}\label{ex:twoimplementations}
The circuit $\circuitUnoMinusX \in \CD[1,1]$ is rewiring-equivalent to the following circuit of~$\SFGform$.
 \begin{equation}\label{eq:extwoimpl1}
  \lower5pt\hbox{$\includegraphics[height=.6cm]{graffles/impl1minusxOne.pdf}$}
 \end{equation}
 The witnessing derivation uses first $\eqref{eq:scalarsum}^{\op}$ and then Lemma~\ref{prop:star=refl} to transform $\circuitUnoMinusX$ into a circuit which is a rewiring of \eqref{eq:extwoimpl1}.
  \begin{equation} \label{eq:extwoimpl2}
 \lower7pt\hbox{$\includegraphics[height=.8cm]{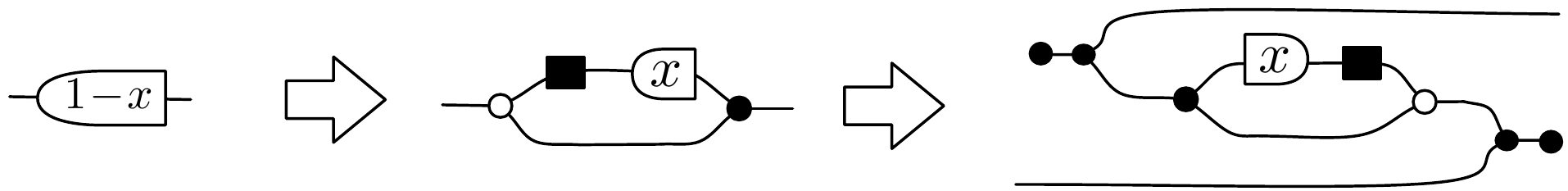}$}
\end{equation}
There is another circuit of $\SFGform$ which is rewiring-equivalent --- in fact, equal in $\IBpoly$ --- to $\circuitUnoMinusX$, resulting from the following derivation in $\IBpoly$:
 \begin{equation} \label{eq:extwoimpl3}
\lower15pt\hbox{$ \includegraphics[height=1.8cm]{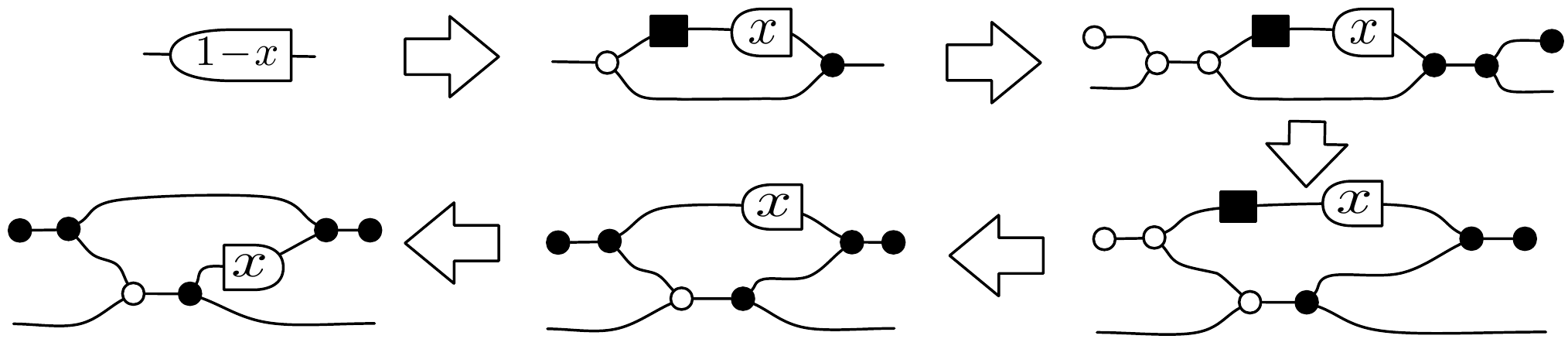}$}
\end{equation}
  Intuitively, the choice between the two different realisations of $\circuitUnoMinusX$ depends on whether one considers the input port to be on the right  --- \eqref{eq:extwoimpl2} --- or on the left  --- \eqref{eq:extwoimpl3}. We will be able to formalise this observation in Example~\ref{ex:implementations}, after that flow directionality has been explicitly introduced.

  As a concluding remark, note that, modulo the notation $\antipodesquare$ adopted for $\antipodeop$ since Example~\ref{ex_fibonacci}, the two circuits of Example~\ref{exm:opsem} also appear in \eqref{eq:extwoimpl2} and \eqref{eq:extwoimpl3}: by Corollary~\ref{cor:completeness}, we just proved that they have the same the semantics.
 \end{example}



We now turn to the proof of Theorem \ref{thm:realisability}, demonstrating how a procedure like the one of Example~\ref{ex:twoimplementations} is possible in general for any circuit of $\CD$. We shall work with matrices (and the corresponding circuits) of a particular shape. We say that a matrix over $\frpoly$ is in \emph{rational form} if all its entries are in fact rationals (in $\ratio$) and:
\begin{enumerate}
 \item for each non-zero row, there is a pivot entry with value $1$.
 \item in the column of a pivot, the pivot is the only non-zero entry.
\end{enumerate}
An example is given below, where $r_1,r_2,r_3\in \ratio$.
$$\scriptsize \left(%
\begin{array}{cccc}
  r_1  & 0 & 1 & 0\\
  r_2   & 1 & 0 & 0\\
  r_3   & 0 & 0 & 1\\
  0     & 0 & 0 & 0
\end{array}\right)$$

The following lemma is the final ingredient for the proof of Theorem~\ref{thm:realisability}---its proof, in Appendix~\ref{appendix:SFG}, is an easy exercise in linear algebra.
\begin{lemma}\label{lemma:rowequivalent}
Every $\frpoly$-matrix is row equivalent to one in rational form.
\end{lemma}
\begin{proof}[Proof of Theorem \ref{thm:realisability}]
Fix a circuit $c \in \CD[n,m]$. In the following, we give a recipe to transform $c$, using the equational theory of $\IBpoly$, into the rewiring of a circuit in $\SFGform$. To improve readability, we shall draw any circuit as if both $n$ and $m$ were $2$. It should be clear how our argument generalizes.
\begin{enumerate}[(i)]
    \item First we transform $c$ into the circuit $c_1$ on the right: the two are equal in $\IBpoly$ by \eqref{eq:gensnake}.
            \begin{equation*}
            \includegraphics[height=1.3cm]{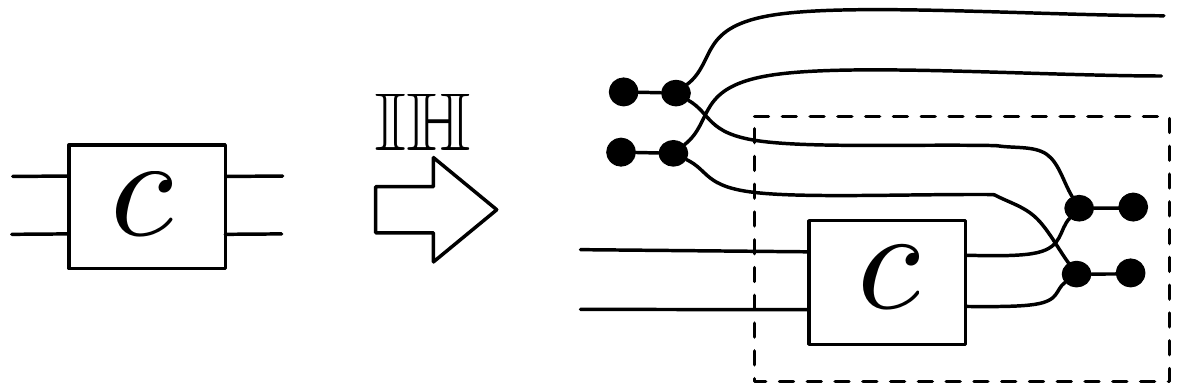}
            \end{equation*}
        Let us call $c_2$ the circuit $n+m \to 0$ delimited by the dotted square in the picture above. Since $c$ is a rewiring of $c_2$, it should be clear that, if $c_2$ can be rearranged as the rewiring of a circuit in $\SFGform$, then so can $c$. Therefore, in the sequel we shift our focus to $c_2$.
\item Theorem~\ref{Th:factIBR} allows us to rewrite $c_2$ in cospan form, as the composition along a middle boundary $z$ of $c_3$ and $c_4$ below, while preserving equality in $\IBpoly$. By definition of cospan form, $c_3$ is an arrow of $\FC$, while $c_4$ is an arrow of $\FCop$. For the sake of readability, we will draw $z$ as if it was $2$.
            \begin{equation*}
            \includegraphics[height=1.2cm]{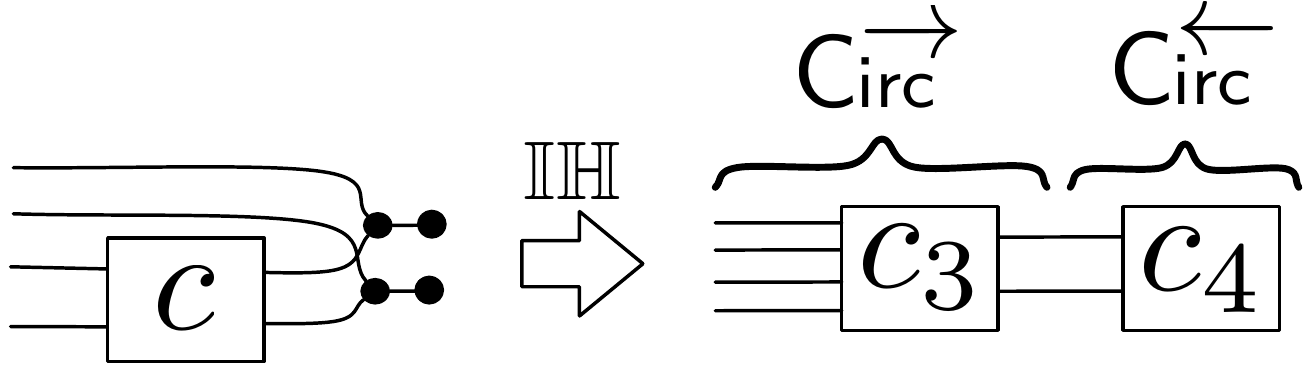}
            \end{equation*}
\item Since we are reasoning in $\IBpoly$, all the equations of $\ABpolyop$ hold. Now, $0$ is both the initial and the terminal object in $\Matpoly$; because $\ABpolyop \cong \Matpolyop$, this means that there is exactly one circuit of $\FCop$, up to equality in $\ABpolyop$, of type $z \to 0$. Therefore, $c_4$ and the $z$-fold monoidal product of $\Wcounit$ (a circuit that we call $c_5$) are equal in $\ABpolyop$ --- and thus in $\IBpoly$. We are then allowed to make the following rewriting:
            \begin{equation*}
            \includegraphics[height=.7cm]{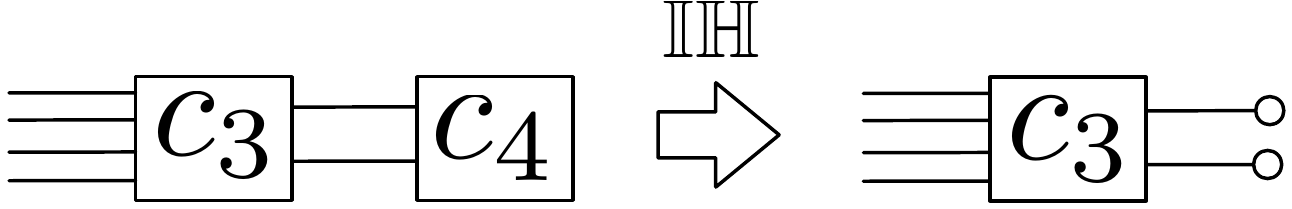}
            \end{equation*}
\item Since $c_3$ is in $\FC$, its polynomial semantics $\dsemHAO{} \: \FC \to \Matpoly$ gives us a $z \times (m+n)$ matrix $M$. This means that, when considered as a circuit of $\CD$, $c_3$ has semantics $\dsem{c_3} = \{(\vv, M \vv \mid \vv \in \frpoly^{m+n}\}$. Conversely, as discussed in \S~\ref{sec:theorymatr}, there is a canonical way of representing $M$ as a circuit $c_6$ of $\FC$ in matrix form, below right:

\[
M =\left(%
\begin{array}{cccc}
 \!\! p_{11}\!\! &\! p_{21}\!\! &\! p_{31}\!\! &\! p_{41} \!\!\\
 \!\! p_{12}\!\! &\! p_{22}\!\! &\! p_{32}\!\!& \!p_{42} \!\!
\end{array}\right)
\qquad \qquad
c_6 = \lower40pt\hbox{$\includegraphics[height=3.2cm]{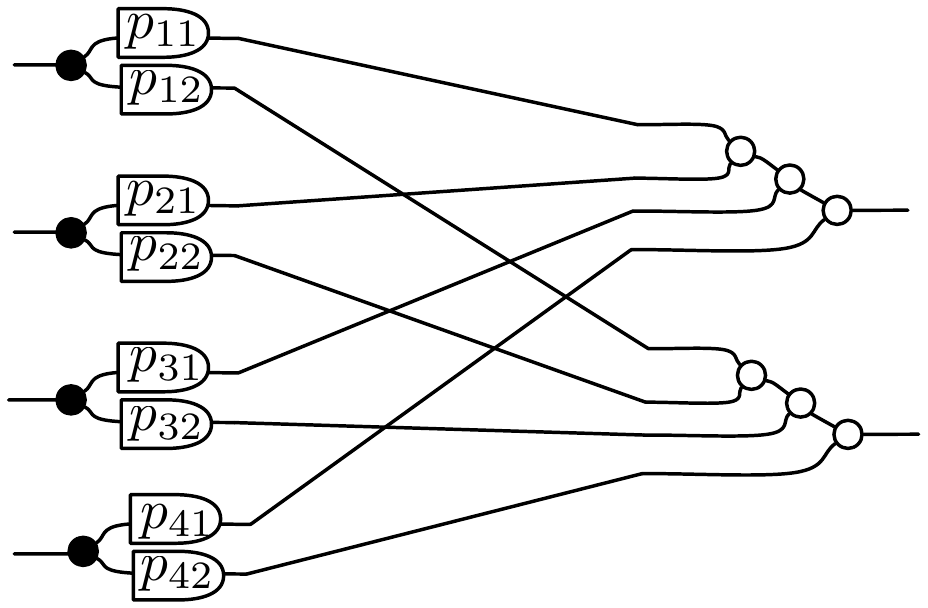}$}
\]
     We have that $\dsem{c_6} = \{(\vv, M \vv \mid \vv \in \frpoly^{m+n}\} = \dsem{c_3}$ and thus, by Corollary~\ref{prop:isoIHsoundcomplete}, $c_6 \eqIH c_3$. Therefore, we can rewrite our circuit $c_3 \poi c_5$ as follows:
 \begin{equation*}
   \includegraphics[height=2.8cm]{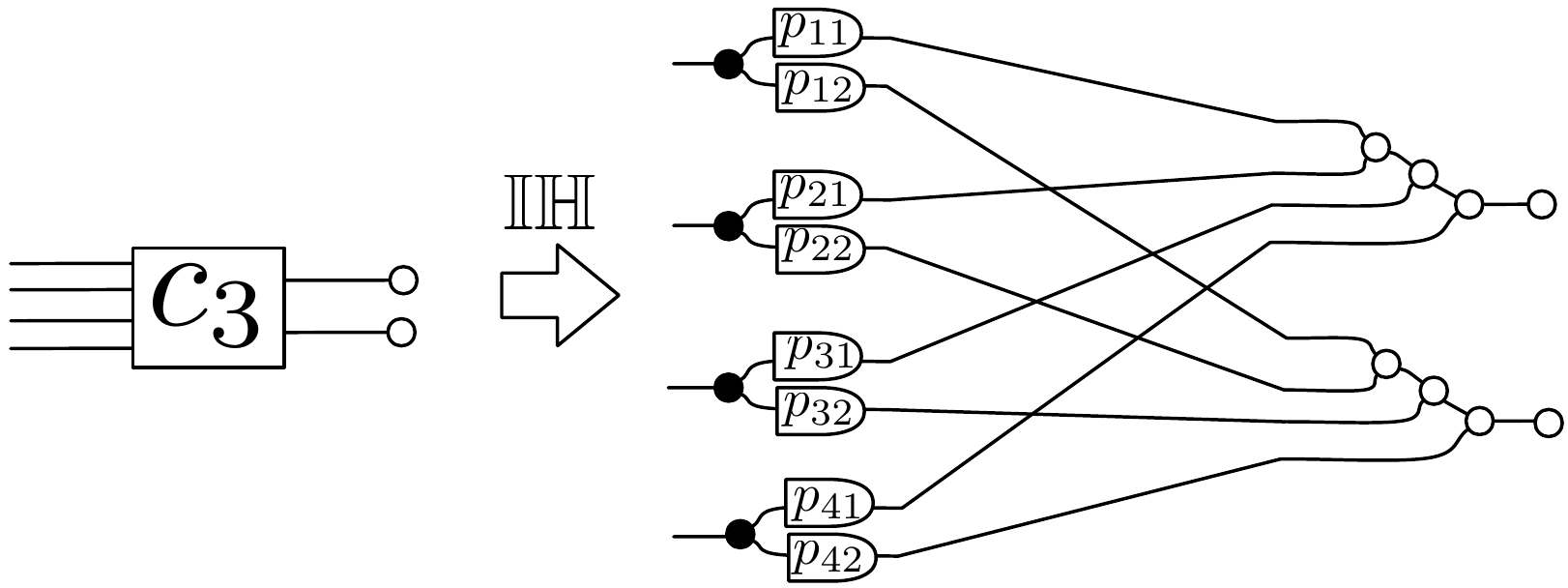}
 \end{equation*}
 \item Using Lemma \ref{lemma:rowequivalent}, we can then transform $M$ into a matrix $\widehat{M}$ in rational form --- for instance, the one on the left below. Since $\widehat{M}$ is a matrix over $\frpoly$, as observed in Remark~\ref{rmk:matrixformFractions}, there is a canonical circuit $c_7$ of $\IBpoly$, below right, representing it.

\[
\widehat{M} = \left(%
\begin{array}{cccc}
 \!\! 1\!\! &\! p_1 / q_1\!\! &\! 0\!\! &\! p_3 / q_3 \!\!\\
 \!\! 0\!\! & p_2 / q_2\! \!\! &\! 1\!\!& \! p_4 / q_4 \!\!
\end{array}\right)
\quad \qquad  \qquad \quad
c_7 = \lower30pt\hbox{$\includegraphics[height=2.2cm]{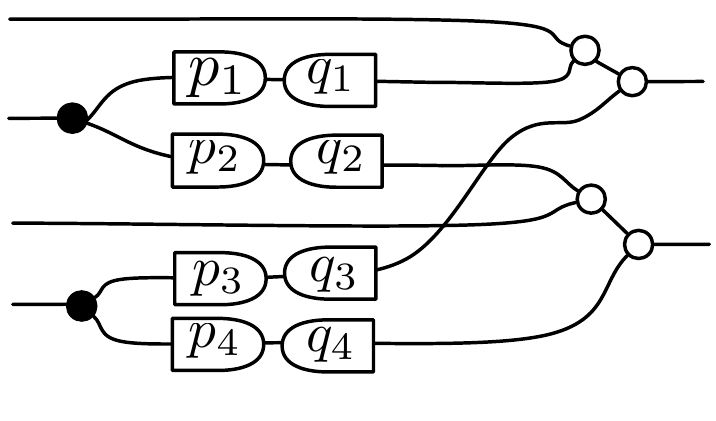}$}
\]
     By definition of rational form, each non-zero row $R$ in $\widehat{M}$ is associated with a pivot column $C$ with the only non-zero value $1$ at the intersection of $R$ and $C$. In order to graphically represent this property in $c_7$, we assume the following choice of pivots: the first and the third column for the first and second row respectively. Observe that an entry with value $0$ corresponds to the circuit $\zeroscalar$, which in $\IBpoly$ is equal to $\zeroscalarr$: therefore we can avoid drawing the corresponding link in the circuit $c_7$.

    We now claim that $c_6\poi c_5 \eqIH c_7 \poi c_5$. By Theorem \ref{th:IBR=SVR}, to check this it suffices to show that $\dsem{c_6\poi c_5} = \dsem{c_7 \poi c_5}$. We can compute:
    \begin{eqnarray*}
    \dsem{c_6\poi c_5} = \dsem{c_6}\poi\dsem{c_5} = \{(\vv, M \vv) \mid \vv \in \frpoly^{m+n}\} \poi \{ (\zerov , \matrixNull) \} = \{(\vv, \matrixNull) \mid M \vv = \zerov \} \\
    \dsem{c_7 \poi c_5} = \dsem{c_7}\poi\dsem{c_5} = \{(\vv, \widehat{M} \vv) \mid \vv \in \frpoly^{m+n}\} \poi \{ (\zerov , \matrixNull) \} = \{(\vv, \matrixNull) \mid \widehat{M} \vv = \zerov \}.
    \end{eqnarray*}
    Thus $\dsem{c_6\poi c_5} = \dsem{c_7 \poi c_5}$ amounts to the statement that $M$ and $\widehat{M}$ have the same kernel, which is true because they are row-equivalent. Therefore, $c_6\poi c_5 \eqIH c_7 \poi c_5$ and we can make the following rewriting:
\begin{equation*}
   \includegraphics[height=3cm]{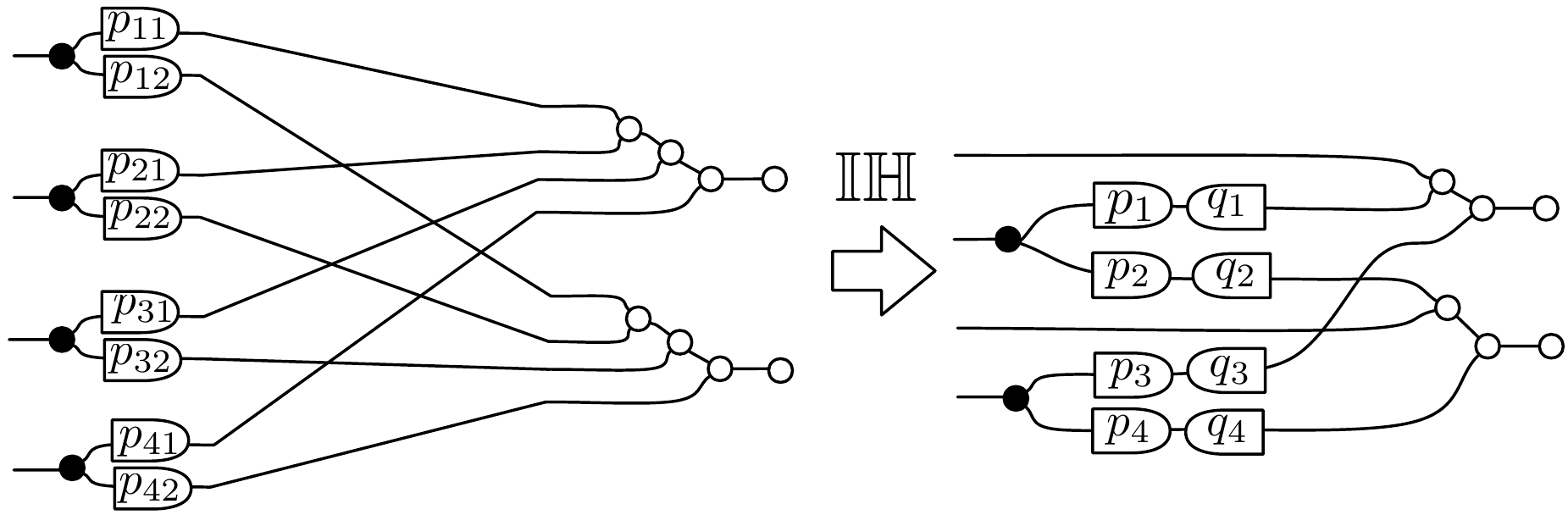}.
 \end{equation*}
 \begin{remark}\label{remark:graphicalProof}
  Instead of passing through the semantics, one can argue that $c_6\poi c_5 \eqIH c_7 \poi c_5$ in a more direct fashion by performing the linear algebraic manipulations involved in the proof of Lemma~\ref{lemma:rowequivalent} graphically. Indeed, similarly to what we did to prove Lemma~\ref{lemma:invertiblestar}, we can mimic the row operations used to transform $M$ into $\widehat{M}$ at the circuit level, using the equational theory of $\IBpoly$. This procedure involves a sequence of row-equivalent matrices $M_0,M_1,\dots,M_h$ represented by circuits $d_0,d_1,\dots,d_h$, where $M_0 = M$, $d_0 = c_6$ and $M_h = \widehat{M}$, $d_h = c_7$. At each step, two kinds of operation can be applied to $M_i$ in order to obtain $M_{i+1}$: the first is multiplying a row by an element $\frac{p_1}{p_2} \in \frpoly$, the second is replacing a row $R_1$ by $R_1 + \frac{p_1}{p_2} R_2$, where $R_2$ is another row. Bearing in mind that rows correspond to entries on the right boundary of $d_i$, the application of these two operations can be mimicked graphically as on the left and on the right below respectively.
  \begin{equation*}
\includegraphics[height=1cm]{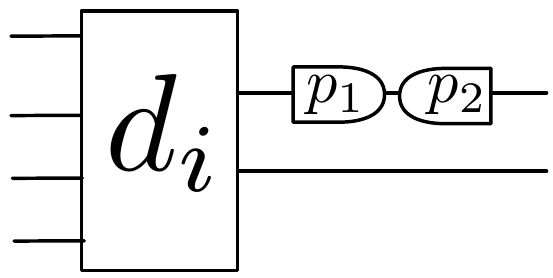}
\qquad \qquad
\includegraphics[height=1cm]{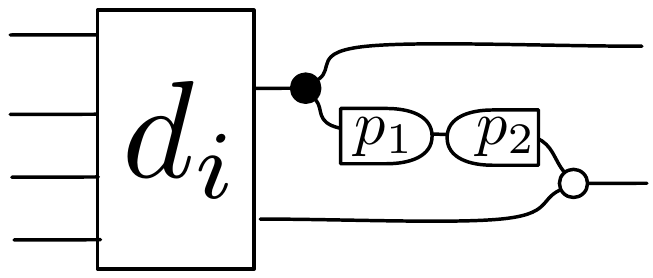}
\end{equation*}
On the left, we represent the first row being multiplied by $\frac{p_1}{p_2}$. On the right, we have the second row being summed with the first one multiplied by $\frac{p_1}{p_2}$: the semantics of $\Wmult$ and $\Bcomult$ confirm our description. Since these are row operations, the resulting circuit $d_{i+1}$ will still correspond to a matrix, namely $M_{i+1}$. An equational derivation can show that, modulo composition with $c_5$, the transformation of $d_i$ into $d_{i+1}$ is sound in $\IBpoly$:
 \begin{equation*}
\lower10pt\hbox{$\includegraphics[height=1cm]{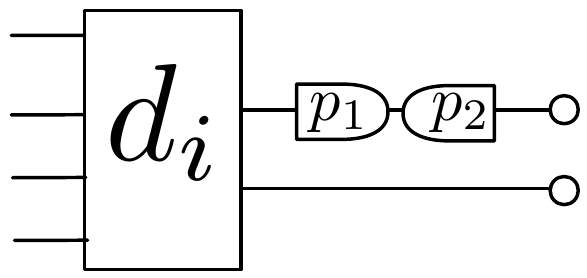}$}
\ \ \eqIH\ \
\lower10pt\hbox{$\includegraphics[height=1cm]{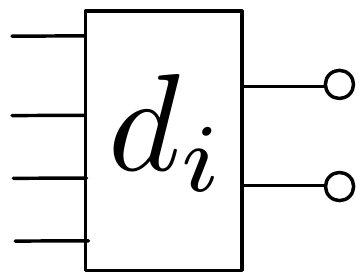}$}
\ \ \eqIH\ \
\lower10pt\hbox{$\includegraphics[height=1cm]{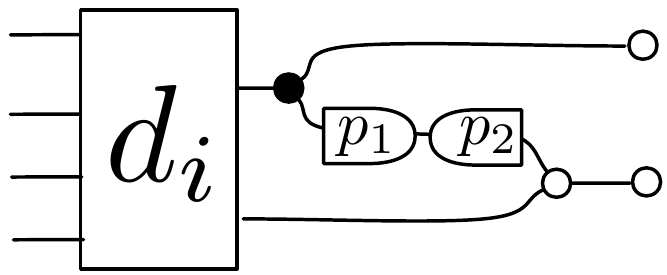}$}
\end{equation*}
\end{remark}
\end{enumerate}
 \noindent \paragraph{Proof of Theorem \ref{thm:realisability}, continued}
\begin{enumerate}[resume,label=({\roman*})]
\item We now focus on circuit $c_7 \poi c_5$. Our next step is to use associativity and commutativity of $\Wmult$ to make one of the two legs of each component $\Wccl$ be always attached to the pivot-wire of the corresponding row. Also, we use the laws~\eqref{eq:sliding1}-\eqref{eq:sliding2} of SMCs to push the pivot-wires towards the top of the circuit, as follows:
    \begin{equation*}
   \includegraphics[height=2.5cm]{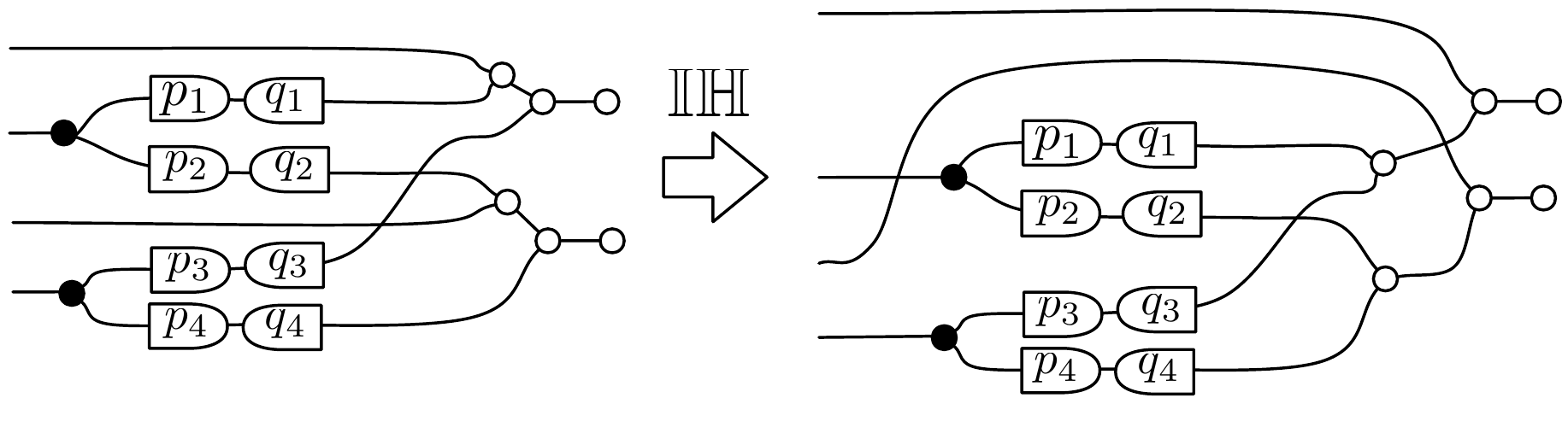}
 \end{equation*}
\item We can now remove the components of shape $\Wccl$ by turning them into rewiring structure. This can be done by using axiom $\eqref{eq:rcc}$ of $\IBpoly$:
    \begin{equation*}
   \includegraphics[height=2.5cm]{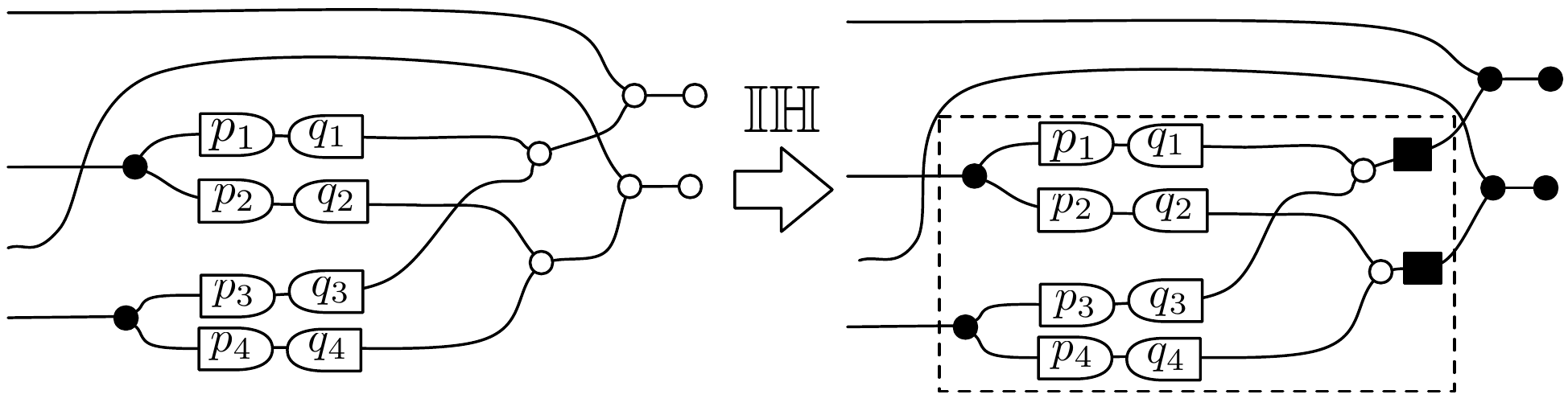}
 \end{equation*}
\item Let us call $c_8$ the right-hand circuit above: it is a rewiring of the circuit inscribed into the dotted square, which we call $c_9$. Since $c_7$ was constructed starting by a matrix in rational form, for all the components $\rationalcircuit$ in $c_8$, $\frac{p}{q}$ is a rational. Thus, using that $\SFG \cong \Matratio$ (Theorem~\ref{th:SFGcharactRationals}), we can rewrite in $\IBpoly$ each such component as a circuit $\tilde{c}$ in $\SFGform$:
    \begin{equation*}
   \includegraphics[height=2.5cm]{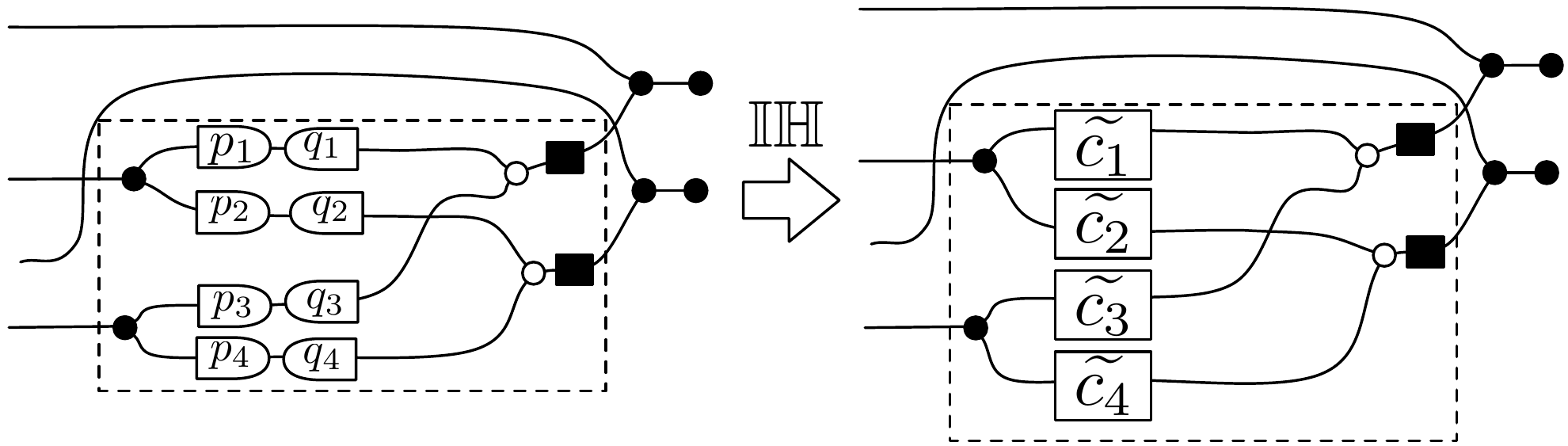}
 \end{equation*}
Now, observe that $c_9$ can be seen as the composition of circuits in $\SFGform$.
    \begin{equation*}
   \includegraphics[height=2.5cm]{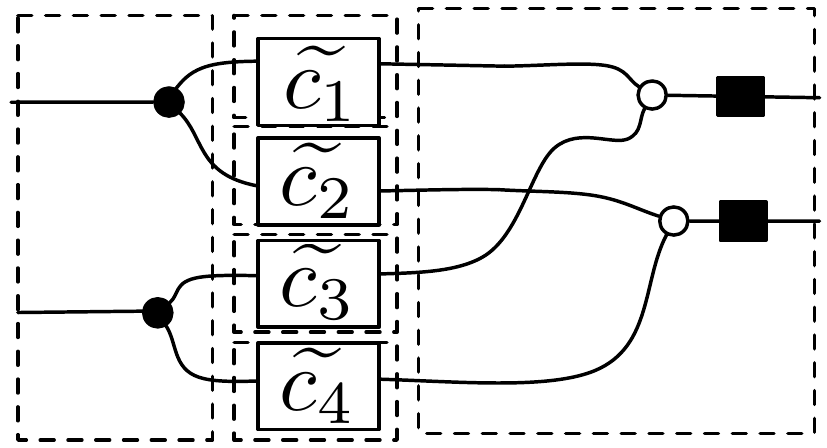}
 \end{equation*}
It follows that $c_9$ is also in $\SFGform$ and thus $c_8$ is the rewiring of a circuit in $\SFGform$. Since $c_8$ was obtained by $c_2$ by only using rewriting steps allowed by the equational theory of $\IBpoly$, the statement of the theorem follows.
\end{enumerate}\vspace{-.6cm}
\end{proof}

Combined with Proposition \ref{prop:SFGfree} and Lemma \ref{lemma:rew}, Theorem~\ref{thm:realisability} allows us to state the following result.

\begin{corollary}[Realisability]\label{cor:rea}
Every circuit $c$ of $\CD$ is equal in $\IBpoly$ to some circuit $d$ of $\CD$ that is deadlock and initialisation free.
\end{corollary}

Corollary~\ref{cor:rea} enables us to crystallise the main achievement of this section: any denoted behaviour $\strsem{c}$ can be properly \emph{realised} by a circuit $d$ for which, by full abstraction, the denotational and the operational perspective coincide. Moreover, $d$ being a rewiring of a signal flow graph, we can think of it as a properly \emph{executable} circuit --- we shall return to this point in \S\ref{sec:types} (Remark~\ref{rmk:executableOpsem}), when flow directionality is formally introduced in our theory. 

\medskip

We conclude this section with another observation stemming from Theorem~\ref{thm:realisability}, about the number of inputs of the realisation of a circuit.
\begin{proposition}\label{prop:dimension}
Let $c$ be a circuit in $\CD$ and suppose that $c'\: n\to m$ is a rewiring-equivalent circuit in $\SFGform$.
Then the dimension of $\dsem{c}$ is $n$.
\end{proposition}
\begin{proof}
By Lemma~\ref{lemma:rewiringspaces}, $\dsem{c}\cong\dsem{c'}$ as vector spaces. Since $c'$ is in $\SFGform$,
$\dsem{c'}$ is a functional subspace by Theorem~\ref{th:SFGcharactRationals}, whence its dimension is $n$.
\end{proof}
\begin{corollary}\label{corollary:dimension}
Let $c$ be a circuit in $\CD$ and suppose that $c_1\: n_1\to m_1$ and $c_2\: n_2\to m_2$, in $\SFGform$, are
rewiring-equivalent to $c$. Then $n_1=n_2$ and $m_1=m_2$.
\end{corollary}
Intuitively, for a circuit $c' \in \SFGform[n,m]$, inputs are always on the left and outputs are always on the right boundary. Although this partition may be different in a rewiring $c$ of $c'$, Proposition~\ref{prop:dimension} tells us that the information about the number of inputs of $c$ can be retrieved as the dimension of $\dsem{c}$. Moreover, Corollary \ref{corollary:dimension} guarantees that it remains constant for different realisations of $c$, despite of them exhibiting distinct input/ouput partitions (as in Example~\ref{ex:twoimplementations}).



\section{Directing the Flow}
\label{sec:types}

In the traditional presentation of signal flow graphs (see e.g.~\cite{mason1953feedback}), wires are directed, signifying the direction of signal flow. Throughout the previous sections, we have been referring to flow direction only on an intuitive level. We now introduce directionality explicitly, claiming that it can be really treated as a \emph{derivative} notion of our theory of circuits. We then present some applications and examples supporting our statement.

In order to model orthodox signal flow graphs we first need to introduce an alternative syntax, which we call the \emph{directed} signal flow calculus. We will need components that resemble those of $\FC$, but which are explicitly oriented from left to right.
\[
e\  :: = \
\lower8pt\hbox{$\includegraphics[height=.8cm]{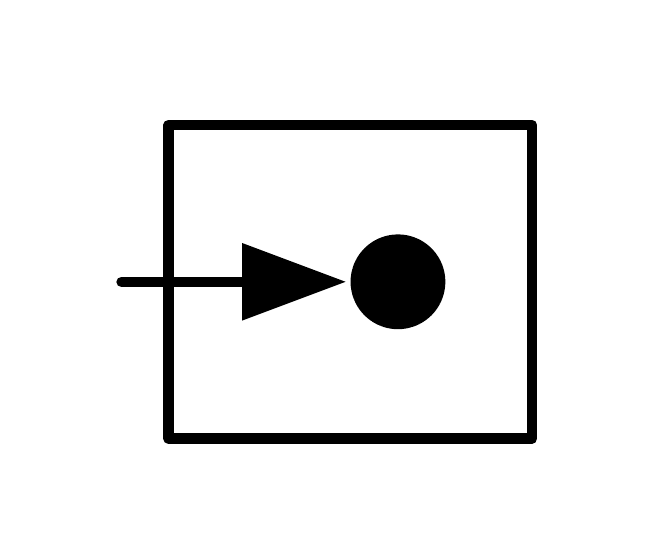}$}
\ |\
\lower8pt\hbox{$\includegraphics[height=.8cm]{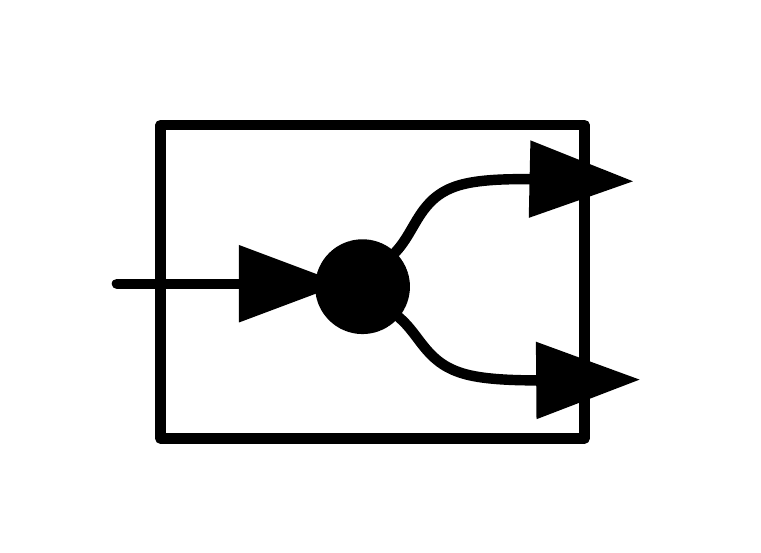}$}
\ |\
\lower8pt\hbox{$\includegraphics[height=.8cm]{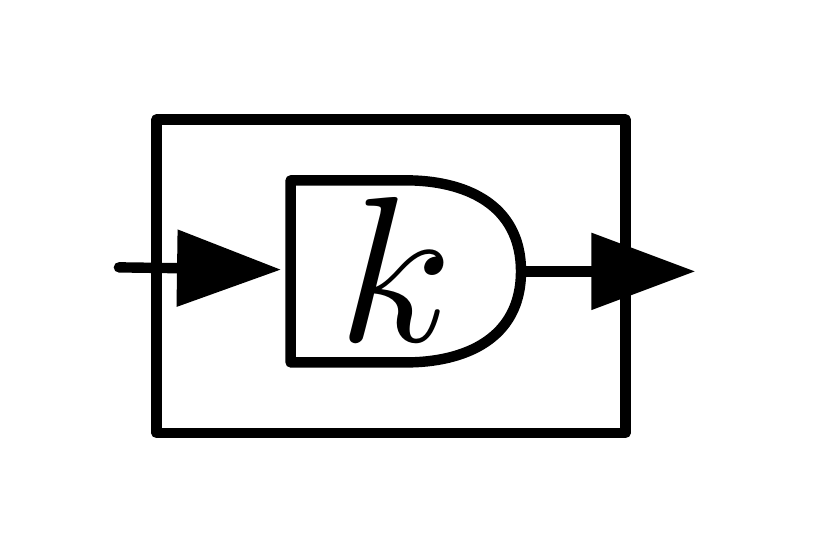}$}
\ |\
\lower8pt\hbox{$\includegraphics[height=.8cm]{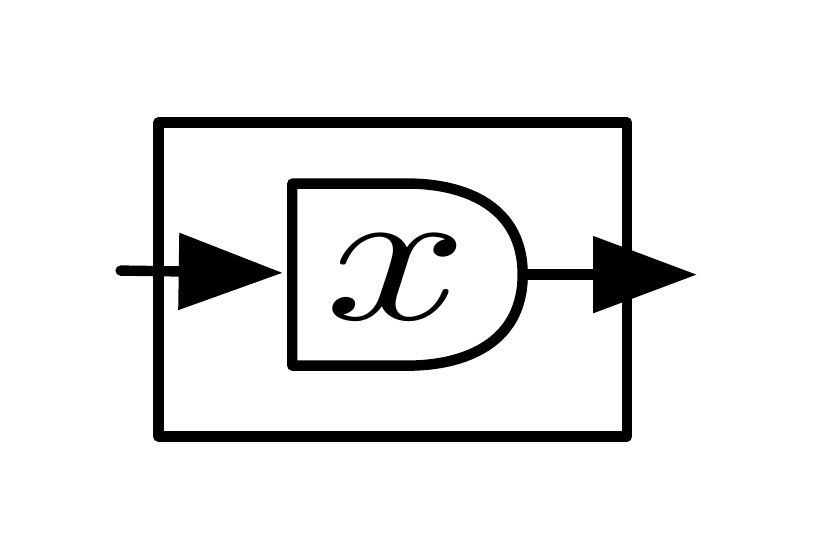}$}
\ |\
\lower8pt\hbox{$\includegraphics[height=.8cm]{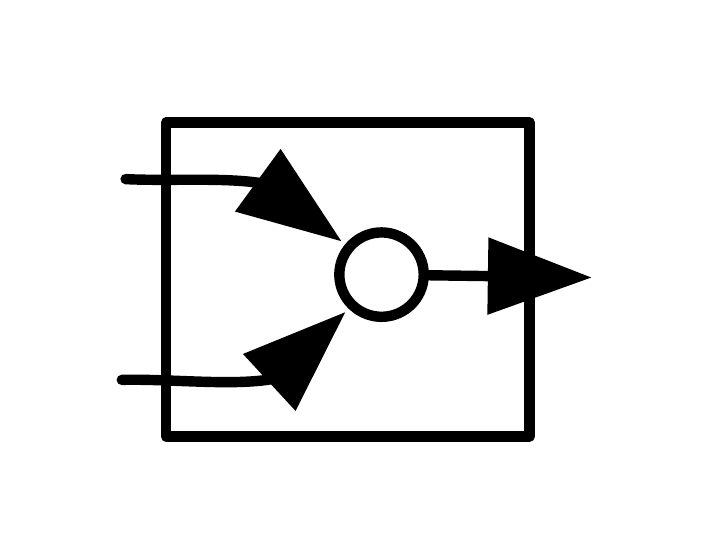}$}
\ |\
\lower8pt\hbox{$\includegraphics[height=.8cm]{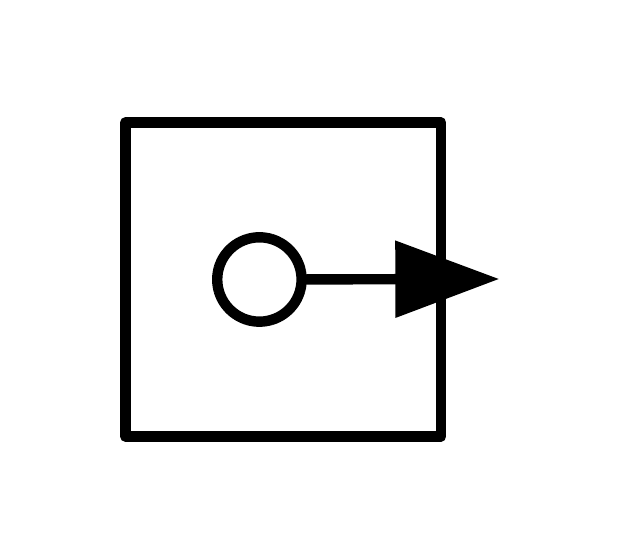}$}
\]
We also require some ``pure'' wiring: since signal flow is explicit, we include two versions of the identity wire and four of the twist:
\[
w \ ::= \
\lower8pt\hbox{$\includegraphics[height=.8cm]{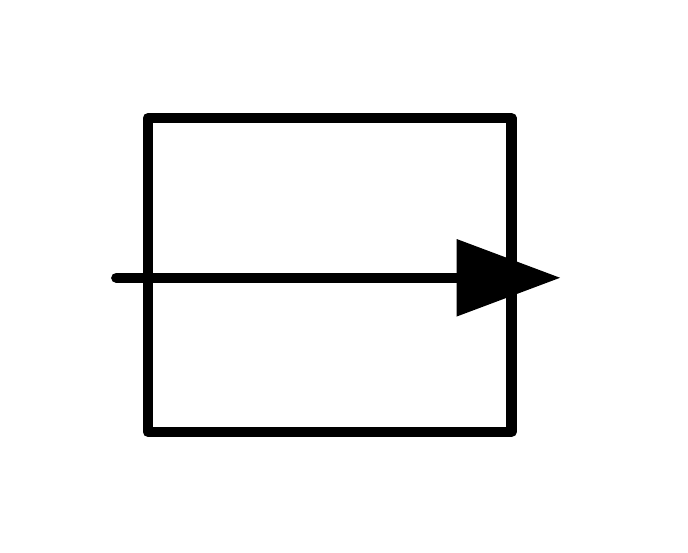}$}
\ |\
\lower8pt\hbox{$\includegraphics[height=.8cm]{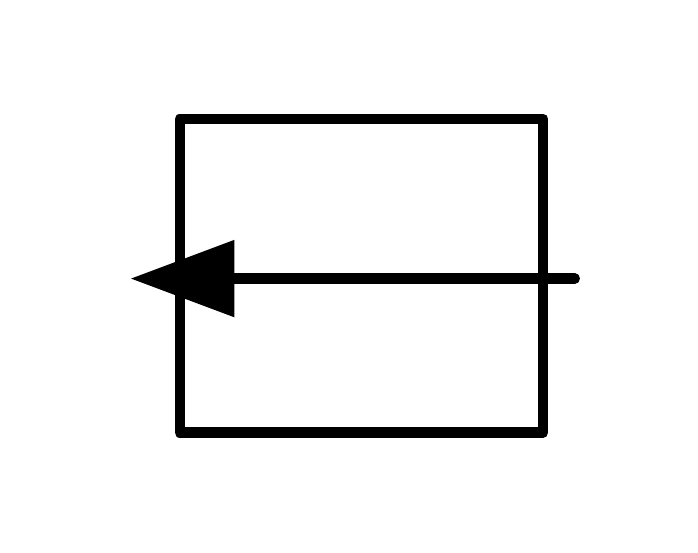}$}
\ |\
\lower8pt\hbox{$\includegraphics[height=.8cm]{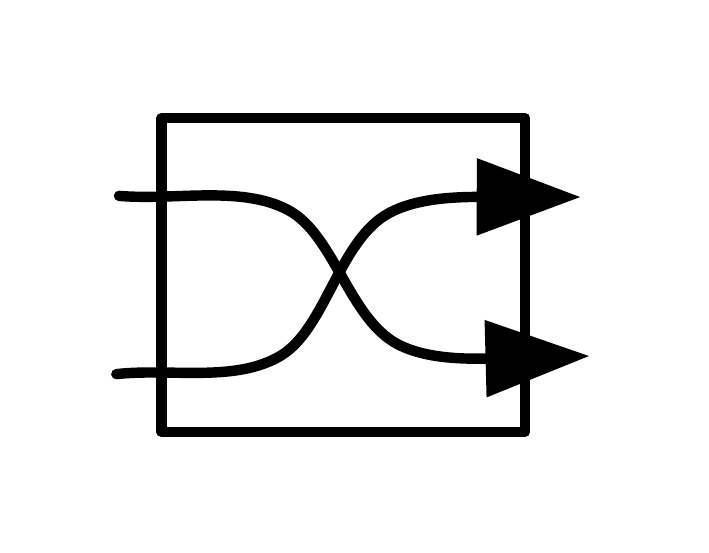}$}
\ |\
\lower8pt\hbox{$\includegraphics[height=.8cm]{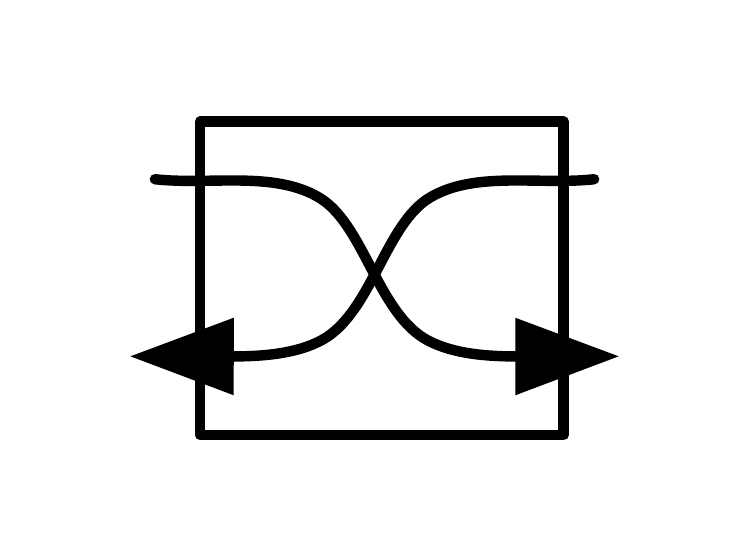}$}
\ |\
\lower8pt\hbox{$\includegraphics[height=.8cm]{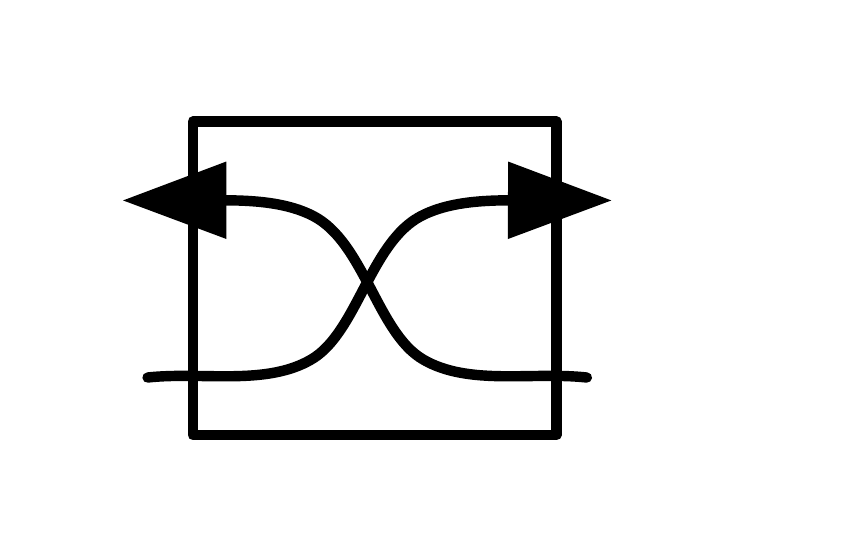}$}
\!\!\! |
\lower8pt\hbox{$\includegraphics[height=.8cm]{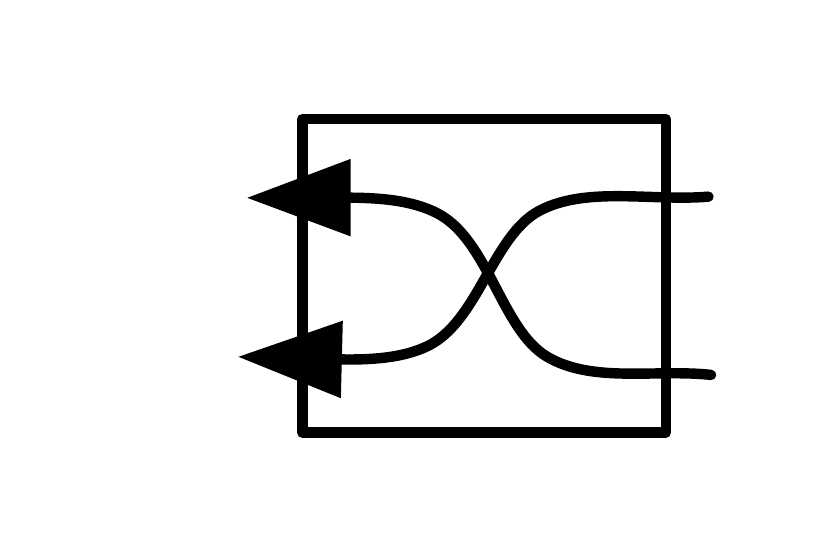}$}
\]
These basic components above are given a sorting $(u,v)$ where
 $u,v\in \{\LHD,\,\RHD\}^*$; for instance:
 \[
 \!\!\lower10pt\hbox{$\includegraphics[height=.9cm]{graffles/dBcomult}$}\!\!\!\!\mathrel{:}\sort{\RHD}{\RHD\RHD} \text{ and }
 \!\!\lower10pt\hbox{$\includegraphics[height=.9cm]{graffles/dtwisttwo}$}\!\!\!\!\!\!\mathrel{:}\sort{\LHD\RHD}{\RHD\LHD}.\]
Traditional signal flow graphs are obtained by composing components $e$ and $w$ using the operations $\poi$ and $\tns$, for which we reuse the sorting rules of Fig.~\ref{fig:sortInferenceRules}, together with guarded feedback operations $\Tr{\RHD}(\cdot)$ that take a circuit of sort
$(\RHD^{1+m},\RHD^{1+n})$ and yield a circuit of sort
$(\RHD^{m},\RHD^{n})$. The associated sorting rule is thus:
\[
\reductionRule{\typeJudgment{}{c}{\sort{\RHD^{1+n}}{\RHD^{1+m}}}}{ \typeJudgment{}{\Tr{\RHD}(c)}{\sort{\RHD^n}{\RHD^m}} }
\]
This is represented graphically as follows:
\[
\lower12pt\hbox{$\includegraphics[height=1cm]{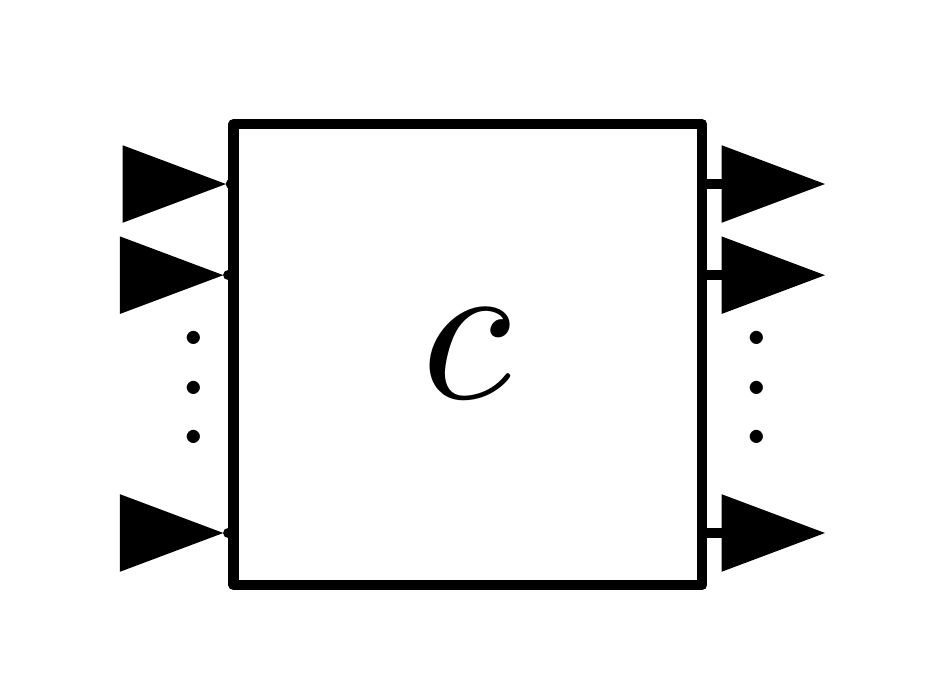}$}
\quad
\longmapsto
\quad
\lower15pt\hbox{$\includegraphics[height=1.3cm]{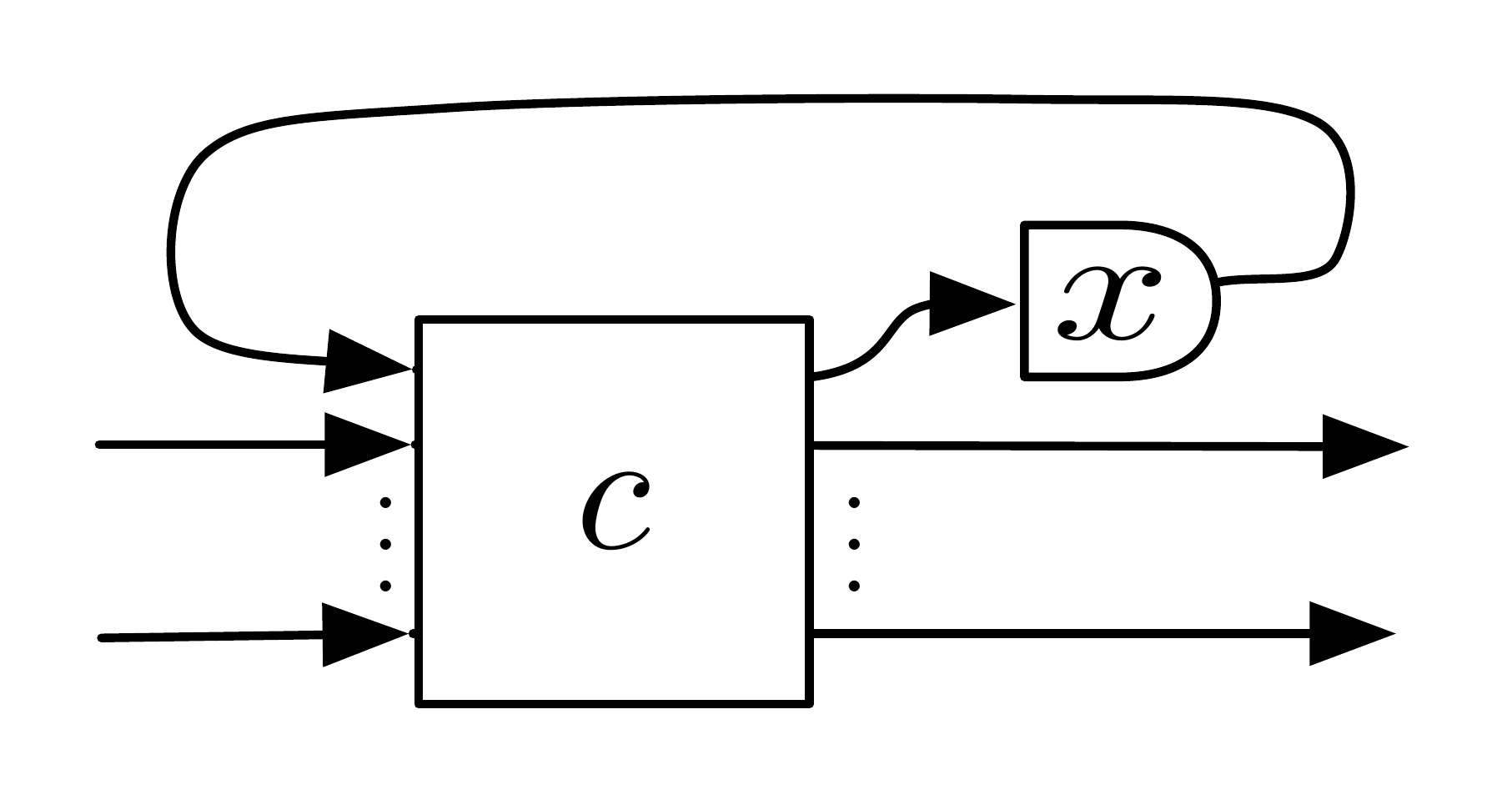}$}
\]
The syntax for directed signal flow graphs is thus:
\[
sf\ ::=\ e\ |\ w\ |\ sf\poi sf\ |\ sf\tns sf\ |\ \Tr{\RHD}(sf)
\]
Finally, we include top-level operations reminiscent of the rewiring in \S\ref{sec:realisability}: $L^{\RHD}$, $L^{\LHD}$, $R^{\RHD}$ and $R^{\LHD}$, with sorting rules:
\[
\reductionRule
{\typeJudgment{}{c}{\sort{\RHD u}{v}}}
{\typeJudgment{}{L^{\RHD}(c)}{\sort{u}{\LHD v}}}
\ \ \
\reductionRule
{\typeJudgment{}{c}{\sort{\LHD u}{v}}}
{\typeJudgment{}{L^{\LHD} (c)}{\sort{u}{\RHD v}}}
\ \ \
\reductionRule
{\typeJudgment{}{c}{\sort{u}{\RHD v}}}
{\typeJudgment{}{R^{\RHD}(c)}{\sort{\LHD u}{v}}}
\ \ \
\reductionRule
{\typeJudgment{}{c}{\sort{u}{\LHD v}}}
{\typeJudgment{}{R^{\LHD} (c)}{\sort{\RHD u}{v}}}
\]
In the graphical rendering below we leave out the arrowheads on wires where direction is arbitrary:
\[
\includegraphics[height=1cm]{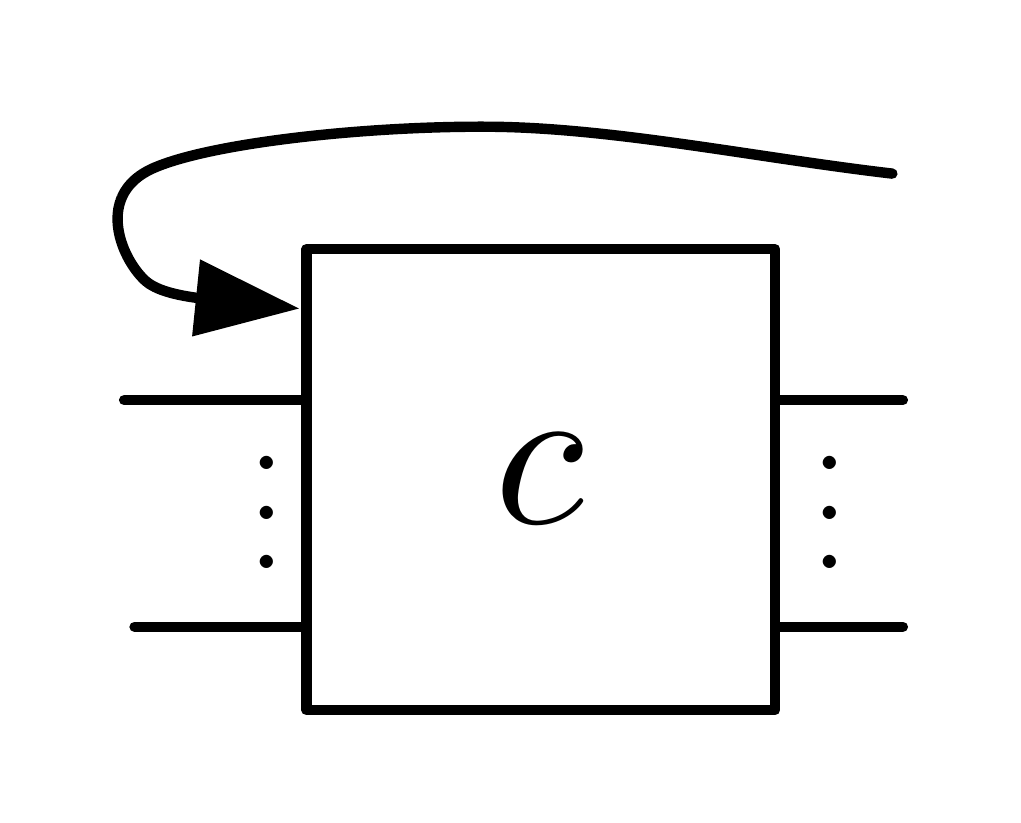}
\quad\quad\quad
\includegraphics[height=1cm]{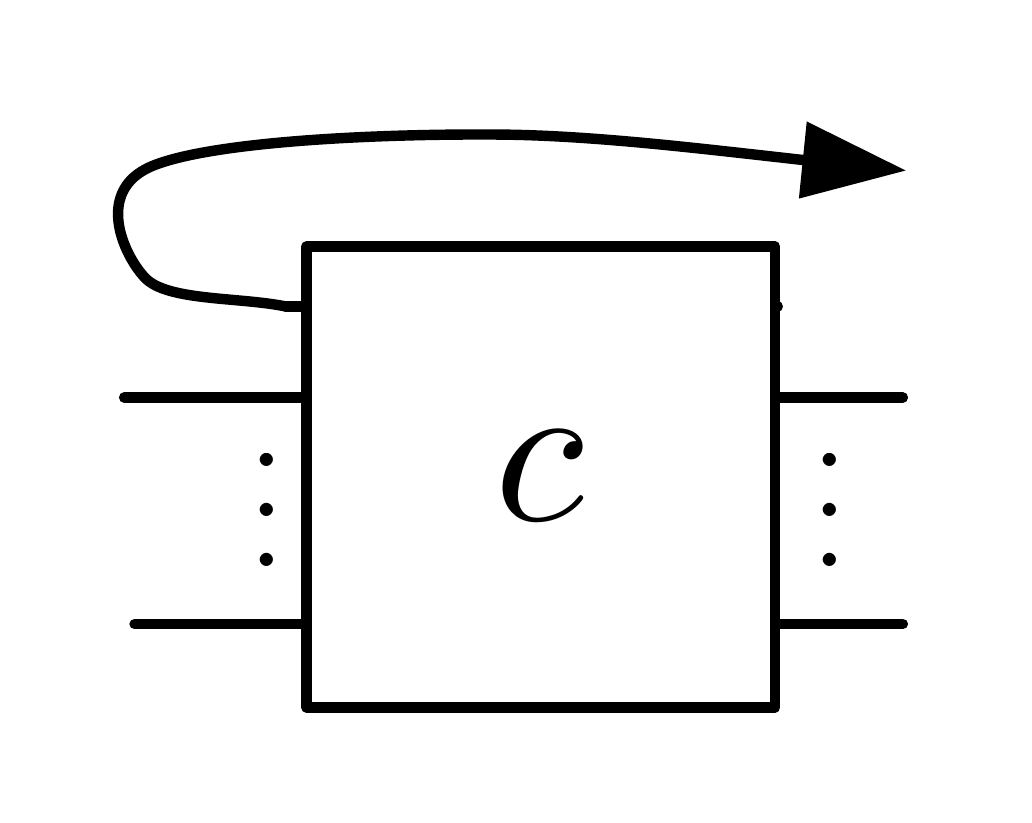}
\quad\quad\quad
\includegraphics[height=1cm]{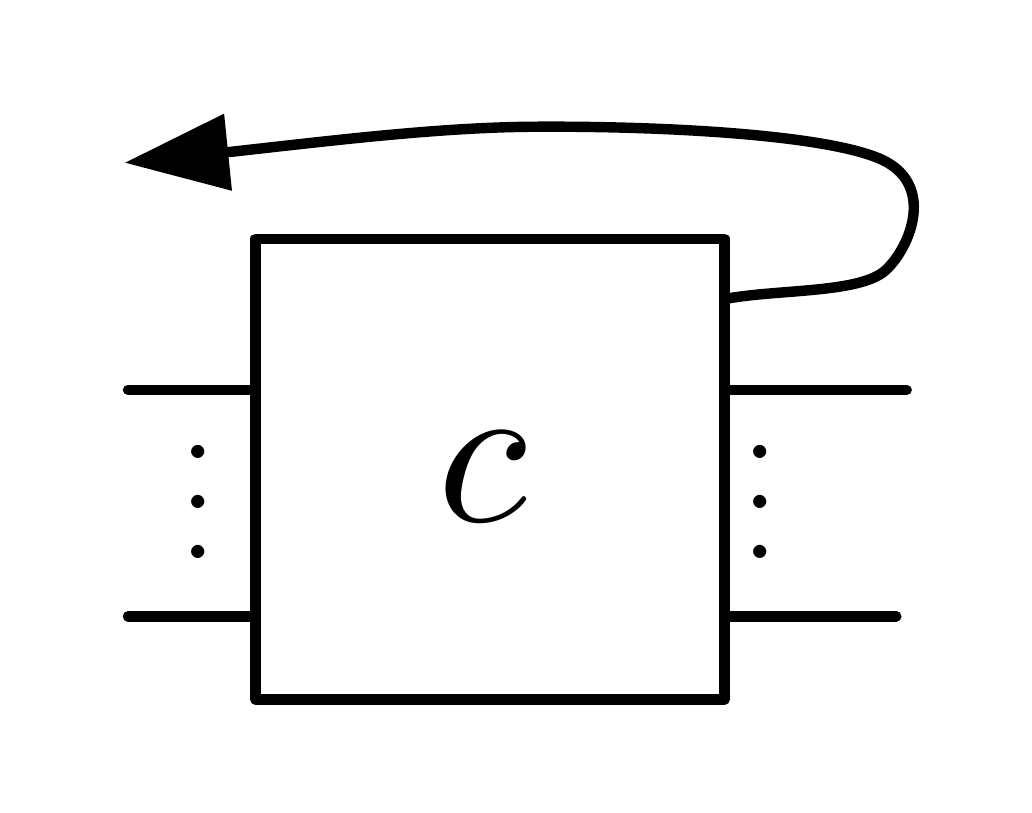}
\quad\quad\quad
\includegraphics[height=1cm]{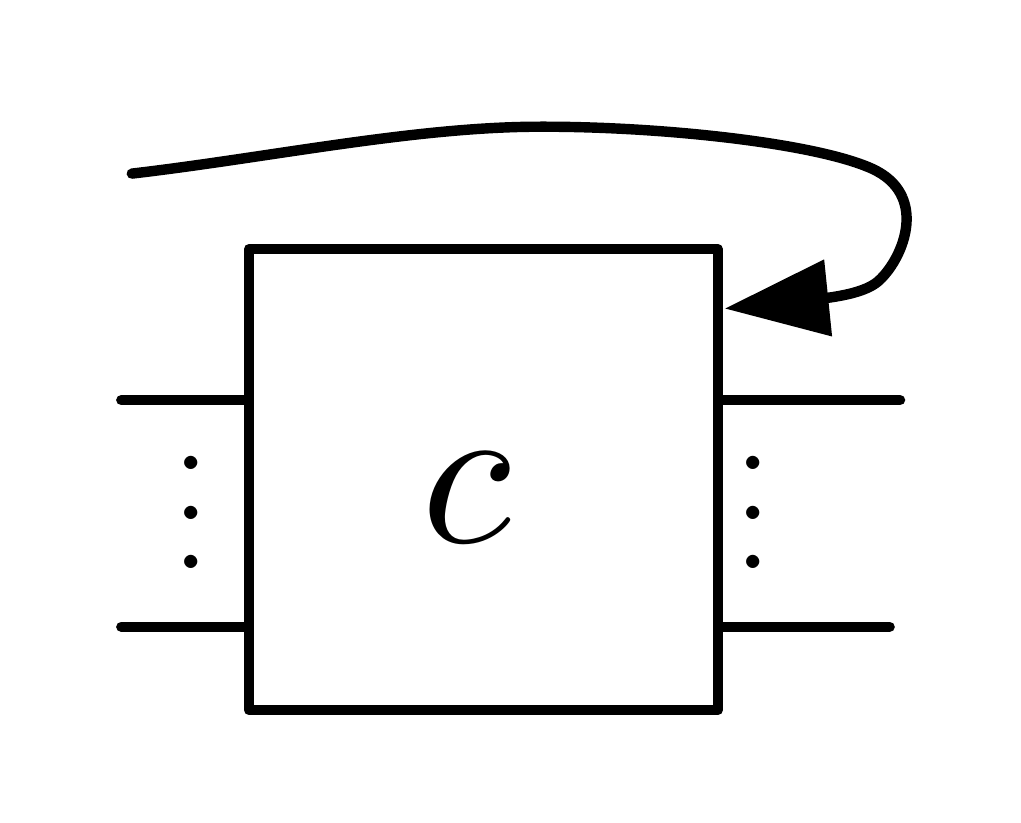}
\]
Circuits of the directed signal flow calculus are thus specified by the following grammar:
\[
d \ ::=\ sf\ |\ L^{\RHD}d\ |\ L^{\LHD} d\ |\ R^{\RHD}d\ |\ R^{\LHD}d\ |\ d \poi w\ |\ w\poi d
\]
Note that the composition at the top level is restricted to disallow the introduction of unguarded feedback.


Rather than defining the operational semantics directly, we can obtain the expected behaviour by first translating directed terms to the signal flow calculus. Intuitively, the inductively defined translation $\E$ ``erases directions'' from the wires:
\[
\!\!\!\lower8pt\hbox{$\includegraphics[height=.7cm]{graffles/dBcounit}$}
\!\!\mapsto \BcounitT,
\,
\lower8pt\hbox{$\includegraphics[height=.7cm]{graffles/dBcomult}$}
\!\!\mapsto \BcomultT
\cdots
\!\!
\lower7pt\hbox{$\includegraphics[height=.7cm]{graffles/dtwisttwo}$}
\!\!\!\!\mapsto \symNetT,
\,
\!\!
\lower7pt\hbox{$\includegraphics[height=.7cm]{graffles/dtwistop}$}
\!\!\mapsto \symNetT,
\]
\[
sf_1 \poi sd_2 \mapsto \E(sf_1) \poi \E(sf_2),\quad
sf_1 \tns sf_2 \mapsto \E(sd_1) \tns \E(sf_2),
\]
\[
\Tr{\RHD}(sf) \mapsto \Tr{}(\E(sf)),
\quad
\mathsf{L}^{\star}(d) \mapsto \Rwl{}(\E(d))
\quad
\mathsf{R}^{\star}(d) \mapsto \Rwr{}(\E(d)).
\]
where $\star\in\{\LHD,\RHD\}$ and $\Tr{}{}$, $\Rwl{}{}$ and $\Rwr{}{}$ are defined as in \S\ref{sec:SFcalculus} and \S\ref{sec:realisability}. 


\begin{remark}\label{rmk:executableOpsem} For any circuit $c$ in the image of $\E$, the rules of the operational semantics (Fig.~\ref{fig:operationalSemantics}) can be really thought as describing the step-by-step \emph{execution} of a state machine. Indeed, for any transition $c_1 \dtrans{\xx}{\yy} c_2$ between $c$-states, directed sort discipline on the formation of $c$ allows to identify certain values in vectors $\xx$, $\yy$ as \emph{inputs}, and all the others as \emph{outputs} which can be computed from the inputs and the state $c_1$. In particular this solves, for circuits in the image of $\E$, the unbounded non-determinism implicit in the operational rule for sequential composition that we observed in Remark~\ref{rmk:signalflowOpSem}. 
\end{remark}

A related observation is that directed sort discipline prevents us from writing problematic circuits where signal flow is incompatible, like in the examples in \S\ref{sec:fullabstract}. Thus, using Proposition~\ref{prop:SFGfree} and Lemma~\ref{lemma:rew}, we are able to state:
\begin{proposition}\label{pro:directedopsem}
For any circuit $d$ of the directed signal flow calculus, $\E(d)$ is deadlock and initialisation free.
\end{proposition}
\begin{proof}
The proof is by induction on the structure of terms. One first proves that for any directed signal flow graph $sf$, $\E(sf)$ is in $\SFGform$. Then, that for any directed circuit $d$, $\E(d)$ is a rewiring (Def.~\ref{def:rewiring}) of a circuit in $\SFGform$.
Proposition~\ref{prop:SFGfree} and Lemma~\ref{lemma:rew} conclude the proof.
\end{proof}

Proposition~\ref{pro:directedopsem} explains why we do not pursue an axiomatisation for the operational semantics: whenever it is executable --- i.e., for circuits in the image of $\E$, it is also fully abstract. Thus it agrees with the denotational semantics and we can reason about it using the equations of $\IBpoly$.

Moreover, we emphasise that the syntactic restrictions on the directed calculus do not affect the expressiveness since, thanks to Theorem \ref{thm:realisability}, rewirings of signal flow graphs denote all the possible behaviours. Thus, informally speaking, all circuits in $\CD$ can be directed (modulo~$\IBpoly$).

\begin{proposition}\label{pro:directedexpressivity}
For any circuit $c$ of $\CD$, there exists a directed circuit $d$ such that $\E(d)\eqIH c$
\end{proposition}
\begin{proof}
 Let $c$ be a circuit in $\CD$. Then, by Theorem~\ref{thm:realisability} there exists a $c'\eqIH c$ that is the rewiring of some circuit in $\SFGform$. By induction on the structure of $c'$, one can easily check that there exists a directed circuit $d$ such that $\E(d)=c'$.
\end{proof}

Propositions~\ref{pro:directedopsem} and \ref{pro:directedexpressivity} have two interesting consequences.
First, Proposition~\ref{pro:directedopsem} and the full abstraction result mean that we can use the equational theory of $\IBpoly$ to safely reason about orthodox signal flow graphs and extensions---indeed, all the circuits in the directed signal flow calculus. Roughly speaking, the procedure is: forget the directions and then apply any rewriting within $\IBpoly$. This confirms the intuition that, like for electrical circuits, also for signal flow graphs directionality is \emph{not} a primitive notion as originally advocated in \cite{mason1953feedback}. 

Second, Proposition~\ref{pro:directedopsem}, Proposition~\ref{pro:directedexpressivity} and full abstraction tell us that the denotational semantics of any circuit of the signal flow calculus can be properly realised by some directed circuit. We can therefore use the ``more liberal'' signal flow calculus to specify circuits and the ``more restrictive'' directed calculus to implement them. One can then check that an implementation $d$ adheres to a specification $c$ by mean of the graphical reasoning supported by $\IBpoly$. Indeed $\E(d)\eqIH c$, means that $d$ implements, without deadlocks or initialisation, the behaviour denoted by $c$. Note that, while an implementation is a directed circuit---typically featuring feedbacks---we are being deliberately vague about what kind of circuit in $\CD$ constitutes a specification: in examples that we consider these are typically generating functions that can be obtained in a standard way (see e.g.~\cite{Wilf2006}) from recurrence formulas, like the Fibonacci function of Example~\ref{ex_fibonacci}. We illustrate these ideas by resuming Example~\ref{ex:twoimplementations} below. 

\begin{example}\label{ex:implementations}
The circuit $\circuitUnoMinusX$, studied in Example~\ref{ex:twoimplementations}, can be seen as the string diagrammatic specification of the generating function $\frac{1}{1-x}$, which yields the constant stream $1,1,1,1,\dots$. Indeed,
\[
\strsem{\lower3pt\hbox{$\includegraphics[height=.4cm]{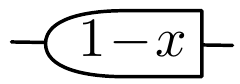}$}} = [(1-x, 1) ] = [(1, \frac{1}{1-x}) ]
\]
We now claim that the following directed circuits realise the specification.
\begin{equation}\label{eq:directedcircuits}
 \lower9pt\hbox{$\includegraphics[height=.9cm]{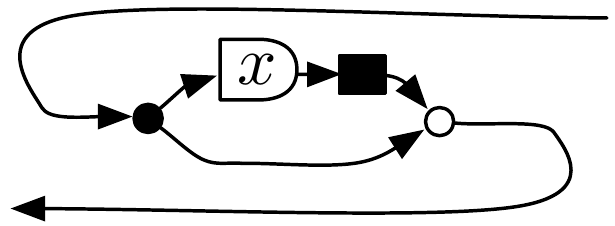}$} \qquad \qquad
 \lower9pt\hbox{$\includegraphics[height=.9cm]{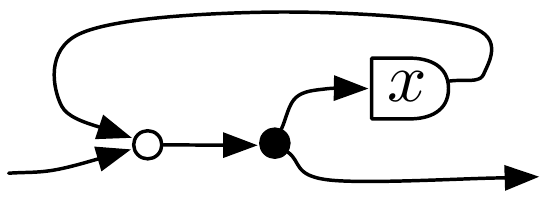}$}
 \end{equation}
To prove it, we first throw away all the directions from the wires, transforming the directed circuits into circuits of $\CD$. Then we need to show that the resulting circuits are equal to the specification in $\IBpoly$:
\begin{equation*}
\lower7pt\hbox{$\includegraphics[height=.8cm]{graffles/rewiring5.pdf}$}
\qquad \eqIH \qquad
\lower3pt\hbox{$\includegraphics[height=.4cm]{graffles/circuit1-xPicture.pdf}$}
\qquad \eqIH \qquad
\lower5pt\hbox{$\includegraphics[height=.6cm]{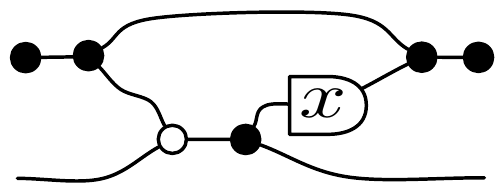}$}
\end{equation*}
The graphical derivations witnessing this statement are~\eqref{eq:extwoimpl2} and \eqref{eq:extwoimpl3}. The fact that the two implementations have flow directionality as in \eqref{eq:directedcircuits} substantiate our claim, made in Example~\ref{ex:twoimplementations}, that they yield different input/ouput partitions of the specification $\circuitUnoMinusX$.
\end{example}

\begin{example} A similar procedure can be used to check the observational equivalence of directed signal flow graphs. For instance, take:
\begin{eqnarray}\label{impl}
\lower10pt\hbox{$\includegraphics[height=1cm,width=3.4cm]{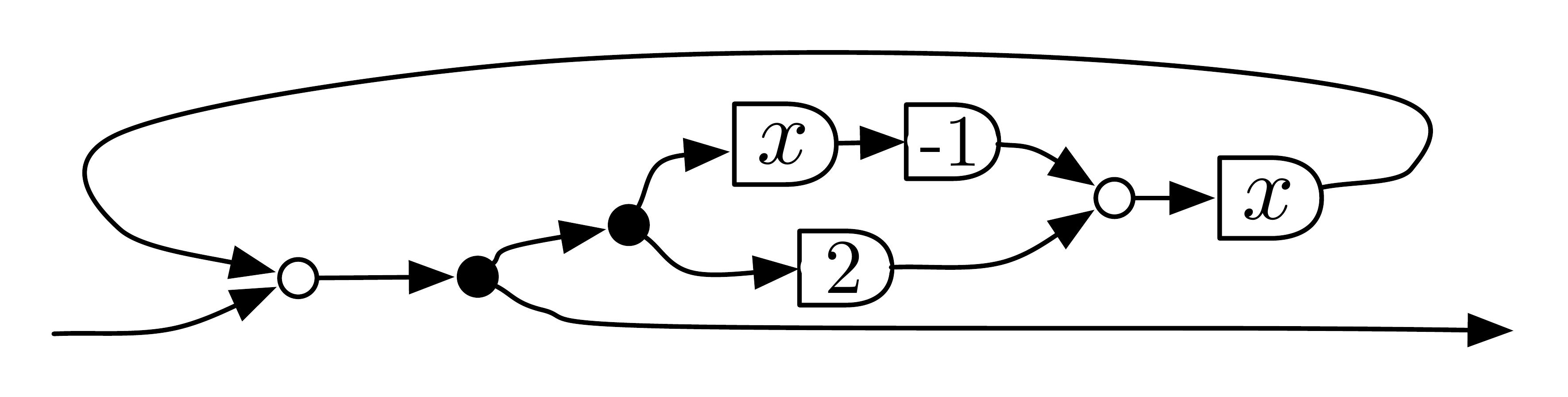} \qquad\qquad
\includegraphics[height=1cm,width=4.4cm]{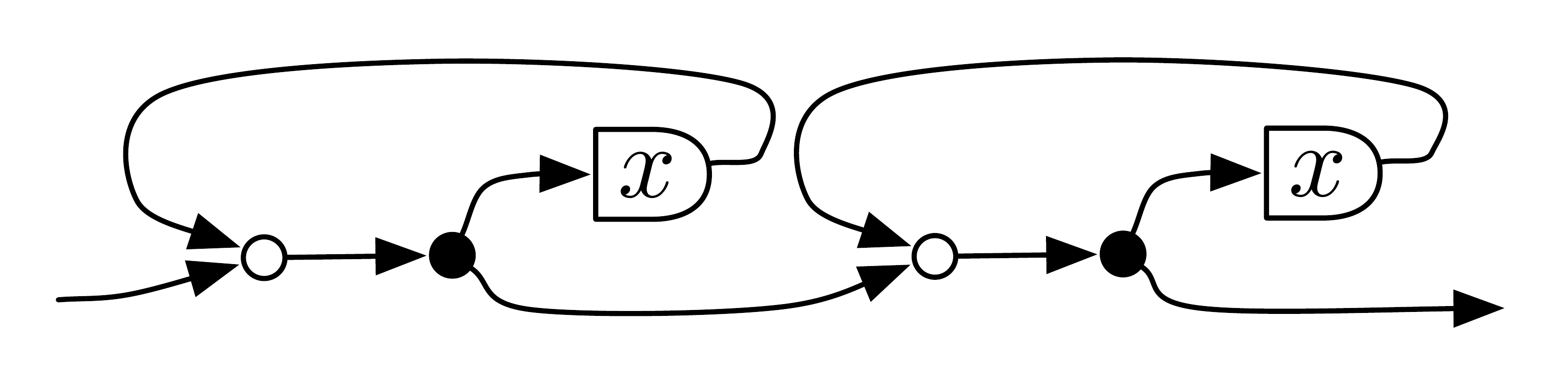}$}
\end{eqnarray}
First, we forget the direction of the flow and we obtain the circuits $c_3$ and $c_4$ depicted below.
\[
c_3 \df \lower8pt\hbox{$\includegraphics[height=.8cm,width=3.6cm]{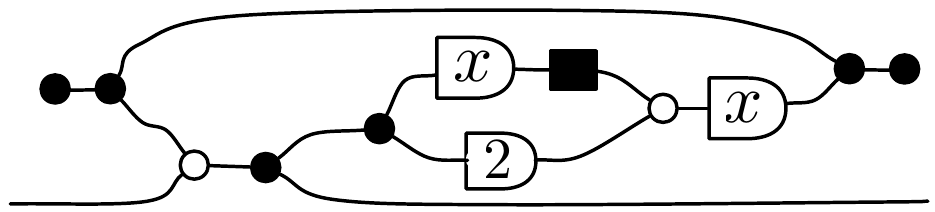}$} \qquad\qquad
c_4 \df \lower11pt\hbox{$\includegraphics[height=1cm,width=4.6cm]{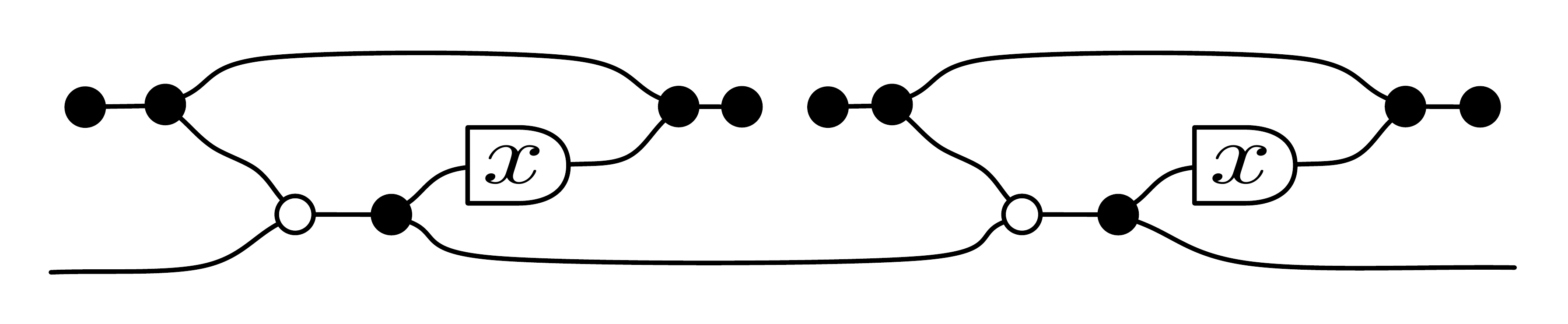}$}
\]
Then, by virtue of Proposition~\ref{pro:directedopsem} and full abstraction, we can safely use $\eqIH$ to check $\osem{c_3}=\osem{c_4}$.
Observe that $c_3$ is like in Example~\ref{exm:opsem}. For $c_4$, we can immediately derive that
\[\lower9pt\hbox{$\includegraphics[height=.9cm,width=4.4cm]{graffles/123456a.pdf}$}\eql{\eqref{eq:extwoimpl3}} \lower3pt\hbox{$\includegraphics[height=.4cm]{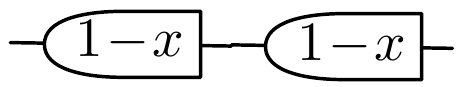}$}\eql{\eqref{eq:scalarmult}$^{\op}$} \lower4pt\hbox{$\includegraphics[height=.5cm]{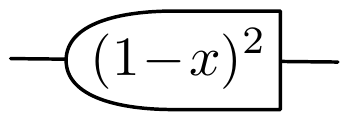}$}.\]

To conclude that $c_3 \eqIH c_4$, we only have to check that $c_3$ is equal in $\IBpoly$ to the rightmost circuit above. This is shown as follows, along the same lines of derivation~\eqref{eq:extwoimpl3}:
\[\includegraphics[height=2.3cm]{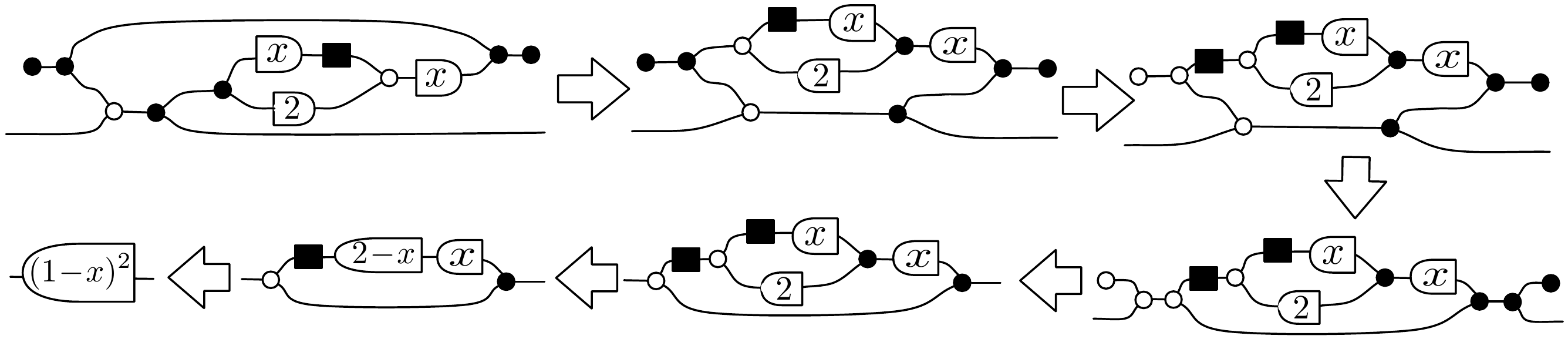}.\]

The circuits in \eqref{impl} can then be thought of as two different implementations of $\circuitUnoMinusXSquare$, which specifies the generating function $\frac{1}{(1-x)^2}$ of the sequence $1,2,3,4,\dots$.
\end{example}

\chapter{Conclusions and Future Directions}\label{chapter:conclusion}

In this thesis we focused on a characterisation of \emph{linear} dynamical systems. However, we would like to stress the importance of the \emph{methodology} grounding our investigation, as we believe it is fruitfully applicable to other classes of systems. The first principle underlying our approach is \emph{modularity}: similarly to $\IH{}$, most of the interesting graphical languages are quite complex; we contend that the right mix of PROP operations --- sum, fibered sum and especially composition --- enables formal analyses that are hampered by monolithic approaches. The second principle is the conviction that a formal theory of circuit diagrams should not be endowed with a primitive notion of causality. We believe that signal flow graphs provide a strong case of how discarding flow directionality is not only harmless, but beneficial for a compositional understanding of computing devices.

In the rest of this section, we mention some of the most promising directions for future research. 


\paragraph{Computing with distributive laws} Distributive laws yield rewriting systems: when we read off the equations from the graph of the law, these are \emph{oriented}, thus it is natural to view them as rewriting rules. Unfortunately, a distributive law does not tell much about the algorithmic aspects of the associated rewriting: for instance, it guarantees that a factorisation exists for each diagram, but it does not produce a terminating sequence of rewriting steps to reach it. An interesting question is thus what principles and general assumptions guarantee \emph{termination} and other useful properties of the rewriting associated with the law. More generally, we are interested in exploring the computational content of distributive laws: for instance, the development of new algorithms for linear algebraic manipulations in $\IH{}$ would greatly improve the appeal of our graphical syntax.

Exploring these aspects would also clarify the connection between our work and existing approaches to string diagram rewriting, such as the research surrounding the \emph{Quantomatic} software package~\cite{Kissinger_quantomatic} and the line of works on higher-dimensional rewriting systems~\cite{Burroni1993,Lafont95-equationalReasoningTwoDimDiagrams,Lafont2003,Mimram14}. 


\paragraph{Categorical control theory} Our approach to signal flow graphs only scratched the surface of the rich body of work devoted to these structures in control theory. A research thread that we find particularly worth investigating is Willems' behavioural approach --- see e.g.~\cite{Willems-linearsystems,Willems2007}. Our proximity to Willems' vision of systems is efficaciously synthesised by the following passage in~\cite{Willems-linearsystems}.

\begin{quote}
\emph{Adding a signal flow direction is often a figment of one's imagination, and when something is not real, it will turn out to be cumbersome sooner or later. [...] The input/output framework is totally inappropriate for dealing with all but the most special system interconnections. [The input/output representation] often needlessly complicates matters, mathematically and conceptually. A good theory of systems takes the behavior as the basic notion and the reference point for concepts and definitions, and switches back and forth between a wide variety of convenient representations.}
\end{quote}
Within this perspective, Willems' analysis treats behaviour as a \emph{relation} rather than a function, which makes it close to our approach and opens the way for studying the system properties considered in~\cite{Willems2007,Willems-linearsystems}, such as controllability, stabilisability and autonomy. Understanding their string diagrammatic significance may yield syntactic characterisations, similar to those that we proposed for deadlock and initialisation freedom, and result in new methods for static analysis.

Another question stemming from traditional control theory is how to extend our approach to include an operational semantics for the continuous-time interpretation of signal flow graphs, where the registers stand for differentiation/integration of a signal. This perspective has been pursued in~\cite{BaezErbele-CategoriesInControl}, using the Laplace transform to reduce the continuous case to the discrete case. Another relevant work is~\cite{Mimram_hyperreals}, which studies the denotational semantics of a class of network diagrams in continuous time using hyperreals. However, these authors focus only on the denotational picture, whereas we would also like to understand how to give an operational account of the continuous-time interpretation. The approach described in~\cite{PavlovicEscardo:1998}, which gives a coinductive formalism for elementary calculus not far removed from our stream semantics, may be helpful in this respect.

\paragraph{Beyond signal flow graphs: other models of computation} It is a natural research question to extend our methodology to other classes of systems: we mention below some of the most promising directions.

First, we believe our approach may be fruitful for graphical formalisms that, like signal flow graphs, are typically translated to more traditional mathematics and seldom reasoned about directly. We think in particular of electrical circuits --- which already attracted categorical modeling~\cite{Fong2013} --- and Kahn process networks~\cite{Brookes98-KahnPrinciple} --- whose formalisation in monoidal categories has been studied in~\cite{Panangaden_dataflow,Mimram_hyperreals}. The key is to understand how to generalise our picture to the modeling of other kinds of behaviour, such as non-linearity, asynchrony and non-determinism.

Another important class of formalisms that we would like to tackle are automata. A promising starting point is the observation that signal flow graphs can be thought as weighted automata on a singleton alphabet --- see e.g.~\cite{Rutten08_rationalstreamscoalgebraically}. Can we characterise the string diagrammatic analogue of other classes of automata? There are at least two ways of extending the case of signal flow diagrams: the first is to enlarge the alphabet, which yields circuits processing infinite trees in place of streams; the second is to generalise the elements of streams from fields to semirings --- this is motivated by the case of non-deterministic automata, which are weighted automata for the boolean semiring.

Yet another interesting research direction concerns quantitative models of computation. In~\cite{Fritz_stochasticmatrices}, Fritz gave a presentation by generators and relations for the PROP of stochastic matrices. By analogy with the case studied in this thesis, we may hope to use this result to achieve a presentation for the PROP whose arrows are the probabilistic analogue of relations, and indeed there are notions of this kind appearing in the literature~\cite{Panangaden98probabilisticrelations}. The technical challenge here is that stochastic matrices have an inherently \emph{non-local} character, because values in a column must always sum to $1$. This prevents the PROP of stochastic matrices from having the limits and colimits necessary to form our cube construction.

Our approach was inspired by graphical formalisms developed in quantum information and concurrency theory: we argued that $\IH{}$ lies somehow at the intersection, as it captures the interplay of Frobenius and Hopf algebras that appear in both research threads. It is a natural question to ask whether we can apply our methodology to model the features which are specific of each of these areas. For instance, we observed that $\IH{\Z_{\scriptscriptstyle 2}}$ only describes the phase-free fragment of the ZX-calculus: it would be interesting to extend our characterisation to include in the picture properly quantum aspects such as the presence of phases. We believe that a modular characterisation for the full ZX-calculus would be particularly valuable for the quantum community: completeness for ZX and its variants have been intensively studied in recent years~\cite{Backens-ZXcompleteness1,Backens-ZXcompleteness2,ZXIncomplete14}.

Another motivating example for our developments were the algebras for Petri nets~\cite{Soboci'nski2010,Bruni2013} and for connectors~\cite{Bruni2006} developed in concurrency theory. In the former case, we believe that our modular understanding may bring important insights towards a complete axiomatisation for the proposed calculi: preliminary investigations in this direction~\cite{Sobocinski2013a} are promising as the identified structure is quite close to $\IH{\Z_{\scriptscriptstyle 2}}$. In the latter case, the formalism for stateless connectors proposed in~\cite{Bruni2006} is extremely close to $\IH{\Z_{\scriptscriptstyle 2}}$, the only difference being that $\Wmult$ is undefined on input $\matrixOneOne$. Nonetheless, this difference is crucial as it prevents us from modeling properly concurrent features such as mutual exclusion. We leave as future work to explore a modular account of the calculus of stateless connectors, which may also give useful insights for an equational theory of calculi of \emph{stateful} connectors like Reo~\cite{Arbab2004}.
\appendix
\chapter{Omitted Proofs}



\section{Proofs of Chapter~2}\label{App:proofsABR}
\label{app:Lawvere}
 \label{app:Lawvere}

\begin{proof}[Proof of Lemma~\ref{lemma:LawNoEquations}] 


The proof given in the main text is completed by verifying that $\lambda$ is compatible with $E$. 
For this purpose, suppose that $\tr{d_1 \in \T}$ and $\tr{d_2 \in \T}$ are such that $d_1 = d_2$ modulo $E$. We need to check that $\lambda(\tr{d_1 \in \T}\tr{c \in \Com}) = \lambda(\tr{d_2 \in \T}\tr{c \in \Com})$ modulo $E$.

We can suppose without loss of generality that $d_1 = d_2$ holds because there are sub-diagrams $a_1$ of $d_1$ and $a_2$ of $d_2$ such that $a_1 = a_2$ is an equation of $E$ --- which also means, because $(\Sigma,E)$ is cartesian, that $a_1$ and $a_2$ have target $1$. Thus, for $i \in \{1,2\}$, we can depict $\tr{d_i}\tr{c}$ as follows.
\begin{eqnarray*}
\lower11pt\hbox{$\includegraphics[height=2.5cm]{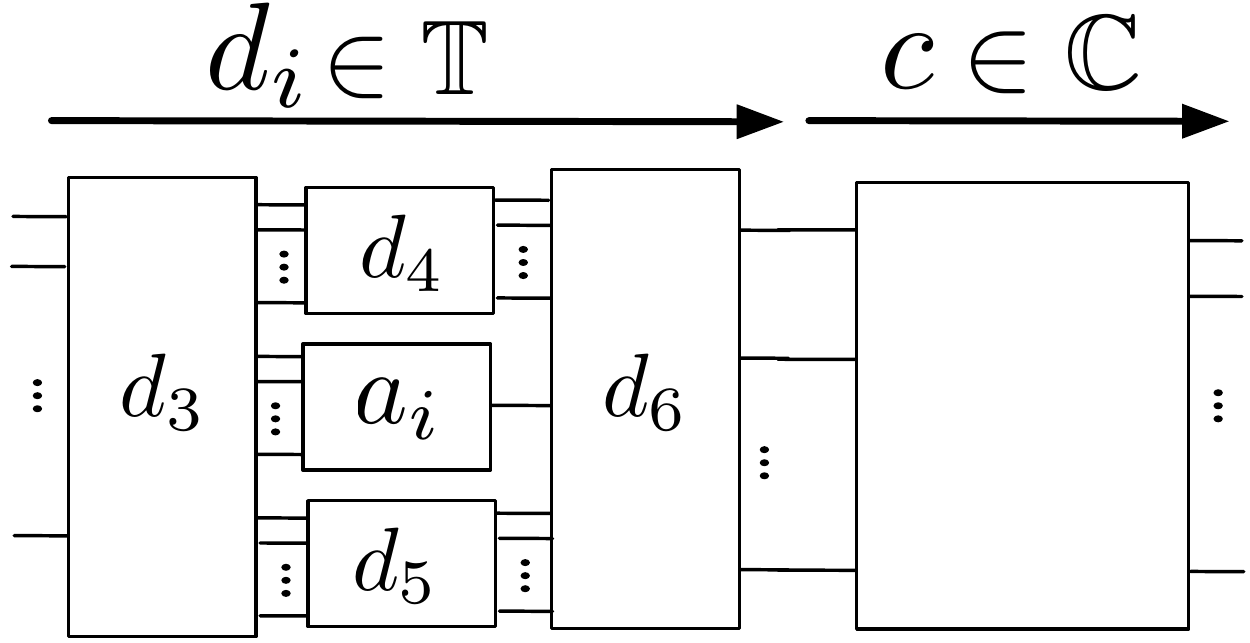}$}
\end{eqnarray*}
We now describe how, using compatibility of $\lambda$ with the structure of $\T$ and $\Com$, the proof reduces to checking that $\lambda$ behaves the same on $a_1$ and on $a_2$. First, we turn our attention to $\tr{c \in \Com}$. Using the same reasoning as in the proof of Theorem~\ref{Th:LawvereCompositePROP}, we can move all the permutation in $c$ to the rightmost part of the diagram --- note that this transformation does not invalidate our argument because $\lambda$ respects the equations of $\Com$.
\begin{eqnarray*}
\lower11pt\hbox{$\includegraphics[height=2.5cm]{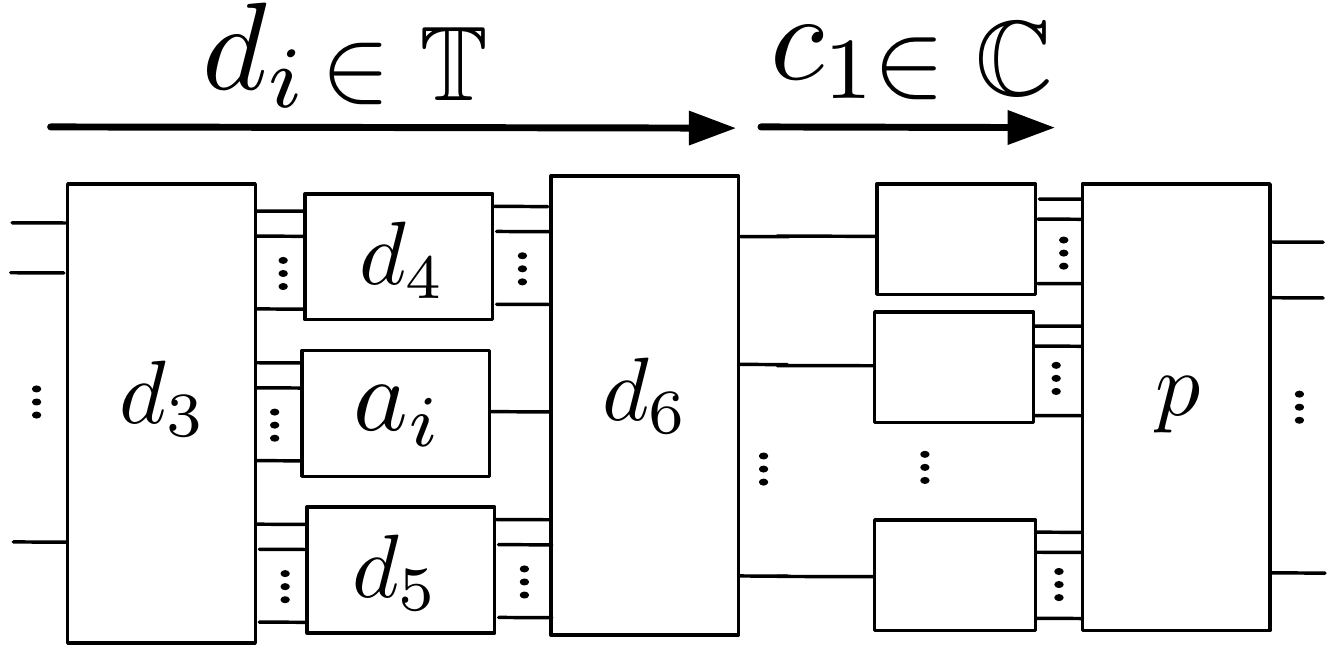}$}
\end{eqnarray*}
In the above diagram, $p$ is an arrow of $\Perm$ and $c_1$ is a $\tns$-product of string diagrams of $\Com$ with source $1$ where no $\symNet$ appears. Being a distributive law, $\lambda$ is compatible with the permutations of $\Com$ and $\T$: we thus have that $$\lambda(\tr{d_i}\tr{c}) = \lambda(\tr{d_i}\tr{c_1 \poi p}) = \lambda(\tr{d_i}\tr{c_1})\tr{p}.$$
Therefore, to prove our statement it suffices to check that $\lambda(\tr{d_1}\tr{c_1}) = \lambda(\tr{d_2}\tr{c_1})$.
\begin{eqnarray*}
\lower11pt\hbox{$\includegraphics[height=2.5cm]{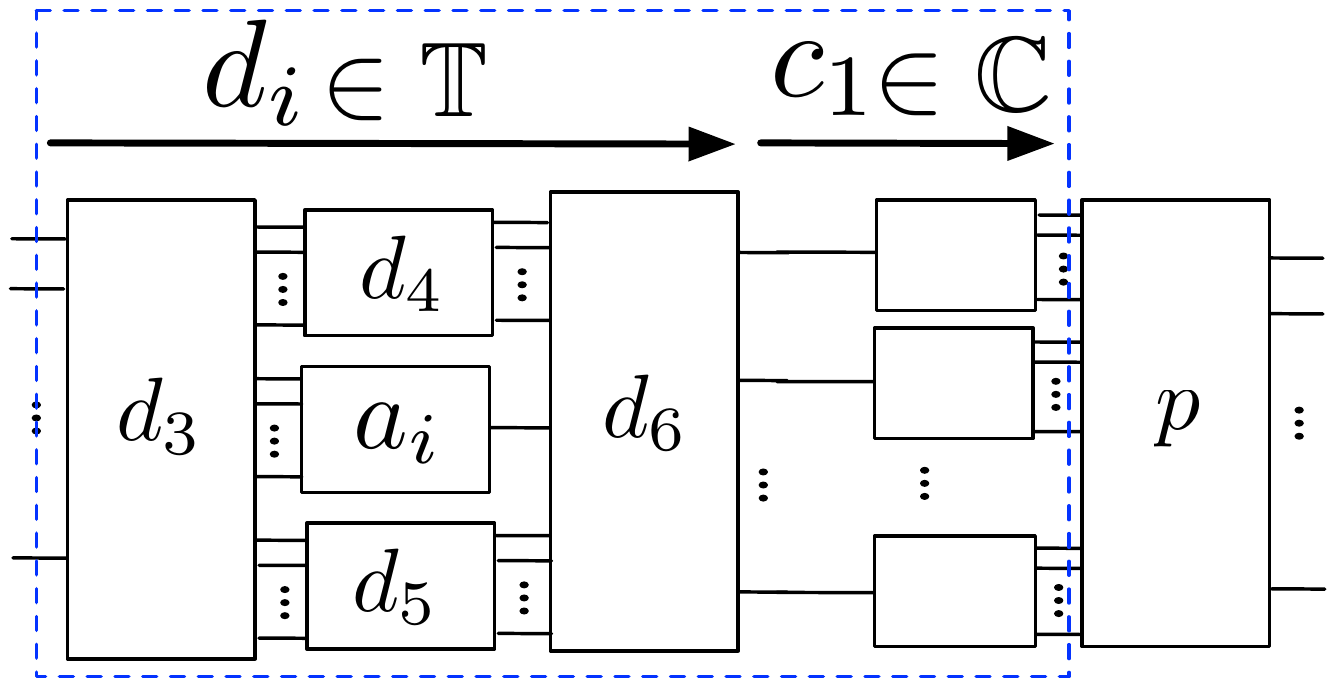}$}
\end{eqnarray*}
 We can now use the fact that $\lambda$ is compatible with composition in $\T$, meaning that
$$\lambda(\tr{d_i}\tr{c_1}) = \lambda(\tr{d_3\poi (d_4 \tns a_i \tns d_5)\poi d_6}\tr{c_1}) = \lambda(\tr{d_3\poi (d_4 \tns a_i \tns d_5)}\tr{c_2 \in \Com})\poi \tr{d_7 \in \T}$$
where $\tr{c_2 \in \Com}\tr{d_7 \in \T} = \lambda(\tr{d_6 \in \T}\tr{c_1 \in \Com})$. This passage can be represented diagrammatically as follows:
\begin{eqnarray*}
\lower25pt\hbox{$\includegraphics[height=2.5cm]{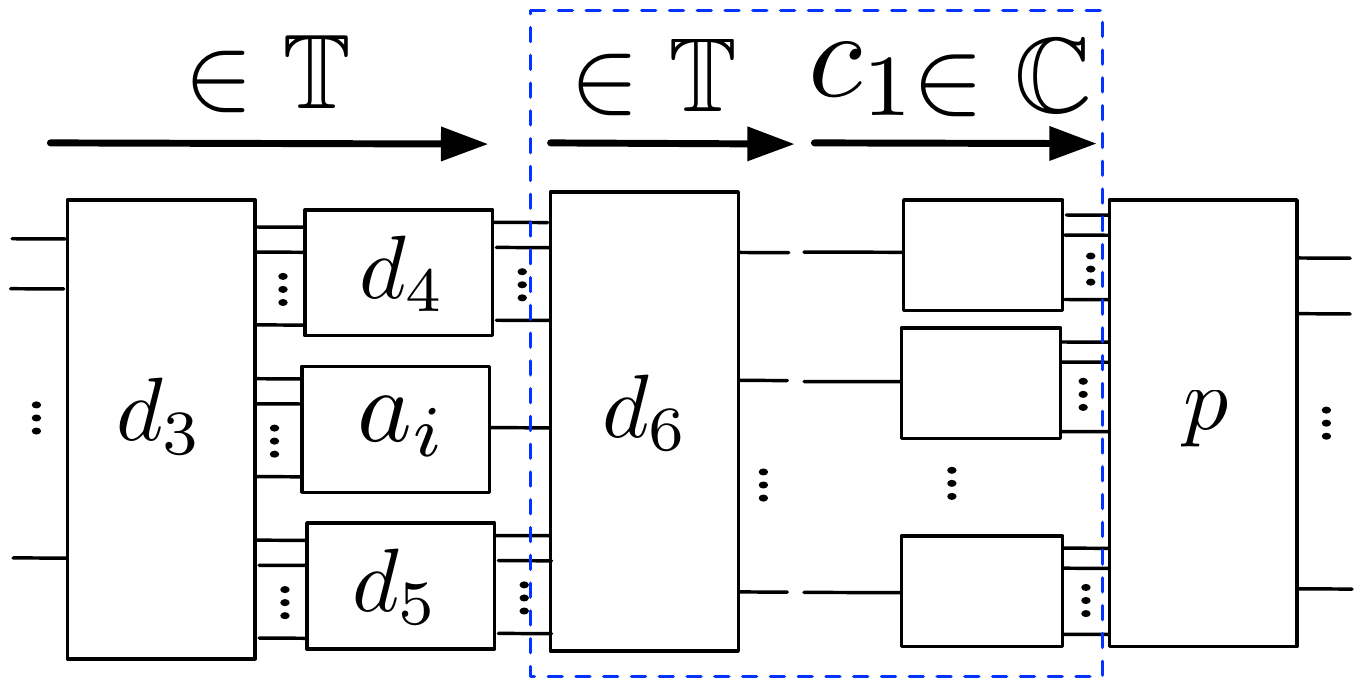}$} \qquad \xymatrix{ \ar@{|->}[rrr]^{\text{apply }\lambda \text{ on }\tr{d_6}\tr{c_1}}&&& } \qquad \lower25pt\hbox{$\includegraphics[height=2.5cm]{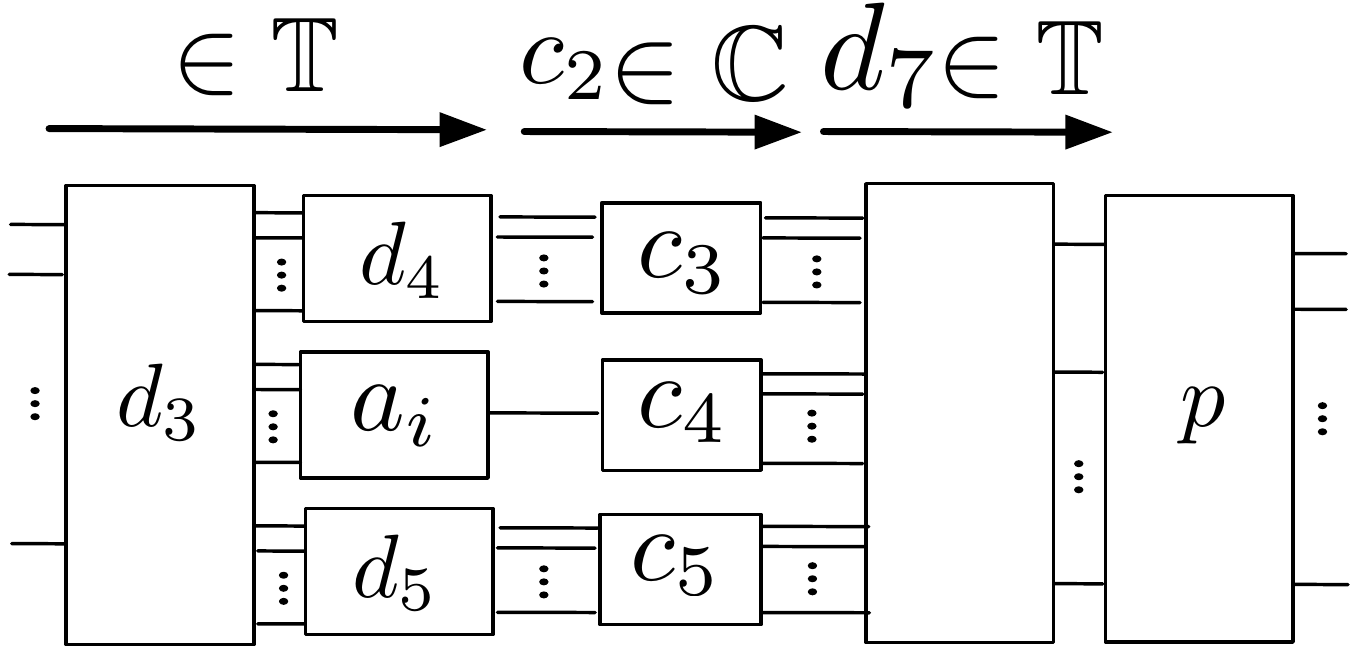}$}.
\end{eqnarray*}
Thus we reduced further our proof to verifying that $\lambda(\tr{d_3\poi (d_4 \tns a_1 \tns d_5)}\tr{c_2 \in \Com}) = {\lambda(\tr{d_3\poi (d_4 \tns a_2 \tns d_5)}\tr{c_2 \in \Com})}$.
\begin{eqnarray*}
\lower11pt\hbox{$\includegraphics[height=2.5cm]{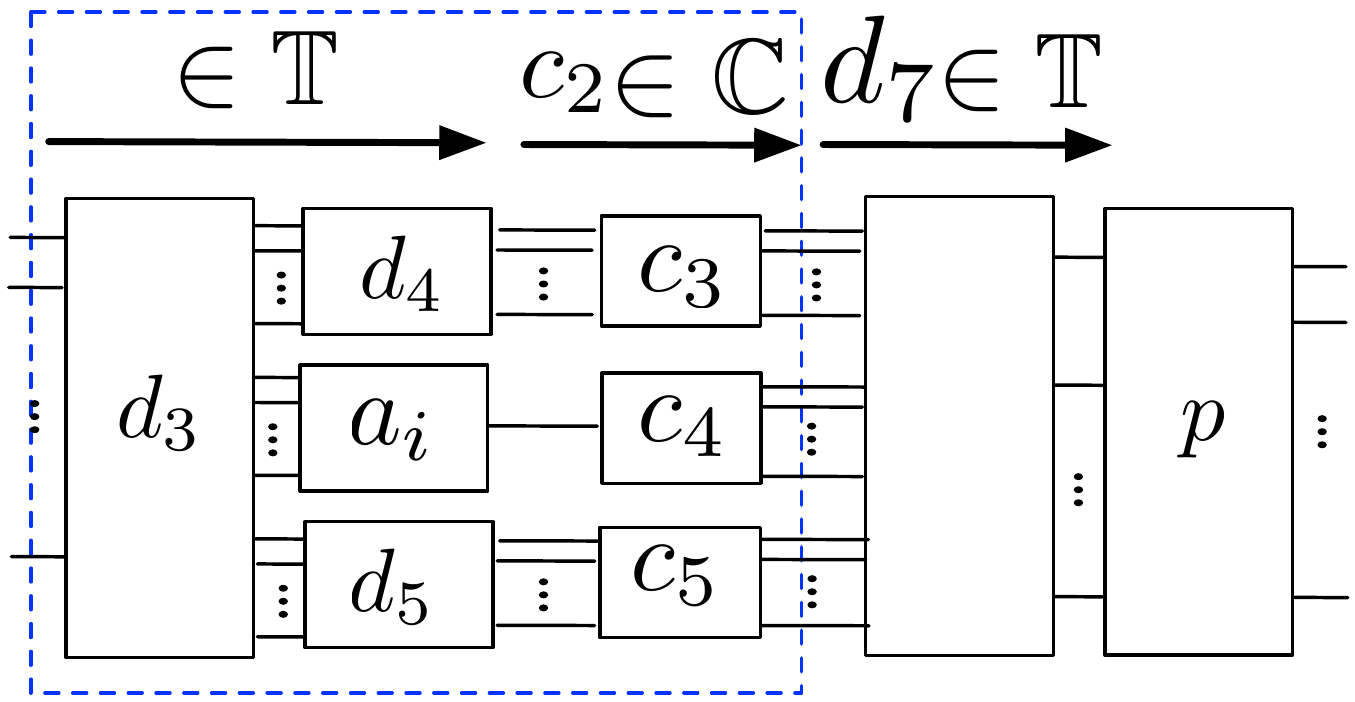}$}
\end{eqnarray*}
We can then reason analogously, using compatibility of $\lambda$ with composition and monoidal product of $\T$. We have that
$$\lambda(\tr{d_3\poi (d_4 \tns a_i \tns d_5)}\tr{c_2 \in \Com}) = \lambda(\tr{d_3}\tr{c_6 \in \Com}) \tr{d_8 \in \T}$$
where $\tr{c_6 \in \Com}\tr{d_8 \in \T} = \lambda(\tr{d_4 \tns a_i \tns d_5}\tr{c_2})$. Thus it suffices to prove that $\lambda(\tr{d_4 \tns a_1 \tns d_5}\tr{c_2}) = {\lambda(\tr{d_4 \tns a_2 \tns d_5}\tr{c_2})}$. But this equality reduces to proving that $\lambda(\tr{a_1}\tr{c_4}) = \lambda(\tr{a_2}\tr{c_4})$, because
$$\lambda(\tr{d_4 \tns a_i \tns d_5}\tr{c_2}) = \lambda(\tr{d_4 \tns a_i \tns d_5}\tr{c_3 \tns c_4 \tns c_5}) = \lambda(\tr{d_4}\tr{c_3}) \tns \lambda(\tr{a_i}\tr{c_4})\tns \lambda(\tr{d_5}\tr{c_5}).$$

Hence it remains to verify that $\lambda(\tr{a_1}\tr{c_4}) = \lambda(\tr{a_2}\tr{c_4})$. Since it has source $1$, $c_4$ must be equal in $\Com$ to
\[   \lower10pt\hbox{$\includegraphics[height=1cm]{graffles/BcomultMultipla.pdf}$} \qquad \text{ or to } \qquad \lower0pt\hbox{$\includegraphics[height=.2cm]{graffles/BcounitNoFrame.pdf}$}\ .\]
Because $\lambda$ preserves the equations of $\Com$, we can compute $\lambda(\tr{a_i}\tr{c_4})$ assuming that $c_4$ has one of the two shapes above. By assumption, $\lambda$ is calculated by factorising $\tr{a_i}\tr{c_4}$ as an arrow of $\LwA{\Sigma}$. This yields, depending on the shape of $c_4$,
\begin{eqnarray*}
\lower14pt\hbox{$\includegraphics[height=1.4cm]{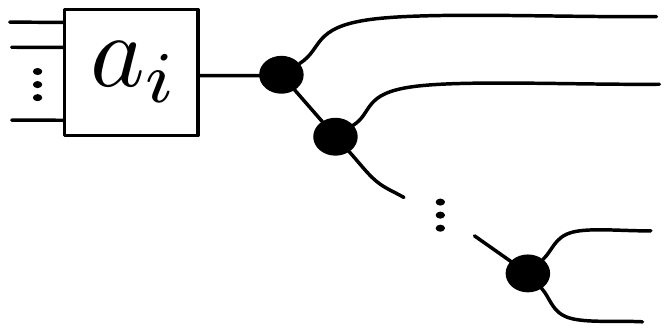}$}\quad \xmapsto{\lambda} \quad \lower29pt\hbox{$\includegraphics[height=3.2cm]{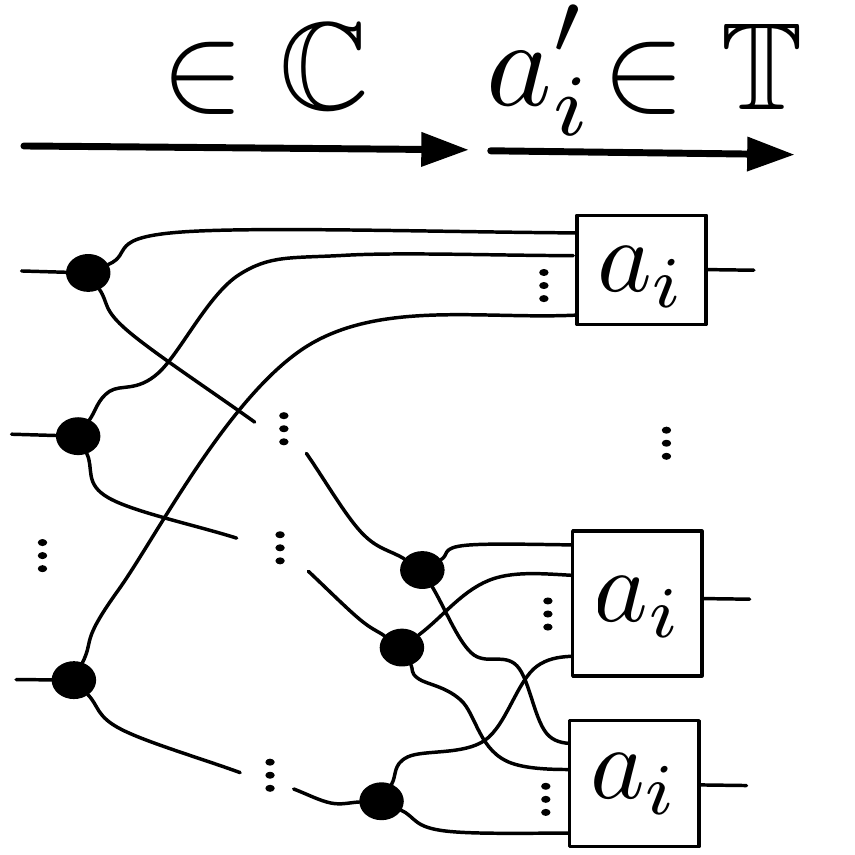}$}\qquad \text{ or} \qquad
\lower8pt\hbox{$\includegraphics[height=.8cm]{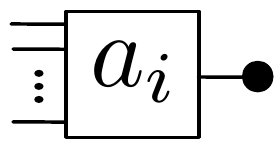}$} \quad \xmapsto{\lambda} \quad \lower8pt\hbox{$\includegraphics[height=1.3cm]{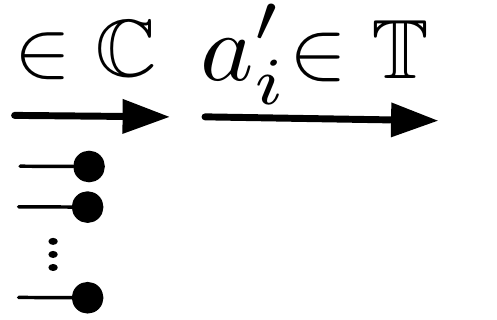}$}.
\end{eqnarray*}
Since $a_1'=a_2'$ is provable from the equations of $E$, we conclude that $\lambda(\tr{a_1}\tr{c_4}) = \lambda(\tr{a_2}\tr{c_4})$.
\end{proof}

\begin{proof}[Proof of Lemma~\ref{lemma:LawvereFactorisation}]
Our starting point is a factorisation for cartesian terms. Suppose $t_i$ is a cartesian term in $\tpl{t_1,\dots,t_m}$: we can see $t_i$ itself as an arrow $n \tr{\tpl{t_i}} 1$. Moreover, we can decompose it as follows.
  \[ \xymatrix@C=40pt{n \ar@/^2pc/[rr]^{\tpl{t_i}} \ar[r]^{\var{t_i}} & \size{\var{t_i}} \ar[r]^{\tpl{\tilde{t}_i}} & 1}.\]
  Here $\var{t_i}$ is regarded as an arrow of type $n \to \size{\var{t_i}}$: this is well-typed because it is a list of length $\size{\var{t_i}}$ where only variables among $x_1,\dots,x_n$ may occur. Also, it belongs to the sub-PROP $\Com$ because all cartesian terms appearing in it are variables. For the arrow $\size{\var{t_i}} \to 1$, the cartesian term $\tilde{t}_i$ is obtained from $t_i$ by replacing the $j$th variable of $t_i$ with $x_j$, for $1 \leq j \leq \size{\var{t_i}}$. Then $\var{\tilde{t}_i} = \tpl{x_1,\dots,x_{\size{\var{t_i}}}}$ by construction, meaning that $\tpl{\tilde{t}_i}$ is a linear list and it is well-typed as an arrow $\size{\var{t_i}} \to 1$. To see that $\var{t_i}\poi \tilde{t}_i = t_i$, note that precomposition with $\var{t_i}$ has the effect of replacing variable $x_j$ in $\tilde{t}_i$ with the $j$th variable in the list $\var{t_i}$, that is, the $j$th variable occurring in $t_i$. Therefore we obtain $t_i$ as the result of composition.

 It is useful to show that our decomposition for $t_i$ is unique up-to permutation, because this will guarantee the same property for $\tpl{t_1,\dots,t_m}$. For the sake of readability, let us use the alias $z$ for $\size{\var{t_i}} \in \N$. Suppose now $n \tr{c\in \Com}z' \tr{\tpl{s_i} \in \T}1$ is such that $c \poi \tpl{s_i} = \tpl{t_i}$. Because $\tpl{s_i}$ is a linear list, $\size{\var{s_i}} = z'$. Because $c \poi \tpl{s_i} = \tpl{t_i}$, this also implies that $z' = \size{\var{s_i}} = \size{\var{t_i}} = z$. Now, since $c \poi \tpl{s_i} = \var{t_i} \poi \tpl{\tilde{t}_i}$, the only difference between $s_i$ and $\tilde{t}_i$ is the order in which the variables $x_1,\dots,x_z$ appear. By construction, $\var{\tilde{t}_i} = \tpl{x_1,\dots,x_z}$ and thus there exists a permutation $p = \tpl{x_{p^{-1}(1)},\dots,x_{p^{-1}(z)}} \: z \to z$ such that $p \poi \tpl{\tilde{t}_i} = \tpl{s_i}$. Since $c \poi \tpl{s_i} = \tpl{t_i}$, precomposing with $c$ the list $\var{s_i} \: z \to z$ of variables occurring in $s_i$ must yield $\var{t_i}$. As we saw, $\var{s_i}$ is just $p = \tpl{x_{p^{-1}(1)},\dots,x_{p^{-1}(z)}} \: z \to z$, meaning that $c \poi p = \var{t_i}$. We can thereby conclude, by commutativity of the following diagram, that the factorisation of $t_i$ is unique up-to permutation.
 \[\vcenter{
            \xymatrix@R=10pt@C=50pt{
                & z \ar@/^/[dr]^-{\tpl{t_i}} & \\
               n\ar@/^/[ur]^-{\var{t_i}} \ar[r]_{c'} & z\ar[u]_-{p} \ar[r]_{\tpl{s_i}} &1
            }
}\]

 It is now straightforward to extend the factorisation to the list $\tpl{t_1,\dots,t_m}$. First, we can form the $m$-fold product of arrows $\tpl{\tilde{t}_i} \: \size{\var{t_i}} \to 1$.
 \begin{eqnarray*}
 \hat{d} \df \tpl{\tilde{t}_1}\times \dots \times \tpl{\tilde{t}_m} & \: & \size{\var{t_1}} \times \dots \times \size{\var{t_m}} \to 1 \times \dots \times 1 \\
 & \: & \size{\var{t_1,\dots,t_m}} \to m.
 \end{eqnarray*}
By its type, $\hat{d}$ is a list $\tpl{\hat{t}_1,\dots,\hat{t}_m}$ of $m$ elements on variables $\size{\var{t_1,\dots,t_m}}$. By definition of $\times$, a term $\hat{t}_i$ is obtained from $\tilde{t}_i$ by replacing the variable $x_j$ with $x_{j+k}$, where $k = \size{\var{\tilde{t}_1,\dots,\tilde{t}_{i-1}}}$. Thus $\tpl{\hat{t}_1,\dots,\hat{t}_m}$ is a linear list for which $\var{\hat{t}_1,\dots,\hat{t}_m} = [x_1,\dots,x_{\size{\var{t_1,\dots,t_m}}}]$.

 We now construct a second arrow $\hat{c} \: n \to \size{\var{t_1,\dots,t_m}}$. The existence of arrows $\var{t_i} \: n \to \size{\var{t_i}}$ triggers the universal property of the product $\size{\var{t_1,\dots,t_m}}$, which yields
  \[ \hat{c} \df \var{t_1,\dots,t_m}  \:  n \to \size{\var{t_1,\dots,t_m}}. \]
 By construction, $\hat{d}$ and $\hat{c}$ make the following diagram commute, where the unlabeled arrows are the product projections.
 \begin{eqnarray}\label{eq:Lawvere_factGen}
 \vcenter{
 \xymatrix@C=20pt{
 && \size{\var{t_1}} \ar@/^2pc/[rrr]^{\tpl{\tilde{t}_1}}& \!\! \!\! \!\! \dots\!\! \!\! & \size{\var{t_m}}\ar@/^2pc/[rrr]^{\tpl{\tilde{t}_m}}  & 1 & \dots & 1 \\
 n \ar@/^/[urr]^{\var{t_1}} \ar@/^/[urrrr]^<<<<<<<<<<<{\tpl{\var{t_m}}}  \ar@/_/@{-->}[rrr]_{\hat{c}} &&& \size{\var{t_1,\dots,t_m}} \ar[ul] \ar[ur] \ar@/_/@{-->}[rrr]_{\hat{d}} &&& m \ar[ul] \ar[ur] &
 }
 }
 \end{eqnarray}
To see that $\hat{c} \poi \hat{d} = \tpl{t_1,\dots,t_m}$, observe that $\tpl{t_1,\dots,t_m}$ is given by the universal property of the product $m$:
 \begin{eqnarray}\label{eq:Lawvere_factComm}
 \vcenter{
 \xymatrix@C=30pt{
   & 1 & \dots & 1 \\
 n \ar@/^/[ur]^{\tpl{t_1}} \ar@/^/[urrr]^<<<<<<<<{\tpl{t_m}}  \ar@/_/@{-->}[rr]_{\tpl{t_1,\dots,t_m}} && m \ar[ul] \ar[ur] &
 }
 }
 \end{eqnarray}
Now, by construction, $\var{t_i}\poi \tpl{\tilde{t}_i} = \tpl{t_i}$ for each $1 \leq i \leq m$ and $\hat{c}\poi \hat{d}$ makes \eqref{eq:Lawvere_factGen} commute. This implies that $\hat{c}\poi \hat{d}$ also makes \eqref{eq:Lawvere_factComm} commute, meaning by uniqueness that it is equal to $\tpl{t_1,\dots,t_m}$.

The proof is concluded by observing that the proposed factorisation is unique up-to permutation: this immediately follows by uniqueness of the factorisation for each $\tpl{t_i}$, proved above, and universal property of the product.
\end{proof} 

\bigskip

\label{App:proofsER} We conclude this section with the proofs that $\SpanToER \: \Injop \bicomp{\Perm} {\Inj} \to \ER$ and $\CospanToER \: \F \bicomp{\Perm} \Fop \to \ER$ --- introduced in \S~\ref{sec:ER} --- are PROP morphisms.

\begin{proposition}\label{prop:SpanToERFunctor} $\SpanToER$ is a PROP morphism.	
\end{proposition}
\begin{proof}
The hard part is to show that $\SpanToER$ preserves composition. For this purpose, consider the following diagram in $\Inj$, expressing the composition of arrows $\tl{f_1}\tr{g_1}$ and $\tl{f_2}\tr{g_2}$ in $\Injop \bicomp{\Perm} {\Inj}$.
\begin{equation}\label{eq:pbSpanToER}
\vcenter{
\xymatrix@=12pt{
& & {r} \pushoutcorner \ar[dl]_{h_1} \ar[dr]^{h_2} \\
& {z_1} \ar[dl]_{f_1} \ar[dr]^{g_1} & & {z_2} \ar[dl]_{f_2} \ar[dr]^{g_2} \\
{n} & & {z} & & {m}
}
}
\end{equation}
Pullbacks in $\Inj$ are constructed as in $\Set$. Therefore,
\begin{align}
g_1(v) = f_2(w) &\Leftrightarrow \text{ there is }u \in \ord{r}\text{ such that }h_1(u) = v, h_2(u) = w \nonumber \\
&\Leftrightarrow \text{(since $h_1, h_2$ are injections) }\inv{h_1}(v) = \inv{h_2}(w). \label{eq:pbinj}\tag{$\star$}
 \end{align}
 Also, by definition of composition in $\Injop \bicomp{\Perm} \Inj$, $(\tl{f_1} \tr{g_1}) \poi (\tl{f_2}\tr{g_2})
=\ \tl{f_1 h_1}\tr{g_2 h_2}$.

We now split the statement to prove in two distinct lemmas.

\begin{lemma}\label{SpanToERLemma1} If  $(v,w)$ is in $\SpanToER(\tl{f_1 h_1}\tr{g_2 h_2})$ then $(v,w)$ is in $\SpanToER(\tl{f_1} \tr{g_1}) \poi \SpanToER(\tl{f_2}\tr{g_2})$.\end{lemma}
\begin{lemma}\label{SpanToERLemma2}  If $(v,w)$ is in $\SpanToER(\tl{f_1} \tr{g_1}) \poi \SpanToER(\tl{f_2}\tr{g_2})$ then $(v,w)$ is in $\SpanToER(\tl{f_1 h_1}\tr{g_2 h_2})$.\end{lemma}
\begin{proof}[Proof of Lemma~\ref{SpanToERLemma1}]
Suppose that $(v,w)\in \SpanToER(\tl{f_1 h_1}\tr{g_2 h_2})$. By definition of $\SpanToER$, the following three can be the case.
\begin{enumerate}[label=(\Roman{*})]
  \item If $v=w$, then $(v,w) \in \SpanToER(\tl{f_1} \tr{g_1}) \poi \SpanToER(\tl{f_2}\tr{g_2})$ because an equivalence relation is reflexive.
  \item Otherwise, suppose that $\inv{(f_1 h_1)}(v) = \inv{(g_2 h_2)}(w)$. It follows, by observation \eqref{eq:pbinj}, that $g_1(\inv{f_1}(v)) = f_2(\inv{g_2}(w))$. Call this element $u$. Then $\inv{g_1}(u) = \inv{f_1}(v)$ and $\inv{f_2}(u) = \inv{g_2}(w)$, witnessing $(v,u) \in \SpanToER(\tl{f_1} \tr{g_1})$ and $(u,w) \in \SpanToER(\tl{f_2}\tr{g_2})$. Therefore $(v,w) \in \SpanToER(\tl{f_1} \tr{g_1}) \poi \SpanToER(\tl{f_2}\tr{g_2})$.
  \item The case in which $\inv{(g_2 h_2)}(v) = \inv{(f_1 h_1)}(w)$ is handled symmetrically with respect to (II).
\end{enumerate}
\end{proof}
In order to prove Lemma~~\ref{SpanToERLemma2}, it is first convenient to makes the following observation. Intuitively, it states that there is no need of transitivity when reasoning about equivalence relations generated by injective functions.
\begin{lemma}\label{lemma:EqRel=SymReflforInj}
Let $e^{=S}$ denote the symmetric and reflexive closure of a relation $e$.
Then:
$$\big(\SpanToER(\tl{f_1} \tr{g_1}) \relcomp \SpanToER(\tl{f_2}\tr{g_2})\big)^{=S} = \eqr{\SpanToER(\tl{f_1} \tr{g_1}) \relcomp \SpanToER(\tl{f_2}\tr{g_2})}$$
\end{lemma}
\begin{proof}
First, it is not hard to check that $\big(\SpanToER(\tl{f_1} \tr{g_1}) \relcomp \SpanToER(\tl{f_2}\tr{g_2})\big)^{=S}$ is an equivalence relation: reflexivity and symmetry are given by construction and transitivity can be verified by case analysis on pairs in $\big(\SpanToER(\tl{f_1} \tr{g_1}) \relcomp \SpanToER(\tl{f_2}\tr{g_2})\big)^{=S}$.

Also, $\big(\SpanToER(\tl{f_1} \tr{g_1}) \relcomp \SpanToER(\tl{f_2}\tr{g_2})\big)^{=S}$ contains $\SpanToER(\tl{f_1} \tr{g_1}) \relcomp \SpanToER(\tl{f_2}\tr{g_2})$ and is included in $\eqr{\SpanToER(\tl{f_1} \tr{g_1}) \relcomp \SpanToER(\tl{f_2}\tr{g_2})}$. Since by construction $\eqr{\SpanToER(\tl{f_1} \tr{g_1}) \relcomp \SpanToER(\tl{f_2}\tr{g_2})}$ is the smallest equivalence relation containing $\SpanToER(\tl{f_1} \tr{g_1}) \relcomp \SpanToER(\tl{f_2}\tr{g_2})$, it follows that
$$\big(\SpanToER(\tl{f_1} \tr{g_1}) \relcomp \SpanToER(\tl{f_2}\tr{g_2})\big)^{=S} = \eqr{\SpanToER(\tl{f_1} \tr{g_1}) \relcomp \SpanToER(\tl{f_2}\tr{g_2})}.$$
\end{proof}
We are now ready to supply a proof of Lemma~\ref{SpanToERLemma2}.
\begin{proof}[Proof of Lemma~\ref{SpanToERLemma2}]
Our assumption is that $(v,w) \in \SpanToER(\tl{f_1} \tr{g_1}) \poi \SpanToER(\tl{f_2}\tr{g_2})$. By Lemma~\ref{lemma:EqRel=SymReflforInj}, $(v,w) \in \big(\SpanToER(\tl{f_1} \tr{g_1}) \relcomp \SpanToER(\tl{f_2}\tr{g_2})\big)^{=S}$, meaning that we only have three cases to consider.
\begin{itemize}
\item If $v =w$, then clearly $(v,w)$ is in the equivalence relation $\SpanToER(\tl{f_1 h_1}\tr{g_2 h_2})$.
\item Otherwise, suppose that $(v,w)$ is in $\SpanToER(\tl{f_1}\tr{g_1}) \relcomp\SpanToER(\tl{f_2} \tr{g_2})$. Since by assumption $v$ and $w$ are distinct elements of $\ord{n}$ or $\ord{m}$, we are in the situation in which $v$ is in $\ord{n}$ and $w$ is in $\ord{m}$. Also, there exists $u \in \ord{z}$ such that $(v,u)$ is in $\SpanToER(\tl{f_1}\tr{g_1})$ and $(u,w)$ is in $\SpanToER(\tl{f_2} \tr{g_2})$. By definition, this means that there are elements $q_1 \in \ord{z_1}$ and $q_2 \in \ord{z_2}$ such that $f_1(q_1) = v$, $g_1(q_1) = u$, $f_2(q_2) = u$ and $g_2(q_2) = w$. Then, by property $(\star)$ of the pullback square in~\eqref{eq:pbSpanToER}, there exists $q_3 \in \ord{r}$ such that $h_1(q_3) = q_1$ and $h_2(q_3) = q_2$. Therefore, $f_1h_1(q_3) = v$ and $g_2h_2 (q_3) = w$, meaning that $(v,w)$ is in $\SpanToER(\tl{f_1 h_1}\tr{g_2 h_2})$.
\item In the remaining case, $(v,w)$ is in $\SpanToER(\tl{g_2}\tr{f_2}) \relcomp\SpanToER(\tl{g_1} \tr{f_1})$, meaning that $v$ is in $\ord{m}$ and $w$ is in $\ord{n}$. By a symmetric argument, we can conclude that $\inv{f_1 h_1}(w) = \inv{g_2 h_2} (v)$ and thus $(v,w)$ is in $\SpanToER(\tl{f_1 h_1}\tr{g_2 h_2})$.
\end{itemize}
In both cases we were able to prove that $(v,w)$ is in $\SpanToER(\tl{f_1 h_1}\tr{g_2 h_2})$. This concludes the proof of Lemma~\ref{SpanToERLemma2}.
\end{proof}
The proof of Lemma~\ref{SpanToERLemma2} concludes the proof of Proposition~\ref{prop:SpanToERFunctor}.
\end{proof}
\begin{proposition}\label{lemma:Psifunctor}
$\CospanToER$ is a PROP morphism from $\F \bicomp{\Perm} \Fop$ to $\ER$.	
\end{proposition}
\begin{proof}
The hard part is to show that $\CospanToER$ preserves composition. In the diagram below let the centre square be a pushout diagram.
\begin{eqnarray}\label{eq:poComp}
\vcenter{
\xymatrix{
n \ar[dr]_{p_1} & & \ar[dl]_{q_1} z \ar[dr]^{p_2} & & m \ar[dl]^{q_2}\\
& z_1 \ar[dr]_{p_3} & & z_3 \ar[dl]^{q_3} & \\
 & & z_3 & & &
}
}
\end{eqnarray}
The equivalence to check is the following, for $u,v \in \ord{n+m}$:
\begin{itemize}
\item[(i)] $(u,v) \in \CospanToER(\tr{p_1}\tl{q_1} \poi \tr{p_2}\tl{q_2}) = \CospanToER(\tr{p_3 p_1}\tl{q_3 q_2}).$
\item[(ii)] $(u,v) \in \CospanToER(\tr{p_1}\tl{q_1})\poi \CospanToER(\tr{p_2}\tl{q_2}).$
\end{itemize}

For the direction $(i) \To (ii)$, we split the proof into the following cases:
\begin{enumerate}
 \item First, suppose that $u \in \ord{n}$ and $v \in \ord{m}$. Then by assumption $p_3 p_1 (u) = q_3 q_2 (v)$. By the way pushouts are constructed in $\F$, there are elements $w_1, \dots, w_n$ of $z$ such that
     \begin{align*}
     p_1(u) & = q_1(w_1) \\
     p_2(w_1) & = q_2(w_2) \\
     q_1(w_2) & = q_1(w_3) \\
     \dots & \dots \\
     p_2(w_n) & = q_2(v).
     \end{align*}
     We now show that all pairs in the chain $u,w_1,\dots,w_n,v$ are in the equivalence relation generated by $\CospanToER(\tr{p_1}\tl{q_1}) \relcomp \CospanToER(\tr{p_2}\tl{q_2})$:
     \begin{itemize}
     \item since $p_1(u) = q_1(w_1)$ then $(u,w_1) \in \CospanToER(\tr{p_1}\tl{q_1})$ and since $(w_1,w_1) \in \CospanToER(\tr{p_2}\tl{q_2})$ then $(u, w_1) \in  \CospanToER(\tr{p_1}\tl{q_1}) \relcomp \CospanToER(\tr{p_2}\tl{q_2})$.
     \item For $i \ls n$ odd, $p_2(w_i) = p_2(w_{i+1})$ implies $(w_{i+1},w_i) \in \CospanToER(\tr{p_2}\tl{q_2})$ and since $(w_{i+1},w_{i+1}) \in \tr{p_1}\tl{q_1}$ then $(w_{i+1},w_i) \in \CospanToER(\tr{p_1}\tl{q_1}) \relcomp \CospanToER(\tr{p_2}\tl{q_2})$ and so $(w_i,w_{i+1}) \in
         \eqr{\CospanToER(\tr{p_1}\tl{q_1}) \relcomp \CospanToER(\tr{p_2}\tl{q_2})}$.
      \item For $i \ls n$ even, $q_1(w_i) = q_1(w_{i+1})$ implies $(w_i,w_{i+1}) \in \CospanToER(\tr{p_1}\tl{q_1})$ and since $(w_{i+1},w_{i+1}) \in \tr{p_2}\tl{q_2}$ then $(w_i,w_{i+1}) \in
         \CospanToER(\tr{p_1}\tl{q_1}) \relcomp \CospanToER(\tr{p_2}\tl{q_2})$.
      \item Finally, $p_2(w_n) = q_2(v)$ implies $(w_n,v) \in \CospanToER(\tr{p_2}\tl{q_2})$ and since $(w_n,w_n) \in \CospanToER(\tr{p_2}\tl{q_2})$ then $(w_n,v) \in \CospanToER(\tr{p_1}\tl{q_1}) \relcomp \CospanToER(\tr{p_2}\tl{q_2})$.
\end{itemize}
It follows by transitivity that $(u,v) \in \eqr{\CospanToER(\tr{p_1}\tl{q_1}) \relcomp \CospanToER(\tr{p_2}\tl{q_2})}$. Since $u,w \in \ord{n+m}$, then also $(u,v) \in \CospanToER(\tr{p_1}\tl{q_1}) \poi \CospanToER(\tr{p_2}\tl{q_2})$.
\item We now focus on the case in which $u,v \in i_1(n)$. By assumption, $p_3 p_1 (u) = p_3 p_2 (v)$. Therefore, by the way pushouts are constructed in $\F$, there are elements $w_1, \dots, w_n$ of $\ord{z}$ such that
     \begin{align*}
     p_1(u) & = q_1(w_1) \\
     p_2(w_1) & = q_2(w_2) \\
     q_1(w_2) & = q_1(w_3) \\
     p_2(w_3) & = p_2(w_4) \\
     \dots & \dots \\
     q_1(w_n) & = p_2(v).
     \end{align*}
     Analogously to case 1 considered above, we can show that all pairs in the chain $u,w_1,\dots,w_n,v$ are in the equivalence relation generated by $\CospanToER(\tr{p_1}\tl{q_1}) \relcomp \CospanToER(\tr{p_2}\tl{q_2})$, implying that $(u,v) \in \CospanToER(\tr{p_1}\tl{q_1}) \poi \CospanToER(\tr{p_2}\tl{q_2})$.
     \item The case in which $u \in \ord{m}$ and $v \in \ord{n}$, and the one in which $u,v \in \ord{m}$, are handled symmetrically to the first and the third case above, respectively.
\end{enumerate}

\noindent We now provide a proof for direction $(ii) \To (i)$. To this aim, it will be useful to first show the following statement:
 \begin{itemize}
 \item[(*)] Consider diagram \eqref{eq:poComp} and $(u,v) \in \CospanToER(\tr{p_1}\tl{q_1}) \relcomp \CospanToER(\tr{p_2}\tl{q_2})$. Then $u$ and $v$ are mapped (either by $p_3p_1$, $q_3p_2$, $q_3q_2$ of $p_3q_1$, according to the set they belong to) to the same element of $z_3$.
\end{itemize}
For the proof of such statement, let $(u,v) \in \CospanToER(\tr{p_1}\tl{q_1}) \relcomp \CospanToER(\tr{p_2}\tl{q_2})$, be witnessed by $w$ (necessarily an element of $z$) such that $(u,w) \in \CospanToER(\tr{p_1}\tl{q_1})$ and $(w,v) \in \CospanToER(\tr{p_2}\tl{q_2})$. We reason by cases:
\begin{itemize}
    \item if $u \in \ord{n}$ and $v\in \ord{z}$ then, by definition of $\CospanToER$, $p_1(u) = q_1(w)$ and $p_2(w) = p_2(v)$. By commutativity of \eqref{eq:poComp}, $p_3 q_1 (w) = q_3 p_2(w)$ and thus $p_3 p_1 (u) = q_3 p_2 (v)$. 
    \item If $u \in \ord{z}$ and $v\in \ord{m}$ then $q_1(u) = q_1(w)$ and $p_2(w) = q_2(v)$. Since $p_3 q_1 (w) = q_3 p_2(w)$, we have that $p_3 q_1 (u) = q_3 q_2 (v)$.
    \item If $u \in \ord{n}$ and $v\in \ord{m}$ then $p_1(u) = q_1(w)$ and $p_2(w) = q_2(v)$. Since $p_3 q_1 (w) = q_3 p_2(w)$, we have that $p_3 p_1 (u) = q_3 q_2 (v)$.
    \item Finally, if $u,v \in \ord{z}$ then $q_1(u) = q_1(w)$ and $p_2(w) = q_1(u)$. Now, $p_3 q_1 (w) = q_3 p_2(w)$ implies that $p_3 q_1 (u) = q_3 q_1 (v)$.
    \end{itemize}
It is not hard to show that statement (*) extends to any pair $(u,v)$ in $\eqr{\CospanToER(\tr{p_1}\tl{q_1}) \relcomp \CospanToER(\tr{p_2}\tl{q_2})}$. Indeed, such membership is witnessed by a chain $w_1,\dots,w_n$ where $w_1 = u$, $w_n = w$ and for all $i \ls n$ either $(w_i,w_{i+1})$ or $(w_{i+1},w_i)$ is in $\eqr{\CospanToER(\tr{p_1}\tl{q_1}) \relcomp \CospanToER(\tr{p_2}\tl{q_2})}$. By (*) this means that $w_i$ and $w_{i+1}$ are mapped into the same element of $z_3$. It follow that also the first element $u$ and the last element $w$ of the chain are mapped into the same element of $z_3$, that is:
\begin{align*}
p_3 p_1 (u) &= q_3 q_1 (v) && \text{if } u \in \ord{n}, v \in \ord{m} \\
p_3 p_1 (u) &= p_3 p_1 (v) && \text{if } u,v \in \ord{n} \\
q_3 q_1 (u) &= p_3 p_1 (v) && \text{if } u \in \ord{m}, v \in \ord{n} \\
q_3 q_1 (u) &= q_3 q_1 (v) && \text{otherwise, } u, v \in \ord{m}.
\end{align*}
In any of the cases above, $(u,v) \in \CospanToER(\tr{p_3 p_1} \tl{q_3 q_1})$.
\end{proof}

\section{Proofs of Chapter~3}

This appendix gives more details on the equational theories of $\IBRw$ and $\IBR$ which have not been included in the main text for space reasons.

\subsection{The Frobenius Laws in $\IBRw$}\label{AppFrob} The presence of Frobenius laws both for the white --- \eqref{eq:WFrob} --- and for the black structure --- \eqref{eq:BFrob} --- make valid various deformations of the internal topology of circuits of $\IBRw$, as long as the connections between boundaries are preserved. We list here some useful laws of that kind. In describing the derivation steps, we occasionally use the notation $(n)^{op}$, which means the counterpart in $\ABRop$ of a valid equation $(n)$ in $\ABR$.
\begin{equation}\label{eq:Bfrobcomult}\tag{Fr1}
\lower11pt\hbox{$\includegraphics[height=1cm]{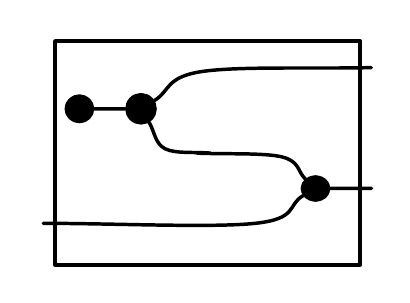}$}
\eql{\eqref{eq:BFrob}}
\lower11pt\hbox{$\includegraphics[height=1cm]{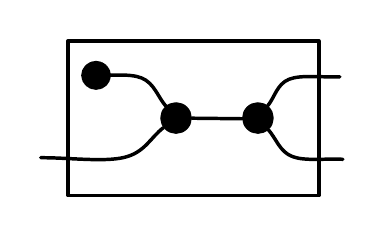}$}
\eql{\eqref{eq:bcomoncomm},\eqref{eq:bcomonunitlaw}}
\lower6pt\hbox{$\includegraphics[height=.6cm]{graffles/Bcomult.pdf}$}
\eql{\eqref{eq:bcomonunitlaw}}
\lower11pt\hbox{$\includegraphics[height=1cm]{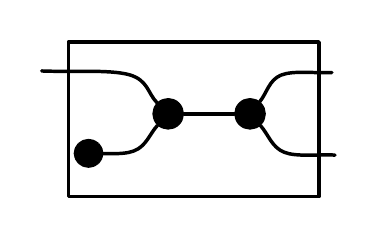}$}
\eql{\eqref{eq:BFrob}}
\lower11pt\hbox{$\includegraphics[height=1cm]{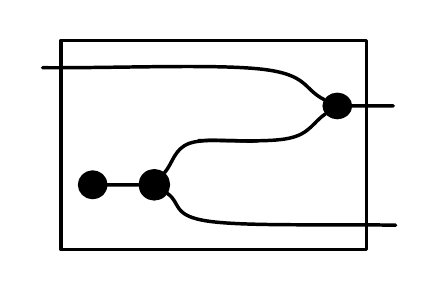}$}
\end{equation}
\begin{equation}\label{eq:Bsnake}\tag{Fr2}
\lower11pt\hbox{$\includegraphics[height=1cm]{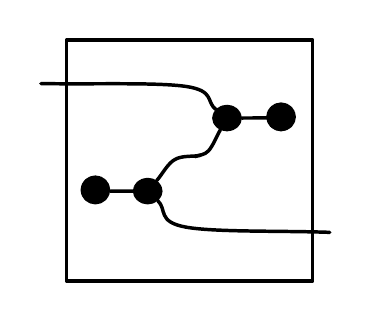}$}
\!\!\!\!\eql{\eqref{eq:BFrob}}\!\!\!\!
\lower10pt\hbox{$\includegraphics[height=1cm]{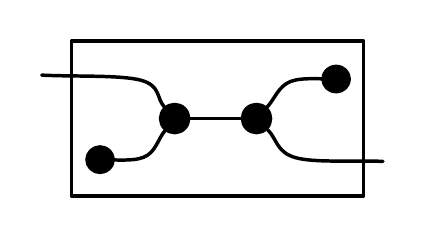}$}
\!\!\!\!\eql{\eqref{eq:bcomoncomm},\eqref{eq:bcomonunitlaw},\eqref{eq:bcomonunitlaw}$^{op}$}\!\!\!\!
\lower4pt\hbox{$\includegraphics[height=.6cm]{graffles/idcircuit.pdf}$}
\!\!\!\!\eql{\eqref{eq:bcomonunitlaw},\eqref{eq:bcomoncomm}$^{op}$,\eqref{eq:bcomonunitlaw}$^{op}$}\!\!\!\!
\lower10pt\hbox{$\includegraphics[height=1cm]{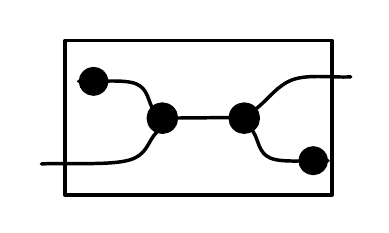}$}
\!\!\!\!\eql{\eqref{eq:BFrob}}\!\!\!\!
\lower11pt\hbox{$\includegraphics[height=1cm]{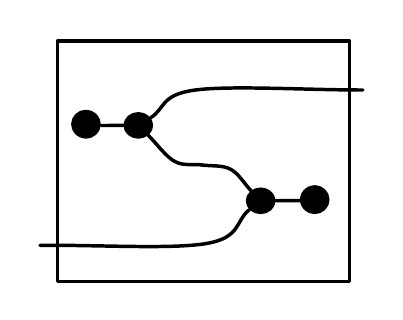}$}
\end{equation}

The following laws are derived analogously. The ones involving the white structure use the white Frobenius axiom~\eqref{eq:WFrob}.

\begin{multicols}{2}\noindent
\begin{equation}\label{eq:Bfrobmult}\tag{Fr3}
\lower11pt\hbox{$\includegraphics[height=1cm]{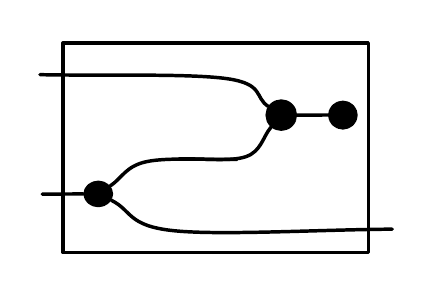}$} =
\lower6pt\hbox{$\includegraphics[height=.6cm]{graffles/Bmult.pdf}$} =
\lower11pt\hbox{$\includegraphics[height=1cm]{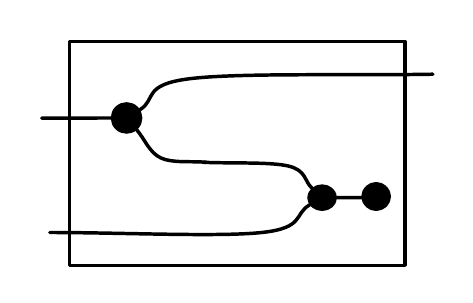}$}
\end{equation}
\begin{equation}\label{eq:Wsnake}\tag{Fr4}
\lower11pt\hbox{$\includegraphics[height=1cm]{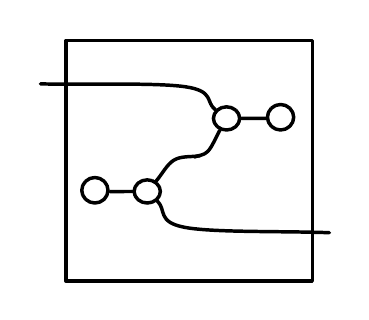}$} =
\lower6pt\hbox{$\includegraphics[height=.6cm]{graffles/idcircuit.pdf}$} =
\lower11pt\hbox{$\includegraphics[height=1cm]{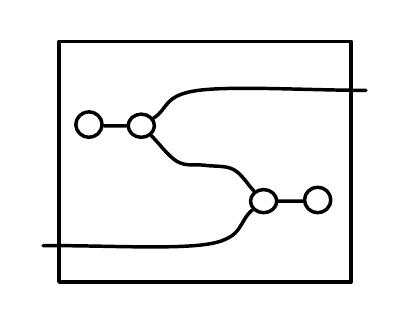}$}
\end{equation}
\end{multicols}
\begin{multicols}{2}\noindent
\begin{equation}\label{eq:Wfrobcomult}\tag{Fr5}
\lower11pt\hbox{$\includegraphics[height=1cm]{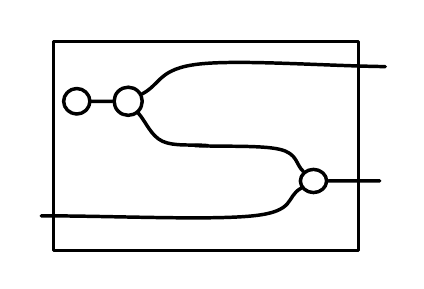}$} =
\lower6pt\hbox{$\includegraphics[height=.6cm]{graffles/Wcomult.pdf}$} =
\lower11pt\hbox{$\includegraphics[height=1cm]{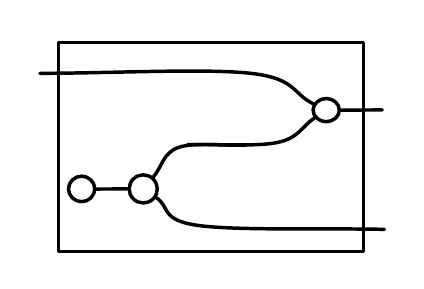}$}
\end{equation}
\begin{equation}\label{eq:Wfrobmult}\tag{Fr6}
\lower11pt\hbox{$\includegraphics[height=1cm]{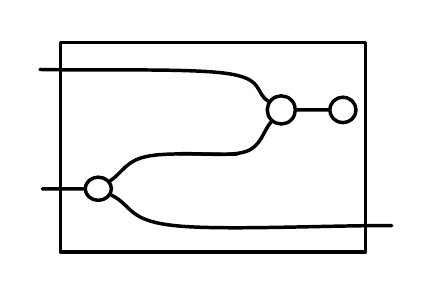}$} =
\lower11pt\hbox{$\includegraphics[height=1cm]{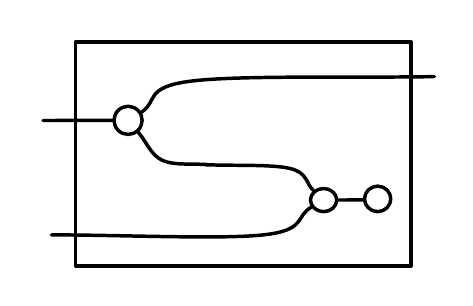}$} =
\lower6pt\hbox{$\includegraphics[height=.6cm]{graffles/Wmult.pdf}$}
\end{equation}
\end{multicols}

\noindent For later reference, we also record the following derivation.
\begin{equation}\label{eq:wccantipodesquare}\tag{Fr7}
\lower10pt\hbox{$\includegraphics[height=1cm]{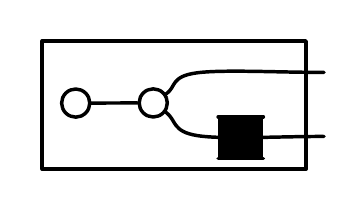}$} \eql{\eqref{eq:scalarwunit}}
\lower10pt\hbox{$\includegraphics[height=1cm]{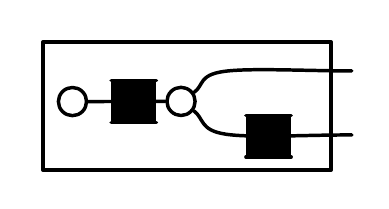}$} \eql{\eqref{eq:scalarwmult}$^{\op}$}
\lower10pt\hbox{$\includegraphics[height=1cm]{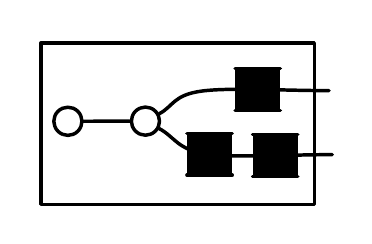}$}
\eql{\eqref{eq:scalarmult}}
\lower10pt\hbox{$\includegraphics[height=1cm]{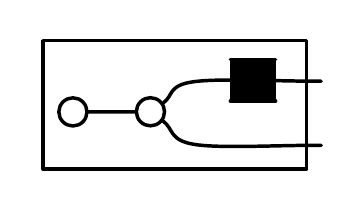}$}
\end{equation}
The same equation reflected about the $y$-axis and the black counterparts are proven analogously.
\begin{multicols}{3}\noindent
\begin{equation}\label{eq:lwccantipodesquare}\tag{Fr8}
\lower8pt\hbox{$\includegraphics[height=.8cm]{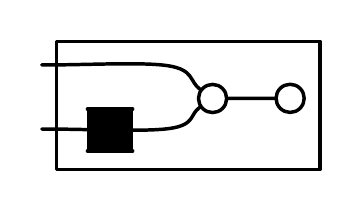}$} \!\!\!=\!\!\!
\lower8pt\hbox{$\includegraphics[height=.8cm]{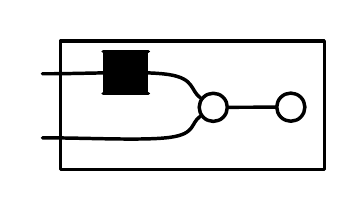}$}
\end{equation}
\begin{equation}\label{eq:bccantipodesquare}\tag{Fr9}
\lower8pt\hbox{$\includegraphics[height=.8cm]{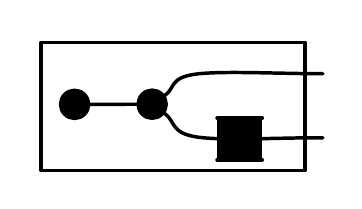}$} \!\!\!=\!\!\!
\lower8pt\hbox{$\includegraphics[height=.8cm]{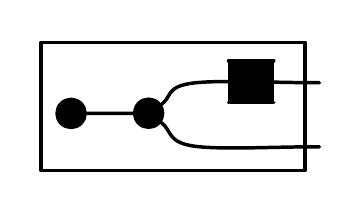}$}
\end{equation}
\begin{equation}\label{eq:lbccantipodesquare}\tag{Fr10}
\lower8pt\hbox{$\includegraphics[height=.8cm]{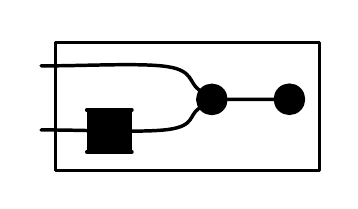}$} \!\!\!=\!\!\!
\lower8pt\hbox{$\includegraphics[height=.8cm]{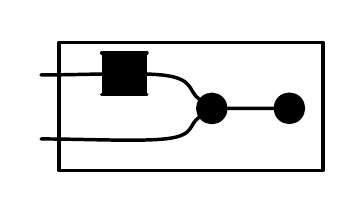}$}
\end{equation}
\end{multicols} 

\subsection{Derived Laws of $\IBRw$}\label{AppDerLaws} In this section we supply the equational proofs of the laws stated in Section~\ref{sec:ibrw}. We begin with the derivations of \eqref{eq:lccb} and \eqref{eq:uniqueantipode}.
\begin{equation*}
\lower9pt\hbox{$\includegraphics[height=.8cm]{graffles/wccantipodel.pdf}$} \eql{\eqref{eq:rcc}}
\lower9pt\hbox{$\includegraphics[height=.8cm]{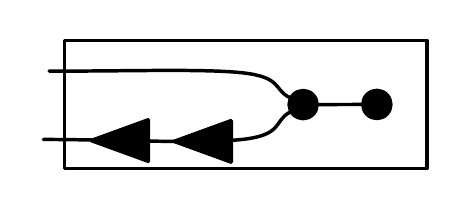}$} \eql{\eqref{eq:scalarmult}}
\lower8pt\hbox{$\includegraphics[height=.7cm]{graffles/bccl.pdf}$}
\end{equation*}
\begin{equation*}
\lower4pt\hbox{$\includegraphics[height=.5cm]{graffles/antipode.pdf}$} \eql{\eqref{eq:lcm}}
\lower9pt\hbox{$\includegraphics[height=.8cm]{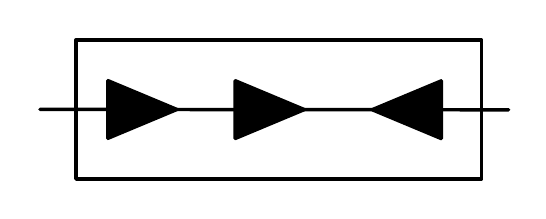}$} \eql{\eqref{eq:scalarmult}}
\lower4pt\hbox{$\includegraphics[height=.5cm]{graffles/antipodeop.pdf}$}
\end{equation*}
The derivation of \eqref{eq:rccb} is analogous to the one of \eqref{eq:lccb}, with \eqref{eq:lcc} used in place of \eqref{eq:rcc}. Now that \eqref{eq:uniqueantipode} has been proven, we follow the convention to write $\antipodesquare$ for both $\antipode$ and $\antipodeop$. We give next the derivation for~\eqref{eq:QFrob}:
\begin{equation*}
\lower60pt\hbox{$\includegraphics[height=4.8cm]{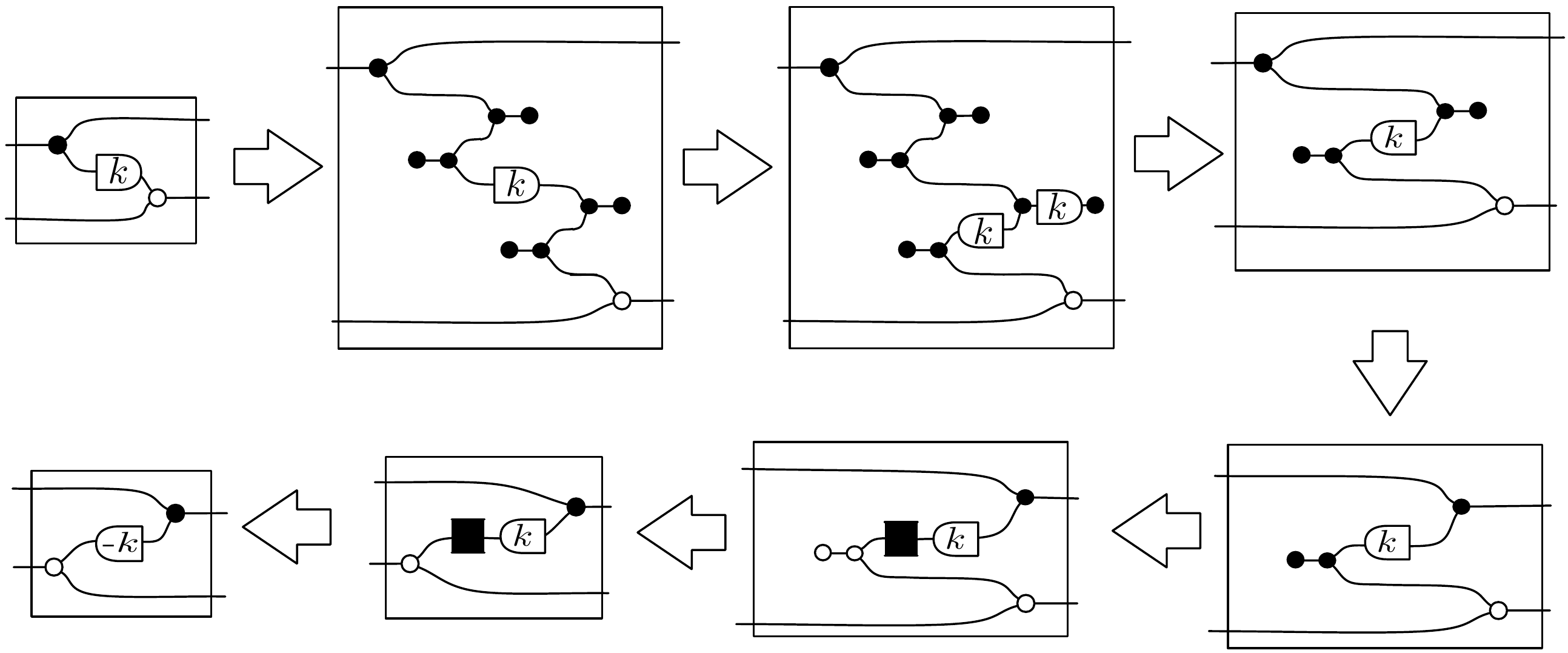}$}
\end{equation*}
The first step uses twice \eqref{eq:Bsnake}. The successive steps use: \eqref{eq:BccscalarAxiomOne}, \eqref{eq:scalarbcounit}, \eqref{eq:Bfrobmult}, \eqref{eq:lccb} and \eqref{eq:wccantipodesquare}, \eqref{eq:Wfrobcomult}, \eqref{eq:scalarmult}.

We show below the proof of \eqref{eq:coscalarwunit}, where $l \neq 0$. The ones for \eqref{eq:scalarwcounit} is symmetric.
\begin{equation*}
\lower6pt\hbox{$\includegraphics[height=.7cm]{graffles/Wunitcoscalarl.pdf}$} \eql{\eqref{eq:scalarwunit}}
\lower7pt\hbox{$\includegraphics[height=.8cm]{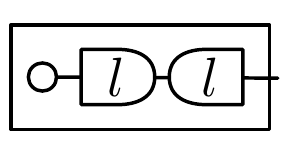}$}
 \eql{\eqref{eq:lcm}}
\lower6pt\hbox{$\includegraphics[height=.7cm]{graffles/Wunit.pdf}$}
\end{equation*}
Next we give the derivation of \eqref{eq:coscalarbcomult}, where $l \neq 0$. The one of \eqref{eq:scalarbmult} is analogous.
\begin{equation*}
\lower9pt\hbox{$\includegraphics[height=.9cm]{graffles/coscalarBcomult_der1.pdf}$} \!\!\eql{\eqref{eq:Bfrobcomult}}\!\!
\lower10pt\hbox{$\includegraphics[height=1cm]{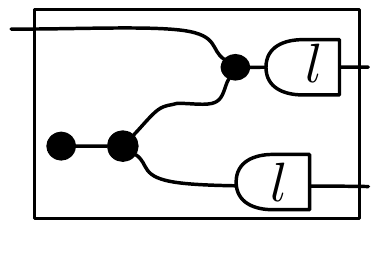}$} \!\!\eql{\eqref{eq:BccscalarAxiomTwo}}\!\!
\lower10pt\hbox{$\includegraphics[height=1cm]{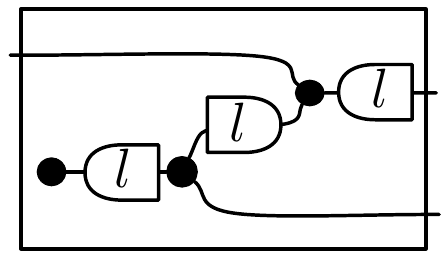}$} \!\!\eql{\eqref{eq:scalarbcomult}$^{\op}$}\!\!
\lower10pt\hbox{$\includegraphics[height=1cm]{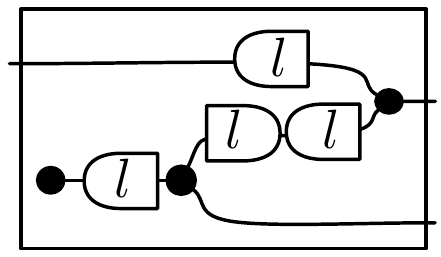}$} \!\!\eql{\eqref{eq:lcm}}\!\!
\lower10pt\hbox{$\includegraphics[height=1cm]{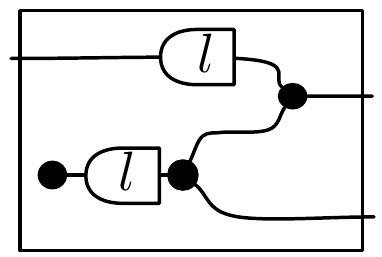}$} \!\!\eql{\eqref{eq:scalarbcounit}$^{\op}$}\!\!
\lower10pt\hbox{$\includegraphics[height=1cm]{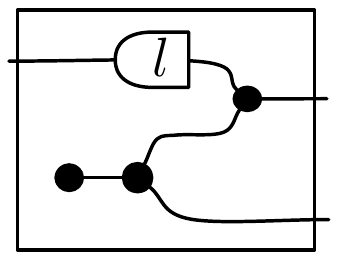}$} \!\!\eql{\eqref{eq:Bfrobcomult}}\!\!
\lower9pt\hbox{$\includegraphics[height=.9cm]{graffles/coscalarBcomult_der7.pdf}$}
\end{equation*}
We now consider the task of deriving law \eqref{eq:papillon}. For the first half:
\begin{center}
\includegraphics[height=4cm]{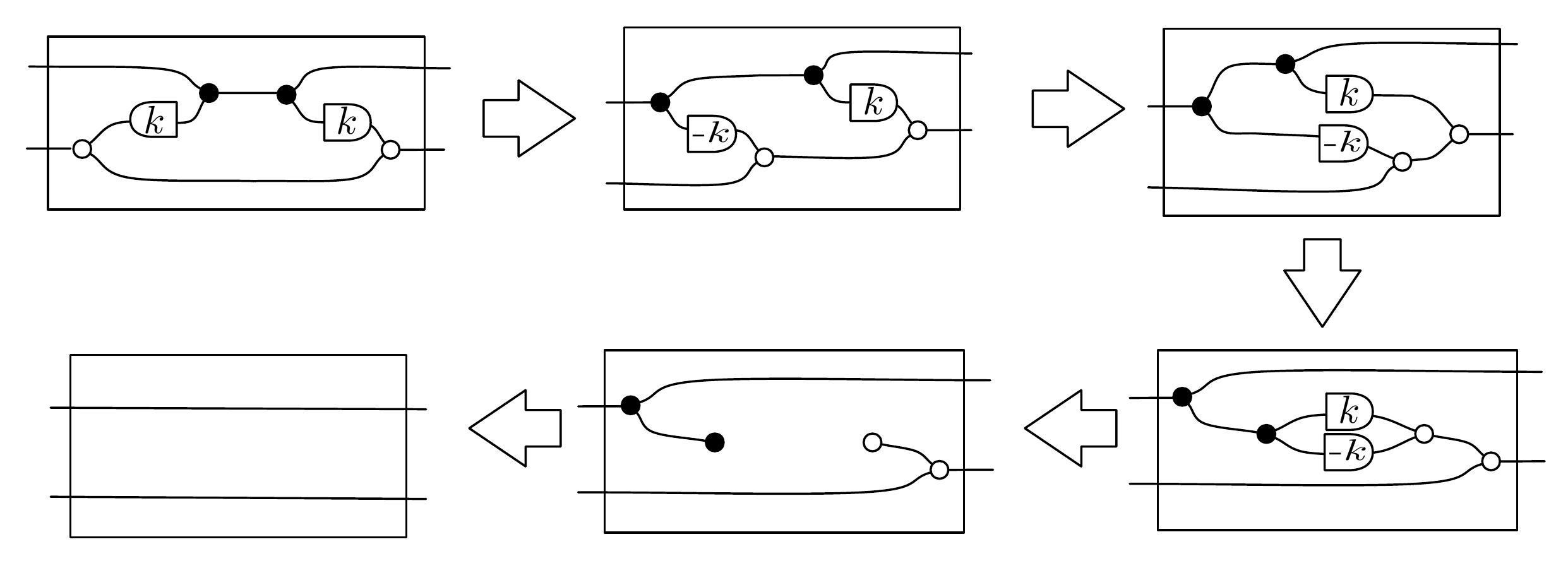}
\end{center}
The sequence of equations that are used is the following: \eqref{eq:QFrob}, \eqref{eq:sliding1}, \eqref{eq:bcomonassoc} and \eqref{eq:wmonassoc}, \eqref{eq:scalarsum} and \eqref{eq:zeroscalar}, \eqref{eq:bcomonunitlaw} and \eqref{eq:wmonunitlaw}.
The second half of \eqref{eq:papillon} is derived analogously as follows.
\begin{center}
\includegraphics[height=4cm]{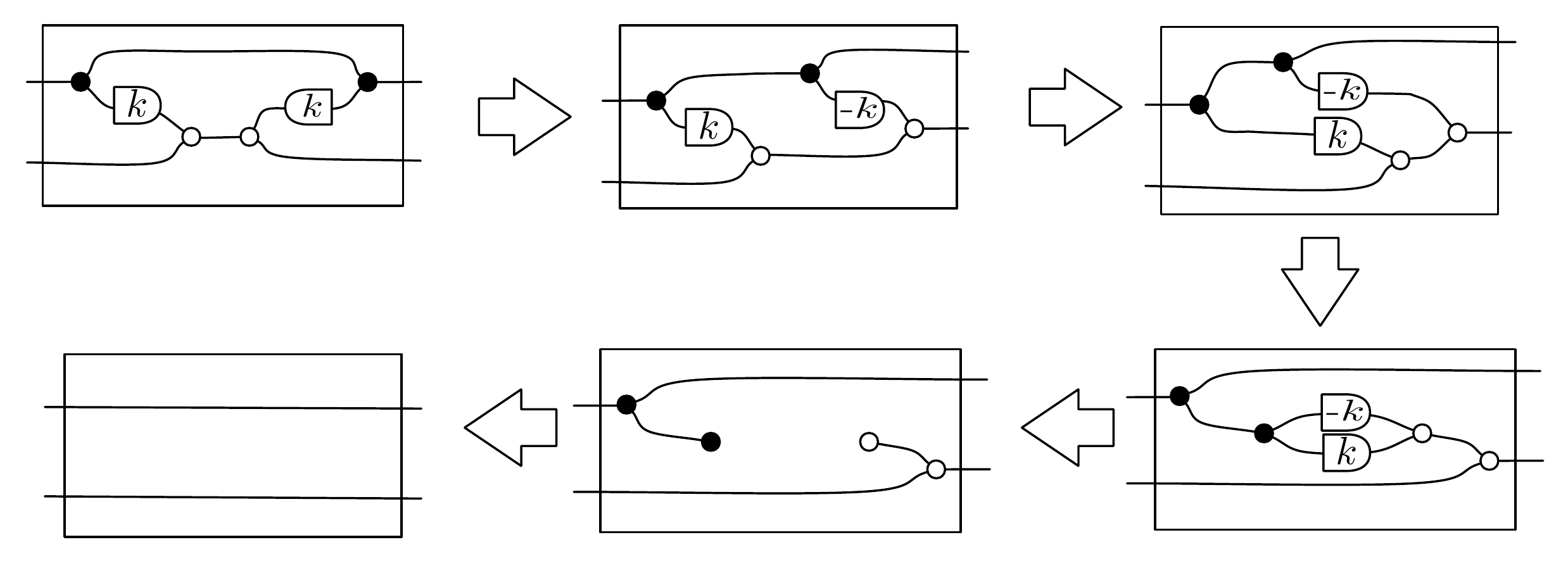}
\end{center}
In order to show the validity of \eqref{eq:wunitcancelbcomult}, we proceed by induction on the coarity $n \geq 1$ of the circuit, i.e., the number of gates on the right boundary. For the case $n = 1$, we have the following derivation, where $l \neq 0$.
\begin{equation}\label{eq:wcounitcancelbcomultder}
\lower40pt\hbox{$\includegraphics[height=3cm]{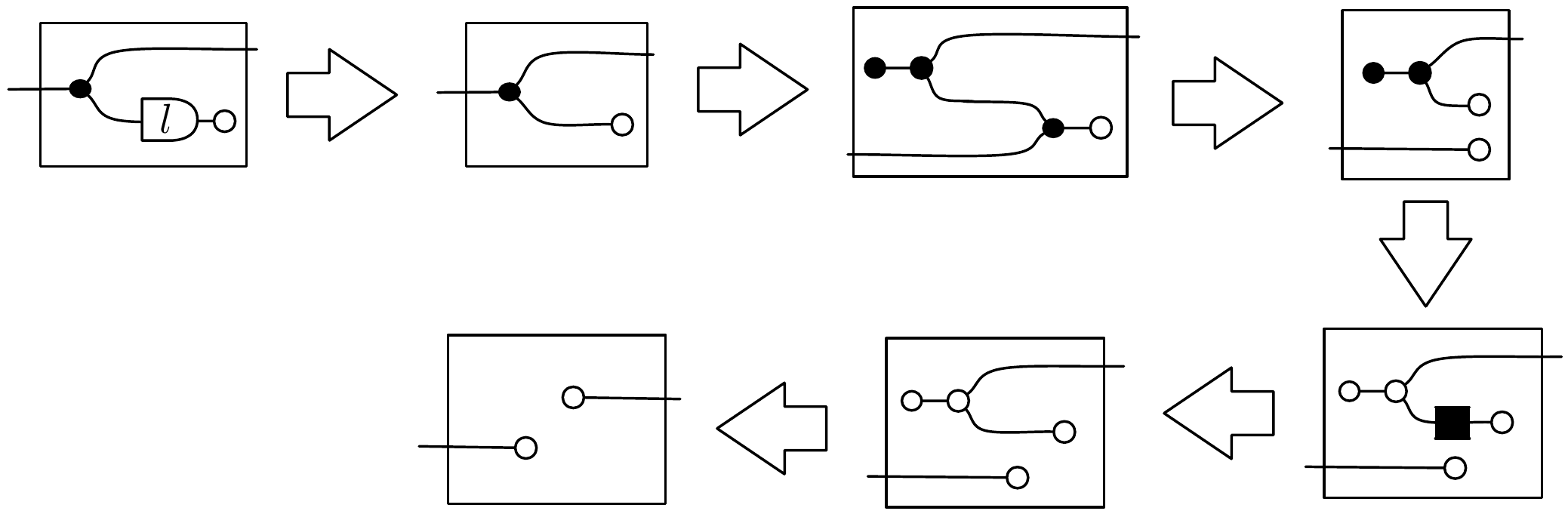}$}
\end{equation}
The sequence of applied laws is: \eqref{eq:scalarwcounit}, \eqref{eq:Bfrobcomult}, \eqref{eq:unitsr}$^{op}$, \eqref{eq:lccb}, \eqref{eq:scalarwunit}, \eqref{eq:wmonunitlaw}$^{op}$. The inductive case is handled as follows.
\begin{equation*}
\lower12pt\hbox{$\includegraphics[height=1.1cm]{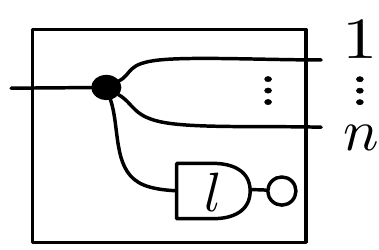}$} \eql{\eqref{eq:scalarwcounit}}
\lower12pt\hbox{$\includegraphics[height=1.1cm]{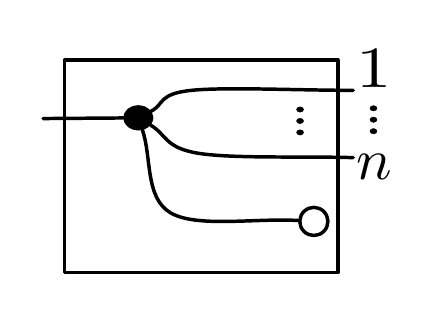}$} \eql{\eqref{eq:bcomonassoc}}
\lower12pt\hbox{$\includegraphics[height=1.1cm]{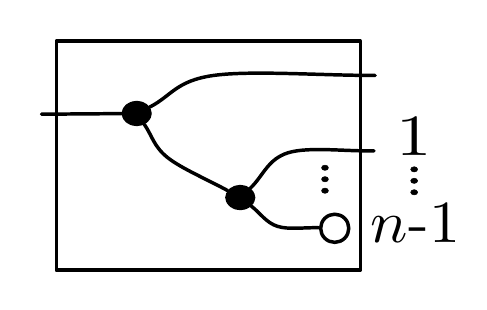}$}
\eql{Ind. hyp.}
\lower12pt\hbox{$\includegraphics[height=1.1cm]{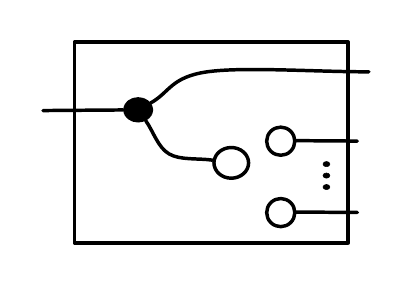}$}
\eql{\eqref{eq:wcounitcancelbcomultder}}
\lower10pt\hbox{$\includegraphics[height=1cm]{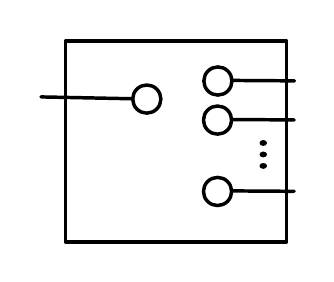}$}
\end{equation*}

\noindent Finally, we show the derivation for \eqref{eq:Bsep}. The sequence of applied laws is \eqref{eq:wbone}, \eqref{eq:bcomonunitlaw}+\eqref{eq:bcomonunitlaw}$^{\op}$, \eqref{eq:scalarsum}+\eqref{eq:scalarsum}$^{\op}$, \eqref{eq:bcomonassoc}+\eqref{eq:bcomonassoc}$^{\op}$, \eqref{eq:papillon}.
\begin{equation*}
\includegraphics[height=2cm]{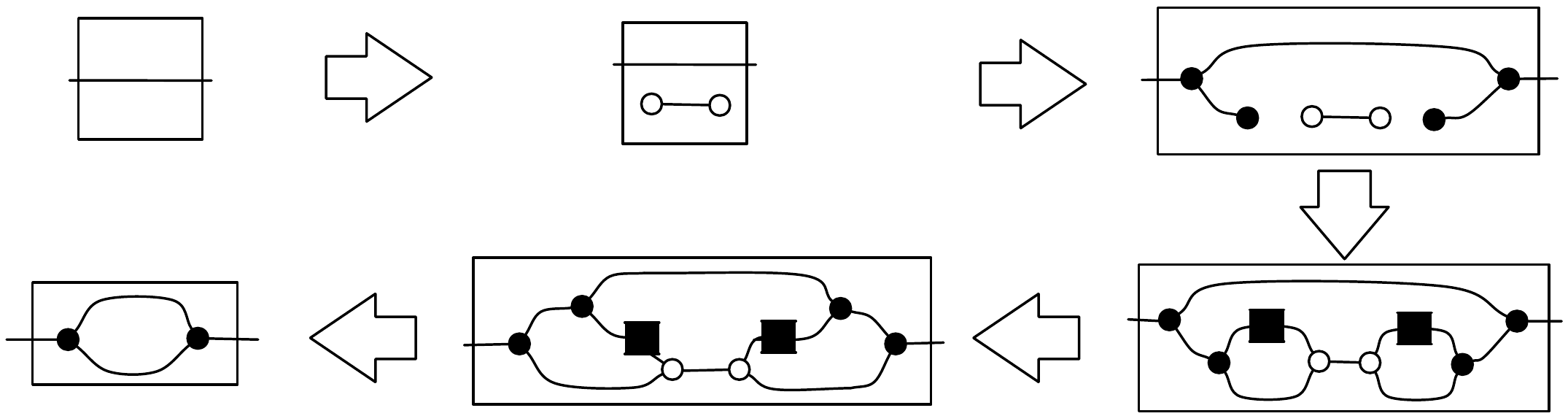}
\end{equation*}

\subsection{The Compact Closed Structure of $\IBRw$}\label{AppCC} We give more detailed proofs to the statements of Section~\ref{sec:cc}.

\begin{proof}[Proof of Proposition \ref{prop:snakecc}] We give the argument proving the left side of \eqref{eq:gensnake} --- the proof for the right side is completely symmetric. We proceed by induction on $n$. For the case $n = 1$, the statement is given by \eqref{eq:Bsnake}. For the inductive step, let $n = i+1$. In the sequel we show the equality
\begin{eqnarray}\label{eq:ccsnakeInd}
\lower15pt\hbox{$\includegraphics[height=1.3cm]{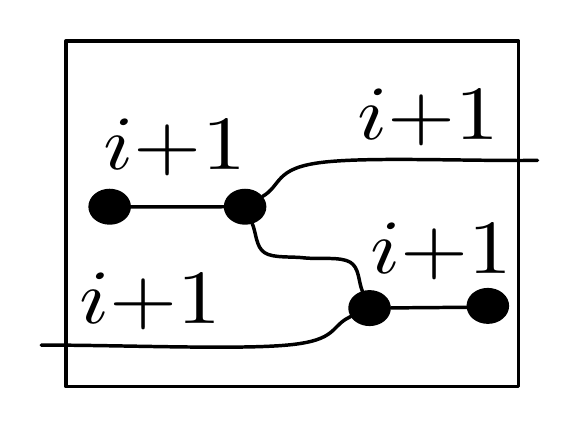}$}
&=&
\lower15pt\hbox{$\includegraphics[height=1.3cm]{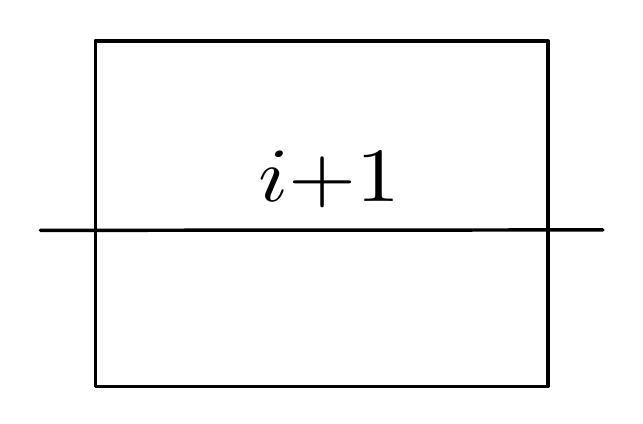}$} \end{eqnarray}
yielding the left side of \eqref{eq:gensnake}. For this purpose, we record the following derived law, allowing to``move'' the compact closed structure past the symmetries of $\IBRw$.
\begin{eqnarray}\label{eq:moveccpastsym}
\lower11pt\hbox{$\includegraphics[height=1cm]{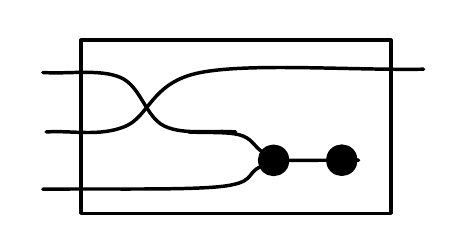}$}
&\eql{\eqref{eq:sliding2}}&
\lower11pt\hbox{$\includegraphics[height=1cm]{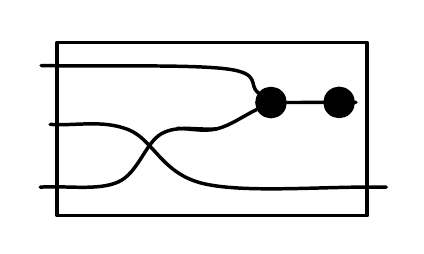}$} \end{eqnarray}
We can now proceed with the derivation of \eqref{eq:ccsnakeInd}. The diagram on the left side of \eqref{eq:ccsnakeInd} has the following shape.
\begin{center}
\includegraphics[height=2.6cm]{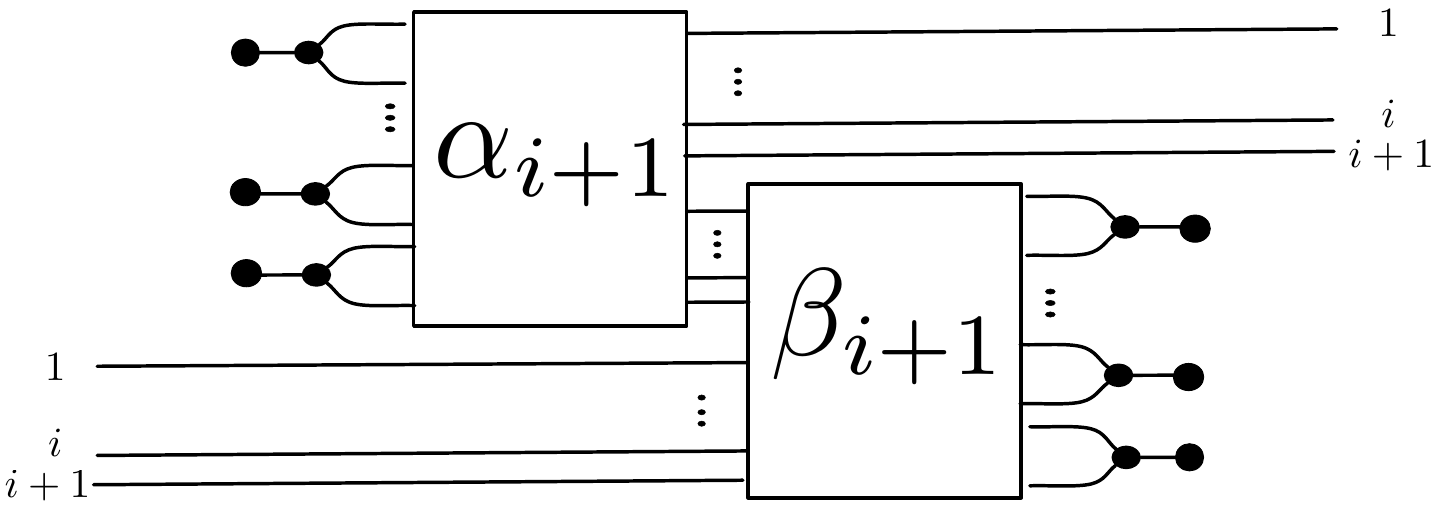}
\end{center}
By definition, port $1$ of the bottommost diagram $\rccB$ (call it $c_l$) connects to port $i+1$ on the right boundary and port $2$ connects to port $1$ of the bottommost diagram $\lccB$ (call it $c_r$). The other port of $c_r$ connects to port $i+1$ on the left boundary. By iteratively applying \eqref{eq:moveccpastsym} to $c_r$, we can move it towards the middle of the diagram, past all the symmetries in $\beta_{i+1}$. The resulting diagram is the following:
\begin{center}
\includegraphics[height=3cm]{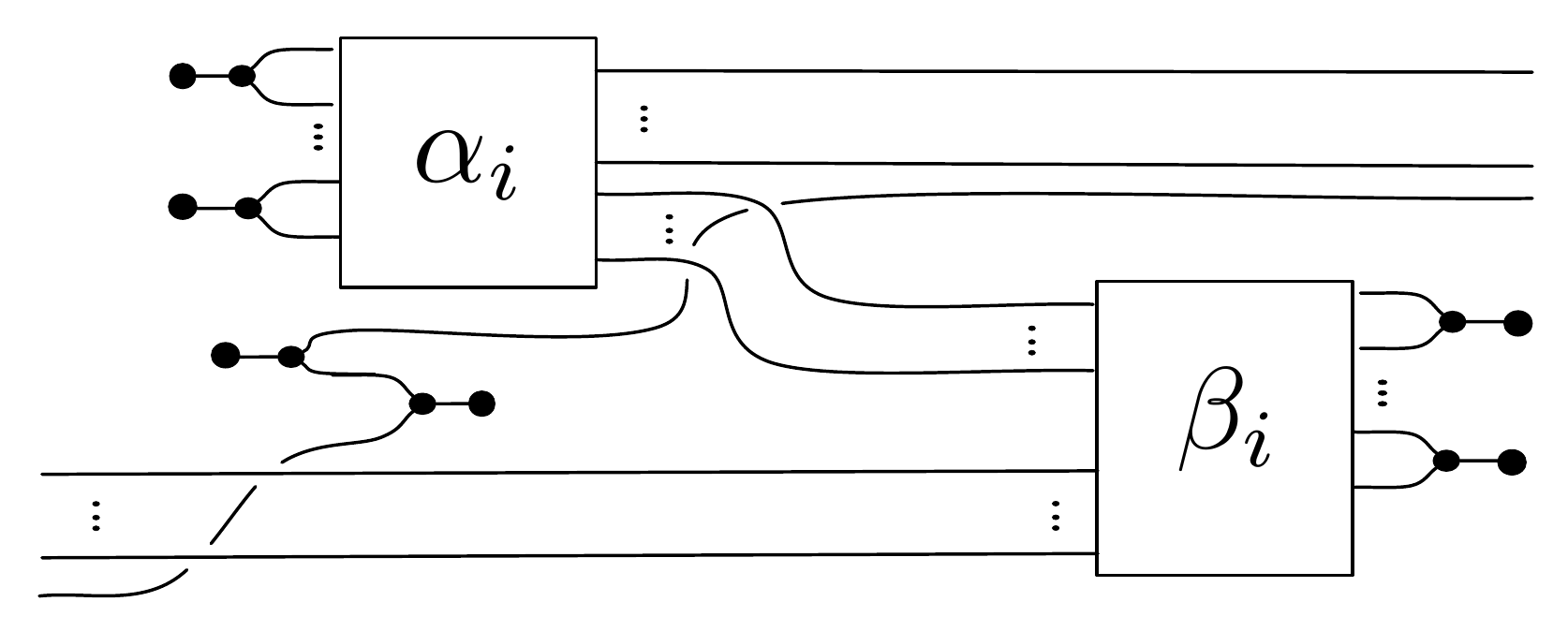}
\end{center}
Note that, now that we isolated $c_l$ and $c_r$, the diagrams $\alpha_{i+1}$ and $\beta_{i+1}$ become by definition $\alpha_{i}$ and $\beta_i$ --- also, the application of \eqref{eq:moveccpastsym} does not affect the arity of the symmetries in the diagram. We are now in position to apply \eqref{eq:Bsnake}:
\begin{center}
\includegraphics[height=3cm]{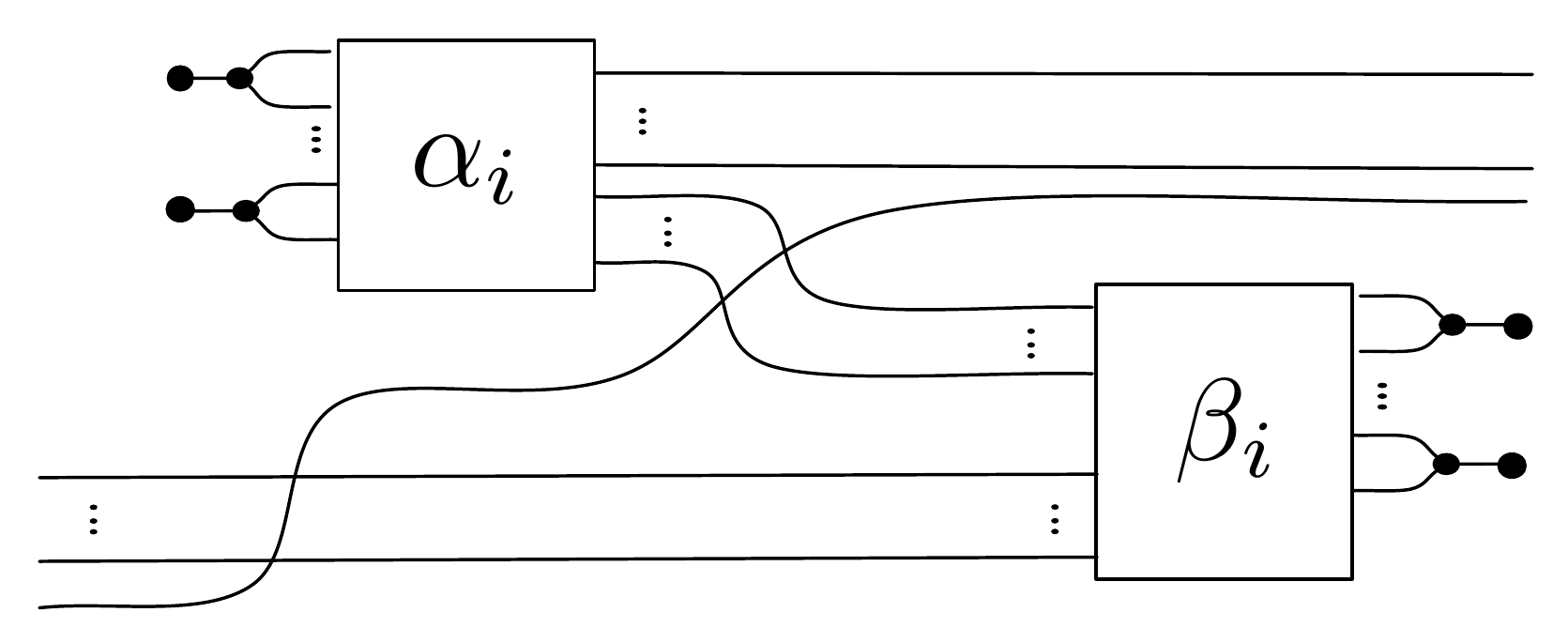}
\end{center}
The next step is to iteratively apply \eqref{eq:moveccpastsym} to move the middle identity wire towards the bottom.
\begin{center}
\includegraphics[height=3cm]{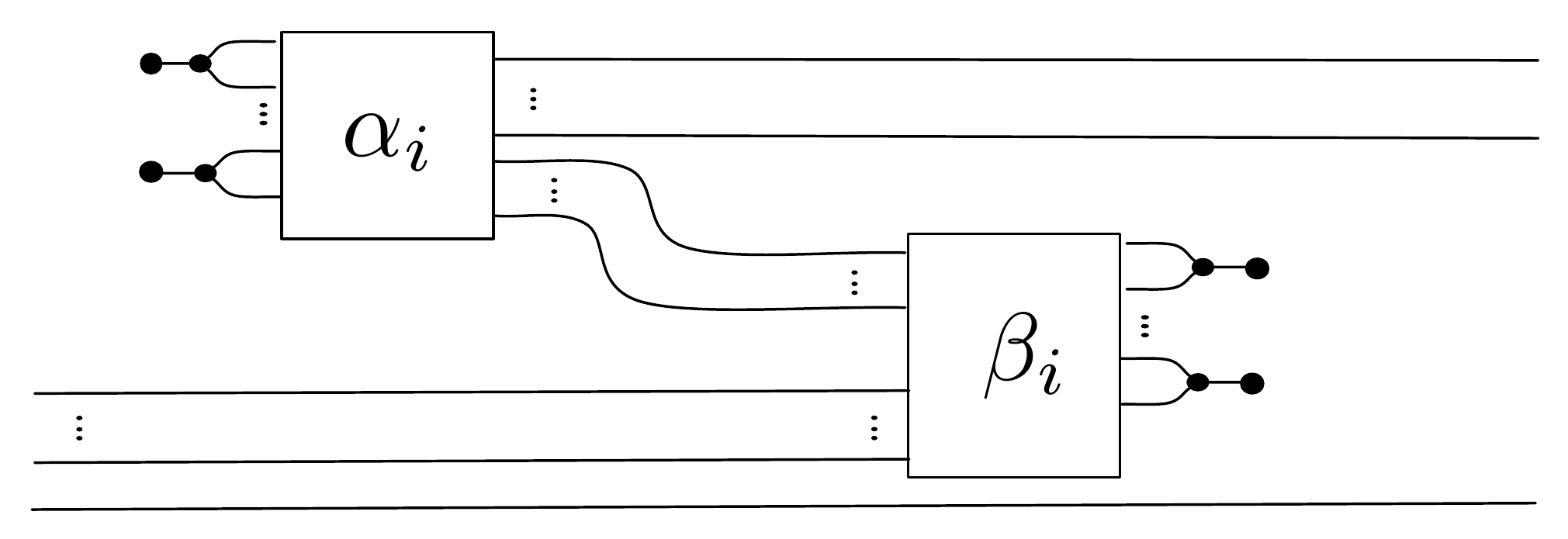}
\end{center}
It is now possible to apply the inductive hypothesis on $i$, obtaining as a result the desired identity diagram as on the right side of \eqref{eq:ccsnakeInd}.
\end{proof}

\begin{proof}[Proof of Proposition \ref{prop:star=refl}]
The proof is by induction on a string diagram $c$ of $\IBRw$. First we give the derivations for the four base cases of white/black unit/counit.
\begin{equation*}
\lower16pt\hbox{$\includegraphics[height=1.3cm]{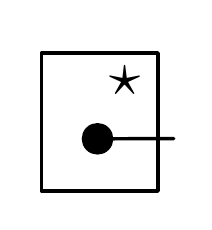}$}
\eql{Def. $\coc{\cdot}$}
\lower15pt\hbox{$\includegraphics[height=1.2cm]{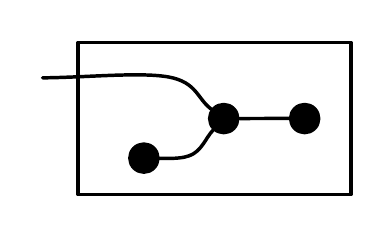}$}
\eql{\eqref{eq:bcomonunitlaw}$^{\op}$}
\ \lower8pt\hbox{$\includegraphics[height=.7cm]{graffles/Bcounit.pdf}$} \end{equation*}
\begin{equation*}
\lower15pt\hbox{$\includegraphics[height=1.3cm]{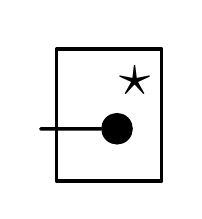}$}
\eql{Def. $\coc{\cdot}$}
\lower15pt\hbox{$\includegraphics[height=1.2cm]{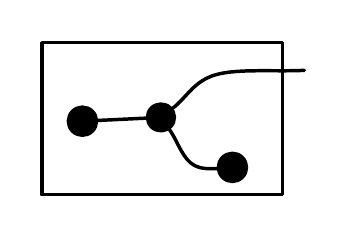}$}
\eql{\eqref{eq:bcomonunitlaw}}
\ \lower8pt\hbox{$\includegraphics[height=.7cm]{graffles/Bunit.pdf}$} \end{equation*}
\begin{equation*}
\lower15pt\hbox{$\includegraphics[height=1.3cm]{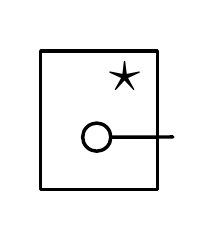}$}
\!\eql{Def. $\coc{\cdot}$}\!
\lower15pt\hbox{$\includegraphics[height=1.2cm]{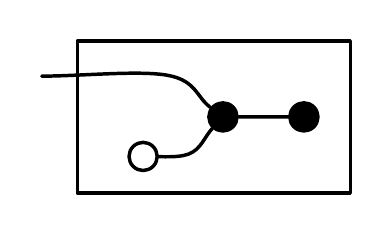}$}
\!\eql{\eqref{eq:lccb},\eqref{eq:lwccantipodesquare}}\!
\lower15pt\hbox{$\includegraphics[height=1.2cm]{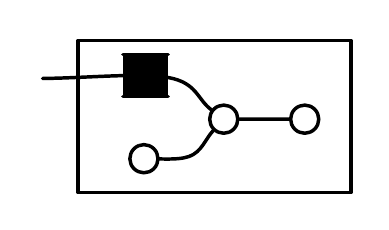}$}
\!\eql{\eqref{eq:wmonunitlaw}}\!
\lower15pt\hbox{$\includegraphics[height=1.2cm]{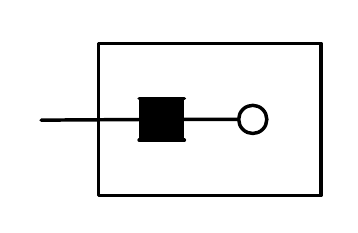}$}
\!\eql{\eqref{eq:scalarwunit}}\!
\ \lower8pt\hbox{$\includegraphics[height=.7cm]{graffles/Wcounit.pdf}$} \end{equation*}
\begin{equation*}
\lower15pt\hbox{$\includegraphics[height=1.3cm]{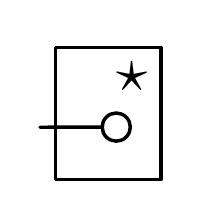}$}
\!\eql{Def. $\coc{\cdot}$}\!
\lower15pt\hbox{$\includegraphics[height=1.2cm]{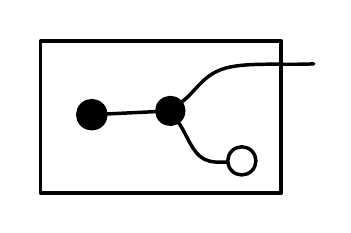}$}
\!\eql{\eqref{eq:rccb}}\!
\lower15pt\hbox{$\includegraphics[height=1.2cm]{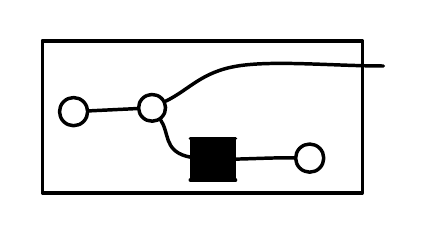}$}
\!\eql{\eqref{eq:scalarwunit}$^{\op}$}\!
\lower15pt\hbox{$\includegraphics[height=1.2cm]{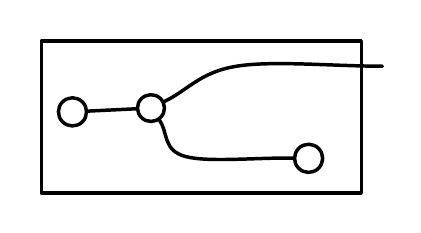}$}
\!\eql{\eqref{eq:wmonunitlaw}$^{\op}$}\!
\ \lower8pt\hbox{$\includegraphics[height=.7cm]{graffles/Wunit.pdf}$} \end{equation*}

We now consider the base cases $\scalar$ and $\coscalar$, for $k \in \PID$.
\begin{equation*}
\lower9pt\hbox{$\includegraphics[height=.8cm]{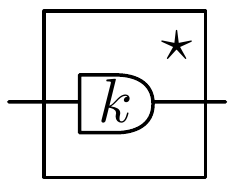}$}
\!\eql{Def. $\coc{\cdot}$}\!
\lower14pt\hbox{$\includegraphics[height=1.3cm]{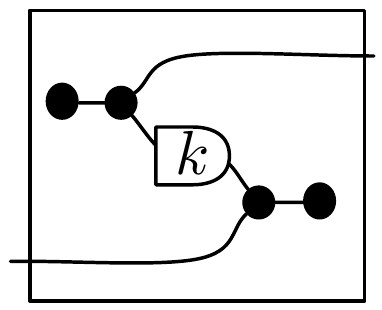}$}
\!\eql{\eqref{eq:BccscalarAxiomOne}}\!
\lower14pt\hbox{$\includegraphics[height=1.3cm]{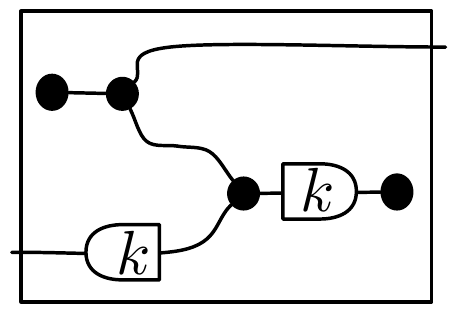}$}
\!\eql{\eqref{eq:scalarbcounit}}\!
\lower14pt\hbox{$\includegraphics[height=1.3cm]{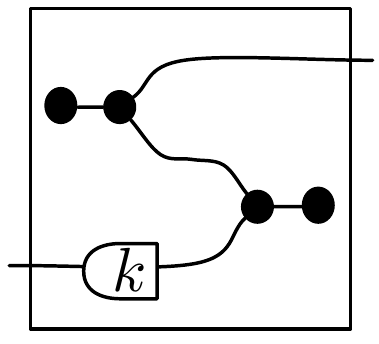}$}
\!\eql{\eqref{eq:Bsnake}}\!
\ \lower9pt\hbox{$\includegraphics[height=.9cm]{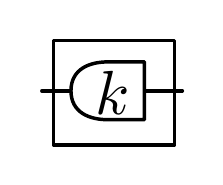}$}
\end{equation*}
\begin{equation*}
\lower9pt\hbox{$\includegraphics[height=.8cm]{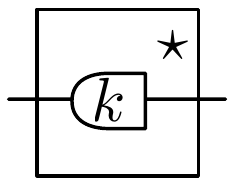}$}
\!\eql{Def. $\coc{\cdot}$}\!
\lower14pt\hbox{$\includegraphics[height=1.3cm]{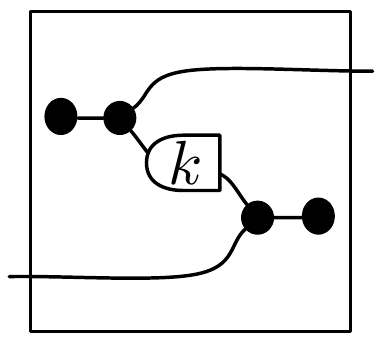}$}
\!\eql{\eqref{eq:BccscalarAxiomTwo}}\!
\lower14pt\hbox{$\includegraphics[height=1.3cm]{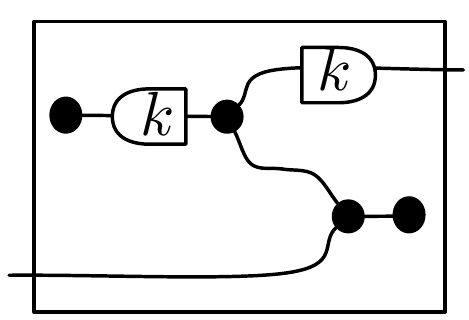}$}
\!\eql{\eqref{eq:scalarbcounit}$^{\op}$}\!
\lower14pt\hbox{$\includegraphics[height=1.3cm]{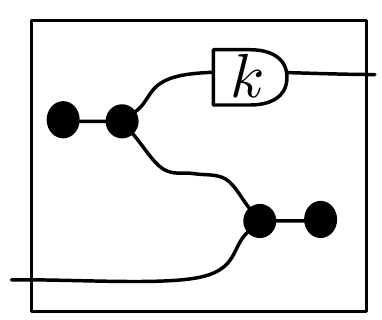}$}
\!\eql{\eqref{eq:Bsnake}}\!
\ \lower9pt\hbox{$\includegraphics[height=.9cm]{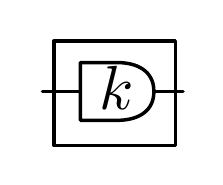}$}
\end{equation*}

We also provide the derivation for the base case $\Bcomult$.
\begin{center}
\includegraphics[height=5.5cm]{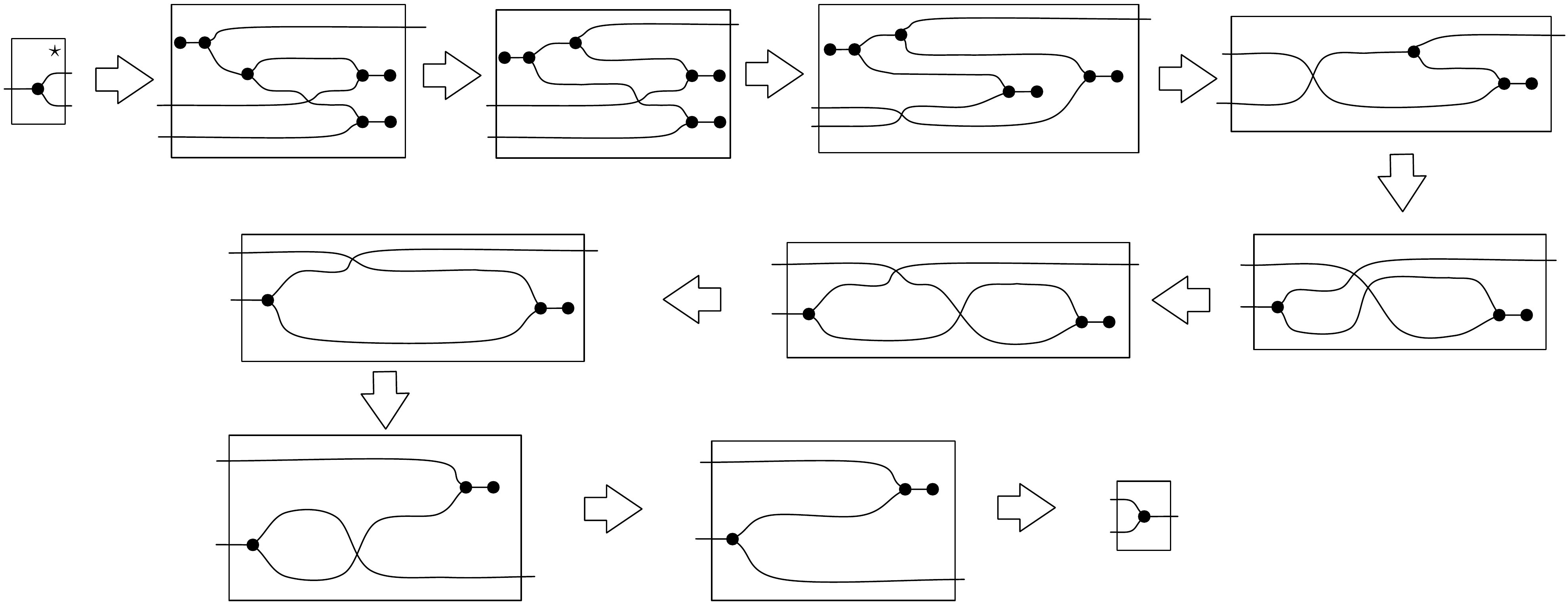}
\end{center}
The sequence of applied laws is: definition of $\coc{\cdot}$, \eqref{eq:bcomonassoc}, \eqref{eq:moveccpastsym}, \eqref{eq:Bsnake}, \eqref{eq:sliding2}, \eqref{eq:sliding1}, \eqref{eq:bcomoncomm}$^{\op}$, \eqref{eq:moveccpastsym}, \eqref{eq:bcomoncomm}, \eqref{eq:Bfrobmult}.

The remaining base cases of generators $\Bcomult$, $\Wmult$ and $\Wcomult$ are handled in an analogous way by using the Frobenius laws derived in \S~\ref{AppFrob}. The proof is concluded by examining the two inductive cases. For sequential composition:
   \begin{equation*}    \lower10pt\hbox{$\includegraphics[height=1cm]{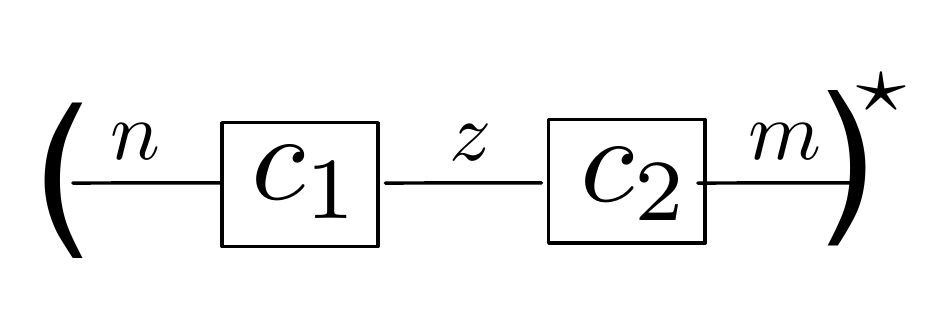}$}
  \!\eql{}\!
   \lower10pt\hbox{$\includegraphics[height=.9cm]{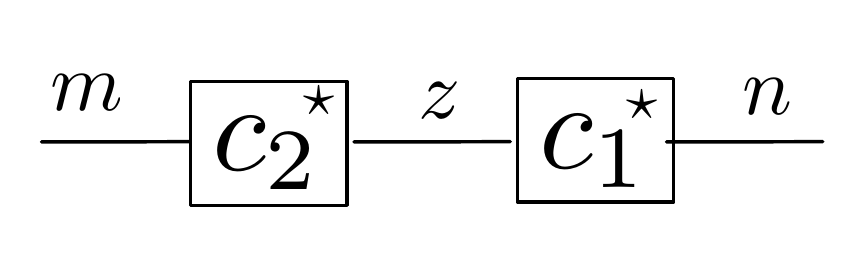}$}
  \!\eql{Ind. hyp.}\!
   \lower10pt\hbox{$\includegraphics[height=.9cm]{graffles/reflcompr.pdf}$}
  \!\eql{}\!
   \lower10pt\hbox{$\includegraphics[height=1cm]{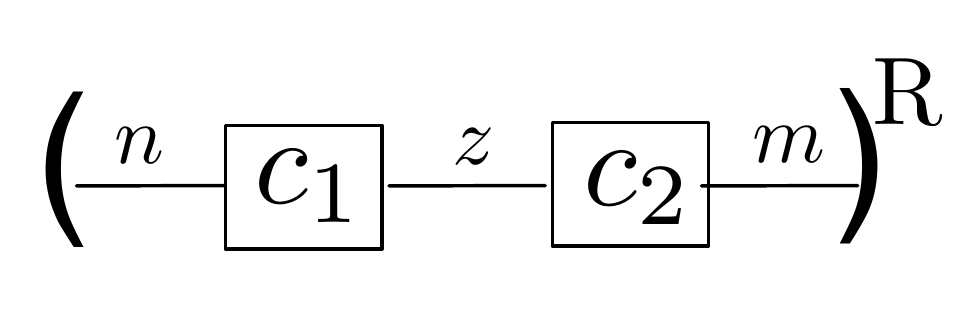}$}
   \end{equation*}
The derivation for the case of parallel composition $\tns$ is analogous.\end{proof} 

\subsection{Derived Laws of $\IBR$}\label{AppDerLawsIH} We conclude the proof of Proposition~\ref{prop:topfacepushout} by verifying that \eqref{eq:wbone}, \eqref{eq:BccscalarAxiomOne}, \eqref{eq:BccscalarAxiomTwo}, \eqref{eq:bbone}, \eqref{eq:WcccoscalarAxiomOne} and \eqref{eq:WcccoscalarAxiomTwo} are all derivable in $\IBR$. The following is the derivation of \eqref{eq:wbone}.
\begin{equation*}
\lower7pt\hbox{$\includegraphics[height=.7cm]{graffles/idzerocircuit.pdf}$}
\eql{\eqref{eq:bwbone}}
\lower8pt\hbox{$\includegraphics[height=.8cm]{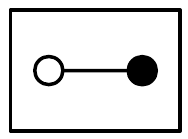}$}
\eql{\eqref{eq:BSepIBR}}
\lower9pt\hbox{$\includegraphics[height=.9cm]{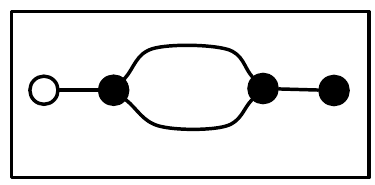}$}
\eql{\eqref{eq:lccb}}
\lower9pt\hbox{$\includegraphics[height=.9cm]{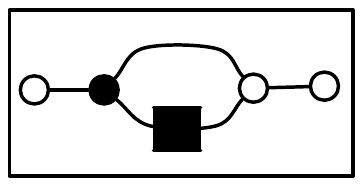}$}
\eql{\eqref{eq:scalarsum}}
  \lower8pt\hbox{$\includegraphics[height=.8cm]{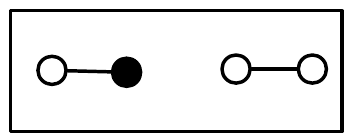}$}
\eql{\eqref{eq:bwbone}}
  \lower6pt\hbox{$\includegraphics[height=.6cm]{graffles/WBone.pdf}$}
\end{equation*}
The derivation of~\eqref{eq:bbone} is the ``photografic negative'' of the one of \eqref{eq:wbone}. We now show the derivations for \eqref{eq:BccscalarAxiomOne} and \eqref{eq:WcccoscalarAxiomOne}. For $l \neq 0$:
\begin{equation*}
\lower9pt\hbox{$\includegraphics[height=.8cm]{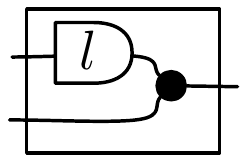}$}
\eql{\eqref{eq:lcmopIH}}
\lower9pt\hbox{$\includegraphics[height=.8cm]{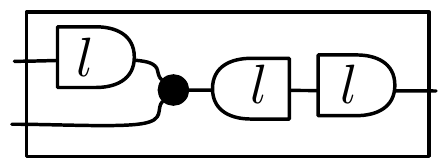}$}
\eql{\eqref{eq:scalarbcomult}$^{\op}$}
\lower9pt\hbox{$\includegraphics[height=.8cm]{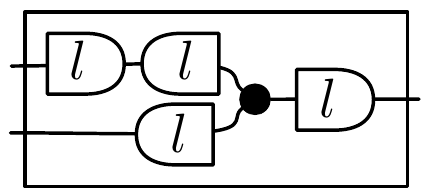}$}
\eql{\eqref{eq:lcmIH}}
\lower9pt\hbox{$\includegraphics[height=.8cm]{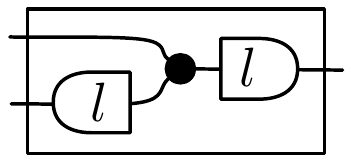}$}
\end{equation*}\noindent
\begin{equation*}
\lower9pt\hbox{$\includegraphics[height=.8cm]{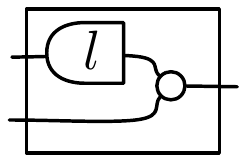}$}
\eql{\eqref{eq:lcmopIH}}
\lower9pt\hbox{$\includegraphics[height=.8cm]{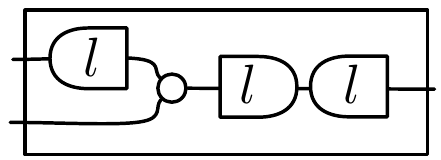}$}
\eql{\eqref{eq:scalarwmult}}
\lower9pt\hbox{$\includegraphics[height=.8cm]{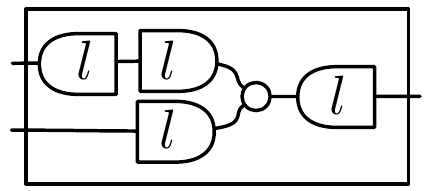}$}
\eql{\eqref{eq:lcmopIH}}
\lower9pt\hbox{$\includegraphics[height=.8cm]{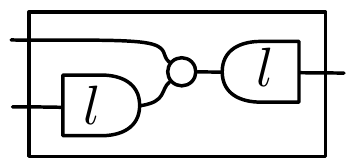}$}.
\end{equation*}
The zero cases:
\begin{equation*}
\lower9pt\hbox{$\includegraphics[height=.8cm]{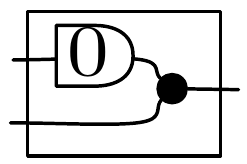}$}
\!\!\eql{\eqref{eq:zeroscalar}}\!\!
\lower9pt\hbox{$\includegraphics[height=.8cm]{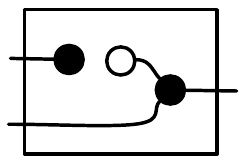}$}
\!\!\eql{\eqref{eq:Bfrobmult},\eqref{eq:lccb}}\!\!
\lower12pt\hbox{$\includegraphics[height=1.1cm]{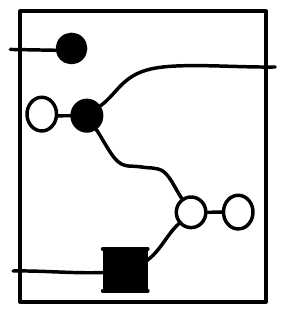}$}
\!\!\eql{\eqref{eq:unitsr}}\!\!
\lower10pt\hbox{$\includegraphics[height=.9cm]{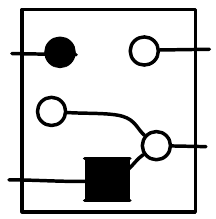}$}
\!\!\eql{\eqref{eq:wmonunitlaw},\eqref{eq:scalarwunit}}\!\!
\lower9pt\hbox{$\includegraphics[height=.8cm]{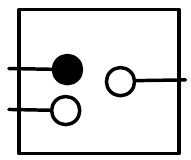}$}
\!\!\eql{\eqref{eq:bcomonunitlaw}}\!\!
\lower9pt\hbox{$\includegraphics[height=.8cm]{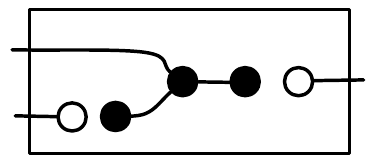}$}
\!\!\eql{\eqref{eq:zeroscalar},\eqref{eq:zeroscalar}$^{\op}$}\!\!
\lower9pt\hbox{$\includegraphics[height=.8cm]{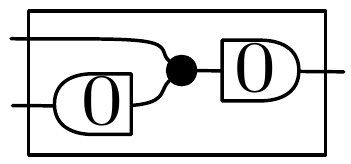}$}
\end{equation*}\noindent
\begin{equation*}
\lower9pt\hbox{$\includegraphics[height=.8cm]{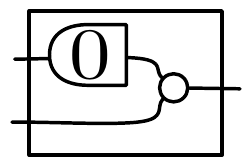}$}
\!\!\eql{\eqref{eq:zeroscalar}$^{\op}$}\!\!
\lower9pt\hbox{$\includegraphics[height=.8cm]{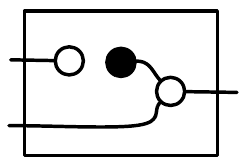}$}
\!\!\eql{\eqref{eq:Wfrobmult},\eqref{eq:rcc}}\!\!
\lower12pt\hbox{$\includegraphics[height=1.1cm]{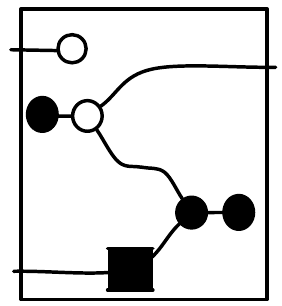}$}
\!\!\eql{\eqref{eq:unitsl}$^{\op}$}\!\!
\lower10pt\hbox{$\includegraphics[height=.9cm]{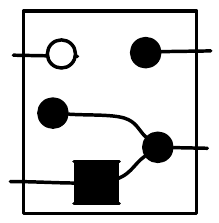}$}
\!\!\eql{\eqref{eq:bcomonunitlaw},\eqref{eq:scalarbcounit}$^{\op}$}\!\!
\lower9pt\hbox{$\includegraphics[height=.8cm]{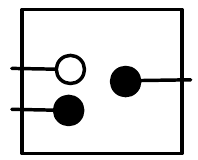}$}
\!\!\eql{\eqref{eq:wmonunitlaw}}\!\!
\lower9pt\hbox{$\includegraphics[height=.8cm]{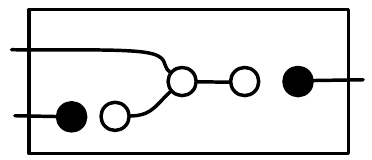}$}
\!\!\eql{\eqref{eq:zeroscalar},\eqref{eq:zeroscalar}$^{\op}$}\!\!
\lower9pt\hbox{$\includegraphics[height=.8cm]{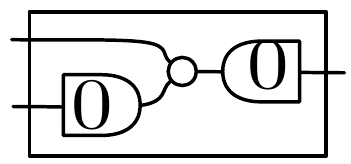}$}
\end{equation*}
The other two equations \eqref{eq:BccscalarAxiomTwo} and \eqref{eq:WcccoscalarAxiomTwo} are proven symmetrically. 

\section{Proofs of Chapter~4}\label{app:appendixStream}

\begin{proof}[Proof of Lemma \ref{lemma:pullbackpreservation}]
Recall from \S\ref{sec:colimitsMat} how pullbacks in $\Matpoly$ are calculated: the diagram below left is a pullback in $\Matpoly$ if and only if, in the diagram in $\RMod{\poly}$ below right, we have
\[ C= \KerAB \poi \pi_1 \quad \text{and} \quad D= \KerAB \poi \pi_2\]
 where $\KerAB$ is the kernel of the map $( A \mid \minusB ) \colon \poly^{n} \tns \poly^{m}\to \poly^z$.
\[
\xymatrix@C=25pt@R=30pt{
&\ar[dl]_{C} r \ar[dr]^{D} & \\
n \ar[dr]_{A}  & & m \ar[dl]^{B} \\
& z  & }
\quad \quad
\xymatrix@C=40pt@R=25pt{
&\ar[dl]_{C} \poly^r \ar[d]|{\KerAB} \ar[dr]^{D} & \\
\poly^n \ar[dr]_{A}  & \poly^{n} \tns \poly^{m} \ar[l]_{\pi_1} \ar[r]^{\pi_2}  & \poly^m   \ar[dl]^{B} \\
& \poly^z &}
\]
The same holds for pullbacks in $\Matfps$, since $\fps$ is also a PID. Therefore, our proof reduces to check that, for arbitrary $M$ in $\Matpoly$, $\Ker{\wstream{M}}= \wstream{\Ker{M}}$.


Now, given a matrix $M$ in $\Matpoly$, its kernel can be obtained following the recipe of Proposition~\ref{prop:kernelfromHNF}: given an invertible matrix $U$ such that $MU$ is in Hermite Normal Form (HNF), the initial $0$-columns of $U$ yield $\Ker{M}$.

By definition, in order to check that a matrix is in HNF, it suffices to verify the position of the $0$-entries. The embedding $\wstream{\cdot}$ preserves $0$: therefore, since $MU$ is in HNF then also $\wstream{M}\wstream{U}$ must be in HNF. To conclude, let $\vv_1,\dots,\vv_r$ be the initial $0$-columns of $U$, yielding $\Ker{M}$. Since $\wstream{M}\wstream{U}$ is in HNF, by Proposition~\ref{prop:kernelfromHNF} the same vectors $\vv_1,\dots,\vv_r$, now considered as the first $r$ columns of $\wstream{U}$, yield the matrix $\Ker{\wstream{M}}$. Therefore  $\wstream{\Ker{M}} = \Ker{\wstream{M}}$. \end{proof}

For the next proof, it is useful to first fix some notation. The embeddings between $\poly$, $\fps$ $\frpoly$ and $\laur$, defined in \S \ref{sec:stream}, lift to the faithful morphisms of the corresponding PROPs of matrices, as summarised below.
%
%
\begin{equation}\label{eq:squarering}
\vcenter{
\xymatrix@R=0.3cm@C=0.3cm{
\Matfps \ar@{^{(}->}[rr]^{\embfpsfls} & & \Matlaur\\
\\
\Matpoly \ar@{^{(}->}[rr]_{\embpolyfrpoly} \ar@{^{(}->}[uu]^{\embpolyfps} && \Matfrpoly \ar@{_{(}->}[uu]_{\embfrpolyfls}
}
}
\end{equation}

\begin{proof}[Proof of Lemma \ref{lemma:concreteDescrUnivPropFLS}]
By Lemma~\ref{lemma:phipsifull}, for every $H\in \SVpoly[n,m]$ there exists a span $n\tl{V}k\tr{W}m$ in $\Matpoly$ such that
$\Phi(n\tl{V}k\tr{W}m)=H$, i.e., $$H=\{\,(\uu,\vv)\ |\ \uu\in \frpoly^n,\, \vv\in \frpoly^m,\, \exists \ww\in \frpoly^k.\; \embpolyfrpoly(V)\ww=\uu \wedge \embpolyfrpoly(W)\ww = \vv \,\}.$$
For $1\leq i \leq k$, let $\vv_i\in \poly^n$ and $\ww_i\in \poly^m$ be the $i$-th column vectors of $V$ and $W$, respectively.
Then, $\{(\embpolyfrpoly(\vv_i),\embpolyfrpoly(\ww_i)) \ | \ 1\leq i \leq k \}$ spans $H$.

Since $[\embfrpolyfls]$ makes the rightmost front face of \eqref{eq:cube2} commute, it maps $H$ into $\Phi'\circ \Theta(H)$ that is
$$\{\,(\uu,\vv)\ |\ \uu\in \laur^n,\, \vv\in \laur^m,\, \exists \ww\in \laur^k.\; \embfpsfls(\wstream{V})\ww=\uu \wedge \embfpsfls(\wstream{W})\ww = \vv \,\}$$ which, by \eqref{eq:squarering}, is
$$\{\,(\uu,\vv)\ |\ \uu\in \laur^n,\, \vv\in \laur^m,\, \exists \ww\in \laur^k.\; \wlrn{\embpolyfrpoly(V)}\ww=\uu \wedge \wlrn{\embpolyfrpoly(W)}\ww = \vv \,\}\text{.}$$
This space is spanned by $\{(\wlrn{\embpolyfrpoly(\vv_i)},\wlrn{\embpolyfrpoly(\ww_i)}) \ | \ 1\leq i \leq k \}$. This set is obtained by embedding via $\embfrpolyfls$ the generators of $H$ into $\laur$. Therefore $[\embfrpolyfls]$ applied to $H$ yields $[\wlrn{H}]$.\end{proof}

\begin{proof}[Proof of Proposition~\ref{prop:circuittomatrix}]
The proof is by structural induction on $c$, following the inductive definition of $\SFGform$ in \S\ref{sec:SFcalculus}. If $c$ is in $\FC$ then
$\dsem{c}=[({\ee}_i, A{\ee}_i)]_{i\leq n}$ for some matrix $A\in \Matpoly [n,m]$
and clearly any (ordinary) polynomial is rational.

Inductively, suppose that $\dsem{c\: n+1\to m+1}$ is $[({\ee}_i, A{\ee}_i)]_{i\leq n+1}$ for some $A\in {\Matratio [n+1,m+1]}$. We need to show that $\Tr{}(c)$ is
$[({\ee}_i, A'{\ee}_i)]_{i\leq n}$ for some $A'\in \Matratio [n,m]$.

For this purpose, suppose that $\vv= \left(%
			 {\scriptsize   \begin{array}{c}
				\!\! \sigma_1 \!\!\\
				 \!\!\vdots \!\!\\
				 \!\!\sigma_{n+1}\!\!
			   \end{array} }\right)$
and $\ww = \left(%
			  {\scriptsize  \begin{array}{c}
				 \!\!\tau_1 \!\!\\
				 \!\!\vdots \!\!\\
				\!\! \tau_{m+1}\!\!
			   \end{array} }\right)$ are $\frpoly$-vectors such that $A \vv=\ww$.
This means that
\[
\begin{array}{rcl}
 \tau_1 & = & A_{1,1} \sigma_1 + A_{1,2} \sigma_2 + \dots + A_{1,{n+1}} \sigma_{n+1}\\
 \vdots \\
 \tau_i & = & A_{i,1} \sigma_1 + A_{i,2} \sigma_2 + \dots + A_{i,{n+1}} \sigma_{n+1}\\
\vdots \\
 \tau_m & = & A_{m,1} \sigma_1 + A_{m,2} \sigma_2 + \dots + A_{{m+1},{n+1}} \sigma_{n+1}\\
 \end{array}
\]
The semantics of $\Tr{}(c)$ is the subspace corresponding to the solution of the above system of equations plus
$$\sigma_1 = x \cdot \tau_1\text{.}$$
By replacing $\sigma_1$ with $x \cdot \tau_1$ in the first equation, one can deduce that $\tau_1(1 - A_{1,1} \cdot x) = \sum_{j=2}^{n+1}A_{1,j}\sigma_j$.
Note that $1 - A_{1,1} \cdot x \neq 0$ since, by assumption, $A_{1,1} \neq \frac{1}{x}$. Therefore we can safely conclude that
 $$\tau_1 = \sum_{j=2}^{n+1}\left( \frac{A_{1,j}}{1 - A_{1,1} \cdot x} \right)\sigma_j$$
We can now replace $\sigma_1$ by $x \cdot \sum_{j=2}^{n+1}\left( \frac{A_{1,j}}{1 - A_{1,1} \cdot x} \right)\sigma_j$ in the above system of equations and obtain
$$\tau_i = A_{i,1} x \cdot \sum_{j=2}^{n+1}\left( \frac{A_{1,j}}{1 - A_{1,1} \cdot x} \right)\sigma_j + \sum_{j=2}^{n+1}A_{i,j}\sigma_j$$
for all $2 \leq i \leq m+1$.
We thus have $m$ equations with $n$ variables (namely $\sigma_j$ for $2 \leq j \leq n+1$). These form a matrix $A'$ with $m$ columns and $n$ rows. In order to conclude, we have to show that all the entries of this matrix are rationals.

Since $A_{1,1}$ is a rational we can write it as $\frac{p}{k+q \cdot x}$ for some polynomials $p,q$ and scalar $k\neq 0$. So  $1 - A_{1,1} \cdot x = \frac{k + (q-p) \cdot x}{k+ q \cdot x}$ and $\frac{1}{1 - A_{1,1} \cdot x}= \frac{k+q \cdot x}{k+(q-p)\cdot x}$ which is a rational since $k\neq 0$. Since rationals form a ring, i.e., they are closed under $+$ and $\cdot$, all the entries of $A'$ are rationals.

The remaining inductive cases are the ones in which $c = c_1 \poi c_2$ and $c = c_1 \tns c_2$ for circuits $c_1, c_2$ of $\SFGform$. The statement is easily verified by functoriality of $\dsem{\cdot}$ and definition of $\tns$ and $\poi$ in $\SVpoly$.
\end{proof} \label{appendix:SFG}

\begin{proof}[Proof of Proposition \ref{prop:OsemFunctor}] It is routine to check that $\osemO$ preserves identities and the symmetric monoidal structure. Thus we focus on showing that $\osemO$ preserves composition. To keep notation simple, we present the argument for $1$-to-$1$ circuits: the general case does not present any further challenge. Let $c_1 \in \CD[1,1]$ and $c_2 \in \CD[1,1]$ be circuits and fix fps $\alpha = k_0k_1\dots$ and $\beta = l_0l_1\dots$. It is immediate by definition that the following statements are equivalent.

\begin{enumerate}
\item $(\alpha,\beta)$ is in $\itr{c_1}\poi\itr{c_2}$.
\item There exists a fps $\gamma = r_0r_1\dots$ such that $(\alpha,\gamma)$ is in $\itr{c_1}$ and $(\gamma,\beta)$ is in $\itr{c_2}$.
\item For each $i \in \N$, there are $r_i \in \field$, $c_1$-states $s_1$, $t_1$ and $c_2$-states $s_2$, $t_2$ such that $s_1\dtrans{k_i}{r_i} t_1$ and $s_2 \dtrans{r_i}{l_i} t_2$.
\item For each $i \in \N$, there are $c_1\poi c_2$-states $s_1\poi s_2$ and $t_1 \poi t_2$ such that $s_1\poi s_2 \dtrans{k_i}{l_i} t_1 \poi t_2$.
\item $(\alpha,\beta)$ is in $\itr{c_1\poi c_2}$.
\end{enumerate}
This shows that $\itr{c_1}\poi \itr{c_2} = \itr{c_1 \poi c_2}$. In order to show that $\ftr{c_1}\poi \ftr{c_2} = \ftr{c_1 \poi c_2}$ we can use an analogous chain of equivalences, where $\alpha$, $\beta$ and $\gamma$ will be finite instead of infinite traces. We can thereby conclude that $\osem{c_1}\poi \osem{c_2} = \osem{c_1\poi c_2}$.
\end{proof}

For the proof of Theorem \ref{thm:deadinit}, it is useful to fix the following two lemmas.

\begin{lemma}\label{lemma:initallowed} For $c \in \CD[n,m]$, let $s$ be its initial state. Then $s \dtrans{\zerov}{\zerov} s$.
\end{lemma}
\begin{proof} The statement is easily verified by induction on $c$. \end{proof}

\begin{lemma}\label{lemma:SatOsemFunctor} $\phantom{bla}$
\begin{itemize}
\item Given $c_1 \in \CD[n,z]$ and $c_2 \in \CD[z,m]$, $\Sat\osem{c_1\poi c_2} = \Sat\osem{c_1} \poi \Sat\osem{c_2}$.
 \item Given $c_1 \in \CD[n,m]$ and $c_2 \in\CD[z,r]$, $\Sat\osem{c_1\tns c_2} = \Sat\osem{c_1} \tns \Sat\osem{c_2}$.
 \end{itemize}
    \end{lemma}
\begin{proof} Let us focus on the first statement. Analogously to Proposition~\ref{prop:OsemFunctor}, we confine ourselves to presenting the argument for the $1$-to-$1$ case, that is, we take $c_1$ and $c_2$ both in $\CD[1,1]$. Fix notation $\pair{f_1}{g_1} \df \osem{c_1} \in \Relpoly \times \Relfps[1,1]$ and $\pair{f_2}{g_2} \df \osem{c_2} \in \Relpoly \times \Relfps[1,1]$ and suppose that $(\sigma, \tau)$ is in $\Sat\big(\pair{f_1}{g_1} \poi \pair{f_2}{g_2}\big)$. This means that there is $(\alpha,\beta) \in g_1 \poi g_2$ and $z \in \Z$ generating $(\sigma, \tau)$. It follows that there is a fps $\gamma$ such that $(\alpha,\gamma)$ is in $g_1$ and $(\gamma,\beta)$ is in $g_2$. Define the fls $\rho$ by $\rho(i + z) \df \gamma(i)$ if $i \geq 0$ and $\rho(i+z) = 0$ otherwise. Then $(\alpha,\gamma)$ with $z$ generates $(\sigma,\rho)$ and $(\gamma,\beta)$ with $z$ generates $(\rho,\tau)$. By definition this means that $(\sigma,\rho)$ is in $\Sat\pair{f_1}{g_1}$ and $(\rho,\tau)$ is in $\Sat\pair{f_2}{g_2}$. Therefore $(\sigma,\tau)$ is in $\Sat\pair{f_1}{g_1} \poi \Sat\pair{f_2}{g_2}$.

For the converse direction, suppose that $(\sigma,\tau)$ is an element of $\Sat\pair{f_1}{g_1} \poi \Sat\pair{f_2}{g_2}$. Then there is a fls $\rho$ such that $(\sigma,\rho)$ is in $\Sat\pair{f_1}{g_1}$ and $(\rho,\tau)$ is in $\Sat\pair{f_2}{g_2}$. By definition, this means that there are $(\alpha,\gamma) \in g_1$ and $z \in \Z$ generating $(\sigma,\rho)$, and $(\gamma',\beta) \in g_2$ and $z' \in \Z$ generating $(\rho,\tau)$. We now distinguish two cases.
\begin{itemize}
\item Suppose first that $z' \leq z$. For $j < (z \mn z')$, $\gamma'(j) =\rho(j+z') = 0$ because $j + z' < (z\mn z')+z' = z$ and $z$ is assumed to be smaller or equal than the degree of $\rho$. For $j \geq (z\mn z')$, $\gamma'(j) = \rho(j+z') = \rho(j \mn (z\mn z') + z) = \gamma ( j \mn (z\mn z'))$. Intuitively, this means that $\gamma'$ is given by prefixing $\gamma$ with $z\mn z'$ elements with value $0$. In other words, $\gamma' = \gamma \cdot x^{z\mn z'}$. Now, define the fps $\alpha'$ by $\alpha'(j) = \alpha(j \mn (z\mn z'))$ for $j \geq (z\mn z')$ and $\alpha(j) = 0$ otherwise. By construction, $\alpha' = \alpha \cdot x^{z\mn z'}$. 
    We then want to show that $(\alpha' ,\gamma' )$ is in $g_2$. For this purpose, let $s_0 \dtrans{k_0}{l_0} s_1 \dtrans{k_1}{l_0} \dots$ be the computation of $c$ from initial state $s_0$ yielding the trace $(\alpha,\gamma) \in \itr{c}$. As observed in Lemma~\ref{lemma:initallowed}, a transition $s_0 \dtrans{0}{0} s_0$ is always possible. Therefore the following is a computation of $c$, where $s_j = s_0$ for all $j \leq z'\mn z$.
   $$s_0 \dtrans{0}{0} s_1 \dtrans{0}{0} \dots \dtrans{0}{0} s_{(z'\mn z)} \dtrans{k_0}{l_0} s_{(z'\mn z)+1} \dtrans{k_1}{l_1} \dots$$
   By definition of $\alpha'$ and $\gamma'$, such computation yields the trace $(\alpha' ,\gamma' )$ as an element of $g_2$. Since by assumption $(\gamma',\beta)$ is in $g_2$, it follows that $(\alpha',\beta)$ is in $g_1 \poi g_2$. Also notice that $\alpha'(j) = \alpha(j \mn (z\mn z')) = \sigma(j \mn (z\mn z') + z) = \sigma(j + z')$ for $j \geq (z\mn z')$, and otherwise, for $j \ls (z \mn z')$, $\alpha'(j) = 0 = \sigma(j+z')$ because $j+z' \ls (z \mn z')+z' = z$ and $z$ is smaller or equal to the degree of $\sigma$. Thus $(\alpha',\beta)$ and $z'$ generate $(\sigma,\tau)$. It follows that $(\sigma,\tau)$ is in $\Sat(\pair{f_1}{g_1}\poi \pair{f_2}{g_2})$.
\item The case in which $z < z'$ is handled symmetrically: instead of constructing $\alpha'$, we let $\beta'$ be $\beta \cdot x^{z'\mn z}$. Then one can check that $(\alpha,\beta')$ is in $g_1 \poi g_2$ and, with the choice of the instant $z$, it generates $(\sigma,\tau)$. It follows that $(\sigma,\tau)$ is in $\Sat(\pair{f_1}{g_1}\poi \pair{f_2}{g_2})$.
  \end{itemize}
  The proof of the second statement in the lemma is conceptually straightforward, following a similar line of reasoning as the argument for the first statement.
\end{proof}

\noindent We split the proof of Theorem \ref{thm:deadinit} in three parts.

\begin{proof}[Proof of Theorem~\ref{thm:deadinit}]~\phantom{bla}

\begin{proof}[Proof of point \ref{point:Posem=Dsem}] The proof is by induction on $c$. We just show two representative base cases:
\begin{itemize}
\item if $c = \Wmult$, then let $\sigma$ and $\tau$ be fls, say with degree $a$ and $b$ respectively. Pick $z$ smaller than $a$ and $b$ and let fps $\alpha$ and $\beta$ be defined by $\alpha(i) \df \sigma(i+z)$ and $\beta(i) \df \tau(i+z)$ respectively, for $i \in \N$. Then the pair $(\left(\begin{array}{c}
				 \!\!\! \sigma\!\!\! \\
				 \!\!\! \tau\!\!\!
				\end{array}\right),\sigma+\tau  )$, which is in $\strsem{c}$, is generated by $(\left(\begin{array}{c}
				 \!\!\! \alpha\!\!\! \\
				 \!\!\! \beta\!\!\!
				\end{array}\right), \alpha+\beta )$ and $z$. In particular, observe that by definition $\alpha+\beta(i) = \sigma+\tau(i+z)$. Since $(\left(\begin{array}{c}
				 \!\!\! \alpha\!\!\! \\
				 \!\!\! \beta\!\!\!
				\end{array}\right), \alpha+\beta )$ is in $\itr{c}$ then $(\left(\begin{array}{c}
				 \!\!\! \sigma\!\!\! \\
				 \!\!\! \tau\!\!\!
				\end{array}\right),\sigma+\tau  )$ is in $\Sat\osem{c}$. This show that $\strsem{c} \subseteq \Sat\osem{c}$. Conversely, suppose that $\rho$, $\phi_1$ and $\phi_2$ are fls such that  $(\left(\begin{array}{c}
				 \!\!\! \phi_1\!\!\! \\
				 \!\!\! \phi_2\!\!\!
				\end{array}\right),\rho  )$ is in $\Sat\osem{c}$. Then there are fps $\gamma$, $\delta_1$, $\delta_2$ and $z \in \Z$ such that $(\left(\begin{array}{c}
				 \!\!\! \delta_1\!\!\! \\
				 \!\!\! \delta_2\!\!\!
				\end{array}\right),\gamma  )$ is in $\osem{c}$ and together with $z$ generates $(\left(\begin{array}{c}
				 \!\!\! \phi_1\!\!\! \\
				 \!\!\! \phi_2\!\!\!
				\end{array}\right),\rho  )$. This means in particular that $\gamma = \delta_1 + \delta_2$. Also, for $i \geq 0$, $\rho(i+z) = \gamma(i) = \delta_1(i) + \delta_2(i) = \phi_1(i+z) + \phi_2(i+z)$. For $i < 0$, $\rho(i+z) = 0 = \phi_1(i+z) = \phi_2(i+z)$, because $z$ is smaller or equal to the minimum among the degrees of $\rho$, $\phi_1$ and $\phi_2$. This shows that  $\rho = \phi_1+\phi_2$ and therefore $(\left(\begin{array}{c}
				 \!\!\! \phi_1\!\!\! \\
				 \!\!\! \phi_2\!\!\!
				\end{array}\right),\rho  )$ is in $\strsem{c}$.
\item For $c = \circuitXT$, let $\sigma$ be a fls with degree $d$. Any element of $\strsem{c}$ has shape $(\sigma,\sigma \cdot x)$. Define the fps $\alpha$ by putting $\alpha(i) \df \sigma(i+d)$. Clearly, $(\alpha , \alpha \cdot x)$ together with the choice of $d$ for $z \in \Z$ generate $(\sigma,\sigma \cdot x)$. Also, by definition $(\alpha , \alpha \cdot x)$ is in $\osem{c}$, meaning that $(\sigma,\sigma \cdot x)$ is in $\Sat\osem{c}$. Therefore $\strsem{c} \subseteq \Sat\osem{c}$. For the converse direction, let $(\sigma,\tau)$ be in $\Sat\osem{c}$. This means that it is generated by some $(\alpha,\beta) \in \osem{c}$ and $z \in \Z$. By definition of $\osem{c}$, $\alpha$ and $\beta$ are of shape $k_0k_1k_2\dots$ and $0k_0k_1\dots$ respectively, that is, $\beta = \alpha \cdot x$. For $i \geq 0$ we have that $\sigma(i+z) = \alpha(i) = \beta(i+1) = \tau(i+1+z)$. For $i < 0$ we have that $\sigma(i+z) = 0 = \tau(i+1+z)$. To see this last point, observe that $z$ is by assumption smaller or equal to the degree of $\sigma$ and smaller than the degree of $\tau$, because $\tau(z) = \beta(0) = 0$. We can thereby conclude that $\tau = \sigma \cdot x$ and thus $(\sigma,\tau)$ is in $\strsem{c}$. This shows that $\Sat\osem{c} \subseteq \strsem{c}$.
\end{itemize}

The inductive cases of composition by $\tns$ and $\poi$ are just given by application of the inductive hypothesis, Lemma~\ref{lemma:SatOsemFunctor} and functoriality of $\strsemO$. For instance, for $c = c_1 \poi c_2$ we have:
\[\Sat\osem{c_1 \poi c_2} = \Sat\osem{c_1} \poi \Sat\osem{c_2} \eql{IH} \strsem{c_1} \poi \strsem{c_2} = \strsem{c_1\poi c_2} .\]
\end{proof}

\begin{proof}[Proof of point \ref{point:deadIncl}] Let $\pair{f}{g} = \Forget\Sat\osem{c}$. For the infinite traces, we now show that $\itr{c} \subseteq g$. Suppose that $(\vlist{\alpha},\vlist{\beta})$ is an element of $\itr{c}$. For each fps $\alpha_j$ in $\vlist{\alpha}$ and $\beta_k$ in $\vlist{\beta}$, define fls $\sigma_j$ and $\tau_k$ as $\sigma_j(i) = \alpha(i)$ and $\tau_k(i) = \beta_k(i)$ for $i \geq 0$, and as $0$ for $i <0$. This gives vectors $\vlist{\sigma}$ and $\vlist{\tau}$ such that $(\vlist{\sigma},\vlist{\tau})$ is generated by $(\vlist{\alpha},\vlist{\beta})$ for $z = 0$. It follows that $(\vlist{\sigma},\vlist{\tau}) \in \Sat\osem{c}$ and $(\vlist{\alpha},\vlist{\beta}) \in g$.

For the finite traces, we need to show that $\ftr{c} \subseteq f$. Because $c$ is deadlock free, any trace $(\vlist{\alpha},\vlist{\beta})$ in $\ftr{c}$ is the prefix of some infinite trace in $\itr{c}$. Since we just proved that $\itr{c} \subseteq g$, then $(\vlist{\alpha},\vlist{\beta})$ is also the prefix of some infinite trace in $g$. It follows by definition of $\Forget$ that $(\vlist{\alpha},\vlist{\beta})$ is an element of $f$. We can thereby conclude that $\osem{c} \subseteq \Forget\Sat\osem{c}$.
\end{proof}

\begin{proof}[Proof of point \ref{point:initIncl}] For the sake of readability, we show the proof only for $1$-to-$1$ circuits $c$ of $\CD$. The general case is just a more involved formulation of the same argument. For this purpose, let $\pair{f}{g} = \Forget\Sat\osem{c}$ and suppose that $(\alpha,\beta)$ is an element of $g$, which we want to show to be in $\itr{c}$. By definition of $\Forget$ there are fls $\sigma$, $\tau$ and $z\in \Z$ such that $(\alpha,\beta)$ generates $(\sigma,\tau)$ and $(\sigma,\tau)$ is in $\Sat\osem{c}$. Then, by definition of $\Sat$, there are also $\alpha'$, $\beta'$ and $z'\in \Z$ such that $(\alpha',\beta')$ is in $\itr{c}$ and generates $(\sigma,\tau)$. From this, we shall derive that also $(\alpha,\beta)$ is in $\itr{c}$: the idea is that $(\alpha',\beta')$ and $(\alpha,\beta)$ are the same, modulo a prefix of $0$s. We now distinguish two cases.
\begin{itemize}
\item First, suppose that $z' \leq z$. Intuitively, this is the case in which $(\alpha',\beta')$ is $(\alpha,\beta)$ with the addition of a prefix of $0$s at the beginning of the two streams. Indeed, for $j < (z \mn z')$, $\alpha'(j) = \sigma(j+z') = 0$ because $j + z' < (z\mn z')+z' = z$ and $z$ is assumed to be smaller or equal to the degree of $\sigma$. Also, for $j \geq (z\mn z')$, $\alpha'(j) = \sigma(j+z') = \sigma(j \mn (z\mn z') + z) = \alpha ( j \mn  (z\mn z'))$. This shows that, given notation $k_0 k_1 \dots$ for $\alpha$, we can write  $\alpha'$ as $0 \dots 0 k_0 k_1 \dots$, where the prefix $0 \dots 0$ is of length $z \mn z'$. In other words, $\alpha' = \alpha \cdot x^{z \mn z'}$. With an analogous reasoning we can calculate that, given notation $l_0 l_1 \dots$ for $\beta$, then $\beta' = 0 \dots 0 l_0 l_1 \dots$, where the prefix $0 \dots 0$ is of length $z \mn z'$ --- that is, $\beta' = \beta \cdot x^{z \mn z'}$.
 Now, since $(\alpha',\beta')$ is in $\itr{c}$ then there is a computation of $c$ as follows from the initial $c$-state $s_0$.
 $$s_0 \dtrans{0}{0} s_1 \dtrans{0}{0} \dots s_{(z\mn z')} \dtrans{k_0}{l_0} s_{(z\mn z')+1} \dtrans{k_1}{l_1} \dots$$
  Since $c$ is initialisation-free, $s_0 = s_1 = \dots = s_{(z\mn z')}$. It follows that $s_{(z\mn z')} \dtrans{k_0}{l_0} s_{(z\mn z')+1} \dtrans{k_1}{l_1} \dots$ is also a computation of $c$ from the initial state, meaning that $(k_0k_1\dots,l_0l_1\dots) = (\alpha,\beta)$ is in $\itr{c}$.
\item The remaining case to consider is the one in which $z' > z$. We can then use a symmetric argument to conclude that $\alpha = \alpha' \cdot x^{z'-z}$ and $\beta = \beta' \cdot x^{z'-z}$. Now, let $s_0 \dtrans{k_0}{l_0} s_1 \dtrans{k_1}{l_1} \dots$ be the computation of $c$ from initial state $s_0$ yielding the trace $(\alpha',\beta') \in \itr{c}$. As observed in Lemma~\ref{lemma:initallowed}, a transition $s_0 \dtrans{0}{0} s_0$ is always possible. Therefore the following is a computation of $c$, where $s_j = s_0$ for all $j \leq z'\mn z$.
   $$s_0 \dtrans{0}{0} s_1 \dtrans{0}{0} \dots \dtrans{0}{0} s_{(z'\mn z)} \dtrans{k_0}{l_0} s_{(z'\mn z)+1} \dtrans{k_1}{l_1} \dots$$
   yielding the trace $(\alpha,\beta)$ as an element of $\itr{c}$.
\end{itemize}
The argument above show that $g \subseteq \itr{c}$. It remains to show that $f \subseteq \ftr{c}$. By definition, a finite trace $(\alpha,\beta)$ in $f$ is always a prefix of an infinite trace in $g$. Since $g \subseteq \itr{c}$ then $(\alpha,\beta)$ is also a prefix of an infinite trace in $\itr{c}$, meaning that it is an element of $\ftr{c}$. This concludes the proof that $\Forget\Sat\osem{c} \subseteq \osem{c}$.
\end{proof}
The proof of point \ref{point:initIncl} concludes the proof of Theorem~\ref{thm:deadinit}. \end{proof}

\begin{proof}[Proof of Lemma~\ref{lemma:rowequivalent}]
We show a procedure similar to Gaussian elimination that, using elementary row operations, transforms $n\times m$ matrices to rational form.

First, we set all the entries of the first row to be polynomials $p_1, \dots, p_m$ (simply by multiplying this row by the product of all denominators). Like for formal Laurent series, we define the \emph{degree} of a polynomial $k_0+k_1x+\dots+k_zx^z$ to be the smallest $k_i\neq 0$; for instance $1+x$ has degree $0$ and $x+x^2$ has degree 1. Amongst $p_1,\dots, p_m$, we pick $p_{\kappa_1}$ with minimal degree and we multiply the first row by $\frac{1}{p_{\kappa_1}}$. In the resulting row, all the entries are rationals, since they are fractions $\frac{p_i}{p_{\kappa_1}}$ where the denominator has degree smaller or equal than the numerator. Moreover in the $\kappa_1$-th position there is $1$. We call $\kappa_1$ the \emph{pivot} of the first row and this sub-procedure \emph{rationalization} of a row.

Second, we bring to $0$ all the entries below the first pivot. Like in Gaussian elimination, this can be done by simply adding to each row a scalar multiple of the first one. This second sub-procedure is the \emph{downward substitution} of a pivot.

Rationalization and downward substitution can be iteratively applied to all the (non-zero) rows in the matrix so to obtain a novel matrix where (a) all the entries are rationals, (b) each (non-zero) row has a pivot with coefficient $1$ and (c) all the entries \emph{below} a pivot are $0$. For instance, one obtains a matrix as the following.
$$ \left(%
{\scriptsize \begin{array}{cccc}
  r_1   & r_2 & 1 & r_3 \\
  r_4   & 1   & 0 & r_5\\
  r_6   & 0   & 0 & 1\\
  0     & 0   & 0 & 0
\end{array}}\right)$$
For having a matrix in rational form, we only need to set to $0$ all the entries \emph{above} a pivot. We start from the last (non-zero) row, which we call $s$ with pivot $\kappa_s$. Like for downward substitution, we can add to each row above $s$ a scalar multiple of $s$, but we have to do it carefully, by checking that the resulting rows are in the good shape. Take a (non-zero) row $j$ above $s$, and call $r_{j\kappa_s}$ the $\kappa_s$-th entry of such row. By virtue of (a), $r_{j\kappa_s}$ is a rational. By adding to the row $j$, the row $s$ multiplied by $-r_{j\kappa_s}$ we obtain a new row where (d) the $\kappa_s$-element is 0; (e) the entry at the pivot $\kappa_j$ is 1 since, by (c), in $s$ the $\kappa_j$-th element is $0$; (f)
all the entries of the row are rationals, since they are obtained by additions and multiplications of rationals (and rationals form a ring).
We can repeat this for all the pivots $\kappa_u$ and for all the rows above $u$ and we will eventually obtain a matrix
where by (f) all the entries are rationals, each row has, by (e), a pivot with entry $1$ and all the entries above and below a pivot are $0$ by (c) and (d).
\end{proof}

\chapter{Works not Included in the Thesis}

During the PhD I also worked on other research topics, which are not detailed in this thesis for reasons of homogeneity. In the following I briefly report on these works and the associated publications.

\paragraph{\S 1 Expressive completeness for fragments of the modal $\mu$-calculus} The modal $\mu$-calculus (MC) is a well-balanced specification language for describing and verifying properties of transition systems. A landmark result in the area is the Janin-Walukiewicz theorem, stating that MC is the bisimulation-invariant fragment of monadic second-order logic (MSO). That means, MC is as expressive as MSO on properties not distinguishing between bisimilar transition systems, like computationally relevant properties usually do.

The Janin-Walukiewicz theorem has been extended also to finite, transitive and tree models. However, very little is known about analogous expressive completeness results for fragments of MC. In my master thesis I answered this question for the alternation-free fragment of the $\mu$-calculus (AFMC), which is particularly appealing in applications for its computational feasibility. By means of an automata theoretic construction, I proved that AFMC is the bisimulation-invariant fragment of well-founded MSO: the latter is the fragment of MSO where, roughly speaking, second-order quantifiers are only allowed to range over well-founded portions of the model. During my PhD, together with my master thesis supervisors Y. Venema and A. Facchini, we prepared a joint paper on these results which appeared at LICS'13~\cite{FacchiniVZ13}. Our collaboration continued in another paper appearing at CSL-LICS'14~\cite{CarreiroFVZ14}, also with F. Carreiro, in which we tailored a fragment of MC where the application of fixpoint operators is restricted by a continuity condition and prove that it is the bisimulation-inviariant fragment of weak MSO.
\begin{itemize}[itemsep=1pt,topsep=1pt,parsep=1pt,partopsep=1pt]
\item F. Carreiro, A. Facchini, Y. Venema, F. Zanasi - \emph{Weak MSO: automata and expressiveness modulo bisimilarity}, CSL-LICS'14.
\item A. Facchini, Y. Venema, F. Zanasi - \emph{A Characterization Theorem for the Alternation-Free Fragment of the Modal $\mu$-Calculus}, LICS'13.
\end{itemize}

\paragraph{\S 2 Bialgebraic semantics for logic programming} Recent work by E. Komendantskaya and J. Power~\cite{KomPowerCSL11} introduces a coalgebraic treatment of logic programming. This approach brings several benefits (for instance, allows for reasoning about coinductive programs) but fails to be compositional: intuitively, their semantics is not compatible with substitutions that may occur in formulas. In categorical terms, the space of formulas is represented by a presheaf \emph{F} on a Lawvere theory; a logic program and its semantics should then be coalgebras on~\emph{F}, but they fail to be natural transformations, forcing the authors to introduce laxness in the picture.
The drawbacks of the coalgebraic approach were the main motivation for F. Bonchi and me to propose a bialgebraic environment for logic programming. We model logic programs as bialgebras on presheaves.
We obtain compositionality with respect to substitutions by applying \emph{saturation techniques}, previously employed in the setting of nominal process calculi. Categorically, saturation amounts to a certain right Kan extension of the presheaf modeling formulas. The bialgebra structure allows us to achieve also a second form of compositionality, with respect to the conjunction of formulas in a goal. This work appeared in a journal article~\cite{BonchiZanasiLPJournal}, which extends a conference paper~\cite{BonchiZ_Calco13} winner of the Best Paper award at CALCO'13.


\begin{itemize}[itemsep=1pt,topsep=1pt,parsep=1pt,partopsep=1pt]
\item F. Bonchi, F. Zanasi - \emph{Bialgebraic Semantics for Logic Programming}, Logical Methods in Computer Science, Vol. 11, Issue 2, 2015.
\item F. Bonchi, F. Zanasi - \emph{Saturated Semantics for Coalgebraic Logic Programming}, CALCO'13 (Best Paper Award).
\end{itemize}
\paragraph{\S 3 A Categorical environment for systems with internal moves} Coalgebras has been successful in modeling a wide range of state-based systems. However, they seem to be inadequate to describe the behaviour of systems where internal transitions play a role --- notable examples include automata with $\epsilon$-transitions and labelled transition systems with $\tau$-transitions. In a collaboration with F. Bonchi, A. Silva and S. Milius, we had the intuition that the gap could be filled by considering a more general form of coinduction --- usually called \emph{parametric corecursion}~\cite{AdamekMV02} --- and the attached body of work on iteration theories and Elgot monads. We developed this idea into a framework where the type of a system with internal moves is an Elgot monad satisfying an initial algebra-final coalgebra coincidence. Our work resulted in a paper which I presented at CMCS'14~\cite{Bonchi-CMCS14} and a subsequent journal publication~\cite{BMSZ-daggerJournal}.
\begin{itemize}[itemsep=1pt,topsep=1pt,parsep=1pt,partopsep=1pt]
\item F. Bonchi, S. Milius, A. Silva, F. Zanasi - \emph{Killing Epsilons with a Dagger - A Coalgebraic Study of Systems with Algebraic Label Structure}, Theoretical Computer Science, accepted for publication, 2015.
\item F. Bonchi, S. Milius, A. Silva, F. Zanasi - \emph{How to Kill Epsilons with a Dagger - A Coalgebraic Take on Systems with Algebraic Label Structure} - \mbox{CMCS'14}.
\end{itemize}    

\printglossary[style = mystyle,title=Glossary of Symbols]
\bibliographystyle{acm}
\bibliography{catBib3}
\end{document}